\RequirePackage{fix-cm}

\documentclass[leqno, fleqn, centertags, 12pt]{article}

\usepackage{latexsym}
\usepackage{amsmath}
\usepackage{amssymb}
\usepackage{verbatim}

\usepackage[T1]{fontenc}

\usepackage[upright, widespace]{fourier}

\usepackage{bm}

\usepackage{fullpage}

\usepackage{tikz-cd}

\tikzcdset{row sep/boom/.initial=5em}
\tikzcdset{column sep/boom/.initial=5.1em}

\tikzcdset{row sep/boomm/.initial=5.6em}
\tikzcdset{column sep/boomm/.initial=4.25em}

\tikzcdset{row sep/sboom/.initial=5.6em}
\tikzcdset{column sep/sboom/.initial=5.712em}

\tikzcdset{row sep/rboom/.initial=5.225em}

\tikzcdset{column sep/sma/.initial=3em}
\tikzcdset{row sep/sma/.initial=4.25em}

\tikzcdset{row sep/tri/.initial=5em}
\tikzcdset{column sep/tri/.initial=2.4em}

\newskip\nineskipamount \nineskipamount=9pt plus 0pt minus 0pt
\newskip\zeroskipamount \zeroskipamount=0pt plus 0pt minus 0pt
\usepackage[noorphans,vskip=\nineskipamount]{quoting}

\newcommand{\dis}{\displaystyle}

\makeatletter
\renewcommand{\@makefntext}[1]{\vspace*{0.5ex}\parindent=0em
\hspace*{-0.4em}
\hbox to 0.4em{\hss\@makefnmark}\hspace*{0.4em}{#1}
}
\makeatother

\newcounter{mysectionnumber}
\newcommand{\mysection}[2]{\setcounter{footnote}{0}
\setcounter{myparnum}{0}
\refstepcounter{mysectionnumber}
\vspace{27pt}{\Large {\themysectionnumber.} {#1}}\label{#2}\vspace*{15pt}}

\newcommand{\myuppar}[1]{\vspace{\medskipamount}\textbf{#1}\hspace*{0.5em}}

\newcommand{\myit}[1]{\textbf{\textit{#1}}\hspace{0.0em}}

\newcounter{myparnum}[mysectionnumber]
\renewcommand{\themyparnum}{\arabic{mysectionnumber}.\arabic{myparnum}}
\newcommand{\mypar}[2]{\refstepcounter{myparnum}{\vspace{\medskipamount}\textbf{{\themyparnum. #1}\label{#2}}\hspace{0.5em}}}

\newcounter{mylemmanum}[myparnum]
\newcommand{\proof}{\vspace{\medskipamount}{\textbf{{\emph{Proof}.}}\hspace*{1em}}}

\newcommand{\eproof}{ $\blacksquare$}

\def\sss{\hspace{0.05em}\ }
\def\dss{\hspace{0.1em}\ }
\def\trs{\hspace{0.15em}\ }
\def\qss{\hspace{0.2em}\ }
\def\pss{\hspace{0.3em}\ }
\def\oss{\hspace{0.4em}\ }

\def\halfff{\hspace*{0.025em}}
\def\fff{\hspace*{0.05em}}
\def\dff{\hspace*{0.1em}}
\def\trf{\hspace*{0.15em}}
\def\qff{\hspace*{0.2em}}
\def\pff{\hspace*{0.3em}}
\def\off{\hspace*{0.4em}}

\newcommand{\hnsp}{\hspace*{-0.05em}}
\newcommand{\nsp}{\hspace*{-0.1em}}
\newcommand{\nnsp}{\hspace*{-0.15em}}
\newcommand{\snsp}{\hspace*{-0.175em}}
\newcommand{\dnsp}{\hspace*{-0.2em}}

\renewcommand{\leq}{\leqslant}
\renewcommand{\geq}{\geqslant}

\newcommand{\zzz}{\mathbf{Z}}

\newcommand{\ccc}{\mathbf{C}}
\newcommand{\rrr}{\mathbf{R}}

\newcommand{\image}{\operatorname{Im}\trf}

\newcommand{\kernel}{\operatorname{Ker}\trf}

\newcommand{\pr}{\operatorname{p{\fff}r}}

\newcommand{\id}{\operatorname{id}}

\newcommand{\gr}{\operatorname{G{\fff}r}}

\newcommand{\fin}{\operatorname{fin}}

\newcommand{\mult}{\operatorname{Mult}\trf}
\newcommand{\fred}{^{\dff \operatorname{Fred}}}

\newcommand{\re}{\operatorname{R{\fff}e}\trf}
\newcommand{\im}{\operatorname{I{\fff}m}\trf}

\newcommand{\ti}{\operatorname{t-ind}}
\newcommand{\ai}{\operatorname{a-ind}}

\def\omult{{%
    \setbox0\hbox{${=\mathrel{\mkern-5mu}=}$}%
    \rlap{\hbox to \wd0{\hss$|\mathrel{\mkern-7.5mu}|$\hss}}\box0
    }}

\newcommand{\num}[1]{|\qff #1 \qff|}

\newcommand{\fclass}[1]{[\snsp [\dff #1\dff]\snsp ]}
\newcommand{\norm}[1]{\|\qff #1 \qff\|}
\newcommand{\bnorm}[1]{\left\|\qff #1 \qff\right\|}
\newcommand{\sco}[1]{\langle\trf #1 \trf\rangle}
\newcommand{\bsco}[1]{\left\langle\trf #1 \trf\right\rangle}
\newcommand{\basco}[1]{\bigl\langle\trf #1 \trf\bigr\rangle}

\newcommand{\sa}{\mathrm{sa}}

\newcommand{\ffin}{^{\dff \mathrm{fin}}}
\newcommand{\comp}{^{\dff \mathrm{comp}}}
\newcommand{\phg}{_{\mathrm{phg}}}

\newcommand{\ttoo}{\hspace*{0.2em}\longrightarrow\hspace*{0.2em}}

\begin{document}

\setlength{\baselineskip}{12pt plus 0pt minus 0pt}
\setlength{\parskip}{12pt plus 0pt minus 0pt}
\setlength{\abovedisplayskip}{12pt plus 0pt minus 0pt}
\setlength{\belowdisplayskip}{12pt plus 0pt minus 0pt}

\newskip\smallskipamount \smallskipamount=3pt plus 0pt minus 0pt
\newskip\medskipamount   \medskipamount  =6pt plus 0pt minus 0pt
\newskip\bigskipamount   \bigskipamount =12pt plus 0pt minus 0pt

\author{Nikolai\qss V.\qss Ivanov}
\title{The\qss index\qss of\qss 
self-adjoint\qss Shapiro-Lopatinskii\qss
boundary\qss problems of\qss order\qss one} 
\date{}

\footnotetext{\hspace*{-0.65em}\copyright\qss 
Nikolai\qss V.\qss Ivanov,\qss 2022.\trs 
The work reported in\sss the first\sss version of\dss this paper was done 
during my\dss Consultantship year at\sss Michigan\sss State\sss University.\qss
I am grateful\dss to\dss Keith\dss Promislow,\oss the Chair of\qss Department\sss of\trs Mathematics,\qss
for offering me this opportunity.\oss\vspace{3pt}}

\footnotetext{\hspace*{-0.65em}I am grateful\sss to\trs G.\dss Nagy,\pss M.\dss Prokhorova,\pss
and\trs J.\dss Schenker\sss for suggestions and stimulating discussions.\oss}

\maketitle

\renewcommand{\baselinestretch}{1}
\selectfont

\vspace*{6ex}


\hspace*{0.44\textwidth}\myit{\hspace*{0em}\large Contents}  \vspace*{3ex} \vspace*{\bigskipamount}\\ 
\hspace*{0.1\textwidth}\hbox to 0.8\textwidth{\myit{\phantom{1}1.}\hspace*{0.5em} Introduction \hfil 2}\hspace*{0.5em} \vspace*{0.25ex}\\
\hspace*{0.1\textwidth}\hbox to 0.8\textwidth{\myit{\phantom{1}2.}\hspace*{0.5em} Geometric algebra at\dss the boundary\hfil 7}  \hspace*{0.5em} \vspace*{0.25ex}\\
\hspace*{0.1\textwidth}\hbox to 0.8\textwidth{\myit{\phantom{1}3.}\hspace*{0.5em} Self-adjoint\sss symbols and\dss boundary conditions \hfil 16} \hspace*{0.5em} \vspace*{0.25ex}\\
\hspace*{0.1\textwidth}\hbox to 0.8\textwidth{\myit{\phantom{1}4.}\hspace*{0.5em} Multiplicative properties of\dss symbols \hfil 24} \hspace*{0.5em} \vspace*{0.25ex}\\
\hspace*{0.1\textwidth}\hbox to 0.8\textwidth{\myit{\phantom{1}5.}\hspace*{0.5em} Abstract\dss boundary\sss problems \hfil 41} \hspace*{0.5em} \vspace*{0.25ex}\\  
\hspace*{0.1\textwidth}\hbox to 0.8\textwidth{\myit{\phantom{1}6.}\hspace*{0.5em} Multiplicative properties of\dss abstract\sss operators \hfil 47} \hspace*{0.5em} \vspace*{0.25ex}\\
\hspace*{0.1\textwidth}\hbox to 0.8\textwidth{\myit{\phantom{1}7.}\hspace*{0.5em} Pseudo-differential\sss operators and\dss 
boundary\sss conditions \hfil 57}  \hspace*{0.5em} \vspace*{0.25ex}\\ 
\hspace*{0.1\textwidth}\hbox to 0.8\textwidth{\myit{\phantom{1}8.}\hspace*{0.5em} Multiplicative properties of\dss pseudo-differential\sss operators \hfil 62}  \hspace*{0.5em} \vspace*{0.25ex}\\
\hspace*{0.1\textwidth}\hbox to 0.8\textwidth{\myit{\phantom{1}9.}\hspace*{0.5em} Glueing and cutting \hfil 70}  \hspace*{0.5em} \vspace*{0.25ex}\\
\hspace*{0.1\textwidth}\hbox to 0.8\textwidth{\myit{10.}\hspace*{0.5em} The\sss index\sss theorem \hfil 77}  \hspace*{0.5em} \vspace*{0.25ex}\\
\hspace*{0.1\textwidth}\hbox to 0.8\textwidth{\myit{11.}\hspace*{0.5em} Self-adjoint\trs Fredholm\dss relations \hfil 81}  \hspace*{0.5em} \vspace*{0.25ex}\\
\hspace*{0.1\textwidth}\hbox to 0.8\textwidth{\myit{12.}\hspace*{0.5em} Boundary\dss triplets and\dss the analytical\dss index \hfil 87}  \hspace*{0.5em} \vspace*{0.25ex}\\
\hspace*{0.1\textwidth}\hbox to 0.8\textwidth{\myit{13.}\hspace*{0.5em} Realizations of\dss boundary\sss symbols \hfil 96}  \hspace*{0.25em} \vspace*{0.25ex}\\
\hspace*{0.1\textwidth}\hbox to 0.8\textwidth{\myit{14.}\hspace*{0.5em} Special\dss boundary conditions \hfil 104} 
\hspace*{0.5em} \vspace*{0.25ex}\\
\hspace*{0.1\textwidth}\hbox to 0.8\textwidth{\myit{15.}\hspace*{0.5em} Dirac-like boundary problems \hfil 112}  \hspace*{0.5em}\vspace*{1.5ex}\\
\hspace*{0.1\textwidth}\hbox to 0.8\textwidth{\myit{References} \hfil 118}\hspace*{0.5em}  \hspace*{0.5em}  \vspace*{0.25ex}

\renewcommand{\baselinestretch}{1}
\selectfont

\newpage
\mysection{Introduction}{introduction}

\myuppar{Elliptic\sss local\dss boundary problems.}
Let\sss $X$\sss be a compact\sss manifold\sss with\sss boundary\sss 
$Y\off =\off \partial\dff X$\dss
and\dss let\sss $E$\sss be a\dss Hermitian\dss vector bundle over $X$\nnsp.\oss
In order\sss to state a\qss \emph{boundary problem}\qss for sections of\sss $E$\sss
one needs a differential\sss or pseudo-differential\sss operator $P$ 
acting on sections of\sss $E$\sss and\sss some\qss 
\emph{boundary conditions}\pss on sections of\sss $E$\nnsp.\oss
In\sss this paper we are interested\sss in\trs \emph{local}\pss boundary conditions
described\sss by a\qss \emph{boundary operator}\qss $B$\sss
taking\sss sections of\sss $E$ over $X$\sss to sections of\dss a vector bundle $G$\sss over $Y$\dnsp.\oss
The boundary problem defined\dss by\sss $P\fff,\pff B$\sss
is\dss called an\qss \emph{elliptic\dss boundary\sss problem}\pss if\dss $P$\sss is\dss elliptic
and\sss $B$\sss satisfies\sss the\qss \emph{Shapiro-Lopatinskii\dss condition}.\oss

In\sss this paper we consider only\sss operators $P$ of\dss order $1$\nnsp.\oss
In\sss this case one can assume\sss that\sss $B$ has\sss the form\sss
$B\off =\off B_{\trf Y}\trf \circ\trf \gamma$\nnsp,\oss
where $\gamma$ takes\sss sections
of\sss $E$ over $X$\sss to\sss their restrictions\sss to $Y$\dnsp,\oss
and\sss $B_{\trf Y}$\sss is\dss a differential\sss or pseudo-differential\sss operator  
from sections of\dss $E\trf |\trf Y$\sss to sections of\sss $G$\nnsp.\oss
A crucial\sss role\dss is\dss played\sss by\sss the kernel $N$ of\dss the symbol\sss
of\sss $B_{\trf Y}$ over\sss the unit\sss cotangent\sss sphere bundle $S\fff Y$ of\sss $Y$\dnsp.\oss
It\dss is\dss a subbundle of\dss the\sss lift\sss to $S\fff Y$ of\dss  the restriction $E\trf |\trf Y$\dnsp.\oss
We will\sss call\dss the subbundle $N$ the\qss \emph{kernel-symbol}\qss of\dss $B_{\trf Y}$ and\sss $B$\nnsp.\oss
If\sss $N$\sss is\dss equal\dss to\sss the\sss lift\sss to $S\fff Y$ of\dss a subbundle $N_{\trf Y}$ of\sss
$E\trf |\trf Y$\dnsp,\oss we will\sss say\sss that $N$ is\qss
\emph{bundle-like},\oss borrowing\sss a\sss term of\qss Bandara,\qss Goffeng,\oss 
and\dss Saratchandran\qss \cite{bgs}.\oss
In\sss this case we will\sss usually write $N$\sss instead of\sss $N_{\trf Y}$\nnsp.\oss

\myuppar{The\dss Atiyah--Bott--Singer\dss index\sss theorem.}
If\dss $P\fff,\pff B$\sss is\dss elliptic,\oss
then\sss $P\dff \oplus\dff B$\sss defines\dss Fredholm\dss operators between
appropriate\dss Sobolev\dss spaces and\dss hence\sss the\sss index of\dss
$P\dff \oplus\dff B$\sss is\dss defined.\oss
By\sss the\dss Atiyah--Bott--Singer\dss index\sss theorem\qss \cite{ab1}\qss
this index\dss is\dss equal\sss to\sss the\qss \emph{topological\dss index}\pss
of\dss $P\fff,\pff B$\sss defined\sss in\sss terms of\dss 
the principal\sss symbols of\dss $P$ and $B_{\trf Y}$\nsp.\oss
According\sss to\sss the\sss introduction\sss to\dss Atiyah\dss and\dss Singer\qss \cite{as1},\oss
the\dss Atiyah--Bott--Singer\dss proof\dss of\dss this index\sss theorem\sss
was planned\dss to be included\sss in\sss the series of\dss papers\qss
\cite{as1}\qss --\qss \cite{as5}\qss by\trs Atiyah,\oss Segal,\oss and\trs Singer.\oss
But\sss this proof\dss was never published.\oss
An\sss improved\sss version of\dss the analytic\sss part\sss
was written down\dss by\trs H\"{o}rmander\qss \cite{h}.\oss
According\sss to\dss Atiyah\qss \cite{a1},\qss \cite{a2},\oss
the crux of\dss the matter\dss is\dss the definition of\dss the\sss topological\dss index.\oss
The definition\dss is\dss based on\qss (and actually\sss led\sss to)\qss
the elementary\sss proof\dss of\trs Bott\dss periodicity\dss theorem of\qss
Atiyah\dss and\dss Bott\qss \cite{ab2}.\oss
Another\sss key element\dss is\dss a deformation of\dss the symbol\sss of\dss $P$
to a special\sss form,\oss
which we will\sss call\qss \emph{bundle-like}.\oss
This deformation,\oss in\sss turn,\oss is\dss based on a deformation of\dss
the kernel-symbol\sss $N$\sss to a bundle-like one.\oss

\myuppar{Self-adjoint\dss elliptic operators on closed\sss manifolds.}
For\sss individual\sss self-adjoint\dss Fredholm\dss operators\sss
there\dss is\dss no analogue of\dss index.\oss
But\sss if\dss one considers\sss families of\dss operators 
instead of\dss individual\sss operators,\oss
the\sss index\sss theory for self-adjoint\trs  Fredholm\dss operators\sss turns out\sss
to be as rich as\sss the one for\dss Fredholm\dss operators.\oss
Such a\sss theory was developed\sss by\trs
Atiyah--Singer\qss \cite{as}\qss
and\sss applied\dss by\trs Atiyah--Patodi--Singer\qss \cite{aps}\qss
to prove an analogue of\dss the\dss Atiyah--Singer\dss theorem\sss
for families of\dss self-adjoint\sss elliptic operators on closed\sss manifolds.\oss

\myuppar{Self-adjoint\sss local\sss elliptic boundary\sss problems.}
It\dss is\dss only\sss natural\sss to ask\sss for an analogue of\dss the\dss
Atiyah--Bott--Singer\dss index\sss theorem\sss for
families of\dss self-adjoint\sss local\sss boundary problems,\oss
but,\oss apparently,\oss this question hardly attracted any serious attention
until\sss the\sss last\sss decade.\oss
About\sss ten years ago\sss the physics of\dss graphene and,\oss
in\sss particular,\oss questions posed\dss by\dss M.\dss Kats\-nelson,\oss served as a stimulus\sss
to address\sss this question at\sss least\sss for operators of\dss order $1$ and\dss families
parameterized\dss by circle,\oss i.e.\qss for\sss the spectral\sss flow.\oss
The relevant\sss physics\dss is\dss discussed\dss in\qss \cite{kn}.\oss
The case when\sss $\dim\dff X\off =\off 2$\sss was addressed\dss by\dss
Prokhorova\qss \cite{p1},\oss Katsnelson\sss and\dss Nazaikinskii\qss \cite{kn},\oss
and\sss then\sss by\dss Prokhorova\qss \cite{p2},\qss \cite{p3}\qss in\sss greater generality.\oss
The spectral\dss flow of\dss linear\sss families
of\trs Dirac\dss operators was considered\dss by\dss
Gorokhovsky\sss and\dss Lesch\qss \cite{gl}.\oss

In\sss the present\sss paper\sss this question\dss is\dss 
addressed\sss for operators of\dss order $1$ in arbitrary dimension.\oss
Many\sss features specific\sss to\sss the self-adjoint\sss case
are already\sss present\sss in\sss this case.\oss
Following\dss H\"{o}rmander\qss \cite{h},\oss
we do not\sss strive for\sss the greatest\sss possible\sss generality
and consider only\sss pseudo-differential\sss operators belonging\sss
to a class introduced\dss in\qss \cite{h},\oss Chapter\dss 20.\oss
Strictly speaking,\oss these operators are not\sss pseudo-differential\dss
near\sss the boundary,\oss except\sss when\sss they are actually differential\oss operators.\oss 
But\sss these operators can be approximated\dss by\sss
pseudo-differential\sss operators in a very strong sense.\oss

We will\sss consider families of\dss boundary problems\sss $P\fff,\pff B$\nnsp,\oss
parameterized\dss by a\sss topological\sss space $Z$\nnsp,\oss
such\sss that\sss the operators\sss $P$ are self-adjoint\sss and\sss induce unbounded
self-adjoint\sss operators in appropriate\dss Sobolev\dss spaces
when restricted\sss to\sss the kernel\sss $\kernel\dff B$\nnsp.\oss
The\sss latter property requires\sss the kernel-symbol\sss $N$\sss of\sss $B$\sss
to be a\sss lagrangian\sss subbundle with respect\sss to 
an\sss indefinite metric defined\sss 
in\sss terms of\dss the symbol\sss of\dss $P$ at\sss $Y$\dnsp.\oss
It\sss turns out\sss that\sss one also need\dss to assume\sss that\sss 
the boundary operators\sss $B$\sss are\qss \emph{bundle-like}\pss 
in\sss the sense\sss that\sss they are\dss
induced\dss by a morphisms of\dss bundles\sss $E\trf |\trf Y\qff \ttoo\qff G$\nnsp.\oss
Then\sss the kernel-symbols $N$\sss are also bundle-like.\oss

For such families one can define both\sss the\sss analytical\dss index\sss
reflecting\sss behavior of\dss the induced unbounded self-adjoint\sss operators 
and\dss the\sss topological\dss index depending only on\sss the families of\dss the symbols of\dss $P$
and\sss the kernel-symbols $N$\nnsp.\oss
Both\sss indices are elements of\dss $K^{\dff 1}\dff (\trf Z\trf)$\nnsp,\oss
as\sss they should\dss be.\oss
Our\sss main\sss result\dss is\trs Theorem\qss \ref{index-theorem-full}\qss
to\sss the effect\sss these\dss indices are equal.\oss

\myuppar{The\sss topological\dss index.}
Similarly\sss to\trs Atiyah--Bott--Singer\trs theorem\qss \cite{ab1},\oss
defining\sss the\sss topological\dss index\dss is\dss a key\sss part.\oss
As in\qss \cite{ab1},\oss this\dss is\dss done by\sss using\sss the 
boundary conditions\sss to deform\sss the symbol\sss $\sigma$ of\sss $P$\sss
to a standard\sss form on\sss the boundary.\oss

Similarly\sss to\qss \cite{ab1},\oss
we first\sss deform\sss the boundary conditions over\sss 
individual\dss points\sss $u\qff \in\qff S\dff Y$\dnsp.\oss
This involves some nice geometry of\dss self-adjoint\sss operators and\dss lagrangian subspaces
in finite dimensions and\dss is\dss the subject\sss
of\trs Section\qss \ref{boundary-algebra}.\oss

In\dss Section\qss \ref{symbols-conditions}\qss we perform\sss
deformations from\dss Section\qss \ref{boundary-algebra}\qss simultaneously\sss
for all\sss $u\qff \in\qff S\dff Y$\dnsp.\oss
This brings $\sigma$ over\sss the boundary\sss $Y$\sss 
into a standard\sss form depending almost\sss only on $N$\nnsp,\oss
but\dss in order\sss to get\sss an analogue of\trs Atiyah--Bott\trs definition
one needs also\sss to make\sss $N$ and\sss $\sigma$ bundle-like.\oss
In\qss \cite{ab1}\qss this\dss is\dss done  
by\sss the\qss \emph{Fifth\dss Homotopy}.\oss
Somewhat\sss disappointingly,\oss this\sss homotopy\dss does not\sss work
in\sss the self-adjoint\sss case.\oss
Moreover,\oss there\dss is\dss a non-trivial\sss obstruction\sss 
to making $\sigma$ and\sss $N$\sss bundle-like by a deformation.\oss
See\dss Section\qss \ref{symbols-conditions}.\oss
By\sss these reasons one needs\sss to require at\sss least\sss that\sss
$N$\sss can\sss be deformed\sss to a bundle-like kernel-symbol\sss after some stabilization.\oss
When such a deformation exists,\oss the\sss remaining\sss part\sss
of\dss the definition\dss is\dss similar\sss to\sss the\dss 
Atiyah--Singer\dss definition of\dss 
the\sss topological\dss index for families\qss \cite{as4}.\oss\vspace{-0.125pt}

\myuppar{The\sss analytical\dss index.}
Defining\sss the analytical\dss index also encounters new difficulties.\oss
The analytical\dss index of\dss families of\dss self-adjoint\sss elliptic operators on closed\sss
manifolds was defined\dss by\trs Atiyah--Patodi--Singer\qss \cite{aps}.\oss
This definition\dss is\dss based on\dss Atiyah--Singer\qss \cite{as}\qss 
theory of\dss self-adjoint\qss (actually,\oss skew-adjoint\fff)\pss \emph{bounded}\pss 
operators in a\dss Hilbert\dss space.\oss 
Applying\sss this definition\sss to differential\sss and\sss 
pseudo-differential\sss operators requires,\oss
as\sss the first\sss step,\oss 
replacing\dss given operators by operators of\dss order $0$\nnsp.\oss
With\sss the\sss theory of\dss pseudo-differential\sss operators at\sss hand,\oss
this\dss is\dss very simple on closed\sss manifolds,\oss but\sss this\dss is\dss not\sss so
on\sss manifolds with boundary.\oss
Author's approach\sss to\sss the analytical\sss index\qss \cite{i2}\qss
works equally\sss well\sss for\sss families of\dss bounded and\sss unbounded operators
in a\dss Hilbert\dss space with\sss
minimal\sss continuity\sss properties.\oss

Section\qss \ref{abstract-index}\qss complements\qss \cite{i2}\qss 
by a\dss little\sss theory of\dss abstract\dss boundary\sss problems
serving\sss as a bridge\sss between\sss pseudo-differential\dss boundary\sss problems
and\sss the\sss theory developed\sss in\qss \cite{i2}.\oss
This\sss theory\sss may be considered as an axiomatic\sss parameterized\sss version of\dss some aspects of\dss
the classical\dss theory of\dss boundary\sss problems as presented,\oss for example,\oss
in\qss \cite{w},\oss Section\qss 15.\oss

The framework of\trs Section\qss \ref{abstract-index}\qss resembles\sss that\sss of\dss 
the\sss theory of\pss \emph{boundary\dss triplets}\pss
(see,\oss for example,\oss Schm\"{u}dgen\qss \cite{s},\oss 
Chapter\qss 14\qss for\sss the\sss latter).\oss
While\sss the\sss theory of\dss boundary\sss triplets\dss
is\dss devoted\sss to a classification of\dss self-adjoint\sss extensions
of\dss an arbitrary\sss symmetric operator,\oss Section\qss \ref{abstract-index}\qss is\dss concerned\sss
with\sss proving\sss that\sss some\dss Fredholm\dss operators are self-adjoint,\oss
and\sss with\sss the dependence of\dss the operators on\sss parameters.\oss
Nevertheless,\oss the\sss two\sss theories are related,\oss
and we will\sss use\sss the\sss theory of\dss boundary\sss triplets\sss 
in\dss Section\qss \ref{realizations}\qss to explain why\sss it\dss is\dss 
necessary\sss to require\sss kernel-symbols\sss to be bundle-like.\oss\vspace{-0.125pt}

\myuppar{The proof\dss of\trs the index\sss theorem.}
The proof\dss of\dss our\sss index\sss theorem,\oss 
Theorem\qss \ref{index-theorem-full},\oss
consists of\dss two steps.\oss
The values of\dss the symbol\sss $\sigma$ of\dss $P$ 
on\sss (co)vectors normal\dss to\sss the boundary\sss $Y$ define 
a splitting\sss of\dss the bundle\sss $E\trf|\trf Y$\sss into a direct\sss sum\sss
$E\trf|\trf Y
\off =\off
E^{\dff +}_{\trf Y}\qff \oplus\qff E^{\dff -}_{\trf Y}$\nsp.\oss
First\sss we prove\sss the index\sss theorem\sss under\sss the assumption\sss
that\sss the bundle\sss $E^{\dff +}_{\trf Y}$ can\sss be extended\sss from\sss
the boundary\sss $Y$\sss to\sss the whole manifold\sss $X$\dss
(in a manner continuously depending\sss on\sss the parameter\sss $z\qff \in\qff Z$\nsp).\oss
This\dss is\dss done in\dss Theorem\qss \ref{index-theorem-ext}.\oss
Then\sss we reduce\sss the general\sss case\sss to\sss this one.\oss

In order\sss to prove\dss  Theorem\qss \ref{index-theorem-ext}\qss
we adapt\sss to\sss the self-adjoint\dss boundary problems\dss
H\"{o}rmander's\dss approach\sss to\dss Atiyah--Bott-Singer\dss
index\sss theorem.\oss
See\qss \cite{h},\oss Chapter\qss 20.\oss
The first\sss step\dss is\dss to construct\sss sufficiently\sss many standard\dss boundary problems
guaranteed\dss to have vanishing\sss analytical\sss index.\oss
This\dss is\dss done in\dss Section\qss \ref{pdo}\qss using\trs
Proposition\qss 20.3.1\qss from\qss \cite{h}.\oss
The next\sss step\dss is\dss to\sss take\qss (families of{\dff})\qss boundary problems
standard\dss near\sss the boundary and\sss glue\sss them\sss to standard ones\sss
to get\sss self-adjoint\sss operators on\sss the double of\dss $X$\nnsp,\oss
which\dss is\dss a closed\sss manifold.\oss
This\dss is\dss the\sss topic of\dss Section\qss \ref{glueing}.\oss
This constructions preserves both\sss the analytical\sss
and\sss the\sss topological\sss indices.\oss
The proof\dss for\sss the analytical\sss index\dss is\dss based on\sss the proof\dss of\trs
Proposition\qss 20.3.2\qss in\qss \cite{h}.\oss
The proof\dss for\sss the\sss topological\sss index\dss is\dss fairly\sss routine.\oss
In\dss Section\qss \ref{index-theorems}\qss
we apply\sss the results of\trs Sections\qss \ref{pdo}\qss and\qss \ref{glueing}\qss
to reduce\sss the index\sss problem\sss to\sss the case of\dss closed\sss manifolds,\oss
and\sss then apply\sss the\dss Atiyah--Patodi--Singer\dss index\sss theorem\qss \cite{aps}.\oss
This proves\dss Theorems\qss \ref{index-theorem-ext}.\vspace{1pt}

The additional\sss assumption\sss in\dss Theorem\qss \ref{index-theorem-ext}\qss
is\dss caused\sss by\sss the\sss limitations of\dss the construction of\dss
boundary\sss problems with\sss index\sss zero used\sss in\dss Section\qss \ref{pdo}.\oss
A simple doubling\sss trick\sss shows\sss that\sss this assumption always holds
for\sss the direct\sss sum of\dss a boundary problem\sss with a copy of\dss it.\oss
This proves\dss Corollary\qss \ref{index-theorem-two},\oss an\sss index\sss theorem 
modulo elements of\dss order $2$ in\sss $K^{\dff 1}\dff (\trf Z\trf)$\nnsp.\oss\vspace{1pt}

In order\sss to prove\sss the\sss full\dss index\sss theorem,\oss
Theorem\qss \ref{index-theorem-full},\oss
we use a much\sss more sophisticated doubling\sss tool.\oss
Namely,\oss we\sss take\sss the $\omult_{\fff 1}$ product\sss
of\dss our boundary\sss problems\sss with an elliptic\qss (not\sss self-adjoint\fff)\qss
operator $Q$ of\dss index $1$ in a vector\sss bundle $E\fff'$ 
on an auxiliary\sss closed\sss manifold\sss $V$\dnsp.\oss 
The $\omult_{\fff 1}$ product\dss is\dss the self-adjoint\sss version
of\dss the product\sss used\dss by\dss Atiyah--Singer\qss \cite{as1}.\oss
Both versions were used already\sss in\dss Palais\dss seminar\qss \cite{pa}.\oss
It\sss turns out\sss that\sss the $\omult_{\fff 1}$ product\sss with $Q$\sss
satisfies\sss the assumption of\trs Theorem\qss \ref{index-theorem-ext}.\oss
At\sss the same\sss time\sss taking\sss the $\omult_{\fff 1}$ product\sss with $Q$\sss
preserves both\sss the\sss topological\sss and\dss the analytical\dss indices.\oss
These facts reduce\sss the full\dss index\sss theorem\sss to\sss the special\sss
case considered\sss in\dss Theorem\qss \ref{index-theorem-ext}.\oss\vspace{1pt}

Proving\sss these properties of\dss the $\omult_{\fff 1}$ products\dss  
is\dss most\sss technical\sss part\sss of\dss the proof\dss
of\trs Theorem\qss \ref{index-theorem-full}.\oss
We discuss\sss the $\omult_{\fff 1}$ products of\dss symbols,\oss
operators in\dss Hilbert\dss spaces,\oss
and\sss pseudo-dif\-fer\-en\-tial\sss operators\sss
in\dss Sections\qss \ref{mult},\qss \ref{mult-operators},\oss 
and\pss \ref{mult-pdo}\qss respectively.\oss
The proof\dss of\trs Theorem\qss \ref{index-theorem-ext},\oss
as also\sss the results of\trs Sections\qss \ref{relations}\dss --\dss \ref{odd},\oss
do not\sss depend on\sss properties of\dss $\omult_{\fff 1}$ products.\oss\vspace{1pt}

\myuppar{Realization of\dss kernel-symbols.}
For\dss Atiyah--Singer\qss \cite{as1},\qss \cite{as4}\qss approach\sss to\sss index\sss theorems
one needs\sss to be able\sss to define not\sss only\sss the analytical\sss index of\dss
operators and\sss families of\dss operators,\oss but\sss also\sss the analytical\sss index of\dss
symbols and\dss families of\dss symbols.\oss
This\dss is\dss done by using\sss pseudo-differential\sss operators\sss to construct\qss
\emph{realizations}\qss of\dss symbols by operators in spaces of\dss sections
of\dss vector bundles,\oss mostly\sss in\dss Sobolev\dss spaces.\oss\vspace{1pt}

It\dss turns out\sss that\sss for self-adjoint\dss boundary problems\sss
there\dss is\dss a non-trivial\sss obstruction\sss to\sss the existence of\pss 
\emph{self-adjoint}\pss realizations.\oss
It\dss is\dss similar\sss to\sss the obstruction\sss to defining\sss the\sss topological\dss index\qss
(actually,\oss it\dss is\dss weaker\sss than\sss the\sss latter).\oss
Moreover,\oss when\sss this obstruction vanish,\oss the realization\dss is\dss not\sss unique,\oss
even up\sss to homotopy.\oss
While\sss the\sss topological\dss index does not\sss depend\sss on\sss the realization,\oss
the analytical\dss index does,\oss and very dramatically.\oss
This means\sss that,\oss in\sss general,\oss 
the analytical\dss index of\dss families of\dss self-adjoint\dss
pseudo-differential\dss 
boundary\sss problems\dss is\dss not\sss determined\dss by\sss
their principal\sss symbols.\oss
But\sss if\dss the boundary conditions are bundle-like,\oss
then\sss there\dss is\dss a canonical\sss realization depending only on\sss the symbol.\oss
This explains why\sss the proof\dss of\dss the index\sss theorem\sss 
works only\sss for\sss bundle-like boundary conditions.\vspace{1pt}

These\sss issues are discussed\sss in\dss Section\qss \ref{realizations}.\oss
In\sss the\sss last\sss subsection of\dss this section\sss we 
construct\sss  families of\dss operators
having\sss the\sss topological\dss index\sss zero and,\oss 
at\dss least\sss for compact\sss space of\dss parameters,\oss 
all\sss elements of\dss $K^{\dff 1}\dff (\trf Z\trf)$
as\sss their analytical\dss indices.\oss
Except\sss of\dss the\sss last\sss subsection,\oss   
Section\qss \ref{realizations}\qss depends only on\dss
Sections\qss \ref{boundary-algebra},\oss \ref{symbols-conditions},\oss 
\ref{abstract-index},\oss and\pss \ref{pdo}.\oss
The\sss last\sss subsection depends on\sss the\sss theory of\qss
\emph{self-adjoint\dss relations}\pss and\qss 
\emph{boundary\trs triplets},\oss
the\sss topic of\trs Sections\qss \ref{relations}\qss 
and\qss \ref{boundary-triplets}.\oss

A similar phenomenon of\dss the non-uniqueness of\dss realizations
in\sss the context\sss of\trs Atiyah--Patodi--Singer\trs boundary conditions\sss
was encountered and\dss 
thoroughly\sss investigated\dss by\trs
Melrose\dss and\dss Piazza\qss \cite{mp1},\qss \cite{mp2}.\oss
These boundary conditions\sss define canonical\sss realizations of\dss some ker\-nel-symbols,\oss
but\sss do not\sss continuously depend on parameters.\oss

\myuppar{Self-adjoint\sss relations and\dss boundary\sss triplets.}
We review\sss the\sss theory of\dss self-adjoint\sss relations and\dss boundary\sss triplets
in\dss Sections\dss \ref{relations}\dss and\dss \ref{boundary-triplets},\oss
referring\sss to\sss the\sss textbook of\trs 
Schm\"{u}dgen\qss \cite{s}\qss for\sss the details and\dss the history of\dss this\sss theory.\oss
In\dss Section\qss \ref{relations}\qss we show\sss that\sss there are natural\sss extensions of\dss 
notions of\qss Fredholm\dss families and of\dss the analytical\dss index of\dss such\sss families\sss
to self-adjoint\sss relations.\oss
In\dss Section\qss \ref{boundary-triplets}\qss we use\sss the famous\dss Krein--Naimark\dss
resolvent\sss formula\sss to prove a\qss ``relative\sss index\sss theorem''\qss
in an abstract\sss setting.\oss
See\qss Theorem\qss \ref{spectral-triplet-index}.\oss
In\dss Section\qss \ref{boundary-triplets}\qss we also connect\sss 
the\sss theory of\dss boundary\sss triplets with\sss the\sss theory
developed\sss in\dss Section\qss \ref{abstract-index}.\oss

\myuppar{Special\dss boundary conditions.}
We will\sss say\sss that\sss a boundary condition\dss is\qss \emph{special}\pss
if\dss it\dss satisfies not\sss only\sss the\dss Shapiro-Lopatinskii\dss condition,\oss
but\sss also a natural\qss ``dual''\qss form of\dss it.\oss
For example,\oss if\dss $P$\sss is\dss a\qss \emph{differential\sss operator}\qss
of\dss order $\nsp 1$\nnsp,\oss then every elliptic boundary condition\dss is\dss
automatically\sss special.\oss
Section\qss \ref{special-conditions}\qss is\dss devoted\sss
to such boundary conditions for self-adjoint\sss symbols minimally\sss
resembling\sss at\sss the boundary $Y$\sss symbols of\trs Dirac operators.\oss
We call\sss such symbols\qss \emph{anti-commuting}.\oss
See\dss Section\qss \ref{special-conditions}\qss for\sss the precise definition.\oss
The main\sss result\sss of\dss this section\dss is\dss an\dss
Agranovich--Dynin\dss type\sss theorem\sss
for\sss the\sss topological\dss index,\oss Theorem\qss \ref{t-index-difference}.\oss

Given an anti-commuting\sss self-adjoint\sss symbol $\sigma$ together with 
a special\dss boundary condition,\oss
one can define a natural\dss boundary condition\qss ``dual''\qss to\sss the given one.\oss
One can also do\sss the same for families.\oss
Theorem\qss \ref{t-index-difference}\qss expresses\sss the difference of\dss the\sss
topological\dss indices of\dss $\sigma$ with\sss these\sss two boundary conditions\sss
in\sss terms of\dss the restriction of\dss $\sigma$\sss to\sss the boundary.\oss
The proof\dss naturally\sss leads\sss to the original\trs Bott\trs periodicity\sss map.\oss
See\dss Theorem\qss \ref{basic-difference}\trs and\dss its\sss proof.\oss

\myuppar{Dirac-like boundary\sss problems.}
In\dss Section\qss \ref{odd}\qss we consider an even more narrow class of\dss self-adjoint\sss
operators\dss $P$\dss and\sss boundary conditions.\oss 
We call\dss the corresponding\dss boundary problems\qss \emph{Dirac-like},\oss
although we require only a small\dss part\sss of\dss the rich algebraic structure 
of\trs Dirac\dss operators\sss to be present.\oss
Basically,\oss the operator $P$\sss is\dss required\sss to be odd\sss with respect\sss to a
$\zzz/2$\dnsp-grading of\dss the bundle $E$\nnsp,\oss
its symbol\dss is\dss required\dss to induce a skew-adjoint\sss symbol\sss
on a\qss ``half\dff''\qss of\dss the bundle $E\trf |\trf Y$\dnsp,\oss
and\sss the boundary condition\dss is\dss required\sss to commute with\sss this induced symbol.\oss
See\dss Section\qss \ref{odd}\qss for\sss the precise definitions.\oss

The main\sss result\sss of\trs Section\qss \ref{odd}\qss is\trs Theorem\qss \ref{dirac-flip}
expressing\sss the\sss topological\sss and\dss the analytical\dss
index of\dss families of\dss such boundary problems\sss in\sss terms of\dss
the\sss topological\dss index of\dss some families of\dss 
self-adjoint\dss operators induced on\sss the boundary.\oss
See\dss Theorem\qss \ref{dirac-flip}.\oss
This\sss theorem\dss is\dss motivated\dss by\sss the results of\trs
Gorokhovsky\sss and\dss Lesch\qss \cite{gl},\oss
who considered only\dss Dirac\dss operators\qss
({\fff}but\sss indicated\dss that\sss many of\dss their arguments are valid\sss
in\sss greater\sss generality)\qss
and\sss only\sss the spectral\dss flow
of\trs linear families.\oss
The proof\dss of\trs Theorem\qss \ref{dirac-flip}\qss is\dss almost\sss entirely\sss
topological,\oss in contrast\sss with\sss the heat\sss equation methods used\sss in\qss \cite{gl}.\oss

\newpage
\mysection{Geometric\qss algebra\qss at\qss the\qss boundary}{boundary-algebra}

\myuppar{Elliptic pairs.}
Let\sss $E$\sss be a finitely dimensional\sss vector space over\sss $\ccc$\sss 
equipped\sss with a\dss Hermitian\dss positive definite 
scalar product\sss $\sco{\dff \bullet\dff,\qff \bullet \dff}$\dnsp,\oss
and\dss let\dss $\sigma,\qff \tau\dff \colon\dff E\qff \ttoo\qff E$\sss 
be\sss two\sss linear operators.\oss
We will\sss say\sss that\sss $\sigma,\qff \tau$\sss is\dss an\qss
\emph{elliptic\sss pair}\pss if\dss the operator\sss
$a\trf \sigma\qff +\qff b\trf \tau$\sss is\dss invertible\sss
for every\sss $a\fff,\qff b\qff \in\qff \rrr$\sss such\sss that
$(\dff a\fff,\qff b\trf)\off \neq\off 0$\nnsp.\oss
Then,\oss in\sss particular,\qss $\sigma$ and $\tau$ are invertible.\oss
Moreover,\pss $\lambda\qff +\qff \sigma^{\dff -\dff 1}\dff \tau$\sss
is\dss invertible for every\sss 
$\lambda\qff \in\qff \rrr$\nnsp.\oss 
In other words,\sss 
$\sigma^{\dff -\dff 1}\dff \tau$\sss has no real\sss eigenvalues.\oss
For\sss the rest\sss of\dss this section we will\sss assume\sss that\sss
$\sigma,\qff \tau$\sss is\dss an elliptic pair.\oss 
Let\sss $\rho\off =\off \sigma^{\dff -\dff 1}\dff \tau$\nnsp.\oss

\myuppar{The ordinary differential\sss equation associated\sss
with\sss $\sigma,\qff \tau$\dnsp.}
Let $t$\sss be a variable running over $\rrr$\sss and\dss let\sss $\partial\off =\off d/d t$\sss
be\sss the usual\sss differentiation operator.\oss
Let\sss $D\off =\off -\qff i\qff \partial$\nnsp.\oss
We are interested\sss in\sss the differential\sss equation
$(\trf \sigma\trf D\qff +\qff \tau\trf)\dff(\trf f\trf)
\off =\off 
0$
for functions $f\dff \colon\dff \rrr\qff \ttoo\qff E$\nnsp.\oss
Obviously,\oss it\dss is\dss equivalent\sss to\sss the equation\sss
$D\trf(\trf f\trf)\qff +\qff \rho\trf (\trf f\trf)\off =\off 0$\nnsp.\oss
It\dss is\dss well\dss known\sss that\sss every\sss solution of\dss such an equation\dss is\dss a\sss
linear\sss combination of\dss solutions of\dss the form\sss
$f\dff(\dff t\trf)\off =\off p\dff(\dff t\trf)\qff e^{\dff -\dff i\dff \lambda\fff t}$\nnsp,\oss
where\sss $\lambda\qff \in\qff \ccc$\sss is\dss an eigenvalue of\sss $\rho$\sss and 
$p\dff(\dff t\trf)$ is\dss a polynomial\sss
with coefficients in\sss the generalized eigenspace\sss 
$\mathcal{E}_{\trf \lambda}\dff(\trf \rho\dff)\qff \subset\pff E$ of\dss
$\rho$\sss corresponding\dss to $\lambda$\nnsp.\oss\vspace{-0.4pt}

A solution of\dss the form\sss $p\dff(\dff t\trf)\qff e^{\dff -\dff i\dff \lambda\fff t}$\sss
is\dss bounded on\sss $\rrr_{\qff \geq\dff 0}$\sss
if\dss and\sss only\dss if\dss either $\re -\qff i\dff \lambda\qff <\qff 0$\sss
or\sss $\re -\qff i\dff \lambda\off =\off 0$\sss and\sss 
$p\dff(\dff t\trf)$ is\dss a constant.\oss
Since\sss $\rho$\sss has no real\sss eigenvalues,\oss this\sss condition\dss
is\dss equivalent\sss to\sss $\im\dff \lambda\qff <\qff 0$\sss and such solutions are actually\sss
exponentially decreasing when\sss $t\qff \ttoo\qff \infty$\nnsp.\oss
The solutions of\dss the form $p\dff(\dff t\trf)\qff e^{\dff -\dff i\dff \lambda\fff t}$\sss
with\sss $\im\dff \lambda\qff >\qff 0$\sss are\sss
exponentially increasing when\sss $t\qff \ttoo\qff \infty$\nnsp.\oss
Let\sss us\sss denote by\sss $\mathcal{M}_{\dff +}\dff(\trf \rho\dff)$\sss 
and\sss $\mathcal{M}_{\dff -}\dff(\trf \rho\dff)$\sss
the spaces of\dss solutions exponentially decreasing and,\oss respectively,\oss
increasing\sss when\sss $t\qff \ttoo\qff \infty$\nnsp.\oss\vspace{-0.4pt}

At\sss the same\sss time a solution $f\dff(\dff t\trf)$\sss is\dss determined\sss by\dss its
initial\sss value\sss $f\dff(\dff 0\dff)$\nnsp.\oss
It\sss follows\sss that\sss $\mathcal{M}_{\dff +}\dff(\trf \rho\dff)$\sss is\dss canonically\sss isomorphic\sss
to\sss the space\sss $\mathcal{L}_{\dff -}\dff(\trf \rho\dff)\qff \subset\pff E$\nnsp,\oss
the sum of\dss the generalized eigenspaces of\dss $\rho$\sss corresponding\sss
to\sss the eigenvalues $\lambda$\sss with\sss $\im\dff \lambda\qff <\qff 0$\nnsp.\oss
Similarly,\pss $\mathcal{M}_{\dff -}\dff(\trf \rho\dff)$\sss is\dss canonically\sss isomorphic\sss
to\sss the space\sss $\mathcal{L}_{\dff +}\dff(\trf \rho\dff)\qff \subset\pff E$\nnsp,\oss
the sum of\dss the generalized eigenspaces of\dss $\rho$\sss corresponding\sss
to\sss the eigenvalues $\lambda$\sss with\sss $\im\dff \lambda\qff >\qff 0$\nnsp.\oss\vspace{-0.4pt}

\myuppar{Self-adjoint\sss elliptic pairs and\dss the associated\dss Pontrjagin\sss spaces.}
From\sss now on we will\sss assume\sss that\sss the elliptic pair\sss
$\sigma,\qff \tau$\sss is\qss \emph{self-adjoint}\pss
in\sss the sense\sss that\sss the operators\sss $\sigma,\qff \tau$\sss are
self-adjoint\sss with respect\sss to\sss the scalar\sss product
$\sco{\dff \bullet\dff,\qff \bullet \dff}$\dnsp.\oss
As usual,\oss for a\sss linear operator 
$\varphi\dff \colon\fff E\qff \ttoo\qff E$\sss
we will\sss denote\sss by\sss $\varphi^{\fff *}$\sss 
its adjoint\sss with respect\sss to\sss $\sco{\dff \bullet\dff,\qff \bullet \dff}$\nnsp.\oss
Let\sss us define a new scalar\sss product $[\trf \bullet\fff,\qff \bullet\trf]$
on\sss $E$\nnsp,\oss
which\sss in\sss general\sss will\dss be indefinite.\oss
Namely,\oss for\sss $u\fff,\qff v\qff \in\pff E$\sss let\sss\vspace{-0.4pt} 
\[
\quad
[\trf u\fff,\qff v\trf]
\off =\off 
\sco{\dff \sigma\trf(\trf u\trf)\fff,\qff v\dff}
\qff.
\]

\vspace{-12pt}\vspace{-0.4pt}
Since\sss the operator $\sigma$\sss is\dss self-adjoint\sss and\sss $\sco{\dff \bullet\dff,\qff \bullet \dff}$\sss
is\dss Hermitian,\pss 
$\sco{\dff \sigma\trf(\trf u\trf)\fff,\qff v\dff}
\off =\off
\overline{\sco{\dff \sigma\trf(\trf v\trf)\fff,\qff u\dff}}$\trs
and\sss hence\sss 
$[\trf u\fff,\qff v\trf]
\off =\off
\overline{[\trf v\fff,\qff u\trf]}$\sss
for every\sss $u\fff,\qff v\qff \in\pff E$\nnsp,\oss
i.e.\qss the scalar\sss product\sss $[\trf \bullet\fff,\qff \bullet\trf]$\sss is\dss Hermitian.\oss
Since $\sigma$\sss is\dss invertible,\oss
the scalar\sss product\sss $[\trf \bullet\fff,\qff \bullet\trf]$\sss is\dss non-degenerate.\oss
Therefore\sss $E$\sss together with\sss the product\sss $[\trf \bullet\fff,\qff \bullet\trf]$\sss
is\dss a\qss \emph{Pontrjagin\dss space},\oss 
or an\qss \emph{indefinite\sss inner\sss product\sss space}.\oss

The operator\sss $\rho$\sss turns out\sss to be\qss
\emph{self-adjoint\sss with\sss respect\sss to}\pss
$[\trf \bullet\fff,\qff \bullet\trf]$\nnsp,\oss
i.e.\qss such\sss that\vspace{1.5pt}
\[
\quad
[\qff \rho\dff(\trf u\trf)\fff,\qff v\trf]
\off =\off
[\trf u\fff,\qff \rho\dff(\trf v\trf)\trf]
\pff
\]

\vspace{-12pt}\vspace{1.5pt}
for every\sss $u\fff,\qff v\qff \in\pff E$\nnsp.\oss
Indeed,\oss since\sss $\tau$\sss is\dss self-adjoint,\pss
$\sco{\dff \tau\trf(\trf u\trf)\fff,\qff v\dff}
\off =\off
\sco{\dff u\fff,\qff \tau\trf(\trf v\trf)\dff}$\nnsp,\pss
or,\oss what\dss is\dss the same,\pss
$\sco{\dff \sigma\dff \circ\dff \rho\trf(\trf u\trf)\fff,\qff v\dff}
\off =\off
\sco{\dff u\fff,\qff \sigma\dff \circ\dff \rho\trf(\trf v\trf)\dff}$\dss
for every\sss $u\fff,\qff v\qff \in\pff E$\nnsp.\oss
This implies\sss that\sss\vspace{1.75pt}
\[
\quad
[\qff \rho\dff(\trf u\trf)\fff,\qff v\trf]
\off =\off
\sco{\trf \sigma\dff \circ\dff \rho\dff(\trf u\trf)\fff,\qff v\trf}
\off =\off
\sco{\dff u\fff,\qff \sigma\dff \circ\dff \rho\trf(\trf v\trf)\dff}
\off =\off
\sco{\dff \sigma\trf(\dff u\trf)\fff,\qff \rho\trf(\trf v\trf)\dff}
\off =\off
[\trf u\fff,\qff \rho\dff(\trf v\trf)\trf]
\pff,
\]

\vspace{-12pt}\vspace{1.75pt}
as claimed.\oss
The operator\sss $\rho$\sss cannot\sss be self-adjoint\sss
with respect\sss to\sss $\sco{\dff \bullet\dff,\qff \bullet \dff}$\nnsp.\oss
Indeed,\vspace{1.75pt}
\[
\quad
\rho^{\fff *}
\off =\off 
(\trf \sigma^{\dff -\dff 1}\dff \tau\trf)^{\fff *}
\off =\off
\tau^{\dff *}\dff (\trf \sigma^{\dff -\dff 1}\trf)^{\fff *}
\off =\off
\tau\trf \sigma^{\dff -\dff 1}
\off =\off
\sigma\trf \sigma^{\dff -\dff 1}\trf \tau\trf \sigma^{\dff -\dff 1}
\off =\off
\sigma\trf \rho\trf \sigma^{\dff -\dff 1}
\pff.
\]

\vspace{-12pt}\vspace{1.75pt}
Hence\sss $\rho^{\fff *}\off =\off \rho$\dss 
if\dss and\sss only\dss if\dss
$\rho$ commutes\sss with\sss $\sigma$\nnsp,\oss
or,\oss equivalently,\pss $\tau$\sss commutes with\sss $\sigma$\dnsp.\oss
But,\oss if\dss $\tau$\sss commutes with\sss $\sigma$\dnsp,\oss
then\sss there exists a common eigenvector\sss of\sss $\sigma$\sss and\sss $\tau$\dnsp,\oss
and\sss this contradicts\sss to\sss the ellipticity assumption.\oss
At\sss the same\sss time\sss $\rho$\sss is\dss often\sss skew-adjoint\sss with respect\sss to\sss
$\sco{\dff \bullet\dff,\qff \bullet \dff}$\nnsp,\oss
i.e.\qss such\sss that\sss $\rho^{\fff *}\off =\off -\qff \rho$\nnsp.\oss 
Clearly,\oss this happens\sss if\dss and\sss only\dss if\sss
$\rho$ anti-commutes\sss with\sss $\sigma$\nnsp,\oss
or,\oss equivalently,\pss $\tau$\sss anti-commutes with\sss $\sigma$\dnsp.\oss
Such elliptic pairs\sss $\sigma\fff,\qff \tau$\sss usually arise from\dss Clifford\dss modules,\oss 
but,\sss as we will\sss see,\oss they\sss are quite ubiquitous. 

By\sss introducing\sss the factor\sss $i\off =\off \sqrt{\dff -\qff 1}$\sss
at\sss various places one can rephrase our discussion\sss in\sss the\sss
language of\qss \emph{complex symplectic spaces},\oss
but\sss the\sss language of\qss Pontrjagin spaces seems\sss to be more natural\dss
for our purposes.\oss

\mypar{Lemma.}{orthogonality}
\emph{If\pss $\lambda\off \neq\off \overline{\mu}$\nsp,\oss
then\sss the generalized eigenspaces\sss
$\mathcal{E}_{\trf \lambda}\dff(\trf \rho\dff)$\sss and\sss
$\mathcal{E}_{\trf \mu}\dff(\trf \rho\dff)$\sss
are orthogonal\sss with respect\sss to\sss
the scalar\sss product\dss $[\trf \bullet\fff,\qff \bullet\trf]$\nnsp.\oss}

\proof
The proof\trs is\dss similar\sss to\sss the one in\sss the positive definite case.\oss
See\dss Pontrjagin\qss \cite{p},\oss Section\qss 2\dff({\fff}A{\fff}),\oss 
or\dss Bong\'{a}r\qss \cite{b},\oss Theorem\qss 3.3.\oss  \eproof

\mypar{Corollary.}{plus-minus-orth}
\emph{The restrictions of\qss $[\trf \bullet\fff,\qff \bullet\trf]$\sss
to\sss $\mathcal{L}_{\dff +}\dff(\trf \rho\dff)$
and\dss to\sss $\mathcal{L}_{\dff -}\dff(\trf \rho\dff)$\sss
are equal\dss to zero.\oss}  \eproof

\mypar{Corollary.}{dimension}
\emph{The dimensions of\trs the spaces\sss
\sss 
$\mathcal{L}_{\dff +}\dff(\trf \rho\dff)$
and\dss $\mathcal{L}_{\dff -}\dff(\trf \rho\dff)$\sss
are equal.\oss}

\proof
Since\sss $[\trf \bullet\fff,\qff \bullet\trf]$\sss
is\dss non-degenerate,\oss for every\sss linear\sss map\sss
$l\dff \colon\dff \mathcal{L}_{\dff +}\dff(\trf \rho\dff)\qff \ttoo\qff \ccc$\sss
there exists\sss $v\qff \in\qff E$\sss such\sss that\sss
$l\trf(\dff u\trf)\off =\off [\trf u\fff,\qff v\trf]$\sss
for\sss $u\qff \in\qff \mathcal{L}_{\dff +}\dff(\trf \rho\dff)$\nnsp.\oss
Since\sss $[\trf u\fff,\qff w\trf]\off =\off 0$\sss
for\sss $u\fff,\qff w\qff \in\qff \mathcal{L}_{\dff +}\dff(\trf \rho\dff)$\nnsp,\oss
there exists\sss $v\qff \in\qff \mathcal{L}_{\dff -}\dff(\trf \rho\dff)$\sss 
with\sss this property.\oss
It\sss follows\sss that\sss the natural\sss map from\sss
$\mathcal{L}_{\dff -}\dff(\trf \rho\dff)$\sss to\sss the dual\sss space of\dss
$\mathcal{L}_{\dff +}\dff(\trf \rho\dff)$\sss defined\sss by\sss
$[\trf \bullet\fff,\qff \bullet\trf]$\sss is\dss surjective.\oss 
Hence\sss
$\dim\dff \mathcal{L}_{\dff +}\dff(\trf \rho\dff)
\qff \leq\qff
\dim\dff \mathcal{L}_{\dff -}\dff(\trf \rho\dff)$\nnsp.\oss
By\sss the same reasons\sss
$\dim\dff \mathcal{L}_{\dff -}\dff(\trf \rho\dff)
\qff \leq\qff
\dim\dff \mathcal{L}_{\dff +}\dff(\trf \rho\dff)$\nnsp.\oss  \eproof

\myuppar{Lagrangian subspaces.}
Here we\sss largely\sss follow\qss F.\dss Latour\qss \cite{l},\oss Section\qss I.1.\oss
A subspace\sss $L\qff \subset\qff E$\sss is\dss said\sss to be\qss
\emph{lagrangian}\qss if\dss $L$\sss is\dss equal\sss to its orthogonal\sss complement\sss
with respect\sss to\sss $[\trf \bullet\fff,\qff \bullet\trf]$\nnsp.\oss
Since\sss $[\trf \bullet\fff,\qff \bullet\trf]$\sss is\dss non-degenerate,\oss
the dimension of\dss a\sss lagrangian subspace\dss is\dss equal\sss to\sss the
half\dss of\dss the dimension of\dss $E$\nnsp.\oss
In\sss particular,\oss if\dss a\sss lagrangian subspace exists,\oss
then\sss $\dim\dff E$\sss is\dss even.\oss
Two\sss lagrangian subspaces\sss $L\fff,\qff L'$\sss are said\dss to be\qss
\emph{transverse}\pss if\dss $L\qff +\qff L'\off =\off E$\nnsp,\oss
or,\oss equivalently,\qss $L\dff \cap\dff L'\off =\off 0$\nnsp.\oss 
Our main example of\dss transverse\sss lagrangian subspaces\dss is\sss 
$\mathcal{L}_{\dff +}\dff(\trf \rho\dff)\fff,\off 
\mathcal{L}_{\dff -}\dff(\trf \rho\dff)$\nnsp.\oss

For a\sss lagrangian subspace\sss $L\qff \subset\qff E$\sss let\sss us\sss
denote\sss by\sss $\Lambda_{\fff L}$\sss the set\sss of\dss all\dss
lagrangian subspaces\sss transverse\sss to\sss $L$\nnsp.\oss
Let\sss $M\qff \in\qff \Lambda_{\fff L}$\sss and\dss let\sss
$q\dff \colon\dff E\qff \ttoo\qff L$\sss
be\sss the projection along\sss $M$\nnsp.\oss
Then\vspace{1.5pt}
\[
\quad
q\dff(\dff x\trf)\off =\off x
\off\off
\mbox{for every}\qff\off
x\qff \in\qff L
\off\qff
\mbox{and}\off\off
\]

\vspace{-36pt}
\[
\quad
[\trf q\dff(\trf u\trf)\fff,\qff v\trf]
\qff +\qff
[\trf u\fff,\qff q\dff(\trf v\trf)\trf]
\off =\off
[\dff u\fff,\qff v\trf]
\off\off
\mbox{for every}\qff\off
u\fff,\qff v\qff \in\qff E
\qff.
\]

\vspace{-12pt}\vspace{1.5pt}
Indeed,\oss the first\sss property\dss is\dss trivial,\oss and\dss it\dss is\dss
sufficient\sss to check\sss the second one in\sss the case when each of\sss $u\fff,\qff v$\sss
belongs\sss to either\sss $M$\sss or\sss $L$\nnsp.\oss
But\sss in\sss this case\sss it\dss is\dss trivial.\oss
Conversely,\oss if\dss $q\dff \colon\dff E\qff \ttoo\qff L$\sss has\sss these properties,\oss
then an\sss immediate verification shows\sss that\sss 
$M\off =\off \kernel q$\sss is\dss a\sss lagrangian\sss subspace\sss transverse\sss to\sss $L$\nnsp.\oss
If\sss $q\fff'$\sss is\dss another map with\sss these properties,\oss
then\sss $\delta\off =\off q\qff -\qff q\fff'$\sss is\dss a\sss linear map equal\sss to zero on\sss
$L$\sss and such\sss that\vspace{1.5pt}
\begin{equation}
\label{delta}
\quad
[\trf \delta\dff(\trf u\trf)\fff,\qff v\trf]
\qff +\qff
[\trf u\fff,\qff \delta\dff(\trf v\trf)\trf]
\off =\off
0
\end{equation}

\vspace{-12pt}\vspace{1.5pt}
for every\sss $u\fff,\qff v\qff \in\qff E$\nnsp.\oss
Therefore\sss $\Lambda_{\fff L}$\sss is\dss an affine space over\sss the
vector space of\dss such\sss maps $\delta$\nnsp.\oss
In\sss particular,\qss $\Lambda_{\fff L}$ is\dss contractible.\oss
The structure of\dss an affine space on\sss $\Lambda_{\fff L}$\sss
does not\sss depend on\sss $M$\nnsp.\oss
If\sss a base point\sss $M\qff \in\qff \Lambda_{\fff L}$\sss is\dss fixed,\oss
then we can\sss identify $\Lambda_{\fff L}$ with\sss the vector space of\dss
such maps $\delta$\nnsp.\oss
The\sss lagrangian subspace corresponding\sss 
to\sss $\delta$\sss is\dss equal\dss to\sss
the kernel\sss $\kernel (\dff q\qff -\qff \delta\trf)$\sss
and\sss hence\sss is\dss equal\sss to\sss the graph\sss
$\{\trf u\qff +\qff \delta\dff(\dff u\trf)\qff |\qff u\qff \in\qff M\qff\}$\sss
of\dss the restriction\sss $\delta\dff |\qff M$\sss 
of\dss $\delta$\sss to\sss $M$\nnsp.\oss\vspace{0.375pt}

\myuppar{Boundary conditions.}
A\qss \emph{boundary condition}\pss
for\sss the elliptic\sss pair\sss $\sigma\fff,\qss \tau$\sss
is\dss defined as a\sss subspace\sss $N$\sss complementary\sss to\sss
$\mathcal{L}_{\dff -}\dff(\trf \rho\dff)$\nnsp,\oss
i.e.\qss such\sss that\sss
$N\dff \cap\dff \mathcal{L}_{\dff -}\dff(\trf \rho\dff)\off =\off 0$\sss
and\sss
$N\qff +\qff \mathcal{L}_{\dff -}\dff(\trf \rho\dff)\off =\off E$\nnsp.\oss
A boundary condition\sss $N$\sss is\dss said\sss to be\qss
\emph{self-adjoint}\pss if\dss $N$\sss is\dss lagrangian.\oss
A\sss lagrangian subspace $N$\sss is\dss a boundary condition\sss if\dss 
either $N\qff +\qff \mathcal{L}_{\dff -}\dff(\trf \rho\dff)\off =\off E$
or\sss $N\dff \cap\dff \mathcal{L}_{\dff -}\dff(\trf \rho\dff)\off =\off 0$\nnsp.\oss
For\sss the rest\sss of\dss this section we will\sss assume\sss that\sss $N$\sss
is\dss a self-adjoint\sss boundary condition\sss for\sss 
$\sigma\fff,\qff \tau$\dnsp.\vspace{0.375pt}

\myuppar{Parameters and\sss deformations.}
In our applications\sss the elliptic pair\sss and\sss the boundary condition\sss
will\sss continuously\sss depend on a parameter\sss $x\qff \in\qff X$\nnsp,\oss
where $X$\sss is\dss a\sss topological\sss space.\oss
Moreover,\oss the vector space\sss $E$\sss will\sss depend on\sss $x\qff \in\qff X$\nnsp.\oss
More precisely,\oss we will\dss be given a family\sss
$E_{\dff x}\dff,\pff x\qff \in\qff X$\sss of\dss vector spaces forming a
vector bundle over\sss $X$\nnsp,\oss
scalar\sss products\sss $\sco{\dff \bullet\dff,\qff \bullet \dff}_{\dff x}\dff,\pff x\qff \in\qff X$\sss
on\sss the vector spaces\sss $E_{\dff x}$\sss defining a\dss Hermitian\sss structure
on\sss this vector bundle,\oss 
self-adjoint\sss operators\sss
$\sigma_{\dff x}\dff,\qff \tau_{\dff x}\dff \colon\dff
E_{\dff x}\qff \ttoo\qff E_{\dff x}$\nnsp,\oss
and subspaces\sss $N_{\dff x}\qff \subset\qff E_{\dff x}$\sss
such\sss that\dss\vspace{1.5pt} 
\[
\quad
E\off =\off E_{\dff x}\dff,\off 
\sigma\off =\off \sigma_{\dff x}\dff,\off 
\tau\off =\off \tau_{\dff x}
\quad
\mbox{and}\quad
N\off =\off N_{\dff x}
\]

\vspace{-12pt}\vspace{1.5pt}
satisfy our assumptions 
for every\sss $x\qff \in\qff X$\sss and continuously depend on $x$\nnsp.\oss
Usually\sss $X$\sss is\dss a compact\sss manifold,\oss
but\sss in\sss this section\sss it\dss is\dss sufficient\sss to assume\sss
that\sss $X$\sss is\dss paracompact.\oss

The rest\sss of\dss this section\dss is\dss devoted\sss to deforming\sss
the elliptic pair\sss $\sigma\fff,\qss \tau$\sss together with\sss the 
boundary condition\sss $N$\sss into a more or\sss less canonical\dss form.\oss
The deformation\sss itself\dss will\sss be sufficiently canonical\sss 
to continuously depend on\sss $x$\sss when a parameter\sss $x\qff \in\qff X$\sss
as above\dss is\dss present.\oss
The deformation will\sss be carried out\sss in several\sss stages,\oss
and\sss this continuity\sss property\sss will\sss be more or\sss less
obvious at\sss each stage.\oss

\myuppar{The first\sss deformation.}
Our first\sss goal\dss is\dss to make\sss
the operator\sss $\sigma$\sss not\sss only self-adjoint,\oss but\sss also unitary,\oss
i.e.\qss such\sss that\sss 
$\sigma\off =\off \sigma^{\fff *}\off =\off \sigma^{\dff -\dff 1}$\dnsp.\oss
Let\sss $\num{\sigma}\off =\off \sqrt{\sigma^{\fff *}\dff \sigma}$\dnsp,\oss
and\dss for\sss 
$\alpha\qff \in\qff [\dff 0\fff,\qff 1/2\trf]$\sss 
let\vspace{1.5pt}
\[
\quad
\sigma_{\dff \alpha}
\off =\off
\num{\sigma}^{\dff -\dff \alpha}\qff \sigma\qff \num{\sigma}^{\dff -\dff \alpha}
\qff,\quad
\tau_{\dff \alpha}
\off =\off
\num{\sigma}^{\dff -\dff \alpha}\qff \tau\qff \num{\sigma}^{\dff -\dff \alpha}
\qff,\quad
\mbox{and}
\]

\vspace{-36pt}
\[
\quad
\rho_{\dff \alpha}
\off =\off
\sigma_{\dff \alpha}^{\dff -\dff 1}\qff \tau_{\dff \alpha}
\off =\off
\num{\sigma}^{\dff \alpha}\qff \rho\qff \num{\sigma}^{\dff -\dff \alpha}
\qff.
\]

\vspace{-12pt}\vspace{1.5pt}
Then\sss the ellipticity assumption\sss holds for\sss
$\sigma_{\dff \alpha}\dff,\qff \tau_{\dff \alpha}$\dss in\sss the role
of\sss $\sigma\fff,\qff \tau$\nnsp.\oss
Indeed,\oss\vspace{1.5pt}
\[
\quad
a\trf \sigma_{\dff \alpha}\qff +\qff b\trf \tau_{\dff \alpha}
\off =\off
\num{\sigma}^{\dff -\dff \alpha}\dff 
(\dff a\trf \sigma\qff +\qff b\trf \tau\trf)\trf 
\num{\sigma}^{\dff -\dff \alpha}\qff 
\]

\vspace{-12pt}\vspace{1.5pt}
is\dss invertible\sss together\sss with\sss
$a\trf \sigma\qff +\qff b\trf \tau$\sss
for every\sss $a\fff,\qff b\qff \in\qff \rrr$\sss such\sss that
$(\trf a\fff,\qff b\trf)\off \neq\off 0$\nnsp.\oss
Let\sss\vspace{1.5pt}
\[
\quad
[\trf u\fff,\qff v\trf]_{\dff \alpha}
\off =\off 
\sco{\dff \sigma_{\dff \alpha}\trf(\trf u\trf)\fff,\qff v\dff}
\qff.
\]

\vspace{-12pt}\vspace{1.5pt}
Then\sss $[\trf \bullet\fff,\qff \bullet\trf]_{\dff \alpha}$\sss
is\dss a non-degenerate\dss Hermitian\dss form on\sss $E$\nnsp,\oss
and\sss $\rho_{\dff \alpha}$\sss is\dss self-adjoint\sss with respect\sss to
$[\trf \bullet\fff,\qff \bullet\trf]_{\dff \alpha}$\nsp.\oss
Clearly,\oss the operators\sss
$\sigma_{\dff \alpha}\dff,\qff \tau_{\dff \alpha}$\nsp,\oss and\sss
$\rho_{\dff \alpha}$\sss continuously depend on\sss $\alpha$\nnsp,\oss
as also\sss the form\sss $[\trf \bullet\fff,\qff \bullet\trf]_{\dff \alpha}$\nsp.\oss
The operators\sss $\rho_{\dff \alpha}$\sss are conjugate\sss to\sss $\rho$\sss
and\sss hence have\sss the same eigenvalues,\oss
and\sss the spaces\sss
$\mathcal{L}_{\dff +}\dff(\trf \rho_{\dff \alpha}\trf)$\sss and\sss
$\mathcal{L}_{\dff -}\dff(\trf \rho_{\dff \alpha}\trf)$\sss
also continuously depend on\sss $\alpha$\nnsp.\oss
For\sss $\alpha\off =\off 1/2$\sss the operator\sss
$\sigma_{\dff \alpha}
\off =\off 
\sigma_{\dff 1/2}
\off =\off
\num{\sigma}^{\dff -\dff 1}\dff \sigma$\sss
is\dss self-adjoint\sss and\sss has only\sss 
$1$\sss and\sss $-\qff 1$\sss as eigenvalues.\oss
It\sss follows\sss that\sss $\sigma_{\dff 1/2}$\sss 
is\dss an\sss involution,\oss
i.e.\dss $\sigma_{\dff 1/2}^{\dff 2}\off =\off 1$\nnsp.\oss
In\sss turn,\oss this implies\sss that\sss 
$\sigma_{\dff 1/2}^{\fff *}\off =\off \sigma_{\dff 1/2}\off =\off \sigma_{\dff 1/2}^{\dff -\dff 1}$\sss
and\sss hence\sss $\sigma_{\dff 1/2}$\sss is\dss a\sss unitary operator.\oss

We need also deform\sss the boundary condition\sss $N$\nnsp.\oss
Let\dss
$N_{\dff \alpha}
\off =\off 
\num{\sigma}^{\dff \alpha}\dff(\trf N\trf)$\nnsp,\qss
$\alpha\qff \in\qff [\trf 0\fff,\qff 1/2\trf]$\nnsp.\oss 
Then\vspace{3pt}
\[
\quad
\bsco{\dff \sigma_{\dff \alpha}^{\vphantom{\alpha}}\trf(\trf N_{\dff \alpha}\trf)\fff,\qff N_{\dff \alpha}\dff}
\off =\off
\bsco{\dff \num{\sigma}^{\dff -\dff \alpha}\qff \sigma\qff \num{\sigma}^{\dff -\dff \alpha}\qff 
\num{\sigma}^{\dff \alpha}\dff(\trf N\trf)\fff,\qff \num{\sigma}^{\dff \alpha}\dff(\trf N\trf)\dff}
\]

\vspace{-35.5pt}
\[
\quad
\phantom{\bsco{\dff \sigma_{\dff \alpha}^{\vphantom{\alpha}}\trf(\trf N_{\dff \alpha}\trf)\fff,\qff 
N_{\dff \alpha}\dff}
\off }
=\off
\bsco{\dff \num{\sigma}^{\dff -\dff \alpha}\qff \sigma\qff (\trf N\trf)\fff,\qff 
\num{\sigma}^{\dff \alpha}\dff(\trf N\trf)\dff}
\off =\off
\bsco{\dff \num{\sigma}^{\dff \alpha}\qff \num{\sigma}^{\dff -\dff \alpha}\qff \sigma\qff 
(\trf N\trf)\fff,\qff 
N\dff}
\]

\vspace{-35.5pt}
\[
\quad
\phantom{\bsco{\dff \sigma_{\dff \alpha}^{\vphantom{\alpha}}\trf(\trf N_{\dff \alpha}\trf)\fff,\qff 
N_{\dff \alpha}\dff}
\off =\off
\bsco{\dff \num{\sigma}^{\dff -\dff \alpha}\qff \sigma\qff (\trf N\trf)\fff,\qff 
\num{\sigma}^{\dff \alpha}\dff(\trf N\trf)\dff}
\off }
=\off
\bsco{\dff \sigma\qff 
(\trf N\trf)\fff,\qff 
N\dff}
\off =\off
0
\qff,
\]

\vspace{-12pt}\vspace{1.5pt}
i.e.\qss
$\sigma_{\dff \alpha}\trf(\trf N_{\dff \alpha}\trf)$\sss
is\dss orthogonal\sss to\sss $N_{\dff \alpha}$\sss
and\sss hence\sss $N_{\dff \alpha}$\sss is\dss lagrangian subspace\sss
with respect\sss to\sss $[\trf \bullet\fff,\qff \bullet\trf]_{\dff \alpha}$\nsp.\oss
A standard\sss verification shows\sss that\sss the subspace\sss $N_{\dff \alpha}$\sss
is\dss transverse\sss to\sss $\mathcal{L}_{\dff -}\dff(\trf \rho_{\dff \alpha}\trf)$\sss
and\sss hence\dss is\dss a self-adjoint\sss boundary condition\sss
for\sss $\sigma_{\dff \alpha}\dff,\qff \tau_{\dff \alpha}$\nsp.\oss

\myuppar{The second deformation.}
Now we would\dss like\sss to deform\sss the operator\sss $\rho$\nnsp,\oss
while keeping\sss its properties,\oss to an operator\sss having only\sss
$i$\sss and\sss $-\qff i$\sss as\sss the eigenvalues.\oss
This deformation\dss is\dss independent\sss of\dss the first\sss one 
and can\sss be done either after or before\sss it.\oss
Let\sss $p_{\dff +}$\sss be\sss the projection of\dss $E$\sss onto\sss
$\mathcal{L}_{\dff +}\dff(\trf \rho\dff)$\sss along\sss
$\mathcal{L}_{\dff -}\dff(\trf \rho\dff)$\nnsp,\oss
and\dss let\sss $p_{\dff -}$\sss be\sss the projection of\dss $E$\sss onto\sss
$\mathcal{L}_{\dff -}\dff(\trf \rho\dff)$\sss along\sss
$\mathcal{L}_{\dff +}\dff(\trf \rho\dff)$\nnsp.\oss
Let\sss $\rho_{\trf 0}\off =\off i\dff p_{\dff +}\qff -\qff i\dff p_{\dff -}$\nsp.\oss
We claim\sss that\sss $\rho_{\trf 0}\dff \colon\dff E\qff \ttoo\qff E$\sss is\dss
self-adjoint\sss with respect\sss to\sss $[\trf \bullet\fff,\qff \bullet\trf]$\nnsp.\oss
Let\sss $u\off =\off u_{\dff +}\qff +\qff u_{\dff -}$\sss
and\sss $v\off =\off v_{\dff +}\qff +\qff v_{\dff -}$\sss
with\vspace{3.5pt}
\[
\quad
u_{\dff +}\dff,\qff v_{\dff +}
\qff \in\qff
\mathcal{L}_{\dff +}\dff(\trf \rho\dff)
\quad
\mbox{and}\quad
u_{\dff -}\dff,\qff v_{\dff -}
\qff \in\qff
\mathcal{L}_{\dff -}\dff(\trf \rho\dff)
\qff.
\]

\vspace{-12pt}\vspace{3.75pt}
Then\quad 
$[\trf \rho_{\trf 0}\trf(\dff u\trf)\fff,\qff v\trf]
\off =\off
[\trf i\dff u_{\dff +}\qff -\qff i\dff u_{\dff -}\fff,\qff v\trf]
\off =\off
[\trf i\dff u_{\dff +}\fff,\qff v\trf]
\qff -\qff
[\trf i\dff u_{\dff -}\fff,\qff v\trf]$\vspace{3.75pt}
\[
\quad
=\off
[\trf i\dff u_{\dff +}\fff,\qff v_{\dff +}\qff +\qff v_{\dff -}\trf]
\qff -\qff
[\trf i\dff u_{\dff -}\fff,\qff v_{\dff +}\qff +\qff v_{\dff -}\trf]
\off =\off
[\trf i\dff u_{\dff +}\fff,\qff v_{\dff -}\trf]
\qff -\qff
[\trf i\dff u_{\dff -}\fff,\qff v_{\dff +}\trf]
\]

\vspace{-33pt}\vspace{0.5pt}
\[
\quad
=\off
[\trf u_{\dff +}\fff,\qff -\qff i\dff v_{\dff -}\trf]
\qff +\qff
[\trf u_{\dff -}\fff,\qff i\dff v_{\dff +}\trf]
\off =\off
[\trf u_{\dff +}\qff +\qff u_{\dff -}\fff,\qff i\dff v_{\dff +}\qff -\qff i\dff v_{\dff -}\trf]
\off =\off
[\trf u\fff,\qff \rho_{\trf 0}\trf(\dff v\trf)\trf]
\qff.
\]

\vspace{-12pt}\vspace{3.5pt}
This proves our\sss claim.\oss
Clearly,\oss the only eigenvalues of\dss $\rho_{\trf 0}$\sss
are\sss $i$\sss and\sss $-\qff i$\nnsp.\oss
Let\sss $\tau_{\dff 0}\off =\off \sigma\dff \rho_{\trf 0}$\nsp.\oss
Then\sss
$\sco{\dff \tau_{\dff 0}\trf(\trf u\trf)\fff,\qff v\dff}
\off =\off
\sco{\dff \sigma\dff \rho_{\trf 0}\trf(\trf u\trf)\fff,\qff v\dff}
\off =\off
[\trf \rho_{\trf 0}\trf(\dff u\trf)\fff,\qff v\trf]$\sss
and\dss
$\sco{\dff u\fff,\qff \tau_{\dff 0}\trf(\trf v\trf)\dff}
\off =\off
[\trf u\fff,\qff \rho_{\trf 0}\trf(\dff v\trf)\trf]$\nnsp.\oss
It\sss follows\sss that\sss $\tau_{\dff 0}$\sss is\dss self-adjoint\sss
with respect\sss to\sss
$\sco{\dff \bullet\fff,\qff \bullet\dff}$\nnsp.\oss
Let\sss us\dss connect\sss $\rho\fff,\qff \tau$\sss with\sss
$\rho_{\trf 0}\dff,\qff \tau_{\dff 0}$\sss respectively\sss
by\sss linear\sss homotopies\vspace{1.5pt}\vspace{1.5pt}
\[
\quad
\rho_{\trf \beta}
\off =\off
(\dff 1\qff -\qff \beta\trf)\trf \rho
\qff +\qff
\beta\trf \rho_{\trf 0}
\quad
\mbox{and}\quad
\tau_{\trf \beta}
\off =\off
(\dff 1\qff -\qff \beta\trf)\trf \tau
\qff +\qff
\beta\trf \tau_{\dff 0}
\qff,
\]

\vspace{-12pt}\vspace{1.5pt}\vspace{1.5pt}
where\sss $\beta\qff \in\qff [\dff 0\fff,\qff 1\trf]$\nnsp.\oss
Then\sss 
$\rho_{\trf \beta}\off =\off \sigma^{\dff -\dff 1}\dff \tau_{\trf \beta}$\sss
for every\sss $\beta$\nnsp.\oss
Clearly,\vspace{1.5pt}\vspace{1.5pt}
\[
\quad
\mathcal{L}_{\dff +}\dff(\trf \rho_{\trf \beta}\dff)
\off =\off
\mathcal{L}_{\dff +}\dff(\trf \rho\dff)
\quad
\mbox{and}\quad
\mathcal{L}_{\dff -}\dff(\trf \rho_{\trf \beta}\dff)
\off =\off
\mathcal{L}_{\dff -}\dff(\trf \rho\dff)
\]

\vspace{-12pt}\vspace{1.5pt}\vspace{1.5pt}
for every $\beta$\nnsp.\oss
The ellipticity assumption\sss holds for\sss 
$\sigma\fff,\qff \tau_{\trf \beta}$\sss and every $\beta$\nnsp.\oss
Indeed,\oss\vspace{1.5pt}
\[
\quad
a\trf \sigma\qff +\qff b\trf \tau_{\dff \beta}
\off =\off
\sigma\dff 
(\dff a\qff +\qff b\trf \rho_{\trf \beta}\trf)
\]

\vspace{-12pt}\vspace{1.5pt}
is\dss invertible if\sss $a\fff,\qff b\qff \in\qff \rrr^{\dff 2}\qff \smallsetminus\qff \{\trf 0\trf\}$\sss 
because $\sigma$ is\dss invertible and\sss $\rho_{\trf \beta}$\sss
has no real\sss eigenvalues.\oss

This deformation does not\sss affect\sss the scalar product\sss
$[\trf \bullet\fff,\qff \bullet\trf]$\sss and\sss the space
$\mathcal{L}_{\dff -}\dff(\trf \rho\dff)$\nnsp.\oss
Therefore $N$\sss remains a boundary condition during\sss the second deformation\sss
$\sigma_{\dff \beta}$\nsp,\dss $\tau_{\dff \beta}$\nsp,\dss $\rho_{\trf \beta}$\nsp,\qss
$\beta\qff \in\qff [\dff 0\fff,\qff 1\trf]$\nnsp.\oss

\myuppar{Lagrangian subspaces when\sss the operator\sss $\sigma$\sss is\dss unitary.}
Suppose\sss that\sss the operator\sss $\sigma$\sss is\sss self-adjoint\sss and\sss unitary.\oss 
Let\sss $L$\sss be a\sss lagrangian subspace.\oss
Since\sss $\sigma$\sss is\dss unitary,\pss $\sigma\dff(\trf L\trf)$\sss is\dss also\sss lagrangian.\oss
Indeed,\oss the\sss lagrangian property of\dss $L$\sss means\sss that\sss
$\sigma\dff(\trf L\trf)$\sss is\dss orthogonal\dss to $L$ with respect\sss to
$\sco{\dff \bullet\dff,\qff \bullet \dff}$\nnsp.\oss
Since\sss the operator\dss $\sigma$\sss is\dss unitary,\oss this implies\sss that\sss
$\sigma\dff(\trf \sigma\dff(\trf L\trf)\trf)\off =\off L$\sss
is\dss orthogonal\sss to\sss $\sigma\dff(\trf L\trf)$\sss with respect\sss to
$\sco{\dff \bullet\dff,\qff \bullet \dff}$\nnsp,\oss
and\sss hence\sss the space\sss $\sigma\dff(\trf L\trf)$\sss is\dss lagrangian.\oss

Since\sss $\sigma$\sss is\dss self-adjoint\sss and unitary,\pss
$E$ can be decomposed\sss into a direct\sss sum\sss
$E\off =\off E^{\dff +}\dff \oplus\qff E^{\dff -}$\nsp,\oss
where\sss $E^{\dff +}$\nsp,\oss and\sss $E^{\dff -}$ are\sss
the eigenspaces of\sss $\sigma$\sss corresponding\sss to\sss
the eigenvalues $1$ and $-\qff 1$ respectively.\oss
Since $\sigma$\sss is\dss self-adjoint,\oss the spaces\sss
$E^{\dff +}$\nsp,\oss and\sss $E^{\dff -}$ are\sss
orthogonal\sss with respect\sss to\sss
$\sco{\dff \bullet\dff,\qff \bullet \dff}$\nnsp.\oss
Clearly,\oss the scalar\sss product\sss $[\trf \bullet\fff,\qff \bullet\trf]$\sss
is\dss equal\sss to\sss
$\sco{\dff \bullet\dff,\qff \bullet \dff}$\sss
on\sss $E^{\dff +}$\sss and\sss to\sss
$-\qff \sco{\dff \bullet\dff,\qff \bullet \dff}$\sss on\sss $E^{\dff -}$\dnsp.\oss
In\sss the orthogonal\sss decomposition\sss
$E\off =\off E^{\dff +}\dff \oplus\dff E^{\dff -}$\sss
the operator\sss $\sigma$\sss takes\sss the matrix\sss form\vspace{1.5pt}
\begin{equation}
\label{sigma-st}
\quad
\sigma
\off =\off\dff
\begin{pmatrix}
\off\dff 1 &
0 \off
\vspace{4.5pt} \\
\off\dff 0 &
-\qff 1 \off 
\end{pmatrix}
\off.
\end{equation}

\vspace{-12pt}\vspace{1.5pt}
The\sss lagrangian subspace\sss $L$\sss 
is\dss transverse\sss to both 
$E^{\dff +}$\sss and\sss $E^{\dff -}$\nsp,\oss
i.e.\dss 
$L\dff \cap\dff E^{\dff +}
\off =\off
L\dff \cap\dff E^{\dff -}
\off =\off
0$\nnsp.\oss
This implies\sss that\sss the projection of\dss $L$\sss onto $E^{\dff +}$
along\sss $E^{\dff -}$\sss is\dss injective,\oss 
as also\sss the projection onto $E^{\dff -}$
along\sss $E^{\dff +}$\dnsp.\oss
Since\sss the dimension of\dss $L$\sss is\dss equal\sss to\sss the half\dss
of\dss the dimension of\dss $E$\nnsp,\oss it\sss follows\sss that\sss the 
dimensions of\dss $E^{\dff +}$ and\sss $E^{\dff -}$ are also equal\sss to\sss the half\dss
of\dss the dimension of\dss $E$\nnsp.\oss
Also,\qss $L$\sss is\dss equal\sss to\sss the graph\sss
$\{\trf u\qff +\qff \varphi\dff(\dff u\trf)\qff |\qff u\qff \in\qff E^{\dff +}\qff\}$\sss
of\dss a unique\sss linear\sss map\sss
$\varphi\dff \colon\dff E^{\dff +}\qff \ttoo\qff E^{\dff -}$\dnsp.\oss
We claim\sss that\sss $\varphi$\sss is\dss an\dss isometry\sss
with respect\sss to\sss $\sco{\dff \bullet\dff,\qff \bullet \dff}$\nnsp.\oss
Indeed,\oss if\dss $u\fff,\qff v\qff \in\qff E^{\dff +}$\dnsp,\oss then\vspace{3pt}\vspace{-0.25pt}
\[
\quad
0
\off =\off
[\trf u\qff +\qff \varphi\dff(\dff u\trf)\fff,\qff v\qff +\qff \varphi\dff(\dff v\trf)\trf]
\off =\off
\sco{\dff \sigma\trf(\dff u\qff +\qff \varphi\dff(\dff u\trf)\trf)\fff,\qff 
v\qff +\qff \varphi\dff(\dff v\trf) \dff}
\]

\vspace{-34.75pt}\vspace{-0.25pt}
\[
\quad
=\off
\sco{\dff u\qff -\qff \varphi\dff(\dff u\trf)\fff,\qff 
v\qff +\qff \varphi\dff(\dff v\trf) \dff}
\off =\off
\sco{\dff u\fff,\qff v \dff}
\qff -\qff
\sco{\dff \varphi\dff(\dff u\trf)\fff,\qff v \dff}
\qff +\qff
\sco{\dff u\fff,\qff \varphi\dff(\dff v\trf) \dff}
\qff -\qff
\sco{\dff \varphi\dff(\dff u\trf)\fff,\qff \varphi\dff(\dff v\trf) \dff}
\qff.
\]

\vspace{-12pt}\vspace{3pt}\vspace{-0.25pt}
The spaces\sss $E^{\dff +}$ and\sss $E^{\dff -}$ are orthogonal\sss
with respect\sss to\sss the product $\sco{\dff \bullet\dff,\qff \bullet \dff}$\sss 
and\sss therefore\sss
$\sco{\dff \varphi\dff(\dff u\trf)\fff,\qff v \dff}
\off =\off
\sco{\dff u\fff,\qff \varphi\dff(\dff v\trf) \dff}
\off =\off
0$\nnsp.\oss
Hence\sss 
$\sco{\dff \varphi\dff(\dff u\trf)\fff,\qff \varphi\dff(\dff v\trf) \dff}
\off =\off
\sco{\dff u\fff,\qff v \dff} $\sss
for every\sss $u\fff,\qff v\qff \in\qff E^{\dff +}$\nsp,\oss
i.e.\dss $\varphi$\sss is\dss an\sss isometry.\oss
If\dss $L'$\sss is\dss another\sss lagrangian subspace and\sss
$\varphi'$\sss is\dss the corresponding\sss isometry,\oss
then\sss $L'$\sss is\dss transverse\sss to\sss $L$\sss
if\dss and\sss only\dss if\dss
$\varphi\dff(\dff u\trf)
\off \neq\off
\varphi'\dff(\dff u\trf)$\sss
for every\sss nonzero\sss $u\qff \in\qff E^{\dff +}$\dnsp.\oss

Let\sss us use\sss the isometry\sss $\varphi$\sss to identify $E^{\dff -}$
with\sss $E^{\dff +}$\dnsp.\oss
This identification\sss turns\sss $L$\sss into\sss the diagonal\sss
$\Delta
\off =\off 
\{\trf (\dff u\fff,\qff u\trf)\qff |\qff u\qff \in\qff E^{\dff +}\qff\}$\nnsp.\oss
Let\sss us denote by $F^{\trf \perp}$
the orthogonal\sss complement\sss of\dss $F\qff \subset\qff E$ with respect\sss to\sss 
$\sco{\dff \bullet\dff,\qff \bullet \dff}$\dnsp.\oss
Then\sss
$\Delta^{\fff \perp}
\off =\off 
\{\trf (\dff u\fff,\qff -\qff u\trf)\qff |\qff u\qff \in\qff E^{\dff +}\qff\}$\sss
is\dss also a\sss lagrangian subspace.\oss

\myuppar{The\sss third deformation.}
After\sss two deformations\sss the operator\sss $\sigma$\sss
is\dss self-adjoint\sss and\sss unitary,\oss and\sss $\rho$\sss has only\sss
$i$\sss and\sss $-\qff i$\sss as eigenvalues.\oss
In\sss particular,\pss $\rho$\sss is\dss uniquely determined\sss by\sss
the spaces\sss
$\mathcal{L}_{\dff +}\dff(\trf \rho\trf)$\sss and\sss
$\mathcal{L}_{\dff -}\dff(\trf \rho\trf)$\nnsp.\oss
In\sss the\sss third deformation\sss we will\sss use\sss the boundary condition\sss $N$\sss
to bring\sss $\mathcal{L}_{\dff -}\dff(\trf \rho\trf)$\sss
into a standard\sss form.\oss
To begin with,\oss let\sss 
$\varphi\dff \colon\dff E^{\dff +}\qff \ttoo\qff E^{\dff -}$\sss
be\sss the isometry corresponding\sss to\sss the\sss lagrangian subspace\sss $N$\nnsp.\oss
Let\sss us\sss identify $E^{\dff +}$ with\sss $E^{\dff -}$\sss by\sss the isometry\sss
$\varphi$\nnsp,\oss
and\sss let\sss $F\off =\off E_{\dff +}\off =\off E_{\dff -}$\nsp.\oss
With\sss these identifications\sss $N$\sss turns into\sss the diagonal\sss
$\Delta\off =\off \{\trf (\dff u\fff,\qff u\trf) \qff |\qff u\qff \in\qff F \trf\}$\sss
and $\varphi$\sss into\sss the identity map\sss
$\id_{\trf F}\dff \colon\dff F\qff \ttoo\qff F$\nnsp.\oss

The\sss lagrangian subspace\sss $\mathcal{L}_{\dff -}\dff(\trf \rho\dff)$\sss
corresponds\sss to\sss some\sss isometry\sss 
$\psi_{\dff -}\dff \colon\dff F\qff \ttoo\qff F$\dnsp.\oss
As we saw,\oss the\sss transversality of\sss $\mathcal{L}_{\dff -}\dff(\trf \rho\dff)$\sss
and\sss $N$\sss implies\sss that\sss 
$\psi_{\dff -}\dff(\dff u\trf)\off \neq\off u$\sss
for every\sss $u\qff \in\qff F$\dnsp.\oss
In other\sss terms,\pss $1$\sss is\dss not\sss an eigenvalue of\sss $\psi_{\dff -}$\nsp.\oss
Therefore\sss there\dss is\dss a canonical\sss spectral\sss deformation 
$\psi_{\dff -}\dff(\trf \alpha\trf)$\nsp,\qss
$\alpha\qff \in\qff [\dff 0\fff,\qff 1\dff]$\sss in\sss the class of\dss isometries
not\sss having $1$ as an eigenvalue,\oss such\sss that\sss
$\psi_{\dff -}\off =\off \psi_{\dff -}\dff(\dff 0\dff)$\sss and\sss
$\psi_{\dff -}\dff(\dff 1\dff)$\sss has only\sss $-\qff 1$\sss as an eigenvalue,\oss
i.e.\qss $\psi_{\dff -}\dff(\dff 1\dff)\off =\off -\qff \id_{\trf F}$\nsp.\oss
The corresponding\sss lagrangian subspace\dss is\dss the anti-diagonal\sss 
$\Delta^{\fff \perp}
\off =\off
\{\trf (\dff u\fff,\qff -\qff u\trf)\qff |\qff u\qff \in\qff F\qff\}$\nnsp.\oss

Let\sss $L_{\dff -}\dff(\trf \alpha\trf)$\sss be\sss the\sss lagrangian subspace
corresponding\sss to\sss $\psi_{\dff \alpha}$\nsp.\qss
Then\sss $L_{\dff -}\dff(\trf \alpha\trf)$\sss 
is\dss transverse\sss to\sss $N$\sss for every $\alpha$\nnsp.\oss 
Also,\pss 
$L_{\dff +}\dff(\trf 0\trf)\off =\off \mathcal{L}_{\dff +}\dff(\trf \rho\dff)$\sss 
is\dss transverse\sss to\sss 
$L_{\dff -}\dff(\trf 0\trf)\off =\off \mathcal{L}_{\dff -}\dff(\trf \rho\dff)$\nnsp.\oss
By\sss the discussion of\dss the spaces\sss
$\Lambda_{\fff L}$\sss\sss above,\oss 
for each\sss $\alpha\qff \in\qff [\dff 0\fff,\qff 1\trf]$\sss
the space of\dss lagrangian subspaces\sss transverse\sss to\sss
$L_{\dff -}\dff(\trf \alpha\trf)$\sss 
is\dss an affine space over\sss the vector\sss space of\dss
linear\sss maps\sss 
$\delta\dff \colon\dff 
E\qff \ttoo\qff L_{\dff -}\dff(\trf \alpha\trf)$\sss
such\sss that\qss (\ref{delta})\qss holds.\oss 
These\sss affine spaces form a\sss locally\sss trivial\dss bundle 
over\sss $[\dff 0\fff,\qff 1\trf]$\dss
(since\sss $[\dff 0\fff,\qff 1\trf]$\sss is\dss contractible,\oss
this bundle\dss is\dss actually\sss trivial\fff).\oss
By\sss the homotopy\sss lifting\sss property\sss the deformation\sss
$L_{\dff -}\dff(\trf \alpha\trf)$\nsp,\dss
$\alpha\qff \in\qff [\dff 0\fff,\qff 1\trf]$\sss
can\sss be accompanied\dss by\sss a deformation\sss $L_{\dff +}\dff(\trf \alpha\trf)$\nsp,\qss
$\alpha\qff \in\qff [\dff 0\fff,\qff 1\trf]$\sss starting\sss
with\sss $L_{\dff +}\dff(\trf 0\trf)$\sss and such\sss that\sss
$L_{\dff +}\dff(\trf \alpha\trf)$\sss
is\dss transverse\sss to\sss $L_{\dff -}\dff(\trf \alpha\trf)$\sss
for every\sss $\alpha$\nnsp.\oss
Moreover,\oss the deformation\sss $L_{\dff +}\dff(\trf \alpha\trf)$\nsp,\qss
$\alpha\qff \in\qff [\dff 0\fff,\qff 1\trf]$\sss
is\dss unique up\sss to\sss homotopy\qss
(between\sss homotopies)\qss when\sss parameters are present.\oss
Cf.\dss the discussion of\dss loops and\dss bundles in\qss \cite{i1},\oss Section\qss 16.\oss
It\sss follows\sss that\sss the $L_{\dff +}\dff(\trf 1\trf)$\sss
is\dss unique up\sss to a homotopy\sss when\sss parameters are present.\oss
It\dss is\dss worth\sss to point\sss out\sss
that\sss the deformation $L_{\dff +}\dff(\trf \alpha\trf)$\nsp,\qss
$\alpha\qff \in\qff [\dff 0\fff,\qff 1\trf]$\sss
is\qss \emph{not\sss a spectral\sss deformation}\pss of\dss the isometry corresponding\sss
$\mathcal{L}_{\dff +}\dff(\trf \rho\dff)$\nnsp.\oss
Let\sss
$\rho_{\dff \alpha}\dff \colon\dff E\qff \ttoo\qff E$\sss
be\sss the\sss linear\sss map such\sss that\sss\vspace{1.5pt}\vspace{0.125pt}
\[
\quad
\rho_{\dff \alpha}\dff(\dff u\trf)\off =\off i\dff u
\quad
\mbox{for}\quad
u\qff \in\qff L_{\dff +}\dff(\trf \alpha\trf)
\quad
\mbox{and}\quad
\]

\vspace{-36pt}
\[
\quad
\rho_{\dff \alpha}\dff(\dff u\trf)\off =\off -\qff i\dff u
\quad
\mbox{for}\quad
u\qff \in\qff L_{\dff -}\dff(\trf \alpha\trf)
\qff.
\]

\vspace{-12pt}\vspace{1.5pt}\vspace{0.125pt}
Let\sss $\tau_{\dff \alpha}\off =\off \sigma\dff \rho_{\dff \alpha}$\nsp.\oss
As in\sss the discussion of\dss the second\sss deformation,\oss
we see\sss that\sss $\rho_{\dff \alpha}$\sss is\dss self-adjoint\sss 
with respect\sss to\sss $[\trf \bullet\fff,\qff \bullet\trf]$\sss
and\sss hence\sss $\tau_{\dff \alpha}$\sss is\dss self-adjoint\sss 
with respect\sss to\sss $\sco{\dff \bullet\dff,\qff \bullet \dff}$\sss 
for every $\alpha$\nnsp.\oss
Clearly,\pss
$\mathcal{L}_{\dff +}\dff(\trf \rho_{\dff \alpha}\dff)
\off =\off
L_{\dff +}\dff(\trf \alpha\trf)$\sss
and\dss
$\mathcal{L}_{\dff -}\dff(\trf \rho_{\dff \alpha}\dff)
\off =\off
L_{\dff -}\dff(\trf \alpha\trf)$\sss
and\sss that\sss $\rho_{\dff \alpha}$\sss has only\sss
$i$\sss and\sss $i-\qff $\sss as eigenvalues\sss for every $\alpha$\nnsp.\oss
Since\sss $N$\sss is\dss transverse\sss to\sss
$\mathcal{L}_{\dff -}\dff(\trf \rho_{\dff \alpha}\dff)$\sss
by\sss the construction,\pss
$N$\sss is\dss a boundary condition\sss for\sss
$\sigma\fff,\qff \tau_{\dff \alpha}$\sss for every $\alpha$\sss
and stays intact\sss during\sss this deformation.\oss
To sum up,\oss at\sss the end of\dss the\sss third deformation\sss
$\sigma$\sss has\sss the standard\sss form\qss (\ref{sigma-st}),\oss
the boundary condition\sss $N$\sss is\dss equal\sss to\sss the diagonal\sss $\Delta$\nnsp,\oss
and\sss the\sss eigenspace\sss $\mathcal{L}_{\dff -}\dff(\trf \rho\dff)$\sss
is\dss equal\sss to\sss $\Delta^{\fff \perp}$\dnsp.\oss

\myuppar{The fourth deformation.}
It\dss is\dss similar\sss to\sss the\sss third one.\oss
Now we will\dss bring\sss the space\sss 
$\mathcal{L}_{\dff +}\dff(\trf \rho\dff)$\sss
to a standard\sss form.\oss
After\sss the\sss third deformation\sss it\dss is\dss 
a\sss lagrangian subspace\sss transverse\sss to\sss 
$\mathcal{L}_{\dff -}\dff(\trf \rho\dff)\off =\off \Delta^{\fff \perp}$\dnsp.\oss
The\sss lagrangian subspace\sss
$N\off =\off \Delta$\sss is\dss also\sss transverse\sss to\sss
$\mathcal{L}_{\dff -}\dff(\trf \rho\dff)\off =\off \Delta^{\fff \perp}$\dnsp.\oss
Since\sss the space of\dss lagrangian subspaces\sss
transverse\sss to a given one,\oss say\sss to\sss $\Delta^{\fff \perp}$\dnsp,\oss
is\dss an affine space,\oss
we can deform\sss $\mathcal{L}_{\dff +}\dff(\trf \rho\dff)$\sss to\sss $N$\sss
by a\sss linear deformation\sss in\sss this affine space.\oss
Moreover,\oss this\dss is\dss a canonical\sss deformation.\oss
As in\sss the\sss third\sss deformation,\oss
this deformation of\sss $\mathcal{L}_{\dff +}\dff(\trf \rho\dff)$\sss
defines a deformation of\sss $\rho$ and\sss $\tau$\dnsp,\oss
and\sss $N$\sss remains a boundary condition during\sss this deformation.\oss
The fourth deformation\sss brings\sss $\sigma\fff,\off \rho\fff,\off \tau$\sss
and\sss $N$\sss into a normal\sss form.\oss
As we already saw,\pss $\sigma$\sss has\sss the\sss form\qss (\ref{sigma-st}).\oss
Let\sss us\sss treat\sss $\varphi$ as\sss the identification\sss map\sss
$E^{\dff +}\qff \ttoo\qff E^{\dff -}$\nsp.\oss
Then\sss $N\off =\off \Delta$\nnsp.\oss
Since\sss $\Delta$\sss and\sss $\Delta^{\fff \perp}$\sss
are\sss the eigenspaces of\sss $\rho$\sss with\sss
eigenvalues\sss $i$\sss and\sss $-\qff i$\sss respectively,\oss
one can easily compute\sss $\rho$ and\sss $\tau$\nnsp.\oss
We see\sss that\vspace{3pt}\vspace{-0.75pt}
\begin{equation}
\label{sigma-rho-tau-st}
\quad
\sigma
\off =\off\dff
\begin{pmatrix}
\off\dff 1 &
0 \off
\vspace{4.5pt} \\
\off\dff 0 &
-\qff 1 \off 
\end{pmatrix}
\qff,\quad
\rho
\off =\off\dff
\begin{pmatrix}
\off 0 &
i \qff\off
\vspace{4.5pt} \\
\off\dff i &
0 \qff\off 
\end{pmatrix}
\quad
\mbox{and}\quad
\tau
\off =\off\dff
\begin{pmatrix}
\off 0 &
i \off
\vspace{4.5pt} \\
\off -\qff i &
0 \off 
\end{pmatrix}
\off.
\end{equation}

\vspace{-12pt}\vspace{3pt}\vspace{-0.75pt}
\myuppar{The normal\dss form.}
When\sss $\sigma\fff,\off \rho\fff,\off \tau$\sss have\sss the form\qss 
(\ref{sigma-rho-tau-st})\qss and\sss $N\off =\off \Delta$\nnsp,\oss
we will\sss say\sss that\sss
$\sigma\fff,\off \rho\fff,\off \tau$ and\sss $N$\sss
are\qss \emph{in\dss the standard\dss form}.\oss 
When everything depends on parameters,\oss
it\dss may\sss happen\sss that\sss the vector spaces\sss
$E^{\dff +}\dff,\off E^{\dff -}$\sss
depend on a parameter,\oss
or\sss that\sss 
$E^{\dff +}\dff,\off E^{\dff -}$\sss
do not\sss depend on a parameter,\oss
but\sss 
$\varphi$\sss does.\oss
In\sss this case\sss it\dss is\dss
not\sss possible\sss to identify\dss 
$E^{\dff +}$\sss with\sss $E^{\dff -}$\sss
by\sss $\varphi$\dss and one has\sss to use
a normal\dss form explicitly\sss involving\sss $\varphi$\nnsp.\oss
For every\sss isometry\sss
$\varphi\dff \colon\dff E^{\dff +}\qff \ttoo\qff E^{\dff -}$\dss
let\vspace{3pt}
\begin{equation}
\label{diagonals}
\quad
\Delta\dff(\trf \varphi\trf)
\off =\off
\bigl\{\qff
(\trf u\fff,\qff \varphi\dff(\dff u\trf)\qff)
\qff \bigl|\qff
u\qff \in\qff E^{\dff +}
\trf\bigr\}
\quad
\mbox{and}\quad
\Delta^{\fff \perp}\dff(\trf \varphi\trf)
\off =\off
\bigl\{\qff
(\trf u\fff,\qff -\qff \varphi\dff(\dff u\trf)\qff)
\qff \bigl|\qff
u\qff \in\qff E^{\dff +}
\trf\bigr\}
\qff.\quad
\end{equation}

\vspace{-12pt}\vspace{3pt}
If\dss $\sigma$ has\sss the form\qss (\ref{sigma-st})\qss
with respect\sss to\sss the decomposition\sss
$E\off =\off E^{\dff +}\dff \oplus\dff E^{\dff -}$\sss
and\sss 
$\Delta\dff(\trf \varphi\trf)\dff,\off 
\Delta^{\fff \perp}\dff(\trf \varphi\trf)$\sss
are\sss the eigenspaces of\dss $\rho$\sss with\sss 
the eigenvalues\sss $i\fff,\off -\qff i$\sss respectively,\oss
then\vspace{3pt}
\begin{equation}
\label{sigma-rho-tau-phi}
\quad
\sigma
\off =\off\dff
\begin{pmatrix}
\off\dff 1 &
0 \off
\vspace{4.5pt} \\
\off\dff 0 &
-\qff 1 \off 
\end{pmatrix}
\qff,\quad
\rho
\off =\off\dff
\begin{pmatrix}
\off 0 &
i\trf \varphi^{\fff *} \off
\vspace{4.5pt} \\
\off\dff i\trf \varphi &
0 \off 
\end{pmatrix}
\qff,\quad
\mbox{and}\quad
\tau
\off =\off\dff
\begin{pmatrix}
\off 0 &
i\trf \varphi^{\fff *} \off
\vspace{4.5pt} \\
\off -\qff i\trf \varphi &
0 \off 
\end{pmatrix}
\off,
\end{equation}

\vspace{-12pt}\vspace{3pt}
as a simple computation,\oss taking\sss into account\sss that\sss
$\varphi^{\fff *}\off =\off \varphi^{\dff -\dff 1}$\dnsp,\oss
shows.\oss
If\dss also\vspace{3pt}
\begin{equation}
\label{delta-delta}
N
\off =\off
\Delta\dff(\trf \varphi\trf)
\off =\off
\mathcal{L}_{\dff -}\dff(\trf \rho\trf)
\quad
\mbox{and\dss hence}\quad
N^{\dff \perp}
\off =\off
\Delta^{\fff \perp}\dff(\trf \varphi\trf)
\off =\off
\mathcal{L}_{\dff +}\dff(\trf \rho\trf)
\qff,
\end{equation}

\vspace{-12pt}\vspace{3pt}
we will\sss say\sss that\sss
$\sigma\fff,\off \rho\fff,\off \tau$ and\sss $N$\sss
are\qss \emph{in\sss the normal\dss form},\oss or are\qss \emph{normalized}.\oss
Equivalently,\pss
$\sigma\fff,\off \rho\fff,\off \tau$ and\sss $N$\sss
are in\sss the normal\dss form\dss if\dss
$\sigma$\sss is\dss a unitary,\pss
$N\off =\off \mathcal{L}_{\dff +}\dff(\trf \rho\dff)$\nnsp,\pss
$N^{\dff \perp}\off =\off \mathcal{L}_{\dff -}\dff(\trf \rho\dff)$\nnsp,\oss
and\sss the spaces\sss 
$\mathcal{L}_{\dff +}\dff(\trf \rho\dff)\fff,\off
\mathcal{L}_{\dff -}\dff(\trf \rho\dff)$\sss
are\sss the eigenspaces of\dss $\rho$\sss with\sss 
the eigenvalues\sss $i\fff,\off -\qff i$\sss respectively.\oss
Clearly,\oss in\sss this case $\sigma$ anti-commutes with $\tau$ and\sss $\rho$\nnsp.\oss

\myuppar{Diagonal\dss standard\dss form.}
Clearly,\vspace{3pt}
\[
\quad
\sigma\dff(\trf u\fff,\qff \varphi\dff(\dff u\trf)\trf)
\off =\off
(\trf u\fff,\qff -\qff \varphi\dff(\dff u\trf)\trf)
\qff,\qquad
\sigma\dff(\trf u\fff,\qff -\qff \varphi\dff(\dff u\trf)\trf)
\off =\off
(\trf u\fff,\qff \varphi\dff(\dff u\trf)\trf)
\qff,
\]

\vspace{-34.5pt}
\[
\quad
\rho\dff(\trf u\fff,\qff \varphi\dff(\dff u\trf)\trf)
\off =\off
(\trf i\trf \varphi^{\dff -\dff 1}\fff \circ\qff \varphi\dff(\dff u\trf)\fff,\qff
i\trf \varphi(\dff u\trf)\qff)
\off =\off
i\trf (\trf u\fff,\qff \varphi\dff(\dff u\trf)\qff)
\qff,\quad
\mbox{and}\quad
\]

\vspace{-34.5pt}
\[
\quad
\rho\dff(\trf u\fff,\qff -\qff \varphi\dff(\dff u\trf)\trf)
\off =\off
(\trf -\qff i\trf \varphi^{\dff -\dff 1}\fff \circ\qff \varphi\dff(\dff u\trf)\fff,\qff
i\trf \varphi(\dff u\trf)\qff)
\off =\off
-\qff i\trf (\trf u\fff,\qff -\qff \varphi\dff(\dff u\trf)\qff)
\qff.
\]

\vspace{-12pt}\vspace{3pt}
One can always\sss identify\sss 
$\Delta\dff(\trf \varphi\trf)$ with $\Delta^{\fff \perp}(\trf \varphi\trf)$\sss
by\sss the map\sss
$(\dff u\fff,\qff v\trf)
\off \longmapsto\off
(\dff u\fff,\qff -\qff v\trf)$\nnsp.\oss
In\sss terms of\dss this identification 
and\sss the resulting\sss decomposition\sss
$E
\off =\off 
\Delta\dff(\trf \varphi\trf)
\dff \oplus\dff 
\Delta\dff(\trf \varphi\trf)$\sss
we have\vspace{3pt}
\begin{equation}
\label{sigma-rho-tau-alt}
\quad
\sigma
\off =\off\dff
\begin{pmatrix}
\off 0 &
1 \qff\off
\vspace{4.5pt} \\
\off\dff 1 &
0 \qff\off 
\end{pmatrix}
\qff,\quad
\rho
\off =\off\dff
\begin{pmatrix}
\off i &
0 \off
\vspace{4.5pt} \\
\off 0 &
-\qff i \off 
\end{pmatrix}
\qff,\quad
\mbox{and}\quad
\tau
\off =\off\dff
\begin{pmatrix}
\off 0 &
-\qff i \off
\vspace{4.5pt} \\
\off i &
0 \off 
\end{pmatrix}
\off,
\end{equation}

\vspace{-12pt}\vspace{3pt}
which\dss is\dss often more convenient\sss than\qss (\ref{sigma-rho-tau-st}).\oss
We will\sss call\pss (\ref{sigma-rho-tau-alt})\qss 
\emph{the diagonal\dss standard\dss form}.\oss

\myuppar{Odd\sss graded\dss elliptic pairs.}
Suppose\sss that\sss for some orthogonal\sss decomposition\sss
$E
\off =\off 
F\dff \oplus\trf F$\sss
and a\sss linear\sss map\sss
$\bm{\tau}\dff \colon\dff
F\qff \ttoo\qff F$\sss
the operators\sss \sss $\sigma\fff,\off \tau\fff,\off \rho$\sss
have\sss the form\vspace{3pt}
\begin{equation}
\label{sigma-rho-tau-phi-alt}
\quad
\sigma
\off =\off\dff
\begin{pmatrix}
\off 0 &
1 \qff\off
\vspace{4.5pt} \\
\off\dff 1 &
0 \qff\off 
\end{pmatrix}
\qff,\quad
\rho
\off =\off\dff
\begin{pmatrix}
\off \bm{\tau} &
0 \off
\vspace{4.5pt} \\
\off 0 &
\bm{\tau}^{\dff *} \off 
\end{pmatrix}
\quad
\mbox{and}\quad
\tau
\off =\off\dff
\begin{pmatrix}
\off 0 &
\bm{\tau}^{\dff *} \off
\vspace{4.5pt} \\
\off \bm{\tau} &
0 \off 
\end{pmatrix}
\off.
\end{equation}

\vspace{-12pt}\vspace{3pt}
In\sss this case we will\sss say\sss that\sss elliptic pairs
$\sigma\fff,\off \tau$\sss
is\pss \emph{graded}\pss and\qss \emph{odd}.\oss
In\sss this case\sss\vspace{3pt}
\[
\quad
\mathcal{L}_{\dff +}\dff(\trf \rho\dff)
\off =\off
\mathcal{L}_{\dff +}\dff(\trf \bm{\tau}\dff)
\dff \oplus\dff
\mathcal{L}_{\dff -}\dff(\trf \bm{\tau}^{\dff *}\dff)
\quad
\mbox{and}\quad
\mathcal{L}_{\dff -}\dff(\trf \rho\dff)
\off =\off
\mathcal{L}_{\dff -}\dff(\trf \bm{\tau}\dff)
\dff \oplus\dff
\mathcal{L}_{\dff +}\dff(\trf \bm{\tau}^{\dff *}\dff)
\pff.
\]

\vspace{-12pt}\vspace{3pt}
If\dss $\bm{\tau}$\sss is\dss a skew-adjoint\sss isometry,\oss
i.e.\qss 
$-\qff \bm{\tau}\off =\off \bm{\tau}^{\dff *}\off =\off \bm{\tau}^{\dff -\dff 1}$\nsp\dnsp,\oss
then 
$\mathcal{L}_{\dff +}\dff(\trf \rho\dff)\fff,\off
\mathcal{L}_{\dff -}\dff(\trf \rho\dff)$\sss
are\sss the eigen\-spaces of\dss $\rho$\sss with\sss 
the eigenvalues\sss $i\fff,\off -\qff i$\sss respectively.\oss
If\dss also\sss
$N\off =\off \mathcal{L}_{\dff +}\dff(\trf \rho\dff)$\nnsp,\oss
then\sss
$\sigma\fff,\off \tau$\nnsp,\qss $\rho$\sss and\sss $N$\sss are\sss in\sss the
normal\dss form.\oss
The corresponding decomposition\dss  
is
$E
\pff =\pff 
\Delta_{\trf F}\qff \oplus\qff \Delta_{\trf F}^{\fff \perp}$\nsp,\oss
where\sss
$\Delta_{\trf F}
\off =\off 
\{\trf (\dff u\fff,\qff u\trf)\qff |\qff u\qff \in\qff F\qff\}$\sss
and\dss
$\Delta_{\trf F}^{\fff \perp}
\off =\off 
\{\trf (\dff u\fff,\qff -\qff u\trf)\qff |\qff u\qff \in\qff F\qff\}$\nnsp.\oss

\myuppar{The positive eigenspaces of\sss 
$\sigma\dff \cos\dff \theta\qff +\qff \tau\dff \sin\dff \theta$\nnsp.}
Let\sss $\theta$\sss run over\sss the interval\sss $[\dff 0\fff,\qff \pi\trf]$\nnsp.\oss
Then\sss the points\sss $(\dff \cos\dff \theta\fff,\off \sin\dff \theta \trf)$\sss
run over a standard\sss half-circle.\oss
If\dss $\sigma\fff,\qff \tau$\sss is\dss a self-adjoint\sss elliptic pair,\oss then\sss
$\sigma\dff \cos\dff \theta\qff +\qff \tau\dff \sin\dff \theta$\nnsp,\qss
$\theta\qff \in\qff [\dff 0\fff,\qff \pi\trf]$\sss
is\dss a continuous family of\dss self-adjoint\sss invertible operators.\oss
Let\sss
$E^{\dff +}\fff(\trf \theta\trf)\fff,\off\off
E^{\dff -}\fff(\trf \theta\trf)
\qff \subset\qff E$\sss
be\sss the sum of\dss eigenspaces of\dss
$\sigma\dff \cos\dff \theta\qff +\qff \tau\dff \sin\dff \theta$\sss
corresponding\sss to\sss the positive and\sss negative eigenvalues respectively.\oss
Clearly,\pss 
$E^{\dff +}\fff(\dff 0\dff)
\off =\off
E^{\dff +}$\dnsp,\qss
$E^{\dff +}\fff(\dff \pi\dff)
\off =\off
E^{\dff -}$\nsp,\oss
and\dss $E^{\dff +}\fff(\trf \theta\trf)$\nnsp,\qss
$\theta\qff \in\qff [\dff 0\fff,\qff \pi\trf]$\sss
is\dss a path connecting\sss $E^{\dff +}$\sss with\sss $E^{\dff -}$\nsp.\oss
The subspaces\sss $E^{\dff +}$\sss and\sss $E^{\dff -}$\sss
do not\sss change during\sss our deformations,\oss
and\sss $E^{\dff +}\fff(\trf \theta\trf)$\nnsp,\qss
$\theta\qff \in\qff [\dff 0\fff,\qff \pi\trf]$\sss remains\sss to be a path
connecting\sss $E^{\dff +}$\sss with\sss $E^{\dff -}$\nsp.\oss
Similarly,\pss
$E^{\dff -}\fff(\trf \theta\trf)$\nnsp,\qss
$\theta\qff \in\qff [\dff 0\fff,\qff \pi\trf]$\sss
is\dss a continuous path connecting\sss $E^{\dff -}$\sss with\sss $E^{\dff +}$\dnsp.\oss
These paths play a crucial\dss role in\sss the definition of\dss the\sss topological\dss index.\oss

It\dss is\dss instructive\sss to compute\sss the path\sss
$E^{\dff +}\fff(\trf \theta\trf)\fff,\off
\theta\qff \in\qff [\dff 0\fff,\qff \pi\trf]$\sss
for\sss $\sigma\fff,\off \rho\fff,\off \tau$\sss 
as\sss in\qss (\ref{sigma-rho-tau-phi}).\oss
After\sss passing\sss to\sss the decomposition\sss
$E
\off =\off 
\Delta\dff(\trf \varphi\trf)\dff \oplus\dff \Delta^{\fff \perp}\dff(\trf \varphi\trf)$\sss
and\sss the diagonal\sss standard\sss form\qss (\ref{sigma-rho-tau-alt})\qss
we\sss get\vspace{3pt}
\[
\quad
\sigma\dff \cos\dff \theta\qff +\qff \tau\dff \sin\dff \theta
\off =\off\dff
\begin{pmatrix}
\off 0 &
\pff \cos\dff \theta\qff -\qff i\dff \sin\dff \theta \off
\vspace{6pt} \\
\off \cos\dff \theta\qff +\qff i\dff \sin\dff \theta &
0 \off 
\end{pmatrix}
\off\dff =\off\dff
\begin{pmatrix}
\off 0 &
\off \overline{z} \off
\vspace{6pt} \\
\off z &
0 \off 
\end{pmatrix}
\off,
\]

\vspace{-12pt}\vspace{3pt}
where\sss $z\off =\off \cos\dff \theta\qff +\qff i\dff \sin\dff \theta$\nnsp.\oss
It\sss follows\sss that\sss\vspace{3pt}
\[
\quad
E^{\dff +}\fff(\trf \theta\trf)
\off =\off
\bigl\{\qff 
(\qff a\fff,\qff z\dff a \qff)
\qff\fff \bigl|\qff\fff
a\qff \in\qff \Delta\dff(\trf \varphi\trf)
\qff\bigr\}
\off \subset\off\dff
\Delta\dff(\trf \varphi\trf)\dff \oplus\dff \Delta^{\fff \perp}\dff(\trf \varphi\trf)
\qff.
\]

\vspace{-12pt}\vspace{3pt}
By\sss passing\sss to\sss the standard\sss form\qss (\ref{sigma-rho-tau-phi})\qss
and\sss taking\sss into account\qss (\ref{diagonals}),\oss
we see\sss that\vspace{4.5pt}
\[
\quad
E^{\dff +}\fff(\trf \theta\trf)
\off =\off
\bigl\{\pff 
\bigl(\trf u\fff,\qff \varphi\dff(\dff u\trf)\qff\bigr)
\qff +\qff
z\qff \bigl(\trf u\fff,\qff -\qff \varphi\dff(\dff u\trf)\qff\bigr)
\pff\fff \bigl|\pff\fff
u\qff \in\qff E^{\dff +}
\pff\bigr\}
\]

\vspace{-31.5pt}
\[
\quad
\phantom{E^{\dff +}\fff(\trf \theta\trf)
\off }
=\off
\bigl\{\pff 
\bigl(\qff (\trf 1\qff +\qff z\trf)\dff u\fff,\qff 
(\trf 1\qff -\qff z\trf)\qff \varphi\dff(\dff u\trf)\qff\bigr)
\pff\fff \bigl|\pff\fff
u\qff \in\qff E^{\dff +}
\pff\bigr\}
\qff
\]

\vspace{-12pt}\vspace{3pt}
and\dss hence\vspace{3.0pt}
\begin{equation}
\label{plus-path}
\quad
E^{\dff +}\fff(\trf \theta\trf)
\off =\off
\left\{\pff \left.
\left(\qff u\fff,\off 
\frac{1\qff -\qff z}{1\qff +\qff z}\qff \varphi\trf(\trf u\trf)\qff\right)
\off\fff \right|\off\fff
u\qff \in\qff E^{\dff +}
\pff\right\}
\qff,
\quad
\end{equation}

\vspace{-12pt}\vspace{3.0pt}
where\sss for\sss $z\off =\off -\qff 1$\sss the\sss last\sss expression\dss
is\dss interpreted as\sss $E^{\dff -}$\dnsp.\oss

\newpage
\mysection{Self-adjoint\qss symbols\qss and\qss boundary\qss conditions}{symbols-conditions}

\myuppar{Self-adjoint\sss bundle automorphisms.}
Let\sss $F$\sss be a\dss Hermitian\dss vector bundle over a\sss topological\sss space\sss $X$\nnsp,\oss
and\sss $a\dff \colon\dff F\qff \ttoo\qff F$\sss be an endomorphisms of\dss $F$
covering\sss $\id_{\trf X}$\nsp.\oss
For\dss $Y\qff \subset\qff X$\sss we denote by\sss $F\trf |\qff Y$\sss
the restriction of\dss $F$\sss to\sss $Y$ and\dss by\sss $a\trf |\trf Y$\sss
the induced endomorphism\sss $F\trf |\qff Y\qff \ttoo\qff F\trf |\qff Y$\dnsp.\oss
If\dss $a$\sss is\dss a self-adjoint\sss automorphism,\oss then $a$ defines a decomposition\sss 
$F
\off =\off 
\mathbb{P}^{\dff +}\fff (\trf a\trf)
\trf \oplus\trf 
\mathbb{P}^{\dff -}\fff (\trf a\trf)$\sss
into\sss the sums 
$\mathbb{P}^{\dff +}\fff (\trf a\trf)$\sss and\sss $\mathbb{P}^{\dff -}\fff (\trf a\trf)$\sss
of\dss  eigen\-spaces of\sss $a$\sss corresponding\sss to\sss the positive and\sss
negative eigenvalues of\sss $a$\sss respectively.\oss
Up\sss to homotopy\sss $a$\sss
is\dss determined\dss by\sss  
$\mathbb{P}^{\dff +}\fff (\trf a\trf)\qff \subset\pff F$\dnsp.\oss

\myuppar{Self-adjoint\sss symbols.}
Let\sss $X$\sss be a\sss compact\sss manifold.\oss
For convenience,\oss we will\sss assume\sss that\sss $X$\sss
is\dss equipped\sss with a riemannian\sss metric.\oss 
Let\sss $T\dff X$\sss be\sss the\sss tangent\sss bundle of\dss $X$\sss 
and\sss let\sss $B\dff X$\sss and\sss $S\dff X$\sss
be\sss the bundles of\dss unit\sss balls and\sss 
unit\sss spheres in\sss $T\dff X$\sss respectively.\oss
Let\sss $\pi\dff \colon\dff T\dff X\qff \ttoo\qff X$\sss
be\sss the bundle projection.\oss
When\sss there\dss is\dss no danger of\dss confusion,\oss we will\sss denote\sss by\sss $\pi$\sss
also various\sss related\sss projections,\oss such as\sss
$S\dff X\qff \ttoo\qff X$\sss
and\sss
$B\dff X\qff \ttoo\qff X$\nnsp.\oss
We will\sss denote\sss the fiber of\dss $T\dff X$\sss at\sss $x\qff \in\qff X$\sss
by\sss $T_x\qff X$\sss and use similar notations for other bundles.\oss

Let\sss $E$\sss be a complex vector bundle on\sss $X$\sss
equipped\sss with a\sss Hermitian\sss metric.\oss 
In contrast\sss with\sss the riemannian\sss metric on $X$\sss
the\dss Hermitian\dss metric on $E$ will\sss play a key role.\oss
An\qss \emph{(elliptic self-adjoint)\qss symbol}\oss 
is\dss a self-adjoint\sss automorphism\sss $\sigma$\sss
of\dss the restriction\sss $\pi^{\fff *}\dff E\trf |\trf S\dff X$\nnsp.\oss
For such a symbol\sss $\sigma$\sss we will\sss use\sss the notations\sss
$E^{\dff +}\fff (\trf \sigma\trf)\off =\off \mathbb{P}^{\dff +}\fff (\trf \sigma\trf)$\sss
and\sss
$E^{\dff -}\fff (\trf \sigma\trf)\off =\off \mathbb{P}^{\dff -}\fff (\trf \sigma\trf)$\nnsp.\oss

\myuppar{Self-adjoint\sss symbols of\trs order \nsp$1$\dnsp.}
Suppose\sss that\sss $X$\sss has non-empty\sss boundary\sss
$Y\off =\off \partial\trf X$\nnsp.\oss
Let\sss $S\fff Y$\sss be\sss the bundle of\dss
unit\sss spheres in\sss the\sss tangent\dss bundle $T\dff Y$\dnsp,\oss
and\dss let\sss $B\dff X_{\trf Y}$\sss and\sss $S\dff X_{\trf Y}$ 
be\sss the restrictions\sss to $Y$ of\dss the bundles 
$B\dff X$\sss and\sss $S\dff X$\sss respectively.\oss
Let\dss $D\dff X\off =\off S\fff X\qff \cup\qff B\dff X_{\trf Y}$\nsp.\oss

For $y\qff \in\qff Y$\dss let\sss $\nu_y$\sss
be\sss unit\sss normal\sss to\sss the\sss tangent\sss space $T_y\dff Y$\sss 
in\sss the\sss tangent\sss space $T_y\trf X$\sss
pointing\sss into $X$\nnsp.\oss
We will\sss say\sss that\sss a symbol\sss $\sigma$\sss
is\dss a\qss \emph{self-adjoint\sss symbol\sss of\trs order\dss $1$}\qss
if\dss for every\sss 
$y\qff \in\qff Y$\dnsp,\qss $u\qff \in\qff S_y\dff Y$\sss and\sss
$\theta\qff \in\qff [\trf 0\fff,\qff \pi\trf]$\sss we have\vspace{3pt} 
\begin{equation}
\label{sin-cos}
\quad
\sigma\dff(\trf
\nu_y\dff \cos\dff \theta
\pff +\qff
u\dff \sin\dff \theta
\trf)
\off =\off
\sigma_y\trf \cos\dff \theta
\pff +\qff
\tau_{u}\trf \sin\dff \theta
\end{equation}

\vspace{-12pt}\vspace{3pt}
for some self-adjoint\sss endomorphisms\sss $\sigma_y$\sss and\sss $\tau_u$\sss
of\dss the fiber of\dss $E_{\dff y}$\sss of\dss $E$ over $y$\nnsp.\oss
Such a symbol $\sigma$ may\dss be\sss thought\sss
as\sss the symbol\sss of\dss a pseudo-differential\sss operator\sss
of\dss order\sss $1$\nnsp.\oss
Clearly,\pss $\sigma_y\off =\off \sigma\dff(\trf \nu_y\trf)$\nnsp.\oss
Since\sss $\sigma\dff(\trf u\trf)$\sss is\dss invertible for every\sss
$u\qff \in\qff S_y\qff X_{\trf Y}$\nsp,\oss
the pair\sss $\sigma_y\dff,\qff \tau_u$\dss is\dss an elliptic\sss pair
in\sss the sense of\trs Section\qss \ref{boundary-algebra}.\oss
For\sss $u\qff \in\qff S_y\qff X_{\trf Y}$\sss let\vspace{3pt}
\[
\quad
\rho_{\fff u}
\off =\off 
\sigma_y^{\dff -\dff 1}\dff \tau_u
\qff.
\]

\vspace{-12pt}\vspace{3pt}
Clearly,\oss the point\sss $y\off =\off \pi\trf(\dff u\trf)$\sss
is\dss determined\sss by\sss the\sss tangent\sss vector\sss $u$\nnsp.\oss
We will\sss always assume\sss that\sss $y\fff,\qff u$\sss are related\sss
in\sss this way.\oss
For $y\qff \in\qff Y$\dss let\sss
$E_{\dff y}^{\dff +}$\sss and\sss $E_{\dff y}^{\dff -}$\sss
be\sss the fibers at $\nu_y\qff \in\qff S\fff X_{\trf Y}\qff \subset\qff S\dff X$\sss of\dss the bundles\sss
$E^{\dff +}\fff (\trf \sigma\trf)$\sss and\sss $E^{\dff -}\fff (\trf \sigma\trf)$\sss
respectively.\oss
Then\sss
$E_{\dff y}
\off =\off
E_{\dff y}^{\dff +}\qff \oplus\qff E_{\dff y}^{\dff -}$\nsp.\oss

\myuppar{Bundle-like symbols and extensions of\dss the bundles\sss
$E^{\dff +}\fff (\trf \sigma\trf)$\sss and\sss
$E^{\dff -}\fff (\trf \sigma\trf)$\nnsp.}
We will\sss say\sss that\sss a self-adjoint\sss symbol\sss $\sigma$\sss of\dss order\sss $1$\sss
is\pss \emph{bundle-like}\oss if\sss 
$\tau_{\fff u}$\sss depends only\sss on\sss $y\off =\off \pi\trf(\dff u\trf)$\nnsp.\oss
This property depends only on\sss the restriction of\sss $\sigma$\sss to\sss $S\dff X_{\trf Y}$.\oss
If\dss $\sigma$\sss is\dss bundle-like,\oss then\sss\vspace{2.5pt}
\[
\quad
\sigma\trf(\trf
\nu_y\dff \cos\dff \theta
\pff +\qff
u\dff \sin\dff \theta
\trf)
\off =\off
\sigma_y\trf \cos\dff \theta
\pff +\qff
\tau_{u}\trf \sin\dff \theta
\]

\vspace{-12pt}\vspace{2.5pt}
depends only on $y$ and\sss $\theta$\nnsp,\oss
and\dss hence\sss the rule\vspace{2.5pt}
\[
\quad
\bm{\sigma}\trf(\trf
\nu_y\dff \cos\dff \theta
\pff +\qff
r\halfff u\dff \sin\dff \theta
\trf)
\off =\off
\sigma\trf(\trf
\nu_y\dff \cos\dff \theta
\pff +\qff
u\dff \sin\dff \theta
\trf)
\qff,
\]

\vspace{-12pt}\vspace{2.5pt}
where\sss $r\qff \in\qff [\trf 0\fff,\qff 1\trf]$\nnsp,\oss
defines a\qss \emph{canonical\sss extension}\qss
of\dss $\sigma$\sss to a self-adjoint\sss automorphism $\bm{\sigma}$ 
of\dss the restriction of\dss the bundle\sss $\pi^{\fff *}\dff E$\sss to\sss
$D\dff X\off =\off S\fff X\qff \cup\qff B\dff X_{\trf Y}$\nsp.\oss
Of\dss course,\oss the bundle-like property of\dss $\sigma$\sss 
matters only\sss for\sss $r\off =\off 0$\nnsp.\oss
If\dss $\sigma$\sss is\dss bundle-like,\oss the bundles\vspace{2.5pt}
\[
\quad
\mathbb{E}^{\dff +}\fff (\trf \sigma\trf)
\off =\off
\mathbb{P}^{\dff +}\fff (\trf \bm{\sigma}\trf)
\quad
\mbox{and}\quad
\mathbb{E}^{\dff -}\fff (\trf \sigma\trf)
\off =\off
\mathbb{P}^{\dff -}\fff (\trf \bm{\sigma}\trf)
\]

\vspace{-12pt}\vspace{2.5pt}
extend\sss from\sss $S\fff X$\sss to\sss
$D\dff X\off =\off S\fff X\qff \cup\qff B\dff X_{\trf Y}$\sss the bundles\sss
$E^{\dff +}\fff (\trf \sigma\trf)$\sss and\sss $E^{\dff -}\fff (\trf \sigma\trf)$\sss
respectively.\oss
The extended\dss bundles\sss
$\mathbb{E}^{\dff +}\fff (\trf \sigma\trf)\fff,\off \mathbb{E}^{\dff -}\fff (\trf \sigma\trf)$\sss 
lead\sss to\sss their $K$\dnsp-theory\sss classes\vspace{2.5pt}
\[
\quad
e^{\dff +}\fff (\trf \sigma\trf)\dff,\off\off
e^{\dff -}\fff (\trf \sigma\trf)
\pff \in\pff
K^{\dff 0}\dff (\trf D\dff X\trf)
\off =\off
K^{\dff 0}\dff (\trf S\dff X\dff \cup\dff B\dff X_{\trf Y}\trf)
\qff.
\]

\vspace{-12pt}\vspace{2.5pt}
Let\sss
$\varepsilon^{\dff +}\fff (\trf \sigma\trf)\dff,\off\off
\varepsilon^{\dff -}\fff (\trf \sigma\trf)$\sss
be\sss the\sss images of\dss these classes under\sss the coboundary\sss map\vspace{2.5pt}
\[
\quad
K^{\dff 0}\dff (\trf D\dff X\trf)
\qff \ttoo\qff
K^{\dff 1}\dff (\trf B\dff X\fff,\qff  D\dff X\trf)
\qff.
\]

\vspace{-12pt}\vspace{2.5pt}
The kernel\sss of\dss the coboundary\sss map\dss is\dss equal\sss to\sss the image
of\dss the natural\sss map\vspace{2.5pt}
\[
\quad
K^{\dff 0}\dff (\trf B\dff X\trf)
\qff \ttoo\qff
K^{\dff 0}\dff (\trf D\dff X\trf)
\qff.
\]

\vspace{-12pt}\vspace{2.5pt}
It\sss follows\sss that\sss $\varepsilon^{\dff +}\fff (\trf \sigma\trf)\off =\off 0$\sss
if\trs $\mathbb{E}^{\dff +}\fff (\trf \sigma\trf)$\sss
can\sss be extended\sss to\sss $B\dff X$\nnsp.\oss
Clearly,\pss  
$\mathbb{E}^{\dff +}\fff (\trf \sigma\trf)\dff \oplus\trf \mathbb{E}^{\dff -}\fff (\trf \sigma\trf)$\sss
is\dss equal\dss to\sss the restriction\sss
$\pi^{\fff *}\dff E\trf |\trf D\dff X$\nnsp,\oss
and\sss hence
$\varepsilon^{\dff +}\fff (\trf \sigma\trf)
\qff +\qff
\varepsilon^{\dff -}\fff (\trf \sigma\trf)
\off =\off
0$\nnsp.\oss

\myuppar{Boundary conditions.}
A\qss \emph{boundary\sss condition}\pss for a self-adjoint\sss symbol\sss $\sigma$\sss
of\dss degree $1$\sss is\dss a subbundle\sss $N$\sss 
of\dss the restriction of\dss the bundle\sss 
$\pi^{\fff *}\dff E$\sss to\sss $S\dff Y$\dnsp.\oss
Let\sss $\sigma$\sss be a self-adjoint\sss symbol\sss of\dss order $1$ and\sss $N$\sss 
be a boundary condition for $\sigma$\nnsp.\oss 

\dnsp$N$\sss
is\dss said\sss to be\qss \emph{self-adjoint}\qss if\trs for every
$u\qff \in\qff S\fff Y$\sss the fiber $N_{\fff u}$ 
is\dss a\sss lagrangian\sss subspace
of\dss $E_{\dff y}$\sss with respect\sss to\sss the form\sss
$[\trf a\fff,\qff b\trf]_{\dff y}
\off =\off
\sco{\dff \sigma_y\qff a\fff,\qff b\dff}$\nnsp,\oss
where $\sco{\dff \bullet\fff,\qff \bullet\dff}$ is\dss
the\sss Hermitian\sss metric\sss in\sss $E$\nnsp.\oss

\dnsp$N$\sss
is\dss said\sss to be an\qss \emph{elliptic}\qss or\sss satisfying\sss the\qss 
\emph{Shapiro--Lopatinskii\dss condition}\pss
if\trs for every
$u\qff \in\qff S\fff Y$\sss the fiber\sss $N_{\fff u}$\sss 
is\dss complementary\sss to\sss
$\mathcal{L}_{\dff -}\dff(\trf \rho_{\fff u}\dff)$\nnsp,\oss
i.e.\qss\vspace{3pt}
\[
\quad
N_{\fff u}\dff \cap\dff \mathcal{L}_{\dff -}\dff(\trf \rho_{\fff u}\dff)
\off =\off
0\quad
\mbox{and}\quad
N_{\fff u}\qff +\qff \mathcal{L}_{\dff -}\dff(\trf \rho_{\fff u}\dff)
\off =\off
E_{\dff y}
\qff,
\]

\vspace{-12pt}\vspace{3pt}
where\sss $y\off =\off \pi\trf(\dff u\trf)$\nnsp.\oss 
If\dss $N$\sss is\dss elliptic,\oss then\sss $N$\sss is\dss self-adjoint\dss
if\dss and\sss only\trs if\trs for every\sss
$u\qff \in\qff S\fff Y$\sss the fiber\sss $N_{\fff u}$\sss
is\dss a self-adjoint\sss boundary condition for\sss 
$\sigma_y\dff,\qff \tau_u$\dss
in\dss the sense of\trs Section\qss \ref{boundary-algebra}.\oss

\dnsp$N$\dss 
is\dss said\sss to be\qss \emph{bundle-like}\pss
if\dss the fibers $N_{\fff u}$ depend only on 
$y\off =\off \pi\trf(\dff u\trf)$\nnsp,\oss
i.e.\qss if\dss $N$\sss is\dss
equal\dss to\sss the\sss lift\sss to\sss $S\dff Y$\sss of\dss a subbundle
of\dss the bundle\sss 
$E_{\trf Y}$\nsp,\oss 
the restriction of\dss the bundle\sss $E$\sss to\sss $Y$\dnsp.\oss

\myuppar{Self-adjoint\dss boundary conditions.}
The following\sss facts immediately\sss follow\sss from\trs Section\qss \ref{boundary-algebra}.\oss
If\dss $N$\sss is\dss a self-adjoint\sss boundary condition,\oss
then\sss the dimension of\dss $E$\sss is\dss even,\oss
the dimension of\dss $N$\sss is\dss equal\dss to\sss the half\dss
of\dss the dimension of\dss $E$\nnsp,\oss
and\vspace{3pt}
\[
\quad 
N_{\fff u}\trf \cap\qff E_{\dff y}^{\dff +}
\off =\off
N_{\fff u}\trf \cap\qff E_{\dff y}^{\dff -}
\off =\off
0
\]

\vspace{-12pt}\vspace{3pt}
for every\sss $u\qff \in\qff S\dff Y$\dnsp.\oss
It\sss follows\sss that\sss the dimensions of\dss
$E_{\dff y}^{\dff +}$ and\sss $E_{\dff y}^{\dff -}$\sss
are also equal\dss to\sss the half\dss
of\dss $\dim\dff E$\nnsp,\oss
and\sss $N_{\fff u}$\sss is\dss equal\sss to\sss the graph\sss
of\dss an\sss isomorphism\sss
$\varphi_{\fff u}\dff \colon\dff
E_{\dff y}^{\dff +}
\qff \ttoo\qff
E_{\dff y}^{\dff -}$\nsp,\oss
i.e.\vspace{1.5pt}
\[
\quad
N_{\fff u}
\off =\off
\left\{\pff
a\qff +\qff \varphi_{\fff u}\dff(\trf a\trf)\qff \left|\pff
a\qff \in\qff E_{\dff y}^{\dff +}
\pff \right.\right\}
\qff.
\]

\vspace{-12pt}\vspace{1.5pt}
If\dss $\sigma$\sss is\dss self-adjoint\sss and\sss unitary,\oss 
i.e.\sss $\sigma\off =\off \sigma^{\dff *}\off =\off \sigma^{\dff -\dff 1}$\dnsp,\oss
then\sss the isomorphisms\sss $\varphi_{\fff u}$\sss are isometries,\oss
and every\sss family of\dss such\sss isometries defines a self-adjoint\sss
boundary condition.\oss

\myuppar{Normalized\sss symbols and\dss boundary conditions.}
Suppose\sss that\sss $\sigma$\sss
is\dss a self-adjoint\sss symbol\sss of\dss order $1$
and\dss $N$\sss is\dss an elliptic self-adjoint\dss boundary condition\sss for $\sigma$\nnsp.\oss
We will\sss say\sss that\sss the pair\sss $\sigma\dff,\off N$\sss
is\qss \emph{normalized}\oss if\dss  
$\sigma_y\dff,\off \rho_{\fff u}\dff,\off \tau_u$\sss and\sss
$N_{\fff u}$\sss are normalized\sss in\sss 
the sense of\trs Section\qss \ref{boundary-algebra}\qss
for every\sss $u\qff \in\qff S\dff Y$\dnsp.\oss
If\dss $\sigma\dff,\off N$\sss is\dss normalized,\oss
then\sss the operators\sss $\rho_{\fff u}$\sss and\sss $\tau_{\fff u}$\sss are determined\dss
by\sss $\sigma_y$\sss and\sss the boundary condition\sss $N_{\fff u}$\nsp.\oss
In\sss particular,\oss if\dss $\sigma\dff,\off N$\sss is\dss normalized
and\sss $N$\sss is\dss bundle-like,\oss then\sss $\sigma$\sss is\dss
also bundle-like.\oss

\myuppar{Deformations of\dss symbols and\dss boundary conditions.}
In\trs Section\qss \ref{boundary-algebra}\qss
we con\-struct\-ed a deformation of\dss an arbitrary elliptic pair\sss together\sss
with a self-adjoint\sss boundary condition\sss to 
a normalized\sss elliptic pair and\sss boundary condition.\oss 
Moreover,\oss we saw\sss that\sss this construction\dss is\sss nearly canonical\sss
and can be applied\sss in continuous
families of\dss elliptic pairs and self-adjoint\sss boundary conditions.\oss
The resulting deformations of\dss families are unique up\sss to homotopy.\oss
In\sss particular,\oss such deformation can\sss be applied\sss to\sss the family of\dss
elliptic pairs\sss $\sigma_y\dff,\qff \tau_u$\dss and\sss
boundary conditions\sss $N_{\fff u}$\sss
parameterized\sss by\sss $u\qff \in\qff S\fff Y$\dnsp.\oss
By\sss the usual\dss homotopy extension\sss property\sss 
this deformation can\sss be extended\sss to a deformation of\dss
the symbol\sss $\sigma$\nnsp.\oss
This deformation ends with\sss
new\sss pair $\sigma\fff,\off N$\nnsp,\oss 
and\dss this new\sss pair\sss $\sigma\fff,\off N$\sss is\dss normalized.\oss 

If\dss $N$\sss was a bundle-like boundary condition\sss
before\sss the deformation,\oss
it\sss remains bundle-like during\sss the deformation.\oss
Indeed,\oss the boundary condition\sss may\sss be changed only\sss during\sss the first\sss
stage of\dss the deformation.\oss
This stage\dss is\dss canonical\sss and depends only\sss on\sss
$\sigma_y\dff,\off y\qff \in\qff Y$\dnsp,\oss and\sss hence keeps\sss 
$N$\sss being\sss bundle-like.\oss
Also,\oss if\trs both\sss $\sigma$\sss and\sss $N$\sss
were bundle-like before\sss the deformation,\oss
then,\oss obviously,\oss the deformation can be arranged\dss to\sss
keep\sss this property.\oss

\myuppar{Elementary\sss symbols and\sss boundary conditions.}
Let\sss us construct\sss some standard symbols and\sss boundary conditions
which we will\sss call\qss \emph{elementary}\pss ones.\oss

Suppose\sss that\sss a collar neighborhood of\dss $Y$\sss in\sss $X$\sss is\dss identified\sss
with\sss $Y\dff \times\dff [\trf 0\fff,\qff 1\dff)$\nnsp,\oss
and\sss points\sss in\sss this collar are denoted\sss by\sss $(\trf y\fff,\qff x_{\dff n}\trf)$\nnsp,\oss
where\sss $y\qff \in\qff Y$\sss and\sss $x_{\dff n}\qff \in\qff [\trf 0\fff,\qff 1\dff)$\nnsp.\oss
We will\sss treat\sss $x_{\dff n}$\sss as\sss the function\sss
$(\trf y\fff,\qff x_{\dff n}\trf)\off \longmapsto\off x_{\dff n}$\sss
on\sss the collar.\oss
Let\sss 
$\varphi\dff \colon\dff 
(\dff -\qff 1\fff,\qff 1\dff)
\qff \ttoo\qff 
[\trf 0\fff,\qff 1\trf]$\sss
be a smooth function with compact\sss support\sss
such\sss that\sss $\varphi$\sss is\dss equal\sss to $1$ in a neighborhood of\sss $0$\dss
(in\sss this section we will\sss use $\varphi$ only on $[\trf 0\fff,\qff 1\dff)$\nnsp).\oss
Suppose\sss that\sss $F$\sss is\dss a\dss Hermitian\dss vector\sss bundle on\sss $X$\sss
and\sss that\sss the restriction\sss 
$F\trf |\trf Y\dff \times\dff [\trf 0\fff,\qff 1\dff)$\sss 
is\dss presented as an orthogonal\sss
direct\sss sum\sss $F^{\dff +}\dff \oplus\dff F^{\dff -}$\nsp.\oss
Let\sss $\lambda\fff,\off \lambda^{\fff +}\nsp,\off \lambda^{\fff -}$\sss
be some positive real\sss numbers,\oss
and\sss let\sss us use\sss the same notations for the automorphisms
\vspace{1.5pt}
\[
\quad
\lambda\dff \colon\dff F\qff \ttoo\qff F\dff,\quad
\lambda^{\fff +}\dff \colon\dff F^{\dff +}\qff \ttoo\qff F^{\dff +}\dff,\quad
\lambda^{\fff -}\dff \colon\dff F^{\dff -}\qff \ttoo\qff F^{\dff -}
\qff,\quad
\]

\vspace{-12pt}\vspace{1.5pt}
defined as\sss the multiplications\sss by\sss the numbers\sss
$\lambda\fff,\off \lambda^{\fff +}\nsp,\off \lambda^{\fff -}$\sss
respectively.\oss
Let\sss us\dss identify\sss the\sss tangent\sss bundle over\sss the collar\sss
with\sss the\sss pull-back\sss of\dss its restriction\sss $T\dff X\trf |\trf Y$\sss
to\sss the boundary\sss by\sss the map\sss
$(\trf y\fff,\qff x_{\dff n}\trf)\off \longmapsto\off y$\nnsp.\oss
Then we can define endomorphisms\vspace{1.5pt}\vspace{0.5pt}
\[
\quad
\sigma^{\dff \mathrm{i}}
\off =\off
i\trf (\dff 1\qff -\qff \varphi\dff(\trf x_{\dff n}\trf)\trf)\qff 
\lambda\qff
(\dff 1\qff -\qff \varphi\dff(\trf x_{\dff n}\trf)\trf)
\off =\off
i\trf (\dff 1\qff -\qff \varphi\dff(\trf x_{\dff n}\trf)\trf)^{\dff 2}\qff 
\lambda
\qff,
\]

\vspace{-33pt}
\[
\quad
\sigma^{\dff \mathrm{b}}_{\dff +}\qff
\left(\trf
\nu_y\dff \cos\dff \theta
\pff +\qff
u\dff \sin\dff \theta
\trf\right)
\off =\off
\varphi\dff(\trf x_{\dff n}\trf)\qff
\left(\trf
\cos\dff \theta\qff +\qff i\trf \lambda^{\fff +}\dff \sin\dff \theta
\qff\right)
\qff,
\]

\vspace{-33pt}
\[
\quad
\sigma^{\dff \mathrm{b}}_{\dff -}\qff
\left(\trf
\nu_y\dff \cos\dff \theta
\pff +\qff
u\dff \sin\dff \theta
\trf\right)
\off =\off
\varphi\dff(\trf x_{\dff n}\trf)\qff
\left(\trf
-\qff \cos\dff \theta\qff +\qff i\trf \lambda^{\fff -}\dff \sin\dff \theta \vphantom{\lambda^{\fff +}}
\qff\right)
\qff,
\]

\vspace{-12pt}\vspace{1.5pt}\vspace{1.5pt}
of\dss the bundles\sss
$F\fff,\off \pi^{\fff *}\dff F^{\dff +},\off \pi^{\fff *}\dff F^{\dff -}$\nsp.\oss
Let\vspace{1.5pt}\vspace{1.5pt}
\[
\quad
\sigma^{\dff \mathrm{b}}
\off =\off
\sigma^{\dff \mathrm{b}}_{\dff +}
\qff \oplus\qff
\sigma^{\dff \mathrm{b}}_{\dff -}
\qff \colon\qff
\pi^{\fff *}\dff F^{\dff +}\qff \oplus\qff \pi^{\fff *}\dff F^{\dff -}
\off \ttoo\off
\pi^{\fff *}\dff F^{\dff +}\qff \oplus\qff \pi^{\fff *}\dff F^{\dff -}
\off.
\]

\vspace{-12pt}\vspace{1.5pt}\vspace{1.5pt}
Since\sss $\varphi\dff(\trf x_{\dff n}\trf)\off =\off 0$\sss
for $x_{\dff n}$ close\sss to $1$\nsp,\oss
the endomorphisms\sss $\sigma^{\dff \mathrm{b}}$
canonically extends\sss to an endomorphism\sss
$\sigma^{\dff \mathrm{b}}\dff \colon\dff F\qff \ttoo\qff F$\dnsp.\oss
Let\sss\vspace{1.5pt}
\[
\quad
\sigma
\off =\off
\sigma^{\dff \mathrm{b}}
\qff +\qff
\sigma^{\dff \mathrm{i}}
\dff \colon\dff
\pi^{\fff *}\dff F\qff \ttoo\qff\pi^{\fff *}\dff  F
\qff,
\]

\vspace{-12pt}\vspace{1.5pt}
and\sss let\sss us define an endomorphism\sss
$\sigma^{\dff \sa}$\sss of\dss the bundle\sss 
$\pi^{\fff *}\dff E\off =\off \pi^{\fff *}\dff F\qff \oplus\qff \pi^{\fff *}\dff F$\sss by\sss
the matrix\vspace{1.5pt}
\[
\quad
\sigma^{\dff \sa}
\off =\off\dff
\begin{pmatrix}
\qff 0\qff &
\sigma\off
\vspace{4.5pt} \\
\off \sigma^{\fff *} &
0 \dff\off 
\end{pmatrix}
\qff.
\]

\vspace{-12pt}\vspace{1.5pt}
A\sss trivial\sss verification shows\sss that\sss $\sigma^{\dff \sa}$\sss is\dss
a self-adjoint\sss symbol\sss of\trs order\dss $1$\nnsp.\oss

By\sss the very definition\sss the symbol\sss $\sigma^{\dff \sa}$\sss is\dss
bundle-like.\oss
Also,\oss the operator\sss $\sigma^{\dff \mathrm{i}}$\sss over\sss $u\qff \in\qff S\dff X$\sss 
depends only on\sss $y\off =\off \pi\trf(\dff u\trf)$\sss
over\sss the whole $X$\nnsp.\oss 
It\sss follows\sss that\sss $\sigma^{\dff \sa}$\sss
can\sss be extended\sss to an\sss isomorphism\sss
$\pi^{\fff *}\dff E\qff \ttoo\qff \pi^{\fff *}\dff E$\sss
over\sss $B\dff X$\nnsp.\oss
It\sss follows\sss that\sss the bundles\sss
$\mathbb{E}^{\dff +}\fff (\trf \sigma^{\dff \sa}\trf)$\sss
and\sss
$\mathbb{E}^{\dff -}\fff (\trf \sigma^{\dff \sa}\trf)$\sss
can\sss be extended\sss to\sss $B\dff X$\nnsp.\oss
In\sss turn,\oss this\sss implies\sss that\sss 
$\varepsilon^{\dff +}\fff \left(\trf \sigma^{\dff \sa}\trf\right)
\off =\off
\varepsilon^{\dff -}\fff \left(\trf \sigma^{\dff \sa}\trf\right)
\off =\off
0$\nnsp.\oss

There are natural\dss boundary conditions for\sss $\sigma^{\dff \sa}$\dnsp.\oss
For every\sss $y\qff \in\qff Y$\sss 
the endomorphism\sss $\sigma_y^{\dff \sa}$\sss of\vspace{1.5pt}
\[
\quad
E_{\dff y}
\off =\off
F_{\fff y}\dff \oplus\dff F_{\fff y}
\off =\off
\left(\trf F_{\fff y}^{\dff +}\dff \oplus\dff F_{\fff y}^{\dff -}\trf\right)
\qff \oplus\qff
\left(\trf F_{\fff y}^{\dff +}\dff \oplus\dff F_{\fff y}^{\dff -}\trf\right)
\]

\vspace{-12pt}\vspace{1.5pt}
is\dss given\sss by\sss the formula\sss\vspace{1.5pt}
\[
\quad
 (\trf u^{\dff +}\fff,\qff u^{\dff -}\fff,\qff v^{\dff +}\fff,\qff v^{\dff -}\trf)
\off \longmapsto\off
(\trf v^{\dff +}\fff,\qff -\qff v^{\dff -}\fff,\qff u^{\dff +}\fff,\qff -\qff u^{\dff -}\trf)
\qff,
\]

\vspace{-12pt}\vspace{1.5pt}
and\dss the corresponding\sss indefinite\dss Hermitian\dss  
product\sss $[\trf \bullet\fff,\qff \bullet\trf]_{\dff y}$\dss 
has\sss the form\vspace{1.5pt}
\[
\quad
\bigl[\qff
(\trf u^{\dff +}\fff,\qff u^{\dff -}\fff,\qff v^{\dff +}\fff,\qff v^{\dff -}\trf)\fff,\off\off
(\trf a^{\dff +}\fff,\qff a^{\dff -}\fff,\qff b^{\dff +}\fff,\qff b^{\dff -}\trf)
\qff\bigr]
\]

\vspace{-34.5pt}\vspace{-0.25pt}
\[
\quad
=\off
\sco{\dff u^{\dff +}\dff,\qff b^{\dff +} \dff}
\qff -\qff
\sco{\dff u^{\dff -}\dff,\qff b^{\dff -} \dff}
\qff +\qff
\sco{\dff v^{\dff +}\dff,\qff a^{\dff +} \dff}
\qff -\qff
\sco{\dff v^{\dff -}\dff,\qff a^{\dff -} \dff}
\qff.
\]

\vspace{-15.75pt}
Let\sss \vspace{-6.25pt}
\[
\quad
N_{\fff u}^{\dff \sa}
\off =\off 
\mathcal{L}_{\dff +}\dff(\trf \rho_{\fff u}\trf)
\off =\off
\left(\trf F_{\fff y}^{\dff +}\dff \oplus\dff 0\trf\right)
\qff \oplus\qff
\left(\trf 0\dff \oplus\dff F_{\fff y}^{\dff -}\trf\right)
\qff,
\]

\vspace{-12pt}\vspace{1.5pt}
where\sss $y\off =\off \pi\trf(\dff u\trf)$\nnsp.\oss
Then\sss 
$\mathcal{L}_{\dff -}\dff(\trf \rho_{\fff u}\trf)
\off =\off
(\trf N_{\fff u}^{\dff \sa}\trf)^{\dff \perp}$\trs
and\sss since\vspace{1.5pt}
\[
\quad
\bigl[\qff
(\trf u^{\dff +}\fff,\qff 0\fff,\qff 0\fff,\qff v^{\dff -}\trf)\fff,\off\off
(\trf a^{\dff +}\fff,\qff 0\fff,\qff 0\fff,\qff b^{\dff -}\trf)
\qff\bigr]
\]

\vspace{-36pt}
\[
\quad
=\off
\sco{\dff u^{\dff +}\dff,\qff 0 \dff}
\qff -\qff
\sco{\dff 0\dff,\qff b^{\dff -} \dff}
\qff +\qff
\sco{\dff 0\dff,\qff a^{\dff +} \dff}
\qff -\qff
\sco{\dff v^{\dff -}\dff,\qff 0 \dff}
\off =\off
0
\qff,
\]

\vspace{-12pt}\vspace{1.5pt}
the subspaces\sss 
$N_{\fff u}^{\dff \sa}\pff \subset\pff E_{\dff y}$ 
are\sss lagrangian and\sss define a self-adjoint\dss 
boundary condition $N^{\dff \sa}$ for $\sigma^{\dff \sa}$\nsp\dnsp.\oss 
Clearly,\qss $\sigma^{\dff \sa},\off N^{\dff \sa}$\sss are\dss bundle-like and\dss if\dss
$\lambda^{\fff +}\off =\off \lambda^{\fff -}\off =\off 1$\sss 
are\sss in\sss the diagonal\dss standard\dss form.\oss  
We will\sss call\sss symbols and\sss boundary conditions of\dss the form\sss
$\sigma^{\dff \sa},\off N^{\dff \sa}$\dss
\emph{elementary}.\oss

\myuppar{An obstruction\sss to deforming boundary conditions\sss to bundle-like ones.}
Let\sss $\sigma$\sss be a self-adjoint\sss symbol\sss of\dss the order $1$\nnsp,\oss
and\dss let\dss $E^{\dff +}_{\trf Y}$\sss and\sss $E^{\dff -}_{\trf Y}$\dss 
be\sss the bundles over\sss $Y$\sss having as fibers over\sss $y\qff \in\qff Y$\sss
the spaces $E_{\dff y}^{\dff +}$\sss and\sss $E_{\dff y}^{\dff -}$\sss respectively.\oss

Let\sss $N$\sss be a self-adjoint\sss elliptic boundary condition for\sss $\sigma$\nnsp,\oss
and suppose\sss that\sss $\sigma\fff,\off N$\sss are normalized.\oss
Then\sss $N$\sss defines for every\sss $u\qff \in\qff S\fff Y$\sss
an\sss isometry
$\varphi_{\fff u}\dff \colon\dff
E_{\dff y}^{\dff +}
\qff \ttoo\qff
E_{\dff y}^{\dff -}$,\oss
and\sss these isometries define an\sss isometric\sss isomorphism of\dss bundles\sss\vspace{1.5pt}
\[
\quad
\bm{\varphi}\dff \colon\dff 
\pi^{\fff *}\dff E^{\dff +}_{\trf Y}
\qff \ttoo\qff 
\pi^{\fff *}\dff E^{\dff -}_{\trf Y}
\off.
\]

\vspace{-12pt}\vspace{1.5pt}
These bundles are well\sss defined also over $B\fff Y$\sss
and\dss the isometry\sss $\bm{\varphi}$\sss leads\sss to a class\vspace{1.5pt}
\[
\quad
\mathcal{I}\dff(\trf N\trf)
\off \in\off
K^{\dff 0}\dff (\trf B\fff Y,\qff S\fff Y\trf)
\off =\off
K^{\dff 0}\dff (\trf T\dff Y\trf)
\qff,
\]

\vspace{-12pt}\vspace{1.5pt}
where\sss the\sss last\sss group\dss is\dss the $K$\dnsp-theory with compact\sss supports.\oss
Clearly,\oss if\dss $N$\sss is\dss bundle-like,\oss then $\bm{\varphi}$ extends\sss to $B\fff Y$\sss
and\sss hence\sss $\mathcal{I}\dff(\trf N\trf)\off =\off 0$\nnsp.\oss
In\sss fact,\oss a converse\dss is\dss also\sss true.\oss

\mypar{Proposition.}{obstruction}
\emph{$\mathcal{I}\dff(\trf N\trf)\off =\off 0$\sss
if\dss and\dss only\dss if\trs
there exist\dss elementary\sss
$\sigma^{\dff \sa},\off N^{\dff \sa}$\dss
such\sss that\sss $N\dff \oplus\dff N^{\dff \sa}$\sss
can\sss be deformed\sss to a bundle-like\sss boundary condition\sss for\sss 
$\sigma\dff \oplus\dff \sigma^{\dff \sa}$\dnsp.\oss}

\proof
If\dss $\sigma^{\dff \sa},\off N^{\dff \sa}$\dss
are elementary,\oss then\sss $N^{\dff \sa}$\sss is\dss bundle-like and\sss hence\sss
$\mathcal{I}\dff(\trf N^{\dff \sa}\trf)\off =\off 0$\nnsp.\oss
It\sss follows\sss that\sss
$\mathcal{I}\dff(\trf N\trf)
\off =\off 
\mathcal{I}\dff(\trf N\dff \oplus\dff N^{\dff \sa}\trf)$\nnsp.\oss
This implies\sss the\qss ``if''\qss
part\sss of\dss the proposition.\oss

Suppose now\sss that\dss $\mathcal{I}\dff(\trf N\trf)\off =\off 0$\nnsp.\oss
Then\sss there exists a bundle\sss $F$\sss over $Y$ such\sss that\sss
the direct\sss sum of\sss $\bm{\varphi}$\sss and\sss 
the identity\sss isomorphism\vspace{1.5pt}\vspace{-0.375pt}
\[
\quad
\id_{\qff \pi^{\fff *}\dff F} \dff \colon\dff 
\pi^{\fff *}\dff F
\qff \ttoo\qff 
\pi^{\fff *}\dff F
\off
\]

\vspace{-12pt}\vspace{1.5pt}\vspace{-0.375pt}
can be extended\sss to $B\fff Y$\dnsp.\oss
Moreover,\oss by adding\sss to\sss $F$\sss another bundle $F\fff'$ over $Y$\dnsp,\oss
if\dss necessary,\oss we can ensure\sss that\sss $F$\sss is\dss a\sss trivial\sss bundle,\oss
and,\oss in\sss particular,\oss can\sss be extended\sss to $X$\nnsp.\oss
Let\sss us apply\sss the construction of\dss elementary symbols and\sss
boundary conditions\sss to\sss the bundles\sss $F$\dnsp,\qss
$F^{\dff +}\off =\off F$\dnsp,\qss
$F^{\dff -}\off =\off 0$\nnsp,\oss
and\sss
$\lambda\off =\off \lambda^{\fff +}\off =\off \lambda^{\fff -}\off =\off 1$\dss
({\fff}the value of\dss $\lambda^{\fff -}$\sss is\dss actually\sss irrelevant\sss because\sss
the corresponding\sss bundle\dss is\dss zero).\oss
Let\sss $E\fff'\off =\off F\dff \oplus\dff F$\dnsp,\oss
and\dss let\sss $\sigma\fff'\off =\off \sigma^{\dff \sa}\fff,\off 
N\fff'\off =\off N^{\dff \sa}$\sss be\sss
the resulting\sss elementary symbol\sss and\sss the boundary condition\sss $E\fff'$\dnsp.\oss
The above construction\sss provides $\sigma\fff'\dff,\pff N\fff'$\sss in\sss
the diagonal\dss normal\dss form.\oss
Passing\sss to\sss the\qss ``non-diagonal''\qss normal\dss form shows\sss that\sss
taking\sss the direct\sss sum of\dss $\sigma\dff,\pff N$\sss with\sss
$\sigma\fff'\dff,\pff N\fff'$\sss
results in\sss adding\sss $F$\sss to both\sss
$E_{\dff y}^{\dff +}$\sss and\sss $E_{\dff y}^{\dff -}$\sss
and\sss adding\sss the identity\sss morphism\sss
$\id_{\qff \pi^{\fff *}\dff F}$\sss of\dss 
$\pi^{\fff *}\dff F$\sss to\sss $\bm{\varphi}$\nnsp.\oss
Therefore\dss it\dss is\dss sufficient\sss to show\sss that\sss $N$\sss
can\sss be deformed\sss to a bundle-like boundary condition\dss if\dss
$\bm{\varphi}$\sss can\sss be extended\sss to\sss an\sss isomorphism\sss 
$\bm{\widetilde{\varphi}}$\sss over\sss $B\fff Y$\dnsp.\oss
But,\oss clearly,\sss the restriction of\dss $\bm{\widetilde{\varphi}}$\sss
to\sss the zero section of\dss $B\fff Y$\sss
defines a bundle-like boundary condition\sss 
homotopic\sss to\sss $N$\nnsp.\oss  \eproof

\myuppar{Boundary conditions in\sss the classical\dss form.}
Classically,\oss the boundary conditions for\sss $\sigma$\sss
are described\sss by\sss a vector bundle $G$\sss over\sss $Y$\sss
and\sss a surjective bundle map\sss
$b\dff \colon\dff
\pi^{\fff *}\dff E\trf |\trf S\fff Y
\qff \ttoo\qff
G$\sss
over\sss the projection\sss $S\fff Y\qff \ttoo\qff Y$\dnsp.\oss
The corresponding\sss subbundle $N$\sss is\dss 
simply\sss the kernel\sss $\kernel\dff b$\nnsp.\oss
By\sss the definition,\pss
$b$\sss satisfies\sss the\dss Shapiro--Lopatinskii\dss condition\dss
if\dss $N$\sss has\sss this\sss property,\oss
and\dss is\dss self-adjoint\sss if\dss $N$\sss is.\oss
In\sss the non-self-adjoint\sss case\sss the existence of\dss
such $G\fff,\pff b$\sss ensures\sss that\sss 
after adding a symbol\sss of\dss the\sss form\sss
$\sigma^{\dff \mathrm{b}}
\qff +\qff
\sigma^{\dff \mathrm{i}}$\sss
one can deform\sss the symbol\sss to a bundle-like one.\oss
See\sss the\qss \emph{Fifth\dss homotopy}\qss of\trs
Atiyah\dss and\dss Bott\qss \cite{ab1},\oss
or\sss the\trs \emph{Step\dss II}\pss in\dss H\"{o}rmander's\trs proof\dss
of\trs Proposition\qss 20.3.3\qss in\qss \cite{h}.\oss
In\sss its\sss turn,\oss deforming\sss the symbol\dss to a bundle-like one
is\dss the key\sss step\sss in
defining\dss the\qss \emph{topological\dss index}.\oss
See\trs Atiyah\qss \cite{a1},\oss \cite{a2}.\oss

In\sss the self-adjoint\sss case\sss the situation\dss is\dss different.\oss
Namely,\oss every self-adjoint\sss boundary condition\sss $N$\sss
can\sss be realized as\sss the kernel\sss of\dss 
a bundle map\sss $b$\sss as above.\oss
Indeed,\oss for every\sss $u\qff \in\qff S\fff Y$\sss
the subspace\sss $N_{\fff u}\qff \subset\qff E_{\dff y}$\sss
is\dss lagrangian and\dss hence\dss transverse\sss to $E^{\dff +}_{\dff y}$.\oss
Therefore\sss there\dss is\dss a unique projection\sss
$b_{\dff u}\dff \colon\dff E_{\dff y}\qff \ttoo\qff E^{\dff +}_{\dff y}$\sss
with\sss the kernel\sss $N_{\fff u}$\nsp.\oss
These projections define a surjective bundle map\sss
$b\dff \colon\dff
\pi^{\fff *}\dff E\trf |\trf S\fff Y
\qff \ttoo\qff
E^{\dff +}$\sss
with\sss the kernel\sss $N$\sss
covering\sss $S\fff Y\qff \ttoo\qff Y$\dnsp.\oss
At\sss the same\sss time every element\sss of\sss
$K^{\dff 0}\dff (\trf B\fff Y,\qff S\fff Y\trf)$
can be realized as $\mathcal{I}\dff(\trf N\trf)$ for some self-adjoint\sss  
$N\dff,\pff \sigma$\sss
satisfying\dss the\dss Shapiro--Lopatin\-skii\dss condition.\oss
Namely,\oss one can define\sss the operators\sss 
$\rho_{\fff u}\dff,\pff \tau_{\dff u}$\sss
by\qss (\ref{sigma-rho-tau-phi}).\oss
By\sss these reasons\sss in\sss the self-adjoint\sss case\sss
the classical\dss form of\dss boundary conditions does not\dss
help\sss to define\sss the\sss topological\dss index.\oss

\myuppar{The class\sss
$\varepsilon^{\dff +}\fff (\trf \sigma\fff,\qff N\qff)$\nnsp.}
We would\dss like\sss to define an analogue of\dss the class\sss
$\varepsilon^{\dff +}\fff (\trf \sigma\trf)$\sss 
when\sss $\sigma$\sss is\dss not\sss necessarily\sss a\sss
bundle-like symbol.\oss
Such an analogue exist\dss if\dss a self-adjoint\sss elliptic\sss boundary condition\sss $N$\sss 
for $\sigma$ is\dss given,\oss
and,\oss moreover,\pss $N$\sss is\dss bundle-like or at\dss least\sss 
$\mathcal{I}\dff(\trf N\trf)\off =\off 0$\nnsp.\oss

Proposition\qss \ref{obstruction}\qss implies\sss that\sss
$N\dff \oplus\dff N^{\dff \sa}$\sss
can\sss be deformed\sss to a bundle-like\sss boundary condition\sss $N\fff'$\sss for\sss 
$\sigma\dff \oplus\dff \sigma^{\dff \sa}$\sss
for some elementary\sss $\sigma^{\dff \sa},\off N^{\dff \sa}$\dnsp.\oss
After\sss this we can use\sss$N\fff'$\sss to deform\sss 
$\sigma\dff \oplus\dff \sigma^{\dff \sa},\off N\fff'$\sss
to bundle-like pair\sss $\sigma\fff'\fff,\off N\fff'$\dss
(while keeping\sss $N\fff'$\sss intact\fff).\oss
As we saw,\pss
$\varepsilon^{\dff +}\fff (\trf \sigma^{\dff \sa}\trf)
\off =\off 
0$\nnsp.\oss
Since\sss the invariant\sss
$\varepsilon^{\dff +}\fff (\trf \bullet\trf)$\sss
is\dss obviously additive with respect\sss to\sss the direct\sss sums,\oss
this suggests\sss to introduce\sss the class\sss 
$\varepsilon^{\dff +}\fff (\trf \sigma\fff,\qff N\qff)
\off =\off
\varepsilon^{\dff +}\fff (\trf \sigma\fff'\qff)$\nnsp.\oss
Let\sss us check\sss that\sss 
$\varepsilon^{\dff +}\fff (\trf \sigma\fff,\qff N\qff)$\sss
does not\sss depend on\sss the choices made.\oss
Since\sss the normalizing\sss deformations can be applied\sss to families,\oss
in\sss particular,\oss to homotopies,\oss
$\varepsilon^{\dff +}\fff (\trf \sigma\fff,\qff N\qff)$\sss
does not\sss depend on\sss the choice of\dss deformation of\dss
$N\dff \oplus\dff N^{\dff \sa}$\sss
to a bundle-like\sss boundary condition\sss for fixed\sss 
$\sigma^{\dff \sa},\off N^{\dff \sa}$\dnsp.\oss
Different\sss $\sigma^{\dff \sa},\off N^{\dff \sa}$ lead\sss to\sss the same class because
one can simultaneously add\sss two of\dss them and\dss then use
either one of\dss them.\oss
If\dss $\sigma\fff,\pff N$\sss are already normalized and\dss bundle-like,\oss
then\sss there\dss is\dss no need\sss to deform\sss
$N\dff \oplus\dff N^{\dff \sa}$\sss
and\sss
$\sigma\dff \oplus\dff \sigma^{\dff \sa}$\dss
(and even\sss to add\sss $\sigma^{\dff \sa},\off N^{\dff \sa}$)\qss
and\dss hence\vspace{1.5pt}
\[
\quad
\varepsilon^{\dff +}\fff (\trf \sigma\fff,\qff N\qff)
\off =\off
\varepsilon^{\dff +}\fff (\trf \sigma\dff \oplus\dff \sigma^{\dff \sa}\trf)
\off =\off
\varepsilon^{\dff +}\fff (\trf \sigma\trf)
\qff +\qff
\varepsilon^{\dff +}\fff (\trf \sigma^{\dff \sa}\trf)
\off =\off
\varepsilon^{\dff +}\fff (\trf \sigma\trf)
\qff.
\]

\vspace{-12pt}\vspace{1.5pt}
If\trs both $\sigma$ and $N$ are bundle-like,\oss
but\sss the pair $\sigma\fff,\off N$\dss is\dss not\dss normalized,\oss
we can\sss deform\sss $\sigma\fff,\off N$\sss to a normalized\sss pair in\sss the class of\dss
bundle-like pairs.\oss
This deformation\sss will\sss induce a deformation of\dss the subbundle\sss 
$\mathbb{E}^{\dff +}\fff (\trf \sigma\trf)$\nnsp.\oss
It\sss follows\sss that\sss in\sss this case\sss
$\varepsilon^{\dff +}\fff (\trf \sigma\fff,\qff N\qff)
\off =\off
\varepsilon^{\dff +}\fff (\trf \sigma\trf)$\nnsp.\oss

\myuppar{The\sss topological\dss index.}
Let\sss us embed\sss $X$\sss into\sss
$\rrr_{\qff \geq\trf 0}^{\dff n}
\off =\off
\rrr^{\dff n\dff -\dff 1}\dff \times\qff \rrr_{\qff \geq\trf 0}$\sss
for some $n$ in such a way\sss that\sss
$Y\off =\off X\qff \cap\qff (\trf \rrr^{\dff n\dff -\dff 1}\dff \times\dff 0\trf)$\sss
and\sss $X$\sss is\dss transverse\sss to\sss
$\rrr^{\dff n\dff -\dff 1}\dff \times\dff 0$\sss
in\sss $\rrr^{\dff n}$\dnsp.\oss
Let\sss $\mathbb{N}$\sss be\sss the normal\dss bundle\sss 
to\sss $X$\sss in\sss $\rrr^{\dff n}$\dnsp.\oss
The normal\dss bundle of\dss $T\dff X$\sss 
in\sss 
$T\dff \rrr_{\qff \geq\trf 0}^{\dff n}
\off =\off
\rrr_{\qff \geq\trf 0}^{\dff n}
\dff \times\dff
\rrr^{\dff n}$\sss
can\sss be identified\sss with\sss the\sss lift\sss of\dss the bundle\sss
$\mathbb{N}\dff \oplus\dff \mathbb{N}$\sss 
to\sss $T\dff X$\nnsp,\oss
where\sss $\mathbb{N}\dff \oplus\trf 0$\sss is\dss simply\sss
the normal\dss bundle\sss $\mathbb{N}$\nnsp,\qss and\sss
$0\trf \oplus\trf \mathbb{N}$\sss consists of\dss
vectors\sss tangent\sss to\sss the fibers of\dss $\mathbb{N}$\nnsp.\oss
The bundle\sss $\mathbb{N}\dff \oplus\dff \mathbb{N}$\sss
has a natural\sss complex structure resulting\sss from\sss identifying\sss it\dss with\sss
$\mathbb{N}\trf \oplus\trf i\trf \mathbb{N}
\off =\off
\mathbb{N}\trf \otimes\trf \ccc$\nnsp.\oss

The previous paragraph\sss implies\sss that\sss a\sss tubular\sss
neighborhood of\dss $T\dff X$\sss 
in\sss $T\dff \rrr_{\qff \geq\trf 0}^{\dff n}$\sss
can\sss be identified\sss with\sss the bundle\sss $\mathbb{U}$\sss of\dss 
unit\sss balls in\sss 
$\mathbb{N}\dff \oplus\dff \mathbb{N}$\nnsp.\oss
Let\sss $\mathbb{S}$\sss be\sss the bundle of\dss unit\sss spheres in\sss
$\mathbb{N}\dff \oplus\dff \mathbb{N}$\nnsp.\oss
Let\sss us consider\sss the restrictions\sss
$\mathbb{U}\dff |\trf B\dff X$\nnsp,\qss
$\mathbb{S}\dff |\trf B\dff X$\nnsp,\qss
$\mathbb{U}\dff |\trf S\dff X$\nnsp,\oss and\sss
$\mathbb{U}\dff |\trf D\dff X$\sss
to\sss $B\dff X$\nnsp,\qss $B\dff X$\nnsp,\qss $S\dff X$\nnsp,\oss and\sss $D\dff X$\sss
respectively.\oss
The complex structure on\sss $\mathbb{N}\dff \oplus\dff \mathbb{N}$\sss
leads\sss to\sss the\trs Thom\dss isomorphism\vspace{3pt}
\[
\quad
\mathrm{Th}\dff \colon\dff
K^{\dff 1}\dff (\trf B\dff X\fff,\off  D\dff X\trf)
\qff \ttoo\qff
K^{\dff 1}\dff \bigl(\qff \mathbb{U}\dff |\trf B\dff X\fff,\off 
(\trf \mathbb{S}\dff |\trf B\dff X\trf)
\qff \cup\qff 
(\trf \mathbb{U}\dff |\trf D\dff X\trf)
\qff\bigr)
\qff.
\]

\vspace{-12pt}\vspace{3pt}
Clearly,\pss
$(\trf \mathbb{S}\dff |\trf B\dff X\trf)
\qff \cup\qff 
(\trf \mathbb{U}\dff |\trf D\dff X\trf)$\sss
is\dss the boundary of\dss  
$\mathbb{U}\dff |\trf B\dff X 
\qff \subset\qff
T\dff \rrr^{\dff n}
\off =\off
\rrr^{\dff n}
\dff \times\dff
\rrr^{\dff n}$\dnsp,\oss
and\sss $\mathbb{U}\dff |\trf B\dff X_{\trf Y}$\sss
is\dss equal\dss to\sss the intersection of\dss this boundary
with\sss $(\trf \rrr^{\dff n\dff -\dff 1}\dff \times\dff 0\trf)\dff \times\dff \rrr^{\dff n}$\dnsp.\oss
Let\sss $S^{\dff 2 n}$\sss be\sss the compactification of\trs
$\rrr^{\dff n}
\fff \times\dff
\rrr^{\dff n}$\sss by a point $p_{\dff 0}$\nsp.\oss
There\dss is\dss a canonical\dss map\vspace{3pt}
\[
\quad
c_{\trf X}\dff \colon\dff
S^{\dff 2 n}
\qff \ttoo\qff
\mathbb{U}\dff |\trf B\dff X\qff\bigl/\dff
\bigl(\qff 
(\trf \mathbb{S}\dff |\trf B\dff X\trf)
\qff \cup\qff 
(\trf \mathbb{U}\dff |\trf D\dff X\trf)
\qff\bigr)
\qff
\]

\vspace{-12pt}\vspace{3pt}
from\sss $S^{\dff 2 n}$\sss to\sss the quotient\sss of\dss the subspace\dss
$\mathbb{U}\dff |\trf B\dff X 
\qff \subset\qff
\rrr^{\dff n}
\dff \times\dff
\rrr^{\dff n}$\sss
by\sss its boundary,\oss
and\dss this map induces\sss 
the pull-back\sss homomorphism\vspace{3pt}\vspace{-0.125pt}
\[
\quad
c_{\trf X}^{\dff *}\qff \colon\dff
K^{\dff 1}\dff \bigl(\qff \mathbb{U}\dff |\trf B\dff X\fff,\off 
(\trf \mathbb{S}\dff |\trf B\dff X\trf)
\qff \cup\qff 
(\trf \mathbb{U}\dff |\trf D\dff X\trf)
\qff\bigr)
\qff \ttoo\qff
K^{\dff 1}\dff (\trf S^{\dff 2 n}\fff,\qff p_{\dff 0}\trf)
\pff.
\]

\vspace{-12pt}\vspace{3pt}
Of\dss course,\pss
$K^{\dff 1}\dff (\trf S^{\dff 2 n}\fff,\qff p_{\dff 0}\trf)\off =\off 0$\nnsp,\oss
but\sss let\sss us pretend\sss that\sss we don't\dss know\sss this
and define\vspace{1.5pt}
\[
\quad
t\qff(\trf \sigma\fff,\qff N\qff)
\off \in\off
K^{\dff 1}\dff (\trf S^{\dff 2 n}\fff,\qff p_{\dff 0}\trf)
\]

\vspace{-12pt}\vspace{1.5pt}
as\sss the image of\dss the class\sss 
$\varepsilon^{\dff +}\fff (\trf \sigma\fff,\qff N\qff)
\qff \in\qff
K^{\dff 1}\dff (\trf B\dff X\fff,\off  D\dff X\trf)$\sss
under\sss the composition\sss $c_{\trf X}^{\dff *}\dff \circ\trf \mathrm{Th}$\nnsp.\oss
This\dss is\qss (almost\fff)\qss the\qss \emph{topological\dss index}\pss of\dss
$\sigma\fff,\off N$\nnsp,\oss but\sss
in order\sss to get\sss something non-zero,\oss
one needs\sss to allow symbols and\sss 
boundary conditions depending on a parameter.\oss

\myuppar{Families of\dss self-adjoint\sss symbols and\dss boundary conditions.}
Suppose\sss that\sss our symbols and\dss boundary conditions
continuously depend on a parameter\sss $z\qff \in\qff Z$\nnsp,\oss
where\sss $Z$\sss is\dss a\sss topological\sss space.\oss
Then everything\sss in\sss this section can\sss be arranged\sss to be
continuously depending on $z$\nnsp.\oss
In\sss fact,\oss our constructions are either canonical\sss and\sss hence
continuously depend on $z$\nnsp,\oss
or are\sss based on\sss the extension of\dss homotopies,\oss
which can be arranged\sss to be continuously depending on $z$\nnsp.\oss
Moreover,\oss one can allow\sss the manifold\sss $X$\sss
to continuously depend\sss on $z$\nnsp,\oss in\sss the sense\sss that\sss
$X\off =\off X\trf(\trf z\trf)$\sss is\dss the fiber over $z$ of\dss a\sss locally\sss trivial\dss bundle\sss
$p\dff \colon\dff W\qff \ttoo\qff Z$\nnsp,\oss
as\sss in\qss \cite{as4}.\oss
Then\sss the boundaries\sss $\partial\dff X\trf(\trf z\trf)$\sss form another\sss bundle\sss
$V\qff \ttoo\qff Z$\nnsp,\oss
and one can define\sss bundles\sss $T\dff W$\nnsp,\qss $S\dff W$\dnsp,\qss
$S\dff W_{\dff V}$\nsp,\qss etc.,\oss involving,\oss of\dss course,\oss
only\sss vectors\sss tangent\sss to\sss the fibers.\oss
Moreover,\oss for a family\sss 
$\sigma\trf(\trf z\trf)\fff,\off N\trf(\trf z\trf)\dff,\off z\qff \in\qff Z$\sss
of\dss self-adjoint\sss symbols and\dss boundary conditions,\oss
one can define\sss the class\sss 
$\mathcal{I}\dff(\trf N\trf)
\off \in\off
K^{\dff 0}\dff (\trf B\fff V,\qff S\fff V\trf)$\sss
and\dss if\dss $\mathcal{I}\dff(\trf N\trf)\off =\off 0$\nnsp,\oss
define\sss the class\vspace{1.5pt}
\[
\quad
\varepsilon^{\dff +}\fff (\trf \sigma\fff,\qff N\qff)
\pff \in\off 
K^{\dff 1}\dff (\trf B\dff W\fff,\qff  S\dff W\dff \cup\dff B\dff W_{\trf V}\trf)
\qff.
\]

\vspace{-12pt}\vspace{1.5pt}
Following\dss Atiyah\dss and\dss Singer\qss \cite{as4},\oss
 we will\sss now construct\sss the forward\sss image
of\dss $\varepsilon^{\dff +}\fff (\trf \sigma\fff,\qff N\qff)$ in\sss 
$K^{\dff 1}\dff (\trf Z\trf)$\nnsp,\oss
assuming\sss that\sss $Z$\sss is\dss compact.\oss
There\dss is\dss a natural\sss number $n$\sss and a continuous map\sss
$f\dff \colon\dff W\qff \ttoo\qff \rrr^{\dff n}$\sss
such\sss that\sss the restriction $f_{\dff z}$ of $f$ to every\sss fiber\sss
$X\trf(\trf z\trf)\off =\off p^{\dff -\dff 1}\dff(\trf z\trf)$\sss
is\dss a smooth embedding.\oss
Moreover,\oss we can assume\sss that\sss the image of $f$\sss is\dss contained\sss in\sss
the positive half-space\sss
$\rrr_{\qff \geq\trf 0}^{\dff n}
\off =\off
\rrr^{\dff n\dff -\dff 1}\dff \times\qff \rrr_{\qff \geq\trf 0}$\nsp,\oss
and\dss that\sss for every\sss $z\qff \in\qff Z$\vspace{1.5pt}
\[
\quad
\partial\trf X\trf(\trf z\trf)
\off =\off
f^{\dff -\dff 1}_{\dff z}\dff\left(\qff \rrr^{\dff n\dff -\dff 1}\dff \times\qff 0 \qff\right)
\qff.
\]

\vspace{-12pt}\vspace{1.5pt}
and\sss $f_{\dff z}\dff(\qff X\trf(\trf z\trf)\trf)$\sss is\dss transverse\sss to\sss
$\rrr^{\dff n\dff -\dff 1}\dff \times\qff 0$\sss in\sss $\rrr^{\dff n}$\dnsp.\oss
Let\sss us identify\sss $W$\sss with\sss its image in\sss
$Z\dff \times\dff \rrr^{\dff n}$\sss under\sss the map\sss
$x
\off \longmapsto\off 
(\trf p\dff(\dff x\trf)\fff,\qff f\dff(\dff x\trf)\trf)$\nnsp.\oss
Now\sss we can carry\sss the constructions of\dss the previous subsection\sss
with\sss the parameter\sss $z\qff \in\qff Z$\sss and\sss get\sss an element\sss
$t\trf(\trf \sigma,\pff N\qff)
\qff \in\qff
K^{\dff 1}\dff (\trf S^{\dff 2 n}\dff \times\dff Z\fff,\qff p_{\dff 0}\dff \times\dff Z\trf)$\nnsp,\oss
the\qss \emph{topological\dss pre-index},\oss
as\sss the image of\dss
$\varepsilon^{\dff +}\fff (\trf \sigma\fff,\qff N\qff)$\sss
under\sss the\trs Thom\dss isomorphism and\dss the parametrized
analogue of\dss the homomorphism\sss $c_{\trf X}^{\dff *}$\nsp.\oss
Finally,\oss we define\sss the\qss
\emph{topological\dss index}\vspace{3pt}
\[
\quad
\ti\qff (\trf \sigma,\pff N\qff)
\off \in\off
K^{\dff 1}\dff (\trf Z\trf)
\]

\vspace{-12pt}\vspace{3pt}
of\dss the family 
$\sigma\trf(\trf z\trf)\fff,\off N\trf(\trf z\trf)\dff,\off z\qff \in\qff Z$\sss
as\sss the image of\dss
$t\qff( \sigma,\pff N\qff)$\sss
under\sss the\dss Bott\trs periodicity\sss map.\oss

\newpage
\mysection{Multiplicative\qss properties\qss of\pss symbols}{mult}

\myuppar{Coboundary\sss maps.}
For a space $X$ with a base point $x_{\dff 0}$\sss
the reduced\sss $K$\dnsp-group $\widetilde{K}^{\trf 0}\fff(\trf X\trf)$\sss is\dss defined
as\sss the kernel\sss of\dss the canonical\sss map\sss
$K^{\dff 0}\fff(\trf X\trf)
\qff \ttoo\qff
K^{\dff 0}\fff(\trf \{\trf x_{\dff 0}\qff\}\trf)$\nnsp.\oss
There\dss is\dss a canonical\dss isomorphism\sss
$K^{\dff 0}\fff(\trf X\trf)
\qff \ttoo\qff
\widetilde{K}^{\trf 0}\fff(\trf X\trf)
\dff \oplus\trf
\zzz$\nnsp.\oss
For a\qss (nice)\qss subspace\sss $Y\qff \subset\qff X$\sss
the relative group\sss
$K^{\dff 0}\fff(\trf X\fff,\qff Y\trf)$\sss is\dss defined as\sss
$\widetilde{K}^{\trf 0}\fff(\trf X/Y\trf)$\nnsp,\oss
where $Y/Y$\sss is\dss taken as\sss the base point\sss of\sss $X/Y$\dnsp.\oss

For a space $X$\sss let\sss $X^{\dff +}$\sss be\sss the disjoint\sss union of\sss $X$\sss and a point,\oss 
taken as\sss the base point\sss of\sss $X^{\dff +}$\dnsp.\oss
We will\sss denote by $\bm{\Sigma}\dff X$\sss
the reduced suspension of\sss $X^{\dff +}$\dnsp,\oss
and\dss by\sss $\bm{\Sigma}^{\fff n}\trf X$\sss 
the \nsp$n$\dnsp-fold\sss reduced suspension of\sss $X^{\dff +}$\dnsp.\oss
Then\sss 
$\bm{\Sigma}^{\fff n}\trf X
\off =\off
X\dff \times\dff S^n\nsp/\fff X\dff \times\dff s_{\dff 0}$\nsp,\oss
where $s_{\dff 0}$\sss is\dss the base point\sss of\sss $S^n$\dnsp,\oss
and\sss $\bm{\Sigma}^{\dff 0}\trf X\off =\off X^{\dff +}$\dnsp.\oss
Recall\sss that\sss the negative $K$\dnsp-groups are defined as\sss
$K^{\dff -\dff n}\fff(\trf X\trf)
\off =\off
\widetilde{K}^{\trf 0}\fff(\trf \bm{\Sigma}^{\fff n}\dff X\trf)$\nnsp,\oss
and\sss that\sss
$K^{\dff n}\fff(\trf X\trf)
\off =\off
\widetilde{K}^{\dff n\dff +\qff 1}\fff(\trf \bm{\Sigma}\dff X\trf)$\sss
for every\sss $n\qff \in\qff \zzz$\nnsp.\oss

Let\sss us denote by\sss $\mathbf{C}\dff X$\sss the cone of\sss $X$\nnsp,\oss
i.e.\dss $X\dff \times\dff [\trf 0\fff,\qff 1\trf]$\sss
with\sss $X\dff \times\dff 0$\sss 
collapsed\sss to a point.\oss
We\sss identify $X$ with 
$X\dff \times\dff 1\qff \subset\pff \mathbf{C}\dff X$ 
and\dss take\sss the\sss latter point\sss as\sss the base point\sss of\sss 
$\mathbf{C}\dff X$\nnsp.\oss
Suppose\sss that\sss $X$\sss is\dss compact\sss and\sss
$Y\qff \subset\qff X$\sss is\dss a closed subset.\oss
The coboundary\sss maps\sss
\[
\partial\dff \colon\dff
K^{\dff n}\fff(\trf Y\trf)
\qff \ttoo\qff
K^{\dff n\dff +\qff 1}\fff(\trf X\fff,\qff Y\trf)
\]

\vspace{-12pt}
can be defined as\sss the maps
\[
\quad
K^{\dff n}\fff(\trf Y\trf)
\off =\off
\widetilde{K}^{\dff n\dff +\qff 1}\fff(\trf \bm{\Sigma}\dff Y\trf)
\qff \ttoo\qff
\widetilde{K}^{\dff n\dff +\qff 1}\fff(\trf X\qff \cup\qff \mathbf{C}\dff Y\trf)
\off =\off
\widetilde{K}^{\dff n\dff +\qff 1}\fff(\trf X/Y\trf)
\off =\off
K^{\dff n\dff +\qff 1}\fff(\trf X\fff,\qff Y\trf)
\qff
\]

\vspace{-12pt}
induced\dss by\sss the obvious 
quotient\sss map\sss
$X\qff \cup\qff \mathbf{C}\dff Y\qff \ttoo\qff \bm{\Sigma}\dff Y$\nnsp.\oss
The coboundary maps commute with\sss the\dss Bott\dss periodicity\sss map,\oss
i.e.\qss for every $n\qff \in\qff \zzz$\sss the diagram\vspace{-3.75pt}
\[
\quad
\begin{tikzcd}[column sep=boom, row sep=boom]
K^{\dff n}\fff(\trf Y\trf)
\arrow[d, "\dis \qff \beta"]
\arrow[r, "\dis \partial"]
&
K^{\dff n\dff +\qff 1}\fff(\trf X\fff,\qff Y\trf)
\arrow[d, "\dis \qff \beta"]
\\
K^{\dff n\dff -\dff 2}\fff(\trf Y\trf)
\arrow[r, "\dis \partial"]
&
K^{\dff n\dff -\dff 1}\fff(\trf X\fff,\qff Y\trf)
\end{tikzcd}
\]

\vspace{-12pt}\vspace{-0.75pt}
is\dss commutative,\oss where\sss the vertical\sss arrows $\beta$ 
are\sss the\dss Bott\dss periodicity\sss maps.\oss

The coboundary map\sss
$\partial\dff \colon\dff
K^{\dff -\dff 1}\fff(\trf Y\trf)
\qff \ttoo\qff
K^{\dff 0}\fff(\trf X\fff,\qff Y\trf)$\sss
admits an alternative description.\oss
Every\sss element\sss of\dss
$K^{\dff -\dff 1}\fff(\trf X\trf)
\off =\off
\widetilde{K}^{\trf 0}\fff(\trf \bm{\Sigma}\dff X\trf)$\dss
has\sss the form\sss
$[\trf E\trf]\qff -\qff [\dff e\trf]$\nnsp,\oss
where\sss $E$\sss is\dss a vector bundle over\sss $\bm{\Sigma}\dff Y$
and\sss $[\dff e\trf]$\sss is\dss the\sss trivial\dss bundle of\dss
dimension $e\off =\off \dim\dff E$\nnsp.\oss
The bundle $E$ can be obtained by glueing two
trivial\dss bundles of dimension ${\nsp}e{\nsp}$ over two halves of\sss $\bm{\Sigma}\dff Y$\sss
by a map 
$\varphi\fff \colon\fff 
Y\dff \ttoo\dff GL\dff(\dff e\trf)$\nnsp.\oss
Then\sss the\sss lift\sss of\sss $E$\sss to\sss 
$X\qff \cup\qff \mathbf{C}\dff Y$\sss also results from glueing\sss
trivial\dss bundles of\dss dimension $e$ over\sss
$X\qff \cup\qff (\qff Y\qff \times\qff [\trf 1/2\fff,\qff 1\trf] \trf)$\sss
and\sss 
$Y\qff \times\qff [\trf 0\fff,\qff 1/2\trf]\fff /(\trf Y\dff \times\dff 0\trf)$\sss
by\sss the map
$y\dff \times\dff 1/2
\off \longmapsto\off
\varphi\dff(\trf y\trf)$\nnsp.\oss
It\dss is\dss easy\sss to see\sss that\sss
the image of\dss $[\trf E\trf]\qff -\qff [\dff e\trf]$\sss
in\sss the group $K^{\dff 0}\fff(\trf X\fff,\qff Y\trf)$\sss
is\dss equal\dss to\sss the difference construction
$d\trf(\trf [\dff e\trf]\fff,\qff [\dff e\trf]\fff,\qff \varphi\trf)$\nnsp,\qss
where now $[\dff e\trf]$\sss is\sss the\sss trivial\dss bundle of\dss
dimension $e$ over $X$\nnsp.

\myuppar{Self-adjoint\sss endomorphisms.}
Suppose\sss that\sss $E$\sss is\dss a bundle over $X$ equipped with a\dss Hermitian\dss metric and
$p\dff \colon\dff E\qff \ttoo\qff E$\sss
is\dss a self-adjoint\sss endomorphism of\sss $E$ such\sss that $p$\sss
is\dss an\sss isomorphism over $Y$\dnsp.\oss
Recall\dss that\sss $\mathbb{P}^{\dff +}\fff (\trf p\trf |\trf Y\trf)$\sss is\dss the subbundle of\sss $E\trf|\trf Y$
having as\sss fibers\sss the sums of\dss eigenspaces of\sss $p$\sss 
corresponding\sss to\sss the positive eigenvalues of\sss $p$\nnsp.\oss 
By an abuse of\dss notations we will\sss denote by\sss 
$\mathbb{P}^{\dff +}\fff (\trf p\trf |\trf Y\trf)$\sss
also\sss the corresponding element\sss of\sss $K^{\dff 0}\fff(\trf Y\trf)$\nnsp.\oss
We are interested\sss in\sss\vspace{1.5pt}
\[
\quad
\partial\qff \mathbb{P}^{\dff +}\fff (\trf p\trf |\trf Y\trf)
\off \in\off
K^{\dff 1}\fff(\trf X\fff,\qff Y\trf)
\qff.
\]

\vspace{-12pt}\vspace{1.5pt}
Let\sss us consider\sss the family of\dss endomorphisms\sss
$p_{\trf \eta}\dff \colon\dff E\qff \ttoo\qff E$\nnsp,\qss 
$\eta\qff \in\qff [\trf 0\fff,\qff 2\dff \pi\trf]$
defined\dss by\vspace{3pt}
\[
\quad
p_{\trf \eta}
\off =\off
\id\dff \cos\dff \eta
\qff +\qff
p\trf i\dff \sin\dff \eta
\quad\pff\halfff
\mbox{for}\quad
\eta\qff \in\qff [\trf 0\fff,\qff \pi\trf]
\qff,
\]

\vspace{-34.5pt}
\[
\quad
p_{\trf \eta}
\off =\off
\id\qff
(\trf \cos\dff \eta
\qff +\qff
i\dff \sin\dff \eta\trf)
\quad
\mbox{for}\quad
\eta\qff \in\qff [\trf \pi\fff,\qff 2\dff \pi\trf]
\qff.
\]

\vspace{-12pt}\vspace{3pt}
Clearly,\pss $p_{\trf \eta}$\sss are isomorphisms over $Y$\dnsp,\oss
and\sss $p_{\trf \eta}$\sss is\dss the identity\sss isomorphism of\sss $E$ over $X$\sss
for\sss $\eta\off =\off 0$\sss or\sss $2\dff \pi$\nnsp,\oss
and\dss hence\sss this family defines an endomorphism\vspace{1.5pt}
\[
\quad
\overline{p}\trf \colon\dff
E\dff \times\dff [\trf 0\fff,\qff 2\dff \pi\trf]
\qff \ttoo\qff
E\dff \times\dff [\trf 0\fff,\qff 2\dff \pi\trf]
\]

\vspace{-12pt}\vspace{1.5pt}
of\dss the bundle\sss 
$\overline{E}
\off =\off
E\dff \times\dff [\trf 0\fff,\qff 2\dff \pi\trf]$\sss
over\sss
$X\dff \times\dff [\trf 0\fff,\qff 2\dff \pi\trf]$\sss
such\sss that\sss $\overline{p}$\sss is\dss an\sss isomorphism over\sss\vspace{1.5pt}
\[
\quad
Y\dff \times\dff [\trf 0\fff,\qff 2\dff \pi\trf]
\qff \cup\qff
X\dff \times\dff \{\trf 0\fff,\qff 2\dff \pi\trf\}
\qff.
\]

\vspace{-12pt}\vspace{1.5pt}
Therefore\sss $\overline{p}$\sss defines\sss the difference class\vspace{3pt}
\[
\quad
d\trf\left(\qff \overline{p}\pff\right)
\off =\off
d\trf\left(\pff \overline{E}\fff,\off \overline{E}\fff,\off \overline{p} \pff\right)
\off \in\off
K^{\dff 0}\fff\bigl(\qff
X\dff \times\dff [\trf 0\fff,\qff 2\dff \pi\trf]\fff,\off
Y\dff \times\dff [\trf 0\fff,\qff 2\dff \pi\trf]
\qff \cup\qff
X\dff \times\dff \{\trf 0\fff,\qff 2\dff \pi\trf\}
\qff\bigr)
\]

\vspace{-33pt}
\[
\quad
\phantom{d\trf\left(\qff \overline{p}\pff\right)
\off =\off
d\trf\left(\pff \overline{E}\fff,\off \overline{E}\fff,\off \overline{p} \pff\right)
\off }
=\off
K^{\dff 0}\fff(\trf
\bm{\Sigma}\dff X\fff,\off
\bm{\Sigma}\dff Y
\trf)
\off =\off
K^{\dff -\dff 1}\fff(\trf X\fff,\qff Y\trf)
\qff.
\]

\vspace{-12pt}\vspace{3pt}
\mypar{Lemma.}{cobound-p-plus}
\emph{Let\qss
$\beta\dff \colon\dff
K^{\dff 1}\fff(\trf X\fff,\qff Y\trf)
\qff \ttoo\qff
K^{\dff -\dff 1}\fff(\trf X\fff,\qff Y\trf)$\sss
be\sss the\dss Bott\dss periodicity\sss map.\oss
Then}
\vspace{1.5pt}
\[
\quad
\beta\trf\left(\qff
\partial\qff\fff \mathbb{P}^{\dff +}\fff (\trf p\trf |\trf Y\trf)
\qff\right)
\off =\off\dff
d\trf\left(\qff \overline{p}\pff\right)
\qff.
\]

\vspace{-12pt}\vspace{1.5pt}
\proof
Let $E\fff'$\sss be a bundle over $X$ such\sss that $E\dff \oplus\dff E\fff'$\sss
is\dss a\sss trivial\dss bundle.\oss
Let\sss $p\fff'\dff \colon\dff E\fff'\qff \ttoo\qff E\fff'$\sss
be equal\dss to\sss $-\qff \id$\nnsp.\oss
Replacing\sss $E$\sss by $E\dff \oplus\dff E\fff'$
and $p$\sss by\sss $p\dff \oplus\dff p\fff'$\sss obviously
does not\sss change\sss $\mathbb{P}^{\dff +}\fff (\trf p\trf |\trf Y\trf)$\nnsp.\oss
Clearly,\pss
$p\fff'_{\dff \pi\dff +\dff \eta}
\off =\off
p\fff'_{\dff \pi\dff -\dff \eta}$\sss
and\dss hence\sss the\sss loop
$p\fff'_{\dff \eta}$,\qss $\eta\qff \in\qff [\trf 0\fff,\qff 2\dff \pi\trf]$\sss
is\dss homotopic\sss with\sss the endpoints fixed\dss to\sss the constant\dss loop
in\sss the class of\dss loops of\dss automorphisms of\dss $E\fff'$\dnsp.\oss
It\sss follows\sss that\dss replacing\sss $E$\sss by $E\dff \oplus\dff E\fff'$
and $p$\sss by\sss $p\dff \oplus\dff p\fff'$\sss
does not\sss changes\sss the class\vspace{0pt}
\[
\quad
d\trf\bigl(\trf \overline{p}\pff\bigr)
\off =\off
d\trf\bigl(\pff \overline{E}\fff,\off \overline{E}\fff,\off \overline{p} \pff\bigr)
\]

\vspace{-12pt}\vspace{0pt}
either.\oss
It\sss follows\sss that\sss we may assume\sss that\sss $E$\sss is\dss a\sss trivial\dss bundle.\oss

Deforming\sss $p$\sss to its unitary\sss part\sss $p\trf \num{p}^{\dff -\dff 1}$\sss also
does not\sss affects our $K$\dnsp-theory classes.\oss
Therefore we can assume\sss that\sss $p$\sss over $Y$\sss is\dss unitary.\oss
Then\sss $p\off =\off p_{\dff +}\qff -\qff p_{\dff -}$ over $Y$\dnsp,\oss
where $p_{\dff +}$\sss is\dss the orthogonal\sss projection of\dss $E$\sss onto\sss
the subbundle\sss $\mathbb{P}^{\dff +}\fff (\trf p\trf |\trf Y\trf)$\nnsp,\oss
and $p_{\dff -}$\sss is\dss the orthogonal\sss projection onto its orthogonal\sss complement.\oss
Also,\qss 
$p_{\dff +}\qff +\qff p_{\dff -}\off =\off \id$ and\vspace{1.5pt}
\[
\quad
p_{\trf \eta}
\off =\off
\id\dff \cos\dff \eta
\pff +\pff
p\dff i\dff \sin\dff \eta
\]

\vspace{-36pt}
\[
\quad
\phantom{p_{\trf \eta}
\off }
=\off
(\trf p_{\dff +}\qff +\qff p_{\dff -}\trf)\trf \cos\dff \eta
\pff +\pff
(\trf p_{\dff +}\qff -\qff p_{\dff -}\trf)\qff i\dff \sin\dff \eta
\quad
\]

\vspace{-36pt}
\[
\quad
\phantom{p_{\trf \eta}
\off }
=\off
p_{\dff +}\qff (\trf \cos\dff \eta\qff +\qff i\dff \sin\dff \eta\trf)
\pff +\pff
p_{\dff -}\qff (\trf \cos\dff \eta\qff -\qff i\dff \sin\dff \eta\trf)
\quad
\]

\vspace{-36pt}
\[
\quad
\phantom{p_{\trf \eta}
\off }
=\off
p_{\dff +}\qff e^{\dff i\trf \eta}
\pff +\pff
p_{\dff -}\qff e^{\dff -\dff i\trf \eta}
\quad
\]

\vspace{-12pt}\vspace{1.5pt}
for $\eta\qff \in\qff [\trf 0\fff,\qff \pi\trf]$\nnsp.\oss
It\sss follows\sss that\sss in every\dss fiber\sss $E_{\dff y}$\nsp,\dss $y\qff \in\qff Y$\sss the path
$p_{\trf \eta}$\nsp,\dss $\eta\qff \in\qff [\trf 0\fff,\qff 2\dff \pi\trf]$\sss
is\dss essentially\dss the path of\oss isomorphisms\dss
$E_{\dff y}\qff \ttoo\qff E_{\dff y}$\dss assigned\dss by\qss Bott\qss \cite{bott}\qss
to\sss the subspace\dss $\mathbb{P}^{\dff +}\fff (\trf p\trf)_{\dff y}$\sss of\dss $E_{\dff y}$\nsp.\oss
The part\sss $p_{\trf \eta}$\nsp,\dss $\eta\qff \in\qff [\trf \pi\fff,\qff 2\dff \pi\trf]$\sss
differs from\dss Bott's\dss part,\oss but\sss the only\sss thing\sss that\sss matters\dss is\dss
that\dss this\dss is\dss a standard\sss path connecting\sss $-\qff \id$ with $\id$\nnsp.\oss

Let\sss us consider\sss the family of\dss subbundles\sss
$E_{\trf \eta\fff,\qff \theta}\qff \subset\qff E\dff \oplus\dff E$\nnsp,\qss 
$\eta\fff,\qff \theta\qff \in\qff [\trf 0\fff,\qff 2\dff \pi\trf]$
defined\dss by\vspace{3pt}
\[
\quad
E_{\trf \eta\fff,\qff \theta}
\off =\off
\bigl\{\pff 
\bigl(\qff u\dff \cos\dff \theta/2\fff,\off 
p_{\trf \eta}\trf(\dff u\trf)\trf \sin\dff \theta/2\qff\bigr)
\pff\fff \bigl|\pff\fff
u\qff \in\qff E
\pff\bigr\}
\quad
\mbox{for}\quad
\theta\qff \in\qff [\trf 0\fff,\qff \pi\trf]
\qff,
\]

\vspace{-34.5pt}
\[
\quad
E_{\trf \eta\fff,\qff \theta}
\off =\off
\bigl\{\pff 
\bigl(\qff u\dff \cos\dff \theta/2\fff,\off 
u\trf \sin\dff \theta/2\qff\bigr)
\pff\fff \bigl|\pff\fff
u\qff \in\qff E
\pff\bigr\}
\quad
\mbox{for}\quad
\theta\qff \in\qff [\trf \pi\fff,\qff 2\dff \pi\trf]
\qff.
\]

\vspace{-12pt}\vspace{3pt}
This family defines a subbundle $\mathbf{E}$ of\sss
$(\trf E\dff \oplus\dff E\trf)\dff \times\dff I$
over\sss
$Y\dff \times\dff I$\nnsp,\oss
where\sss 
$I\off =\off [\trf 0\fff,\qff 2\dff \pi\trf]\dff \times\dff [\trf 0\fff,\qff 2\dff \pi\trf]$\nnsp.\oss
Let\sss $\mathbf{E_{\trf 0}}$ be defined\sss by\sss the same formulas,\oss but\sss with
$p_{\trf \eta}\trf(\dff u\trf)$ replaced\sss by $u$\nnsp.\oss
Clearly,\pss $\mathbf{E}$\sss is\dss equal\dss to\sss 
$\mathbf{E_{\trf 0}}$\sss 
over\sss the boundary\sss $\partial\trf I$\nnsp.\oss
The difference class
$d\trf(\trf \mathbf{E}\fff,\qff \mathbf{E_{\trf 0}}\dff,\qff \id\trf)$\sss 
is\dss an element\sss\vspace{3pt}
\[
\quad
\mathbf{E_{\trf \bm{\beta}}}
\off \in\off
K^{\dff 0}\fff(\trf
Y\dff \times\dff I\fff,\off
Y\dff \times\dff \partial\trf I
\trf)
\off =\off
K^{\dff -\dff 2}\fff(\trf Y\trf)
\off =\off
K^{\dff -\dff 1}\fff(\trf \bm{\Sigma}\dff Y\trf)
\qff.
\]

\vspace{-12pt}\vspace{3pt}
By\sss the construction,\oss
this\dss is\dss the image $\beta\qff \mathbb{P}^{\dff +}\fff (\trf p\trf |\trf Y\trf)$ 
of\sss $\mathbb{P}^{\dff +}\fff (\trf p\trf |\trf Y\trf)$
under\sss the\dss Bott\dss periodicity\sss map\sss $\beta$\nnsp.\oss
There are obvious isomorphisms of\dss $\mathbf{E}$\sss and\sss 
$(\trf E\dff \oplus\dff 0\trf)\dff \times\dff I$\sss over 
$Y\dff \times\dff [\trf 0\fff,\qff 2\dff \pi\trf]\dff \times\dff [\trf 0\fff,\qff \pi\trf]$\sss
and over\sss 
$Y\dff \times\dff [\trf 0\fff,\qff 2\dff \pi\trf]\dff \times\dff [\trf \pi\fff,\qff 2\dff \pi\trf]$\nnsp.\oss
They differ over\sss 
$Y\dff \times\dff [\trf 0\fff,\qff 2\dff \pi\trf]\dff \times\dff \pi$\sss
by\sss the endomorphism\sss $\overline{p}$\nnsp,\oss
and\dss hence\sss $\mathbf{E}$\sss is\dss isomorphic\sss to\sss the result\sss of\dss
glueing of\dss two\sss trivial\dss bundles by\sss $\overline{p}$\nnsp.\oss
The above alternative description of\dss the maps\sss $\partial$\sss with
$(\trf \bm{\Sigma}\dff X\fff,\off \bm{\Sigma}\dff Y \trf)$
in\sss the role of\dss $(\trf X\fff,\off Y \trf)$
implies\sss that\vspace{3pt}\vspace{-0.4pt}
\[
\quad
\partial\trf \mathbf{E_{\trf \bm{\beta}}}
\off =\off
d\trf\bigl(\pff \overline{E}\fff,\off \overline{E}\fff,\off \overline{p} \pff\bigr)
\off \in\off
K^{\dff 0}\fff(\trf
\bm{\Sigma}\dff X\fff,\off
\bm{\Sigma}\dff Y
\trf)
\off =\off
K^{\dff -\dff 1}\fff(\trf X\fff,\qff Y\trf)
\qff.
\]

\vspace{-12pt}\vspace{3pt}\vspace{-0.4pt}
Since\sss the\dss Bott\dss periodicity commutes with\sss the coboundary\sss maps,\oss
it\sss follows\sss that\sss\vspace{3pt}\vspace{-0.4pt}
\[
\quad
\beta\dff\left(\qff \partial\qff \mathbb{P}^{\dff +}\fff (\trf p\trf |\trf Y\trf)\qff\right)
\off =\off
\partial\dff\left(\qff \beta\qff \mathbb{P}^{\dff +}\fff (\trf p\trf |\trf Y\trf)\qff\right)
\off =\off
\partial\trf \mathbf{E_{\trf \bm{\beta}}}
\off =\off
d\trf\bigl(\pff \overline{E}\fff,\off \overline{E}\fff,\off \overline{p} \pff\bigr)
\qff.
\]

\vspace{-12pt}\vspace{3pt}\vspace{-0.4pt}
This completes\sss the proof\dss of\dss the\sss lemma.\oss  \eproof

\myuppar{Redefining\sss $\overline{p}$\nnsp.}
The endomorphism\sss $\overline{p}$\sss is\dss an\dss isomorphism over\sss 
$X\dff \times\dff [\trf \pi\fff,\qff 2\dff \pi\trf]$\nnsp.\oss
Hence\sss
$d\trf(\trf \overline{p}\pff)$\sss is\dss equal\dss to $0$ over
$X\dff \times\dff [\trf \pi\fff,\qff 2\dff \pi\trf]$
and after\sss identifying
$[\trf 0\fff,\qff 2\dff \pi\trf]$
with $[\trf 0\fff,\qff \pi\trf]$
we may\sss replace\sss the endomorphism\sss $\overline{p}$\sss
by\sss its restriction\sss to\sss
$X\dff \times\dff [\trf 0\fff,\qff \pi\trf]$\sss
without\sss affecting\sss the difference class\vspace{3pt}
\[
\quad
d\trf\left(\trf \overline{p}\pff\right)
\off \in\off
K^{\dff -\dff 1}\fff(\trf X\fff,\qff Y\trf)
\qff.
\]

\vspace{-12pt}\vspace{3pt}
Let\sss us\sss redefine\sss $\overline{E}$\sss
as\sss $E\dff \times\dff [\trf 0\fff,\qff \pi\trf]$\sss
and\sss redefine\sss $\overline{p}$\sss
as\sss the endomorphism\sss
$\overline{E}\qff \ttoo\qff \overline{E}$\sss
corresponding\sss to\sss the family\sss
$p_{\trf \eta}\dff,\off \eta\qff \in\qff [\trf 0\fff,\qff \pi\trf]$\nnsp.\oss
Clearly,\oss Lemma\qss \ref{cobound-p-plus}\qss remains valid.\oss

\myuppar{Multiplication\sss in\sss relative $K$\dnsp-groups.}
Suppose first\sss that\sss
$p\dff \colon\dff E\qff \ttoo\qff F$\sss
and\sss
$q\dff \colon\dff E\fff'\qff \ttoo\qff F\fff'$\sss
are\sss linear maps of\trs Hermitian\dss vector spaces.\oss
Let\sss $r\qff \in\qff \zzz$\nnsp.\oss
Then\sss the matrix\vspace{3pt}
\begin{equation}
\label{mult-matrix}
\quad
\begin{pmatrix}
\off\dff p\dff \otimes\dff 1 &
(\trf -\qff 1\trf)^{\dff r\dff +\dff 1}\dff 1\dff \otimes\trf q^{\fff *} \off
\vspace{6pt} \\
\off\dff 1\dff \otimes\trf q &
(\trf -\qff 1\trf)^{\dff r}\dff p^{\fff *}\dff \otimes\dff 1 \off 
\end{pmatrix}
\qff 
\end{equation}

\vspace{-12pt}\vspace{3pt}
defines a\sss linear map\vspace{3pt}
\[
\quad
p\pff \omult_{\fff r}\qff q\qff \colon\qff 
\left(\trf E\dff \otimes\dff E\fff' \trf\right)
\qff \oplus\qff
\left(\trf F\dff \otimes\dff F\fff' \trf\right)
\qff \ttoo\qff
\left(\trf F\dff \boxtimes\dff E\fff' \trf\right)
\qff \oplus\qff
\left(\trf E\dff \boxtimes\dff F\fff' \trf\right)
\qff.
\]

\vspace{-12pt}\vspace{3pt}
One immediately\sss checks\sss that\vspace{3pt}
\[
\quad
\left(\trf p\qff +\qff p\fff' \trf\right)\pff \omult_{\fff r}\qff q
\off =\off
\left(\trf p\pff \omult_{\fff r}\qff 0 \trf\right)
\off +\off
\left(\trf p\fff'\pff \omult_{\fff r}\qff q \trf\right)
\quad
\mbox{and}
\qff
\]

\vspace{-34.5pt}
\[
\quad
\left(\trf \lambda\dff p\trf\right)\pff \omult_{\fff r}\qff \left(\trf \lambda\dff q\trf\right)
\off =\off
\lambda\dff \left(\trf p\pff \omult_{\fff r}\qff q \trf\right)
\qff
\]

\vspace{-12pt}\vspace{3pt}
for every\sss $\lambda\qff \in\qff \rrr$\nnsp.\oss
We will\sss not\sss use\sss these\qss ``linearity''\qss 
properties explicitly,\oss preferring direct\sss calculations.\oss
It\dss is\dss well\dss known\sss that\sss
$p\pff \omult_{\fff r}\qff q$\sss
is\dss invertible\sss if\dss either $p$ or $q$\sss is.\oss\vspace{-0pt}

Suppose now\sss that\sss $E\fff,\off F$\sss are vector bundles over a\sss
topological\sss space $X$\nnsp,\oss and\sss
$E\fff'\fff,\off F\fff'$\sss are vector bundles over a\sss
topological\sss space $X\fff'$\dnsp.\oss
Recall\dss that\sss the external\dss tensor product\sss
$E\dff \boxtimes\dff E\fff'$\sss is\dss a vector bundle over\sss $X\dff \times\dff X\fff'$\sss
having\sss the\sss tensor product\sss $E_{\dff x}\dff \otimes\dff E\fff'_{\fff y}$\sss as\sss the fiber
over\sss $(\trf x\fff,\qff y\trf)\qff \in\qff X\dff \times\dff X\fff'$\dnsp.\oss 
Suppose\sss that\sss $E\fff,\qff E\fff'$\dnsp,\oss etc.\qss are\dss Hermitian\sss vector bundles
and\sss \sss
$p\dff \colon\dff E\qff \ttoo\qff F$\sss
and\sss
$q\dff \colon\dff E\fff'\qff \ttoo\qff F\fff'$\sss
are vector\sss bundle maps.\oss
Then\sss the matrix\qss (\ref{mult-matrix})\qss defines a vector\sss bundle map\vspace{3.5pt}
\[
\quad
p\pff \omult_{\fff r}\qff q\qff \colon\qff 
\bigl(\trf E\dff \boxtimes\dff E\fff' \trf\bigr)
\qff \oplus\qff
\left(\trf F\dff \boxtimes\dff F\fff' \trf\right)
\qff \ttoo\qff
\left(\trf F\dff \boxtimes\dff E\fff' \trf\right)
\qff \oplus\qff
\left(\trf E\dff \boxtimes\dff F\fff' \trf\right)
\qff.
\]

\vspace{-12pt}\vspace{3.5pt}
Recall\dss that\sss every element\sss of\dss 
$K^{\dff 0}\fff(\trf X\fff,\qff Y\trf)$ can\sss be represented
as\sss the difference construction
$d\trf(\trf E\fff,\qff F\fff,\qff p\trf)$\nnsp,\oss
where\sss $E\fff,\qff F$\sss are vector bundles on\sss $X$\sss
and\sss $p\dff \colon\dff E\qff \ttoo\qff F$\sss
is\dss map of\dss vector bundles inducing an\sss isomorphism\sss
$E\trf|\trf Y\qff \ttoo\qff F\trf|\trf Y$\dnsp.\oss
For a subspace\sss $Y\fff'\qff \subset\qff X\fff'$\sss
the product\sss of\dss pairs\sss
$(\trf X\fff,\qff Y\trf)\dff \times\dff (\trf X\fff'\fff,\qff Y\fff'\qff)$\sss
is\dss defined as\sss the pair\sss
$(\trf X\dff \times\dff X\fff',\off 
X\dff \times\dff Y\fff'\trf \cup\qff Y\dff \times\dff X\fff'\qff)$\nnsp.\oss
The external\dss tensor product\vspace{3.5pt}
\[
\quad
\widehat{\otimes}\dff \colon\dff
K^{\dff 0}\fff(\trf X\fff,\qff Y\trf)
\dff \otimes\dff
K^{\dff 0}\fff(\trf X\fff'\fff,\qff Y\fff'\trf)
\qff \ttoo\qff
K^{\dff 0}\fff(\trf X\dff \times\dff X\fff',\off 
X\dff \times\dff Y\fff'\trf \cup\qff Y\dff \times\dff X\fff'\trf)
\qff
\]

\vspace{-12pt}\vspace{3.5pt}
can\sss be described\sss as follows.\oss
Let\sss 
$d\trf(\trf E\fff,\qff F\fff,\qff p\trf)
\qff \in\qff K^{\dff 0}\fff(\trf X\fff,\qff Y\trf)$
and\sss
$d\trf(\trf E\fff'\fff,\qff F\fff'\fff,\qff q\trf)
\in\qff K^{\dff 0}\fff(\trf X\fff'\fff,\qff Y\fff'\trf)$\nnsp.\oss
Let\sss us equip each of\dss the bundles\sss 
$E\fff,\pff E\fff'\fff,\off F\fff,\off F\fff'$\sss 
with a\dss Hermitian\sss metic.\oss
Then\sss the map\sss $p\pff \omult_{\fff r}\qff q$\sss is\dss defined and\dss is\dss an isomorphism over\sss
$X\dff \times\dff Y\fff'\trf \cup\qff X\fff'\dff \times\dff Y$\dnsp.\oss
The difference construction applied\sss to
$p\pff \omult_{\fff r}\qff q$ defines an element\sss
$d\trf(\trf p\pff \omult_{\fff r}\qff q\trf)
\qff \in\pff
K^{\dff 0}\fff(\trf X\dff \times\dff X\fff',\off 
X\dff \times\dff Y\fff'\trf \cup\qff Y\dff \times\dff X\fff'\trf)$\sss
and\vspace{4.5pt}
\[
\quad
d\trf(\trf E\fff,\qff F\fff,\qff p\trf)
\qff \widehat{\otimes}\qff
d\trf(\trf E\fff'\fff,\qff F\fff'\fff,\qff q\trf)
\off =\off
d\trf(\trf p\pff \omult_{\fff r}\qff q\trf)
\]

\vspace{-12pt}\vspace{4.5pt}
for every $r$\nnsp.\oss
See,\oss for example,\oss \cite{pa},\qss Section\qss II.5.\oss

\myuppar{The case of\dss self-adjoint $p$\nnsp.}
Suppose now\sss that\sss $E\off =\off F$\sss
and $p$\sss is\dss self-adjoint.\oss
Then\sss the matrix\qss (\ref{mult-matrix})\qss is\dss self-adjoint\sss for $r\off =\off 1$\nnsp.\oss
In\sss this case we will\sss relate\sss $p\pff \omult_1\qff q$\sss
with\sss the product\vspace{4.5pt}
\[
\quad
K^{\dff 1}\fff(\trf X\fff,\qff Y\trf)
\dff \otimes\dff
K^{\dff 0}\fff(\trf X\fff'\fff,\qff Y\fff'\trf)
\qff \ttoo\qff
K^{\dff 1}\fff(\trf X\dff \times\dff X\fff',\off 
X\dff \times\dff Y\fff'\trf \cup\qff Y\dff \times\dff X\fff'\trf)
\qff.
\]

\vspace{-12pt}\vspace{4.5pt}
\mypar{Lemma.}{coboundary-mult}
$\dis
\partial\qff\fff 
\mathbb{P}^{\dff +}\fff (\trf p\trf |\trf Y\trf)
\off \widehat{\otimes}\off\dff
d\trf(\trf E\fff'\fff,\qff F\fff'\fff,\qff q\trf)
\off\dff =\off\qff
\partial\qff\fff 
\mathbb{P}^{\dff +}\fff 
\left(\qff 
p\pff \omult_1\qff q
\qff\fff\bigl|\pff
X\dff \times\dff Y\fff'\trf \cup\qff Y\dff \times\dff X\fff'
\qff\right)$\nnsp.\oss\vspace{1.5pt}

\proof
It\dss is\dss sufficient\sss to prove\sss that\sss 
the\dss Bott\dss periodicity\sss map\sss\vspace{4.5pt}
\[
\quad
\beta\dff \colon\dff
K^{\dff 1}\fff(\trf X\dff \times\dff X\fff',\off 
X\dff \times\dff Y\fff'\trf \cup\qff Y\dff \times\dff X\fff'\qff)
\qff \ttoo\qff
K^{\trf -\dff 1}\fff(\trf X\dff \times\dff X\fff',\off 
X\dff \times\dff Y\fff'\trf \cup\qff Y\dff \times\dff X\fff'\qff)
\qff
\]

\vspace{-12pt}\vspace{4.5pt}
takes both sides of\dss this equality\sss to\sss the same element.\oss
The map $\beta$\sss can\sss be defined as\sss the external\dss product\sss
with\sss an element\sss of\sss
$K^{\dff 0}\fff(\trf S^{\dff 2}\fff,\qff \{\trf s_{\dff 0}\qff\}\trf)$\nnsp,\oss
where\sss $s_{\dff 0}\qff \in\qff S^{\dff 2}$\sss is\dss the base point,\oss
namely\sss with\sss the\dss Bott\dss element.\oss
Therefore\sss the associativity of\dss 
the external\dss product\sss implies\sss that\vspace{4.5pt}
\[
\quad
\beta\qff
\left(\qff
\partial\qff 
\left(\qff
\mathbb{P}^{\dff +}\fff (\trf p\trf |\trf Y\trf)
\qff\right)
\off \widehat{\otimes}\off
d\trf(\trf E\fff'\fff,\qff F\fff'\fff,\qff q\trf)
\qff\right)
\off =\off
\beta\dff \circ\trf \partial\qff
\left(\qff
\mathbb{P}^{\dff +}\fff (\trf p\trf |\trf Y\trf)
\qff\right)
\off \widehat{\otimes}\off
d\trf(\trf E\fff'\fff,\qff F\fff'\fff,\qff q\trf)
\qff.
\]

\vspace{-12pt}\vspace{4.5pt}
Together with\dss Lemma\qss \ref{cobound-p-plus}\qss this implies\sss that\sss it\dss is\dss
sufficient\sss to prove\sss that\vspace{4.5pt}
\begin{equation}
\label{tensor-mult}
\quad
d\trf\left(\qff \overline{p}\pff\right)
\off \widehat{\otimes}\off
d\trf(\trf E\fff'\fff,\qff F\fff'\fff,\qff q\trf)
\off =\off
d\trf\left(\off \overline{p\pff \omult_{\fff 1}\qff q} \off\right)
\qff.
\end{equation}

\vspace{-12pt}\vspace{4.5pt}
The external\sss product\dss
$d\trf\left(\qff \overline{p}\pff\right)
\off \widehat{\otimes}\off
d\trf(\trf E\fff'\fff,\qff F\fff'\fff,\qff q\trf)$\nnsp,\oss 
considered as an element of\vspace{4.5pt}
\[
\quad
K^{\dff 0}\fff\left(\qff X\dff \times\dff [\trf 0\fff,\qff \pi\trf]\dff \times\dff X\fff',\off\off 
X\dff \times\dff [\trf 0\fff,\qff \pi\trf]\dff \times\dff Y\fff'\qff \cup\off 
Y\dff \times\dff [\trf 0\fff,\qff \pi\trf]\dff \times\dff X\fff'
\off \cup\off
X\dff \times\dff \{\trf 0\fff,\qff \pi\trf\}\dff \times\dff X\fff'
\qff\right)
\qff,
\]

\vspace{-12pt}\vspace{4.5pt}
is\dss equal\dss to\sss the difference class 
$d\trf\left(\qff \overline{p}\pff \omult_{\dff 0}\qff q\trf\right)$\nnsp.\oss
Clearly,\pss\vspace{4.5pt}
\[
\quad
d\trf(\trf E\fff'\fff,\qff F\fff'\fff,\qff q\trf)
\off =\off
d\trf(\trf E\fff'\fff,\qff F\fff'\fff,\qff i\dff q\trf)
\]

\vspace{-12pt}\vspace{4.5pt}
and\dss hence\sss this external\dss product\sss also equal\dss to\sss
$d\trf\left(\qff \overline{p}\pff \omult_{\dff 0}\qff (\trf i\dff q\trf)\trf\right)$\nnsp.\oss

By\sss the definition\vspace{3pt}
\[
\quad
p_{\trf \eta}\pff \omult_{\dff 0}\qff (\trf i\dff q\trf)
\off =\off
(\trf \id\dff \cos\dff \eta\pff +\pff i\dff p\dff \sin\dff \eta \trf)
\pff \omult_{\dff 0}\qff (\trf i\dff q\trf)
\]

\vspace{-30pt}
\[
\quad
\phantom{p_{\trf \eta}\pff \omult_{\dff 0}\qff (\trf i\dff q\trf)
\off }
=\off
\begin{pmatrix}
\off\dff (\trf \id\dff \cos\dff \eta\pff +\pff i\dff p\dff \sin\dff \eta \trf)\dff \otimes\dff 1 &
i\dff \otimes\trf q^{\fff *} \off
\vspace{6pt} \\
\off\dff i\dff \otimes\trf q &
(\trf \id\dff \cos\dff \eta\pff -\pff i\dff p\dff \sin\dff \eta \trf)\dff \otimes\dff 1 \off 
\end{pmatrix}
\]

\vspace{-15pt}
\[
\quad
\phantom{p_{\trf \eta}\pff \omult_{\dff 0}\qff (\trf i\dff q\trf)
\off }
=\off
\cos\dff \eta\qff 
\begin{pmatrix}
\off\dff 1\dff \otimes\dff 1 &
0 \off
\vspace{6pt} \\
\off\dff 0 &
1\dff \otimes\dff 1 \off 
\end{pmatrix}
\off +\off
i\pff
\begin{pmatrix}
\off\dff p\dff \sin\dff \eta\dff \otimes\dff 1 &
1\dff \otimes\trf q^{\fff *} \off
\vspace{6pt} \\
\off\dff 1\dff \otimes\trf q &
-\qff p\dff \sin\dff \eta\dff \otimes\dff 1 \off 
\end{pmatrix}
\off.
\]

\vspace{-12pt}
At\sss the same\sss time\vspace{3pt}
\[
\quad
\left(\qff \overline{p\pff \omult_{\fff 1}\qff q} \qff\right)_{\trf \eta}
\off =\off
\cos\dff \eta\qff 
\begin{pmatrix}
\off\dff 1\dff \otimes\dff 1 &
0 \off
\vspace{6pt} \\
\off\dff 0 &
1\dff \otimes\dff 1 \off 
\end{pmatrix}
\off +\off
i\dff \sin\dff \eta\pff
\begin{pmatrix}
\off\dff p\dff \otimes\dff 1 &
1\dff \otimes\trf q^{\fff *} \off
\vspace{6pt} \\
\off\dff 1\dff \otimes\trf q &
-\qff p\dff \otimes\dff 1 \off 
\end{pmatrix}
\]

\vspace{-15pt}
\[
\quad
\phantom{\left(\qff \widetilde{p\pff \omult_{\fff 1}\qff q} \qff\right)_{\trf \eta}
\off }
=\off
\cos\dff \eta\qff 
\begin{pmatrix}
\off\dff 1\dff \otimes\dff 1 &
0 \off
\vspace{6pt} \\
\off\dff 0 &
1\dff \otimes\dff 1 \off 
\end{pmatrix}
\off +\off
i\pff
\begin{pmatrix}
\off\dff p\dff \sin\dff \eta\dff \otimes\dff 1 &
\sin\dff \eta\dff \otimes\trf q^{\fff *} \off
\vspace{6pt} \\
\off\dff \sin\dff \eta\dff \otimes\trf q &
-\qff p\dff \sin\dff \eta\dff \otimes\dff 1 \off 
\end{pmatrix}
\off.
\]

\vspace{-12pt}\vspace{6pt}
For\sss $t\qff \in\qff [\trf 0\fff,\qff 1\trf]$\dss
let\sss
$\sin_{\dff t}\dff \eta
\off =\off
(\trf 1\qff -\qff t\trf)\qff +\qff t\dff \sin\dff \eta$\dss
and\vspace{6pt}
\[
\quad
\left(\qff \overline{p\pff \omult_{\fff 1}\qff q} \qff\right)_{\trf \eta\dff,\dff t} 
\off =\off
\cos\dff \eta\qff 
\begin{pmatrix}
\off\dff 1\dff \otimes\dff 1 &
0 \off
\vspace{6pt} \\
\off\dff 0 &
1\dff \otimes\dff 1 \off 
\end{pmatrix}
\off +\off
i\pff
\begin{pmatrix}
\off\dff p\dff \sin\dff \eta\dff \otimes\dff 1 &
\sin_{\dff t}\dff \eta\dff \otimes\trf q^{\fff *} \off
\vspace{6pt} \\
\off\dff \sin_{\dff t}\dff \eta\dff \otimes\trf q &
-\qff p\dff \sin\dff \eta\dff \otimes\dff 1 \off 
\end{pmatrix}
\qff.
\]

\vspace{-12pt}
Then\vspace{3pt}
\[
\quad
\left(\qff \overline{p\pff \omult_{\fff 1}\qff q} \qff\right)_{\trf \eta\dff,\dff 1} 
\off =\off
\left(\qff \overline{p\pff \omult_{\fff 1}\qff q} \qff\right)_{\trf \eta}
\quad
\mbox{and}\quad
\left(\qff \overline{p\pff \omult_{\fff 1}\qff q} \qff\right)_{\trf \eta\dff,\dff 0} 
\off =\off
p_{\trf \eta}\pff \omult_{\dff 0}\qff (\trf i\dff q\trf)
\qff.
\]

\vspace{-12pt}\vspace{4.5pt}
Clearly,\oss the endomorphism\vspace{4.5pt}
\[
\quad
\begin{pmatrix}
\off\dff p\dff \sin\dff \eta\dff \otimes\dff 1 &
\sin_{\dff t}\dff \eta\dff \otimes\trf q^{\fff *} \off
\vspace{6pt} \\
\off\dff \sin_{\dff t}\dff \eta\dff \otimes\trf q &
-\qff p\dff \sin\dff \eta\dff \otimes\dff 1 \off 
\end{pmatrix}
\off =\off
(\trf p\dff \sin\dff \eta\trf)
\pff \omult_{\fff 1}\qff
(\trf q\trf \sin_{\dff t}\dff \eta \qff)
\]

\vspace{-12pt}\vspace{4.5pt}
is\dss self-adjoint\sss and\dss is\dss an\sss
isomorphism\sss when\sss either\sss
$p\dff \sin\dff \eta$\sss
or\sss $q\trf \sin_{\dff t}\dff \eta$\sss
is\dss an\dss isomorphism.\oss
Hence\sss it\dss is\dss a
self-adjoint\dss isomorphism\sss when\sss
$\eta\qff \in\qff (\trf 0\fff,\qff \pi\trf)$\sss
and\sss $p$\sss or\sss $q$\sss is\dss an\dss isomorphism.\oss
It\sss follows\sss that\sss endomorphisms\sss\vspace{3pt}
\[
\quad
\left(\qff \overline{p\pff \omult_{\fff 1}\qff q} \qff\right)_{\trf \eta\dff,\dff t}
\]

\vspace{-12pt}\vspace{3pt}
are isomorphisms under\sss the same assumptions.\oss 
Also,\oss these endomorphisms are\sss isomorphisms\sss
when\sss $\eta\off =\off 0$\sss or\sss $\pi$\dss 
(and\sss $t\qff \in\qff [\trf 0\fff,\qff 1\trf]$\sss is\dss arbitrary).\oss
Therefore\sss these endomorphisms\sss
define a path\sss
of\dss endomorphisms\sss
$\overline{E}\qff \ttoo\qff \overline{E}$\sss
which are isomorphisms over\vspace{3pt}
\[
\quad
X\dff \times\dff [\trf 0\fff,\qff \pi\trf]\dff \times\dff Y\fff'\qff \cup\off 
Y\dff \times\dff [\trf 0\fff,\qff \pi\trf]\dff \times\dff X\fff'
\off \cup\off
X\dff \times\dff \{\trf 0\fff,\qff \pi\trf\}\dff \times\dff X\fff'
\qff.
\]

\vspace{-12pt}\vspace{3pt} 
This path connects\sss
$\overline{p}\pff \omult_{\dff 0}\qff (\trf i\dff q\trf)$\sss
with\dss
$\overline{p\pff \omult_{\fff 1}\qff q}$\nnsp.\oss
It\sss follows\sss that\sss\vspace{2.5pt}
\[
\quad
d\trf\left(\pff \overline{p\pff \omult_{\fff 1}\qff q} \pff\right)
\off =\off
d\trf\left(\trf \overline{p}\pff \omult_{\dff 0}\qff (\trf i\dff q\trf)\trf\right)
\qff
\]

\vspace{-36pt}
\[
\quad
\phantom{d\trf\left(\pff \widetilde{p\pff \omult_{\fff 1}\qff q} \pff\right)
\off }
=\off
d\trf\left(\qff \overline{p}\pff \omult_{\dff 0}\qff q\trf\right)
\off =\off
d\trf\left(\qff \overline{p} \pff\right)
\off \widehat{\otimes}\off
d\trf(\trf E\fff'\fff,\qff F\fff'\fff,\qff q\trf)
\qff.
\]

\vspace{-12pt}\vspace{2.5pt}
This proves\qss (\ref{tensor-mult})\qss and\dss hence 
completes\sss the proof\dss of\dss lemma.\oss  \eproof

\myuppar{Multiplication of\dss symbols.}
Let\sss $M$\sss be a compact\sss manifold\sss with\sss the boundary\sss $L\off =\off \partial\dff M$\nnsp.\oss
Let\sss $V$\sss be a closed\sss manifold.\oss
As in\dss Section\qss \ref{symbols-conditions}\qss we will\sss denote by\sss $\pi$\sss
the projections\sss $B\fff M\qff \ttoo\qff M$\sss and\sss $B\fff V\qff \ttoo\qff V$\dnsp.\oss
Let\sss $E$\sss be a vector bundle over\sss $M$\sss and\dss let\sss $E\fff'\fff,\off F\fff'$\dss
be vector bundles over\sss $V$\dnsp.\oss
We will\sss assume\sss that\sss all\dss these bundles are equipped\sss with\dss Hermitian\dss metrics.\oss
Let\sss $\sigma$\sss be a self-adjoint\sss automorphism of\dss 
$\pi^{\fff *}\dff E$\sss over\sss $S\fff M$\nnsp,\oss
and\dss let\sss $q$\sss be an\sss isomorphism\sss
$\pi^{\fff *}\dff E\fff'\qff \ttoo\qff \pi^{\fff *}\dff F\fff'$\sss
over\sss $S\fff V$\dnsp.\oss
We will\sss assume\sss that\sss $\sigma$\sss is\dss a self-adjoint\sss symbol\sss of\dss order $1$\nnsp.\oss
For\sss $u\qff \in\qff S\fff M$ and\sss $r\qff \geq\qff 0$\sss
let\sss $\widehat{\sigma}\trf(\trf r\halfff u\trf)\off =\off r\dff \sigma\trf(\trf u\trf)$\nnsp,\oss
and\dss let\sss us define\sss $\widehat{q}$\sss similarly.\oss
Then\sss $\widehat{\sigma}$\sss and\sss $\widehat{q}$\sss extend $\sigma$ and $q$\sss
to\sss $T\dff M$ and\sss $T\dff V$\sss respectively.\oss

Similarly\sss to\dss Section\qss \ref{symbols-conditions},\oss for $y\qff \in\qff L$\sss
we will\sss denote by\sss $\nu_y$\sss the unit\sss normal\sss to\sss $T_y\dff L$\sss 
in\sss $T_y\trf M$\sss pointing\sss into $M$\nnsp,\oss
and\sss will\sss set\sss 
$\Sigma_{\dff y}\off =\off \sigma\dff(\trf \nu_y\trf)$\nnsp.\oss
The notation\sss $\sigma_y$\sss used\sss for\sss $\sigma\dff(\trf \nu_y\trf)$ 
in\dss Section\qss \ref{symbols-conditions}\qss
would\dss be inconvenient\sss because we will\sss usually omit\sss the subscript $y$\nnsp.\oss 
For\sss $u\qff \in\qff S\fff M_{\trf L}$\sss
we will\sss set\sss $\tau_u\off =\off \sigma\dff(\trf u\trf)$\sss
and\sss
$\rho_{\fff u}
\off =\off 
\Sigma_{\dff y}^{\dff -\dff 1}\dff \tau_u$\nsp.\oss
The extension\sss $\widehat{\sigma}$\sss is\dss 
partially\sss linear in\sss the sense\sss that\vspace{1.5pt}
\[
\quad
\widehat{\sigma}\trf(\trf a\dff \nu_y\qff +\qff b\dff u\trf)
\off =\off
a\trf \widehat{\sigma}\trf(\trf \nu_y\trf)
\qff +\qff
b\trf \widehat{\sigma}\trf(\trf u\trf)
\]

\vspace{-12pt}\vspace{1.5pt}
for every $y\qff \in\qff L$ and\sss $u\qff \in\qff T_y\trf L$\nnsp.\oss
This immediately\sss follows\sss from\sss the facts\sss that\sss $\widehat{\sigma}$\sss
is\dss homogeneous of\dss degree $1$\sss by\sss the definitions
and\sss that\sss $\sigma$\sss is\dss a self-adjoint\sss symbol\sss of\dss order $1$\nnsp.\oss 

The product\sss $\widehat{\sigma}\pff \omult_{\fff 1}\qff \widehat{q}$\sss is\dss
defined over\sss 
$T\dff(\trf M\dff \times\dff V\trf)
\off =\off
T\dff M\dff \times\dff T\dff V$\dnsp.\oss
Let\sss us\sss consider\sss this product\sss over\sss
$S\dff(\trf M\dff \times\dff V\trf)$\nnsp.\oss
Every\sss vector\sss
$w\qff \in\qff S\dff(\trf M\dff \times\dff V\trf)$\sss 
has\sss the form
$w\off =\off (\trf u\dff \cos\dff \theta\fff,\qff v\dff \sin\dff \theta\trf)$\sss
with\sss $u\qff \in\qff S\dff M$\nnsp,\qss
$v\qff \in\qff S\fff V$\dnsp,\oss and\sss $\theta\qff \in\qff [\trf 0\fff,\qff \pi/2 \trf]$\nnsp.\oss
Clearly,\vspace{1.5pt}
\[
\quad
\left(\qff \widehat{\sigma}\pff \omult_{\fff 1}\qff \widehat{q} \pff\right)\trf(\trf w\trf)
\off =\off
\sigma\trf(\trf u\trf)\dff \cos\dff \theta
\off \omult_{\fff 1}\pff
q\trf(\trf v\trf)\dff \sin\dff \theta
\pff.
\]

\vspace{-12pt}\vspace{1.5pt}
This formula expresses\sss the restriction of\dss
$\widehat{\sigma}\pff \omult_{\fff 1}\qff \widehat{q}$\sss
to\sss $S\dff(\trf M\dff \times\dff V\trf)$\sss 
directly\sss in\sss terms of\dss $\sigma\fff,\qff q$\nnsp.\oss
By an abuse of\dss notations we will\sss denote\sss this restriction\sss
by\sss $\sigma\pff \omult_{\fff 1}\qff q$\nnsp.\oss

\myuppar{Further assumptions.}
Let\sss $N$\sss be a self-adjoint\sss elliptic 
and\dss bundle-like boundary condition\sss for $\sigma$\nnsp.\oss
We will\sss assume\sss that\sss $\sigma$\sss is\dss bundle-like and\sss
the pair\sss $\sigma\dff,\off N$\sss is\dss normalized.\oss
Then\sss $\tau_u$\sss and\sss
the fiber\sss $N_{\dff u}$\sss of\dss $N$ over\sss $u\qff \in\qff S\fff L$\sss
depend only on\sss $y\off =\off \pi\trf(\dff u\trf)$\nnsp,\oss
and we may denote\sss them\sss by\sss $\tau_y$ and\sss $N_{\dff y}$ respectively.\oss
Also,\qss
$N_{\dff u}\off =\off \mathcal{L}_{\dff +}\dff(\trf \rho_{\fff u}\dff)$\sss
for every\sss $u\qff \in\qff S\fff L$\nnsp,\oss
and\dss $\Sigma_{\dff y}\trf(\trf N_{\dff y}\trf)\off =\off N_{\dff y}^{\dff \perp}$\nsp.\oss
We will\sss also assume\sss that\sss
$q$\sss is\dss an\sss isometric\sss isomorphism.\oss
Then\sss $q^{\dff *} q\off =\off 1$\sss and\sss $q\fff q^{\dff *}\off =\off 1$\nnsp.\oss

\mypar{Lemma.}{symbol-one}
\emph{$\sigma\pff \omult_{\fff 1}\qff q$\sss
is\dss a self-adjoint\sss symbol of\dss order\dss $1$ on\sss $M\dff \times\dff V$\dnsp.\oss}

\proof
Clearly,\oss vectors\sss 
$\nu_y\off =\off(\trf \nu_y\dff,\qff 0\trf)$\sss are\sss
the unit\sss normals\sss to\sss $T\dff(\trf L\dff \times\dff V\trf)$\sss
in\sss $T\dff(\trf M\dff \times\dff V\trf)$\qss
pointing\sss into\sss $M\dff \times\dff V$\dnsp.\oss
Let\sss $(\trf y\fff,\qff z\trf)\qff \in\qff L\dff \times\dff V$\dnsp,\oss
and\sss $w$\sss be a unit\sss tangent\sss vector\sss to\sss $M\dff \times\dff V$
at\sss the point $(\trf y\fff,\qff z\trf)$\nnsp.\oss
Then 
$w\off =\off (\trf u\dff \cos\dff \theta\fff,\qff v\dff \sin\dff \theta\trf)$\sss
for some $u\qff \in\qff S_y\trf M$\nnsp,\qss $v\qff \in\qff S_{\fff z}\fff V$\nsp\dnsp,\pss
and\sss $\theta\qff \in\qff [\trf 0\fff,\qff \pi/2\trf]$\dnsp.\oss
In\sss turn,\oss
$u\off =\off \nu\dff \cos\dff \eta\qff +\qff u\fff'\dff \sin\dff \eta$\dss
for some\sss $u'\qff \in\qff S_y\trf L$\sss
and\sss $\eta\qff \in\qff [\trf 0\fff,\qff \pi\trf]$\nnsp,\oss
where\sss $\nu\off =\off \nu_y$\nsp.\oss
Hence\vspace{3pt}
\[
\quad
w
\off =\off
\nu\dff \cos\dff \eta\dff \cos\dff \theta
\qff +\qff
u\fff'\dff \sin\dff \eta\dff \cos\dff \theta
\qff +\qff
v\dff \sin\dff \theta
\pff.
\]

\vspace{-12pt}\vspace{3pt}
By similar\sss reasons\sss $w$\sss has\sss the form\vspace{3pt}
\[
\quad
w
\off =\off
\nu\dff \cos\dff \xi
\qff +\qff
u\fff'\dff \cos\dff \zeta\dff \sin\dff \xi
\qff +\qff
v\dff \sin\dff \zeta\dff \sin\dff \xi
\]

\vspace{-33pt}
\[
\quad
\phantom{w}
\off =\off
\nu\dff \cos\dff \xi
\qff +\qff
\bigl(\dff u\fff'\dff \cos\dff \zeta
\qff +\qff
v\dff \sin\dff \zeta\qff\bigr)\dff \sin\dff \xi
\pff
\]

\vspace{-12pt}\vspace{3pt}
for some\sss $\xi\qff \in\qff [\trf 0\fff,\qff \pi\trf]$
and\sss $\zeta\qff \in\qff [\trf 0\fff,\qff \pi/2\trf]$\nnsp.\oss
Since\sss
$\nu\fff,\off u\fff'\fff,\off v$\sss are\sss linearly\sss independent,\oss
we have\vspace{4.5pt}
\[
\quad
\cos\dff \xi
\off =\off 
\cos\dff \eta\dff \cos\dff \theta\qff,
\quad
\cos\dff \zeta\dff \sin\dff \xi
\off =\off
\sin\dff \eta\dff \cos\dff \theta\qff,
\quad
\mbox{and}\quad
\sin\dff \zeta\dff \sin\dff \xi
\off =\off
\sin\dff \theta
\qff.
\]

\vspace{-12pt}\vspace{4.5pt}
It\sss follows\sss that\sss
$\left(\qff \sigma\pff \omult_{\fff 1}\qff q \pff\right)\trf(\trf w\trf)$\sss
is\dss equal\dss to\vspace{4.5pt}
\[
\quad
\widehat{\sigma}\trf\left(\qff
\nu\dff \cos\dff \eta\dff \cos\dff \theta
\qff +\qff
u\fff'\dff \sin\dff \eta\dff \cos\dff \theta
\qff\right)
\off \omult_{\fff 1}\off
\left(\qff
\widehat{q}\qff(\trf v\dff \sin\dff \theta\qff)
\qff\right)
\]

\vspace{-31.5pt}
\[
\quad
=\off
\widehat{\sigma}\trf
\left(\qff
\nu\dff \cos\dff \xi
\qff +\qff
u\fff'\dff \cos\dff \zeta\dff \sin\dff \xi
\qff\right)
\off \omult_{\fff 1}\off
\left(\qff
\widehat{q}\qff(\trf v\dff \sin\dff \zeta\dff \sin\dff \xi\qff)
\qff\right)
\]

\vspace{-31.5pt}
\[
\quad
=\off
\bigl(\qff
\widehat{\sigma}\trf\left(\qff
\nu\dff \cos\dff \xi\qff\right)
\pff \omult_{\fff 1}\pff
0
\qff\bigr)
\off +\off
\left(\qff
\widehat{\sigma}\trf (\trf u\fff'\dff \cos\dff \zeta\dff \sin\dff \xi\qff)
\off \omult_{\fff 1}\off
\widehat{q}\qff(\trf v\dff \sin\dff \zeta\dff \sin\dff \xi\qff)
\qff\right)
\]

\vspace{-31.5pt}
\[
\quad
=\off
\bigl(\qff
\widehat{\sigma}\trf\left(\qff
\nu\qff\right)
\pff \omult_{\fff 1}\pff
0
\qff\bigr)\dff \cos\dff \xi
\off +\off
\left(\qff
\widehat{\sigma}\trf (\trf u\fff'\dff \cos\dff \zeta\qff)
\off \omult_{\fff 1}\off
\widehat{q}\qff(\trf v\dff \sin\dff \zeta\qff)
\qff\right)\dff \sin\dff \xi
\]

\vspace{-31.5pt}
\[
\quad
=\off
\left(\qff
\widehat{\sigma}
\pff \omult_{\fff 1}\pff
\widehat{q}
\pff \right)\trf(\trf \nu\fff,\qff 0\trf)\dff \cos\dff \xi
\off +\off
\left(\qff
\widehat{\sigma}
\off \omult_{\fff 1}\off
\widehat{q}
\pff \right)\trf
\left(\qff u\fff'\dff \cos\dff \zeta\qff +\qff v\dff \sin\dff \zeta\qff\right)\dff \sin\dff \xi
\pff.
\]

\vspace{-12pt}\vspace{4.5pt}
Since\sss $(\trf \nu\fff,\qff 0\trf)$ and\sss
$w\fff'
\off =\off
u\fff'\dff \cos\dff \zeta\qff +\qff v\dff \sin\dff \zeta$\sss
are unit\sss vectors,\oss we see\sss that\vspace{4.5pt}
\[
\quad
\left(\qff \sigma\pff \omult_{\fff 1}\qff q \pff\right)\trf(\trf w\trf)
\off =\off
\left(\qff
\sigma\pff \omult_{\fff 1}\qff q
\qff\right)\trf(\trf \nu\fff,\qff 0\trf)\dff \cos\dff \xi
\off +\off
\left(\qff
\sigma\pff \omult_{\fff 1}\qff q
\qff\right)\trf
(\dff w\fff'\qff)\dff \sin\dff \xi
\pff.
\]

\vspace{-12pt}\vspace{4.5pt}
Since\sss
$w\off =\off (\trf \nu\fff,\qff 0\trf)\dff \cos\dff \xi\qff +\qff w\fff'\dff \sin\dff \xi$\nnsp,\oss
this proves\sss the\sss lemma.\oss  \eproof

\mypar{Lemma.}{boundary-product}
\emph{Under\sss the above assumptions\sss
$\mathcal{N}
\pff =\off
\left(\trf N\dff \boxtimes\dff E\fff' \trf\right)
\qff \oplus\qff
\left(\trf N\dff \boxtimes\dff F\fff' \trf\right)$\dss
is\dss a self-adjoint\sss elliptic 
and\dss bundle-like boundary condition\sss for\sss
$\sigma\pff \omult_{\fff 1}\qff q$\nnsp.\oss}

\proof
The values of\dss the symbol\sss
$\sigma\pff \omult_{\fff 1}\qff q$\sss
on\sss the unit\sss normals\sss
$(\trf \nu_y\dff,\qff 0\trf)$\sss
to\sss $T\dff(\trf L\dff \times\dff V\trf)$\sss
are\vspace{3pt}
\[
\quad
\Sigma_y\pff \omult_{\fff 1}\qff 0
\off\dff =\off
\begin{pmatrix}
\off\dff \Sigma_{\dff y}\dff \otimes\dff 1 &
0 \off
\vspace{6pt} \\
\off\dff 0 &
-\qff \Sigma_{\dff y}\dff \otimes\dff 1 \off 
\end{pmatrix}
\qff.
\]

\vspace{-12pt}\vspace{3pt}
Since $N_{\dff y}$\sss 
is\dss a\sss lagrangian subspace,\pss\vspace{3pt}
\[
\quad
\left(\qff \Sigma_{\dff y}\pff \omult_{\fff 1}\qff 0 \pff\right)\trf
\left(\qff N_{\dff y}\qff \boxtimes\qff E\fff'_{\fff z}\qff\right)
\off\qff =\off\qff
\Sigma_{\dff y}\dff\trf\left(\trf N_{\dff y}\trf\right)\qff \boxtimes\qff E\fff'_{\fff z}
\]

\vspace{-12pt}\vspace{3pt}
is\dss orthogonal\dss to\sss
$N_{\dff y}\qff \boxtimes\qff E\fff'_{\fff z}$\dss for every\sss $z\qff \in\qff V$\dnsp.\oss
Similarly,\vspace{3pt}
\[
\quad
\left(\qff \Sigma_{\dff y}\pff \omult_{\fff 1}\qff 0 \pff\right)\trf
\left(\qff N_{\dff y}\qff \boxtimes\qff F\fff'_{\fff z}\qff\right)
\off\qff =\off\qff
-\qff \sigma_y\dff\trf\left(\trf N_{\dff y}\trf\right)\qff \boxtimes\pff E\fff'_{\fff z}
\off\qff =\off\qff
\Sigma_{\dff y}\dff\trf\left(\trf N_{\dff y}\trf\right)\qff \boxtimes\pff E\fff'_{\fff z}
\]

\vspace{-12pt}\vspace{3pt}
is\dss orthogonal\dss to\sss
$N_{\dff y}\qff \boxtimes\qff F\fff'_{\fff z}$\nsp.\oss
Since\sss the dimension of\trs
$\left(\trf N\dff \boxtimes\dff E\fff' \trf\right)
\qff \oplus\qff
\left(\trf N\dff \boxtimes\dff F\fff' \trf\right)$\dss
is\dss equal\sss to\sss the half\dss of\dss the dimension of\trs
$\left(\trf E\dff \boxtimes\dff E\fff' \trf\right)
\qff \oplus\qff
\left(\trf E\dff \boxtimes\dff F\fff' \trf\right)$\nnsp,\oss
it\sss follows\sss that\sss our subbundle\dss is\dss lagrangian.\oss
Therefore\sss it\dss is\dss a self-adjoint\dss boundary condition.\oss
Since $N$\sss is\dss bundle-like,\oss it\dss is\dss bundle-like.

It\sss remains\sss to prove\sss that\dss the boundary condition\sss
$\mathcal{N}$\dss
is\dss elliptic.\oss
We need\sss to prove\sss that\sss the fibers of\dss the subbundle\sss $\mathcal{N}$\sss
are\sss transverse\sss to\sss the subspace\vspace{3pt}
\begin{equation}
\label{l-minus-product}
\quad
\mathcal{L}_{\trf -}\qff
\left(\qff
\left(\qff \Sigma_{\dff y}\pff \omult_{\fff 1}\qff 0 \pff\right)^{\dff -\dff 1}\qff
\left(\qff \sigma\pff \omult_{\fff 1}\qff q \pff\right)\trf(\trf w\trf)
\qff\right)
\qff
\end{equation}

\vspace{-12pt}\vspace{3pt}
for every\sss unit\sss tangent\sss vector\sss 
$w\qff \in\qff S_{y\fff,\trf z}\trf(\trf L\dff \times\dff V\trf)$\nnsp.\oss
Every such unit\sss vector\sss $w$\sss has\sss the form\vspace{3pt}
\[
\quad
w
\off =\off
(\trf u\dff \cos\dff \theta\fff,\qff v\dff \sin\dff \theta\trf)
\qff,
\]

\vspace{-12pt}\vspace{3pt}
where\sss $u\qff \in\qff S_y\trf L$\nsp,\qss
$v\qff \in\qff S_{\fff z}\dff V$\dnsp,\oss
and\sss $\theta\qff \in\qff [\trf 0\fff,\qff \pi/2 \trf]$\nnsp.\oss
By\sss the definition\vspace{3pt}
\[
\quad
\left(\qff \sigma\pff \omult_{\fff 1}\qff q \pff\right)\trf(\trf w\trf)
\off =\off
\begin{pmatrix}
\off\dff \tau_u\trf \cos\dff \theta\dff \otimes\dff 1 &
1\dff \otimes\trf q\trf(\dff v\trf)^{\fff *}\trf \sin\dff \theta \off
\vspace{6pt} \\
\off\dff 1\dff \otimes\trf q\trf(\dff v\trf)\trf \sin\dff \theta &
-\qff \tau_u\trf \cos\dff \theta\dff \otimes\dff 1 \off 
\end{pmatrix}
\qff.
\]

\vspace{-12pt}\vspace{3pt}
Since\sss $\Sigma_{\dff y}$\sss is\dss self-adjoint\sss and\sss unitary,\pss
$\Sigma_{\dff y}^{\dff -\dff 1}\off =\off \Sigma_{\dff y}$\sss and\dss hence\vspace{3pt}
\[
\quad
\left(\qff \Sigma_{\dff y}\pff \omult_{\fff 1}\qff 0 \pff\right)^{\dff -\dff 1}\qff
\left(\qff \sigma\pff \omult_{\fff 1}\qff q \pff\right)\trf(\trf w\trf)
\off =\off
\begin{pmatrix}
\off\dff \rho\trf \cos\dff \theta\dff \otimes\dff 1 &
\Sigma\dff \otimes\trf q^{\fff *}\fff \sin\dff \theta \off
\vspace{6pt} \\
\off\dff -\qff \Sigma\dff \otimes\trf q\dff \sin\dff \theta &
\rho\trf \cos\dff \theta\dff \otimes\dff 1 \off 
\end{pmatrix}
\qff,
\]

\vspace{-12pt}\vspace{3pt}
where we omit\sss subscripts and write  
$q$\sss for\sss $q\trf(\dff v\trf)$\nnsp,\oss
as we will\sss do also in\sss the following calculations.\oss
Let\sss 
$(\dff a\fff,\qff b\trf)\off \neq\off 0$\sss 
be an\sss eigenvector 
of\dss this matrix\sss with\sss the eigenvalue $\lambda$\nnsp.\oss
Then\vspace{3pt}\vspace{-0.375pt}
\[
\quad
(\trf \rho\trf \cos\dff \theta\dff \otimes\dff 1\trf)\trf
(\dff a\trf)
\qff +\qff
(\trf \Sigma\dff \otimes\trf q^{\fff *} \sin\dff \theta\trf)\trf
(\dff b\trf)
\off =\off
\lambda\dff a
\quad
\mbox{and}
\]

\vspace{-34.5pt}
\[
\quad
(\trf -\qff \Sigma\dff \otimes\trf q\dff \sin\dff \theta\trf)\trf
(\dff a\trf)
\qff +\qff
(\trf \rho\trf \cos\dff \theta\dff \otimes\dff 1\trf)\trf
(\dff b\trf)
\off =\off
\lambda\dff b
\qff,
\]

\vspace{-12pt}\vspace{3pt}\vspace{-1.5pt}\vspace{-0.375pt}
or,\oss equivalently,\vspace{3pt}\vspace{-1.5pt}\vspace{-0.375pt}
\[
\quad
(\trf (\trf \rho\trf \cos\dff \theta\qff -\qff \lambda\trf)\dff \otimes\dff 1\trf)\trf
(\dff a\trf)
\off =\off
(\trf -\qff \Sigma\dff \otimes\trf q^{\fff *} \sin\dff \theta\trf)\trf
(\dff b\trf)
\quad
\mbox{and}
\]

\vspace*{-50pt}
\begin{equation}
\label{eigenvector}
{}
\end{equation}

\vspace*{-50pt}
\[
\quad
(\trf (\trf \rho\trf \cos\dff \theta\qff -\qff \lambda\trf)\dff \otimes\dff 1\trf)\trf
(\dff b\trf)
\off =\off
(\trf \Sigma\dff \otimes\trf q\dff \sin\dff \theta\trf)\trf
(\dff a\trf)
\qff.
\]

\vspace{-12pt}\vspace{3pt}\vspace{-0.375pt}
Since we are in a normalized situation,\pss
$\Sigma^{\dff 2}\off =\off 1$\sss and\sss
$\Sigma\dff \rho\off =\off -\qff \rho\dff \Sigma$\nnsp.\oss
It\sss follows\sss that\vspace{3pt}
\[
\quad
(\trf \rho^{\dff 2}\trf \cos^{\dff 2}\dff \theta\qff -\qff \lambda^{\dff 2}\trf)\dff \otimes\dff 1\trf
(\dff a\trf)
\]

\vspace{-33pt}\vspace{-1.5pt}
\[
\quad
=\off
(\qff (\trf \rho\trf \cos\dff \theta\qff +\qff \lambda\trf)\dff \otimes\dff 1\qff)\trf
(\qff (\trf \rho\trf \cos\dff \theta\qff -\qff \lambda\trf)\dff \otimes\dff 1\qff)\trf
(\dff a\trf)
\]

\vspace{-33pt}\vspace{-1.5pt}
\[
\quad
=\off
(\qff (\trf \rho\trf \cos\dff \theta\qff +\qff \lambda\trf)\dff \otimes\dff 1\qff)\trf
(\qff -\qff \Sigma\dff \otimes\trf q^{\fff *} \sin\dff \theta\qff)\trf
(\dff b\trf)
\]

\vspace{-33pt}\vspace{-1.5pt}
\[
\quad
=\off
(\qff \Sigma\dff \otimes\trf q^{\fff *} \sin\dff \theta\qff)
(\qff (\trf \rho\trf \cos\dff \theta\qff -\qff \lambda\trf)\dff \otimes\dff 1\qff)\trf
(\dff b\trf)
\]

\vspace{-33pt}\vspace{-1.5pt}
\[
\quad
=\off
(\qff \Sigma\dff \otimes\trf q^{\fff *} \sin\dff \theta\qff)
(\qff \Sigma\dff \otimes\trf q\dff \sin\dff \theta\qff)\trf
(\dff a\trf)
\]

\vspace{-33pt}\vspace{-1.5pt}
\[
\quad
=\off
(\qff 1\dff \otimes\trf q^{\fff *}\nsp q\dff \sin^{\dff 2}\dff \theta\qff)\trf
(\dff a\trf)
\]

\vspace{-12pt}\vspace{3pt}
Since we are in a normalized situation,\pss $\rho^{\dff 2}\off =\off -\qff 1$\nnsp.\oss
Since\sss $q^{\dff *}\nsp q\off =\off 1$\nnsp,\oss
it\sss follows\sss that\vspace{3pt}
\[
\quad
(\trf -\qff \cos^{\dff 2}\dff \theta\qff -\qff \lambda^{\dff 2}\trf)\trf
a
\off =\off
(\trf \sin^{\dff 2}\dff \theta\trf)\trf
a
\pff
\]

\vspace{-12pt}\vspace{3pt}
and\dss hence\sss
$\lambda^{\dff 2}
\off =\off
-\qff \cos^{\dff 2}\dff \theta\qff -\qff \sin^{\dff 2}\dff \theta
\off =\off
-\qff 1$\nnsp.\oss
Therefore\sss $\lambda\off =\off i$\sss or\sss $-\qff i$\nnsp,\oss
and\sss the\sss last\sss equation\sss imposes no conditions on $a$\nnsp.\oss
Moreover,\oss if\dss $\lambda\off =\off i$\sss or\sss $-\qff i$\nnsp,\oss
then,\oss under our assumptions,\oss each of\dss the equations\qss (\ref{eigenvector})\qss
implies\sss the other.\oss
It\sss follows\sss that\sss the eigenvalues of\vspace{3pt}
\[
\quad
\left(\qff \Sigma_{\dff y}\pff \omult_{\fff 1}\qff 0 \pff\right)^{\dff -\dff 1}\qff
\left(\qff \sigma\pff \omult_{\fff 1}\qff q \pff\right)\trf(\trf w\trf)
\]

\vspace{-12pt}\vspace{3pt}
are\sss $\lambda\off =\off i$\sss or\sss $-\qff i$\nnsp,\oss
and\sss that\sss the eigenspace corresponding\sss to $\lambda$ 
can\sss be described\dss by either of\dss the equations\qss (\ref{eigenvector})\qss 
and\dss hence\dss is\dss equal\dss to\sss the graph of\dss the map\vspace{3pt}
\[
\quad
\varphi_{\dff \lambda}
\off =\off
(\qff -\qff \Sigma\dff \otimes\trf q^{\fff *} \sin\dff \theta \qff)^{\dff -\dff 1}
\trf \circ\qff\halfff
(\qff (\trf \rho\trf \cos\dff \theta\qff -\qff \lambda\trf)\dff \otimes\dff 1 \qff)
\]

\vspace{-34.5pt}
\[
\quad
\phantom{\varphi
\off }
=\off
-\qff \Sigma\dff \circ\dff (\qff (\trf \rho\trf \cos\dff \theta\qff -\qff \lambda\trf)
\dff \otimes\dff 
q\dff \sin^{\dff -\dff 1}\dff \theta \qff)
\qff \colon\qff
E\dff \boxtimes\dff E\fff'
\qff \ttoo\qff
E\dff \boxtimes\dff F\fff'
\qff.
\]

\vspace{-12pt}\vspace{3pt}
Here we continue\sss to omit\sss subscripts corresponding\sss to\sss
the points $y\fff,\qff z\fff,\qff \ldots$ etc.\oss
Since we are in a normalized\sss situation,\qss $\rho$\sss leaves 
$N\off =\off N_{\dff y}$ invariant,\oss
and\dss hence\sss $\rho\trf \cos\dff \theta\qff -\qff \lambda$\sss leaves $N$\sss invariant\sss also.\oss
Since $\Sigma\trf(\trf N\trf)\off =\off N^{\dff \perp}$\dnsp,\oss
it\sss follows\sss that\sss $\varphi_{\dff \lambda}$\sss takes\sss
$N\dff \boxtimes\dff E\fff'$\sss to\sss
$N^{\dff \perp}\dff \boxtimes\dff F\fff'$\nnsp.\oss
In\sss turn,\oss this implies\sss that\sss the graph of\dss $\varphi_{\dff \lambda}$\sss
intersects\sss the subspace\dss
$\left(\trf N\dff \boxtimes\dff E\fff' \trf\right)
\qff \oplus\qff
\left(\trf N\dff \boxtimes\dff F\fff' \trf\right)$\dss
only\sss by $0$\nnsp.\oss
It\sss follows\sss that\sss this subspace\dss is\dss transverse\sss to\sss
both eigenspaces and,\oss in\sss particular,\oss to\sss
the eigenspace\qss (\ref{l-minus-product}).\oss
Hence\dss 
$\mathcal{N}$\dss
is\dss an elliptic boundary condition.\oss  \eproof

\myuppar{Deforming\dss $\sigma\pff \omult_{\fff 1}\qff q$\dss to a bundle-like symbol.}
Let\sss $w$\sss be a unit\sss tangent\sss vector\sss to\sss $M\dff \times\dff V$
at\sss a point\sss $(\trf y\fff,\qff z\trf)\qff \in\qff L\dff \times\dff V$\dnsp.\oss
Then\sss the vector\sss $w$\sss has\sss the form\sss\vspace{3pt}
\[
\quad
w
\off =\off
\nu\dff \cos\dff \xi
\qff +\qff
\bigl(\dff u\fff'\dff \cos\dff \zeta
\qff +\qff
v\dff \sin\dff \zeta\qff\bigr)\dff \sin\dff \xi
\pff
\]

\vspace{-12pt}\vspace{3pt}
for some\sss $u\fff'\qff \in\qff S_y\trf L$\nnsp,\qss
$v\qff \in\qff S_{\fff z}\dff V$\dnsp,\qss
$\xi\qff \in\qff [\trf 0\fff,\qff \pi\trf]$
and\sss $\zeta\qff \in\qff [\trf 0\fff,\qff \pi/2\trf]$\nnsp.\oss
For\sss $t\qff \in\qff [\trf 0\fff,\qff 1\trf]$\dss let\vspace{3pt}
\[
\quad
w\trf(\dff t\trf)
\off =\off
\nu\dff \cos\dff \xi
\off +\off
\bigl(\dff u\fff'\dff \cos\trf (\dff t\trf \zeta\qff)
\qff +\qff
v\dff \sin\trf (\dff t\trf \zeta\qff)\qff\bigr)\trf \sin\dff \xi
\quad
\mbox{and}\quad
\pff
\]

\vspace{-33pt}
\[
\quad
\left(\qff \sigma\pff \omult_{\fff 1}\qff q \pff\right)_{\dff t}\trf(\dff w\qff)
\off =\off
\left(\qff \sigma\pff \omult_{\fff 1}\qff q \pff\right)\trf\bigl(\dff w\trf(\dff t\trf)\qff\bigr)
\]

\vspace{-33pt}
\[
\quad
\phantom{\left(\qff \sigma\pff \omult_{\fff 1}\qff q \pff\right)_{\dff t}\trf(\dff w\qff)
\off }
=\off
\sigma\trf
\left(\trf
\nu\dff \cos\dff \xi
\qff +\qff
u\fff'\dff \cos\trf (\dff t\trf \zeta\qff)\trf \sin\dff \xi
\qff\right)\dff 
\off \omult_{\fff 1}\off
q\trf(\trf v\trf)\dff \sin\trf (\dff t\trf \zeta\qff)\trf \sin\dff \xi
\off
\]

\vspace{-33pt}
\[
\quad
\phantom{\left(\qff \sigma\pff \omult_{\fff 1}\qff q \pff\right)_{\dff t}\trf(\dff w\qff)
\off }
=\off
\bigl(\qff
\Sigma\dff \cos\dff \xi
\qff +\qff
\tau\dff \cos\trf (\dff t\trf \zeta\qff)\trf \sin\dff \xi
\qff\bigr)
\off \omult_{\fff 1}\off
q\trf(\trf v\trf)\dff \sin\trf (\dff t\trf \zeta\qff)\trf \sin\dff \xi
\off
\]

\vspace{-33pt}
\[
\quad
\phantom{\left(\qff \sigma\pff \omult_{\fff 1}\qff q \pff\right)_{\dff t}\trf(\dff w\qff)
\off }
=\off
\bigl(\qff
\bigl(\qff
\Sigma\dff \cos\dff \xi\qff\bigr)
\off \omult_{\fff 1}\off
0\qff\bigr)
\off +\off
\bigl(\trf
\tau\dff \cos\trf (\dff t\trf \zeta\qff)\off \omult_{\fff 1}\off
q\trf(\trf v\trf)\dff \sin\trf (\dff t\trf \zeta\qff)
\qff\bigr)\trf \sin\dff \xi
\off,
\]

\vspace{-12pt}\vspace{3pt}
where\sss \sss $\tau\off =\off \tau_y\off =\off \tau_{u\fff'}$\sss
and\sss
$\Sigma\off =\off \Sigma_{\dff y}$\nsp.\oss
While\sss $u\fff'$\sss is\dss not\sss uniquely determined\dss by $w$\sss
when\sss $\zeta\off =\off \pi/2$ 
and\sss $v$\sss is\dss not\sss uniquely determined\dss by $w$\sss 
when\sss $\zeta\off =\off 0$\nnsp,\oss\vspace{1.5pt}
\[
\quad
\left(\qff \sigma\pff \omult_{\fff 1}\qff q \pff\right)_{\dff t}\trf(\dff w\qff)
\]

\vspace{-12pt}\vspace{1.5pt}
is\dss well\sss defined\dss because\sss $\tau$\sss is\dss independent\sss of\dss $u\fff'$
and\sss $\sin\dff 0\off =\off 0$\nnsp.\oss
Therefore\vspace{2.5pt}
\begin{equation}
\label{product-deformation}
\quad
\left(\qff \sigma\pff \omult_{\fff 1}\qff q \pff\right)_{\dff t}\dff,\off\off
t\qff \in\qff [\trf 0\fff,\qff 1\trf]
\qff
\end{equation}

\vspace{-12pt}\vspace{2.5pt}
is\dss a deformation of\dss the symbol\sss 
$\sigma\pff \omult_{\fff 1}\qff q
\off\dff =\off
\left(\qff \sigma\pff \omult_{\fff 1}\qff q \pff\right)_{\dff 1}$\sss
over\sss $L\dff \times\dff V$\dnsp,\oss
and\sss this deformation can be extended\sss to a deformation of\dss the symbol\sss
$\sigma\pff \omult_{\fff 1}\qff q$\dss over\sss the whole\sss product\sss
$M\dff \times\dff V$\dnsp.\oss
Clearly,\pss
$\left(\qff \sigma\pff \omult_{\fff 1}\qff q \pff\right)_{\dff t}$\dss
is\dss a self-adjoint\sss symbol\sss of\dss order $1$ 
and\trs Lemma\qss \ref{boundary-product}\qss implies\sss that\vspace{1.5pt}
\[
\quad
\mathcal{N}
\off =\off
\left(\trf N\dff \boxtimes\dff E\fff' \trf\right)
\qff \oplus\qff
\left(\trf N\dff \boxtimes\dff F\fff' \trf\right)
\]

\vspace{-12pt}\vspace{1.5pt}
is\dss an elliptic self-adjoint\dss boundary condition\sss for\sss
this symbol\sss for every\sss $t$\nnsp.\oss
Clearly,\oss the symbol\sss
$\left(\qff \sigma\pff \omult_{\fff 1}\qff q \pff\right)_{\trf 0}$\sss
is\dss equal\dss to\sss
$\sigma\pff \omult_{\fff 1}\qff 0$\sss
over\sss $L\dff \times\dff V$\sss
and\sss hence\dss is\dss bundle-like.\oss
Also,\oss this symbol\sss
together\sss with\sss the boundary condition\sss
$\left(\trf N\dff \boxtimes\dff E\fff' \trf\right)
\qff \oplus\qff
\left(\trf N\dff \boxtimes\dff F\fff' \trf\right)$\sss
form a normalized\sss pair.\oss

The deformation\qss (\ref{product-deformation})\qss
can\sss be considered as a special\sss case of\dss the canonical\qss
(up\sss to a homotopy)\qss deformation\sss to a normalized\dss pair.\oss
Indeed,\oss since\sss $\sigma\fff,\off N$\sss is\dss normalized,\pss
the operators\sss $\Sigma_{\dff y}\pff \omult_{\fff 1}\qff 0$\sss are unitary\sss
and\dss hence\sss the first\sss deformation\sss from\dss Section\qss \ref{boundary-algebra}\qss
is\dss not\sss needed.\oss
As we saw\sss in\sss the proof\dss of\qss Lemma\qss \ref{boundary-product},\oss
the operators\sss $\rho$ corresponding\sss to\sss
$\sigma\pff \omult_{\fff 1}\qff q$\sss
have only\sss $i$ and $-\qff i$\sss as eigenvalues.\oss
Therefore\sss the second\sss homotopy\dss is\dss also not\sss needed.\oss
The\sss third\sss deformation\sss moves\sss the eigenspaces corresponding\sss to\sss $-\qff i$\sss
to\sss the orthogonal\sss complements\sss of\dss the fibers of\dss the boundary condition.\oss
For\sss $\left(\qff \sigma\pff \omult_{\fff 1}\qff q \pff\right)_{\dff 0}$\sss
these eigenspaces are\sss the fibers of\trs
$\left(\trf N^{\trf \perp}\dff \boxtimes\dff E\fff' \trf\right)
\qff \oplus\qff
\left(\trf N^{\trf \perp}\dff \boxtimes\dff F\fff' \trf\right)$\nnsp,\oss
and\sss hence our deformation\sss can\sss play\sss the role of\dss
the\sss third deformation\qss (recall\dss that\sss it\dss is\dss unique up\sss to homotopy).\oss
Since our deformation also\sss brings\sss the eigenspaces corresponding\sss to\sss $i$\sss
to\sss the fibers of\trs 
$\mathcal{N}
\off =\off
\left(\trf N\dff \boxtimes\dff E\fff' \trf\right)
\qff \oplus\qff
\left(\trf N\dff \boxtimes\dff F\fff' \trf\right)$\nnsp,\oss
the fourth deformation\dss is\dss also not\sss needed.\oss

\myuppar{Extending\dss symbols\sss to\sss $B\trf(\trf M\dff \times\dff V\trf)_{\trf L\dff \times\dff V}$\nsp.}
Let\sss
$w\qff \in\qff B_{\dff y\fff,\dff z}\trf(\trf M\dff \times\dff V\trf)$\sss
for some\sss $(\trf y\fff,\qff z\trf)\qff \in\qff L\dff \times\dff V$\dnsp,\oss
and\sss let\sss $\nu\off =\off \nu_y$\nsp.\oss
Then\sss the vector\sss $w$\sss has\sss the form\sss\vspace{4.5pt}
\begin{equation}
\label{w-form}
\quad
w
\off =\off
\nu\dff \cos\dff \xi
\qff +\qff
\bigl(\dff r\halfff u\fff'\dff \cos\dff \zeta
\qff +\qff
v\dff \sin\dff \zeta\qff\bigr)\dff \sin\dff \xi
\pff
\end{equation}

\vspace{-12pt}\vspace{4.5pt}
for some\sss $r\qff \in\qff [\trf 0\fff,\qff 1\trf]$\sss 
and\sss
$u\fff'\fff,\pff
v\fff,\off
\xi$
and\sss $\zeta$\sss as above.\oss
For\sss $t\qff \in\qff [\trf 0\fff,\qff 1\trf]$\dss let\vspace{4.5pt}
\[
\quad
\left(\qff \mathbf{\bm{\sigma}\pff \omult_{\fff 1}\qff q} \pff\right)_{\dff t}\trf(\dff w\qff)
\off =\off
\sigma\trf\left(\trf
\nu\dff \cos\dff \xi
\qff +\qff
u\fff'\dff \cos\trf (\dff t\trf \zeta\qff)\dff 
\sin\dff \xi
\qff\right)
\off \omult_{\fff 1}\off
q\trf(\trf v\trf)\dff \sin\trf (\dff t\trf \zeta\qff)\dff \sin\dff \xi\off
\]

\vspace{-12pt}\vspace{4.5pt}
Like above,\oss this endomorphism\dss is\dss well\sss defined even when
$u\fff'$ or $v$\sss is\dss not\sss uniquely determined\dss by $w$\nnsp.\oss
This definition\sss leads\sss to an extension\sss
$\left(\qff \mathbf{\bm{\sigma}\pff \omult_{\fff 1}\qff q} \pff\right)_{\dff t}$\dss
of\dss
$\left(\qff \sigma\pff \omult_{\fff 1}\qff q \pff\right)_{\dff t}$\dss
to\sss\vspace{3pt}\vspace{-0.75pt}
\[
\quad
D\trf(\trf M\dff \times\dff V\trf)
\off =\off
S\trf(\trf M\dff \times\dff V\trf)
\qff \cup\qff
B\trf(\trf M\dff \times\dff V\trf)_{\trf L\dff \times\dff V}
\]

\vspace{-12pt}\vspace{3pt}\vspace{-0.75pt}
for every $t$\nnsp.\oss
Clearly,\pss
$\left(\qff \mathbf{\bm{\sigma}\pff \omult_{\fff 1}\qff q} \pff\right)_{\trf 0}$\dss
is\dss the canonical\sss extension of\dss 
$\left(\qff \sigma\pff \omult_{\fff 1}\qff q \pff\right)_{\trf 0}$\nsp.\oss
By an abuse of\dss notations we will\sss denote\sss
$\left(\qff \mathbf{\bm{\sigma}\pff \omult_{\fff 1}\qff q} \pff\right)_{\dff 1}$\dss
simply\sss by\sss
$\mathbf{\bm{\sigma}\pff \omult_{\fff 1}\qff q}$\nnsp.\oss

\myuppar{Extended\sss symbols and\sss products.}
Despite our notations,\oss the extension\sss 
$\mathbf{\bm{\sigma}\pff \omult_{\fff 1}\qff q}$\sss
is\dss not\sss really a\sss $\omult_{\fff 1}$\sss product.\oss
Nevertheless,\oss it\sss closely\sss related\dss to products\sss
$p\pff \omult_{\fff 1}\qff \widehat{q}$\dss for extensions $p$
of\sss $\sigma$ and,\oss especially,\oss of\dss $\bm{\sigma}$ 
to self-adjoint\sss endomorphisms of\dss $\pi^{\dff *}\fff E\qff |\trf B\dff M$\nnsp.\oss 
To begin\sss with,\oss let\sss $p$ and $p\fff'$\sss be\sss two extensions of\sss $\sigma$\sss
from\sss $S\dff M$\sss to\sss $B\dff M$\nnsp,\oss
and\dss let\vspace{1.5pt}
\[
\quad
p_{\fff s}
\off =\off
s\dff p\qff +\qff (\trf 1\qff -\qff s\trf)\dff p\fff'\fff,\off
s\qff \in\qff [\trf 0\fff,\qff 1\trf]
\]

\vspace{-12pt}\vspace{1.5pt}
be\sss the\sss linear\sss homotopy connecting\sss $p$ and $p\fff'$\dnsp.\oss
Since $p$ and $p\fff'$\sss extend $\sigma$\nnsp,\oss
this homotopy\dss is\dss fixed on $S\dff M$\nnsp.\oss
The products\sss $p_{\fff s}\pff \omult_{\fff 1}\qff \widehat{q}$\dss
form a homotopy connecting\sss $p\pff \omult_{\fff 1}\qff \widehat{q}$\dss
with\sss $p\fff'\pff \omult_{\fff 1}\qff \widehat{q}$\nnsp.\oss
We claim\sss that\sss $p_{\fff s}\pff \omult_{\fff 1}\qff \widehat{q}$\dss
is\dss an\sss isomorphism over\sss $S\trf(\trf M\dff \times\dff V\trf)$ for every $s$\nnsp.\oss
Indeed,\oss let\sss 
$w\off =\off (\trf u\dff \cos\dff \theta\fff,\qff v\dff \sin\dff \theta\trf)$\sss
for some $u\qff \in\qff S\dff M$\nnsp,\qss $v\qff \in\qff S\dff V$\nsp\dnsp,\oss
and\sss $\theta\qff \in\qff [\trf 0\fff,\qff \pi/2\trf]$\dnsp.\oss
Then\vspace{3pt}
\[
\quad
\left(\qff p_{\fff s}\pff \omult_{\fff 1}\pff \widehat{q} \pff\right)\trf (\dff w\trf)
\off =\off
\left(\qff s\dff p\qff +\qff (\trf 1\qff -\qff s\trf)\dff p\fff' \qff\right)\trf (\dff u\dff \cos\dff \theta\trf)
\off \omult_{\fff 1}\off
\widehat{q}\trf(\trf v\dff \sin\trf \theta\trf)
\pff.
\]

\vspace{-12pt}\vspace{3pt}
This endomorphism\dss is\dss invertible\sss if\trs
$\widehat{q}\trf(\trf v\dff \sin\trf \theta\trf)$\sss is\dss invertible,\oss
i.e.\qss if\dss $\sin\dff \theta\off \neq\off 0$\nnsp.\oss
If\dss $\sin\dff \theta\off =\off 0$\nnsp,\oss
then\sss $\theta\off =\off 0$\sss and\dss hence\sss $\cos\dff \theta\off =\off 1$\nnsp.\oss
Since\sss $p\dff(\dff u\trf)\off =\off p\fff'\dff(\dff u\trf)\off =\off \sigma\trf(\dff u\trf)$\nnsp,\oss
in\sss this case our endomorphism\dss is\dss equal\dss to\sss
$\sigma\trf(\dff u\trf)\pff \omult_{\fff 1}\qff 0$ and\dss hence\dss is\dss invertible.\oss
This proves our claim.\oss

\mypar{Lemma.}{extensions}
\emph{Let $p$ be an extension of\dss $\bm{\sigma}$ to $B\dff M$\nnsp.\oss
Then\sss the restriction of\dss $p\pff \omult_{\fff 1}\qff \widehat{q}$\dss to}\vspace{3pt}
\[
\quad
D\trf(\trf M\dff \times\dff V\trf)
\off =\off
S\trf(\trf M\dff \times\dff V\trf)
\qff \cup\qff
B\trf(\trf M\dff \times\dff V\trf)_{\trf L\dff \times\dff V}
\]

\vspace{-12pt}\vspace{3pt}
\emph{is\dss homotopic\sss to\sss
$\mathbf{\bm{\sigma}\pff \omult_{\fff 1}\qff q}$\dss
by a homotopy consisting of\qss bundle\sss isomorphisms.\oss
This homotopy can be extended\dss to\dss $B\trf(\trf M\dff \times\dff V\trf)$
as a homotopy of\trs self-adjoint\dss bundle endomorphisms.\oss}

\proof
Over\sss $S\trf(\trf M\dff \times\dff V\trf)$\sss the endomorphism\sss
$\mathbf{\bm{\sigma}\pff \omult_{\fff 1}\qff q}$\sss
is\dss equal\dss to\sss $\sigma\pff \omult_{\fff 1}\qff q$\nnsp,\oss
i.e.\qss to\sss the restriction of\dss
$\widehat{\sigma}\pff \omult_{\fff 1}\qff \widehat{q}$\sss
to\sss $S\trf(\trf M\dff \times\dff V\trf)$\nnsp.\oss
We will\sss use over\sss $S\trf(\trf M\dff \times\dff V\trf)$\sss
the\sss linear\sss homotopy constructed\dss just\dss before\sss the\sss lemma\sss with\sss 
$p\fff'\off =\off \widehat{\sigma}$\nnsp.\oss
Let\sss us consider\sss what\sss happens 
over\sss $L\dff \times\dff V$\dnsp.\oss
Let\sss
$w\qff \in\qff B\dff(\trf M\dff \times\dff V\trf)_{\trf L\dff \times\dff V}$\nsp,\oss
and\dss let\sss us represent\sss $w$\sss in\sss the form\qss (\ref{w-form}).\oss
By\sss the definitions\vspace{4.5pt}
\[
\quad
\left(\qff p\pff \omult_{\fff 1}\pff \widehat{q} \pff\right)\trf (\dff w\trf)
\off =\off
\bm{\sigma}\trf\left(\trf
\nu\dff \cos\dff \xi
\qff +\qff
r\halfff u\fff'\dff \cos\trf \zeta\trf 
\sin\dff \xi
\qff\right)
\off \omult_{\fff 1}\off
q\trf(\trf v\trf)\dff \sin\trf \zeta\trf \sin\dff \xi\off
\]

\vspace{-31.5pt}
\[
\quad
\phantom{\left(\qff p\pff \omult_{\fff 1}\pff \widehat{q} \pff\right)\trf (\dff w\trf)
\off }
=\off
\sigma\trf\left(\trf
\nu\dff \cos\dff \xi
\qff +\qff
u\fff'\dff  
\sin\dff \xi
\qff\right)
\off \omult_{\fff 1}\off
q\trf(\trf v\trf)\dff \sin\trf \zeta\trf \sin\dff \xi
\off,
\]

\vspace{-31.5pt}
\[
\quad
\left(\qff \mathbf{\bm{\sigma}\pff \omult_{\fff 1}\qff q} \pff\right)\trf(\dff w\qff)
\off =\off
\sigma\trf\left(\trf
\nu\dff \cos\dff \xi
\qff +\qff
u\fff'\dff \cos\trf \zeta\trf 
\sin\dff \xi
\qff\right)
\off \omult_{\fff 1}\off
q\trf(\trf v\trf)\dff \sin\trf \zeta\trf \sin\dff \xi\off
\]

\vspace{-12pt}\vspace{4.5pt}
Let\sss us use a\sss linear\sss homotopy again.\oss
Let\sss
$\cos_{\dff t}\dff \zeta\off =\off (\trf 1\qff -\qff t\trf)\qff +\qff t\dff \cos\trf \zeta$\sss
and\dss let\sss us connect\vspace{3.3pt}
\[
\quad
\left(\qff p\pff \omult_{\fff 1}\pff \widehat{q} \pff\right)\trf (\dff w\trf)
\quad
\mbox{with}\quad
\left(\qff \mathbf{\bm{\sigma}\pff \omult_{\fff 1}\qff q} \pff\right)\trf(\dff w\qff)
\]

\vspace{-12pt}\vspace{3.3pt}
by\sss the path of\dss bundle endomorphisms\vspace{3.3pt}
\[
\quad
\sigma\trf\left(\trf
\nu\dff \cos\dff \xi
\qff +\qff
u\fff'\dff \cos_{\dff t}\dff \zeta\trf 
\sin\dff \xi
\qff\right)
\off \omult_{\fff 1}\off
q\trf(\trf v\trf)\dff \sin\trf \zeta\trf \sin\dff \xi
\]

\vspace{-12pt}\vspace{3.3pt}
with\sss $t\qff \in\qff [\trf 0\fff,\qff 1\trf]$\nnsp.\oss
Like above,\oss such an endomorphism\dss is\dss invertible when\sss 
$\zeta\off \neq\off 0$\sss and\sss $\xi\off \neq\off 0$\nnsp.\oss
If\dss $\xi\off =\off 0$\nnsp,\oss then\sss it\dss is\dss equal\dss to\sss 
$\sigma\trf(\trf \nu\trf)$\sss and\dss hence\dss is\dss invertible.\oss
If\dss $\zeta\off =\off 0$\nnsp,\oss then\sss $\cos_{\dff t}\dff \zeta\off =\off 1$\dss
({\fff}for every $t$\nnsp)\qss and\dss this endomorphism\dss is\dss equal\dss to\sss
$\sigma\trf(\trf
\nu\dff \cos\dff \xi
\qff +\qff
u\fff'\trf 
\sin\dff \xi
\qff)
\off \omult_{\fff 1}\off
0$\sss
and\dss hence\dss is\dss invertible.\oss
If\dss
$w$\sss is\dss a\sss vector in\sss $S\trf(\trf M\dff \times\dff V\trf)$\nnsp,\oss
then\sss this homotopy 
agrees\sss with\sss the\sss linear\sss homotopy
connecting\sss $p\pff \omult_{\fff 1}\qff \widehat{q}$\dss
with\sss $\widehat{\sigma}\pff \omult_{\fff 1}\qff \widehat{q}$\nnsp.\oss
Therefore we can use\sss the\sss latter\sss homotopy over\sss $S\trf(\trf M\dff \times\dff V\trf)$\sss
and\sss the paths we\sss just\sss constructed over\sss
$B\dff(\trf M\dff \times\dff V\trf)_{\trf L\dff \times\dff V}$\nsp.\oss
This proves\sss the\sss first\sss statement\sss of\dss the\sss lemma.\oss
The second one follows from\sss
the homotopy extension\sss theorem.\oss  \eproof

\myuppar{Passing\dss from\sss $B\dff (\trf M\dff \times\dff V\trf)$\sss
to\sss $B\dff M\qff \times\qff B\dff V$\dnsp.}
Let\sss us define a map\vspace{3pt}
\[
\quad
f\dff \colon\dff
B\dff (\trf M\dff \times\dff V\trf)
\qff \ttoo\qff
B\dff M\qff \times\qff B\dff V
\pff
\]

\vspace{-12pt}\vspace{3pt}
as follows.\oss
Every vector\sss
$w\qff \in\qff B\dff (\trf M\dff \times\dff V\trf)$\sss
has\sss the form\sss
$w\off =\off r\trf(\trf u\fff \cos\dff \theta\fff,\qff v\dff \sin\dff \theta\trf)$\sss
for some\sss $r\qff \in\qff [\trf 0\fff,\qff 1 \trf]$\nnsp,\qss $u\qff \in\qff S\dff M$\nnsp,\qss
$v\qff \in\qff S\dff V$\dnsp,\oss and\sss $\theta\qff \in\qff [\trf 0\fff,\qff \pi/2 \trf]$\nnsp.\oss
Let\vspace{3pt}
\[
\quad
f\dff(\dff w\trf)
\off =\off
r\trf(\dff u\fff,\pff v\dff \sin\dff \theta\dff/\hnsp\cos\dff \theta\trf)
\quad
\mbox{for}\quad
\theta\qff \leq\qff \pi/4
\quad
\mbox{and}\quad
\]

\vspace{-33pt}\vspace{-1.5pt}
\[
\quad
f\dff(\dff w\trf)
\off =\off
r\trf(\dff u\cos\dff \theta\dff/\hnsp \sin\dff \theta\fff,\pff v\trf)
\quad
\mbox{for}\quad
\theta\qff \geq\qff \pi/4
\pff.
\]

\vspace{-12pt}\vspace{3pt}
Then $f$\sss is\dss a\sss homeomorphism\sss and\oss\vspace{3pt}
\[
\quad
f\trf\bigl(\trf S\dff (\trf M\dff \times\dff V\trf) \trf\bigr)
\off =\off
S\dff M\qff \times\qff B\dff V
\qff \cup\pff
B\dff M\dff \times\dff S\dff V
\pff,
\]

\vspace{-33pt}\vspace{-1.5pt}
\[
\quad
f\trf\bigl(\trf B\trf(\trf M\dff \times\dff V\trf)_{\trf L\dff \times\dff V} \trf\bigr)
\off =\off
B\dff M_{\trf L}\qff \times\qff B\dff V
\pff.
\]

\vspace{-12pt}\vspace{3pt}
There\dss is\dss a natural\sss radial\dss homotopy\sss between\sss $f$\sss and\sss
the inclusion\sss\vspace{3pt}
\[
\quad
i\dff \colon\dff
B\dff (\trf M\dff \times\dff V\trf)
\qff \ttoo\qff 
B\dff M\dff \times\dff B\dff V
\pff.
\]

\vspace{-12pt}\vspace{3pt}
It\sss follows\sss that\sss
$\left(\qff p\pff \omult_{\fff 1}\pff \widehat{q} \pff\right)\trf \circ\dff f$\sss
is\dss homotopic\sss to\sss
$\left(\qff p\pff \omult_{\fff 1}\pff \widehat{q} \pff\right)\trf \circ\dff i$\nnsp.\oss
Moreover,\oss during\sss this homotopy\sss the bundle endomorphisms remain\sss
invertible over\sss 
$D\trf(\trf M\dff \times\dff V\trf)$\nnsp.\oss
Therefore\sss\vspace{1.5pt}
\[
\quad
\left(\qff p\pff \omult_{\fff 1}\pff \widehat{q} \pff\right)\trf \circ\dff f^{\dff -\dff 1}
\]

\vspace{-12pt}\vspace{1.5pt}
is\dss homotopic\sss to
$p\pff \omult_{\fff 1}\pff \widehat{q}$\trs
by a homotopy\sss keeping endomorphisms invertible over\sss
$D\trf(\trf M\dff \times\dff V\trf)$\dnsp.\oss
This will\sss allow us\sss to\sss replace\sss
$B\dff (\trf M\dff \times\dff V\trf)$\sss
by\sss $B\dff M\qff \times\qff B\dff V$
and\sss the restriction of\dss the product\dss
$p\pff \omult_{\fff 1}\pff \widehat{q}$\qss to\sss 
$B\dff (\trf M\dff \times\dff V\trf)$\sss
by\sss the product\sss
$p\pff \omult_{\fff 1}\pff \widehat{q}$\qss itself.\oss

\myuppar{The\sss invariants\sss
$\varepsilon^{\dff +}\fff (\trf \bullet\fff,\qff \bullet\qff)$\nnsp.}
Since\sss $\sigma$\sss is\dss a self-adjoint\sss symbol\sss of\dss order $1$
and\sss $N$\sss is\dss a bundle-like boundary condition for $\sigma$\nnsp,\oss
the invariant\vspace{1.5pt}
\[
\quad
\varepsilon^{\dff +}\fff (\trf \sigma\fff,\qff N\qff)
\off \in\off
K^{\dff 1}\dff (\trf B\dff M\fff,\off  S\dff M\dff \cup\dff B\dff M_{\trf L}\trf)
\]

\vspace{-12pt}\vspace{1.5pt}
is\dss defined.\oss
At\sss the same\sss time\sss
the symbol\sss
$\sigma\pff \omult_{\fff 1}\qff q$\sss
is\dss a self-adjoint\sss symbol of\dss order\dss $1$ on $M\dff \times\dff V$
by\dss Lemma\qss \ref{symbol-one},\oss
and\sss $\mathcal{N}$\sss
is\dss a self-adjoint\sss elliptic 
and\dss bundle-like boundary condition\sss for\sss
$\sigma\pff \omult_{\fff 1}\qff q$\dss
by\trs Lemma\qss \ref{boundary-product}.\oss
Therefore\sss there invariant\vspace{4.5pt}
\[
\quad
\varepsilon^{\dff +}\dff 
\left(\qff 
\sigma\pff \omult_{\fff 1}\qff q\fff,\off \mathcal{N}
\qff\right)
\off \in\off
K^{\dff 1}\trf 
\bigl(\qff 
B\dff (\trf M\dff \times\dff V\trf)\fff,\off  
S\trf(\trf M\dff \times\dff V\trf)
\qff \cup\qff
B\trf(\trf M\dff \times\dff V\trf)_{\trf L\dff \times\dff V}
\qff\bigr)
\]

\vspace{-12pt}\vspace{4.5pt}
is\dss defined.\oss
If\dss we identify\sss
$B\dff (\trf M\dff \times\dff V\trf)$\sss
with\sss $B\dff M\qff \times\qff B\dff V$\dss by\sss $f$\nnsp,\oss
then\vspace{4.5pt}
\[
\quad
\varepsilon^{\dff +}\dff 
\left(\qff 
\sigma\pff \omult_{\fff 1}\qff q\fff,\off \mathcal{N}
\qff\right)
\off \in\off
K^{\dff 1}\trf 
\bigl(\qff 
B\dff M\qff \times\qff B\dff V\fff,\off  
(\trf S\dff M\dff \cup\dff B\dff M_{\trf L}\trf)\dff \times\dff B\dff V
\qff \cup\qff
B\dff M\dff \times\dff S\dff V
\qff\bigr)
\pff.
\]

\vspace{-12pt}\vspace{4.5pt}
Also,\pss $q$\dss leads\sss to\sss the class\sss
$d\trf(\qff \widehat{q}\pff)
\off =\off
d\trf(\trf \pi^{\fff *}\dff E\fff'\fff,\qff \pi^{\fff *}\dff F\fff'\fff,\qff \widehat{q}\pff)
\off \in\off
K^{\trf 0}\dff(\trf B\dff V\fff,\pff S\dff V\trf)$\nnsp.\oss

\mypar{Theorem.}{epsilon-product}
$\dis
\varepsilon^{\dff +}\dff 
\left(\qff 
\sigma\pff \omult_{\fff 1}\qff q\fff,\off \mathcal{N}
\qff\right)
\off =\off\dff
\varepsilon^{\dff +}\fff (\trf \sigma\fff,\pff N\qff)
\off \widehat{\otimes}\off\dff
d\trf(\qff \widehat{q}\pff)$\nnsp.\oss

\proof
In order\sss to find\sss
$\varepsilon^{\dff +}\dff 
\left(\qff 
\sigma\pff \omult_{\fff 1}\qff q\fff,\off \mathcal{N}
\qff\right)$\sss
we need\dss to deform\sss first\sss the symbol\sss
$\sigma\pff \omult_{\fff 1}\qff q$\sss 
to a bundle-like symbol.\oss
As we saw,\oss the deformation\sss
$\left(\qff \sigma\pff \omult_{\fff 1}\qff q \pff\right)_{\dff t}$\nsp,\qss
$t\qff \in\qff [\trf 0\fff,\qff 1 \trf]$\sss
connects\sss
$\sigma\pff \omult_{\fff 1}\qff q$\sss
with\sss the bundle-like symbol\sss
$\left(\qff \sigma\pff \omult_{\fff 1}\qff q \pff\right)_{\trf 0}$\nsp.\oss

Next,\oss we need\sss to\sss take\sss the canonical\sss extension of\dss
the symbol\dss
$\left(\qff \sigma\pff \omult_{\fff 1}\qff q \pff\right)_{\trf 0}$\sss 
to\sss
$D\trf(\trf M\dff \times\dff V\trf)$\nnsp.\oss
As we saw,\oss
$\left(\qff \mathbf{\bm{\sigma}\pff \omult_{\fff 1}\qff q} \pff\right)_{\trf 0}$\sss
is\dss equal\dss to\sss the canonical\sss extension.\oss
Therefore by\sss the definition\vspace{3.75pt}
\[
\quad
\varepsilon^{\dff +}\dff 
\left(\qff 
\sigma\pff \omult_{\fff 1}\qff q\fff,\off \mathcal{N}
\qff\right)
\off =\off\dff
\partial\qff\dff
\mathbb{P}^{\dff +}\fff 
\left(\qff
\left(\qff \mathbf{\bm{\sigma}\pff \omult_{\fff 1}\qff q} \pff\right)_{\trf 0}
\qff\right)
\qff,
\]

\vspace{-12pt}\vspace{3pt}
where\dss $\partial$\dss is\dss the coboundary\dss map\vspace{4.5pt}
\[
\quad
K^{\trf 0}\dff (\qff D\trf(\trf M\dff \times\dff V\trf)\qff)
\qff \ttoo\qff
K^{\dff 1}\dff (\trf B\trf(\trf M\dff \times\dff V\trf)\fff,\off  D\trf(\trf M\dff \times\dff V\trf)\trf)
\qff.
\]

\vspace{-12pt}\vspace{4.5pt}
The deformation\sss
$\left(\qff \mathbf{\bm{\sigma}\pff \omult_{\fff 1}\qff q} \pff\right)_{\dff t}$\dss
shows\sss that\sss the bundle\sss
$\mathbb{P}^{\dff +}\fff 
\left(\qff
\left(\qff \mathbf{\bm{\sigma}\pff \omult_{\fff 1}\qff q} \pff\right)_{\trf 0}
\qff\right)$\sss
is\dss isomorphic\sss to\vspace{4.5pt}
\[
\quad
\mathbb{P}^{\dff +}\fff 
\left(\qff
\left(\qff \mathbf{\bm{\sigma}\pff \omult_{\fff 1}\qff q} \pff\right)_{\dff 1}
\qff\right)
\off =\off\dff
\mathbb{P}^{\dff +}\fff 
\left(\qff
\mathbf{\bm{\sigma}\pff \omult_{\fff 1}\qff q}
\qff\right)
\pff.
\]

\vspace{-12pt}\vspace{4.5pt}
Let\sss $p$ be an extension of $\bm{\sigma}\nsp$ to \nsp$B\dff M$\nnsp.\qss
Lemma\qss \ref{extensions}\qss implies\sss that\sss \vspace{3pt}
\[
\quad
\mathbb{P}^{\dff +}\fff 
\left(\qff
\mathbf{\bm{\sigma}\pff \omult_{\fff 1}\qff q}
\qff\right)
\]

\vspace{-12pt}\vspace{3pt}
is\dss isomorphic\sss to\sss the bundle\sss
$\mathbb{P}^{\dff +}\fff 
\left(\qff
p\pff \omult_{\fff 1}\qff \widehat{q}
\off \bigl|\pff
D\trf(\trf M\dff \times\dff V\trf)
\qff\right)$\dss
and\dss hence\vspace{4.5pt}
\begin{equation}
\label{epsilon-dp}
\quad
\varepsilon^{\dff +}\dff 
\left(\qff 
\sigma\pff \omult_{\fff 1}\qff q\fff,\off \mathcal{N}
\qff\right)
\off =\off\dff
\partial\qff\dff
\mathbb{P}^{\dff +}\fff 
\left(\qff
p\pff \omult_{\fff 1}\qff \widehat{q}
\off \bigl|\pff
D\trf(\trf M\dff \times\dff V\trf)
\qff\right)
\pff.
\end{equation}

\vspace{-12pt}\vspace{4.5pt}
Let\sss us now\sss identify\sss
$B\dff (\trf M\dff \times\dff V\trf)$\sss
with\sss $B\dff M\qff \times\qff B\dff V$\dss by\sss $f$\nnsp.\oss
Then we can\sss replace\sss the restriction\vspace{4.5pt}
\[
\quad
p\pff \omult_{\fff 1}\qff \widehat{q}
\off \bigl|\pff
D\trf(\trf M\dff \times\dff V\trf)
\]

\vspace{-12pt}\vspace{4.5pt}
by\sss the restriction\dss
$p\pff \omult_{\fff 1}\qff \widehat{q}
\off \bigl|\pff
D\dff M\dff \times\dff B\dff V
\qff \cup\qff
B\dff M\dff \times\dff S\dff V$\dnsp.\oss
Let\sss us apply\trs Lemma\qss \ref{coboundary-mult}\qss to\sss \vspace{4.5pt}
\[
\quad
B\dff M\dff,\off\qff
D\dff M\dff,\off\qff
B\dff V\dff,\off\qff
S\dff V\dff,\off\qff
p\dff,\off\qff
\widehat{q}
\]

\vspace{-12pt}\vspace{4.5pt}
in\sss the roles of\trs 
$X\dff,\off\qff
Y\dff,\off\qff
X\fff'\dff,\off\qff
Y\fff'\dff,\off\qff
p\dff,\off\qff
q$\dss
respectively.\oss
Lemma\qss \ref{coboundary-mult}\qss implies\sss that\vspace{4.5pt}
\[
\quad
\partial\qff\fff 
\mathbb{P}^{\dff +}\fff 
\left(\qff 
p\pff \omult_1\pff \widehat{q}
\off \bigl|\pff
D\dff M\dff \times\dff B\dff V
\qff \cup\qff
B\dff M\dff \times\dff S\dff V
\qff\right)
\off\dff =\off\qff
\partial\qff\fff 
\mathbb{P}^{\dff +}\fff (\qff p\trf |\trf D\dff M\qff)
\off \widehat{\otimes}\off\dff
d\trf(\qff \widehat{q}\pff)
\pff.
\]

\vspace{-12pt}\vspace{4.5pt}
The\sss theorem\sss immediately\sss follows\sss from\sss this equality\sss together\sss with\sss
the equality\qss (\ref{epsilon-dp})\qss
and\dss the fact\sss that\dss
$\varepsilon^{\dff +}\fff (\trf \sigma\fff,\pff N\qff)
\off =\off
\partial\qff\fff
\mathbb{P}^{\dff +}\fff (\qff p\trf |\trf D\dff M\qff)$\nnsp.\oss  \eproof

\myuppar{The\sss topological\dss index and\dss the product\sss of\dss symbols.}
Let\sss $X\off =\off M\dff \times\dff V$\dnsp.\oss
Let\sss us embed\sss $M$\sss into\sss
$\rrr_{\qff \geq\trf 0}^{\dff m}$\sss
for some $m$ in such a way\sss that\sss
$L\off =\off M\qff \cap\trf (\trf \rrr^{\dff m\dff -\dff 1}\dff \times\dff 0\trf)$\sss
and\sss $M$\sss is\dss transverse\sss to\sss
$\rrr^{\dff m\dff -\dff 1}\dff \times\dff 0$\sss
in\sss $\rrr^{\dff m}$\dnsp.\oss
Let\sss $\mathbb{N}_{\trf M}$\sss be\sss the normal\dss bundle\sss 
to\sss $M$\sss in\sss $\rrr^{\dff m}$\dnsp.\oss
Let\sss us embed $V$ into\sss $\rrr^{\dff n}$\sss for some $n$\sss
and\dss let\sss $\mathbb{N}_{\trf V}$\sss be\sss the normal\dss bundle\sss 
to $V$ in\sss $\rrr^{\dff n}$\dnsp.\oss
These embeddings\sss lead\dss to an embedding of\dss $X\off =\off M\dff \times\dff V$ into\sss 
$\rrr_{\qff \geq\trf 0}^{\dff m\dff +\dff n}
\off =\off
\rrr_{\qff \geq\trf 0}^{\dff m}\dff \times\dff \rrr^{\dff n}$\sss
with\sss the normal\dss bundle\sss 
$\mathbb{N}
\off =\off
\mathbb{N}_{\trf M}\qff \times\pff \mathbb{N}_{\trf V}$\nsp.\oss

The roundness of\dss the unit\dss balls in\sss the\sss tangent\sss spaces\dss
is\dss irrelevant,\oss and we can\sss redefine $B\dff X$
as\sss the product\sss $B\dff M\dff \times\dff B\dff V$\dnsp.\oss
Then $S\fff M$ should\dss be redefined as
$(\trf B\dff M\dff \times\dff S\dff V\trf)
\qff \cup\qff
(\trf S\fff M\dff \times\dff B\dff V\trf)$\nnsp.\oss
We can deal\sss with\sss the normal\dss bundles in\sss the same way.\oss
Namely,\oss
let\sss $\mathbb{U}_{\trf M}$\sss and\sss $\mathbb{U}_{\trf V}$\sss
be\sss the bundles of\dss unit\dss balls in\sss the\sss lifts of\dss
$\mathbb{N}_{\trf M}\dff \oplus\dff \mathbb{N}_{\trf M}$\sss
and\sss
$\mathbb{N}_{\trf V}\dff \oplus\dff \mathbb{N}_{\trf V}$\sss
to\sss $B\fff M$\sss and\sss $B\dff V$ respectively.\oss
Similarly,\oss 
let\sss $\mathbb{S}_{\trf M}$\sss and\sss $\mathbb{S}_{\trf V}$\sss
be\sss the bundles of\dss unit\sss spheres in\sss the\sss lifts of\dss
$\mathbb{N}_{\trf M}\dff \oplus\dff \mathbb{N}_{\trf M}$\sss
and\sss
$\mathbb{N}_{\trf V}\dff \oplus\dff \mathbb{N}_{\trf V}$\sss
to\sss $B\fff M$\sss and\sss $B\dff V$ respectively.\oss
Then\sss 
$\mathbb{U}
\off =\off 
\mathbb{U}_{\trf M}\qff \times\pff \mathbb{U}_{\trf V}$\sss
is\dss a bundle over\sss $B\dff X\off =\off B\fff M\dff \times\dff B\dff V$\sss
which we can use instead of\dss the bundle of\dss the unit\dss balls\sss
in\sss the\sss lift\sss of\dss 
$\mathbb{N}\dff \oplus\dff \mathbb{N}$\sss
to\sss $B\dff X$\nnsp,\oss
and\dss then\vspace{3pt}
\[
\quad
\mathbb{S}
\off =\off
\bigl(\qff \mathbb{U}_{\trf M}\dff \times\trf \mathbb{S}_{\trf V} \qff\bigr)
\qff \cup\qff
\bigl(\qff \mathbb{S}_{\trf M}\dff \times\trf \mathbb{U}_{\trf V} \qff\bigr)
\qff \ttoo\qff
B\fff M\dff \times\dff B\dff V
\off =\off
B\dff X
\]

\vspace{-12pt}\vspace{3pt}
should\sss be used\sss instead of\dss the bundle of\dss unit\sss spheres in this\sss lift.\oss
With\sss these conventions\sss
$D\dff X
\off =\off
(\trf B\dff M\dff \times\dff S\fff V\trf)
\qff \cup\qff
(\trf S\dff M\dff \times\dff B\fff V\trf)
\qff \cup\qff
(\trf B\fff M_{\trf L}\dff \times\dff B\dff V\trf)$\nnsp.\oss
Also,\oss\vspace{3pt}
\[
\quad
\bigl(\qff 
\mathbb{U}\trf |\trf B\dff X\dff,\off
(\trf \mathbb{S}\dff |\trf B\dff X\trf)
\qff \cup\qff 
(\trf \mathbb{U}\dff |\trf D\dff X\trf)
\qff\bigr)
\qff
\]

\vspace{-34.5pt}
\[
\quad
=\off
\bigl(\qff 
\mathbb{U}_{\trf M}\trf |\trf B\dff M\dff,\off
(\trf \mathbb{S}_{\trf M}\trf |\trf B\dff M\trf)
\qff \cup\qff 
(\trf \mathbb{U}_{\trf M}\trf |\trf D\dff M\trf)
\qff\bigr)
\off \times\off
\bigl(\qff 
\mathbb{U}_{\trf V}\trf |\trf B\dff V\dff,\off
(\trf \mathbb{S}_{\trf V}\trf |\trf B\dff V\trf)
\qff \cup\qff 
(\trf \mathbb{U}_{\trf V}\trf |\trf S\dff V\trf)
\qff\bigr)
\qff
\]

\vspace{-12pt}\vspace{3pt}
and\dss hence\vspace{3pt}
\[
\quad
\mathbb{U}\trf |\trf B\dff X\qff\bigl/\qff
\bigl(\qff
(\trf \mathbb{S}\dff |\trf B\dff X\trf)
\qff \cup\qff 
(\trf \mathbb{U}\dff |\trf D\dff X\trf)
\qff\bigr)
\qff
\]

\vspace{-34.5pt}
\[
\quad
=\off
\bigl(\qff 
\mathbb{U}_{\trf M}\trf |\trf B\dff M\qff\bigl/\qff
(\trf \mathbb{S}_{\trf M}\trf |\trf B\dff M\trf)
\qff \cup\qff 
(\trf \mathbb{U}_{\trf M}\trf |\trf D\dff M\trf)
\qff\bigr)
\off \wedge\off
\bigl(\qff 
\mathbb{U}_{\trf V}\trf |\trf B\dff V\qff\bigl/\qff
(\trf \mathbb{S}_{\trf V}\trf |\trf B\dff V\trf)
\qff \cup\qff 
(\trf \mathbb{U}_{\trf V}\trf |\trf S\dff V\trf)
\qff\bigr)
\qff.
\]

\vspace{-12pt}\vspace{3pt}
Therefore\sss the canonical\sss quotient\sss map\vspace{3pt}
\[
\quad
c_{\trf X}\dff \colon\dff
S^{\dff 2 m}\dff \wedge\qff S^{\dff 2 n}
\off =\off
S^{\dff 2 m\qff +\qff 2 n}
\qff \ttoo\qff
\mathbb{U}\dff |\trf B\dff X\qff\bigl/\dff
\bigl(\qff 
(\trf \mathbb{S}\dff |\trf B\dff X\trf)
\qff \cup\qff 
(\trf \mathbb{U}\dff |\trf D\dff X\trf)
\qff\bigr)
\qff
\]

\vspace{-12pt}\vspace{3pt}
is\dss equal\dss to\sss the wedge product\sss $c_{\trf M}\dff \wedge\dff c_{\trf V}$\sss
of\dss the canonical\sss quotient\sss maps\vspace{3pt}
\[
\quad
c_{\trf M}\dff \colon\dff
S^{\dff 2 m}
\qff \ttoo\qff
\mathbb{U}_{\trf M}\trf |\trf B\dff M\qff\bigl/\qff
\bigl(\qff 
(\trf \mathbb{S}_{\trf M}\trf |\trf B\dff M\trf)
\qff \cup\qff 
(\trf \mathbb{U}_{\trf M}\trf |\trf D\dff M\trf)
\qff\bigr)
\quad
\mbox{and}\quad
\]

\vspace{-34.5pt}
\[
\quad
c_{\trf V}\dff \colon\dff
S^{\dff 2 n}
\qff \ttoo\qff
\mathbb{U}_{\trf V}\trf |\trf B\dff V\qff\bigl/\qff
\bigl(\qff 
(\trf \mathbb{S}_{\trf V}\trf |\trf B\dff V\trf)
\qff \cup\qff 
(\trf \mathbb{U}_{\trf V}\trf |\trf S\fff V\trf)
\qff\bigr)
\pff.
\]

\vspace{-12pt}\vspace{3pt}
The\dss Thom\dss isomorphism\vspace{3pt}
\[
\quad
\mathrm{Th}_{\trf X}\dff \colon\dff
K^{\dff 1}\dff (\trf B\dff X\fff,\off  D\dff X\trf)
\qff \ttoo\qff
K^{\dff 1}\dff \bigl(\qff \mathbb{U}\dff |\trf B\dff X\fff,\off 
(\trf \mathbb{S}\dff |\trf B\dff X\trf)
\qff \cup\qff 
(\trf \mathbb{U}\dff |\trf D\dff X\trf)
\qff\bigr)
\qff
\]

\vspace{-12pt}\vspace{3pt}
is\dss the composition of\dss the pull-back\sss map\sss
$K^{\dff 1}\dff (\trf B\dff X\fff,\off  D\dff X\trf)
\qff \ttoo\qff
K^{\dff 1}\dff (\qff \mathbb{U}\dff |\trf B\dff X\fff,\off 
\mathbb{U}\dff |\trf D\dff X \qff)$\sss
with\sss the multiplication\sss by\sss the\dss Thom\sss class\sss
$\lambda_{\trf X}
\qff \in\qff 
K^{\trf 0}\fff(\trf \mathbb{U}_{\trf X}\dff,\pff \mathbb{S}_{\trf X}\trf)$\nnsp.\oss
Similarly,\oss the\dss Thom\dss isomorphisms\vspace{3pt}
\[
\quad
\mathrm{Th}_{\trf M}\dff \colon\dff
K^{\dff 1}\dff (\trf B\fff M\fff,\off  D\dff M\trf)
\qff \ttoo\qff
K^{\dff 1}\dff \bigl(\qff \mathbb{U}_{\trf M}\dff |\trf B\fff M\fff,\off 
(\trf \mathbb{S}_{\trf M}\dff |\trf B\dff M\trf)
\qff \cup\qff 
(\trf \mathbb{U_{\trf M}}\dff |\trf D\dff M\trf)
\qff\bigr)
\qff
\]

\vspace{-34.5pt}
\[
\quad
\mathrm{Th}_{\trf V}\dff \colon\dff
K^{\dff 0}\dff (\trf B\dff V\fff,\off  S\fff V\trf)
\qff \ttoo\qff
K^{\dff 0}\dff \bigl(\qff \mathbb{U_{\trf V}}\dff |\trf B\dff V\fff,\off 
(\trf \mathbb{S}_{\trf V}\dff |\trf B\dff V\trf)
\qff \cup\qff 
(\trf \mathbb{U}_{\trf V}\dff |\trf S\fff V\trf)
\qff\bigr)
\qff
\]

\vspace{-12pt}\vspace{3pt}
are\sss the compositions of\dss the pull-back\sss maps with\sss
the multiplication\sss by\sss the\dss Thom\sss classes\sss
$\lambda_{\trf M}\dff,\off \lambda_{\trf V}$\nsp.\oss
The multiplicative\sss formula
$\lambda_{\trf X}
\off =\off 
\lambda_{\trf M}\dff \boxtimes\dff \lambda_{\trf V}$\sss
implies\sss that\sss\vspace{2.75pt}
\[
\quad
\mathrm{Th}_{\trf X}\dff(\trf a\qff \boxtimes\qff b\qff)
\off =\off
\mathrm{Th}_{\trf M}\dff(\trf a\trf)
\pff \boxtimes\pff\fff
\mathrm{Th}_{\trf V}\dff(\trf b\qff)
\] 

\vspace{-12pt}\vspace{2.75pt}
for every\sss
$a\qff \in\qff K^{\dff 1}\fff(\trf B\dff M\dff,\qff D\dff M\trf)$\sss
and\sss
$b\qff \in\qff K^{\trf 0}\fff(\trf B\dff V\fff,\qff S\fff V\trf)$\nnsp,\oss
where\sss $\boxtimes$\sss denotes\sss the external\dss tensor\sss products.\oss
In\sss turn,\oss this implies\sss that\vspace{3pt}
\begin{equation}
\label{index-product}
\quad
c_{\trf X}^{\dff *}\trf
\left(\qff \mathrm{Th}_{\trf X}\dff(\trf a\qff \boxtimes\qff b\qff) \qff\right)
\off =\off
c_{\trf M}^{\dff *}\trf
\bigl(\qff 
\mathrm{Th}_{\trf M}\dff(\trf a\trf) \qff\bigr)
\pff\fff \boxtimes\pff\fff
c_{\trf V}^{\dff *}\trf
\bigl(\qff 
\mathrm{Th}_{\trf V}\dff(\trf b\qff) \qff\bigr)
\pff,
\end{equation} 

\vspace{-12pt}\vspace{3pt}
where\sss $\boxtimes$\sss on\sss the right\sss denotes\sss 
the external\dss tensor product\vspace{2.75pt}
\[
\quad
\boxtimes\dff \colon\dff
K^{\dff 1}\fff\left(\trf S^{\dff 2 m}\trf\right)
\qff \otimes\qff
K^{\trf 0}\fff\left(\trf S^{\dff 2 n}\trf\right)
\qff \ttoo\qff
K^{\dff 1}\fff\left(\trf S^{\dff 2 m\qff +\qff 2 n}\trf\right)
\pff.
\]

\vspace{-12pt}\vspace{2.75pt}
As in\dss Section\qss \ref{symbols-conditions}\qss we pretend\dss that\sss
we don't\dss know\sss that\sss
$K^{\dff 1}\fff\left(\trf S^{\dff 2 m}\trf\right)
\off =\off
K^{\dff 1}\fff\left(\trf S^{\dff 2 m\qff +\qff 2 n}\trf\right)
\off =\off
0$\nnsp.\oss

\mypar{Theorem.}{t-index-product}
\emph{Suppose\sss that\sss everything continuously depends on a parameter $z\qff \in\qff Z$\nnsp.\oss
Then}\vspace{2.75pt}
\[
\quad
\ti\qff  
\left(\qff 
\sigma\pff \omult_{\fff 1}\qff q\fff,\off \mathcal{N}
\qff\right)
\off =\off
\ti\qff 
\left(\trf  \sigma\fff,\pff N
\qff\right) 
\off \otimes\off
\ti\qff 
(\trf q\qff)
\off
\]

\vspace{-12pt}\vspace{2.75pt}
\emph{as elements of\dss $K^{\dff 1}\fff(\trf Z\trf)$\nnsp.\oss}

\proof
Let\sss us apply\sss the multiplicative property\qss (\ref{index-product})\qss to\sss 
$a\off =\off \varepsilon^{\dff +}\fff (\trf \sigma\fff,\pff N\qff)$
and\sss
$b\off =\off d\trf(\qff \widehat{q}\pff)$\nnsp,\oss
where $\sigma\fff,\pff N$ and $q$\sss are as\sss in\dss Theorem\qss \ref{epsilon-product}.\oss
Together with\trs Theorem\qss \ref{epsilon-product}\qss this
shows\sss that\vspace{3pt}
\[
\quad
c_{\trf X}^{\dff *}\trf
\left(\qff \mathrm{Th}_{\trf X}\dff
\left(\trf \varepsilon^{\dff +}\dff 
\left(\qff 
\sigma\pff \omult_{\fff 1}\qff q\fff,\off \mathcal{N}
\qff\right)
\qff\right) 
\qff\right)
\off =\off
c_{\trf M}^{\dff *}\trf
\bigl(\qff 
\mathrm{Th}_{\trf M}\dff
\left(\trf \varepsilon^{\dff +}\fff (\trf \sigma\fff,\pff N\qff)
\trf\right) 
\qff\bigr)
\off\qff\boxtimes\off\qff
c_{\trf V}^{\dff *}\trf
\left(\qff 
\mathrm{Th}_{\trf V}\dff
\left(\trf
d\trf
\left(\qff \widehat{q}
\pff\right)
\qff\right) 
\qff\right)
\off.
\]

\vspace{-12pt}\vspace{3pt}
In\sss terms of\dss the\sss topological\dss pre-index\sss $t\qff(\trf \bullet\trf)$\sss 
this means\sss that\vspace{3pt}\vspace{-0.125pt}
\[
\quad
t\qff  
\left(\qff 
\sigma\pff \omult_{\fff 1}\qff q\fff,\off \mathcal{N}
\qff\right)
\off =\off
t\qff 
\left(\trf  \sigma\fff,\pff N
\qff\right) 
\off \boxtimes\off
t\qff 
(\trf q\qff)
\off.
\]

\vspace{-12pt}\vspace{3pt}
At\sss this point\sss we merely\sss proved\sss in a complicated\sss way\sss the identity\sss
$0
\off =\off 
0\off \boxtimes\off t\qff (\trf q\qff)$\nnsp.\oss
But,\oss as\sss in\dss Section\qss \ref{symbols-conditions},\oss
this proof\dss works when everything\sss depends on\sss the\sss parameter\sss $z\qff \in\qff Z$\nnsp,\oss
and\sss the\sss last\sss displayed\sss formula\sss should\dss be understood as an equality\sss in\sss
$K^{\dff 1}\fff\left(\trf S^{\dff 2 m\qff +\qff 2 n}\dff \times\dff Z\fff,\off 
p_{\dff 0}\dff \times\dff Z\trf\right)$\nnsp.\oss
It\sss remains\sss to apply\sss the\dss Bott\dss periodicity\sss map.\oss  \eproof

\newpage
\mysection{Abstract\qss boundary\qss problems}{abstract-index}

\myuppar{The framework.}
Let\sss $H_{\dff 0}$\sss be a separable\sss Hilbert\sss space and\sss $H_{\fff 1}$\sss
be a dense subspace of\dss $H_{\dff 0}$\nsp,\oss
which\dss is\dss a\sss Hilbert\sss space in\sss its own\sss right\qss
({\halfff}the\dss Hilbert\sss space structure of\dss $H_{\fff 1}$\sss is\dss
not\sss induced\sss from\sss $H_{\dff 0}$\nsp).\oss
Let\sss $H^{\dff \partial}$ also be a\sss Hilbert\sss space and\dss let\sss
$H_{\dff 1/2}^{\dff \partial}$\sss be a dense subspace of\dss $H^{\dff \partial}$\nsp,\oss
which\dss is\dss a\sss Hilbert\sss space in\sss its own\sss right.\oss
Suppose\sss that\sss the inclusion maps\sss 
$\iota\dff \colon\dff
H_{\fff 1}\qff \ttoo\qff H_{\dff 0}$\sss
and\sss
$H_{\dff 1/2}^{\dff \partial}\qff \ttoo\qff H^{\dff \partial}$\sss
are bounded operators with\sss respect\sss
to\sss these\sss Hilbert\sss space structures.\oss
Let\sss \vspace{0pt}
\[
\quad
\gamma\dff \colon\dff
H_{\fff 1}\qff \ttoo\qff H_{\dff 1/2}^{\dff \partial}
\]

\vspace{-12pt}\vspace{0pt}
be a surjective bounded operator such\sss that
the kernel\dss $\kernel\fff \gamma$\sss is\dss dense in\sss
$H_{\dff 0}$\nsp.\oss
Let\sss
$\sco{\dff \bullet\dff,\qff \bullet \dff}_{\dff 0}$\sss
and\sss
$\sco{\dff \bullet\dff,\qff \bullet \dff}_{\dff \partial}$\sss
be\sss the scalar products in\sss $H_{\dff 0}$\sss 
and\sss $H^{\dff \partial}$\sss respectively.\oss

\myuppar{Projections.}
A\qss \emph{projection}\pss (not\sss necessarily\sss orthogonal\fff)\qss 
in a\dss Hilbert\sss space $H$\sss is a bounded operator
$p\dff \colon\dff H\qff \ttoo\qff H$\sss such\sss that\sss $p\dff \circ\dff p\off =\off p$\sss
and\sss hence\sss 
$p\dff \circ\dff (\trf 1\qff -\qff p\trf)
\off =\off
(\trf 1\qff -\qff p\trf)\dff \circ\dff p
\off =\off 
0$\nnsp.\oss
If\dss $p$\sss is\dss a projection,\oss then\sss $1\qff -\qff p$\sss
is\dss also a projection,\oss and,\oss clearly,\qss
$\image\dff (\trf 1\qff -\qff p\trf)
\qff \subset\qff
\kernel\dff p$\sss
and\sss
$\image\dff p
\qff \subset\qff
\kernel\dff (\trf 1\qff -\qff p\trf)$\nnsp.\oss
The identity\sss $p\qff +\qff (\trf 1\qff -\qff p\trf)\off =\off 1$\sss 
implies\sss that\sss
$\image\dff p\qff +\qff \image\dff (\trf 1\qff -\qff p\trf)\off =\off H$\sss
and\sss
$\kernel\dff p\qff \cap\qff \kernel\dff (\trf 1\qff -\qff p\trf)\off =\off 0$\nnsp.\oss
By\sss taking into account\sss the above inclusions we see\sss that\sss
$\image\dff p\qff \cap\qff \image\dff (\trf 1\qff -\qff p\trf)\off =\off 0$\sss
and\sss
$\kernel\dff p\qff +\qff \kernel\dff (\trf 1\qff -\qff p\trf)\off =\off H$\nnsp.\oss
It\sss follows that\sss
$\image\dff (\trf 1\qff -\qff p\trf)
\off =\off
\kernel\dff p$\sss
and\sss
$\image\dff p
\off =\off
\kernel\dff (\trf 1\qff -\qff p\trf)$\nnsp.\oss
Of\dss course,\oss all\dss this\dss is\dss well\dss known.\oss

\myuppar{Boundary conditions.}
Let\sss $A\dff \colon\dff H_{\fff 1}\qff \ttoo\qff H_{\dff 0}$\sss
be a bounded operator.\oss
A\qss \emph{boundary\sss condition}\pss for $A$\sss is\dss a projection\sss
$\Pi \dff \colon\dff
H^{\dff \partial}\qff \ttoo\qff H^{\dff \partial}$\sss
leaving\sss the subspace\sss $H_{\dff 1/2}^{\dff \partial}$\sss invariant.\oss
We will\sss denote by\sss $\Pi_{\dff 1/2}$\sss the projection\sss
$H_{\dff 1/2}^{\dff \partial}\qff \ttoo\qff H_{\dff 1/2}^{\dff \partial}$\sss
induced\sss by\sss $\Pi$\sss and assume\sss that\sss $\Pi_{\dff 1/2}$\sss is\dss bounded.\oss

The\qss \emph{boundary\dss problem}\qss $A\fff,\qff \Pi $\sss
is\dss the problem of\dss finding solutions $u$ of\dss the equation\sss
$A\dff u\off =\off 0$\sss subject\sss to\sss the condition\sss
$(\trf 1\qff -\qff \Pi_{\dff 1/2}\trf)\dff \circ\dff \gamma\qff (\trf u\trf)\off =\off 0$\nnsp.\oss
Let\sss\vspace{1.5pt}
\[
\quad
G_{\dff 1/2}
\off =\off
\kernel\dff \Pi_{\dff 1/2}
\off =\off
\image\dff (\trf 1\qff -\qff \Pi_{\dff 1/2}\trf)
\quad
\mbox{and}\quad
\]

\vspace{-36pt}
\[
\quad
G 
\off =\off
\kernel\dff \Pi 
\off =\off
\image\dff (\trf 1\qff -\qff \Pi \trf)
\off.
\]

\vspace{-12pt}\vspace{1.5pt}
Clearly,\pss 
$G_{\dff 1/2}
\off =\off
\kernel\dff \Pi_{\dff 1/2}
\off =\off
\kernel\dff \Pi \qff \cap\qff H_{\dff 1/2}
\off =\off
G \qff \cap\qff H_{\dff 1/2}$\nsp.\oss
Let\sss\vspace{0.75pt}
\[
\quad
\Gamma\dff \colon\dff
H_{\fff 1}\qff \ttoo\qff G_{\dff 1/2}
\]

\vspace{-12pt}\vspace{0.75pt} 
be\sss the composition\sss
$\Gamma\off =\off (\trf 1\qff -\qff \Pi_{\dff 1/2}\trf)\dff \circ\dff \gamma$\nnsp,\oss
and\dss let\vspace{1.0pt}
\[
\quad
A\dff \oplus\dff \Gamma\dff \colon\dff
H_{\fff 1}\qff \ttoo\qff H_{\dff 0}\dff \oplus\dff G_{\dff 1/2}
\qff
\]

\vspace{-12pt}\vspace{1.0pt}
be\sss the operator defined\sss by\sss
$u\off \longmapsto\off (\trf A\dff u\fff,\qff \Gamma\dff u\trf)$\nnsp.\oss
In\sss terms of\dss the operator\sss $A\dff \oplus\dff \Gamma$ the boundary\sss problem\sss
$A\fff,\qff \Pi$\sss
is\dss the problem of\dss finding\sss the kernel\sss
$\kernel\dff A\dff \oplus\dff \Gamma$\dnsp.\oss

\myuppar{The\dss Lagrange\sss identity.}
Let\sss  
$\Sigma\dff \colon\dff 
H^{\dff \partial}\qff \ttoo\qff H^{\dff \partial}$\dss
be a self-adjoint\sss invertible bounded operator\sss
leaving\sss the subspace\sss $H_{\dff 1/2}^{\dff \partial}$\sss invariant\sss
and such\sss that\sss the\qss
\emph{abstract\dss Lagrange\dss identity}\vspace{1.5pt}
\begin{equation}
\label{lagrange}
\quad
\sco{\dff A\dff u\fff,\qff v \dff}_{\dff 0}
\qff -\qff
\sco{\dff u\dff,\qff A\dff v \dff}_{\dff 0}
\off =\off
\sco{\dff i\trf \Sigma\dff \gamma\dff u\dff,\qff \gamma\dff v \dff}_{\dff \partial}
\qff
\end{equation}

\vspace{-12pt}\vspace{1.5pt}
holds for every\sss $u\fff,\qff v\qff \in\qff H_{\fff 1}$\nsp.\oss
Such identities are also known as\qss
\emph{Green\dss formulas}.\oss
Note\sss that\dss 
$\Sigma\dff \gamma\dff u$\sss and\sss $\gamma\dff v$\sss
belong\sss to
$H_{\dff 1/2}^{\dff \partial}$\nsp,\oss
but\sss the scalar\sss product\sss on\sss the right\dss
is\dss taken\sss in\sss $H^{\dff \partial}$\nsp.\oss

\myuppar{Self-adjointness and\sss elliptic regulatiry.}
A boundary condition $\Pi $\sss for $A$\sss
and\sss the boundary\sss problem\sss $A\fff,\qff \Pi$\sss are said\sss to be\qss
\emph{self-adjoint}\pss if\dss $\Pi $ is\dss a self-adjoint\sss projection\sss
in\sss $H^{\dff \partial}$\sss and\vspace{0pt}
\[
\quad
\Sigma\trf(\qff \image\dff \Pi  \qff)
\off =\off
\kernel\dff \Pi 
\qff,
\]

\vspace{-12pt}\vspace{0pt}
and\qss \emph{symmetric}\pss if\dss the projection\sss $\Pi $\sss is\dss self-adjoint\sss and\sss
$\Sigma\trf(\qff \image\dff \Pi  \qff)
\off \subset\off
\kernel\dff \Pi $\nsp.\oss
Let\sss us\sss equip\sss $H_{\dff 0}\dff \oplus\dff G$\sss
with\sss the\dss Hermitian\dss scalar product\sss induced\sss from\sss
$H_{\dff 0}\dff \oplus\dff H^{\dff \partial}$\nsp,\oss
and\sss denote by ${\fff\bullet\fff}^{\trf \perp}$\sss orthogonal\sss
complements with respect\sss to\sss this product.\oss
We will\sss say\sss that\sss the boundary problem\sss  $A\fff,\qff \Pi$\sss
is\qss \emph{elliptic regular}\pss if\trs
$\image A\dff \oplus\dff \Gamma
\off =\off
V^{\trf \perp}$\sss
for a finitely dimensional\sss subspace\sss
$V\qff \subset\pff
H_{\fff 1}\dff \oplus\dff G$\nnsp.\oss
In\sss this case $V$\sss is\dss closed\sss in\sss $H_{\dff 0}\dff \oplus\dff G$\sss
and\dss hence\sss
$(\qff \image A\dff \oplus\dff \Gamma\qff)^{\trf \perp}
\off =\off
V$\dnsp.\oss

\mypar{Theorem.}{duality-lemma}
\emph{If\trs $A\fff,\qff \Pi$\sss is\dss a self-adjoint\dss elliptic regular\sss boundary\sss problem,\oss  
then\sss $(\qff \image A\dff \oplus\dff \Gamma\qff)^{\trf \perp}$\sss 
is\dss equal\dss to\sss the image of\pss
$\kernel\dff A\dff \oplus\dff \Gamma$ under\sss the map\sss
$u
\qff \longmapsto\qff 
(\trf u\fff,\qff i\trf \Sigma\dff \circ\trf \gamma\dff u\trf)$\dnsp.}

\proof
Let\sss $h\qff \in\qff H_{\fff 1}$\nsp,\qss
$u\qff \in\qff H_{\fff 1}$\nsp,\oss
and\sss 
$w\qff \in\qff G 
\off =\off
\kernel\dff \Pi $\nsp.\oss
Then\sss 
$A\dff \oplus\dff \Gamma\dff(\trf h\trf)
\off =\off
(\trf A\dff h\fff,\qff \Gamma\dff h\trf)$\sss
is\dss orthogonal\sss to\sss $(\trf u\fff,\qff w\trf)$\sss
if\dss and\dss only\dss if\vspace{1.5pt}\vspace{-0.125pt}
\begin{equation}
\label{ort}
\quad
\sco{\dff A\dff h\fff,\qff u \dff}_{\dff 0}
\off +\off
\sco{\dff \Gamma\dff h\fff,\qff w \dff}_{\dff \partial}
\off =\off
0
\qff.
\end{equation}

\vspace{-12pt}\vspace{1.5pt}
Since\sss $w\qff \in\qff \kernel\dff \Pi $\sss
and\sss $\Pi_{\dff 1/2}$\sss is\dss equal\sss to\sss $\Pi $\sss
on\sss the image of\dss $\gamma$\nnsp,\oss\vspace{3pt}
\[
\quad
\sco{\dff \Gamma\dff h\fff,\qff w \dff}_{\dff \partial}
\off =\off
\sco{\dff \gamma\dff h\fff,\qff w \dff}_{\dff \partial}
\off -\off
\sco{\dff \Pi \dff \circ\dff \gamma\dff h\fff,\qff w \dff}_{\dff \partial}
\]

\vspace{-36pt}
\[
\quad
\phantom{\sco{\dff \Gamma\dff h\fff,\qff w \dff}_{\dff \partial}
\off }
=\off
\sco{\dff \gamma\dff h\fff,\qff w \dff}_{\dff \partial}
\off -\off
\sco{\dff \gamma\dff h\fff,\qff \Pi \dff w \dff}_{\dff \partial}
\off\dff =\off
\sco{\dff \gamma\dff h\fff,\qff w \dff}_{\dff \partial}
\qff.
\]

\vspace{-12pt}\vspace{3pt}
Now\qss (\ref{lagrange})\qss and\dss the self-adjointness\sss of\dss $\Sigma$\sss imply\sss that\sss
the\sss left\sss hand side of\qss (\ref{ort})\qss it\dss equal\dss to\vspace{3pt}
\[
\quad
\sco{\dff A\dff h\fff,\qff u \dff}_{\dff 0}
\off +\off
\sco{\dff \gamma\dff h\fff,\qff w \dff}_{\dff \partial}
\off
\]

\vspace{-33pt}\vspace{1.5pt}
\[
\quad
=\off
\sco{\dff h\fff,\qff A\dff u \dff}_{\dff 0}
\off +\off
\sco{\dff i\trf \Sigma\dff \circ\trf \gamma\dff h\dff,\qff \gamma\dff u \dff}_{\dff \partial}
\off +\off
\sco{\dff \gamma\dff h\fff,\qff w \dff}_{\dff \partial}
\]

\vspace{-36pt}\vspace{1.5pt}
\[
\quad
=\off
\sco{\dff h\fff,\qff A\dff u \dff}_{\dff 0}
\off -\off
\sco{\dff \gamma\dff h\dff,\qff i\trf \Sigma\dff \circ\trf \gamma\dff u \dff}_{\dff \partial}
\off +\off
\sco{\dff \gamma\dff h\fff,\qff w \dff}_{\dff \partial}
\]

\vspace{-36pt}\vspace{1.5pt}
\[
\quad
=\off
\sco{\dff h\fff,\qff A\dff u \dff}_{\dff 0}
\off +\off
\sco{\dff \gamma\dff h\dff,\qff 
w\qff -\qff i\trf \Sigma\dff \circ\dff \gamma\trf u \dff}_{\dff \partial}
\off.
\]

\vspace{-12pt}\vspace{3pt}
Suppose\sss that\qss (\ref{ort})\sss holds for every\sss
$h\qff \in\qff H_{\fff 1}$\dnsp.\oss
Then\sss 
$\sco{\dff h\fff,\qff A\dff u \dff}_{\dff 0}
\off =\off
0$\sss
if\dss $\gamma\dff h\off =\off 0$\nnsp.\oss
By our assumptions such\sss $h$\sss are dense in\sss $H_{\dff 0}$\nsp.\oss
It\sss follows\sss that\sss $A\dff u\off =\off 0$\nnsp.\oss
In\sss turn,\oss this implies\sss that\vspace{3pt}
\[
\quad
\sco{\dff \gamma\dff h\dff,\qff 
w\qff -\qff i\trf \Sigma\dff \circ\dff \gamma\trf u \dff}_{\dff \partial}
\off =\off 
0
\]

\vspace{-12pt}\vspace{3pt}
for every\sss $h\qff \in\qff H_{\fff 1}$\nsp.\oss
Since\sss $\gamma$\sss is\dss a map onto\sss $H_{\dff 1/2}^{\dff \partial}$\sss
and\sss $H_{\dff 1/2}^{\dff \partial}$\sss is\dss
dense in\sss $H^{\dff \partial}$\nsp,\oss
it\sss follows\sss that\sss\vspace{3pt}
\[
\quad
w\qff -\qff i\trf \Sigma\dff \circ\trf \gamma\dff u
\off =\off 
0
\off.
\]

\vspace{-12pt}\vspace{3pt}
Since\sss 
$w
\qff \in\qff 
G 
\off =\off 
\kernel\dff \Pi $\sss
and\sss $\Pi $\sss is\dss a self-adjoint\dss boundary condition,\oss 
it\sss follows\sss that\sss\vspace{3pt}
\[
\quad
\gamma\dff u
\off \in\off 
\Sigma^{\dff -\dff 1}\dff(\qff \kernel\dff \Pi  \qff)
\off =\off
\image\dff \Pi 
\off =\off
\kernel\dff(\trf 1\qff -\qff \Pi \trf)
\off
\]

\vspace{-12pt}\vspace{3pt}
and\dss hence\sss 
$u\qff \in\qff \kernel\dff \Gamma$\dnsp.\oss
Since\sss $A\dff u\off =\off 0$\nnsp,\oss
this\sss implies\sss that\sss
$u\qff \in\qff \kernel\dff A\dff \oplus\dff \Gamma$\dnsp.\oss
It\sss follows\sss that\dss if\sss 
$(\trf u\fff,\qff w\trf)$\sss is\dss orthogonal\dss to\sss
$A\dff \oplus\dff \Gamma\dff(\trf H_{\fff 1}\trf)$\nnsp,\oss
then\sss 
$u\qff \in\qff \kernel\dff A\dff \oplus\dff \Gamma$\sss
and\sss
$w\off =\off i\trf \Sigma\dff \circ\trf \gamma\dff u$\nnsp.\oss

Conversely,\oss let\sss
$u\qff \in\qff \kernel\dff A\dff \oplus\dff \Gamma$\sss
and\sss
$w\off =\off i\trf \Sigma\dff \circ\trf \gamma\dff u$\nnsp.\oss
Then\sss 
$\gamma\dff u
\qff \in\qff 
\kernel\dff (\trf 1\qff -\qff \Pi \trf)
\off =\off
\image\dff \Pi $\sss 
and\vspace{3pt}
\[
\quad
w
\off =\off
i\trf \Sigma\dff \circ\trf \gamma\dff u
\qff \in\qff \Sigma\dff(\qff \image\dff \Pi \qff)
\off =\off
\kernel\dff \Pi 
\off =\off
G 
\off.
\]

\vspace{-12pt}\vspace{3pt}
As we saw above,\oss this implies\sss that\sss
$\sco{\dff A\dff h\fff,\qff u \dff}_{\dff 0}
\off +\off
\sco{\dff \Gamma\dff h\fff,\qff w \dff}_{\dff \partial}$\trs
is\dss equal\sss to\vspace{3pt}
\[
\quad
\sco{\dff h\fff,\qff A\dff u \dff}_{\dff 0}
\off +\off
\sco{\dff \gamma\dff h\dff,\qff 
w\qff -\qff i\trf \Sigma\dff \circ\trf \gamma\dff u \dff}_{\dff \partial}
\off =\off
\sco{\dff h\fff,\qff 0 \dff}_{\dff 0}
\off +\off
\sco{\dff \gamma\dff h\dff,\qff 
0 \dff}_{\dff \partial}
\off =\off
0
\off,
\]

\vspace{-12pt}\vspace{3pt}
i.e.\dss
$(\trf A\dff h\fff,\qff \Gamma\dff h\trf)$\sss
is\dss orthogonal\sss to\sss $(\trf u\fff,\qff w\trf)$\sss
for every\sss $h\qff \in\qff H_{\fff 1}$\nsp.\oss
The\sss theorem\sss follows.\oss  \eproof

\mypar{Corollary.}{in-half}
\emph{If\trs $A\fff,\qff \Pi$\sss is\dss a self-adjoint\dss elliptic regular\sss boundary\sss problem,\oss  
then\sss $(\qff \image A\dff \oplus\dff \Gamma\qff)^{\trf \perp}$\sss
is\dss contained\sss in\sss $H_{\dff 1}\dff \oplus\dff G_{\dff 1/2}$\nsp.\oss}

\proof
Theorem\qss \ref{duality-lemma}\qss implies\sss that\dss if\dss 
$(\dff u\fff,\off w\trf)
\qff \in\qff 
H_{\dff 0}\dff \oplus\dff G $\sss 
is\dss orthogonal\dss to\sss
$A\dff \oplus\dff \Gamma\trf(\trf H_{\fff 1}\trf)$\nnsp,\oss
then\sss
$w\off =\off i\trf \Sigma\dff \circ\trf \gamma\dff u$\nnsp.\oss
Since\sss $\Sigma$\sss leaves\sss $H_{\dff 1/2}$\sss invariant,\oss
it\sss follows\sss that\sss
$w\qff \in\qff H_{\dff 1/2}$\nsp.\oss
In\sss turn,\oss this implies\sss that\sss 
$w
\qff \in\qff 
G \qff \cap\qff H_{\dff 1/2}
\off =\off
G_{\dff 1/2}$\nsp.\oss  \eproof

\myuppar{The induced\sss unbounded operator.}
Suppose\sss that\sss $\Pi $\sss is\dss a symmetric\sss boundary condition.\oss
Let\sss 
$A_{\trf \Gamma}\dff \colon\dff
\kernel\fff \Gamma\qff \ttoo\qff H^{\dff 0}$\sss
be\sss the operator\sss induced\sss by\sss
$A\dff \oplus\dff \Gamma$\dnsp.\oss
It\sss is\dss an unbounded operator\sss in\sss $H_{\dff 0}$\nsp.\oss
Since\sss $\kernel\dff \gamma$\nnsp,\oss and\dss hence\sss
$\kernel\fff \Gamma$\dnsp,\oss are dense in\sss $H_{\dff 0}$\nsp,\oss
the operator\sss $A_{\trf \Gamma}$\sss is\dss densely defined.\oss

Since\sss $\Pi $\sss is\dss a self-adjoint\sss projection,\oss
the kernels\sss
$\kernel\dff (\trf 1\qff -\qff \Pi \trf)$\sss
and\sss $\kernel\dff \Pi $\sss are orthogonal\sss in\sss
$H^{\dff \partial}$\nsp.\oss
Since\sss $\Pi $\sss is\dss a symmetric boundary condition,\pss
$\Sigma$\sss maps\sss
$\kernel\dff (\trf 1\qff -\qff \Pi \trf)
\off =\off
\image\dff \Pi $\sss 
into\sss $\kernel\dff \Pi $\nsp.\oss
It\sss follows\sss that\sss
$\sco{\dff i\trf \Sigma\dff \gamma\dff u\dff,\qff \gamma\dff v \dff}_{\dff \partial}
\off =\off
0$\sss
for every\sss $u\fff,\qff v\qff \in\qff \kernel\fff \Gamma$\dnsp,\oss
and\qss (\ref{lagrange})\qss implies\sss that\sss
$A_{\trf \Gamma}$\sss is\dss a symmetric operator.\oss
But\sss our current\sss assumptions are not\sss sufficient\sss to prove\sss that\sss
$A_{\trf \Gamma}$\sss is\dss a closed operator,\oss
and\dss that $A_{\trf \Gamma}$ is\dss a self-adjoint\sss 
if\dss the boundary condition\sss $\Pi $\sss
is\dss self-adjoint.\oss

\mypar{Theorem.}{reduction}
\emph{If\trs $A\fff,\qff \Pi$\sss is\dss a self-adjoint\dss elliptic regular\sss boundary\sss problem,\oss
then\dss the kernel\dss $\kernel\dff A_{\trf \Gamma}$\sss
is\dss equal\sss to\sss the orthogonal\sss complement\sss
of\dss the image\sss
$\image\fff A_{\trf \Gamma}
\qff \subset\qff 
H_{\dff 0}$\nsp.\oss}

\proof
Clearly,\pss
$\image\fff A_{\trf \Gamma}\dff \oplus\dff 0
\off =\off
(\trf \image\fff A\dff \oplus\dff \Gamma \trf)
\qff \cap\qff
(\trf H_{\dff 0}\dff \oplus\dff 0 \trf)$\nnsp.\oss
If\trs $\image A\dff \oplus\dff \Gamma
\off =\off
V^{\trf \perp}$\sss for\sss $V$ as in\sss the definition of\dss elliptic regularity,\oss
then\sss
$\image\fff A_{\trf \Gamma}\dff \oplus\dff 0
\off =\off
V^{\trf \perp}
\qff \cap\qff
(\trf H_{\dff 0}\dff \oplus\dff 0 \trf)$\nnsp,\oss
and\dss hence\sss
$(\trf \image\fff A_{\trf \Gamma} \trf)^{\trf \perp}
\off =\off
V_{\dff 0}$\nsp,\oss
where\sss $V_{\dff 0}$\sss is\dss the image of\dss $V$\sss
under\sss the projection\sss
$H_{\dff 0}\dff \oplus\dff G
\qff \ttoo\qff
H_{\dff 0}$\nsp.\oss
Since\sss 
$V\qff \subset\pff
H_{\fff 1}\dff \oplus\dff G$\nnsp,\oss
we have\sss
$V_{\dff 0}\qff \subset\qff H_{\dff 1}$\nsp.\oss
Therefore\sss the orthogonal\sss complement\sss
$(\trf \image\fff A_{\trf \Gamma} \trf)^{\trf \perp}$\sss
is\dss contained\sss in\sss $H_{\fff 1}$\nsp.\oss

Suppose\sss that\sss $v\qff \in\qff H_{\fff 1}$\sss
is\dss orthogonal\dss to\sss $\image\fff A_{\trf \Gamma}$\nnsp.\oss
Then\sss
$\sco{\dff A\dff u\fff,\qff v \dff}_{\dff 0}
\off =\off
0$\sss
forsss $u\qff \in\qff \kernel\dff \Gamma$\dnsp,\oss
and,\oss in\sss particular,\oss
for\sss $u\qff \in\qff \kernel\dff \gamma$\nnsp.\oss
For such $u$\sss the identity\qss (\ref{lagrange})\qss implies\sss that\sss
$\sco{\dff u\dff,\qff A\dff v \dff}_{\dff 0}
\off =\off
0$\nnsp.\oss
Since\sss $\kernel\dff \gamma$\sss is\dss dense\sss in $H_{\dff 0}$\nsp,\oss
it\sss follows\sss that\sss $A\dff v\off =\off 0$\nnsp.\oss
Since\sss $\Sigma$\sss is\dss skew-adjoint,\pss
$A\dff v\off =\off 0$\sss
together\sss with\sss the identity\qss (\ref{lagrange})\qss implies\sss that\vspace{3pt}
\[
\quad
\sco{\dff A\dff u\fff,\qff v \dff}_{\dff 0}
\qff +\qff
\sco{\dff \gamma\dff u\dff,\qff i\trf \Sigma\dff \gamma\dff v \dff}_{\dff \partial}
\off =\off
\sco{\dff A\dff u\fff,\qff v \dff}_{\dff 0}
\qff -\qff
\sco{\dff i\trf \Sigma\dff \gamma\dff u\dff,\qff \gamma\dff v \dff}_{\dff \partial}
\off =\off
0
\]

\vspace{-12pt}\vspace{3pt}
for every $u\qff \in\qff H_{\fff 1}$\nsp.\oss
Hence 
$(\trf v\fff,\qff i\trf \Sigma\dff \gamma\dff v\trf)$\sss 
belongs\sss to\sss the orthogonal\sss complement\sss of\qss 
$A\dff \oplus\dff \Gamma\trf(\trf H_{\fff 1}\trf)$
in\dss $H_{\dff 0}\dff \oplus\dff G $\nsp.\oss
By\trs Theorem\qss \ref{duality-lemma}\qss this implies\sss that
$v\qff \in\qff \kernel\dff A\dff \oplus\dff \Gamma
\off =\off
\kernel\dff A_{\trf \Gamma}$\nsp.\oss

Conversely,\oss if\dss $u\qff \in\qff \kernel\dff \Gamma$\sss and\sss
$v\qff \in\qff \kernel\dff A_{\trf \Gamma}$\dnsp,\oss
then\qss (\ref{lagrange})\qss implies\sss that\sss\vspace{3pt}
\[
\quad
\sco{\dff A\dff u\fff,\qff v \dff}_{\dff 0}
\off =\off
\sco{\dff i\trf \Sigma\dff \gamma\dff u\dff,\qff \gamma\dff v \dff}_{\dff \partial}
\off =\off
-\qff
\sco{\dff \gamma\dff u\dff,\qff i\trf \Sigma\dff \gamma\dff v \dff}_{\dff \partial}
\off.
\]

\vspace{-12pt}\vspace{3pt}
At\sss the same\sss time\sss $u\fff,\qff v\qff \in\qff \kernel\dff \Gamma$\sss
and\dss hence\sss
$\gamma\dff u\fff,\off \gamma\dff v
\off \in\off 
\kernel\dff (\trf 1\qff -\qff \Pi \trf)
\off =\off
\image\dff \Pi $\nsp.\oss
Since\sss $\Pi $\sss is\dss a self-adjoint\dss boundary condition,\oss
this implies\sss that\sss
$\Sigma\dff \gamma\dff v
\off \in\off
\kernel\dff \Pi $\nsp.\oss
Since\sss $\Pi $\sss is\dss a self-adjoint\sss projection,\pss
this,\oss in\sss turn,\pss implies\sss that\sss
$\sco{\dff \gamma\dff u\dff,\qff i\trf \Sigma\dff \gamma\dff v \dff}_{\dff \partial}
\off =\off
0$
and\dss hence\sss
$\sco{\dff A\dff u\fff,\qff v \dff}_{\dff 0}\off =\off 0$\nnsp.\oss
It\sss follows\sss that\sss every\sss
$v\qff \in\qff \kernel\dff A_{\trf \Gamma}$\sss
is\dss orthogonal\dss to\sss 
$A_{\trf \Gamma} \trf(\trf H_{\fff 1}\trf)$\nnsp.\oss
This completes\sss the proof.\oss  \eproof

\myuppar{Further\sss assumptions.}
Let\sss $A\fff,\qff \Pi$\sss be a self-adjoint\dss elliptic regular\sss boundary\sss problem.\oss 
Let\sss us assume\sss that\sss the inclusion\sss
$\iota\dff \colon\dff H_{\fff 1}\qff \ttoo\qff H_{\dff 0}$\sss
is\dss a compact\sss operator,\oss
and\dss 
that\sss the bound\-ed operator\sss $A\dff \oplus\dff \Gamma\dff \colon\dff
H_{\fff 1}\qff \ttoo\qff H_{\dff 0}\dff \oplus\dff G_{\dff 1/2}$\sss
is\dss Fredholm.\oss

\mypar{Theorem.}{fredholm}
\emph{\dnsp$A_{\trf \Gamma}$\sss defines an\dss unbounded
self-adjoint\sss operator\sss in\sss $H_{\dff 0}$\nsp.\oss
Moreover,\oss this operator\dss is\dss Fredholm,\oss
has discrete spectrum,\oss and\dss is\dss an operator\sss
with compact\sss resolvent.\oss}

\proof
Since\sss
$A\dff \oplus\dff \Gamma$\sss
is\dss Fredholm,\oss
the image of\sss
$A\dff \oplus\dff \Gamma$\sss
is\dss closed\sss subspace of\dss finite codimension\sss 
in\sss $H_{\dff 0}\dff \oplus\dff G_{\dff 1/2}$\nsp.\oss
It\sss follows\sss that\sss the image\sss
$\image\fff A_{\trf \Gamma}$\sss
is\dss closed\sss subspace of\dss finite codimension\sss in\sss $H_{\dff 0}$\nsp.\oss
By\trs Theorem\qss \ref{reduction}\qss the orthogonal\sss complement\sss of\dss
$\image\fff A_{\trf \Gamma}$\sss in $H_{\dff 0}$\sss is\dss equal\dss to\sss 
$\kernel\dff A_{\trf \Gamma}$\nsp.\oss
Equivalently,\pss
$H_{\dff 0}$\sss is\dss the orthogonal\sss direct\sss sum\sss
$\image\fff A_{\trf \Gamma}
\qff \oplus\qff
\kernel\dff A_{\trf \Gamma}$\nsp.\oss
Let\sss $\pi$\sss be\sss the orthogonal\dss projection of\dss 
$H_{\dff 0}
\qff \ttoo\qff 
\image\fff  A_{\trf \Gamma}$\nsp.\oss 
Then\sss $1\qff -\qff \pi$\sss is\dss the orthogonal\dss projection of\dss $H_{\dff 0}$
onto\sss $\kernel\fff A_{\trf \Gamma}$\nsp.\oss

Let\sss $\lambda\qff \in\qff \ccc$\nnsp,\oss
and\dss let\sss us consider\sss the operator\sss
$A_{\trf \Gamma}\qff -\qff \lambda\trf \iota
\dff \colon\dff
\kernel\fff \Gamma
\qff \ttoo\qff
H_{\dff 0}$\nsp.\oss
Let $\mathcal{H}$ be\sss the orthogonal\sss complement\sss of\dss
$\kernel\dff A_{\trf \Gamma}$ in\sss 
$\kernel\fff \Gamma\qff \subset\qff H_{\fff 1}$\nsp.\oss
Then\sss $A_{\trf \Gamma}$\sss induces an\sss invertible operator\sss
$\mathcal{H}\qff \ttoo\qff \image\fff  A_{\trf \Gamma}$\nsp.\oss
Since\sss the operator\sss $\iota\dff \colon\dff H_{\fff 1}\qff \ttoo\qff H_{\dff 0}$\sss
is\dss bounded,\oss the operator\sss\vspace{2.5pt}
\[
\quad
\pi\trf \circ\trf
\left(\trf
A_{\trf \Gamma}\qff -\qff \lambda\trf \iota
\trf\right)
\off =\off
A_{\trf \Gamma}\qff -\pff \pi\trf \circ\trf \lambda\trf \iota
\dff \colon\dff
\mathcal{H}\qff \ttoo\qff \image\fff  A_{\trf \Gamma}
\]

\vspace{-12pt}\vspace{2.5pt}
is\dss invertible\sss if\dss $\num{\lambda}$\sss is\dss sufficiently small.\oss
With respect\sss to\sss the decompositions\vspace{2.5pt}
\[
\quad
\kernel\fff \Gamma
\off =\off
\mathcal{H}\dff \oplus\dff \kernel\dff A_{\trf \Gamma}
\quad
\mbox{and}\quad
H_{\dff 0}
\off =\off
\image\fff A_{\trf \Gamma}
\qff \oplus\qff
\kernel\dff A_{\trf \Gamma}
\]

\vspace{-12pt}\vspace{2.5pt}
the operator\sss
$A_{\trf \Gamma}\qff -\qff \lambda\trf \iota$\sss
has\sss the matrix\sss form\vspace{1.5pt}
\[
\quad
\begin{pmatrix}
\off\dff \pi\trf \circ\trf
\left(\trf
A_{\trf \Gamma}\qff -\qff \lambda\trf \iota
\trf\right) &
* \off
\vspace{6pt} \\
\off\dff 0 &
-\qff \lambda\trf \iota \off 
\end{pmatrix}
\off,
\]

\vspace{-12pt}\vspace{1.5pt}
where $*$ stands for\sss the operator\sss
$\kernel\dff A_{\trf \Gamma}
\qff \ttoo\qff
\image\fff A_{\trf \Gamma}$\sss
induced\dss by\sss
$\pi\trf \circ\trf
\left(\trf
A_{\trf \Gamma}\qff -\qff \lambda\trf \iota
\trf\right)$\nnsp.\oss
The precise nature of\dss this operator\dss is\dss not\sss important.\oss
The above matrix\sss form shows\sss that\sss
$A_{\trf \Gamma}\qff -\qff \lambda\trf \iota$\sss
is\dss invertible if\dss $\num{\lambda}$\sss is\dss sufficiently small\sss
and\sss $\lambda\off \neq\off 0$\nnsp.\oss 
In\sss particular,\oss the unbounded operators\sss
$A_{\trf \Gamma}\qff -\qff \lambda\trf \iota$\sss
and\sss
$A_{\trf \Gamma}\qff -\qff \overline{\lambda}\trf \iota$\sss
have $H_{\dff 0}$ as\sss their images\dss if\dss
$\num{\lambda}$\sss is\dss sufficiently small\sss
and\sss $\lambda\off \neq\off 0$\nnsp.\oss
Since\sss $A_{\trf \Gamma}$\sss is\dss symmetric,\oss
this implies\sss that\sss $A_{\trf \Gamma}$\sss is\dss self-adjoint\sss
and,\pss in\sss particular,\pss closed.\oss
See,\oss for example,\oss \cite{s},\oss Proposition\qss 3.11.\oss
Since\sss
$\iota\dff \colon\dff H_{\fff 1}\qff \ttoo\qff H_{\dff 0}$\sss 
is\dss compact,\oss for such\sss $\lambda$\sss the operator\sss 
$\iota\dff \circ\trf (\trf A_{\trf \Gamma}\qff -\qff \lambda\trf \iota\trf)^{\dff -\dff 1}
\dff \colon\dff
H_{\dff 0}\qff \ttoo\qff H_{\dff 0}$\sss 
is\dss compact.\oss
This means\sss that\sss $A_{\trf \Gamma}\qff -\qff \lambda\trf \iota$\nnsp,\oss
considered as an\sss unbounded operator\sss in\sss $H_{\dff 0}$\nsp,\oss
has compact\sss resolvent.\oss
This\sss implies\sss that\sss $A_{\trf \Gamma}$\sss has compact\sss resolvent\sss
and\sss hence has discrete spectrum.\oss
Since\sss $A_{\trf \Gamma}$\sss is\dss closed,\oss
$\kernel\dff A_{\trf \Gamma}$\sss is\dss finitely dimensional,\pss
and\sss $\image\fff A_{\trf \Gamma}$\sss is\dss closed
and\sss has finite codimension,\pss $A_{\trf \Gamma}$\sss is\dss Fredholm.\oss \eproof

\myuppar{The continuity of\dss $A_{\trf \Gamma}$\sss
as a function of\dss $\gamma\fff,\off A\fff,\off \Sigma\fff,\off \Pi$\nnsp.}
Let\sss us fix\dss Hilbert\sss spaces\sss
$H_{\dff 0}$\nsp,\qss $H_{\fff 1}$\nsp,\qss $H^{\dff \partial}$,\oss
$H_{\dff 1/2}^{\dff \partial}$\nsp,\pss
as also\sss the inclusions\sss 
$\iota\dff \colon\dff H_{\fff 1}\qff \ttoo\qff H_{\dff 0}$\sss
and\dss
$H_{\dff 1/2}^{\dff \partial}\qff \ttoo\qff H^{\dff \partial}$,\oss
and consider\sss the space of\dss quadruples\sss
$\gamma\fff,\off A\fff,\off \Sigma\fff,\off \Pi $\sss
satisfying\dss all\dss our assumptions.\oss
Let\sss us equip\sss this space with\sss the\sss topology\sss
induced\sss by\sss the norm\sss topologies on\sss the spaces
of\dss bounded\sss operators and\sss  
projections $\Pi_{\dff 1/2}$\nsp,\oss 
and\sss the space of\dss unbounded self-adjoint\sss operators
in\sss $H_{\dff 0}$\sss with\sss the\sss topology of\dss the convergence
in\sss the norm\sss resolvent\sss sense.\oss
See\qss \cite{rs},\oss Section\qss VIII.7,\oss for\sss the\sss latter.\oss

\mypar{Lemma.}{continuity}
\emph{$A_{\trf \Gamma}$ continuously depends on\dss
$\gamma\fff,\off A\fff,\off \Sigma\fff,\off \Pi$\nnsp.\oss
If\trs $A\dff \oplus\dff \Gamma$\dss is\dss invertible,\oss
then\sss $A_{\trf \Gamma}$\sss is\dss invertible as an\sss operator\sss
$\kernel\fff \Gamma \qff \ttoo\qff H_{\dff 0}$\sss
and\sss the inverse\sss 
$A_{\trf \Gamma}^{\dff -\dff 1}\dff \colon\dff H_{\dff 0}\qff \ttoo\qff H_{\dff 1}$\sss
continuously depends on\sss the quadruple.\oss}

\proof
Let\sss us\sss fix\sss some quadruple\sss
$\bm{\gamma}\fff,\off \mathbf{A}\fff,\off \bm{\Sigma}\fff,\off \bm{\Pi }$\sss
and use\sss the boldface font\sss for\sss objects related\sss
to\sss this quadruple.\oss
Let\sss consider\sss quadruples\sss
$\gamma\fff,\off A\fff,\off \Sigma\fff,\off \Pi $\sss
close\sss to\sss this quadruple.\oss
Let\vspace{3pt}
\[
\quad
\widetilde{A\dff \oplus\dff \Gamma}\qff \colon\qff
H_{\fff 1}\dff \oplus\dff \image\fff \bm{\Pi_{\dff 1/2}}
\qff \ttoo\qff 
H_{\dff 0}\dff \oplus\dff H_{\dff 1/2}^{\dff \partial}
\]

\vspace{-12pt}\vspace{3pt}
be\sss the bounded operator defined\sss by\sss the rule\sss 
$\widetilde{A\dff \oplus\dff \Gamma}\qff \colon\qff
(\trf h\fff,\qff w\trf)
\off \longmapsto\off
\left(\qff 
A\dff h\fff,\off 
\Gamma\dff h 
\pff +\pff
w
\qff\right)$\nnsp.\oss
The operator
$\widetilde{\mathbf{A}\dff \oplus\dff \bm{\Gamma}}$\dss
is\dss Fredholm,\oss
being\sss the direct\sss sum of\dss the\dss Fredholm\sss operator\sss
$\mathbf{A}\dff \oplus\dff \bm{\Gamma}$\sss
and\sss the identity\sss 
$\image\fff \bm{\Pi_{\dff 1/2}}\qff \ttoo\qff \image\fff \bm{\Pi_{\dff 1/2}}$.\oss
Hence\sss the operator\dss
$\widetilde{A\dff \oplus\dff \Gamma}$\dss
is\trs Fredholm\dss for quadruples close\sss to\sss the fixed one.\oss
Suppose\sss first\sss that\sss 
$\mathbf{A}\dff \oplus\dff \bm{\Gamma}$\sss is\dss
invertible.\oss
Then\sss
$\widetilde{\mathbf{A}\dff \oplus\dff \bm{\Gamma}}$\sss
is\dss also invertible,\oss
and\sss hence\sss
$\widetilde{A\dff \oplus\dff \Gamma}$\sss
is\dss invertible for quadruples close\sss to\sss the fixed one.\oss
The inverse of\dss the\sss latter operator\sss continuously depends on\sss
the quadruple.\oss
If\trs 
$\Gamma\dff h\off =\off 0$\nnsp,\oss then\vspace{3pt}
\[
\quad
\left(\qff \widetilde{A\dff \oplus\dff \Gamma}\qff\right)^{\dff -\dff 1}\qff
(\trf 
A\dff h\fff,\qff 0
\trf)
\off =\off
\left(\qff \widetilde{A\dff \oplus\dff \Gamma}\qff\right)^{\dff -\dff 1}\qff
(\trf 
A\dff h\fff,\qff \Gamma\dff h\qff +\qff 0
\trf)
\off =\off
(\trf h\fff,\qff 0\trf)
\]

\vspace{-12pt}\vspace{3pt}
and\dss hence\sss the operator 
$A_{\trf \Gamma}\dff \colon\dff
\kernel\fff \Gamma\qff \ttoo\qff H_{\dff 0}$ 
is\dss invertible and\vspace{3pt}
\[
\quad
\left(\qff A_{\trf \Gamma}^{\dff -\dff 1}\dff(\trf u\trf)\fff,\qff 0\qff\right)
\off =\off\qff
\left(\qff \widetilde{A\dff \oplus\dff \Gamma}\qff\right)^{\dff -\dff 1}\dff
(\dff u\fff,\qff 0\trf)
\]

\vspace{-12pt}\vspace{3pt}
for every $u\qff \in\qff H_{\dff 0}$\dnsp.\oss
Since
$\dnsp\left(\qff \widetilde{A\dff \oplus\dff \Gamma}\qff\right)^{\dff -\dff 1}\dnsp$
continuously depends on\sss the quadruple,\oss
we see\sss that\vspace{3pt}
\[
\quad
A_{\trf \Gamma}^{\dff -\dff 1}\dff \colon\dff
H_{\dff 0}
\qff \ttoo\qff
\kernel\dff \Gamma
\qff \subset\qff
H_{\fff 1}
\]

\vspace{-12pt}\vspace{3pt}
continuously depends on\sss the quadruple.\oss
This proves\sss the second claim of\dss the\sss lemma.\oss
By\sss taking\sss the composition of\dss these inverses\sss with\sss the inclusion\sss
$\iota\dff \colon\dff H_{\fff 1}\qff \ttoo\qff H_{\dff 0}$\sss
we conclude\sss that\sss 
$A_{\trf \Gamma}^{\dff -\dff 1}\dff \colon\dff
H_{\dff 0}
\qff \ttoo\qff
H_{\dff 0}$\sss
continuously depends on\sss the quadruple.\oss
But\sss this means exactly\sss that\sss $A_{\trf \Gamma}$\sss
continuously depends on\sss the quadruple in\sss the\sss topology of\dss the convergence
in\sss the norm\sss resolvent\sss sense.\oss
This proves\sss the\sss lemma  
when\sss
$\mathbf{A}\dff \oplus\dff \bm{\Gamma}$\sss is\dss
invertible.\oss
In\sss general,\qss
$(\trf \mathbf{A}\qff -\qff \lambda\trf)\dff \oplus\dff \bm{\Gamma}$\sss
is\dss invertible for some $\lambda$\nnsp.\oss
By\sss the already\sss proved special\sss case\sss
$(\trf A\qff -\qff \lambda\trf)_{\dff \Gamma}$\sss
continuously depends on\sss the quadruple for\sss 
quadruples close\sss to\sss the fixed one.\oss
But\sss
$(\trf A\qff -\qff \lambda\trf)_{\dff \Gamma}$\sss
is\dss equal\dss to\sss
$A_{\trf \Gamma}\qff -\qff \lambda$\nnsp,\oss
and\sss hence\sss $A_{\trf \Gamma}$\sss 
also continuously depends on\sss the quadruple\qss
(see\qss \cite{k},\oss Theorem\qss IV.2.23).\oss  \eproof

\myuppar{Families of\dss abstract\sss self-adjoint\sss boundary\sss problems.}
Let\sss $Z$\sss be a\sss topological\sss space.\oss
Suppose\sss that\sss for every\sss $z\qff \in\qff Z$\sss we are given a quadruple\sss
$\gamma\trf(\trf z\trf)\dff,\off 
A\trf(\trf z\trf)\dff,\off 
\Sigma\trf(\trf z\trf)\dff,\off 
\Pi \trf(\trf z\trf)$\sss
satisfying\sss all\sss our assumptions and continuously depending on\sss $x$\nnsp.\oss
Then for every $z\qff \in\qff Z$\sss a bounded operator
$\Gamma\trf(\trf z\trf)\dff \colon\dff 
H^{\fff 1}
\qff \ttoo\qff 
G_{\dff 1/2}\trf(\trf z\trf)
\off =\off 
\kernel\dff \Pi_{\dff 1/2}\trf(\trf z\trf)$\sss 
and\sss an unbounded\sss operator\vspace{2pt}
\[
\quad
A_{\trf \Gamma}\trf(\trf z\trf)
\off\dff =\off
A\trf(\trf z\trf)_{\qff \Gamma\trf(\trf z\trf)}
\]

\vspace{-12pt}\vspace{2pt}
in\sss $H_{\dff 0}$\sss are defined.\oss
By\trs Lemma\qss \ref{continuity}\qss the map\sss
$z\off \longmapsto\off A_{\trf \Gamma}\trf(\trf z\trf)$\sss
is\dss continuous with respect\sss the\sss topology of\dss $Z$\sss 
and\sss the\sss topology of\dss the convergence
in\sss the norm\sss resolvent\sss sense on unbounded 
self-adjoint\sss operators in\sss $H_{\dff 0}$\nsp.\oss
It\sss follows\sss that\sss if\dss two real\sss numbers\sss $a\qff <\qff b$\sss
do not\sss belong\sss to\sss the spectrum of\dss $A_{\trf \Gamma}\trf(\trf z_{\trf 0}\trf)$\sss
for some\sss $z_{\trf 0}\qff \in\qff Z$\nnsp,\oss then\sss $a\fff,\qff b$\sss
do not\sss belong\sss to\sss the spectrum of\dss $A_{\trf \Gamma}\trf(\trf z\trf)$\sss
for\sss $z\qff \in\qff Z$\sss sufficiently close\sss to\sss $z_{\trf 0}$\nsp,\oss
and\sss the spectral\sss projections\vspace{3pt}
\[
\quad
P_{\dff [\dff a\fff,\qff b\dff]}\qff \bigl(\trf A_{\trf \Gamma}\trf(\trf z\trf)\qff \bigr)
\]

\vspace{-12pt}\vspace{3pt}
continuously depend on $z$ for such $z$\nnsp.\oss
See\qss \cite{rs},\oss Theorem\qss VIII.23.\oss
This\sss immediately\sss implies\sss that\sss
$A_{\trf \Gamma}\trf(\trf z\trf)\dff,\off z\qff \in\qff Z$\sss
is\dss a\dss Fredholm\sss family of\dss self-adjoint\sss operators in\sss
the sense of\qss \cite{i2}\qss and\sss hence\sss the\qss \emph{analytical\dss index}\qss of\dss
this family\dss is\dss defined.\oss

\newpage
\mysection{Multiplicative\qss properties\qss of\pss abstract\qss operators}{mult-operators}

\myuppar{Tensor products of\trs Hilbert\dss spaces and operators.}
Let $H\fff,\off K$ be\sss two\dss Hilbert\sss spaces.\oss
The algebraic\sss tensor product\sss $H\dff \otimes\dff K$\sss admits
a unique scalar product\sss $\sco{\trf \bullet\fff,\qff \bullet\trf}$\sss
such\sss that\vspace{3pt}
\[
\quad
\sco{\dff h\dff \otimes\dff k\fff,\qff h\fff'\dff \otimes\dff k\fff' \trf}
\off =\off
\sco{\dff h\dff,\qff h\fff' \dff}\dff
\sco{\dff k\dff,\qff k\fff' \dff}
\pff
\]

\vspace{-12pt}\vspace{3pt}
where\sss $h\fff,\qff h\fff'\qff \in\qff H$\nnsp,\qss 
$k\fff,\qff k\fff'\qff \in\qff K$\nnsp,\oss and\dss the scalar products on\sss the right\sss
are\sss taken\sss in $H$ and\dss $K$\sss respectively.\oss
The\dss Hilbert\dss space\sss tensor product\sss
$H\qff\fff \widehat{\otimes}\pff K$\sss
is\dss defined as\sss the completion of\dss $H\dff \otimes\dff K$\sss
with\sss this scalar product.\oss
For\sss bounded operators\sss
$A\dff \colon\dff H\qff \ttoo\qff H\fff'$\sss
and\sss
$B\dff \colon\dff K\qff \ttoo\qff K\fff'$\sss
the algebraic\sss tensor product\sss
$A\dff \otimes\dff B\dff \colon\dff 
H\dff \otimes\dff K
\qff \ttoo\qff 
H\fff'\dff \otimes\dff K\fff'$\sss
extends by continuity\sss to a bounded operator\vspace{3pt}
\[
\quad
A\qff\fff \widehat{\otimes}\pff B\dff \colon\dff 
H\qff\fff \widehat{\otimes}\pff K
\qff \ttoo\qff 
H\fff'\qff\fff \widehat{\otimes}\pff K\fff'
\pff.
\]

\vspace{-12pt}\vspace{3pt}
If\dss $A\fff,\qff B$\sss are bounded and\dss injective,\oss
then\sss $A\qff\fff \widehat{\otimes}\pff B$\sss is\dss also injective.\oss
See\qss \cite{pa},\oss Section\qss XIV.1.\oss
We are mostly concerned\sss with unbounded operators.\oss
For an\sss unbounded operator\sss
$A$ from $H$\sss to $H\fff'$
we will\sss denote by\sss $\mathcal{D}\dff(\trf A\trf)$\sss its domain,\oss
so $\mathcal{D}\dff(\trf A\trf)$\sss is\dss a\sss linear subspace of\dss $H$\sss
and $A$\sss is\dss a\sss linear map\sss
$\mathcal{D}\dff(\trf A\trf)\qff \ttoo\qff H\fff'$\nnsp.\oss
By an abuse of\dss notation we may speak in\sss this situation about\sss an
an unbounded operator\sss $A\dff \colon\dff H\qff \ttoo\qff H\fff'$\sss
with\sss the domain\sss $\mathcal{D}\dff(\trf A\trf)$\nnsp.\oss
For unbounded operators
$A\dff \colon\dff H\qff \ttoo\qff H\fff'$\sss
and\sss
$B\dff \colon\dff K\qff \ttoo\qff K\fff'$\sss
the algebraic\sss tensor product\vspace{3pt}
\[
\quad
A\dff \otimes\dff B\dff \colon\dff 
H\qff\fff \widehat{\otimes}\pff K
\qff \ttoo\qff 
H\fff'\qff\fff \widehat{\otimes}\pff K\fff'
\pff.
\]

\vspace{-12pt}\vspace{3pt}
is\dss an unbounded operator\sss with\sss the domain\sss
$\mathcal{D}\dff(\trf A\trf)\dff \otimes\trf \mathcal{D}\dff(\trf B\trf)$\nnsp.\oss
If\dss $A$ and\sss $B$ are densely defined,\oss i.e.\sss if\dss
$\mathcal{D}\dff(\trf A\trf)$ and\sss $\mathcal{D}\dff(\trf B\trf)$
are dense in $H$ and\sss $K$ respectively,\oss
then\sss $A\dff \otimes\dff B$\sss is\dss also densely defined.\oss
Moreover,\oss if\dss $A$ and\sss $B$ are closable,\oss
then\sss $A\dff \otimes\dff B$\sss is\dss closable.\oss
See\qss \cite{s},\oss Lemma\qss 7.21.\oss
In\sss this case\sss the closure of\dss $A\dff \otimes\dff B$\sss
is\dss denoted\dss by\vspace{3pt}
\[
\quad
A\qff\fff \widehat{\otimes}\pff B\dff \colon\dff 
H\qff\fff \widehat{\otimes}\pff K
\qff \ttoo\qff 
H\fff'\qff\fff \widehat{\otimes}\pff K\fff'
\pff
\]

\vspace{-12pt}\vspace{3pt}
and called\dss the\qss \emph{tensor\dss product}\pss of\dss $A$ and\sss $B$\nnsp.\oss
One should\dss keep in\sss mind\dss that\sss the domain\sss
$\mathcal{D}\dff(\trf A\qff\fff \widehat{\otimes}\pff B\trf)$
is\dss not\sss specified directly and\dss
its\sss determination may\sss be a non-trivial\dss task.\oss
If\dss $A\fff,\pff B$ are self-adjoint,\oss then\sss
$A\qff\fff \widehat{\otimes}\pff B$\sss is\dss also self-adjoint,\oss
and\dss if\dss $A\fff,\pff B$ are densely defined and closable,\oss then 
$(\trf A\qff\fff \widehat{\otimes}\pff B\trf)^{\fff *}
\qff =\off
A^{\fff *}\qff \widehat{\otimes}\pff\fff B^{\fff *}$\dnsp.\oss
See\qss \cite{s},\oss Theorem\qss 7.23\qss and\qss Proposition\qss 7.26.

\myuppar{The framework.}
Let\sss $H_{\trf 0}\dff,\off K_{\trf 0}\dff,\off L_{\trf 0}$\sss 
be\sss separable\dss Hilbert\dss spaces.\oss
Let\sss
$P\dff \colon\dff H_{\dff 0}\qff \ttoo\qff H_{\dff 0}$\sss
be an\sss unbounded\sss densely defined operator\sss with\sss the domain
$H_{\dff 1}\off =\off \mathcal{D}\dff(\trf P\trf)$\nnsp.\oss
Let\sss
$Q\dff \colon\dff K_{\trf 0}\qff \ttoo\qff L_{\trf 0}$\sss
and\sss
$Q^{\fff *}\fff \colon\dff L_{\trf 0}\trf \ttoo\trf K_{\trf 0}$\sss 
be closed densely defined operators with\sss the domains
$K_{\dff 1}\pff =\off \mathcal{D}\dff(\trf Q\trf)$
and\sss
$L_{\dff 1}\off =\off \mathcal{D}\dff(\trf Q^{\fff *}\trf)$\nnsp.\oss
In an agreement\sss with\sss the notations,\oss
we assume\sss that\sss $Q$ and $Q^{\fff *}$
are adjoint\sss to each other as unbounded operators.\oss
In\sss particular,\pss $L_{\dff 1}$\sss is\dss the domain of\dss
the operator adjoint\sss to $Q$ and\sss $K_{\dff 1}$\sss is\dss the domain of\dss
the operator adjoint\sss to $Q^{\fff *}$\dnsp.\oss
We will\sss assume\sss that\sss the vector spaces\sss
$H_{\dff 1}\dff,\off K_{\dff 1}\dff,\off L_{\dff 1}$\sss 
are\dss Hilbert\dss spaces in\sss their own right\sss
such\sss that\sss the operators\vspace{3pt}
\[
\quad
P\dff \colon\dff 
H_{\dff 1}\qff \ttoo\qff H_{\trf 0}\dff,\quad 
Q\dff \colon\dff 
K_{\dff 1}\qff \ttoo\qff L_{\trf 0}\dff,\quad
\mbox{and}\quad
Q^{\fff *}\dff \colon\dff 
L_{\dff 1}\qff \ttoo\qff K_{\trf 0}
\]

\vspace{-12pt}\vspace{3pt}
are bounded\sss and\sss the inclusions\sss\vspace{3pt}
\[
\quad
H_{\dff 1}\qff \ttoo\qff H_{\trf 0}\dff,\quad 
K_{\dff 1}\qff \ttoo\qff K_{\trf 0}\dff,\quad
\mbox{and}\quad
L_{\dff 1}\qff \ttoo\qff L_{\trf 0}
\]

\vspace{-12pt}\vspace{3pt}
are bounded and compact.\oss
We will\sss assume\sss that\sss 
$Q\dff\halfff Q^{\fff *}$ and\sss $Q^{\fff *}\fff Q$ are\dss 
Fredholm\dss and\dss self-adjoint\dss unbounded operators with compact\sss resolvent.\oss
This implies\sss that\sss $Q$ and\sss $Q^{\fff *}$ are\dss Fredholm\dss
as unbounded,\oss and\dss hence as bounded,\oss operators.\oss
See\qss \cite{sc},\oss Theorem\qss 7.12.\oss

Let\sss $D_{\dff 1}$ be a closed subspace of\dss $H_{\dff 1}$\sss
dense in\sss $H_{\trf 0}$\sss
and\dss let\sss
$A\dff \colon\dff D_{\dff 1}\qff \ttoo\qff H_{\trf 0}$\sss
be\sss the restriction of\dss $P$\sss to $D_{\dff 1}$\nsp.\oss
The operator $A$ corresponds\sss to\sss the operator\sss
$A_{\trf \Gamma}$\sss from\dss Section\qss \ref{abstract-index},\oss
and\sss the subspace 
$D_{\dff 1}$ corresponds\sss to\sss $\kernel\fff \Gamma$\dnsp.\oss
We will\sss assume\sss that\sss $A$\sss is\trs Fredholm\dss and\dss
is\dss closed and\sss self-adjoint\sss as an unbounded\sss operator\sss 
$H_{\trf 0}\qff \ttoo\qff H_{\trf 0}$\nsp.\oss
Then $A^{\fff 2}$ also has\sss these properties.\oss
Since\sss the inclusion\sss
$H_{\dff 1}\qff \ttoo\qff H_{\trf 0}$\sss
is\dss compact,\pss $A$ 
is\dss an operator with compact\sss resolvent. 

\myuppar{The $\omult_{\fff 1}$ products.}
Let\sss us define a product\sss $A\pff \omult_{\fff 1}\qff Q$
matching\dss the product\sss $\sigma\pff \omult_{\fff 1}\qff q$ of\dss symbols
when\sss $H_{\dff 0}\dff,\pff H_{\dff 1}\dff,\pff K_{\trf 0}$\nsp,\qss etc.\oss
are\dss Sobolev\dss spaces of\dss sections of\dss bundles.\oss
The matrix\vspace{3pt}
\begin{equation}
\label{aq-matrix}
\quad
\begin{pmatrix}
\off\dff A\dff \otimes\dff 1 &
1\dff \otimes\dff Q^{\fff *} \off
\vspace{6pt} \\
\off\dff 1\dff \otimes\dff Q &
-\qff A\dff \otimes\dff 1 \off 
\end{pmatrix}
\off
\end{equation}

\vspace{-12pt}\vspace{3pt}
defines an\sss unbounded operator\sss in\sss
$(\trf H_{\dff 0}\qff\fff \widehat{\otimes}\pff K_{\trf 0} \trf)
\qff \oplus\qff
(\trf H_{\dff 0}\qff\fff \widehat{\otimes}\pff L_{\trf 0} \trf)$\sss
with\sss the domain\vspace{3pt}
\[
\quad
(\trf D_{\dff 1}\dff \otimes\dff K_{\dff 1}\trf)
\qff \oplus\qff
(\trf D_{\dff 1}\dff \otimes\dff L_{\dff 1} \trf)
\]

\vspace{-12pt}\vspace{3pt}
Clearly,\oss this operator\dss is\dss densely defined and symmetric,\oss
and\dss hence\dss is\dss closable.\oss
The product\sss $A\pff \omult_{\fff 1}\qff Q$\sss is\dss defined as\sss its closure.\oss
The square of\dss the matrix\qss (\ref{aq-matrix})\qss
is\dss equal\dss to\vspace{3pt}
\[
\quad
\begin{pmatrix}
\off\dff A^{\dff 2}\dff \otimes\dff 1\qff +\qff 1\dff \otimes\dff Q^{\fff *}\fff Q &
0 \off
\vspace{6pt} \\
\off\dff 0 &
A^{\dff 2}\dff \otimes\dff 1\qff +\qff 1\dff \otimes\dff Q\trf Q^{\fff *} \off 
\end{pmatrix}
\off
\]

\vspace{-12pt}\vspace{3pt}
and defines an unbounded operator between\sss the same\dss
Hilbert\dss spaces as\sss $A\pff \omult_{\fff 1}\qff Q$\nnsp.\oss
This operator\dss is\dss also densely defined and symmetric,\oss
and\dss hence\dss is\dss closable.\oss
Its closure\dss is\dss equal\dss to\sss the direct\sss sum of\dss
closures of\dss unbounded operators\vspace{3pt}
\[
\quad
A^{\dff 2}\dff \otimes\dff 1\qff +\qff 1\dff \otimes\dff Q^{\fff *}\fff Q
\qff \colon\qff
H_{\trf 0}\qff\fff \widehat{\otimes}\pff K_{\trf 0}
\qff \ttoo\qff
H_{\dff 0}\qff\fff \widehat{\otimes}\pff K_{\trf 0}
\qquad
\mbox{and}\quad
\]

\vspace{-34.5pt}
\[
\quad
A^{\dff 2}\dff \otimes\dff 1\qff +\qff 1\dff \otimes\dff Q\dff\halfff Q^{\fff *}
\qff \colon\qff
H_{\trf 0}\qff\fff \widehat{\otimes}\pff L_{\trf 0}
\qff \ttoo\qff
H_{\dff 0}\qff\fff \widehat{\otimes}\pff L_{\trf 0}
\off.
\]

\vspace{-12pt}\vspace{3pt}
\myuppar{Eigenvectors.}
Since $A$\sss is\dss a self-adjoint\sss operator with compact\sss resolvent,\oss
there\dss exists an
orthonormal\dss basis of\dss $H_{\dff 0}$ consisting of\dss
eigenvectors $u_{\dff i}\qff \in\qff D_{\dff 1}$ of\sss $A$\nnsp.\oss
Here $i$ runs over natural\sss numbers,\oss say.\pss
Moreover,\oss the corresponding eigenvalues $\lambda_{\dff i}$\sss tend\sss to plus or minus
infinity\sss when\sss $i\trf \ttoo\qff \infty$\nnsp.\oss
The eigenvectors of\sss $A^{\fff 2}$ are\sss the same,\oss
with eigenvalues\sss $\lambda_{\dff i}^{\fff 2}$\nsp.\oss

Similarly,\oss 
there\dss exists an orthonormal\dss basis of\dss $K_{\trf 0}$ consisting of\dss
eigenvectors $v_a\qff \in\qff K_{\dff 1}$ of\sss $Q^{\fff *}\fff Q$\nnsp.\oss 
Since\sss $Q^{\fff *}\fff Q$\sss is\dss positive,\oss the corresponding eigenvalues are\sss $\geq\qff 0$\sss
and\sss we can write\sss them\sss in\sss the form\sss $\mu_{\fff a}^{\dff 2}$
with\sss $\mu_{\fff a}\qff \geq\qff 0$\nnsp.\oss
Also,\oss there\dss exists an orthonormal\dss basis of\dss $L_{\trf 0}$ consisting of\dss
eigenvectors $w_{\fff b}\qff \in\qff L_{\dff 1}$ of\sss $Q\fff Q^{\fff *}$\sss
with\sss the eigenvalues\sss $\nu_{b}^{\dff 2}$\nsp,\oss where $\nu_{\fff b}\qff \geq\qff 0$\nnsp.\oss
Moreover,\pss $\mu_{\fff a}$\sss tends\sss to
infinity\sss when\sss $a\qff \ttoo\qff \infty$\nnsp,\oss 
and\sss $\nu_{b}$ have\sss the same property.\oss

\myuppar{Further\sss assumptions.}
Let $\mathcal{K}_{\dff 1}$\sss be\sss the closure\sss in $K_{\dff 1}$
of\dss the\sss linear span of\dss vectors $v_a$ with\sss
$\mu_{\fff a}\off \neq\off 0$\nnsp.\oss
We will\sss assume\sss that\sss the identity\sss map\sss 
$\mathcal{K}_{\dff 1}\qff \ttoo\qff \mathcal{K}_{\dff 1}$\sss
is\dss a\sss topological\dss isomorphism\sss between\sss
the structure of\dss a\dss Hilbert\dss space on\sss $\mathcal{K}_{\dff 1}$\sss
induced\sss from\sss $K_{\dff 1}$\sss and\sss the structure having\sss
the vectors\sss $\mu_{\fff a}^{\dff -\dff 1}\dff v_a$\sss with\sss
$\mu_{\fff a}\off \neq\off 0$\sss
as an orthonormal\sss basis.\oss
We define\sss the space $\mathcal{L}_{\dff 1}$\sss similarly
and assume\sss that\sss it\dss has\sss similar property.\oss

While\sss these assumptions appear\sss to be somewhat\sss artificial,\oss
they hold\sss when\sss $K_{\trf 0}\dff,\pff K_{\dff 1}$ and\sss $L_{\trf 0}\dff,\off L_{\dff 1}$ 
are\sss Sobolev\dss spaces of\dss
sections of\dss vector bundles on a closed\sss manifold and\sss $Q$\sss
is\dss induced\dss by an elliptic operator.\oss
See\qss \cite{pa},\oss Theorem\qss XI.14,\oss
applied\dss to\sss $S\off =\off Q^{\fff *}\fff Q$ or\sss $Q\dff\halfff Q^{\fff *}$\dnsp.\oss

\myuppar{Bounded and\sss analytic vectors.}
Let\sss $T\dff \colon\dff H\qff \ttoo\qff H$\sss be an unbounded operator
in a\dss Hilbert\sss space $H$\nnsp.\oss
Let\sss us denote\sss by\sss $\norm{\bullet}$\sss the norm\sss in\sss $H$\nnsp.\oss
Recall\dss that\sss a vector\sss $x\qff \in\qff H$\sss is\dss said\dss to be a\qss
\emph{bounded\sss vector}\pss for $T$\sss if\dss $x$\sss belongs\sss to
the domain\sss of\dss $T^{\dff n}$\sss for every $n\qff \geq\qff 1$\sss
and\sss $\norm{T^{\dff n}\dff x}\qff \leq\qff C^{\dff n}$\sss for every $n$
and some constant\sss $C$\dss (depending on $x$\nsp).\oss
It\dss is\dss said\dss to be an\qss \emph{analytic vector}\pss
if\dss a weaker estimate of\dss the form\sss 
$\norm{T^{\dff n}\dff x}\qff \leq\qff C^{\dff n}\dff n\dff !$\trs 
holds.\oss
See\qss \cite{s},\oss Section\qss 7.4.\oss

\mypar{Lemma.}{bounded-vectors}
\emph{Every vector of\dss the form
$(\dff u_{\dff i}\dff \otimes\dff v_a\trf)
\dff \oplus\dff
(\dff u_{\dff k}\dff \otimes\dff w_{\fff b}\trf)$
is\dss a bounded\sss vector\sss for\sss
$A\pff \omult_{\fff 1}\qff Q$\nnsp.}

\proof
Let\sss us denote by\sss $\norm{\bullet}_{\trf 0}$\sss the norm\sss in\dss
Hilbert\dss spaces involving only subscripts $0$\nnsp,\oss such as\sss
$H_{\trf 0}\qff\fff \widehat{\otimes}\pff K_{\trf 0}$\nsp,\oss
and\dss by\sss $\norm{\bullet}_{\dff 1}$\sss the norms in\sss
$H_{\dff 1}\dff,\off K_{\dff 1}\dff,\off L_{\dff 1}$\nsp.\oss
Clearly,\oss every\sss  
$(\trf u_{\dff i}\dff \otimes\dff v_a\trf)
\dff \oplus\dff
(\trf u_{\dff k}\dff \otimes\dff w_{\fff b}\trf)$\sss
belongs\sss to\sss the domains of\dss 
even powers of\dss $A\pff \omult_{\fff 1}\qff Q$\dnsp,\oss
and\vspace{3pt}
\[
\quad
(\trf A\pff \omult_{\fff 1}\qff Q\trf)^{\dff 2\fff n}
\qff
\left(\qff
(\trf u_{\dff i}\dff \otimes\dff v_a\trf)
\off \oplus\off
(\trf u_{\dff k}\dff \otimes\dff w_{\fff b}\trf)
\qff\right)
\]

\vspace{-33pt}
\[
\quad
=\off
\bigl(\trf \lambda_{\dff i}^{\dff 2} 
\off +\off
\mu_{\fff a}^{\dff 2} \qff\bigr)^{\dff n}
\qff
\bigl(\trf u_{\dff i}\dff \otimes\dff v_a \trf\bigr)
\off \oplus\off
\bigl(\trf \lambda_{\dff k}^{\dff 2} 
\off +\off
\nu_{\fff b}^{\dff 2} \qff\bigr)^{\dff n}\qff
\bigl(\trf u_{\dff k}\dff \otimes\dff w_{\fff b} \trf\bigr)
\off.
\]

\vspace{-12pt}\vspace{1.5pt}
Therefore\vspace{1.5pt}
\[
\quad
\bnorm{(\trf A\pff \omult_{\fff 1}\qff Q\trf)^{\dff 2\fff n}
\qff
\left(\qff
(\trf u_{\dff i}\dff \otimes\dff v_a\trf)
\off \oplus\off
(\trf u_{\dff k}\dff \otimes\dff w_{\fff b}\trf)
\qff\right)}_{\qff 0}^{\qff 2}
\]

\vspace{-33pt}
\[
\quad
=\off
\bigl(\trf \lambda_{\dff i}^{\dff 2} 
\off +\off
\mu_{\fff a}^{\dff 2} \qff\bigr)^{\dff 2\fff n}
\off +\off\qff
\bigl(\trf \lambda_{\dff k}^{\dff 2} 
\off +\off
\nu_{\fff b}^{\dff 2} \qff\bigr)^{\dff 2\fff n}
\pff.
\]

\vspace{-12pt}\vspace{3pt}
Our assumptions about\sss $\mathcal{K}_{\dff 1}$ and $\mathcal{L}_{\dff 1}$\sss
imply\sss that\dss if\dss 
$\mu_{\fff a}\dff,\pff \nu_{\fff b}\off \neq\off 0$\nnsp,\oss
then\vspace{4.5pt}
\[
\quad
\norm{v_a}_{\dff 1}
\off =\off
\mu_{\fff a}\qff \norm{\mu_{\fff a}^{\dff -\dff 1}\dff v_a}_{\dff 1}
\off \leq\off
C\fff'\qff \mu_{\fff a}
\quad
\mbox{and}\quad
\norm{w_{\fff b}}_{\dff 1}
\off \leq\off
C\fff'\qff \nu_{\fff b}
\]

\vspace{-12pt}\vspace{4.5pt}
for some constant $C\fff'$ independent\sss from $a\fff,\qff b$\nnsp.\oss
Since $Q$\sss and $Q^{\fff *}$ are bounded\sss as operators\sss 
$K_{\dff 1}\qff \ttoo\qff L_{\trf 0}$
and\sss
$L_{\dff 1}\qff \ttoo\qff K_{\trf 0}$,\oss
it\sss follows\sss that\vspace{4.5pt}
\[
\quad
\norm{Q\dff(\trf v_a\trf)}_{\trf 0}
\qff \leq\qff
C\fff''\qff \norm{v_a}_{\dff 1}
\off \leq\off
C\fff''\trf C\fff'\qff \mu_{\fff a}
\quad
\mbox{and}\quad
\norm{Q^{\fff *}\dff(\trf w_{\dff b}\trf)}_{\trf 0}
\off \leq\off
C\fff''\trf C\fff'\qff \nu_{\dff b}
\]

\vspace{-12pt}\vspace{4.5pt}
if\dss 
$\mu_{\fff a}\dff,\pff \nu_{\fff b}\off \neq\off 0$\sss
for some constant\sss $C\fff''$ independent\sss from $a\fff,\qff b$\nnsp.\oss
Let\sss
$C\off =\off C\fff''\trf C\fff'\qff +\qff 1$\nnsp.\oss
Then\vspace{4.5pt}
\[
\quad
\norm{Q\dff(\trf v_a\trf)}_{\trf 0}
\qff \leq\qff
C\qff \mu_{\fff a}
\quad
\mbox{and}\quad
\norm{Q^{\fff *}\dff(\trf w_{\dff b}\trf)}_{\trf 0}
\qff \leq\qff
C\qff \nu_{\dff b}
\off.
\]

\vspace{-12pt}\vspace{4.5pt}
for every $a\fff,\qff b$\nnsp.\oss
Let\dss $r\fff,\qff s\qff \in\qff \rrr$\nnsp.\oss
Then\vspace{4.5pt}
\[
\quad
(\trf A\pff \omult_{\fff 1}\qff Q\trf)
\qff
\left(\qff
r\qff (\trf u_{\dff i}\dff \otimes\dff v_a\trf)
\pff \oplus\pff
s\qff (\trf u_{\dff k}\dff \otimes\dff w_{\dff b}\trf)
\qff\right)
\pff
\]

\vspace{-31.5pt}
\[
\quad
=\off
\left(\qff
r\qff \lambda_{\dff i}\qff (\trf u_{\dff i}\dff \otimes\dff v_a\trf)
\qff +\qff
s\qff (\trf u_{\dff k}\dff \otimes\dff Q^{\fff *}\dff(\trf w_{\fff b}\trf)\trf)
\qff\right)
\pff \oplus\pff
\left(\qff
r\qff (\trf u_{\dff i}\dff \otimes\dff Q\dff(\trf v_a\trf)\trf)
\qff -\qff
s\trf \lambda_{\dff k}\qff (\trf u_{\dff k}\dff \otimes\dff w_{\dff b} \trf)\trf)
\qff\right)
\pff.
\]

\vspace{-12pt}\vspace{4.5pt}
The square of\dss the norm of\dss this vector\sss in\sss
$(\trf H_{\dff 0}\qff\fff \widehat{\otimes}\pff K_{\trf 0} \trf)
\qff \oplus\qff
(\trf H_{\dff 0}\qff\fff \widehat{\otimes}\pff L_{\trf 0} \trf)$\dss
is\vspace{4.5pt}
\[
\quad
\leq\off
2\trf (\trf r\qff \lambda_{\dff i}\trf)^{\dff 2}
\qff +\qff
2\trf s^{\dff 2}\trf \norm{Q^{\fff *}\dff(\trf w_{\fff b}\trf)}_{\trf 0}^{\dff 2}
\qff +\qff
2\trf r^{\dff 2}\trf \norm{ Q\dff(\trf v_a\trf)}_{\trf 0}^{\dff 2}
\qff +\qff
2\trf (\trf s\qff \lambda_{\dff k}\trf)^{\dff 2}
\]

\vspace{-30pt}
\[
\quad
\leq\off
2\trf r^{\dff 2}\qff \lambda_{\dff i}^{\dff 2}
\qff +\qff
2\trf s^{\dff 2}\qff C^{\dff 2}\trf \nu_{\dff b}^{\dff 2}
\qff +\qff
2\trf r^{\dff 2}\qff C^{\dff 2}\trf \mu_{\fff a}^{\dff 2}
\qff +\qff
2\trf s^{\dff 2}\qff \lambda_{\dff k}^{\dff 2}
\]

\vspace{-30pt}
\[
\quad
\leq\off
2\trf C^{\dff 2}\qff
\left(\qff 
r^{\dff 2}\qff \lambda_{\dff i}^{\dff 2}
\qff +\qff
s^{\dff 2}\qff \nu_{\dff b}^{\dff 2}
\qff +\qff
r^{\dff 2}\qff \mu_{\fff a}^{\dff 2}
\qff +\qff
s^{\dff 2}\qff \lambda_{\dff k}^{\dff 2}
\qff\right)
\]

\vspace{-30pt}
\[
\quad
=\off
2\trf C^{\dff 2}\qff
\left(\qff 
r^{\dff 2}\qff
\left(\qff \lambda_{\dff i}^{\dff 2}
\qff +\qff
\mu_{\fff a}^{\dff 2}
\qff\right)
\qff +\qff
s^{\dff 2}\qff 
\left(\qff
\lambda_{\dff k}^{\dff 2}
\qff +\qff
\nu_{\dff b}^{\dff 2}
\qff\right)
\qff\right)
\off.
\]

\vspace{-12pt}\vspace{4.5pt}
It\sss follows\sss that\vspace{4.5pt}
\[
\quad
\bnorm{(\trf A\pff \omult_{\fff 1}\qff Q\trf)^{\dff 2\fff n\qff +\qff 1}
\qff
\left(\qff
(\trf u_{\dff i}\dff \otimes\dff v_a\trf)
\off \oplus\off
(\trf u_{\dff k}\dff \otimes\dff w_{\dff b}\trf)
\qff\right)}_{\qff 0}^{\qff 2}
\]

\vspace{-31.5pt}
\[
\quad
\leq\off
2\trf C^{\dff 2}\qff
\left(\qff 
\bigl(\trf \lambda_{\dff i}^{\dff 2} 
\off +\off
\mu_{\fff a}^{\dff 2} \qff\bigr)^{\dff 2\fff n}
\qff
\left(\qff \lambda_{\dff i}^{\dff 2}
\qff +\qff
\mu_{\fff a}^{\dff 2}
\qff\right)
\off +\off
\bigl(\trf \lambda_{\dff k}^{\dff 2} 
\off +\off
\nu_{\fff b}^{\dff 2} \qff\bigr)^{\dff 2\fff n}
\qff
\left(\qff
\lambda_{\dff k}^{\dff 2}
\qff +\qff
\nu_{\dff b}^{\dff 2}
\qff\right)
\qff\right)
\]

\vspace{-31.5pt}
\[
\quad
\leq\off
2\trf C^{\dff 2}\qff
\left(\qff 
\bigl(\trf \lambda_{\dff i}^{\dff 2} 
\off +\off
\mu_{\fff a}^{\dff 2} \qff\bigr)^{\dff 2\fff n\qff +\qff 1}
\off +\off
\bigl(\trf \lambda_{\dff k}^{\dff 2} 
\off +\off
\nu_{\fff b}^{\dff 2} \qff\bigr)^{\dff 2\fff n\qff +\qff 1}
\qff\right)
\off.
\]

\vspace{-12pt}\vspace{4.5pt}
Together with\sss the above calculation of\dss the norm for even\sss powers
of\dss $A\pff \omult_{\fff 1}\qff Q$\sss this estimate implies\sss that\sss
$(\trf u_{\dff i}\dff \otimes\dff v_a\trf)
\off \oplus\off
(\trf u_{\dff k}\dff \otimes\dff w_{\dff b}\trf)$\sss
are bounded\sss vectors.\oss  \eproof

\mypar{Theorem.}{properties-product}
\emph{The operator\sss $A\pff \omult_{\fff 1}\qff Q$\dss
is\dss self-adjoint.\oss
It\dss is\dss a\dss Fredholm\dss operator\sss with discrete spectrum and compact\sss
resolvent.\oss 
Its\sss kernel\dss is}\vspace{3pt}
\[
\quad
\kernel\dff A\pff \omult_{\fff 1}\qff Q
\off =\off 
(\trf \kernel\fff A\dff \otimes\dff \kernel\fff Q \trf)
\qff \oplus\qff
(\trf \kernel\fff A\dff \otimes\dff \kernel\fff Q^{\fff *} \trf)
\pff,
\]

\vspace{-12pt}\vspace{3pt}
\emph{and\dss its\dss image\dss is\dss equal\dss to\sss the orthogonal\dss
complement\dss to\sss the\sss kernel.\oss}

\proof
By\sss the definition\sss $A\pff \omult_{\fff 1}\qff Q$\sss is\dss closed.\oss
By\trs Lemma\qss \ref{bounded-vectors}\qss 
the subspaces of\dss
bounded,\oss and\dss hence of\dss analytic,\oss vectors are dense in\sss 
$(\trf H_{\dff 0}\qff\fff \widehat{\otimes}\pff K_{\trf 0} \trf)
\dff \oplus\dff
(\trf H_{\dff 0}\qff\fff \widehat{\otimes}\pff L_{\trf 0} \trf)$\nnsp.\oss
Therefore\dss Nelson\dss theorem\qss
(see\qss \cite{s},\oss Theorem\qss 7.16)\qss 
implies\sss that\sss
$A\pff \omult_{\fff 1}\qff Q$\sss is\dss self-adjoint.\oss 

Since\sss $A\pff \omult_{\fff 1}\qff Q$\sss is\dss self-adjoint,\pss 
$(\trf A\pff \omult_{\fff 1}\qff Q\trf)^{\trf 2}$\sss
is\dss self-adjoint.\oss
We already\sss used\sss the fact\sss that\vspace{3pt}
\[
\quad
(\trf A\pff \omult_{\fff 1}\qff Q\trf)^{\dff 2}
\qff
\left(\qff
(\trf u_{\dff i}\dff \otimes\dff v_a\trf)
\off \oplus\off
(\trf u_{\dff k}\dff \otimes\dff w_{\fff b}\trf)
\qff\right)
\]

\vspace{-33pt}
\[
\quad
=\off
\bigl(\trf \lambda_{\dff i}^{\dff 2} 
\off +\off
\mu_{\fff a}^{\dff 2} \qff\bigr)
\qff
\bigl(\trf u_{\dff i}\dff \otimes\dff v_a \trf\bigr)
\off \oplus\off
\bigl(\trf \lambda_{\dff k}^{\dff 2} 
\off +\off
\nu_{\fff b}^{\dff 2} \qff\bigr)\qff
\bigl(\trf u_{\dff k}\dff \otimes\dff w_{\fff b} \trf\bigr)
\off.
\]

\vspace{-12pt}\vspace{3pt}
This shows\sss that\sss the vectors\sss $u_{\dff i}\dff \otimes\dff v_a$\sss
and\sss $u_{\dff k}\dff \otimes\dff w_{\fff b}$\sss are eigenvectors of\dss
$(\trf A\pff \omult_{\fff 1}\qff Q\trf)^{\trf 2}$\sss
and\vspace{3pt}
\[
\quad
\lambda_{\dff i}^{\dff 2} 
\off +\off
\mu_{\fff a}^{\dff 2}
\quad
\mbox{and}\quad
\lambda_{\dff k}^{\dff 2} 
\off +\off
\nu_{\fff b}^{\dff 2}
\]

\vspace{-12pt}\vspace{3pt}
are\sss the corresponding eigenvalues.\oss
At\sss the same\sss time\sss these vectors form an orthonormal\dss basis of\dss
$(\trf H_{\dff 0}\qff\fff \widehat{\otimes}\pff K_{\trf 0} \trf)
\dff \oplus\dff
(\trf H_{\dff 0}\qff\fff \widehat{\otimes}\pff L_{\trf 0} \trf)$\nnsp.\oss
Since\sss the eigenvalues\sss tend\dss to infinity\sss when\sss
$(\trf i\fff,\qff a\trf)\qff \ttoo\qff \infty$ or\sss
$(\trf k\fff,\qff b\trf)\qff \ttoo\qff \infty$\nnsp,\oss
the operator\sss $(\trf A\pff \omult_{\fff 1}\qff Q\trf)^{\dff 2}$\sss
is\dss an operator\sss with discrete spectrum and compact\sss resolvent.\oss
See\qss \cite{s},\oss Proposition\qss 5.12.\oss
In\sss particular,\pss the eigenspaces of\sss
$(\trf A\pff \omult_{\fff 1}\qff Q\trf)^{\trf 2}$\sss
are finitely dimensional.\oss
Clearly,\oss these eigenspaces are invariant\sss under\sss
$A\pff \omult_{\fff 1}\qff Q$\nnsp.\oss
This implies\sss that\sss there exists an orthonormal\dss basis of\dss
$(\trf H_{\dff 0}\qff\fff \widehat{\otimes}\pff K_{\trf 0} \trf)
\dff \oplus\dff
(\trf H_{\dff 0}\qff\fff \widehat{\otimes}\pff L_{\trf 0} \trf)$
consisting of\dss eigenvectors of\dss
$A\pff \omult_{\fff 1}\qff Q$\nnsp.\oss
The corresponding eigenvalues\sss tend\sss to infinity\sss
together\sss with\sss the eigenvalues of\dss
$(\trf A\pff \omult_{\fff 1}\qff Q\trf)^{\trf 2}$\dnsp,\oss
and\dss hence\sss $A\pff \omult_{\fff 1}\qff Q$\sss
has discrete spectrum and compact\sss resolvent.\oss

Since\sss the operator $A\pff \omult_{\fff 1}\qff Q$\sss is\dss self-adjoint,\oss
its kernel\dss is\dss equal\dss to\sss the kernel\sss of\dss its square\sss
$(\trf A\pff \omult_{\fff 1}\qff Q\trf)^{\trf 2}$\dnsp.\oss
The\sss latter\dss is\dss spanned\dss by eigenvectors\sss
$u_{\dff i}\dff \otimes\dff v_a$\sss
and\sss $u_{\dff k}\dff \otimes\dff w_{\fff b}$\sss
with\sss the eigenvalues equal\dss to $0$\nnsp.\oss
Together with\sss the above formulas for eigenvalues\sss this implies\sss
the claim about\sss the kernel.\oss
In\sss particular,\oss the kernel\dss is\dss finitely dimensional.\oss
Since\sss $A\pff \omult_{\fff 1}\qff Q$\sss is\dss self-adjoint,\oss
this immediately\sss implies\sss the claim about\sss the image
and\dss the\dss Fredholm\dss property.\oss  \eproof

\mypar{Lemma.}{eigenspaces-square}
\emph{Every eigenspace of\qss
$(\trf A\pff \omult_{\fff 1}\qff Q\trf)^{\trf 2}$\sss
is\dss equal\dss to a finite direct\sss sum of\dss 
algebraic\sss tensor products of\qss the\dss form\dss
$V\dff \otimes\dff W$\dnsp,\oss
where\qss $V\qff \subset\pff H_{\trf 0}$\sss 
is\trs an eigenspace of\dss $A^{\dff 2}$\dnsp,\oss
and either\qss $W\qff \subset\qff K_{\trf 0}$\sss 
is\trs an eigenspace\sss of\trs $Q^{\fff *}\fff Q$\nnsp,\oss
or\qss $W\qff \subset\qff L_{\trf 0}$
is\trs  an eigenspace\sss of\qss $Q\dff Q^{\fff *}$\nsp\dnsp.\oss}

\proof
This follows\sss from\sss the formulas for\sss the eigenvalues of\dss
$(\trf A\pff \omult_{\fff 1}\qff Q\trf)^{\trf 2}$\sss
corresponding\sss to eigenvectors\sss
$u_{\dff i}\dff \otimes\dff v_a$\sss
and\sss $u_{\dff k}\dff \otimes\dff w_{\fff b}$\sss
which were used\sss in\sss the proof\dss of\trs 
Theorem\qss \ref{properties-product}.\oss  \eproof

\myuppar{Tensor analogues of\trs Sobolev\dss spaces.}
In our applications\sss $H_{\trf 0}\dff,\qff H_{\dff 1}\dff,\qff K_{\trf 0}$\sss etc.\dss
will\dss be\dss 
Sobolev\dss spaces of\dss sections of\dss vector bundles.\oss
Their analogues for\sss tensor\sss products are\sss the space\vspace{1.5pt}
\[
\quad
(\trf H\qff\fff \widehat{\otimes}\pff K \trf)_{\dff 0}
\off =\off
H_{\dff 0}\qff\fff \widehat{\otimes}\pff K_{\trf 0}
\off
\]

\vspace{-12pt}\vspace{1.5pt}
and\sss the space 
$(\trf H\qff\fff \widehat{\otimes}\pff K \trf)_{\dff 1}$
which we are going\sss to define.\oss
For subscripts\sss
$a\fff,\qff b\qff \in\qff \{\qff 0\fff,\qff 1 \qff\}$\sss
let\sss $\sco{\trf \bullet\fff,\qff \bullet \trf}_{\dff a\fff b}$\sss
be\sss the scalar\sss product\sss on\sss 
$H_{\dff a}\qff\fff \widehat{\otimes}\pff K_{\dff b}$\nsp.\oss
Then\vspace{1.5pt}
\[
\quad
\sco{\trf \bullet\fff,\qff \bullet \trf}_{\dff 1}
\off =\off
\sco{\trf \bullet\fff,\qff \bullet \trf}_{\dff 1\fff 0}
\pff +\qff
\sco{\trf \bullet\fff,\qff \bullet \trf}_{\dff 0\fff 1}
\pff
\]

\vspace{-12pt}\vspace{1.5pt}
is\dss a scalar product\sss on $H_{\dff 1}\dff \otimes\dff K_{\dff 1}$\nnsp.\oss
Let\sss
$(\trf H\qff\fff \widehat{\otimes}\pff K \trf)_{\dff 1}$
be\sss the completion of\dss 
$H_{\dff 1}\dff \otimes\dff K_{\dff 1}$\sss with respect\sss to\sss
this scalar product.\oss
The inclusion of\sss 
$H_{\dff 1}\dff \otimes\dff K_{\dff 1}$
into
$H_{\dff 0}\qff\fff \widehat{\otimes}\pff K_{\trf 0}$
extends\sss to a bounded\dss map\vspace{1.5pt}
\[
\quad
\iota
\dff \colon\dff
(\trf H\qff\fff \widehat{\otimes}\pff K \trf)_{\dff 1}
\qff \ttoo\qff
H_{\dff 0}\qff\fff \widehat{\otimes}\pff K_{\trf 0}
\off.
\]

\vspace{-12pt}\vspace{1.5pt}
Let\sss
$(\trf D\qff\fff \widehat{\otimes}\pff K \trf)_{\dff 1}$\sss
be\sss the preimage of\dss
$D_{\dff 1}\qff\fff \widehat{\otimes}\pff K_{\trf 0}$\sss
under\sss the inclusion\sss 
$(\trf H\qff\fff \widehat{\otimes}\pff K \trf)_{\dff 1}
\qff \ttoo\qff
H_{\dff 1}\qff\fff \widehat{\otimes}\pff K_{\trf 0}$
Since\sss this inclusion\dss is\dss bounded,\qss
$(\trf D\qff\fff \widehat{\otimes}\pff K \trf)_{\dff 1}$
is\dss a closed\sss subspace of\dss
$(\trf H\qff\fff \widehat{\otimes}\pff K \trf)_{\dff 1}$\nsp.\oss
The space
$(\trf D\qff\fff \widehat{\otimes}\pff K \trf)_{\dff 1}$
can be also defined as\sss the completion of\dss $D_{\dff 1}\dff \otimes\dff K_{\dff 1}$\sss
with respect\sss to\sss 
$\sco{\trf \bullet\fff,\qff \bullet \trf}_{\dff 1}$\nsp.\oss
Let\sss\vspace{1.5pt}
\[
\quad
\iota_{\dff D}
\dff \colon\dff
(\trf D\qff\fff \widehat{\otimes}\pff K \trf)_{\dff 1}
\qff \ttoo\qff
H_{\dff 0}\qff\fff \widehat{\otimes}\pff K_{\trf 0}
\off
\]

\vspace{-12pt}\vspace{1.5pt}
be\sss the restriction of\dss $\iota$\nnsp.\oss
Of\dss course,\oss these constructions apply 
also\sss to $L$ in\sss the role of\dss $K$\nnsp.\oss

\mypar{Lemma.}{scale-tensor}
\emph{\dnsp$\iota$\sss is\dss an\sss injective map\sss with\sss the image\dss 
$(\trf H_{\dff 1}\qff\fff \widehat{\otimes}\pff K_{\trf 0} \trf)
\dff \cap\dff
(\trf H_{\dff 0}\qff\fff \widehat{\otimes}\pff K_{\dff 1} \trf)$\nnsp.\pss
The image of\qss $\iota_{\dff D}$\sss is\dss equal\dss to\sss
$(\trf D_{\dff 1}\qff\fff \widehat{\otimes}\pff K_{\trf 0} \trf)
\dff \cap\dff
(\trf H_{\dff 0}\qff\fff \widehat{\otimes}\pff K_{\dff 1} \trf)$\nnsp.\oss}

\proof
Clearly,\oss the map\sss $\iota$\sss factors\sss through\sss the natural\dss maps\sss 
$H_{\dff 1}\qff\fff \widehat{\otimes}\pff K_{\trf 0}
\qff \ttoo\qff
H_{\trf 0}\qff\fff \widehat{\otimes}\pff K_{\trf 0}$\sss
and\sss
$H_{\trf 0}\qff\fff \widehat{\otimes}\pff K_{\dff 1}
\qff \ttoo\qff
H_{\trf 0}\qff\fff \widehat{\otimes}\pff K_{\trf 0}$
and\dss hence its image\dss is\dss contained\sss in\sss the indicated\sss intersection.\oss
In order\sss to prove\sss the opposite inclusion,\oss let\sss 
$x\qff \in\qff H_{\dff 0}\qff\fff \widehat{\otimes}\pff K_{\trf 0}$\nsp.\oss
We can\sss present\sss $x$\sss as\sss the sum of\dss a series of\dss the form\sss
$\sum\fff c_{\dff i\fff a}\trf (\trf u_{\dff i}\dff \otimes\dff v_a\trf)$\sss
convergent\sss in $H_{\dff 0}\qff\fff \widehat{\otimes}\pff K_{\trf 0}$\nsp.\oss
Then\sss
$x\qff \in\qff H_{\dff 1}\qff\fff \widehat{\otimes}\pff K_{\trf 0}$\sss
if\dss and\sss only\trs if\trs this series converges in\sss
$H_{\dff 1}\qff\fff \widehat{\otimes}\pff K_{\trf 0}$\nsp,\oss
and
$x\qff \in\qff H_{\trf 0}\qff\fff \widehat{\otimes}\pff K_{\dff 1}$\sss
if\dss and\sss only\trs if\trs this series converges in\sss
$H_{\trf 0}\qff\fff \widehat{\otimes}\pff K_{\dff 1}$\nsp.\oss
Therefore $x$ belongs\sss to\sss the intersection\dss
if\dss and\sss only\trs if\trs this series converges in\sss both\sss
$H_{\dff 1}\qff\fff \widehat{\otimes}\pff K_{\trf 0}$\sss
and\sss
$H_{\trf 0}\qff\fff \widehat{\otimes}\pff K_{\dff 1}$\nsp,\oss
i.e.\qss with respect\sss to\sss the scalar products
$\sco{\trf \bullet\fff,\qff \bullet \trf}_{\dff 0\fff 1}$
and\sss
$\sco{\trf \bullet\fff,\qff \bullet \trf}_{\dff 1\fff 0}$\nsp.\oss
In\sss this case\sss this series also converges with respect\sss to\sss
$\sco{\trf \bullet\fff,\qff \bullet \trf}_{\dff 1}$\nsp,\oss
i.e.\qss in\sss $(\trf H\qff\fff \widehat{\otimes}\pff K \trf)_{\dff 1}$\nsp,\oss
and $x$ belongs\sss to\sss the image of\dss
$\iota\dff \colon\dff
(\trf H\qff\fff \widehat{\otimes}\pff K \trf)_{\dff 1}
\qff \ttoo\qff
(\trf H\qff\fff \widehat{\otimes}\pff K \trf)_{\dff 0}$\nsp.\oss
This proves\sss the claim about\sss the image of\dss $\iota$\nnsp.\oss

In order\sss to prove its injectivity,\oss
suppose\sss that\sss
$x_{\dff i}\qff \in\qff H_{\dff 1}\dff \otimes\dff K_{\dff 1}$\sss
is\dss a\dss Cauchy\dss sequence with respect\sss to\sss
$\sco{\trf \bullet\fff,\qff \bullet \trf}_{\dff 1}$\sss
converging\sss to $0$ in\sss
$H_{\dff 0}\qff\fff \widehat{\otimes}\pff K_{\trf 0}$\nsp.\oss
Then\sss $x_{\dff i}$\sss is\dss also a\dss Cauchy\dss sequence with respect\sss to\sss
$\sco{\trf \bullet\fff,\qff \bullet \trf}_{\dff 1\fff 0}$
and\sss
$\sco{\trf \bullet\fff,\qff \bullet \trf}_{\dff 0\fff 1}$\nsp.\oss
Let\sss $x_{\trf 1}\dff,\off x_{\trf 1\fff 0}$\sss and\sss $x_{\trf 0\fff 1}$\sss 
be\sss the\dss limits of\dss this sequence in\sss
$(\trf H\qff\fff \widehat{\otimes}\pff K \trf)_{\dff 1}\dff,\off 
H_{\dff 1}\qff\fff \widehat{\otimes}\pff K_{\trf 0}$\nsp,\oss
and\sss
$H_{\trf 0}\qff\fff \widehat{\otimes}\pff K_{\trf 1}$
respectively.\oss
The injectivity of\dss
$H_{\dff 1}\qff\fff \widehat{\otimes}\pff K_{\trf 0}
\qff \ttoo\qff
H_{\trf 0}\qff\fff \widehat{\otimes}\pff K_{\trf 0}$\sss
implies\sss that\sss $x_{\trf 1\fff 0}\off =\off 0$\nnsp.\oss
Similarly,\pss $x_{\trf 0\fff 1}\off =\off 0$\nnsp.\oss
Hence
$\sco{\trf x_{\dff i}\fff,\qff x_{\dff i} \trf}_{\dff 1\fff 0}$
and
$\sco{\trf x_{\dff i}\fff,\qff x_{\dff i} \trf}_{\dff 0\fff 1}$\sss
tend\sss to $0$\nnsp,\oss
and\dss hence\sss
$\sco{\trf x_{\dff i}\fff,\qff x_{\dff i} \trf}_{\dff 1}$
also\sss tends\sss to $0$\nnsp.\oss
Therefore\sss $x_{\trf 1}\off =\off 0$\nnsp.\oss
This proves\sss the injectivity of\sss $\iota$\nnsp.\oss
Since\sss $\iota_{\dff D}$\sss is\dss the restriction of\dss $\iota$\nnsp,\oss
the claim about\sss $\iota_{\dff D}$\sss 
follows from\sss the claim about\sss $\iota$\nnsp.\oss \eproof

\myuppar{Modifying\dss Hilbert\dss spaces $H_{\dff 1}\dff,\pff K_{\dff 1}\dff,\pff L_{\dff 1}$\nsp.}
Let $\mathcal{D}_{\dff 1}$\sss be\sss the closure\sss in $D_{\dff 1}$
of\dss the\sss linear span of\dss vectors $u_{\dff i}$ with\sss
$\lambda_{\dff i}\off \neq\off 0$\nnsp.\oss
Since $A$\sss is\dss Fredholm,\oss
the identity\sss map\sss 
$\mathcal{D}_{\dff 1}\qff \ttoo\qff \mathcal{D}_{\dff 1}$\sss
is\dss a\sss topological\dss isomorphism\sss between\sss
the structure of\dss a\dss Hilbert\dss space on\sss $\mathcal{D}_{\dff 1}$\sss
induced\sss from\sss $H_{\dff 1}$\sss and\sss the structure having\sss
the vectors\sss $\lambda_{\dff i}^{\dff -\dff 1}\dff u_{\dff i}$\sss with\sss
$\lambda_{\dff i}\off \neq\off 0$\sss
as an orthonormal\sss basis.\oss
Let\sss us replace\sss the original\sss structure of\dss a\dss Hilbert\sss space\sss 
$D_{\dff 1}$\sss by\sss the direct\sss sum of\dss the\sss latter structure on\sss
$\mathcal{D}_{\dff 1}$\sss and\sss the structure on $\kernel\fff A$
induced\sss from\sss $H_{\trf 0}$\nsp.\oss
Similarly,\oss let\sss change\sss the structures of\dss a\dss Hilbert\dss space on
$K_{\dff 1}$ and\sss $L_{\dff 1}$\sss using\sss vectors\sss
$\mu_{\fff a}^{\dff -\dff 1}\dff v_a$\sss and\sss $\nu_{\fff b}^{\dff -\dff 1}\dff w_{\fff b}$\nsp.\oss

Since\sss the new structures of\trs Hilbert\dss spaces on\sss $H_{\dff 1}$\sss and\sss $K_{\dff 1}$
are\sss topologically\sss isomorphic\sss to\sss the original\sss ones,\oss
using\sss these new structures in\sss the definitions of\dss
$(\trf D\qff\fff \widehat{\otimes}\pff K \trf)_{\dff 1}$
and\sss
$(\trf D\qff\fff \widehat{\otimes}\pff L \trf)_{\dff 1}$
replaces\sss these\dss Hilbert\dss spaces by\sss topologically\sss isomorphic ones.\oss

\mypar{Lemma.}{tensor-basis}
\emph{For\sss the new\dss Hilbert\dss space structure\sss the vectors}\vspace{1.5pt}
\[
\quad
\bigl(\trf \lambda_{\dff i}^{\dff 2} 
\off +\off
\mu_{\dff a}^{\dff 2} \qff\bigr)^{\dff -\dff 1/2}\qff
\bigl(\trf u_{\dff i}\dff \otimes\dff v_a \trf\bigr)
\quad
\mbox{\emph{such\sss that}}\quad\dff
\lambda_{\dff i}^{\dff 2}\qff +\qff \mu_{\dff a}^{\dff 2}
\off \neq\off
0
\quad
\mbox{\emph{and}}\quad
\]

\vspace{-34.5pt}
\[
\quad
u_{\dff i}\dff \otimes\dff v_a
\quad
\mbox{\emph{such\sss that}}\quad\dff
\lambda_{\dff i}\off =\off \mu_{\fff a}\off =\off 0
\]

\vspace{-12pt}\vspace{1.5pt}
\emph{form an orthonormal\dss basis of\qss 
$(\trf D\qff\fff \widehat{\otimes}\pff K \trf)_{\dff 1}$\nsp.\oss}

\proof
With\sss the new\dss Hilbert\sss space structures\vspace{3pt}
\[
\quad
\sco{\trf u_{\dff i}\dff \otimes\dff v_a\fff,\qff u_{\dff i}\dff \otimes\dff v_a \trf}_{\dff 1}
\off
=\off
\sco{\trf u_{\dff i}\dff \otimes\dff v_a\fff,\qff u_{\dff i}\dff \otimes\dff v_a \trf}_{\dff 1\fff 0}
\pff +\qff
\sco{\trf u_{\dff i}\dff \otimes\dff v_a\fff,\qff u_{\dff i}\dff \otimes\dff v_a \trf}_{\dff 0\fff 1}
\]

\vspace{-34.5pt}
\[
\quad
=\off
\sco{\trf u_{\dff i}\fff,\qff u_{\dff i} \trf}_{\dff 1}\qff
\sco{\trf v_a\fff,\qff v_a \trf}_{\dff 0}
\pff +\qff
\sco{\trf u_{\dff i}\fff,\qff u_{\dff i} \trf}_{\dff 0}\qff
\sco{\trf v_a\fff,\qff v_a \trf}_{\dff 1}
\off =\off
\lambda_{\dff i}^{\dff 2}
\pff +\qff
\mu_{\dff a}^{\dff 2}
\pff.
\]

\vspace{-12pt}\vspace{3pt}
Similarly,\pss 
$\sco{\trf u_{\dff i}\dff \otimes\dff v_a\fff,\qff u_{\dff k}\dff \otimes\dff v_{\fff b} \trf}_{\dff 1}
\off =\off
0$
unless\sss
$(\trf i\fff,\qff a\trf)\off =\off (\trf k\fff,\qff b\trf)$\nnsp.\oss
It\sss follows\sss that\sss the vectors from\sss the\sss lemma
have\sss the norm $1$ and are pairwise orthogonal\sss 
with respect\sss to\sss the scalar\sss product\sss
$\sco{\trf \bullet\fff,\qff \bullet \trf}_{\dff 1}$\nsp.\oss
Being\sss the completion of\dss $D_{\dff 1}\dff \otimes\dff K_{\dff 1}$\sss
with respect\sss to\sss 
$\sco{\trf \bullet\fff,\qff \bullet \trf}_{\dff 1}$\nsp,\oss
the space\sss $(\trf D\qff\fff \widehat{\otimes}\pff K \trf)_{\dff 1}$\sss
is\dss the\dss Hilbert\sss space having\sss the vectors from\sss the\sss lemma
as its orthonormal\dss basis.\oss  \eproof

\mypar{Theorem.}{domain-of-product}
\emph{The domain of\qss $A\pff \omult_{\fff 1}\qff Q$\sss is\dss equal\dss to\sss
$(\trf D\qff\fff \widehat{\otimes}\pff K \trf)_{\dff 1}
\dff \oplus\dff
(\trf D\qff\fff \widehat{\otimes}\pff L \trf)_{\dff 1}$\nsp,\oss
and\sss the inclusion of\qss the domain\sss into\sss
$(\trf H_{\dff 0}\qff\fff \widehat{\otimes}\pff K_{\trf 0} \trf)
\qff \oplus\qff
(\trf H_{\dff 0}\qff\fff \widehat{\otimes}\pff L_{\trf 0} \trf)$
is\dss compact.\oss}

\proof
By\trs Theorem\qss \ref{properties-product}\qss the operator\sss 
$A\pff \omult_{\fff 1}\qff Q$\sss is\dss self-adjoint\dss and\sss has compact\sss resolvent.\oss
It\sss follows\sss that\sss there exists an orthonormal\dss basis of\dss
$(\trf H_{\dff 0}\qff\fff \widehat{\otimes}\pff K_{\trf 0} \trf)
\qff \oplus\qff
(\trf H_{\dff 0}\qff\fff \widehat{\otimes}\pff L_{\trf 0} \trf)$
consisting of\dss eigenvectors $e_{\dff n}$\sss of\dss
$A\pff \omult_{\fff 1}\qff Q$\nnsp.\oss
Let\sss $\xi_{\trf n}$\sss be\sss the corresponding eigenvalues.\oss
Then\sss the domain of\dss $A\pff \omult_{\fff 1}\qff Q$\sss
can\sss be identified\sss with\sss the\dss Hilbert\dss space of\dss formal\sss series\sss
$\sum\trf c_{\dff n}\qff e_{\dff n}$\sss
such\sss that\sss\vspace{1.5pt}
\[
\quad
\sum\qff \num{c_{\dff n}}^{\dff 2}\off \xi_{\trf n}^{\dff 2}
\]

\vspace{-12pt}\vspace{1.5pt}
converges.\oss 
Every\sss eigenvector $e_{\dff n}$\sss is\dss also an eigenvector of\dss
$(\trf A\pff \omult_{\fff 1}\qff Q\trf)^{\trf 2}$\sss
with\sss the eigenvalue\sss $\xi_{\trf n}^{\dff 2}$\nsp.\oss
As we saw in\sss the proof\dss of\trs Theorem\qss \ref{properties-product},\oss
this implies\sss that\vspace{1.5pt}
\[
\quad
\xi_{\trf n}^{\dff 2}
\off\dff =\off
\lambda_{\dff i}^{\dff 2} 
\off +\off
\mu_{\fff a}^{\dff 2}
\quad\dff
\mbox{or}\quad
\lambda_{\dff k}^{\dff 2} 
\off +\off
\nu_{\fff b}^{\dff 2}
\pff.
\]

\vspace{-12pt}\vspace{1.5pt}
Also,\qss $e_{\dff n}$\sss is\dss a\sss linear\sss combination
of\dss vectors\sss $u_{\dff i}\dff \otimes\dff v_a$\sss
and\sss $u_{\dff k}\dff \otimes\dff w_{\fff b}$\sss
belonging\sss to\sss the same eigenspace of\dss
$(\trf A\pff \omult_{\fff 1}\qff Q\trf)^{\trf 2}$\dnsp,\oss
i.e.\qss with\sss the same values of\dss 
$\lambda_{\dff i}^{\dff 2} 
\off +\off
\mu_{\fff a}^{\dff 2}$\sss
or\sss
$\lambda_{\dff k}^{\dff 2} 
\off +\off
\nu_{\fff b}^{\dff 2}$.\oss 
Rewriting eigenvectors $e_{\dff n}$ as such\sss linear combinations 
shows\sss that\sss the domain of\dss $A\pff \omult_{\fff 1}\qff Q$\sss
can\sss be identified\sss with\sss the\dss Hilbert\dss space of\dss expressions\vspace{3pt}\vspace{-0.5pt}
\[
\quad
\left(\pff
\sum\qff c_{\dff i\dff a}\qff u_{\dff i}\dff \otimes\dff v_a
\qff\right)
\qff \oplus\qff
\left(\pff
\sum\qff d_{\trf k\fff b}\qff u_{\dff k}\dff \otimes\dff w_{\fff b}
\qff\right)
\]

\vspace{-12pt}\vspace{1.5pt}
such\sss that\sss the sums\vspace{1.5pt}
\[
\quad
\sum\qff \num{c_{\dff i\dff a}}^{\dff 2}\qff
\left(\qff
\lambda_{\dff i}^{\dff 2} 
\off +\off
\mu_{\fff a}^{\dff 2}
\qff\right)
\quad
\mbox{and}\quad\dff
\sum\qff \num{d_{\trf k\fff b}}^{\dff 2}\qff
\left(\qff
\lambda_{\dff k}^{\dff 2} 
\off +\off
\nu_{\fff b}^{\dff 2}
\qff\right)
\]

\vspace{-12pt}\vspace{3pt}\vspace{-0.5pt}
converge.\oss
Lemma\qss \ref{tensor-basis}\qss and\sss its analogue for\sss
$(\trf D\qff\fff \widehat{\otimes}\pff L \trf)_{\dff 1}$\sss
imply\sss that\dss this\dss is\dss nothing else but
$(\trf D\qff\fff \widehat{\otimes}\pff K \trf)_{\dff 1}
\dff \oplus\dff
(\trf D\qff\fff \widehat{\otimes}\pff L \trf)_{\dff 1}$.\oss
This proves\sss the first\sss statement\sss of\dss the\sss theorem.\oss

The norms in\sss
$H_{\dff 0}\qff\fff \widehat{\otimes}\pff K_{\trf 0}$\sss 
of\dss vectors\sss from\trs Lemma\qss \ref{tensor-basis}\qss are\sss 
$(\trf \lambda_{\dff i}^{\dff 2} 
\off +\off
\mu_{\dff a}^{\dff 2} \qff)^{\dff -\dff 1/2}$\dnsp.\oss
Since\sss this expression\sss tends\sss to $0$ when 
$(\trf i\fff,\qff a\trf)\trf \ttoo\qff \infty$\nnsp,\oss
the inclusion\sss
$(\trf D\qff\fff \widehat{\otimes}\pff K \trf)_{\dff 1}
\qff \ttoo\qff
H_{\dff 0}\qff\fff \widehat{\otimes}\pff K_{\trf 0}$\sss
is\dss compact.\oss
Similarly,\pss 
$(\trf D\qff\fff \widehat{\otimes}\pff L \trf)_{\dff 1}
\dff \ttoo\qff
H_{\dff 0}\qff\fff \widehat{\otimes}\pff L_{\trf 0}$\sss
is\dss compact.\oss
The second statement\sss follows.\oss  \eproof

\myuppar{Fredholm\dss families.}
Let\sss us recall\sss some notions from\qss \cite{i2}.\oss 
Let\sss $H$\sss be a separable\dss Hilbert\dss space.\oss
A pair $(\trf A\fff,\qff \varepsilon\trf)$\sss is\dss said\dss to be an\qss
\emph{enhanced operator}\pss if\dss 
$A\dff \colon\dff H\qff \ttoo\qff H$\dss 
is\dss an\sss unbounded self-adjoint\sss operator,\oss
$\varepsilon\qff >\qff 0$\nnsp,\oss
the numbers\sss $\varepsilon$ and\sss $-\qff \varepsilon$\sss
do not\sss belong\sss to\sss the spectrum\sss
$\sigma\trf(\trf A\trf)$\nnsp,\oss
and\sss the spectral\dss projection\sss
$P_{\qff [\dff -\dff \varepsilon\fff,\qff \varepsilon\dff]}\trf(\trf A\trf)$
has finite rank.\oss
Let\sss $Z$\sss be a\sss topological\sss space.\oss
A family\sss
$A_{\dff z}\dff,\pff z\qff \in\qff Z$\sss
of\dss self-adjoint\sss unbounded operators\sss
$H\qff \ttoo\qff H$\dss is\dss said\dss to be a\qss
\emph{Fredholm\dss family}\pss
if\dss for every\sss $z\qff \in\qff Z$\sss
the operator\sss $A_{\dff z}$\sss is\dss Fredholm\dss
and\sss there exsist\sss 
$\varepsilon\off =\off \varepsilon_{\dff z}\qff >\qff 0$\sss
and a neighborhood\sss $U_{\dff z}$\sss of\sss $z$
with\sss the following\sss properties.\oss
First,\oss the pair\sss $(\trf A_{\dff y}\fff,\qff \varepsilon\trf)$\sss
is\dss an enhanced operator\sss for for every\sss
$y\qff \in\qff U_{\dff z}$\nsp.\oss
Second,\oss the subspaces\vspace{1.5pt}
\[
\quad
V_{\fff y}
\off =\off
\image P_{\qff [\dff -\dff \varepsilon\fff,\qff \varepsilon\dff]}\trf(\trf A_{\dff y}\trf)
\]

\vspace{-12pt}\vspace{1.5pt}
depend continuously on $y\qff \in\qff U_{\dff z}$.\oss
Finally,\oss the operators\sss
$V_{\fff y}
\qff \ttoo\qff
V_{\fff y}$
induced\sss by $A_{\dff y}$
depend continuously on \nsp$y\qff \in\qff U_{\dff z}$.\qss
For example,\oss a\sss family continuous in\sss the\sss topology of\dss the convergence
in\sss the norm\sss resolvent\sss sense\dss is\dss Fredholm.\oss
The space $Z$\sss is\dss usually assumed\sss to be compactly\sss generated and\sss paracompact.\oss
Then\sss the\sss analytical\dss index of\trs 
$A_{\dff z}\dff,\pff z\qff \in\qff Z$\sss
is\dss defined
and\sss agrees with\sss the classical\dss analytical\dss index\qss \cite{as}\qss
when\sss the\sss latter\dss is\dss defined.\oss
See\qss \cite{i2}.

\myuppar{Families and\dss the $\omult_{\fff 1}$ products.}
Let\sss $Z$\sss be a\sss topological\sss space and 
$A_{\dff z}\dff \colon\dff H_{\trf 0}\qff \ttoo\qff H_{\trf 0}$\nsp,\qss
$z\qff \in\qff Z$\sss
be a\dss Fredholm\dss family of\dss unbounded operators in\sss $H_{\trf 0}$\nsp.\oss
Suppose\sss that\sss $A\off =\off A_{\dff z}$\sss satisfies all\sss our assumptions
with\sss respect\sss to\sss  
$H_{\trf 0}\dff,\off H_{\dff 1}$ and\sss $D_{\dff 1}$\sss
for every $z$\nnsp.\oss
Let\sss $Q\fff,\pff Q^{\fff *}$ be\sss two operators,\oss independent\sss of\dss $z$\sss
and satisfying\sss all\sss our assumptions with\sss respect\sss to\sss 
$K_{\trf 0}\dff,\off L_{\trf 0}\dff,\off K_{\dff 1}$\sss etc.\oss

We are interested\sss in\sss the relation between\sss the analytical\dss indices
of\dss the families\sss $A_{\dff z}\dff,\off z\qff \in\qff Z$\sss and\sss
$A_{\dff z}\pff \omult_{\fff 1}\qff Q\dff,\off z\qff \in\qff Z$\nnsp.\oss
Lemma\qss \ref{product-fredholm}\qss below\sss implies\sss that\sss
the analytical\dss index of\dss the family\sss
$A_{\dff z}\pff \omult_{\fff 1}\qff Q\dff,\off z\qff \in\qff Z$\dss
is\dss defined.\oss
It\sss turns out\sss that\sss under\sss a natural\sss assumption\sss about\sss $Q$\dss
the analytical\dss indices of\dss the\sss two families are equal.\oss

\mypar{Lemma.}{product-fredholm}
\emph{The\sss family\sss
$A_{\dff z}\pff \omult_{\fff 1}\qff Q\dff,\off z\qff \in\qff Z$\dss
is\trs Fredholm.\oss}

\proof
Lemma\qss \ref{eigenspaces-square}\qss 
implies\sss that\sss the family\sss
$(\trf A_{\dff z}\pff \omult_{\fff 1}\qff Q\trf)^{\trf 2}
\dff,\off z\qff \in\qff Z$\dss
is\trs Fredholm.\oss
Let\sss $z\qff \in\qff Z$ and\dss let\sss 
$\varepsilon\off =\off \varepsilon_{\dff z}$ and\sss $U_{\dff z}$ 
be such as in\sss the definition
of\dss the\dss Fredholm\dss property of\dss the family 
$(\trf A_{\dff y}\pff \omult_{\fff 1}\qff Q\trf)^{\trf 2}\fff,\off z\qff \in\qff Z$\nnsp.\oss
Clearly,\pss  $A_{\dff y}\pff \omult_{\fff 1}\qff Q$\sss leaves\vspace{3pt}
\[
\quad
V_{\fff y}
\off =\off
\image P_{\qff [\dff -\dff \varepsilon\fff,\qff \varepsilon\dff]}\qff
\left(\qff (\trf A_{\dff y}\pff \omult_{\fff 1}\qff Q\trf)^{\trf 2}
\qff\right)
\]

\vspace{-12pt}\vspace{3pt}
invariant\sss and\dss for\sss $y\qff \in\qff U_{\dff z}$\sss
the operators
$V_{\fff y}\dff \ttoo\qff V_{\fff y}$
induced\sss by $A_{\dff y}\pff \omult_{\fff 1}\qff Q$
depend continuously on $y$\nnsp.\oss 
This implies\sss that the family\sss
$A_{\dff z}\pff \omult_{\fff 1}\qff Q\dff,\off z\qff \in\qff Z$\trs
is\trs Fredholm.\oss  \eproof

\mypar{Theorem.}{mult-of-index}
\emph{Suppose\sss that\sss 
$\kernel\fff Q^{\fff *}\off =\off 0$\sss
and\dss $\dim\dff \kernel\fff Q\off =\off 1$\nnsp.\oss
Then\dss the analytical\dss index of\dss the family\dss 
$A_{\dff z}\pff \omult_{\fff 1}\qff Q\dff,\off z\qff \in\qff Z$\dss 
is\dss equal\dss to\sss the analytical\dss index of\trs the\sss family\dss
$A_{\dff z}\dff,\off z\qff \in\qff Z$\nnsp.\oss}

\proof
Since $\kernel\fff Q^{\fff *}\off =\off 0$\nnsp,\oss
Theorem\qss \ref{properties-product}\qss implies\sss that\sss
$\kernel\fff A_{\dff z}\pff \omult_{\fff 1}\qff Q
\off =\off
\kernel\fff A_{\dff z}\dff \otimes\dff \kernel\fff Q$\nnsp.\oss
Since\sss $\dim\dff \kernel\fff Q\off =\off 1$\nnsp,\oss
these kernels are canonically\sss isomorphic\sss to\sss $\kernel\fff A_{\dff z}$\nsp.\oss
Let\sss $z\qff \in\qff Z$\sss and\dss let\sss 
$\varepsilon\off =\off \varepsilon_{\dff z}\qff >\qff 0$\sss
be such\sss that\vspace{1.5pt}
\[
\quad
\image\fff
P_{\trf [\dff -\qff \varepsilon\fff,\qff \varepsilon\trf]}
\qff
\bigl(\qff A_{\dff z}\pff \omult_{\fff 1}\qff Q \qff\bigr)
\off =\off
\kernel\fff A_{\dff z}\dff \otimes\dff \kernel\fff Q
\pff.
\]

\vspace{-12pt}\vspace{1.5pt}
Then\sss $(\qff A_{\dff z}\pff \omult_{\fff 1}\qff Q\fff,\pff \varepsilon\qff)$\sss
is\dss an enhanced operator,\oss and,\oss moreover\halfff,\pss
$(\qff A_{\dff y}\pff \omult_{\fff 1}\qff Q\fff,\pff \varepsilon\qff)$\sss
is\dss an enhanced operator for $y\qff \in\qff U_{\dff z}$\sss
for some open\sss neighborhood\sss $U_{\dff z}$\sss of\dss $z$\nnsp.\oss
If\trs $U_{\dff z}$\sss is\dss sufficiently small,\oss then for\sss
$y\qff \in\qff U_{\dff z}$\sss the\sss subspaces\vspace{1.5pt}
\[
\quad
\image\fff
P_{\trf [\dff -\qff \varepsilon\fff,\qff \varepsilon\trf]}
\qff
\bigl(\qff A_{\dff y}\pff \omult_{\fff 1}\qff Q \qff\bigr)
\quad
\mbox{and}\quad
\image\fff
P_{\trf [\dff -\qff \varepsilon\fff,\qff \varepsilon\trf]}
\qff
\bigl(\qff A_{\dff z}\pff \omult_{\fff 1}\qff Q \qff\bigr)
\]

\vspace{-12pt}\vspace{1.5pt}
are close\sss to each other and,\oss in\sss particular\halfff,\oss have\sss the same dimension.\oss
Also,\oss if\trs $U_{\dff z}$\sss is\dss sufficiently small,\oss then\sss the pair\sss
$(\trf A_{\dff y}\dff,\pff \varepsilon\trf)$\sss
is\dss an enhanced operator for $y\qff \in\qff U_{\dff z}$\nsp,\oss
and\sss the subspaces\vspace{1.5pt}
\[
\quad
\image\fff
P_{\trf [\dff -\qff \varepsilon\fff,\qff \varepsilon\trf]}
\qff
\bigl(\qff A_{\dff z} \qff\bigr)
\quad
\mbox{and}\quad
\image\fff
P_{\trf [\dff -\qff \varepsilon\fff,\qff \varepsilon\trf]}
\qff
\bigl(\qff A_{\dff y} \qff\bigr)
\]

\vspace{-12pt}\vspace{1.5pt}
have\sss the same dimension.\oss
The image\dss
$\image\fff
P_{\trf [\dff -\qff \varepsilon\fff,\qff \varepsilon\trf]}
\qff
(\qff A_{\dff y} \qff)$\sss
is\dss equal\sss to\sss the sum of\dss eigenspaces of\dss $A_{\dff y}$\sss
corresponding\sss to eigenvalues in\sss the interval\sss
$[\trf -\qff \varepsilon\dff,\pff \varepsilon\trf]$\nnsp,\oss
or\halfff,\oss equivalently,\oss to\sss the direct\sss sum of\dss kernels\sss
$\kernel\fff (\trf A_{\dff y}\qff -\qff \lambda \trf)$\sss
with\sss
$\lambda
\qff \in\qff 
[\trf -\qff \varepsilon\dff,\pff \varepsilon\trf]$\nnsp.\oss
Clearly,\oss\vspace{3pt}\vspace{-0.5pt}
\[
\quad
\kernel\dff
\bigl(\qff
\bigl(\qff A_{\dff y}\pff \omult_{\fff 1}\qff Q \qff\bigr)
\qff -\qff
\lambda
\qff\bigr)
\off\qff =\off\qff
\kernel\dff
\begin{pmatrix}
\off\dff (\trf A_{\dff y}\qff -\qff \lambda\trf)\qff\fff \widehat{\otimes}\pff 1 &
1\qff\fff \widehat{\otimes}\pff Q^{\fff *} \off
\vspace{6pt} \\
\off\dff 1\qff\fff \widehat{\otimes}\pff Q &
(\trf -\qff A_{\dff y}\qff -\qff \lambda\trf)\qff\fff \widehat{\otimes}\pff 1 \off 
\end{pmatrix}
\off,
\]

\vspace{-12pt}\vspace{3pt}
and\dss hence\vspace{3pt}
\begin{equation}
\label{spectral-projection}
\quad
\kernel\dff
\bigl(\qff
\bigl(\qff A_{\dff y}\pff \omult_{\fff 1}\qff Q \qff\bigr)
\qff -\qff
\lambda
\qff\bigr)
\off\qff 
\supset\off\qff
\bigl(\qff
\kernel\fff (\trf A_{\dff y}\qff -\qff \lambda \trf)\dff \otimes\dff \kernel\fff Q
\qff\bigr)
\qff \oplus\qff
0
\pff.
\end{equation}

\vspace{-12pt}\vspace{3pt}
It\sss follows\sss that\vspace{3pt}
\[
\quad
\image\fff
P_{\trf [\dff -\qff \varepsilon\fff,\qff \varepsilon\trf]}
\qff
\bigl(\qff A_{\dff y}\pff \omult_{\fff 1}\qff Q \qff\bigr)
\off\qff 
\supset\off\qff
\left(\qff
\left(\qff
\image\fff
P_{\trf [\dff -\qff \varepsilon\fff,\qff \varepsilon\trf]}\qff
(\trf A_{\dff y} \trf)
\qff\right)
\dff \otimes\dff \kernel\fff Q
\qff\right)
\qff \oplus\qff
0
\pff.
\]

\vspace{-12pt}\vspace{3pt}
The dimensions of\dss these\sss two spaces are equal\dss to\sss the dimensions of\vspace{4.5pt}
\[
\quad
\image\fff
P_{\trf [\dff -\qff \varepsilon\fff,\qff \varepsilon\trf]}
\qff
(\trf A_{\dff z}\pff \omult_{\fff 1}\qff Q \trf)
\quad
\mbox{and}\quad
\left(\qff
\image\fff
P_{\trf [\dff -\qff \varepsilon\fff,\qff \varepsilon\trf]}
\qff
(\trf A_{\dff z} \trf)
\qff\right)
\dff \otimes\dff \kernel\fff Q
\]

\vspace{-12pt}\vspace{4.5pt}
respectively.\oss
But\sss under our assumptions\sss these\sss two spaces are canonically\sss isomorphic\sss to\sss
the\sss kernel\sss $\kernel\fff A_{\dff z}$
and\dss hence\sss their dimensions are equal.\oss
It\dss follows\sss that\vspace{3.75pt}
\[
\quad
\image\fff
P_{\trf [\dff -\qff \varepsilon\fff,\qff \varepsilon\trf]}
\qff
\bigl(\qff A_{\dff y}\pff \omult_{\fff 1}\qff Q \qff\bigr)
\off\qff 
=\off\qff
\left(\qff
\left(\qff
\image\fff
P_{\trf [\dff -\qff \varepsilon\fff,\qff \varepsilon\trf]}\qff
(\trf A_{\dff y} \trf)
\qff\right)
\dff \otimes\dff \kernel\fff Q
\qff\right)
\qff \oplus\qff
0
\pff.
\]

\vspace{-12pt}\vspace{3.75pt}
Together\sss with\sss the inclusion\qss (\ref{spectral-projection})\qss
this implies\sss that\vspace{3.75pt}
\[
\quad
\kernel\dff
\bigl(\qff
\bigl(\qff A_{\dff y}\pff \omult_{\fff 1}\qff Q \qff\bigr)
\qff -\qff
\lambda
\qff\bigr)
\off\qff =\off\qff
\bigl(\qff
\kernel\fff (\trf A_{\dff y}\qff -\qff \lambda \trf)\dff \otimes\dff \kernel\fff Q
\qff\bigr)
\qff \oplus\qff
0
\pff.
\]

\vspace{-12pt}\vspace{3.75pt}
for\dss
$\lambda\qff \in\qff
[\trf -\qff \varepsilon\dff,\pff \varepsilon\trf]$\nnsp.\oss
We see\sss that\sss for\sss $y\qff \in\qff U_{\dff z}$\sss the eigenspaces of\dss
$A_{\dff y}$\sss and\sss $A_{\dff y}\pff \omult_{\fff 1}\qff Q$\sss
corresponding\sss to sufficiently small\sss eigenvalues $\lambda$
are canonically\sss isomorphic.\oss

Let\sss us consider\sss the covering $U${\nsp} of\dss $Z$\sss by sufficiently small\sss
neighborhoods $U_{\dff z}\dff,\pff z\qff \in\qff Z$\nnsp.\oss
Let\sss $Z_{\dff U}$\sss be\sss Segal's\dss topological\sss category associated\sss
with\sss this covering.\oss
See\qss \cite{s1}\qss and\qss \cite{i2},\oss Section\qss 2.\oss
Finally,\oss let\sss $\hat{\mathcal{S}}_{\qff 1}$\sss and $\hat{\mathcal{S}}_{\qff 2}$\sss 
be\sss the versions of\dss the category $\hat{\mathcal{S}}$ of\dss subspaces of\dss
a\dss Hilbert\dss space\sss $H$\sss from\qss \cite{i1},\oss \cite{i2}\qss
related\dss to\sss the\dss Hilbert\dss spaces\sss
$H
\off =\off 
H_{\trf 0}$
and\sss
$H
\off =\off
(\trf H\qff\fff \widehat{\otimes}\pff K \trf)_{\dff 0}
\dff \oplus\dff
(\trf H\qff\fff \widehat{\otimes}\pff L \trf)_{\dff 0}$
respectively.\oss
The families\sss
$A_{\dff z}\dff,\off z\qff \in\qff Z$\sss and\sss
$A_{\dff z}\pff \omult_{\fff 1}\qff Q\dff,\off z\qff \in\qff Z$\sss
define index\sss functors\vspace{4.5pt}
\[
\quad
\mathbb{A}_{\qff U,\dff \bm{\varepsilon}}\qff \colon\qff
Z_{\dff U}
\qff \ttoo\qff
\hat{\mathcal{S}}_{\qff 1}
\quad
\mbox{and}\quad
\left(\trf \mathbb{A}\pff \omult_{\fff 1}\qff Q \trf\right)_{\qff U,\dff \bm{\varepsilon}}\qff \colon\qff
Z_{\dff U}
\qff \ttoo\qff
\hat{\mathcal{S}}_{\qff 2}
\pff.
\]

\vspace{-12pt}\vspace{4.5pt}
The canonical\sss isomorphisms between\sss the eigenspaces of\dss
$A_{\dff y}$\sss and\sss $A_{\dff y}\pff \omult_{\fff 1}\qff Q$\sss
show\sss that\sss these index\sss functors are\sss isomorphic.\oss
In\sss fact,\oss if\dss we identify\sss
$\hat{\mathcal{S}}_{\qff 1}$\sss with a subcategory of\dss $\hat{\mathcal{S}}_{\qff 2}$\sss 
by\sss the functor\sss 
$V
\off \longmapsto\off 
(\trf V\dff \otimes\dff \kernel\fff Q\trf)\dff \oplus\dff 0$\nnsp,\oss
then\sss these index\sss functors will\dss be equal.\oss
Then\sss their compositions with\sss the canonical\sss functors\qss
from\qss \cite{i1},\oss \cite{i2}\vspace{3pt}
\[
\quad
\hat{\mathcal{S}}_{\qff 1}
\qff \ttoo\qff Q
\quad
\mbox{and}\quad
\hat{\mathcal{S}}_{\qff 2}
\qff \ttoo\qff Q
\pff,
\]

\vspace{-12pt}\vspace{3pt}
also will\dss be equal,\oss
where\sss for a moment\sss we denoted\dss by $Q$ not\sss an operator but\sss the\dss
Quillen--Segal\qss \cite{s2}\qss category.\oss
Therefore,\sss the corresponding\sss index\sss maps in\sss the sense of\pss \cite{i2}\qss
are also equal\sss and,\oss in\sss particular,\oss homotopic.\oss
Hence\sss the indices of\dss our\sss two families are equal.\oss  \eproof

\newpage
\mysection{Pseudo-differential\qss operators\qss and\qss boundary\qss conditions}{pdo}

\myuppar{Operators.}
As in\dss Section\qss \ref{symbols-conditions},\oss
let\sss $X$\sss be a compact\sss riemannian\sss manifold\sss with\sss the non-empty\sss boundary\sss $Y$\dnsp,\oss
and\sss let\sss $E$\sss be a\dss Hermitian\dss vector bundle on\sss $X$\nnsp.\oss 
We will\sss keep\sss the notations and assumptions of\trs 
Section\qss \ref{symbols-conditions}\qss 
related\sss to $X\fff,\off Y$\dnsp.\oss
In\sss particular,\oss we identify\sss
a collar neighborhood of\dss $Y$\sss in\sss $X$\sss 
with\sss $Y\dff \times\dff [\trf 0\fff,\qff 1\dff)$\sss and
denote by\sss $x_{\dff n}$\sss the $[\trf 0\fff,\qff 1\dff)$\dnsp-coordinate in\sss the collar.\oss
This allows us\sss to define\sss the partial\sss derivative\sss
$\partial/\dff\partial\dff x_{\dff n}$\sss in\sss the collar.\oss
Let\sss $D_{\fff n}\off =\off -\qff i\trf \partial/\dff\partial\dff x_{\dff n}$\nsp.\oss
Let\sss us\sss fix a smooth\sss function\sss
$\varphi\dff \colon\dff 
(\trf -\qff 1\fff,\qff 1\dff)
\qff \ttoo\qff 
[\trf 0\fff,\qff 1\trf]$\sss
with compact\sss support\sss equal\sss to $1$ near\sss $0$\nnsp.\oss
We will\sss identify\sss the restriction of\dss $E$\sss to\sss the collar\sss
with\sss the pull-back of\dss $E\trf |\trf Y$\dnsp.\oss
Let\sss $X^{\fff \circ}\off =\off X\dff \smallsetminus\qff Y$\dnsp.\oss

We will\sss consider\sss operators of\dss order $1$ acting on sections of\sss $E$ 
and\sss belonging\sss to\sss the class introduced\dss by\trs 
H\"{o}rmander\qss \cite{h},\oss Chapter\dss 20.\oss
They have\sss the form\dss
$P\off =\off P^{\dff \mathrm{b}}\qff +\qff P^{\dff \mathrm{i}}$\dnsp,\oss
where\vspace{2.5pt}
\[
\quad
P^{\dff \mathrm{i}}
\off \in\off
\Psi\phg^{\dff 1}\trf
\left(\qff 
X^{\fff \circ}\dff;\off 
E\dff,\off 
E
\qff\right)
\qff
\]

\vspace{-12pt}\vspace{2.5pt}
and\dss the kernel\sss of\dss $P^{\dff \mathrm{i}}$ has compact\sss 
support\sss in\sss
$X^{\fff \circ}
\trf \times\trf 
X^{\fff \circ}$\dnsp,\oss
and\sss\vspace{2.5pt}
\[
\quad
P^{\dff \mathrm{b}}
\off =\off
\varphi\dff(\trf x_{\dff n}\trf)\qff
\bigl(\qff \bm{\Sigma}\trf(\trf x_{\dff n}\trf)\dff \circ\dff D_{\fff n}
\pff +\off
\mathbf{T}\trf(\trf x_{\dff n}\trf)
\qff\bigr)
\qff,
\]

\vspace{-12pt}\vspace{2.5pt}
where\sss $\bm{\Sigma}$\sss is\dss smooth map from $(\trf -\qff 1\fff,\qff 1\trf)$\sss to
endomorphisms of\dss $E\dff |\dff Y$\sss 
and\vspace{2.5pt}
\[
\quad
\mathbf{T}\dff \colon\dff
(\trf -\qff 1\fff,\qff 1\trf)
\qff \ttoo\pff
\Psi\phg^{\dff 1}\trf
\left(\qff 
Y\dff;\off 
E\dff |\dff Y\dff,\off 
E\dff |\dff Y
\qff\right)
\qff
\]

\vspace{-12pt}\vspace{2.5pt}
is\dss a smooth\sss function.\oss
Here\sss $\bm{\Sigma}\trf(\trf x_{\dff n}\trf)$\sss and\sss $\mathbf{T}\trf(\trf x_{\dff n}\trf)$\sss 
are understood as\sss acting\sss by\sss 
$\bm{\Sigma}\trf(\trf t\trf)$\sss and\sss $\mathbf{T}\trf(\trf t\trf)$\sss
respectively\sss over\dss 
$Y\dff \times\dff t$\nnsp,\qss $t\qff \in\qff [\trf 0\fff,\qff 1\dff)$\nnsp.\oss
Of\dss course,\pss $P$ depends only on\sss 
$\bm{\Sigma}\trf(\trf t\trf)\fff,\qff \mathbf{T}\trf(\trf t\trf)$\sss with\sss $t\qff \geq\qff 0$\nnsp.\oss
The\sss operators 
$\bm{\Sigma}\trf(\trf t\trf)\fff,\qff \mathbf{T}\trf(\trf t\trf)$\sss with $t\qff <\qff 0$
define an extension $\widehat{P}$ of\dss $P$\sss to\sss the union\sss 
$\widehat{X}\off =\off X\qff \cup\qff Y\dff \times\dff (\dff -\qff 1\fff,\qff 0\trf]$\nnsp,\oss
where\sss $Y\qff \subset\qff X$\sss is\dss identified\sss with\sss
$Y\dff \times\dff 0\qff \subset\qff Y\dff \times\dff (\dff -\qff 1\fff,\qff 0\trf]$\nnsp.\oss

The operator $\mathbf{T}\trf(\trf x_{\dff n}\trf)$ is\dss not\sss pseudo-differential,\oss
but\sss belongs\sss to\sss the closure\sss $\overline{\mathcal{P}^{\trf 1}}$\sss
of\dss the space\sss $\mathcal{P}^{\trf 1}$ 
of\dss pseudo-differential\sss operators of\dss order $1$\nnsp.\oss 
See\qss \cite{as1},\oss Section\qss 5,\oss for\sss the definitions of\dss these spaces.\oss
In\sss particular,\oss the principal\sss symbol\sss of\dss 
$\mathbf{T}\trf(\trf x_{\dff n}\trf)$\sss is\dss well\sss defined
and\sss $\mathbf{T}\trf(\trf x_{\dff n}\trf)$\sss can be approximated\dss in  by pseudo-differential\sss operators
with principal\sss symbols approximating\sss that\sss of\dss $\mathbf{T}\trf(\trf x_{\dff n}\trf)$\nnsp.\oss
In\sss the simplest\sss case,\oss when\sss the function\sss $\mathbf{T}$\sss
is\dss constant,\pss $\mathbf{T}\trf(\trf x_{\dff n}\trf)$\sss is\sss equal\sss to\sss
the\sss lift\sss of\dss $\mathbf{T}\trf(\trf 0\trf)$ from $Y$ to\sss
$Y\dff \times\dff [\trf 0\fff,\qff 1\dff)$\sss in\sss the sense of\pss
\cite{as1},\oss Section\qss 5.\oss
See also\sss the discussion of\trs lifts in\dss Section\qss \ref{mult-pdo}.\oss
For\sss the general\sss case see\qss \cite{h},\oss Sections\qss 19.2\qss and\dss 20.1.\oss

The principal\sss symbol\sss of\dss $P$\sss
is\dss equal\dss to\sss the sum of\dss the principal\sss symbols of\trs
$P^{\dff \mathrm{b}}$ and\sss $P^{\dff \mathrm{i}}$\dnsp.\oss
We need only\sss the restriction of\dss the principal\sss
symbol\sss to\sss $S\dff X$\nnsp,\oss and will\sss call\dss this restriction\sss
simply\qss \emph{the symbol\dss of}\trs $P$\dnsp.\oss
It\sss is\dss a bundle map 
$\sigma\dff \colon\dff \pi^{\fff *}\dff E\qff \ttoo\qff \pi^{\fff *}\dff E$\nnsp.\oss
Let\sss $\Sigma\off =\off \bm{\Sigma}\trf(\trf 0\trf)$\nnsp.\oss
For\sss $u\qff \in\qff S\dff Y$\sss and\sss $y\off =\off \pi\trf(\dff u\trf)$\sss
we denote by\sss $\sigma_y$ and\sss 
$\tau_{\fff u}$\sss the maps\sss 
$E_{\dff y}\qff \ttoo\qff E_{\dff y}$\sss 
induced\sss by\sss $\Sigma$\sss and\sss the symbol\sss of\dss
$\mathbf{T}\trf(\trf 0\trf)$\sss respectively.\oss
Let\sss 
$\rho_{\fff u}
\off =\off 
\sigma_y^{\dff -\dff 1}\dff \tau_u$\nnsp.\oss
As usual,\pss $P$ and $\sigma$ are said\sss be\qss \emph{elliptic}\pss
if\sss $\sigma$\sss is\dss an automorphism of\sss $\pi^{\fff *}\dff E$\nnsp,\oss
and $\sigma$\sss is\dss said\sss to be\qss 
\emph{self-adjoint}\pss if\dss $\sigma$\sss is\dss self-adjoint\sss in every\sss fiber.\oss

We will\sss say\sss that\sss $P$\sss is\qss \emph{self-adjoint}\pss
if\dss $P$\sss is\dss symmetric on sections with support\sss in $X^{\fff \circ}$\dnsp,\pss
$\bm{\Sigma}\trf(\trf t\trf)$\sss is\dss a self-adjoint\dss bundle map,\oss 
and\sss $\mathbf{T}\trf(\trf t\trf)$\sss is\dss a self-adjoint\sss operator\sss 
for every\sss
$t\qff \in\qff [\trf 0\fff,\qff 1\trf)$\nnsp.\oss
Then\sss $\sigma$\sss is\dss self-adjoint\sss
and\dss the action of\dss $P$\sss on\sss
smooth sections satisfies\sss the\dss Lagrange\dss identity\vspace{3pt}
\begin{equation}
\label{green-formula}
\quad
\sco{\dff P\dff u\fff,\qff v \dff}_{\trf X}
\qff -\qff
\sco{\dff u\dff,\qff P\dff v \dff}_{\trf X}
\off =\off
\sco{\dff i\trf \Sigma\trf \gamma\dff u\dff,\qff \gamma\dff v \dff}_{\trf \partial}
\qff,
\end{equation}

\vspace{-12pt}\vspace{3pt}
where\sss $\sco{\dff \bullet\fff,\qff \bullet\dff}_{\trf X}$\sss and\sss
$\sco{\dff \bullet\dff,\qff \gamma\dff \bullet\dff}_{\trf \partial}$\sss
are\sss the\sss $L_{\dff 2}$ scalar products of\dss sections 
defined\sss by\sss the\dss Hermitian\dss structure,\oss
and\sss $\gamma$\sss is\dss the\qss \emph{trace operator}\pss 
taking a section over $X$\sss to its restriction\sss to $Y$\dnsp.\oss
The identity\qss (\ref{green-formula})\qss is\dss known 
as\sss the\qss \emph{Green\dss formula}\qss and
follows from\sss the integration by\sss parts.\oss

\myuppar{Boundary\sss operators for\sss general\sss elliptic operators of\dss order $1$\nnsp.}
Let\sss $G$\sss be a vector bundle over\sss $Y$\dnsp.\oss
A\qss \emph{boundary operator}\pss maps sections of\sss $E$ over $X$\sss 
to sections of\sss $G$ over $Y$\dnsp.\oss 
H\"{o}rmander\qss \cite{h},\oss Chapter\dss 20,\oss considers boundary operators of\dss the form\sss
$B\off =\off B_{\trf Y}\dff \circ\trf \gamma$\nnsp,\oss 
where\sss $B_{\trf Y}$\sss is\dss a pseudo-differential\sss operator of\sss order $0$ 
from\sss $E\dff |\dff Y$\sss to\sss $G$\nsp.\oss

Let\sss $P$\dnsp,\dss $\sigma$ and\sss $B$ be as above,\oss 
and suppose\sss that\sss $P$\dnsp, $\sigma$ are elliptic.\oss
For\sss every\sss $u\qff \in\qff S\dff Y$\sss and\sss $y\off =\off \pi\trf(\dff u\trf)$\sss
the principal\sss symbol\sss of\dss $B_{\trf Y}$\sss defines a\sss linear\sss map\dss
$\beta_{\fff u}\dff \colon\dff
E_{\dff y}\qff \ttoo\qff G_{\dff y}$\nsp.\oss
Recall\dss that\sss 
$\mathcal{M}_{\dff +}\dff(\trf \rho_{\dff u}\dff)
\qff \subset\pff E_{\dff y}$\sss
is\dss the space of\dss bounded solutions\sss 
$f\dff \colon\dff \rrr_{\qff \geq\qff 0}\qff \ttoo\qff E_{\dff y}$ of\dss
the ordinary\sss differential\sss equation\sss
$D\trf(\trf f\trf)\qff +\qff \rho_{\fff u}\trf (\trf f\trf)\off =\off 0$\nnsp,\oss
where\sss $D\off =\off -\qff i\qff \partial$\sss
and\sss $y\off =\off \pi\trf(\dff u\trf)$\nnsp.\oss
The boundary\sss operator $B$\sss
satisfies\qss 
\emph{Shapiro-Lopatinskii\sss condition\dss for}\dss $P$
if\trs for every\sss $u\qff \in\qff S\dff Y$\sss the map\sss $\beta_{\fff u}$
induces an\sss isomorphism of\dss the space of\dss the 
initial\sss values\sss $f\dff(\trf 0\trf)$\sss of\dss solutions in
$\mathcal{M}_{\dff +}\dff(\trf \rho_{\dff u}\dff)$\sss
with\sss $G_{\dff y}$\nsp.\oss
In\sss this case one also says\sss that $B$ is\qss 
\emph{elliptic}\pss for $P$ or\sss that\sss the pair $(\trf P\fff,\qff B\trf)$
is\qss \emph{elliptic}.\oss

As explained\sss in\dss Section\qss \ref{boundary-algebra},\oss
the space of\dss the initial\sss values of\dss solutions in\sss
$\mathcal{M}_{\dff +}\dff(\trf \rho_{\dff u}\dff)$\sss
is\dss equal\sss to\sss 
$\mathcal{L}_{\dff -}\dff(\trf \rho_{\dff u}\dff)\qff \subset\pff E_{\dff y}$\nsp,\oss
the sum of\dss the generalized eigenspaces of\dss $\rho_{\dff u}$\sss corresponding\sss
to\sss eigenvalues $\lambda$\sss with\sss $\im\dff \lambda\qff <\qff 0$\nnsp.\oss
Therefore\sss the\dss Shapiro-Lopatinskii\sss condition\dss holds\sss
if\dss and\dss only\dss if\trs for every\sss $u\qff \in\qff S\dff Y$\sss
the map\sss $\beta_{\fff u}$\sss induces an\sss isomorphism\sss
$\mathcal{L}_{\dff -}\dff(\trf \rho_{\fff u}\dff)
\qff \ttoo\qff
G_{\dff y}$\nsp,\oss
or,\oss equivalently,\pss
the kernel\dss
$N_{\fff u}\off =\off \kernel\dff \beta_{\dff u}$\sss
is\dss a complementary subspace\sss to\sss $\mathcal{L}_{\dff -}\dff(\trf \rho_{\fff u}\dff)$\nnsp.\oss
The kernels\sss $N_{\fff u}$\sss are fibers of\dss a vector bundle\sss $N$\sss over\sss
$S\dff Y$\dnsp,\oss which we will\sss call\dss the\qss 
\emph{kernel-symbol}\pss of\dss $B_{\trf Y}$ and\sss $B$\nnsp.\oss

\myuppar{Fredholm\dss property.}
For every\sss $s\qff \geq\qff 1$\sss the rule\sss
$u\off \longmapsto\off (\trf P\fff u\fff,\qff B\fff u\trf)$
defines a map\vspace{3pt}
\begin{equation*}
\quad
P\dff \oplus\dff B\qff \colon\qff
H_{\dff s}\trf(\trf X^{\fff \circ}\fff,\qff E\trf)
\off \ttoo\off
H_{\dff s\dff -\dff 1}\trf(\trf X^{\fff \circ}\fff,\qff E\trf)
\qff \oplus\qff
H_{\dff s\dff -\dff 1/2}\trf(\trf Y\fff,\qff G\trf)
\end{equation*}

\vspace{-12pt}\vspace{3pt}
of\trs Sobolev spaces.\oss 
If\dss $(\trf P\fff,\qff B\trf)$\sss is\dss elliptic,\oss
then $P\dff \oplus\dff B$\sss is\dss Fredholm,\oss
the kernel\sss of\dss $P\dff \oplus\dff B$\sss consists of\sss
$C^{\dff \infty}$\dnsp-smooth sections,\oss 
and\sss the image of\dss $P\dff \oplus\dff B$\sss is\dss the orthogonal\sss
complement\sss of\dss a subspace consisting of\dss pairs of\sss
$C^{\dff \infty}$\dnsp-smooth sections.\oss
See\dss H\"{o}rmander\qss \cite{h},\oss Theorem\dss $20.1.8\fff'$\nsp\dnsp.\oss
H\"{o}rmander\sss works with\dss Sobolev\dss spaces\sss
$\overline{H}_{\dff s}\trf(\trf X^{\fff \circ}\fff,\qff E\trf)$
of\dss sections extendable\sss to\sss $\widehat{X}$\sss with a quotient\sss norm,\oss 
but\sss they are isomorphic\sss to 
$H_{\dff s}\trf(\trf X^{\fff \circ}\fff,\qff E\trf)$\nnsp.\oss
See,\pss for example,\pss \cite{ag},\oss Theorem\qss 3.1.1.

\myuppar{Bundle-like\sss boundary conditions in\sss the self-adjoint\sss case.}
Suppose\sss that\sss the pair\sss $(\trf P\fff,\qff B\trf)$\sss is\dss elliptic
and\sss $P$\sss is\dss self-adjoint.\oss
Then $\sigma$\sss is\dss an elliptic self-adjoint\sss symbol\sss of\dss order $1$
in\sss the sense of\qss Section\qss \ref{symbols-conditions}.\oss
Suppose\sss that\sss the boundary operator 
$B\off =\off B_{\trf Y}\dff \circ\trf \gamma$\sss 
is\qss \emph{bundle-like}\pss in\sss
the sense\sss that\sss $B_{\trf Y}$\sss is\dss
induced\sss by\sss a\sss bundle morphism\sss 
$E\trf |\trf Y\qff \ttoo\qff G$\nnsp.\oss
Then\sss the kernel-symbol\sss $N$\sss of\dss $B_{\trf Y}$\sss
is\dss a bundle-like boundary conditions for\sss the symbol $\sigma$
and\sss may\sss be considered as a subbundle of\sss $E\trf |\trf Y$\dnsp.\oss
The kernel\sss of\dss $B$\sss on\sss
$H_{\dff s}\dff(\trf X^{\fff \circ}\fff,\qff E\trf)$\sss
consists of\dss sections\sss $u$\sss such\sss that\dss
$\gamma\dff u
\qff \in\qff\fff
H_{\dff s\dff -\dff 1/2}\dff(\trf Y\fff,\qff N\trf)$\nnsp.\oss
Since only\sss the kernel\sss of\dss $B$\sss matters,\oss
we may assume\sss that\sss $G$\sss is\dss equal\dss to\sss the orthogonal\sss
complement\sss $N^{\dff \perp}$ of\dss $N$\sss in\sss $E\trf |\trf Y$\sss
and\dss the operator\sss $B_{\trf Y}$\sss is\dss induced\sss by\sss the orthogonal\sss projection\sss
$\pi_{\dff N^{\dff \perp}}\dff \colon\dff E\trf |\trf Y\qff \ttoo\qff N^{\dff \perp}$\dnsp.\oss 
Let\sss us\sss take\sss $s\off =\off 1$\sss and\dss let\vspace{1.75pt}
\[
\quad
\Pi\dff \colon\dff
H_{\trf 0}\dff(\trf Y\fff,\qff E\trf |\trf Y\trf)
\off \ttoo\off
H_{\trf 0}\dff(\trf Y\fff,\qff E\trf |\trf Y\trf)
\quad
\mbox{and}\quad
\]

\vspace{-36pt}
\[
\quad
\Pi_{\dff 1/2}\dff \colon\dff
H_{\dff 1/2}\dff(\trf Y\fff,\qff E\trf |\trf Y\trf)
\off \ttoo\off
H_{\dff 1/2}\dff(\trf Y\fff,\qff E\trf |\trf Y\trf)
\]

\vspace{-12pt}\vspace{1.75pt}
be\sss the operators\sss in\dss Sobolev\dss spaces induced\sss by\sss the orthogonal\dss projection\sss 
$\pi_{\dff N}\dff \colon\dff E\trf |\trf Y\qff \ttoo\qff N$\nsp.\oss
Then we can replace\sss $B_{\trf Y}$\sss by\sss $1\qff -\qff \Pi_{\dff 1/2}$\sss
and\sss $B$\sss by\sss 
$\Gamma\off =\off (\trf 1\qff -\qff \Pi_{\dff 1/2}\trf)\dff \circ\dff \gamma$\sss
without\sss affecting\sss the kernels.\oss
If\dss $N$\sss is\dss a self-adjoint\sss boundary condition for\sss $\sigma$\nnsp,\oss
then\sss $\Pi$\sss is\dss a self-adjoint\sss boundary condition\sss for\sss
the operator\sss $P$ acting on\sss the above\dss Sobolev\dss spaces.\oss
In\sss more details,\oss let\vspace{1.75pt}
\[
\quad
H_{\trf 0}
\off =\off
H_{\trf 0}\dff(\trf X^{\fff \circ}\fff,\qff E \trf)
\qff,\quad\off
H_{\dff 1}
\off =\off
H_{\dff 1}\dff(\trf X^{\fff \circ}\fff,\qff E \trf)
\qff,
\]

\vspace{-36pt}
\[
\quad
H^{\dff \partial}
\off =\off
H_{\trf 0}\dff(\trf Y\fff,\qff E\trf |\trf Y\trf)
\qff,\quad\off
H_{\dff 1/2}^{\dff \partial}
\off =\off
H_{\dff 1/2}\dff(\trf Y\fff,\qff E\trf |\trf Y\trf)
\]

\vspace{-12pt}\vspace{1.75pt}
and\sss
$\gamma\dff \colon\dff
H_{\fff 1}\qff \ttoo\qff H_{\dff 1/2}^{\dff \partial}$\sss
be\sss the\sss trace operator.\oss
Now we are in\sss the framework of\trs Section\qss \ref{abstract-index}.\oss

The spaces $H_{\dff 0}$ and\sss $H^{\dff \partial}$ are nothing else but\sss the\dss Hilbert\dss spaces
of\dss $L_{\dff 2}$ sections,\oss
and\dss Green\dss formula\qss (\ref{green-formula})\qss
extends by continuity\sss to\sss sections\sss
$u\fff,\qff v\qff \in\qff H_{\fff 1}$\nsp.\oss
Therefore\qss (\ref{green-formula})\qss can\sss be interpreted as an
abstract\dss Lagrange\dss identity\sss with\dss the operator $P$\dss in\sss the role of\dss $A$\nnsp,\oss
the operator\dss $\Sigma\off =\off \bm{\Sigma}\trf(\trf 0\trf)$\dss
in\sss the role of\dss $\Sigma$\sss and\sss the\sss trace operator\sss $\gamma$\sss
in\sss the role of\dss $\gamma$\sss from\qss (\ref{lagrange}).\oss
If\dss $N$\sss is\dss a self-adjoint\sss boundary condition for\sss $\sigma$\nnsp,\oss then\dss 
$\Sigma\trf(\trf N\trf)
\off =\off 
i\trf \Sigma\trf(\trf N\trf)$\sss 
is\dss orthogonal\dss to\sss $N$\sss and\sss hence\sss
$i\trf \Sigma\trf(\trf \image\dff \Pi\trf)
\off =\off
\kernel\dff \Pi$\nnsp.\oss
Therefore,\oss in\sss this case\sss
$\Pi$\sss is\dss a self-adjoint\sss boundary condition for\sss $P$\sss
in\sss the sense of\trs Section\qss \ref{abstract-index}.\oss
Moreover,\oss Theorem\dss $20.1.8\fff'$\dss in\qss \cite{h}\qss
implies\sss that\sss the boundary\sss problem\sss $P\fff,\qff \Pi$\sss is\dss regular\sss
and\dss the corresponding operator\sss $P\dff \oplus\dff \Gamma$\sss is\dss Fredholm.\oss
As\dss is\dss well\dss known,\oss the inclusion\sss
$H_{\fff 1}\qff \ttoo\qff H_{\dff 0}$\sss is\dss compact.\oss
Therefore,\oss all\sss assumptions of\trs Section\qss \ref{abstract-index}\qss hold.\oss
Hence\sss if\dss $N$\sss is\dss a self-adjoint\dss bundle-like boundary condition,\oss
then\dss Theorem\qss \ref{fredholm}\qss applies.\oss
Hence\sss the restriction\sss $P_{\dff \Gamma}\off =\off P_{\dff B}$\sss
of\dss $P$\sss to\sss $\kernel\dff \Gamma\off =\off \kernel\dff B$\sss
is\dss self-adjoint.\oss
Moreover,\oss it\dss is\dss an\sss unbounded\dss Fredholm\dss operator\sss with discrete spectrum
and compact\sss resolvent.

\myuppar{Families.}
Suppose now\sss that\sss our operators and\sss boundary conditions 
continuously depend on a parameter $z\qff \in\qff Z$\nnsp,\oss
where $Z$\sss is\dss a\sss topological\sss space.\oss
Moreover,\oss as at\sss the end of\trs Section\qss \ref{symbols-conditions},\oss
let\sss us allow\sss the manifold\sss $X$\sss to
continuously depend on $z$\nnsp,\oss in\sss the sense\sss that\sss
$X\off =\off X\trf(\trf z\trf)$\sss is\dss the fiber over $z$ 
of\dss a\sss locally\sss trivial\dss bundle over $Z$\nnsp.\oss
Let\sss $Y\off =\off Y\trf(\trf z\trf)$ be\sss the boundary of\dss $X\trf(\trf z\trf)$\nnsp.\oss
In order\sss to be able\sss to speak\sss about\sss families of\dss operators
in\sss the\dss H\"{o}rmander\dss class,\oss we need\dss to assume\sss that\sss
a family\sss of\dss embeddings\sss 
$c\trf(\trf z\trf)\dff \colon\dff
Y\trf(\trf z\trf)\dff \times\dff [\trf 0\fff,\qff 1\trf)
\qff \ttoo\qff
X\trf(\trf z\trf)$\sss
is\dss fixed,\oss such\sss that\sss
$c\trf(\trf z\trf)\dff(\trf y\fff,\qff 0\trf)\off =\off y$\sss if\dss
$y\qff \in\qff Y\trf(\trf z\trf)$\nnsp.\oss
It\dss is\dss well\sss known\sss that\sss for reasonable spaces $Z$ such a family exists and\dss
is\dss unique up\sss to isotopy.\oss
The embeddings $c\trf(\trf z\trf)$ are used\sss for\sss the identification of\dss
collar neighborhoods of\dss boundaries $Y\trf(\trf z\trf)$ with cylinders.\oss
Of\dss course,\oss we assume\sss that\sss the operators 
$P\off =\off P\trf(\trf z\trf)$ and\sss $B\off =\off B\trf(\trf z\trf)$\sss
continuously depend on\sss parameters.\oss
By\sss the continuous dependence of\dss $P$ on $z$ we understand\dss
the continuous dependence of\dss elements defining\sss $P$\dnsp,\oss
i.e.\qss of\dss the operators\sss $P^{\dff \mathrm{i}}$\dnsp,\oss of\dss the maps\sss
$\bm{\Sigma}\fff,\pff \mathbf{T}$\dnsp,\oss and of\dss functions $\varphi$\nnsp.\oss
We will\sss consider only\sss bundle-like boundary conditions with\sss the obvious
notion of\dss continuous dependence on $z$\nnsp.\oss

Locally over $Z$ we may assume\sss that\sss 
manifolds,\pss bundles,\pss and collars are independent\sss of\dss $z$\nnsp.\oss
Since\sss a bundle-like boundary condition $B$\sss is\dss determined\dss by\sss the corresponding\sss
bundle $N$\nnsp,\oss locally we can also assume\sss that\sss $B$\sss is\dss independent\sss of\dss $z$\nnsp.\oss 
As in\qss \cite{as4},\oss this implies\sss that\sss locally over $Z$\sss the operators\sss
$P\dff \oplus\dff B
\off =\off
P\trf(\trf z\trf)\dff \oplus\dff B\trf(\trf z\trf)$\sss
continuously depend on $z$ in\sss the norm\sss topology.\oss
As explained at\dss the end of\trs Section\qss \ref{abstract-index},\oss
Lemma\qss \ref{continuity}\qss implies\sss that\sss the family of\dss unbounded operators\sss
$P_{\dff B}
\off =\off
P\trf(\trf z\trf)_{\dff B\trf(\trf z\trf)}\dff,\pff
z\qff \in\qff Z$\sss
is\dss a\dss Fredholm\dss family of\dss self-adjoint\sss operators

\myuppar{Elementary operators and\dss boundary conditions.}
Let\sss us review\sss the construction of\dss operators from\dss 
Proposition\qss 20.3.1\qss in\qss \cite{h}.\oss
Suppose\sss that\sss $E\dff|\dff \partial\dff X$\sss is\dss presented
as\sss an orthogonal\sss direct\sss sum\sss
$E\dff|\dff \partial\dff X
\off =\off
E^{\dff +}\dff \oplus\dff E^{\dff -}$\dnsp.\oss
This\sss leads\sss to a similar\sss decomposition of\dss $E$\sss over\sss the collar\sss
$\partial\dff X\dff \times\dff [\trf 0\fff,\qff 1\dff)$\nnsp.\oss
By an abuse of\dss notations we denote\sss the summands of\dss the\sss latter also
by\sss $E^{\dff +},\off E^{\dff -}$\dnsp.\oss
Let\sss 
$\varphi$ be a function such as in\dss Sections\qss \ref{symbols-conditions}.\oss
Let\sss $\lambda^{\fff +}\nsp,\off \lambda^{\fff -}\nsp,\off \lambda\qff >\qff 0$\sss
be positive real\sss numbers and\vspace{3pt}
\[
\quad
\Lambda^{\fff +}\qff \in\qff \Psi\phg^{\dff 1}\trf(\trf \partial\dff X\fff,\qff E^{\dff +}\fff,\qff E^{\dff +}\trf)
\dff,\quad
\Lambda^{\fff -}\qff \in\qff \Psi\phg^{\dff 1}\trf(\trf \partial\dff X\fff,\qff E^{\dff -}\fff,\qff E^{\dff -}\trf)
\dff,\quad
\mbox{and}\quad
\Lambda\qff \in\qff \Psi\phg^{\dff 1}\trf(\trf X^{\fff \circ}\fff,\qff E\fff,\qff E\trf)
\]

\vspace{-12pt}\vspace{3pt}
be pseudo-differential\sss operators
with principal\sss symbols\sss 
$\lambda^{\fff +}\nsp,\off \lambda^{\fff -}$\dnsp,\oss and\sss $\lambda$\sss times\sss
the identity\sss respectively.\oss
Suppose\sss that\sss $\Lambda^{\fff +}\fff,\qff \Lambda^{\fff -}$\sss are self-adjoint\sss
and\sss $\Lambda$\sss is\dss formally self-adjoint.\oss
For a section $u$ of\sss $E$\sss we can represent\sss its
restriction\sss to\sss the collar as\sss $(\trf u^{\dff +},\qff u^{\dff -}\trf)$\nnsp,\oss
where $u^{\dff +},\qff u^{\dff -}$ are sections of\dss $E^{\dff +},\off E^{\dff -}$\sss
respectively.\oss
Let\sss $P^{\dff \mathrm{b}}$\sss be\sss the operator
acting on sections of\dss $E$\sss by\vspace{1.5pt}
\[
\quad
P^{\dff \mathrm{b}}\fff u
\off =\off
\varphi\dff(\trf x_{\dff n}\trf)\qff
\left(\qff
(\trf D_{\fff n}\qff +\qff i\dff \Lambda^{\fff +}\trf)\dff u^{\dff +},\off
(\trf -\qff  D_{\fff n}\qff +\qff i\dff \Lambda^{\fff -}\trf)\dff u^{\dff -}
\qff\right)
\qff.
\]

\vspace{-12pt}\vspace{1.5pt}
Let\dss $Q^{\dff \mathrm{i}}$\sss be\sss the operator\sss acting
on sections of\dss $E$\sss by\vspace{1.5pt}
\[
\quad
Q^{\dff \mathrm{i}}\fff u
\off =\off
(\dff 1\qff -\qff \varphi\dff(\trf x_{\dff n}\trf)\trf)\dff 
\Lambda\dff
\left(\qff (\dff 1\qff -\qff \varphi\dff(\trf x_{\dff n}\trf)\trf)\dff u\qff \right)
\qff,
\]

\vspace{-12pt}\vspace{1.5pt}
and\sss let\sss $P^{\dff \mathrm{i}}\off =\off i\trf Q^{\dff \mathrm{i}}$\dnsp.\oss
Finally,\oss let\sss 
$P\off =\off P^{\dff \mathrm{b}}\qff +\qff P^{\dff \mathrm{i}}$\sss 
and\sss
$B\fff u
\off =\off
u^{\dff -}\dff |\dff \partial\dff X$\nnsp.\oss
It\sss is\dss easy\sss to check\sss that\sss $B$\sss
satisfies\sss Shapiro-Lopatinskii\dss condition\sss with respect\sss to\sss $P$\dnsp,\oss
i.e.\qss the pair $(\trf P\fff,\qff B\trf)$\sss is\dss  elliptic.\oss 
By\qss \cite{h},\oss Proposition\qss 20.3.1,\oss
the operator\sss $P\dff \oplus\dff B$\sss induces\dss Fredholm\dss operators
of\dss index zero in Sobolev\sss spaces.\oss
The following\dss lemma\dss
is\dss a part\sss of\dss the proof\dss of\dss this\dss proposition.\oss

\mypar{Lemma.}{injectivity}
\emph{If\qss $t\fff,\qff C\qff \in\qff \rrr$\dss are sufficiently\sss large,\oss
then}\vspace{0.75pt}
\begin{equation}
\label{h-eq}
\quad
\image\trf \sco{\dff P\dff u\dff,\qff u \dff}_{\dff X}
\off\dff \geq\off
-\qff C\dff
\sco{\dff u\dff,\qff u \dff}_{\dff X}
\off,
\end{equation}

\vspace{-12pt}\vspace{0.75pt}
\emph{for every\sss $u$\sss such\sss that\trs
$B\dff u\off =\off \gamma\dff u^{\dff -}\off =\off 0$\nnsp,\oss
and\sss the operator\trs
$(\trf P\qff +\qff i\fff t\trf \id\trf)\dff \oplus\dff B$\dss
is\dss injective.\oss}

\proof
For\sss the proof\dss of\dss the inequality\qss (\ref{h-eq})\qss
for sufficiently\sss large $C$ and every $u$\sss such\sss that\trs
$B\dff u\off =\off \gamma\dff u^{\dff -}\off =\off 0$ 
see\qss \cite{h},\oss the proof\dss of\trs Proposition\qss 20.3.1.\oss
This inequality\sss implies\sss that\sss\vspace{3pt}
\[
\quad 
\image\trf \sco{\dff (\trf P\qff +\qff i\fff t\trf \id\trf)\dff u\dff,\qff u \dff}_{\dff X}
\off \geq\off
\sco{\dff u\dff,\qff u \dff}_{\dff X}
\]

\vspace{-12pt}\vspace{3pt}
if\dss
$t\qff >\qff C\qff +\qff 1$\sss and\sss 
$B\dff u\off =\off \gamma\dff u^{\dff -}\off =\off 0$\nnsp.\oss 
The injectivity\sss follows.\oss  \eproof

\myuppar{Self-adjoint\sss operators.}
The operators\sss
$-\qff P^{\dff \mathrm{i}},\off -\qff \varphi\dff(\trf x_{\dff n}\trf)\trf i\dff \Lambda^{\fff +}$\dnsp,\oss
and\sss $-\qff \varphi\dff(\trf x_{\dff n}\trf)\trf i\dff \Lambda^{\fff -}$\sss
are formally\sss adjoint\sss to\sss 
$P^{\dff \mathrm{i}},\off \varphi\dff(\trf x_{\dff n}\trf)\trf i\dff \Lambda^{\fff +}$\dnsp,\oss
and\sss $\varphi\dff(\trf x_{\dff n}\trf)\trf i\dff \Lambda^{\fff -}$\sss respectively.\oss
The operator\sss formally adjoint\sss to\sss
$\varphi\dff(\trf x_{\dff n}\trf)\trf D_{\fff n}$\sss
is\dss the operator\sss 
$u\off \longmapsto\off D_{\fff n}\dff(\trf \varphi\dff(\trf x_{\dff n}\trf)\dff u\trf)$\nnsp,\oss
by\dss the\dss Leibniz\trs formula equal\dss to\vspace{3pt}
\[
\quad
u\off \longmapsto\off -\qff i\trf \varphi\fff'\dff(\trf x_{\dff n}\trf)\dff u
\qff +\qff \varphi\dff(\trf x_{\dff n}\trf)\trf D_{\fff n}\dff u
\qff.
\]

\vspace{-12pt}\vspace{3pt}
Let\sss us define\sss the operator\sss$\widetilde{P}^{\trf \mathrm{b}}$
by\vspace{3pt} 
\[
\quad
\widetilde{P}^{\trf \mathrm{b}}\dff u
\off =\off 
\varphi\dff(\trf x_{\dff n}\trf)\qff
\left(\qff
(\trf -\qff
D_{\fff n}\qff +\qff i\dff \Lambda^{\fff +}\trf)\dff u^{\dff +},\off
(\trf D_{\fff n}\qff +\qff i\dff \Lambda^{\fff -}\trf)\dff u^{\dff -}
\qff\right)
\qff
\]

\vspace{-12pt}\vspace{3pt}
and define\sss the operators\dss $\widetilde{P}$\sss and\dss $\widetilde{P}\fff'$\dss 
by\vspace{3pt}
\[
\quad
\widetilde{P}
\off =\off
\widetilde{P}^{\trf \mathrm{b}}
\qff +\qff 
P^{\dff \mathrm{i}}
\quad
\mbox{and}\quad
\widetilde{P}\fff'\dff u
\off =\off
\widetilde{P}\dff u
\qff +\qff 
i\trf \varphi\fff'\dff(\trf x_{\dff n}\trf)\dff u^{\dff +}
\qff -\qff 
i\trf \varphi\fff'\dff(\trf x_{\dff n}\trf)\dff u^{\dff -}
\qff.
\]

\vspace{-12pt}\vspace{3pt}
Then\sss $-\qff \widetilde{P}\fff'$\sss is\dss formally\sss adjoint\sss to\sss $P$\nnsp.\oss
The operators\sss $\widetilde{P}^{\trf \mathrm{b}}$\sss and\sss $\widetilde{P}\fff'$\sss
are operators of\dss the same nature as\sss
$P^{\dff \mathrm{b}}$ and\sss $P$\nnsp,\oss but\sss with\sss the roles of\sss
$E^{\dff +}$ and\sss $E^{\dff -}$\sss interchanged.\oss
Therefore\trs Lemma\qss \ref{injectivity}\qss implies\sss that\dss 
the operator\trs
$(\trf \widetilde{P}\fff'\qff +\qff i\fff t\trf \id\trf)\dff \oplus\dff \widetilde{B}$\nnsp,\oss
where\sss
$\widetilde{B}\dff u\off =\off \gamma\dff u^{\dff +}$\dnsp,\oss
is\dss injective.\oss
It\sss follows\sss that\sss
$(\trf -\qff \widetilde{P}\fff'\qff -\qff i\fff t\trf \id\trf)\dff \oplus\dff \widetilde{B}$\sss
is\dss also injective.\oss
Let\sss us consider\sss the operator\vspace{1.5pt}
\[
\quad
P^{\dff \sa}
\off =\off\dff
\begin{pmatrix}
\off\dff 0 &
P\qff +\qff i\fff t\trf \id \qff\off
\vspace{6pt} \\
\off\qff -\qff \widetilde{P}\fff'\qff -\qff i\fff t\trf \id &
0 \dff\off 
\end{pmatrix}
\qff
\]

\vspace{-12pt}\vspace{1.5pt}
acting\sss in\sss $E\dff \oplus\dff E$\nnsp.\oss
It\dss is\dss formally self-adjoint\sss
and\dss together\sss with\sss the boundary operator\sss\vspace{3pt}
\[
\quad
B^{\dff \sa}\dff \colon\dff
(\dff u\fff,\qff v\trf)
\off \longmapsto\off
(\trf \gamma\trf v^{\dff -}\fff,\qff \gamma\trf u^{\dff +} \trf)
\qff,
\]

\vspace{-12pt}\vspace{3pt}
where\sss $(\trf u\fff,\qff v\trf)\qff \in\qff E\dff \oplus\dff E$\nnsp,\oss
satisfies\dss Shapiro--Lopatinskii\dss condition.\oss
The kernel\sss $N^{\dff \sa}$\sss of\sss $B^{\dff \sa}$\sss
is\dss equal\dss to\sss
$(\trf 0\dff \oplus\dff E^{\dff -}\trf)
\qff \oplus\qff
(\trf E^{\dff +}\dff \oplus\dff 0\trf)$\nnsp.\oss
By\sss the construction,\oss the symbol\sss of\dss $P^{\dff \sa}$\sss
is\dss an elementary symbol,\oss and\sss $N^{\dff \sa}$\sss
is\dss the corresponding\sss boundary condition.\oss
Hence\sss the\sss topological\dss index of\dss the boundary\sss problem\sss
$P^{\dff \sa}\fff,\off N^{\dff \sa}$\sss is\dss zero.\oss
The boundary operator\sss $B^{\dff \sa}$\sss is\dss bundle-like,\oss
and\dss the boundary\sss problem\sss
$P^{\dff \sa}\fff,\off B^{\dff \sa}$\sss
is\dss of\dss the\sss type considered\sss above.\oss
Therefore\sss the results of\qss Section\qss \ref{abstract-index}\qss apply.\oss
Since\sss the operator\sss 
$P^{\dff \sa}\dff \oplus\dff B^{\dff \sa}$\sss is\trs Fredholm,\oss
Theorem\qss \ref{duality-lemma}\qss together\sss with\dss Lemma\qss \ref{injectivity}\qss imply\sss 
that\sss $P^{\dff \sa}\dff \oplus\dff B^{\dff \sa}$\sss\sss 
is\dss an\sss isomorphism between appropriate\dss Sobolev\dss spaces
for sufficiently\sss large\sss $t$\nnsp,\oss
and\dss hence\sss the analytical\dss index\dss is\dss also zero.\oss
Moreover,\oss the indices are equal\dss to zero even when parameters are present.\oss

\newpage
\mysection{Multiplicative\qss properties\qss of\pss pseudo-differential\qss operators}{mult-pdo}

\myuppar{Lifts\sss to\sss products of\dss manifolds.}
Let\sss $W\fff,\pff V$\sss be\sss two smooth manifolds.\oss
Let\sss $E_{\trf W}\dff,\off E_{\trf V}$\sss 
be vector bundles over\sss $W\dff,\off V$\sss respectively,\oss
and\dss let\sss $A\fff,\pff Q$ 
be operators acting on sections of\dss 
$E_{\trf W}\dff,\pff E_{\trf V}$ respectively.\oss
The\sss lifts\sss $A\dff \otimes\dff 1$\sss and\sss
$1\dff \otimes\dff Q$\sss
of\dss $A$ and\sss $Q$ 
act\sss on\sss the sections of\dss the external\dss tensor product\sss
$E_{\trf W}\dff \boxtimes\trf E_{\trf V}$\sss
and are characterized\dss by\sss the property\vspace{3pt}
\[
\quad
A\dff \otimes\dff 1\trf
\left(\qff
f\dff \otimes\dff g
\qff\right)
\off =\off
a\trf\left(\qff
f
\qff\right)\qff \otimes\qff g
\quad
\mbox{and}\quad
1\dff \otimes\dff Q\trf
\left(\qff
f\dff \otimes\dff g
\qff\right)
\off =\off
f\qff \otimes\qff Q\trf\left(\qff
g
\qff\right)
\off,
\]

\vspace{-12pt}\vspace{3pt}
where $f$ and $g$ are sections of\dss
$E_{\trf W}$ and\sss $E_{\trf V}$\sss respectively\sss and\sss
$f\dff \otimes\dff g\trf(\trf x\fff,\qff z\trf)
\off =\off
f\dff(\trf x\trf)\dff \otimes\dff g\trf(\trf z\trf)$\nnsp.\oss
When\sss we will\dss need\dss to stress\sss the second\sss factor,\oss
we will\sss denote\sss $A\dff \otimes\dff 1$\sss and\sss
$1\dff \otimes\dff Q$\sss by\vspace{3pt}
\[
\quad
A\dff \otimes\dff 1_{\dff V}
\quad
\mbox{and}\quad
1_{\trf W}\dff \otimes\dff Q\trf
\]

\vspace{-12pt}\vspace{3pt}
respectively.\oss
The\sss lifts of\dss pseudo-differential\sss operators are almost\sss never
pseudo-differential,\oss but\sss can be approximated\dss by\sss
pseudo-differential\sss operators in\sss a very strong sense.\oss
In\sss particular,\oss the\sss lifts of\dss pseudo-differential\sss operators
of\dss order $1$\sss belong\sss to\sss the class $\overline{\mathcal{P}^{\trf 1}}$\dnsp.\oss
Also,\oss for pseudo-differential\sss operators
of\dss order $1$\sss the\sss lifts\sss $A\dff \otimes\dff 1$ and $1\dff \otimes\dff Q$
depend continuously on $A$ and $Q$ respectively as operators between\dss
Sobolev\dss spaces.\oss
See\qss \cite{as1},\oss Section\qss 5,\oss and\qss \cite{h},\oss
Theorem\qss 19.2.5,\oss Lemma\qss 19.2.6,\oss and\sss  
the proof\dss of\trs Theorem\qss 19.2.7.\oss

\myuppar{Approximations and\dss lifts.}
Let $M$ be a smooth\sss manifold
and\sss $E\fff,\pff F$\sss be vector bundles on $M$\nnsp.\oss
By an approximation of\dss an operator $P$ from sections of\dss $E$\sss to sections of\dss $F$ 
and\dss belonging\sss to\sss the class $\overline{\mathcal{P}^{\trf 1}}$ 
we understand a family 
$P_t\dff,\off t\qff \in\qff (\trf 0\fff,\qff 1\trf]$\sss 
of\dss pseudo-differential\sss operators 
from\sss the class\sss $\Psi\phg^{\dff 1}$  
with\sss the following\sss properties.\oss
First,\oss the symbols of\sss $P_t$ continuously depend on $t$ 
and converge\sss to\sss the symbol\sss of\sss $P$ when 
$t\qff \ttoo\qff 0$\nnsp.\oss
Second,\oss for every $s$\sss the family\vspace{3pt}
\[
\quad
P_t\dff \colon\dff
H_{\dff s}\trf(\trf X^{\fff \circ}\fff,\qff E\trf)
\off \ttoo\off
H_{\dff s\dff -\dff 1}\trf(\trf X^{\fff \circ}\fff,\qff E\trf)\dff,
\quad
t\qff \in\qff (\trf 0\fff,\qff 1\trf]
\]

\vspace{-12pt}\vspace{3pt} 
is\dss continuous in\sss the norm\sss topology and 
converge\sss in\sss the norm\sss topology\sss to\sss $P$ when $t\qff \ttoo\qff 0$\nnsp.\oss 
We are interested\sss mostly\sss in\sss the situation where\sss
$M\off =\off
W\dff \times\dff V$\sss
for some closed\sss manifolds $W\dff,\off V$\sss and\sss $P$\sss is\dss a\sss lift\sss of\dss the form 
$A\dff \otimes\dff 1$\sss or\sss $1\dff \otimes\dff Q$\nnsp.\oss

In\sss fact,\oss we will\sss work\sss with approximations not\sss of\dss the\sss lifts\sss 
$A\dff \otimes\dff 1$\sss and\sss $1\dff \otimes\dff Q$\sss themselves,\oss
but\sss of\dss the\sss lifts
$A_{\trf 0}\dff \otimes\dff 1$\sss and $1\dff \otimes\dff Q_{\trf 0}$\sss
of\dss some operators $A_{\trf 0}$ and $Q_{\trf 0}$ closely\sss related\sss 
to $A$ and $Q$\nnsp.\oss
In\sss particular,\oss the operators $A_{\trf 0}$ and $Q_{\trf 0}$ have\sss the same symbols\sss 
as $A$ and $Q$ respectively.\oss
Also,\oss if\dss $A$\sss
is\dss formally\sss self-adjoint,\oss 
then $A_{\trf 0}$\sss is\dss also formally\sss self-adjoint,\oss
and\trs if\dss $Q^{\fff *}$\sss is\dss formally\sss adjoint\sss to $Q$\nnsp,\oss
then $Q_{\trf 0}^{\fff *}$\sss is\dss formally\sss adjoint\sss to $Q_{\trf 0}$\nsp.\oss

The approximations are constructed\sss first\sss for operators $A$
acting\sss on\sss a\sss trivial\dss bundle over a domain\sss
$U\qff \subset\qff \rrr^{\fff n}$\sss for some $n$
and\dss their\sss lifts\sss to\sss $U\dff \times\dff U\fff'$\sss
for a domain\sss $U\fff'\qff \subset\qff \rrr^{\fff n\fff'}$\sss for some $n\fff'$\dnsp.\oss
In\sss this case approximations depend only on\sss a choice of\dss a smooth\sss family\sss
$\rho^{\fff t}\dff(\trf \xi\fff,\qff \eta\trf)\fff,\off t\qff \in\qff (\trf 0\fff,\qff 1\trf]$\sss 
of\dss real-valued\sss function of\dss two real\sss variables\sss $\xi\fff,\qff \eta$\sss
of\dss the form specified\dss in\qss \cite{as1},\oss Section\qss 5.\oss
The approximation amounts\sss to multiplying\sss the full\sss symbol\sss of\dss
$A\dff \otimes\dff 1$\sss by\sss $\rho^{\fff t}\dff(\trf \xi\fff,\qff \eta\trf)$\sss
with\sss $t\qff \ttoo\qff 0$\nnsp.\oss
H\"{o}rmander\qss \cite{h}\qss chooses first\sss a real-valued\sss function\sss
$\chi\trf(\trf \xi\fff,\qff \eta\trf)$\sss and\sss then considers\sss the family\sss
$\chi\trf(\trf \xi\fff,\qff \varepsilon\dff\eta\trf)\fff,\off \varepsilon\qff >\qff 0$\nnsp.\oss
See\qss \cite{h},\oss Lemma\qss 19.2.6.\oss
Since\sss the functions\sss $\rho^{\fff t}\dff(\trf \xi\fff,\qff \eta\trf)$
and\sss $\chi\trf(\trf \xi\fff,\qff \varepsilon\dff\eta\trf)$\sss
are real-valued,\oss these approximations respect\sss
formal\sss adjoints in\sss thes sense\sss that\dss if\dss
$A^{\fff *}$\sss is\dss formally adjoint\dss to $A$\nnsp,\oss
then\sss the approximation\sss of\dss $A^{\fff *}\dff \otimes\dff 1$\sss
defined\sss by some $t$\sss or $\varepsilon$\sss is\dss formally adjoint\sss to\sss
the approximation of\trs $A\dff \otimes\dff 1$
defined\sss by\sss the same $t$ or $\varepsilon$\nnsp.\oss
In\sss particular,\oss if\dss $A$\sss is\dss formally self-adjoint,\oss
then\sss the approximations of\dss $A\dff \otimes\dff 1$\sss
are also formally self-adjoint.\oss
Of\dss course,\oss the same results hold\sss for approximations of\trs
lifts $1\dff \otimes\dff Q$\nnsp.\oss

Next,\oss one deals with\sss manifolds without\dss boundary.\oss
We\sss largely\sss follow\dss H\"{o}rmander\qss \cite{h},\oss
the proof\dss of\pss Theorem\qss 19.2.7.\oss
Let\sss $W\fff,\qff V$\sss be\sss two such\dss manifolds,\pss 
$A$\sss be a pseudo-dif\-fer\-en\-tial\sss operator on\sss $W$\nnsp,\oss
and\sss $A^{\fff *}$\sss be an operator\sss formally adjoint\sss to $A$\nnsp.\oss
Let\sss us\sss choose open coverings\sss
$W_{\fff i}\dff,\qff i\qff \in\qff I$\sss
and\sss
$V_{j}\dff,\qff j\qff \in\pff J$\sss
of\trs $W$ and\sss $V$ respectively,\oss
such\sss that\sss all\sss $W_{\fff i}\dff,\pff V_j$
are diffeomorphic\sss to domains in euclidean spaces
and\sss the bundles in question are\sss trivial\sss over all\sss $W_{\fff i}\dff,\pff V_j$\nsp.\oss
Let\vspace{3pt}
\[
\quad
\sum\nolimits_{\dff i}\qff \varphi_{\dff i}^{\dff 4}
\off =\off
1
\quad
\mbox{and}\quad\dff
\sum\nolimits_{\dff j}\qff \psi_j^{\dff 4}
\off =\off
1
\pff 
\]

\vspace{-12pt}\vspace{3pt}
be partitions of\dss unity subordinated\sss to\sss these coverings,\oss
and\dss let\sss us consider operators\vspace{1.5pt}
\[
\quad
\varphi_{\dff i}\qff \psi_j\qff
(\trf A\dff \otimes\dff 1\trf)\qff
\varphi_{\dff i}\qff \psi_j
\quad\fff
\mbox{and}\quad\dff
\varphi_{\dff i}\qff \psi_j\qff
(\trf A^{\fff *}\dff \otimes\dff 1\trf)\qff
\varphi_{\dff i}\qff \psi_j
\off,
\]

\vspace{-12pt}\vspace{1.5pt}
where\sss $i\qff \in\qff I$\nnsp,\qss $j\qff \in\qff J$\nnsp.\oss
Clearly,\oss these operators are formally adjoint.\oss
The procedure outlined\sss in\sss the previous paragraph\sss
leads\sss to formally adjoint\sss approximations\sss
$A_{\dff i\fff j}$ and\sss $A^{\fff *}_{\dff i\fff j}$\sss
respectively of\dss these\sss two operators.\oss
Then\sss the sum\vspace{1.5pt}
\[
\quad
\widetilde{A}
\off =\off\dff
\sum\nolimits_{\dff i\fff j}\pff
\varphi_{\dff i}\qff \psi_j\qff
\left(\qff 
A_{\dff i\fff j} 
\qff\right)\qff
\varphi_{\dff i}\qff \psi_j
\]

\vspace{-12pt}\vspace{1.5pt}
approximates\sss the operator\vspace{3pt}
\[
\quad
\sum\nolimits_{\dff i\fff j}\pff
\varphi_{\dff i}\qff \psi_j\qff
\left(\qff 
\varphi_{\dff i}\qff \psi_j\qff
(\trf A\dff \otimes\dff 1\trf)\qff
\varphi_{\dff i}\qff \psi_j 
\qff\right)\qff
\varphi_{\dff i}\qff \psi_j
\]

\vspace{-31.5pt}
\[
\quad
=\off\off
\sum\nolimits_{\dff i\fff j}\pff
\varphi_{\dff i}^{\dff 2}\pff \psi_j^{\dff 2}\pff
\left(\trf A\dff \otimes\dff 1\trf\right)\pff
\varphi_{\dff i}^{\dff 2}\pff \psi_j^{\dff 2} 
\off\off =\off\off
\sum\nolimits_{\dff i}\pff
\varphi_{\dff i}^{\dff 2}\pff
\left(\trf A\dff \otimes\dff 1\trf\right)\pff
\varphi_{\dff i}^{\dff 2}
\off\off =\off\off
\Bigl(\off
\sum\nolimits_{\dff i}\pff
\varphi_{\dff i}^{\dff 2}\pff
A\qff
\varphi_{\dff i}^{\dff 2}
\pff\Bigr)
\trf \otimes\dff 1
\pff,
\]

\vspace{-12pt}\vspace{3pt}
where\sss we used\sss the fact\sss that\qss
$\dis
\sum\nolimits_{\dff j}\pff \psi_j^{\dff 4}
\off =\off
1$\nnsp.\off\oss
Similarly,\oss the sum\vspace{3pt}\vspace{-0.125pt}
\[
\quad
\widetilde{A}^{\fff *}
\off =\off\dff
\sum\nolimits_{\dff i\fff j}\pff
\varphi_{\dff i}\qff \psi_j\qff
\left(\qff 
A^{\fff *}_{\dff i\fff j} 
\qff\right)\qff
\varphi_{\dff i}\qff \psi_j
\]

\vspace{-12pt}\vspace{3pt}
approximates\sss the operator\vspace{1.5pt}
\[
\quad
\Bigl(\off
\sum\nolimits_{\dff i}\pff
\varphi_{\dff i}^{\dff 2}\pff
A^{\fff *}\qff
\varphi_{\dff i}^{\dff 2}
\pff\Bigr)
\dff \otimes\dff 1
\off.
\]

\vspace{-12pt}\vspace{1.5pt}
Since\sss
$A_{\dff i\fff j}$ and\sss $A^{\fff *}_{\dff i\fff j}$\sss
are formally adjoint,\oss the operators\sss
$\widetilde{A}$ and\sss $\widetilde{A}^{\fff *}$ are also formally adjoint,\oss
justifying\sss the notations.\oss
Let\sss us consider\sss the sums\vspace{3pt}
\[
\quad
A_{\trf 0}
\off\fff =\off\dff
\sum\nolimits_{\dff i}\pff
\varphi_{\dff i}^{\dff 2}\pff
A\qff
\varphi_{\dff i}^{\dff 2}
\quad\off
\mbox{and}\quad\off
A^{\fff *}_{\trf 0}
\off\fff =\off\dff
\sum\nolimits_{\dff i}\pff
\varphi_{\dff i}^{\dff 2}\off
A^{\fff *}\qff
\varphi_{\dff i}^{\dff 2}
\]

\vspace{-12pt}\vspace{3pt}
Clearly,\oss $A_{\trf 0}$ and\dss $A^{\fff *}_{\trf 0}$
are formally\sss adjoint,\oss
and\sss formally\sss adjoint\sss pseudo-differential\sss operators\sss
$\widetilde{A}$ and\sss $\widetilde{A}^{\fff *}$
approximate\sss $A_{\trf 0}\dff \otimes\dff 1$
and\sss $A^{\fff *}_{\trf 0}\dff \otimes\dff 1$\sss
respectively.\oss
Similarly,\oss if\dss $A$\sss is\dss formally self-adjoint,\oss
then $\widetilde{A}$\sss is\dss a formally self-adjoint\sss pseudo-differential\sss
operator approximating\sss $A_{\trf 0}\dff \otimes\dff 1$\nnsp,\oss
and\sss $A_{\trf 0}$\sss is\dss also formally self-adjoint.\oss
Unfortunately,\oss usually\sss
$A_{\trf 0}\off \neq \off A$
and\sss
$A^{\fff *}_{\trf 0}\off \neq \off A^{\fff *}$\dnsp.\oss
Still,\oss the symbols of\dss
$A_{\trf 0}$ and\sss $A$ are\sss the same,\oss
as also\sss the symbols of\dss $A^{\fff *}_{\trf 0}$ and\sss $A^{\fff *}$\dnsp.\oss

\myuppar{Approximations in\sss the\dss H\"{o}rmander\dss class.}
Let\sss us return\sss to\sss the\sss framework of\trs Section\qss \ref{pdo}.\oss
Let $V$ be a closed\sss manifold and $E\fff'$ 
be a vector bundle over $V$\dnsp,\oss
and\sss let\sss us consider\sss lifts\sss  
to\sss $X\dff \times\dff V$\dnsp.\oss
In\sss this case we are\sss interested\dss in approximations by operators\sss
in\sss the\dss H\"{o}rmander\dss class in\sss the same sense as before,\oss
except\sss that\sss
the class\sss $\Psi\phg^{\dff 1}$\sss is\dss replaced\dss by\sss 
the\dss H\"{o}rmander\dss class.\oss

The boundary of\dss $X\dff \times\dff V$\sss
is\dss equal\dss to\sss $Y\dff \times\dff V$\dnsp,\oss
and $Y\dff \times\dff [\trf 0\fff,\qff 1\dff)\dff \times\dff V$\sss
is\dss its\sss natural\sss collar.\oss
The\sss lift\dss $P\dff \otimes\dff 1$ of\dss $P$\sss to\sss $X\dff \times\dff V$\sss
has\sss the form of\dss an operator\sss in\sss the\dss H\"{o}rmander\dss class,\oss 
but\sss with\vspace{3pt}
\[
\quad
P^{\dff \mathrm{i}}\dff \otimes\dff 1\dff,\quad 
\bm{\Sigma}\dff \otimes\dff 1\dff,\quad
\mbox{and}\quad\dff
\mathbf{T}\trf(\trf t\trf)\dff \otimes\dff 1
\]

\vspace{-12pt}\vspace{3pt}
in\sss the roles of\dss 
$P^{\dff \mathrm{i}}$\dnsp,\dss $\bm{\Sigma}$\nnsp,\oss and\sss $\mathbf{T}\trf(\trf t\trf)$\nnsp.\oss
Since\sss the\dss lifts 
$P^{\dff \mathrm{i}}\dff \otimes\dff 1$ and\sss
$\mathbf{T}\trf(\trf t\trf)\dff \otimes\dff 1$ 
are usually\sss not\sss pseudo-differential,\pss
$P\dff \otimes\dff 1$\sss is\dss usually\sss not\sss in\sss the\dss
H\"{o}rmander\dss class.\oss

For our purposes\sss it\dss is\dss sufficient\sss to consider\sss the case
when $\mathbf{T}$ is\dss constant.\oss
In\sss this case\sss 
the operator\sss $\mathbf{T}\trf(\trf x_{\dff n}\trf)$\sss is\dss equal\dss to\sss
the\sss lift\sss of\dss $\mathbf{T}\trf(\trf 0\trf)$
to\sss $Y\dff \times\dff (\trf -\qff 1\fff,\qff 1\dff)$\nnsp,\oss 
i.e.\dss\vspace{3pt}
\[
\quad
\mathbf{T}\trf(\trf x_{\dff n}\trf)
\off =\off
\mathbf{T}\trf(\trf 0\trf)\dff \otimes\dff 1_{\dff (\trf -\qff 1\fff,\qff 1\dff)}
\pff 
\]

\vspace{-12pt}\vspace{3pt}
with $E_{\trf (\trf -\qff 1\fff,\qff 1\dff)}$ being\dss the\sss trivial\dss
bundle with\sss the fiber $\ccc$\nnsp.\oss
For\sss $P\dff \otimes\dff 1$\sss the double\dss lift\vspace{3pt}
\[
\quad
\left(\qff \mathbf{T}\trf(\trf 0\trf)\dff \otimes\dff 1_{\dff (\trf -\qff 1\fff,\qff 1\dff)} \qff\right)
\trf \otimes\dff 1_{\dff V}
\pff 
\]

\vspace{-12pt}\vspace{3pt}
plays\sss the role of\dss $\mathbf{T}\trf(\trf x_{\dff n}\trf)$\nnsp.\oss 
Clearly,\oss this double\sss lift\dss is\dss equal\dss to\sss the one-step\sss lift\dss
form\sss $Y$\sss to\sss
$Y\dff \times\dff (\trf -\qff 1\fff,\qff 1\dff)\dff \times\dff V$\dnsp.\oss
This\sss immediately\sss implies\sss that\sss the\dss lift\dss 
$P\dff \otimes\dff 1$ 
belongs\sss to\sss the class\sss $\overline{\mathcal{P}^{\trf 1}}$\dnsp.\oss
At\sss the same\sss time\sss this one-step\sss lift\dss is\dss equal\dss to\sss
another double\sss lift,\oss namely,\oss to\vspace{3pt}
\[
\quad
\bigl(\qff \mathbf{T}\trf(\trf 0\trf)\dff \otimes\dff 1_{\dff V} \qff\bigr)
\trf \otimes\dff 1_{\dff (\trf -\qff 1\fff,\qff 1\dff)}
\pff. 
\]

\vspace{-12pt}\vspace{3pt}
This suggests\sss the following\sss natural\sss way\sss to construct\sss an operator\sss $P_{\fff 0}$\sss
in\sss the\dss H\"{o}rmander\dss class with\sss the same symbol\sss as $P$\sss
and\sss at\dss the same\sss time\sss to construct\sss an approximation 
of\dss the\dss lift\dss $P_{\fff 0}\trf \otimes\dff 1$\sss to $X\dff \times\dff V$\sss
by operators\sss from\dss the\dss H\"{o}rmander\dss class.\oss 

To begin\sss with,\oss there exists a pseudo-differential\sss operator\sss
$P_{\fff 0}^{\dff \mathrm{i}}$ on\sss $X^{\fff \circ}$\sss 
with\sss the same symbol\sss as $P^{\dff \mathrm{i}}$
and approximations\sss $\widetilde{P}^{\trf \mathrm{i}}_t$\
of\trs the\dss lift\sss $P_{\fff 0}^{\dff \mathrm{i}}\dff \times\dff 1$\sss
to\sss $X^{\fff \circ}\dff \times\dff V$\dnsp.\oss
We can also assume\sss that\sss the\sss kernel\sss of\dss $P_{\fff 0}^{\dff \mathrm{i}}$\sss
has compact\sss support.\oss
Next,\oss there exist\sss a pseudo-differential\sss operator\sss
$T_{\fff 0}$ on\sss $Y$\sss 
with\sss the same symbol\sss as $\mathbf{T}\trf(\trf 0\trf)$
and approximations\sss $\widetilde{T}_t$\sss
of\dss the\sss lift\sss $T_{\fff 0}\dff \times\dff 1$\sss
to\sss $Y\dff \times\dff V$\dnsp.\oss
Let\sss $\mathbf{T}_{\fff 0}$\sss be\sss the constant\sss map\sss
$x_{\dff n}\off \longmapsto\off T_{\fff 0}$\nsp.\oss
After\sss replacing\sss $P^{\dff \mathrm{i}}$\sss by\sss
$P_{\fff 0}^{\dff \mathrm{i}}$\sss
and\sss $\mathbf{T}$\sss by\sss $\mathbf{T}_{\fff 0}$\sss
in\sss the definition of\dss $P$
we get\sss an operator\sss $P_{\fff 0}$\sss in\sss
the\dss H\"{o}rmander\dss class with\sss the same symbol\sss as $P$\dnsp.\oss
And after replacing\sss $P^{\dff \mathrm{i}}\dff \times\dff 1$\sss
by\sss $\widetilde{P}^{\trf \mathrm{i}}_t$\sss and\sss 
$\mathbf{T}\trf(\trf 0\trf)\dff \otimes\dff 1_{\dff V}$\sss by\sss
$\widetilde{T}_t$\sss in\sss the above description of\dss $P\dff \otimes\dff 1$\sss
we get\sss operators\sss $\widetilde{P}_t$\sss in\sss the\dss
H\"{o}rmander\dss class approximating\sss $P_{\fff 0}\dff \otimes\dff 1$\nnsp.\oss
Clearly,\oss this\sss way of\dss constructing\sss the operators\sss
$P_{\fff 0}$ and\sss $\widetilde{P}_t$\sss ensures\sss that\sss they are
self-adjoint\dss in\sss the sense of\trs Section\qss \ref{pdo}\qss
if\trs $P$\sss is.

Of\dss course,\oss similar construction applies\sss to\sss the\dss lifts\sss
$1\dff \otimes\dff Q$\sss for pseudo-differential\sss operators $Q$\sss
acting on sections of\dss a vector bundle $E\fff'$\sss over $V$\dnsp,\oss
or\sss from\sss $E\fff'$\sss to another\sss bundle $F\fff'$\nsp.\oss
If\dss $Q^{\fff *}$ and\sss $Q$ are adjoint\sss to each other,\sss ,\oss
then\sss the construction\sss leads to operators\sss $Q_{\trf 0}^{\fff *}$ and\sss $Q_{\trf 0}$
adjoint\sss to each other.\oss
The operators\sss $\widetilde{Q}_t$\sss and\dss $\widetilde{Q}^{\fff *}_t$
are adjoint\sss to each other in\sss the sense\sss that\vspace{2pt}
\begin{equation}
\label{green-simple}
\quad
\bsco{\dff \widetilde{Q}_{\dff t}\trf u\fff,\pff v \dff}_{\trf X\dff \times\dff V}
\off =\off
\bsco{\dff u\dff,\pff \widetilde{Q}^{\fff *}_{\dff t}\dff v \dff}_{\trf X\dff \times\dff V}
\end{equation}

\vspace{-12pt}\vspace{2pt}
when\sss $\widetilde{Q}_{\dff t}\trf u$ and\sss $\widetilde{Q}^{\fff *}_{\dff t}\dff v$ are defined,\oss
and,\oss in\sss particular,\oss when $u\fff,\qff v$ are smooth.\oss

\myuppar{Approximations in\sss families.}
Let\sss us return\sss to\sss the discussion of\dss manifolds without\dss boundary
and suppose\sss that\sss the manifold\sss $W\off =\off W\trf(\trf z\trf)$\sss
depends on a parameter\sss $z\qff \in\qff Z$\sss in\sss the same sense as at\sss
the end of\trs Section\qss \ref{symbols-conditions}.\oss 
So,\pss $W\trf(\trf z\trf)$\sss is\dss the fiber\sss over $z$ of\dss a\sss 
locally\sss trivial\dss bundle over $Z$\nnsp.\oss
Suppose\sss that\sss $A\off =\off A\dff(\trf z\trf)$\sss is\dss
a pseudo-differential\sss operator over $W\trf(\trf z\trf)$\sss
continuously depending on $z$ as in\qss \cite{as4}.\oss
One can allow\sss also manifolds\sss $V\off =\off V\trf(\trf z\trf)$ depending on $z$\nnsp,\oss
but\sss for our purposes\sss it\dss is\dss sufficient\sss to deal\sss with\sss the case
when $V$\sss is\dss fixed.\oss

Let\sss us choose an open covering\sss $U_{\dff i}\dff,\pff i\qff \in\qff I$\sss of\dss $Z$
such\sss that\sss the above bundle with\sss the fibers $W\trf(\trf z\trf)$\sss
is\dss trivial\sss over each $U_{\dff i}$\nsp.\oss
Assuming\sss that\sss $Z$\sss is\dss paracompact,\oss
let\sss us\sss choose a partition of\dss unity\sss $\chi_{\dff i}\dff,\off i\qff \in\qff I$\sss
subordinated\dss to\sss this covering.\oss
Over each\sss $U_{\dff i}$ we can assume\sss that\sss our bundle\dss is\dss
identified\sss with\sss the product\sss bundle\sss $W\dff \times\dff U_{\dff i}\qff \ttoo\qff U_{\dff i}$\nsp,\oss
i.e.\qss that\sss $W\trf(\trf z\trf)$\sss does not\sss depend on $z$\nnsp.\oss
Then we can apply\sss the above approximation procedure and\sss construct\sss
families\sss 
$\widetilde{A}_{\dff i\dff t}\qff(\trf z\trf)\fff,\off 
z\qff \in\qff U_{\dff i}\dff,\off 
t\qff \in\qff (\trf 0\fff,\qff 1\trf]$\sss
of\dss pseudo-differential\sss operators approximating\sss
$A_{\dff i\dff 0}\dff(\trf z\trf)\dff \otimes\dff 1$\sss for some operators
$A_{\dff i\dff 0}\dff(\trf z\trf)$\sss having\sss 
the same symbols as\sss the operators $A\dff(\trf z\trf)$\nnsp.\oss
Then\sss the sum\vspace{3pt}
\[
\quad
\widetilde{A}_{\dff t}\qff(\trf z\trf)
\off\fff =\off\dff
\sum\nolimits_{\dff i}\pff 
\chi_{\dff i}\trf(\trf z\trf)\pff
\widetilde{A}_{\dff i\dff t}\qff(\trf z\trf)
\]

\vspace{-12pt}\vspace{3pt}
approximates\sss the family of\dss operators\sss
$A_{\trf 0}\dff(\trf z\trf)\dff \otimes\dff 1$\nnsp,\oss
where\vspace{3pt}\vspace{-0.125pt}
\[
\quad
A_{\trf 0}\dff(\trf z\trf)
\off\fff =\off\dff
\sum\nolimits_{\dff i}\pff 
\chi_{\dff i}\trf(\trf z\trf)\pff
A_{\dff i\dff 0}\dff(\trf z\trf)
\]

\vspace{-12pt}\vspace{3pt}
are operators with\sss the same symbols as\sss $A\dff(\trf z\trf)$\nnsp.\oss
In\sss the\sss framework of\trs Section\qss \ref{pdo}\qss one can also allow parameters,\oss
at\dss least\sss under\sss the assumption\sss that\sss 
the maps\sss $\mathbf{T}$ are constant\dss for all\sss values of\dss the parameter.\oss
In order\sss to construct\sss the operators\sss $P_{\fff 0}\dff(\trf z\trf)$\sss and\sss the approximations\sss
$\widetilde{P}_t\dff(\trf z\trf)$\sss of\dss
$P_{\fff 0}\dff(\trf z\trf)\dff \otimes\dff 1$\sss
one needs simply\sss to apply\sss the above construction with\sss parameters\sss
to\sss the families\sss $P^{\dff \mathrm{i}}\dff(\trf z\trf)$ and\sss
$\mathbf{T}\trf(\trf 0\trf)\dff(\trf z\trf)$\nnsp.\oss
Similar,\oss but\sss simpler,\oss construction applies\sss 
to\sss the\sss lifts\sss $1\dff \otimes\dff Q$\nnsp.\oss

\myuppar{The $\omult_{\fff 1}$ products\sss and\sss boundary\sss problems.}
Let\sss $P$\sss be an operator\sss from\sss the\dss H\"{o}rmander\sss class
in a bundle $E$ over $X$ as in\trs Section\qss \ref{pdo},\oss
and\dss $Q$\sss be a pseudo-differential\sss operator\sss of\dss order $1$ in
a bundle $E\fff'$\sss over a closed\sss manifold\sss $V$\dnsp.\oss
Let\sss $Q^{\fff *}$\sss be an operator\sss formally adjoint\sss to $Q$\nnsp.\oss
Suppose\sss that\sss $P$\sss is\dss self-adjoint\sss and\dss the corresponding\sss map 
$\mathbf{T}$\sss is\dss constant.\oss
The matrix\vspace{3pt}
\[
\quad
\begin{pmatrix}
\off\dff P\dff \otimes\dff 1 &
1\dff \otimes\trf Q^{\fff *} \off
\vspace{6pt} \\
\off\dff 1\dff \otimes\trf Q &
-\qff P\dff \otimes\dff 1 \off 
\end{pmatrix}
\qff 
\]

\vspace{-12pt}\vspace{3pt}
of\trs lifts of\dss pseudo-differential\sss operators
defines an operator $P\pff \omult_{\fff 1}\qff Q$\sss acting\sss in\sss the bundle\vspace{3pt}
\[
\quad
\bigl(\trf E\dff \boxtimes\dff E\fff' \trf\bigr)
\qff \oplus\qff
\left(\trf E\dff \boxtimes\dff E\fff' \trf\right)
\qff
\]

\vspace{-12pt}\vspace{3pt}
and\dss belonging\sss to\sss the class $\overline{\mathcal{P}^{\trf 1}}$\dnsp.\oss
Let\sss $\sigma$ and $q$ be\sss the symbols of\dss $P$ and\sss $Q$\sss respectively.\oss
Then $\sigma\pff \omult_{\fff 1}\qff q$\dss is\dss the symbol\sss of\dss
$P\pff \omult_{\fff 1}\qff Q$\nnsp.\oss
The\dss Green\dss formulas\qss (\ref{green-formula})\qss and\qss (\ref{green-simple})\qss
for\sss $P$ and\sss $Q$
imply\sss a Green\dss formula for\sss the operator $P\pff \omult_{\fff 1}\qff Q$\nnsp.\oss
Namely, \vspace{4.5pt}
\[
\quad
\basco{\dff (\trf P\pff \omult_{\fff 1}\qff Q \trf)\dff u\fff,\pff v \dff}_{\trf X\dff \times\dff V}
\off -\off
\basco{\dff u\dff,\pff (\trf P\pff \omult_{\fff 1}\qff Q \trf)\dff v \dff}_{\trf X\dff \times\dff V}
\off\qff =\off\qff
\basco{\dff i\trf \bigl(\trf 
(\trf \Sigma\dff \oplus\dff -\qff \Sigma\trf)\dff \otimes\dff 1\trf\bigr)\trf 
\gamma\dff u\dff,\qff \gamma\dff v \dff}_{\trf \partial}
\off,
\]

\vspace{-12pt}\vspace{4.5pt}
where\sss $\gamma$\sss is\dss the\sss trace operator\sss taking sections over\sss $X\dff \times\dff V$\sss
to\sss their restrictions\sss to\sss $Y\dff \times\dff V$\dnsp,\oss
and\sss 
$(\trf \Sigma\dff \oplus\dff -\qff \Sigma\trf)\dff \otimes\dff 1$\sss 
is\dss the\sss lift\sss of\dss $\Sigma\dff \oplus\dff -\qff \Sigma$\sss
from\sss $Y$\sss to\sss $Y\dff \times\dff V$\dnsp.\oss

Suppose\sss that\sss $P$\sss is\dss elliptic and\sss we are given a self-adjoint\dss bundle-like boundary condition
for\sss $P$ satisfying\sss the\dss Shapiro-Lopatinskii\sss condition.\oss
Let\sss the subbundle $N$ of\dss $E\trf |\trf Y$\sss be its kernel-symbol,\oss
and\dss let\sss $\Gamma$\sss be\sss the corresponding operator
defined\sss in\dss Section\qss \ref{pdo}.\oss
Then\sss
the restriction\sss 
$P_{\dff \Gamma}$ of\dss $P$\sss to\sss $\kernel\fff \Gamma$\sss 
is\dss an unbounded self-adjoint\sss operator\sss in
$H_{\trf 0}\trf(\trf X^{\fff \circ}\fff,\qff E\trf)$
with\sss the domain\sss $\kernel\fff \Gamma$\dnsp.\oss
The operators $Q$ and\sss $Q^{\fff *}$ define bounded operators\vspace{3pt}
\[
\quad
H_{\dff 1}\trf(\trf V\fff,\qff E\fff'\qff)
\qff \ttoo\qff
H_{\trf 0}\trf(\trf V\fff,\qff E\fff'\qff)
\pff,
\]

\vspace{-12pt}\vspace{3pt}
which can\sss be considered also as unbounded operators in\sss
$H_{\trf 0}\trf(\trf V\fff,\qff E\fff'\qff)$\nnsp.\oss
Let\vspace{4.5pt}
\[
\quad
H_{\trf 0}
\off =\off
H_{\trf 0}\trf(\trf X^{\fff \circ}\fff,\qff E\trf)
\dff,\quad
H_{\dff 1}
\off =\off
H_{\dff 1}\trf(\trf X^{\fff \circ}\fff,\qff E\trf)
\dff,
\]

\vspace{-31.5pt}
\[
\quad
D_{\dff 1}
\off =\off
\kernel\fff \Gamma
\off \subset\off
H_{\dff 1}\trf(\trf X^{\fff \circ}\fff,\qff E\trf)
\dff,\quad
A\off =\off\dff P_{\dff \Gamma}
\dff,\quad
\mbox{and}\quad
\]

\vspace{-31.5pt}
\[
\quad
K_{\trf s}
\off =\off
L_{\trf s}
\off =\off
H_{\dff s}\trf(\trf V\fff,\qff E\fff'\qff)
\quad
\mbox{for}\quad
s\off =\off 0\fff,\qff 1
\dff.
\]

\vspace{-12pt}\vspace{4.5pt}
These spaces and\sss the operators\sss $P\dff,\pff A\dff,\pff Q\dff,\pff Q^{\fff *}$
satisfy all\sss assumptions of\trs Section\qss \ref{mult-operators}.\oss
Therefore we can define\sss the product\sss
$A\pff \omult_{\fff 1}\qff Q
\off =\off
P_{\dff \Gamma}\pff \omult_{\fff 1}\qff Q$\sss
as an unbounded operator\sss in\vspace{3pt}
\[
\quad
\left(\qff 
H_{\trf 0}\trf(\trf X^{\fff \circ}\fff,\qff E\trf)
\off\fff \widehat{\otimes}\off\dff
H_{\trf 0}\trf(\trf V\fff,\qff E\fff'\qff)
\qff\right)
\off \oplus\off
\left(\qff 
H_{\trf 0}\trf(\trf X^{\fff \circ}\fff,\qff E\trf)
\off\fff \widehat{\otimes}\off\dff
H_{\trf 0}\trf(\trf V\fff,\qff E\fff'\qff)
\qff\right)
\pff,
\]

\vspace{-12pt}\vspace{3pt}
and\dss by\trs Theorem\qss \ref{properties-product}\qss it\dss is\dss a self-adjoint\trs
Fredholm\dss operator with discrete spectrum and compact\sss resolvent.\oss
Since\sss the $0${\dnsp}th\dss Sobolev\dss spaces are simply\dss $L_{\dff 2}$ spaces,\oss
taking\dss the above $\widehat{\otimes}${\dnsp}-products amounts\sss to\sss taking\sss
the products of\dss manifolds.\oss
See\qss \cite{s},\oss Example\qss 7.9.\oss
Therefore\sss the product\sss 
$A\pff \omult_{\fff 1}\qff Q
\off =\off
P_{\dff \Gamma}\pff \omult_{\fff 1}\qff  Q$\sss
can\sss be considered\sss as an\sss unbounded operator\dss in\sss\vspace{3pt}
\[
\quad
H_{\trf 0}\trf\left(\qff X^{\fff \circ}\dff \times\dff V\fff,\pff E\dff \boxtimes\dff E\fff'\qff\right)
\off \oplus\off
H_{\trf 0}\trf\left(\qff X^{\fff \circ}\dff \times\dff V\fff,\pff E\dff \boxtimes\dff E\fff'\qff\right)
\]

\vspace{-33pt}
\[
\quad
=\off
H_{\trf 0}\trf\left(\qff X^{\fff \circ}\dff \times\dff V\fff,\pff 
\left(\qff E\dff \boxtimes\dff E\fff'\qff\right)
\qff \oplus\qff
\left(\qff E\dff \boxtimes\dff E\fff'\qff\right)
\qff\right)
\pff.
\]

\vspace{-12pt}\vspace{3pt}
By\trs Theorem\qss \ref{domain-of-product}\qss the domain of\dss
$A\pff \omult_{\fff 1}\qff Q$\sss 
is\dss equal\dss to\sss
$(\trf D\qff\fff \widehat{\otimes}\pff K \trf)_{\dff 1}
\dff \oplus\dff
(\trf D\qff\fff \widehat{\otimes}\pff K \trf)_{\dff 1}$\nsp.\oss
Let\sss us describe\sss this domain\sss more explicitly.\oss
First\sss of\dss all,\pss 
$(\trf H\qff\fff \widehat{\otimes}\pff K \trf)_{\dff 1}$\sss
can\sss be identified\sss with\sss the\dss Sobolev\dss space\vspace{1.5pt}
\[
\quad
H_{\dff 1}\trf\left(\qff X^{\fff \circ}\dff \times\dff V\fff,\pff E\dff \boxtimes\dff E\fff'\qff\right)
\pff.
\]

\vspace{-12pt}\vspace{1.5pt}
See\qss \cite{pa},\oss Section\qss XIV.2.\oss
Let\sss
$\Gamma\dff \boxtimes\dff E\fff'$\sss
be\sss the composition of\dss the restriction operator\vspace{3pt}
\[
\quad
\gamma
\dff \colon\dff
H_{\dff 1}\trf\left(\qff X^{\fff \circ}\dff \times\dff V\fff,\pff E\dff \boxtimes\dff E\fff'\qff\right)
\off \ttoo\off
H_{\dff 1/2}\trf(\trf Y\dff \times\dff V\fff,\pff E\dff \boxtimes\dff E\fff'\trf |\qff Y\dff \times\dff V\qff)
\pff
\]

\vspace{-12pt}\vspace{1.5pt}
with\sss the operator\vspace{1.5pt}
\[
\quad
H_{\dff 1/2}\dff(\trf Y\dff \times\dff V\fff,\qff E\dff \boxtimes\dff E\fff'\trf |\qff Y\dff \times\dff V\qff)
\off \ttoo\off
H_{\dff 1/2}\dff(\trf Y\dff \times\dff V\fff,\qff N^{\dff \perp}\dff \boxtimes\dff E\fff'\qff)
\pff
\]

\vspace{-12pt}\vspace{3pt}
induced\dss by\sss the orthogonal\sss projection\sss
$E\dff \boxtimes\dff E\fff'\trf |\qff Y\dff \times\dff V
\qff \ttoo\qff
N^{\dff \perp}\dff \boxtimes\dff E\fff'$\dnsp.\oss
Recall\dss that\sss $(\trf D\qff\fff \widehat{\otimes}\pff K \trf)_{\dff 1}$\sss
is\dss the preimage of\dss 
$D_{\dff 1}\qff\fff \widehat{\otimes}\pff K_{\trf 0}$\sss
under\sss the inclusion\sss 
$(\trf H\qff\fff \widehat{\otimes}\pff K \trf)_{\dff 1}
\qff \ttoo\qff
H_{\dff 1}\qff\fff \widehat{\otimes}\pff K_{\trf 0}$\nsp.\oss
Since\sss
$D_{\dff 1}
\off =\off
\kernel\fff \Gamma$\dnsp,\oss
it\sss follows\sss that\sss
$(\trf D\qff\fff \widehat{\otimes}\pff K \trf)_{\dff 1}$\sss
is\dss equal\dss to\sss the kernel\sss of\dss\vspace{3pt}
\[
\quad
\Gamma\dff \boxtimes\dff E\fff'
\dff \colon\dff
H_{\dff 1}\trf\left(\qff X^{\fff \circ}\dff \times\dff V\fff,\pff E\dff \boxtimes\dff E\fff'\qff\right)
\off \ttoo\off
H_{\dff 1/2}\dff(\trf Y\dff \times\dff V\fff,\qff N^{\dff \perp}\dff \boxtimes\dff E\fff'\qff)
\pff.
\]

\vspace{-12pt}\vspace{3pt}
Therefore\sss the domain of\dss
$A\pff \omult_{\fff 1}\qff Q$\sss
is\dss equal\dss to\sss the kernel\sss of\dss
$(\qff \Gamma\dff \boxtimes\dff E\fff'\qff)
\dff \oplus\dff
(\qff \Gamma\dff \boxtimes\dff E\fff'\qff)$\nnsp.\oss
The\sss latter\sss operator can be considered as\sss the composition
of\dss the restriction of\dss $H_{\dff 1}${\dnsp}-sections of\dss\vspace{1.5pt}
\[
\quad 
\left(\qff E\dff \boxtimes\dff E\fff'\qff\right)
\qff \oplus\qff
\left(\qff E\dff \boxtimes\dff E\fff'\qff\right)
\]

\vspace{-12pt}\vspace{1.5pt}
to\sss $Y\dff \times\dff V$\sss with\sss the operator\sss in\sss $H_{\dff 1/2}$\sss induced\dss
by\sss the orthogonal\dss projection\sss to\sss\vspace{1.5pt}
\[
\quad 
\left(\qff N^{\dff \perp}\dff \boxtimes\dff E\fff'\qff\right)
\qff \oplus\qff
\left(\qff N^{\dff \perp}\dff \boxtimes\dff E\fff'\qff\right)
\pff.
\]

\vspace{-12pt}\vspace{1.5pt}
We see\sss that\sss the domain of\dss $A\pff \omult_{\fff 1}\qff Q$\sss
is\dss determined\dss by\sss a bundle-like boundary condition of\dss the\sss type
considered\sss in\dss Section\qss \ref{pdo}.\oss
A direct\sss check shows\sss that\sss it\dss is\dss elliptic and self-adjoint.\oss

\myuppar{Approximations of\sss $\omult_{\fff 1}$ products.}
The above approximation\sss procedure\sss leads\sss to\sss
pseudo-dif\-fer\-en\-tial\sss operators 
$P_{\dff 0}\dff,\off\dff Q_{\trf 0}${\nsp}
and $Q_{\trf 0}^{\fff *}$ with\sss 
the same symbols as $P\fff,\off\qff Q$ and\sss $Q^{\fff *}$ respectively,\oss
and\dss families of\dss operators\qss
$\widetilde{P}_t\dff,\off\qff
\widetilde{Q}_{\dff t}$
and\dss
$\widetilde{Q}^{\fff *}_{\dff t}$\sss in\sss the\dss H\"{o}rmander\sss class
approximating\sss the\sss lifts\vspace{3pt}
\[
\quad
P_{\dff 0}\dff \otimes\dff 1\dff,
\quad
1\dff \otimes\dff Q_{\trf 0}\dff,
\quad
\mbox{and}\quad
1\dff \otimes\dff Q_{\trf 0}^{\fff *}
\]

\vspace{-12pt}\vspace{3pt}
respectively.\oss
Moreover,\oss  
$P_{\dff 0}$ and\sss $\widetilde{P}_t$ are self-adjoint,\pss
$Q_{\trf 0}^{\fff *}$\dss is\dss adjoint\sss to\sss $Q_{\trf 0}$\nsp,\pss
and\dss $\widetilde{Q}^{\fff *}_{\dff t}$\sss 
is\dss adjoint\dss to\sss $\widetilde{Q}_{\dff t}$
in\sss the sense explained above.\oss
The operator\sss $P_{\dff 0}\pff \omult_{\fff 1}\qff Q_{\trf 0}$\sss
has\sss the same symbol\sss as $P\pff \omult_{\fff 1}\qff Q$
and satisfies\sss the same\dss Green\dss formula.\oss
The approximations of\trs lifts\sss lead\dss to operators\vspace{3pt}
\[
\quad
[\qff P\pff \omult_{\fff 1}\qff Q \trf]_{\fff t}
\off\dff =\off\qff
\begin{pmatrix}
\off\dff \widetilde{P}_t &
\widetilde{Q}^{\fff *}_{\dff t} \off
\vspace{6pt} \\
\off\dff \widetilde{Q}_{\dff t} &
-\qff \widetilde{P}_t \off 
\end{pmatrix}
\qff 
\]

\vspace{-12pt}\vspace{3pt}
approximating\sss $P_{\dff 0}\pff \omult_{\fff 1}\qff Q_{\trf 0}$\nsp.\oss
The operators\sss 
$[\qff P\pff \omult_{\fff 1}\qff Q \trf]_{\fff t}$\sss
also satisfy\dss the same\dss Green\dss formula\sss as
$P\pff \omult_{\fff 1}\qff Q$\sss because\sss the matrix entries have\sss
the same adjointness properties and\sss the endomorphism\sss
$\bm{\Sigma}$\sss responsible for\sss the right\dss hand side 
of\dss the\dss Green\dss formula
stays\sss the same during\sss the approximation.\oss
It\sss follows\sss that\sss the operator\sss
$(\qff \Gamma\dff \boxtimes\dff E\fff'\qff)
\dff \oplus\dff
(\qff \Gamma\dff \boxtimes\dff E\fff'\qff)$\sss
is\dss an\sss elliptic self-adjoint\dss boundary operator\sss
for\sss $[\qff P\pff \omult_{\fff 1}\qff Q \trf]_{\fff t}$\nsp.\oss
The corresponding\sss unbounded\sss self-adjoint\sss operators,\oss
as we will\sss see in a moment,\oss
are\sss the restrictions
$[\qff A\pff \omult_{\fff 1}\qff Q \trf]_{\fff t}$ of\dss operators\sss
$[\qff P\pff \omult_{\fff 1}\qff Q \trf]_{\fff t}$\sss
to\vspace{3pt}
\[
\quad
\mathcal{D}
\off =\off\dff
\kernel\fff 
\left(\qff \Gamma\dff \boxtimes\dff E\fff'\qff\right)
\dff \oplus\dff
\left(\qff \Gamma\dff \boxtimes\dff E\fff'\qff\right)
\pff. 
\]

\vspace{-12pt}\vspace{3pt}
Let\sss us define
$[\qff A\pff \omult_{\fff 1}\qff Q \trf]_{\trf 0}$
as\sss the restriction of\dss 
$P_{\dff 0}\pff \omult_{\fff 1}\qff Q_{\trf 0}$\sss 
to $\mathcal{D}$\dnsp.\oss

\mypar{Lemma.}{approx-product}
\emph{If\qss $t\qff \geq\qff 0$ is\dss sufficiently\sss small,\oss
then\dss
$[\qff A\pff \omult_{\fff 1}\qff Q \trf]_{\fff t}$\sss
is\dss a self-adjoint\trs Fredholm\dss operator\sss
with\sss discrete\sss spectrum and\sss compact\sss resolvent\qss
(and\sss its\sss domain\dss is\dss equal\dss to $\mathcal{D}$\dnsp).\oss
For sufficiently\sss small\sss $\varepsilon\qff >\qff 0$\sss
the family\sss
$[\qff A\pff \omult_{\fff 1}\qff Q \trf]_{\fff t}\qff,\off 
t\qff \in\qff [\trf 0\fff,\qff \varepsilon\trf]$\dss
is\dss continuous in\sss the norm\sss resolvent\sss sense and\dss hence\dss
is\trs Fredholm\dss in\dss the sense of\oss \textup{\cite{i2}}.\oss}

\proof
Since\sss $P_{\dff 0}\pff \omult_{\fff 1}\qff Q_{\trf 0}$\sss
has\sss the same symbol\sss as $P\pff \omult_{\fff 1}\qff Q$
and satisfies\sss the same\dss Green\dss formula,\oss
we may assume,\oss to simplify\dss the notations,\oss that\sss $P_{\dff 0}\off =\off P$
and\sss $Q_{\trf 0}\off =\off Q$\nnsp.\oss

Since\sss the symbols of\dss $[\qff P\pff \omult_{\fff 1}\qff Q \trf]_{\fff t}$\sss
converge\sss to\sss the symbol\sss of\dss
$P\pff \omult_{\fff 1}\qff Q$\nnsp,\oss
the boundary condition for $P\pff \omult_{\fff 1}\qff Q$\sss discussed above,\oss
namely,\oss the bundle-like boundary condition determined\dss by\sss the bundle\sss
$(\qff N\dff \boxtimes\dff E\fff'\qff)
\qff \oplus\qff
(\qff N\dff \boxtimes\dff E\fff'\qff)$\nnsp,\oss
is\dss an elliptic self-adjoint\dss boundary condition for\sss
$[\qff P\pff \omult_{\fff 1}\qff Q \trf]_{\fff t}$\sss
for sufficiently small\sss $t\qff >\qff 0$\nnsp.\oss
Therefore\dss the discussion in\dss Section\qss \ref{pdo}\qss applies,\oss
and\dss hence\dss Theorem\qss \ref{fredholm}\qss implies\sss that\sss
for sufficiently small\sss $t\qff >\qff 0$\sss the operator\sss
$[\qff A\pff \omult_{\fff 1}\qff Q \trf]_{\fff t}$\sss
is\dss a self-adjoint\trs Fredholm\dss operator\sss
with\sss discrete\sss spectrum and\sss compact\sss resolvent,\pss
and\sss that\sss its domain\dss is\dss equal\dss to\sss $\mathcal{D}$\nsp\dnsp.\oss
For\sss $t\off =\off 0$\sss the same properties follow\trs
Theorems\qss \ref{properties-product}\qss and\qss \ref{domain-of-product}.\oss
The operators\sss 
$[\qff P\pff \omult_{\fff 1}\qff Q \trf]_{\fff t}$\sss
with $t\qff >\qff 0$
define bounded operators\vspace{3pt}
\[
\quad
H_{\dff 1}\trf\left(\qff X^{\fff \circ}\dff \times\dff V\fff,\pff 
\left(\qff E\dff \boxtimes\dff E\fff'\qff\right)
\qff \oplus\qff
\left(\qff E\dff \boxtimes\dff E\fff'\qff\right)
\qff\right)
\off \ttoo\off
H_{\trf 0}\trf\left(\qff X^{\fff \circ}\dff \times\dff V\fff,\pff 
\left(\qff E\dff \boxtimes\dff E\fff'\qff\right)
\qff \oplus\qff
\left(\qff E\dff \boxtimes\dff E\fff'\qff\right)
\qff\right)
\]

\vspace{-12pt}\vspace{3pt}
continuously depending on\sss $t\qff >\qff 0$\nnsp.\oss
Also,\pss $P\pff \omult_{\fff 1}\qff Q$\sss defines a bounded operator between\sss
the same spaces,\oss and\dss the operators 
$[\qff P\pff \omult_{\fff 1}\qff Q \trf]_{\fff t}$\sss
converge\sss in\sss norm\sss to\sss 
$P\pff \omult_{\fff 1}\qff Q$\sss when $t\qff \ttoo\qff 0$\nnsp.\oss
See\qss \cite{h},\oss Lemma\qss 19.2.6\qss and\trs Theorems\qss 19.2.5\qss and\qss 19.2.7.\oss
Therefore\sss when $t\qff \ttoo\qff 0$\sss the operators
$[\qff A\pff \omult_{\fff 1}\qff Q \trf]_{\fff t}$\sss
converge in\sss norm\sss to\sss the restriction of\dss
$P\pff \omult_{\fff 1}\qff Q$\sss
to $\mathcal{D}$\dnsp,\oss
i.e.\qss to $A\pff \omult_{\fff 1}\qff Q$\nnsp.
It\sss follows\sss that\sss
$[\qff A\pff \omult_{\fff 1}\qff Q \trf]_{\fff t}\qff,\off 
t\qff \in\qff [\trf 0\fff,\qff \varepsilon\trf]$\dss
is\dss continuous as a family of\dss bounded operators\vspace{3pt}
\begin{equation}
\label{from-to}
\quad
\mathcal{D}
\dff \ttoo\pff
H_{\trf 0}\trf\left(\qff X^{\fff \circ}\dff \times\dff V\fff,\pff 
\left(\qff E\dff \boxtimes\dff E\fff'\qff\right)
\qff \oplus\qff
\left(\qff E\dff \boxtimes\dff E\fff'\qff\right)
\qff\right)
\end{equation}

\vspace{-12pt}\vspace{3pt}
for sufficiently small\sss $\varepsilon$\nnsp.\oss
Let\sss $t\qff \in\qff [\trf 0\fff,\qff \varepsilon\trf]$\nnsp.\oss
Since\sss $[\qff A\pff \omult_{\fff 1}\qff Q \trf]_{\dff t}$\sss is\trs Fredholm,\oss
there exists\sss $\lambda\qff \in\qff \rrr$\sss such\sss that\sss
$[\qff A\pff \omult_{\fff 1}\qff Q \trf]_{\dff t}\qff -\qff \lambda$\sss
is\dss invertible as a bounded operator\qss (\ref{from-to}).\oss
Then\sss
$[\qff A\pff \omult_{\fff 1}\qff Q \trf]_{\dff u}\qff -\qff \lambda$\sss
is\dss invertible when $u$ 
is\dss close\sss to $t$ and\sss the inverse 
$(\qff
[\qff A\pff \omult_{\fff 1}\qff Q \trf]_{\dff u}\qff -\qff \lambda
\qff)^{\dff -\dff 1}$
continuously depends on $u$\nnsp.\oss
It\sss follows\sss that\sss $[\qff A\pff \omult_{\fff 1}\qff Q \trf]_{\fff t}$
continuously depends on $t$\sss in\sss the norm\sss resolvent\sss sense.\oss
See\qss \cite{k},\oss Theorem\qss IV.2.25.\oss
As in\sss the\sss last\sss subsection of\trs Section\qss \ref{abstract-index},\oss
this implies\sss that\sss
$[\qff A\pff \omult_{\fff 1}\qff Q \trf]_{\fff t}\qff,\off 
t\qff \in\qff [\trf 0\fff,\qff \varepsilon\trf]$\dss
is\dss a\dss Fredholm\dss family\sss 
in\dss the sense of\oss \textup{\cite{i2}}.\oss  \eproof

\myuppar{Connecting\sss $A\pff \omult_{\fff 1}\qff Q$\sss with\sss
$[\qff A\pff \omult_{\fff 1}\qff Q \trf]_{\trf 0}$.}
For\sss $u\qff \in\qff [\trf 0\fff,\qff 1\trf]$\dss let\vspace{3pt}
\[
\quad 
P_u
\off =\off
u\trf P\qff +\qff (\trf 1\qff -\qff u\trf)\qff P_{\dff 0}
\quad
\mbox{and}\quad
Q_{\dff u}
\off =\off
u\trf Q\qff +\qff (\trf 1\qff -\qff u\trf)\qff Q_{\dff 0}
\pff.
\]

\vspace{-12pt}\vspace{3pt}
For every\sss $u$\sss the symbols of\dss $P_u$ and\sss $Q_{\dff u}$\sss
are equal\dss to\sss the symbols of\dss $P$ and\sss $Q$ respectively,\oss
and\dss hence\sss the symbol\sss of\dss
$P_u\pff \omult_{\fff 1}\qff Q_{\dff u}$\sss
is\dss equal\dss to\sss the symbol\sss of\dss
$P\pff \omult_{\fff 1}\qff Q$\nnsp.\oss
It\sss follows\sss that\sss the operators\sss $P_u\pff \omult_{\fff 1}\qff Q_{\dff u}$\sss
satisfy\sss the same\dss
Green\dss formula as $P\pff \omult_{\fff 1}\qff Q$\nnsp.\oss
In\sss turn,\oss this implies\sss that\sss
$(\qff \Gamma\dff \boxtimes\dff E\fff'\qff)
\dff \oplus\dff
(\qff \Gamma\dff \boxtimes\dff E\fff'\qff)$\sss
is\dss an elliptic boundary condition\sss for\sss
$P_u\pff \omult_{\fff 1}\qff Q_{\dff u}$\sss
in\sss the sense of\trs Section\qss \ref{abstract-index}.\oss
The corresponding\sss unbounded self-adjoint\sss operators,\oss
which we will\sss denote by\sss
$\fclass{\trf A\pff \omult_{\fff 1}\qff Q \qff}_{\dff u}$,\oss
are\sss the restrictions of\dss the operators\sss
$P_u\pff \omult_{\fff 1}\qff Q_{\dff u}$\sss
to $\mathcal{D}$\dnsp.\oss
Theorems\qss \ref{properties-product}\qss and\qss \ref{domain-of-product}\qss imply\sss that\sss 
$\fclass{\trf A\pff \omult_{\fff 1}\qff Q \qff}_{\dff u}$\sss
are self-adjoint\trs Fredholm\dss operators with discrete spectrum and compact\sss resolvent\sss
having $\mathcal{D}$ as\sss the domain.\oss
As\sss in\sss the proof\dss of\trs Lem\-ma\qss \ref{approx-product},\oss
this\sss implies\sss that\sss the family\dss
$\fclass{\trf A\pff \omult_{\fff 1}\qff Q \qff}_{\dff u}\qff,\off 
u\qff \in\qff [\trf 0\fff,\qff 1\trf]$\dss
is\dss continuous in\sss the norm\sss resolvent\sss sense and\dss hence\dss
is\trs Fredholm\dss in\dss the sense of\oss \textup{\cite{i2}}.\oss

\myuppar{Approximation of\sss $\omult_{\fff 1}$ products\sss in\sss families.}
Suppose now\sss that\sss the manifolds 
$X\off =\off X\trf(\trf z\trf)$\nnsp,\qss $V\off =\off V\trf(\trf z\trf)$
and\sss the operators 
$P\off =\off P\trf(\trf z\trf)$\nnsp,\qss $Q\off =\off Q\trf(\trf z\trf)$
continuously depend on a parameter $z\qff \in\qff Z$\dss in\sss the sense
of\qss Section\qss \ref{pdo},\oss
where\sss $Z$\dss is\dss compactly generated and\sss paracompact.\oss
The approximations of\trs lifts and of\sss $\omult_{\fff 1}$ products work\sss with parameters,\oss
as also\sss the construction of\dss the family\sss
$\fclass{\trf A\pff \omult_{\fff 1}\qff Q \qff}_{\dff u}\qff,\off 
u\qff \in\qff [\trf 0\fff,\qff 1\trf]$\nnsp.\oss
The\sss latter construction shows\sss that\sss the families\sss
$A\trf(\trf z\trf)\pff \omult_{\fff 1}\qff Q\trf(\trf z\trf)
\dff,\off 
z\qff \in\qff Z$\sss
and\dss
$[\qff A\trf(\trf z\trf)\pff \omult_{\fff 1}\qff Q\trf(\trf z\trf) \trf]_{\trf 0}
\dff,\off
z\qff \in\qff Z$\sss
are\sss Fredholm\dss homotopic.\oss
Lemma\qss \ref{approx-product}\qss also remains valid when\sss the parameter $z$\sss is\dss present,\oss
but\sss with $\varepsilon\qff >\qff 0$ depending on $z$\nnsp.\oss
Using a partition of\dss unity on $Z$ we can choose\sss
$\varepsilon\off =\off  \varepsilon\trf(\trf z\trf)\qff >\qff 0$ continuously depending on $z$\nnsp.\oss
This\sss leads\sss to a\dss Fredholm\dss family\sss
$[\qff A\trf(\trf z\trf)\pff \omult_{\fff 1}\qff Q\trf(\trf z\trf) \trf]_{\fff t}\qff,\off 
t\qff \in\qff [\trf 0\fff,\qff \varepsilon\trf(\trf z\trf)\trf]\dff,\pff z\qff \in\qff Z$\nnsp.\oss
After an obvious reparametrization we get\sss a\dss Fredholm\dss homotopy\sss
between\sss the family\sss
$[\qff A\trf(\trf z\trf)\pff \omult_{\fff 1}\qff Q\trf(\trf z\trf) \trf]_{\trf 0}$,\dss
$z\qff \in\qff Z$ and\sss the family\sss of\dss restrictions of\dss operators\sss
$[\qff P\trf(\trf z\trf)
\pff \omult_{\fff 1}\qff 
Q\trf(\trf z\trf) \trf]_{\trf \varepsilon\trf(\trf z\trf)}$\sss
to\sss $\mathcal{D}$\dnsp.\oss
It\sss follows\sss that\sss the family
$A\trf(\trf z\trf)\pff \omult_{\fff 1}\qff Q\trf(\trf z\trf)
\dff,\off 
z\qff \in\qff Z$\sss
is\dss also\dss Fredholm\dss homotopic\sss to\sss the\sss latter\sss family.\oss

\newpage
\mysection{Glueing\qss and\qss cutting}{glueing}

\myuppar{The framework.}
We will\dss keep\sss the assumptions and notations of\trs Section\qss \ref{pdo}.\oss
Following\dss H\"{o}rmander\qss \cite{h},\oss Section\qss 20.3,\oss
we introduce\sss the double\sss $\widehat{X}$\sss of\dss $X$\sss
consisting of\dss two copies $X_{\dff 1}$ and $X_{\dff 2}$ of\dss $X$\nnsp.\oss
Let\sss $Y$\sss be\sss the common boundary of\dss $X_{\dff 1}$ and $X_{\dff 2}$\nsp.\oss
We will\sss identify\sss the copy of\dss the collar\sss
$\partial\dff X\dff \times\dff [\trf 0\fff,\qff 1\dff)$\sss
in\sss $X_{\dff 1}$\sss with\sss $Y\dff \times\dff [\trf 0\fff,\qff 1\dff)$\sss
and\sss the copy\sss in $X_{\dff 2}$\sss with\sss
$Y\dff \times\dff (\dff -\qff 1\fff,\qff 0\trf]$\nnsp,\oss
with $x_{\dff n}$\sss changing\sss the sign\sss in\sss the second copy.\oss
The smooth structure on\sss $\widehat{X}$\sss is\dss defined\sss by\sss
taking\sss the smooth structures on\sss 
$X_{\dff 1}$ and $X_{\dff 2}$\sss and\sss the product\sss smooth structure on\sss
$Y\dff \times\dff (\trf -\qff 1\fff,\qff 1\trf)$\nnsp.\oss
This smooth structure depends on\sss the identification of\dss the collar
with\sss $\partial\dff X\dff \times\dff [\trf 0\fff,\qff 1\dff)$\nnsp.\oss

Recall\dss that\sss $E$\sss is\dss a bundle over $X$\nnsp.\oss
Let\sss $E_{\dff 1}$ and\sss $E_{\trf 2}$\sss be\sss the copies of\dss
the bundle\sss $E$\sss over $X_{\dff 1}$ and $X_{\dff 2}$ respectively,\oss
and\dss let\sss $\widehat{E}$\sss be\sss the bundle over\sss $\widehat{X}$\sss
obtained\dss by\sss glueing\sss $E_{\dff 1}$ and\sss $E_{\trf 2}$\sss over\sss $Y$\dnsp.\oss
Let\dss $F_{\dff 2}\off =\off E_{\trf 2}\dff \oplus\dff E_{\trf 2}$\sss
and\dss let\dss $F_{\dff 1}$\sss be a\dss Hermitian\dss bundle over $X_{\dff 1}$\sss
which\dss is\dss equal\dss to\sss the bundle\sss $E_{\dff 1}\dff \oplus\dff E_{\dff 1}$\sss
over\sss the collar\sss
$Y\dff \times\dff [\trf 0\fff,\qff 1\dff)\qff \subset\qff X_{\dff 1}$\nsp.\oss
Let\sss $\widehat{F}$\sss be\sss the bundle over\sss $\widehat{X}$\sss
obtained\dss by\sss glueing\sss $F_{\dff 1}$\sss and\sss $F_{\dff 2}$\sss
over\sss $Y$\dnsp.\oss

\myuppar{The operators.}
As in\dss Section\qss \ref{pdo},\oss
let\sss $E\dff|\dff Y$\sss be presented
as\sss a direct\sss sum\sss
$E\dff|\dff Y
\off =\off
E^{\dff +}\dff \oplus\dff E^{\dff -}$\dnsp.\oss
Let\sss us define operators $P^{\dff \mathrm{b}}$ and\dss $\widetilde{P}^{\trf \mathrm{b}}$\sss 
acting on sections of\dss $\widehat{E}$\sss
by\sss the same formulas as in\dss Section\qss \ref{pdo},\oss
and\dss then define\sss the operator\sss
$\widehat{P}^{\dff \mathrm{b}\qff \sa}$\sss acting on sections of\trs
$\widehat{E}\trf \oplus\trf \widehat{E}$\dss by\sss the matrix\vspace{1.25pt}
\[
\quad
\widehat{P}^{\qff \mathrm{b}\qff \sa}\dff u
\off =\off\dff
\begin{pmatrix}
\off\dff 0 &
P^{\dff \mathrm{b}}\dff u \qff\off
\vspace{6pt} \\
\off\dff -\qff \widetilde{P}^{\trf \mathrm{b}}\dff u
\qff -\qff 
i\trf \varphi\fff'\dff(\trf x_{\dff n}\trf)\dff u^{\dff +}
\qff +\qff 
i\trf \varphi\fff'\dff(\trf x_{\dff n}\trf)\dff u^{\dff -} &
0 \qff\off 
\end{pmatrix}
\qff.
\]

\vspace{-12pt}\vspace{1.25pt}
Since\sss the cut-off\dss function\sss $\varphi$\sss used\sss to define\sss
$P^{\dff \mathrm{b}}$ and\dss $\widetilde{P}^{\trf \mathrm{b}}$\sss 
has support\sss in\sss $(\trf -\qff 1\fff,\qff 1\trf)$\sss 
and\sss $\widehat{F}$\sss is\dss equal\dss to\sss 
$\widehat{E}\trf \oplus\trf \widehat{E}$\sss over\sss
$Y\dff \times\dff (\trf -\qff 1\fff,\qff 1\trf)$\nnsp,\oss
we may,\oss and\sss from\sss now on\sss will,\oss 
consider\sss $\widehat{P}^{\qff \mathrm{b}\qff \sa}$\sss
as an operator acting on sections of\trs $\widehat{F}$\nnsp.\oss
Let\dss $P^{\dff \mathrm{b}\qff \sa}_{\fff 1}$\sss
and\dss $P^{\dff \mathrm{b}\qff \sa}_{\trf 2}$\sss 
be\sss the operators defined\sss in\sss the same way,\oss 
but\sss acting on\sss the sections of\dss 
$F_{\dff 1}$\sss and\sss
$F_{\dff 2}\off =\off E_{\trf 2}\dff \oplus\dff E_{\trf 2}$\sss 
over\sss $X_{\dff 1}$\sss and\sss $X_{\dff 2}$\sss respectively.\oss

Let\sss 
$P^{\dff \mathrm{i}\qff \sa}_{\fff 1}$\sss 
be a formally self-adjoint\sss elliptic operator\sss 
belonging\sss to\vspace{0.875pt}
\[
\quad
\Psi\phg^{\dff 1}\trf
\left(\qff 
X_{\dff 1}\dff \smallsetminus\qff Y\dff;\off 
F_{\dff 1}\dff,\off 
F_{\dff 1}
\qff\right)
\qff
\]

\vspace{-12pt}\vspace{0.875pt}
and\sss such\sss that\sss its kernel\dss has compact\sss support\sss in\sss
$(\trf X_{\dff 1}\dff \smallsetminus\qff Y\trf)
\trf \times\trf 
(\trf X_{\dff 1}\dff \smallsetminus\qff Y\trf)$\nnsp.\oss
Let\dss\vspace{1.25pt}
\[
\quad
P^{\trf \sa}_{\dff 1}
\off =\off
P^{\dff \mathrm{i}\qff \sa}_{\fff 1}
\off +\off
P^{\dff \mathrm{b}\qff \sa}_{\dff 1}
\]

\vspace{-12pt}\vspace{1.25pt}
and\dss let\dss $B^{\dff \sa}_{\dff 1}$\sss be\sss the operator\sss $B^{\dff \sa}$\sss from\dss
Section\qss \ref{pdo}.\oss
Then\sss $P^{\trf \sa}_{\fff 1}$\sss together with\sss the boundary operator\sss
$B^{\dff \sa}_{\dff 1}$\sss satisfies\dss Shapiro-Lopatinskii\dss condition.\oss

The main\sss difference between\dss $P^{\trf \sa}_{\dff 1}$\sss
and\sss the standard operator\sss $P^{\trf \sa}$\sss
from\dss Section\qss \ref{pdo}\qss is\dss in allowing an
arbitrary\sss formally self-adjoint\qss ``interior''\qss operator\sss
$P^{\dff \mathrm{i}\qff \sa}_{\fff 1}$\nsp.\oss
We also omitted\dss the\qss ``correction\sss terms''\qss $\pm\off i\fff t\trf \id$\nnsp,\oss
which will\dss reappear soon.\oss

Let\dss 
$P^{\trf \mathrm{i}}_{\dff 2}$\sss
be\sss the copy over\sss $X_{\dff 2}$\sss
of\dss the operator\sss 
$P^{\dff \mathrm{i}}$\sss 
from\dss Section\qss \ref{pdo},\oss
and\dss let\vspace{1.5pt}
\[
\quad
P^{\dff \mathrm{i}\qff \sa}_{\dff 2}
\off =\off\dff
\begin{pmatrix}
\off\dff 0 &
\dff P^{\dff \mathrm{i}}_{\dff 2} \qff\off
\vspace{6pt} \\
\off\dff -\qff P^{\trf \mathrm{i}}_{\dff 2} &
\dff 0 \qff\off 
\end{pmatrix}
\qff
\]

\vspace{-12pt}\vspace{1.5pt}
and\trs
$P^{\dff \sa}_{\trf 2}
\off =\off\dff
P^{\dff \mathrm{b}\qff \sa}_{\dff 2}
\off +\off
P^{\dff \mathrm{i}\qff \sa}_{\dff 2}$\nsp.\qff\oss
Finally,\oss let\qss\vspace{3pt}
\[
\quad
\widehat{P}^{\qff \sa}
\off =\off\dff
P^{\dff \mathrm{i}\qff \sa}_{\fff 1}
\off +\off
\widehat{P}^{\qff \mathrm{b}\qff \sa}
\off +\off
P^{\dff \mathrm{i}\qff \sa}_{\dff 2}
\off. 
\]

\vspace{-12pt}\vspace{3pt}
The operator\dss $\widehat{P}^{\qff \sa}$ should\dss be\sss thought\sss
as\sss the result\sss of\qss ``glueing''\qss to a\sss fairly\sss general\sss operator\dss
$P^{\trf \sa}_{\dff 1}$\dss the standard operator\dss $P^{\dff \sa}_{\trf 2}$\nsp.\oss

The\sss operator\sss $P^{\dff \sa}_{\trf 2}$\sss has almost\sss 
the same form as\sss the operator\sss $P^{\trf \sa}$\dss
from\dss Section\qss \ref{pdo},\oss
but\sss restricting\trs
$\widehat{P}^{\trf \mathrm{b}\qff \sa}$\sss
to\sss $X_{\dff 2}$\sss changes\sss the sign of\dss $x_{\dff n}$\nsp.\oss
Changing\sss the sign of\dss $x_{\dff n}$\sss back\sss
replaces\sss $D_{\fff n}$\sss by\sss $-\qff D_{\fff n}$\nsp.\oss
Therefore\sss $P^{\dff \sa}_{\trf 2}$\sss has\sss 
the same form as\sss $P^{\trf \sa}_{\dff 1}$\sss
with\sss the roles of\dss $E^{\dff +}$\sss and\sss $E^{\dff -}$\sss interchanged,\oss
and\dss the relevant\dss boundary operator\dss is\vspace{3pt}
\[
\quad
B^{\dff \sa}_{\trf 2}\dff \colon\dff
(\dff u\fff,\qff v\trf)
\off \longmapsto\off
(\trf \gamma\trf u^{\dff +}\fff,\qff \gamma\trf v^{\dff -}\trf)
\qff
\]

\vspace{-12pt}\vspace{3pt}
instead of\dss the boundary operator\dss
$B^{\dff \sa}_{\dff 1}\off =\off B^{\dff \sa}$\dnsp.\oss

\myuppar{Cutting\dss $\widehat{P}^{\qff \sa}$\dnsp.}
H\"{o}rmander's\qss Proposition\qss 20.3.2\qss from\qss \cite{h}\qss
suggests\sss that\sss the analytical\dss index of\dss $\widehat{P}^{\qff \sa}$\sss
should\sss be equal\sss to\sss that\sss of\dss 
$P^{\trf \sa}_{\fff 1}$\sss with\sss the boundary conditions
defined\dss by\sss $B^{\dff \sa}_{\dff 1}$\nsp.\oss
Of\dss course,\oss both of\dss them are equal\dss to zero because we
are dealing\sss with self-adjoint\sss operators,\oss
but\sss we need\sss to forget\sss this fact\sss and\sss prove\sss
the equality\sss in a way allowing\sss to add\sss parameters.\oss

Let\sss us begin by adding\sss to\sss the operators\sss
$P^{\trf \sa}_{\dff 1}$\nsp,\oss
$P^{\dff \sa}_{\trf 2}$\nsp,\oss
and\trs
$\widehat{P}^{\qff \sa}$\trs
the\qss ``correction\sss term''\vspace{1.5pt}
\[
\quad
\psi\trf(\trf x\trf)\qff
\begin{pmatrix}
\off\dff 0 &
i\fff t\trf \id \qff\off
\vspace{6pt} \\
\off\dff -\qff i\fff t\trf \id &
0 \qff\off 
\end{pmatrix}
\qff,
\]

\vspace{-12pt}\vspace{1.5pt}
where $\psi\trf(\trf x\trf)$\sss is\dss equal\dss to $\varphi\trf(\trf x_{\dff n}\trf)$ 
for $x\qff \in\qff X_{\dff 1}$ and\sss to $1$ for $x\dff \in\qff X_{\dff 2}$\nsp.\oss
If\sss $t$\sss is\dss sufficiently\sss large,\oss
then $P^{\dff \sa}_{\trf 2}\dff \oplus\dff B^{\dff \sa}_{\dff 2}$
is\dss an\sss isomorphism between appropriate\dss Sobolev\dss spaces.\oss 
Adding\sss this\sss term does not\sss change\sss the analytical\dss index of\dss $\widehat{P}^{\qff \sa}$\sss
even\sss in\sss families,\oss and does not\sss change\sss the\sss topological\sss
index\sss because\sss this does not\sss affect\sss the principal\sss symbol.\oss

Let\sss $X_{\dff 1}\dff \bigsqcup\qff X_{\dff 2}$\sss be\sss the disjoint\sss union
of\dss $X_{\dff 1}$ and $X_{\dff 2}$\nnsp,\oss
and\dss let\sss $F_{\dff 1}\dff \bigsqcup\qff F_{\dff 2}$\sss
be\sss the disjoint\sss union
of\dss the bundles\sss $F_{\dff 1}$ and $F_{\dff 2}$\nsp,\oss
a bundle over\sss $X_{\dff 1}\dff \bigsqcup\qff X_{\dff 2}$\nsp.\oss
Recall\sss that\sss $F_{\dff 1}$\sss is\dss equal\dss to\sss 
$E_{\dff 1}\trf \oplus\trf E_{\dff 1}$\sss over\sss
$Y\dff \times\dff [\trf 0\fff,\qff 1\trf)$\nnsp,\oss
and\sss let\sss us represent\sss sections of\dss 
$F_{\dff 1}\dff \bigsqcup\qff F_{\dff 2}$\sss
over\sss 
$Y\dff \times\dff (\trf -\qff 1\fff,\qff 0\trf]
\off \bigsqcup\off 
Y\dff \times\dff [\trf 0\fff,\qff 1\dff)$\sss
by quadruples\sss 
$(\dff u_{\dff 1}\dff,\off v_{\dff 1}\dff,\off u_{\dff 2}\dff,\off v_{\dff 2} \trf)$\nnsp,\oss
where\sss $(\dff u_{\dff 1}\dff,\off v_{\dff 1}\trf)$\sss
and\sss $(\dff u_{\dff 2}\dff,\off v_{\dff 2} \trf)$\sss are sections of\vspace{3pt}
\[
\quad
E_{\dff 1}\trf \oplus\trf E_{\dff 1}
\quad
\mbox{and}\quad\dff
E_{\trf 2}\trf \oplus\trf E_{\trf 2}
\]

\vspace{-12pt}\vspace{3pt}
respectively.\oss
For\sss every\sss $\tau\qff \in\qff \rrr$\trs let\dss $B_{\qff \tau}$\sss be\sss 
the boundary operator\vspace{3pt}
\[
\quad
B_{\qff \tau}\qff
\bigl(\dff u_{\dff 1}\dff,\off v_{\dff 1}\dff,\off u_{\dff 2}\dff,\off v_{\dff 2} 
\trf\bigr)
\]

\vspace{-33pt}
\[
\quad
=\off\qff
\left(\pff
\gamma\trf u^{\dff -}_{\dff 1}\off -\off \tau\dff \gamma\trf u^{\dff -}_{\dff2}\dff,\quad
\gamma\trf v^{\dff +}_{\dff 1}\off -\off \tau\dff \gamma\trf v^{\dff +}_{\dff2}\dff,\quad
\gamma\trf u^{\dff +}_{\dff 2}\off -\off \tau\qff \gamma\trf u^{\dff +}_{\dff 1}\dff,\quad 
\gamma\trf v^{\dff -}_{\dff 2}\off -\off \tau\qff \gamma\trf v^{\dff -}_{\dff 1}
\off\right)
\qff
\]

\vspace{-12pt}
Then\vspace{-3pt}
\[
\quad
B_{\qff 0}\qff
\bigl(\dff u_{\dff 1}\dff,\off v_{\dff 1}\dff,\off u_{\dff 2}\dff,\off v_{\dff 2} 
\trf\bigr)
\off =\off
\left(\pff \gamma\trf u^{\dff -}_{\dff 1}\dff,\quad
\gamma\trf v^{\dff +}_{\dff 1}\dff,\quad 
\gamma\trf u^{\dff +}_{\dff 2}\dff,\quad 
\gamma\trf v^{\dff -}_{\dff 2} \off\right)
\quad
\mbox{and}
\]

\vspace{-33pt}
\[
\quad
B_{\trf 1}\qff\fff
\bigl(\dff u_{\dff 1}\dff,\off v_{\dff 1}\dff,\off u_{\dff 2}\dff,\off v_{\dff 2} 
\trf\bigr)
\]

\vspace{-33pt}
\[
\quad
=\off\qff
\left(\pff
\gamma\trf u^{\dff -}_{\dff 1}\off -\off \gamma\trf u^{\dff -}_{\dff2}\dff,\quad
\gamma\trf v^{\dff +}_{\dff 1}\off -\off \gamma\trf v^{\dff +}_{\dff2}\dff,\quad
\gamma\trf u^{\dff +}_{\dff 2}\off -\off \gamma\trf u^{\dff +}_{\dff 1}\dff,\quad 
\gamma\trf v^{\dff -}_{\dff 2}\off -\off \gamma\trf v^{\dff -}_{\dff 1}
\off\right)
\qff.
\]

\vspace{-12pt}\vspace{3pt}
The boundary operators\sss
$B_{\qff \tau}$\sss are non-local\dss for\sss $\tau\off \neq\off 0$\nnsp.\oss
This\sss issue\dss is\dss ignored\sss by\trs H\"{o}rmander\dss in\sss the proof\dss
of\qss Proposition\qss 20.3.2\qss in\qss \cite{h}.\oss
But\sss one can\sss interpret\dss 
$P^{\trf \sa}_{\dff 1}\qff \bigsqcup\pff P^{\dff \sa}_{\trf 2}$\dss 
as an operator over\sss $X$\sss acting on sections of\dss the bundle\sss 
$F_{\dff 1}\trf \oplus\trf (\trf E_{\dff 1}\trf \oplus\trf E_{\dff 1} \trf)$\nnsp,\oss
by\sss $P^{\trf \sa}_{\dff 1}$\sss on\sss the summand\dss $F_{\dff 1}\trf \oplus\trf 0$\nnsp,\oss 
and\sss by\sss $P^{\dff \sa}_{\trf 2}$\sss on\sss the summand\dss 
$0\trf \oplus\trf (\trf E_{\dff 1}\trf \oplus\trf E_{\dff 1} \trf)$\nnsp.\oss
With\sss this interpretations\sss the boundary operators\sss
$B_{\qff \tau}$\sss are\sss local.\oss
Moreover,\oss this\sss interpretation of\qss 
$P^{\trf \sa}_{\dff 1}\qff \bigsqcup\pff P^{\dff \sa}_{\trf 2}$\dss
together with\sss the boundary operators\sss $B_{\qff \tau}$\sss
satisfies\trs Shapiro-Lopatinskii\trs condition,\oss
and\sss the corresponding\sss boundary conditions are self-adjoint.\oss
In order\sss to prove\sss this,\oss it\dss is\dss sufficient\sss to examine\sss
the kernel\sss of\dss $B_{\qff \tau}$\sss over an arbitrary\sss $u\qff \in\qff S\dff Y$\dnsp.\oss
In\sss the obvious notations,\oss this kernel\dss is\dss 
described\sss by\sss the equations\vspace{3pt}
\[
\quad
u^{\dff -}_{\dff 1}\off =\off \tau\qff u^{\dff -}_{\dff 2}\dff,\qquad
v^{\dff +}_{\dff 1}\off =\off \tau\qff v^{\dff +}_{\dff 2}\dff,\qquad
u^{\dff +}_{\dff 2}\off =\off \tau\qff u^{\dff +}_{\dff 1}\dff,\qquad 
v^{\dff -}_{\dff 2}\off =\off \tau\qff v^{\dff -}_{\dff 1}
\qff,
\]

\vspace{-12pt}\vspace{3pt}
and\sss the space\sss
$\mathcal{L}_{\dff -}\dff(\trf \rho_{\fff u}\dff)$\sss
is\dss described\sss by\sss the equations\vspace{3pt}
\[
\quad
u^{\dff +}_{\dff 1}
\off =\off
v^{\dff -}_{\dff 1}
\off =\off
0\dff,\qquad
u^{\dff -}_{\dff 2}
\off =\off
v^{\dff +}_{\dff 2}
\off =\off
0
\qff.
\]

\vspace{-12pt}\vspace{3pt}
Together\sss they\sss imply\sss that\sss also\sss
$u^{\dff +}_{\dff 2}
\off =\off
v^{\dff -}_{\dff 2}
\off =\off
0$\sss
and\sss
$u^{\dff -}_{\dff 1}
\off =\off
v^{\dff +}_{\dff 1}
\off =\off
0$\nnsp.\oss
It\sss follows\sss that\sss the\sss kernel\sss of\dss $B_{\qff \tau}$\sss
is\dss transverse\sss to\sss $\mathcal{L}_{\dff -}\dff(\trf \rho_{\fff u}\dff)$\sss
and\sss hence\dss is\dss a boundary condition\sss for\sss the elliptic pair\sss
$\Sigma_y\fff,\qff \tau_{\fff u}$\sss at\sss $u$\nnsp.\oss
Equivalently,\pss $B_{\qff \tau}$\sss satisfies\trs Shapiro-Lopatinskii\trs condition.\oss
In order\sss to prove self-adjointness we need\sss to check\sss that\sss 
the kernel\dss is\dss lagrangian.\oss
The relevant\trs Hermitian\sss product\dss is\dss the direct\sss sum of\dss
the\dss product\sss $[\trf \bullet\fff,\qff \bullet\trf]_{\dff 1}$\sss
and\sss $[\trf \bullet\fff,\qff \bullet\trf]_{\dff 2}$\sss
on\sss the\sss two copies of\dss $E_{\dff y}\dff \oplus\dff E_{\dff y}$\nsp,\oss
and\sss $[\trf \bullet\fff,\qff \bullet\trf]_{\dff 2}$\sss has\sss the same form as\sss
$[\trf \bullet\fff,\qff \bullet\trf]_{\dff 1}$\nsp,\oss but\sss with\sss
the roles of\dss $E^{\dff +}$\sss and\sss $E^{\dff -}$\sss interchanged.\oss
Therefore\sss the description of\dss the product\sss 
$[\trf \bullet\fff,\qff \bullet\trf]_{\dff 1}$\sss
at\sss the end of\trs Section\qss \ref{pdo}\qss 
implies\sss that\vspace{4.5pt}
\[
\quad
\bigl[\qff
(\dff u_{\dff 1}\dff,\off v_{\dff 1}\dff,\off u_{\dff 2}\dff,\off v_{\dff 2} \trf)\fff,\off\off
(\dff a_{\dff 1}\dff,\off b_{\dff 1}\dff,\off a_{\dff 2}\dff,\off b_{\dff 2} \trf)
\qff\bigr]
\]

\vspace{-33pt}
\[
\quad
=\off
\bsco{\dff u^{\dff +}_{\dff 1}\dff,\off b^{\dff +}_{\dff 1} \dff}
\off -\off
\bsco{\dff u^{\dff -}_{\dff 1}\dff,\off b^{\dff -}_{\dff 1} \dff}
\off +\off
\bsco{\dff v^{\dff +}_{\dff 1}\dff,\off a^{\dff +}_{\dff 1} \dff}
\off -\off
\bsco{\dff v^{\dff -}_{\dff 1}\dff,\off a^{\dff -}_{\dff 1} \dff}
\qff
\]

\vspace{-33pt}
\[
\quad
\qff +\off
\bsco{\dff u^{\dff -}_{\dff 2}\dff,\off b^{\dff -}_{\dff 2} \dff}
\off -\off
\bsco{\dff u^{\dff +}_{\dff 2}\dff,\off b^{\dff +}_{\dff 2} \dff}
\off +\off
\bsco{\dff v^{\dff -}_{\dff 2}\dff,\off a^{\dff -}_{\dff 2} \dff}
\off -\off
\bsco{\dff v^{\dff +}_{\dff 2}\dff,\off a^{\dff +}_{\dff 2} \dff}
\qff.
\]

\vspace{-12pt}\vspace{4.5pt}
If\dss
$(\dff u_{\dff 1}\dff,\off v_{\dff 1}\dff,\off u_{\dff 2}\dff,\off v_{\dff 2} \trf)
\fff,\off 
(\dff a_{\dff 1}\dff,\off b_{\dff 1}\dff,\off a_{\dff 2}\dff,\off b_{\dff 2} \trf)
\off \in\off
\kernel\trf B_{\qff \tau}$\nsp,\oss
then\sss this expression\dss is\dss equal\dss to\vspace{4.5pt}
\[
\quad
\phantom{=\off}
\bsco{\dff u^{\dff +}_{\dff 1}\dff,\off \tau\qff b^{\dff +}_{\dff 2} \dff}
\off -\off
\bsco{\dff \tau\qff u^{\dff -}_{\dff 2}\dff,\off b^{\dff -}_{\dff 1} \dff}
\off +\off
\bsco{\dff \tau\dff v^{\dff +}_{\dff 2}\dff,\off a^{\dff +}_{\dff 1} \dff}
\off -\off
\bsco{\dff v^{\dff -}_{\dff 1}\dff,\qff \tau\off a^{\dff -}_{\dff 2} \dff}
\qff
\]

\vspace{-33pt}
\[
\quad
\dff +\off
\bsco{\dff u^{\dff -}_{\dff 2}\dff,\off \tau\qff b^{\dff -}_{\dff 1} \dff}
\off -\off
\bsco{\dff \tau\qff u^{\dff +}_{\dff 1}\dff,\off b^{\dff +}_{\dff 2} \dff}
\off +\off
\bsco{\dff \tau\qff v^{\dff -}_{\dff 1}\dff,\off a^{\dff -}_{\dff 2} \dff}
\off -\off
\bsco{\dff v^{\dff +}_{\dff 2}\dff,\off \tau\qff a^{\dff +}_{\dff 1} \dff}
\qff.
\]

\vspace{-12pt}\vspace{4.5pt}
Since\sss $\tau$\sss is\dss real,\oss the\sss terms\sss in\sss the\sss last\sss
expression cancel\dss pair-wise,\oss and\sss hence\sss it\dss is\dss equal\dss to zero.\oss
Since\sss the dimension of\dss the kernel\sss of\dss $B_{\qff \tau}$\sss is\dss
equal\sss to\sss the half\dss of\dss the dimension of\dss the domain,\oss
it\sss follows\sss that\sss $\kernel\trf B_{\qff \tau}$\sss is\dss lagrangian,\oss
and\sss hence\sss $B_{\qff \tau}$\sss defines self-adjoint\sss
elliptic boundary conditions for\dss 
$P^{\trf \sa}_{\dff 1}\qff \bigsqcup\pff P^{\dff \sa}_{\trf 2}$\nsp.\oss
Clearly,\oss $B_{\qff \tau}$\sss continuously\sss depends on\sss $\tau$\nnsp.

The operator\sss 
$B_{\qff 0}$\sss is\dss equal\sss to\sss the direct\sss sum\sss
$B^{\dff \sa}_{\dff 1}\qff \oplus\qff B^{\dff \sa}_{\trf 2}$\sss
and\dss hence\dss 
$P^{\trf \sa}_{\dff 1}\qff \bigsqcup\pff P^{\dff \sa}_{\trf 2}$\dss
together\sss with\sss the boundary condition\sss $B_{\qff 0}$\sss
is\dss the direct\sss sum of\dss
$P^{\trf \sa}_{\dff 1}$\dss and\dss $P^{\dff \sa}_{\trf 2}$\dss
together\sss with\sss the boundary conditions\dss
$B^{\dff \sa}_{\dff 1}$\sss and\dss $B^{\dff \sa}_{\trf 2}$\sss
respectively.\oss
Since\sss $P^{\dff \sa}_{\trf 2}\qff \oplus\qff B^{\dff \sa}_{\dff 2}$
is\dss an\sss isomorphism,\oss
the analytical\dss index of\dss 
$P^{\trf \sa}_{\dff 1}\qff \bigsqcup\pff P^{\dff \sa}_{\trf 2}$\dss
with\sss the boundary condition\sss $B_{\qff 0}$\sss
is\dss equal\dss to\sss the analytical\dss index of\dss
$P^{\trf \sa}_{\dff 1}$\dss
with\sss the boundary conditions\dss
$B^{\dff \sa}_{\dff 1}$\nsp.\oss
Moreover,\oss this\dss is\dss true and\sss the proof\dss
remains\sss the same when\dss $P^{\trf \sa}_{\dff 1}$\sss 
depends on a parameter.\oss\vspace{-0.25pt}

Let\sss us consider now\sss the boundary operator\sss $B_{\trf 1}$\nsp.\oss
Up\sss to a permutation of\dss summands,\oss\vspace{1.5pt}
\[
\quad
B_{\trf 1}\qff\fff
\bigl(\dff u_{\dff 1}\dff,\off v_{\dff 1}\dff,\off u_{\dff 2}\dff,\off v_{\dff 2} 
\trf\bigr)
\off =\off
\gamma\trf (\dff u_{\dff 1}\dff,\off v_{\dff 1}\trf)
\off -\off
\gamma\trf (\dff u_{\dff 2}\dff,\off v_{\dff 2}\trf) 
\qff.
\]

\vspace{-12pt}\vspace{1.5pt}
Let\sss us\sss return\sss to\sss the original\sss interpretation of\dss 
$P^{\trf \sa}_{\dff 1}\qff \bigsqcup\pff P^{\dff \sa}_{\trf 2}$\dss 
as an operator over\sss $X_{\dff 1}\dff \bigsqcup\qff X_{\dff 2}$\nsp.\oss
Then vanishing of\trs $B_{\trf 1}$ on a section of\dss
$F_{\dff 1}\dff \bigsqcup\qff F_{\dff 2}$\sss
means simply\sss that\sss the restrictions of\dss this section\sss to\sss
$X_{\dff 1}$\sss and\sss $X_{\dff 2}$\sss agree on\sss the boundaries
and\sss hence define a section of\dss $\widehat{F}$\sss over\sss $\widehat{X}$\nnsp.\oss
In\sss particular,\oss restricting\dss 
$P^{\trf \sa}_{\dff 1}\qff \bigsqcup\pff P^{\dff \sa}_{\trf 2}$\dss
to\sss the kernel\dss $\kernel\trf B_{\trf 1}$\sss almost\sss
recovers\sss the operator\sss $\widehat{P}^{\qff \sa}$\dnsp.\oss
The only\sss problem\dss is\dss that\sss the sections over\sss $\widehat{X}$\sss
defined\sss by\sss elements of\dss $\kernel\trf B_{\trf 1}$\sss
are,\oss in\sss general,\oss only\sss continuous on $Y$ in\sss the directions\sss
transverse\sss to $Y$\dnsp,\oss even\sss
if\dss we consider our operators as acting\sss in\trs Sobolev\dss spaces
ensuring\sss existence of\dss several\sss derivatives.\oss

But\sss we are interested only\sss in\sss the analytical\dss index of\dss
$P^{\trf \sa}_{\dff 1}\qff \bigsqcup\pff P^{\dff \sa}_{\trf 2}$\dss
with\sss the boundary conditions\sss $B_{\trf 1}$\nsp.\oss
Of\dss course,\oss we implicitly assume\sss that\sss the operators\dss
$P^{\dff \mathrm{i}\qff \sa}_{\fff 1}$\sss and\dss $\widehat{P}^{\qff \mathrm{b}\qff \sa}$\sss 
and\dss hence\sss $P^{\trf \sa}_{\dff 1}\qff \bigsqcup\pff P^{\dff \sa}_{\trf 2}$\dss 
continuously\sss depend on a parameter.\oss
By\sss the results of\qss \cite{i2}\qss
the analytical\dss index\dss is\dss determined\sss by a somewhat\sss enhanced\sss
family of\dss kernels\vspace{3pt}
\[
\quad
\kernel\qff
\left(\qff
\left(\qff
P^{\trf \sa}_{\dff 1}\qff \bigsqcup\pff P^{\dff \sa}_{\trf 2}
\qff\right)
\qff \oplus\pff
B_{\trf 1}
\qff\right)
\qff.
\]

\vspace{-12pt}\vspace{3pt}
The enhancement\dss is\dss concerned\sss with\sss the eigenvectors of\dss
the operator\dss
$P^{\trf \sa}_{\dff 1}\qff \bigsqcup\pff P^{\dff \sa}_{\trf 2}$\dss
restricted\sss to\sss the kernel\dss $\kernel \trf B_{\trf 1}$\nsp,\oss
i.e.\qss with\sss the kernels\vspace{3pt}
\[
\quad
\kernel\qff
\left(\qff
\left(\qff
P^{\trf \sa}_{\dff 1}\qff \bigsqcup\pff P^{\dff \sa}_{\trf 2}
\off -\off
\lambda
\qff\right)
\qff \oplus\pff
B_{\trf 1}
\qff\right)
\qff
\]

\vspace{-12pt}\vspace{3pt}
with\sss $\lambda\qff \in\qff \rrr$\nnsp.\oss
By\sss the elliptic regularity\sss the sections of\trs 
$\widehat{F}$ over $\widehat{X}$
defined\dss by elements of\dss these kernels 
are $C^{\dff \infty}$\dnsp-smooth.\oss
It\dss follows\sss that\sss the eigenvalues and\sss eigenspaces defining\sss the analytical\dss index of\dss
the family\sss of\dss operators\dss 
$P^{\trf \sa}_{\dff 1}\qff \bigsqcup\pff P^{\dff \sa}_{\trf 2}$\dss
with\sss the boundary conditions $B_{\trf 1}$\sss
are exactly\sss the same as\sss the eigenvalues and\sss eigenspaces defining\sss the analytical\dss index of\dss
the family\sss of\dss operators\dss $\widehat{P}^{\qff \sa}$\dnsp.\oss
Therefore\sss these\sss two families have\sss the same analytical\dss index.\oss

Since\sss $B_{\qff \tau}$\sss continuously depends on\sss $\tau$\nnsp,\oss
the analytical\dss index of\dss
$P^{\trf \sa}_{\dff 1}\qff \bigsqcup\pff P^{\dff \sa}_{\trf 2}$\dss
with\sss the boundary conditions $B_{\qff \tau}$\sss is\dss
independent\sss of\dss $\tau$\nnsp.\oss
In\sss particular,\oss the analytical\dss index for\sss $\tau\off =\off 1$\sss
is\dss the same as for\sss $\tau\off =\off 0$\nnsp.\oss
It\sss follows\sss that\sss the analytical\dss index of\dss
$P^{\trf \sa}_{\dff 1}$\dss
with\sss the boundary conditions\dss
$B^{\dff \sa}_{\dff 1}$\sss is\dss equal\sss to\sss the analytic\sss index
of\trs $\widehat{P}^{\qff \sa}$\dnsp.\oss

\myuppar{The\sss topological\dss index.}
Our next\sss task\dss is\dss to prove\sss that\sss the\sss topological\sss index of\dss
$P^{\trf \sa}_{\dff 1}$\dss
with\sss the boundary conditions\dss
$B^{\dff \sa}_{\dff 1}$\sss is\dss equal\sss to\sss the\sss topological\sss index
of\trs $\widehat{P}^{\qff \sa}$\dnsp.\oss

The symbols\sss $\sigma_{\dff 1}$\sss and\sss $\sigma_{\dff 2}$\sss 
of\dss the operators\sss
$P^{\trf \sa}_{\dff 1}$\sss and\sss $P^{\dff \sa}_{\trf 2}$\sss respectively
are bundle-like in\sss the sense of\trs Section\qss \ref{symbols-conditions},\oss
and\dss hence\sss the extensions\sss
$\mathbb{E}^{\dff +}\fff (\trf \sigma_{\dff 1}\trf)$\sss 
and\sss
$\mathbb{E}^{\dff +}\fff (\trf \sigma_{\dff 2}\trf)$\sss
of\dss the bundles\sss
$E^{\dff +}\fff (\trf \sigma_{\dff 1}\trf)$\sss 
and\sss
$E^{\dff +}\fff (\trf \sigma_{\dff 2}\trf)$\nnsp.\oss
These extensions depend only on symbols\sss $\sigma_{\dff 1}\dff,\off \sigma_{\dff 2}$\nsp,\oss
but\sss not\sss on\sss the boundary conditions.\oss
Since\sss the restrictions of\dss these symbols\sss to\sss\vspace{1.5pt}
\[
\quad
B\dff X_{\dff 1\dff Y}
\off =\off 
B\dff X_{\dff 2\trf Y}
\off =\off 
B\dff X_{\trf Y}
\]

\vspace{-12pt}\vspace{1.5pt}
are\sss the same,\oss the restrictions of\dss the bundles\sss
$\mathbb{E}^{\dff +}\fff (\trf \sigma_{\dff 1}\trf)$\sss 
and\sss
$\mathbb{E}^{\dff +}\fff (\trf \sigma_{\dff 2}\trf)$\sss
are equal.\oss
By\sss identifying\sss these bundles over\sss $B\dff X_{\trf Y}$\sss
we get\sss a bundle\sss $\mathbb{E}^{\dff +}$\sss over\sss\vspace{3pt}\vspace{0.5pt}
\[
\quad
S\dff X_{\dff 1}\trf \cup\trf S\dff X_{\dff 2}\trf \cup\trf B\dff X_{\trf Y}
\off =\off
S\trf \widehat{X}\qff \cup\trf B\dff X_{\trf Y}
\qff.
\]

\vspace{-12pt}\vspace{3pt}\vspace{0.5pt}
The restriction $\mathbb{E}^{\dff +}\qff |\qff S\trf \widehat{X}$\sss
is\dss the bundle $E^{\dff +}\fff(\trf \widehat{\sigma}\trf)$\sss 
associated\sss with\sss the symbol\sss $\widehat{\sigma}$\sss 
of\trs $\widehat{P}^{\qff \sa}$\dnsp.\oss
Let\vspace{3pt}\vspace{0.5pt}
\[
\quad
e^{\dff +}
\pff \in\pff
K^{\dff 0}\dff (\trf S\trf \widehat{X}\qff \cup\trf B\dff X_{\trf Y}\trf)
\]

\vspace{-12pt}\vspace{3pt}\vspace{0.5pt}
be\sss the class of\sss $\mathbb{E}^{\dff +}$\dnsp,\oss
and\dss let\sss $\varepsilon^{\dff +}$\sss be\sss its\sss 
image under\sss the coboundary\sss map\vspace{3pt}\vspace{0.5pt}
\[
\quad
K^{\dff 0}\dff (\trf S\trf \widehat{X}\qff \cup\trf B\dff X_{\trf Y}\trf)
\qff \ttoo\qff
K^{\dff 1}\dff (\trf B\trf \widehat{X}\fff,\off  S\trf \widehat{X}\qff \cup\trf B\dff X_{\trf Y}\trf)
\qff.
\]

\vspace{-12pt}\vspace{3pt}\vspace{0.5pt}
Let\sss us consider\sss the commutative diagram\vspace{1.5pt}
\[
\quad
\begin{tikzcd}[column sep=sma, row sep=sma]
K^{\dff 0}\dff (\trf S\trf \widehat{X}\qff \cup\trf B\dff X_{\trf Y}\trf)
\arrow[r]
\arrow[d]
&
K^{\dff 0}\dff (\trf S\dff X_{\dff 1}\dff \cup\dff B\dff X_{\dff 1\dff Y}\trf)
\qff \oplus\qff
K^{\dff 0}\dff (\trf S\dff X_{\dff 2}\dff \cup\dff B\dff X_{\dff 2\trf Y}\trf)
\arrow[d]
\\
K^{\dff 1}\dff (\trf B\trf \widehat{X}\fff,\off  S\trf \widehat{X}\qff \cup\trf B\dff X_{\trf Y}\trf)
\arrow[r]
&
K^{\dff 1}\dff (\trf B\dff X_{\dff 1}\fff,\off S\dff X_{\dff 1}\dff \cup\dff B\dff X_{\dff 1\dff Y}\trf)
\qff \oplus\qff
K^{\dff 1}\dff (\trf B\dff X_{\dff 2}\fff,\off S\dff X_{\dff 2}\dff \cup\dff B\dff X_{\dff 2\trf Y}\trf)
\qff,
\end{tikzcd}
\]

\vspace{-12pt}\vspace{1.5pt}
where\sss the horizontal\sss arrows are\sss the direct\sss sums of\dss the restriction\sss maps,\oss
and\sss the vertical\sss arrows are\sss the coboundary maps.\oss
The\sss lower\sss horizontal\sss arrow\dss is\dss an\sss isomorphism\sss
because\sss it\dss is\dss placed\dss between\sss the groups\sss
$K^{\dff 0}\dff (\trf B\dff X_{\trf Y}\fff,\qff B\dff X_{\trf Y}\trf)\off =\off 0$\sss
and\sss
$K^{\dff 1}\dff (\trf B\dff X_{\trf Y}\fff,\qff B\dff X_{\trf Y}\trf)\off =\off 0$\sss
in a\dss Mayer-Vietoris\dss sequence.\oss
The classes $e^{\dff +}\fff,\off 
e^{\dff +}\fff (\trf \sigma_{\dff 1}\trf)\fff,\off 
e^{\dff +}\fff (\trf \sigma_{\dff 2}\trf)$\nnsp,\oss etc. 
are mapped\dss by\sss the arrows of\dss this diagram as\sss in\sss the following diagram.\oss\vspace{1.5pt}
\[
\quad
\begin{tikzcd}[column sep=boomm, row sep=sma]
e^{\dff +}
\arrow[r, mapsto]
\arrow[d, mapsto]
&
e^{\dff +}\fff (\trf \sigma_{\dff 1}\trf)
\qff \oplus\qff
e^{\dff +}\fff (\trf \sigma_{\dff 2}\trf)
\arrow[d, mapsto]
\\
\varepsilon^{\dff +}
\arrow[r, mapsto]
&
\varepsilon^{\dff +}\fff (\trf \sigma_{\dff 1}\trf)
\qff \oplus\qff
\varepsilon^{\dff +}\fff (\trf \sigma_{\dff 2}\trf)
\qff.
\end{tikzcd}
\]

\vspace{-12pt}\vspace{1.5pt}
Since\sss the\sss lower\sss horizontal\sss arrow\dss is\dss an\sss isomorphism,\oss
the class\sss $\varepsilon^{\dff +}$\sss is\dss uniquely determined\sss by\sss
$\varepsilon^{\dff +}\fff (\trf \sigma_{\dff 1}\trf)$\sss
and\sss
$\varepsilon^{\dff +}\fff (\trf \sigma_{\dff 2}\trf)$\nnsp.\oss
Clearly,\oss the image of\dss $\varepsilon^{\dff +}$\sss under\sss the map\vspace{3pt}
\[
\quad
K^{\dff 1}\dff 
\left(\trf
B\trf \widehat{X}\fff,\off  S\trf \widehat{X}\qff \cup\trf B\dff X_{\trf Y}
\trf\right)
\off \ttoo\off
K^{\dff 1}\dff 
\left(\trf 
B\trf \widehat{X}\fff,\off  S\trf \widehat{X}
\trf\right)
\qff
\]

\vspace{-12pt}\vspace{3pt}
is\dss nothing else but\sss $\varepsilon^{\dff +}\fff (\trf \widehat{\sigma}\trf)$\nnsp.\oss
It\sss follows\sss that\sss the class\sss $\varepsilon^{\dff +}\fff (\trf \widehat{\sigma}\trf)$\sss
is\dss uniquely determined\dss by\sss the classes\sss
$\varepsilon^{\dff +}\fff (\trf \sigma_{\dff 1}\trf)$\sss
and\sss
$\varepsilon^{\dff +}\fff (\trf \sigma_{\dff 2}\trf)$\nnsp.\oss

In order\sss to pass from\sss these classes\sss to\sss topological\dss indices,\oss
let\sss us embed\dss the manifold\sss $\widehat{X}$\sss into\sss
$\rrr^{\dff n}$\sss for some $n$ in such a way\sss that\sss\vspace{3pt}
\[
\quad
X_{\dff 1}
\off =\off
\widehat{X}\qff \cap\qff \rrr_{\qff \geq\trf 0}^{\dff n}
\quad
\mbox{and}\quad
X_{\dff 2}
\off =\off
\widehat{X}\qff \cap\qff \rrr_{\qff \leq\trf 0}^{\dff n}
\qff
\]

\vspace{-12pt}\vspace{3pt}
and\sss $X_{\dff 2}$\sss is\dss the reflection of\dss $X_{\dff 1}$\sss in\sss the hyperplane\sss
$\rrr^{\dff n\dff -\dff 1}\dff \times\dff 0$\nnsp.\oss
Let\sss $\mathbb{N}$\sss be\sss the normal\dss bundle\sss 
to\sss $\widehat{X}$\sss in\sss $\rrr^{\dff n}$\dnsp.\oss
The normal\dss bundle of\dss $T\dff \widehat{X}$\sss 
in\sss 
$T\dff \rrr^{\dff n}
\off =\off
\rrr^{\dff n}
\dff \times\dff
\rrr^{\dff n}$\sss
can\sss be identified\sss with\sss the\sss lift\sss of\dss the bundle\sss
$\mathbb{N}\trf \oplus\trf \mathbb{N}$\sss 
to\sss $T\dff \widehat{X}$\nnsp.\oss
As in\trs Section\qss \ref{symbols-conditions},\oss 
the bundle\sss $\mathbb{N}\trf \oplus\trf \mathbb{N}$\sss
has a natural\sss complex structure,\oss
and a\sss tubular\sss
neighborhood of\dss $T\dff \widehat{X}$\sss 
in\sss $T\dff \rrr^{\dff n}$\sss
can\sss be identified\sss with\sss the bundle\sss $\mathbb{U}$\sss of\dss 
unit\sss balls in\sss 
$\mathbb{N}\trf \oplus\trf \mathbb{N}$\nnsp.\oss
Let\sss $\mathbb{S}$\sss be\sss the bundle of\dss unit\sss spheres in\sss
$\mathbb{N}\trf \oplus\trf \mathbb{N}$\nnsp.\oss

The complex structure on\sss 
$\mathbb{N}\trf \oplus\trf \mathbb{N}$\dss
leads\sss to\trs Thom\dss isomorphisms\vspace{4.5pt}
\[
\quad
K^{\dff 1}\dff 
\bigl(\trf 
B\dff X_{\dff 1}\fff,\off  S\dff X_{\dff 1}\dff \cup\dff B\dff X_{\trf Y}
\trf\bigr)
\qff \ttoo\qff
K^{\dff 1}\dff 
\bigl(\qff \mathbb{U}\dff |\trf B\dff X_{\dff 1}\fff,\off 
(\trf \mathbb{S}\dff |\trf B\dff X_{\dff 1}\trf)
\qff \cup\qff 
(\trf \mathbb{U}\dff |\trf S\dff X_{\dff 1}\trf)
\qff \cup\qff 
(\trf \mathbb{U}\dff |\trf B\dff X_{\trf Y}\trf)\qff\bigr)
\qff,
\]

\vspace{-30pt}
\[
\quad
K^{\dff 1}\dff 
\bigl(\trf B\dff X_{\dff 2}\fff,\off  S\dff X_{\dff 2}\dff \cup\dff B\dff X_{\trf Y}
\trf\bigr)
\qff \ttoo\qff
K^{\dff 1}\dff \bigl(\qff \mathbb{U}\dff |\trf B\dff X_{\dff 2}\fff,\off 
(\trf \mathbb{S}\dff |\trf B\dff X_{\dff 2}\trf)
\qff \cup\qff 
(\trf \mathbb{U}\dff |\trf S\dff X_{\dff 2}\trf)
\qff \cup\qff 
(\trf \mathbb{U}\dff |\trf B\dff X_{\trf Y}\trf)\qff\bigr)
\qff,
\]

\vspace{-30pt}
\[
\quad
K^{\dff 1}\dff 
\bigl(\trf 
B\dff \widehat{X}\fff,\off  S\dff \widehat{X}\dff \cup\dff B\dff X_{\trf Y}
\trf\bigr)
\qff \ttoo\qff
K^{\dff 1}\dff \bigl(\qff 
\mathbb{U}\dff |\trf B\dff \widehat{X}\fff,\off 
(\trf \mathbb{S}\dff |\trf B\dff \widehat{X}\trf)
\qff \cup\qff 
(\trf \mathbb{U}\dff |\trf S\dff \widehat{X}\trf)
\qff \cup\qff 
(\trf \mathbb{U}\dff |\trf B\dff X_{\trf Y}\trf)\qff\bigr)
\qff,
\]

\vspace{-30pt}
\[
\quad
K^{\dff 1}\dff 
\bigl(\trf B\dff \widehat{X}\fff,\off  S\dff \widehat{X}
\qff\bigr)
\qff \ttoo\qff
K^{\dff 1}\dff \bigl(\qff \mathbb{U}\dff |\trf B\dff \widehat{X}\fff,\off 
(\trf \mathbb{S}\dff |\trf B\dff \widehat{X}\trf)
\qff \cup\qff 
(\trf \mathbb{U}\dff |\trf S\dff \widehat{X}\trf)
\qff\bigr)
\qff.
\]

\vspace{-12pt}\vspace{4.5pt}
These\trs Thom\dss isomorphisms commute with various restriction\sss maps,\oss
and\dss the analogue of\dss the\sss lower\sss horizontal\sss arrow\sss in\sss the first\sss
diagram\dss is\dss an\sss isomorphism\sss by\sss the same reason as\sss that\sss arrow.\oss
It\sss follows\sss that\sss the image\sss
$\mathrm{Th}\trf(\trf \varepsilon^{\dff +}\trf)$\sss 
of\dss the class\sss 
$\varepsilon^{\dff +}$\sss
under\sss the\sss third\trs Thom\dss isomorphism above\dss
is\dss uniquely determined\dss by\sss the images of\dss the classes\sss
$\varepsilon^{\dff +}\fff (\trf \sigma_{\dff 1}\trf)$\sss
and\sss
$\varepsilon^{\dff +}\fff (\trf \sigma_{\dff 2}\trf)$\sss
under\sss the first\sss two\trs Thom\dss isomorphisms.\oss
Similarly\sss to\sss the class
$\varepsilon^{\dff +}\fff (\trf \widehat{\sigma}\trf)$ itself,\oss
the image of\sss $\varepsilon^{\dff +}\fff (\trf \widehat{\sigma}\trf)$ 
under\sss the\sss last\trs Thom\dss isomorphism\dss
is\dss equal\dss to\sss the image of\dss
$\mathrm{Th}\trf(\trf \varepsilon^{\dff +}\trf)$\sss
under\dss the map\vspace{3pt}
\[
\quad
K^{\dff 1}\dff 
\bigl(\qff 
\mathbb{U}\dff |\trf B\dff \widehat{X}\fff,\off 
(\trf \mathbb{S}\dff |\trf B\dff \widehat{X}\trf)
\qff \cup\qff 
(\trf \mathbb{U}\dff |\trf S\dff \widehat{X}\trf)
\qff \cup\qff 
(\trf \mathbb{U}\dff |\trf B\dff X_{\trf Y}\trf)\qff\bigr)
\qff \ttoo\qff
K^{\dff 1}\dff \bigl(\qff \mathbb{U}\dff |\trf B\dff \widehat{X}\fff,\off 
(\trf \mathbb{S}\dff |\trf B\dff \widehat{X}\trf)
\qff \cup\qff 
(\trf \mathbb{U}\dff |\trf S\dff \widehat{X}\trf)
\qff\bigr)
\qff
\]

\vspace{-12pt}\vspace{3pt}
induced\dss by\sss the quotient\sss map\vspace{3pt}
\[
\quad
q\dff \colon\dff
\mathbb{U}\dff |\trf B\dff \widehat{X}\qff\bigl/\dff
\bigl(\qff 
(\trf \mathbb{S}\dff |\trf B\dff \widehat{X}\trf)
\qff \cup\qff 
(\trf \mathbb{U}\dff |\trf S\dff \widehat{X}\trf)
\qff\bigr)
\qff \ttoo\qff
\mathbb{U}\dff |\trf B\dff \widehat{X}\qff\bigl/\dff
\bigl(\qff 
(\trf \mathbb{S}\dff |\trf B\dff \widehat{X}\trf)
\qff \cup\qff 
(\trf \mathbb{U}\dff |\trf S\dff \widehat{X}\trf)
\qff \cup\qff 
(\trf \mathbb{U}\dff |\trf B\dff X_{\trf Y}\trf)\qff\bigr)
\qff.
\]

\vspace{-12pt}\vspace{3pt}
The composition of\dss the canonical\dss map\vspace{3pt}
\[
\quad
S^{\dff 2 n}
\qff \ttoo\qff
\mathbb{U}\dff |\trf B\dff \widehat{X}\qff\bigl/\trf
\bigl(\qff 
\mathbb{S}\qff \cup\qff 
(\trf \mathbb{U}\dff |\trf S\dff \widehat{X}\trf)
\qff\bigr)
\qff
\]

\vspace{-12pt}\vspace{3pt}
with\sss $q$\sss
is\dss equal\dss to\sss the sum\qss
(as in\sss the definition of\trs homotopy\sss groups)\qss 
of\dss the maps\vspace{3pt}
\[
\quad
S^{\dff 2 n}
\qff \ttoo\qff
\mathbb{U}\dff |\trf B\dff X_{\dff 1}\qff\bigl/\trf
\bigl(\qff 
\mathbb{S}\qff \cup\qff 
(\trf \mathbb{U}\dff |\trf S\dff X_{\dff 1}\trf)
\qff \cup\qff 
(\trf \mathbb{U}\dff |\trf B\dff X_{\trf Y}\trf)\qff\bigr)
\quad
\mbox{and}\quad
\qff
\]

\vspace{-33pt}\vspace{-1.2pt}
\[
\quad
S^{\dff 2 n}
\qff \ttoo\qff
\mathbb{U}\dff |\trf B\dff X_{\dff 2}\qff\bigl/\trf
\bigl(\qff 
\mathbb{S}\qff \cup\qff 
(\trf \mathbb{U}\dff |\trf S\dff X_{\dff 2}\trf)
\qff \cup\qff 
(\trf \mathbb{U}\dff |\trf B\dff X_{\trf Y}\trf)\qff\bigr)
\qff.
\]

\vspace{-12pt}\vspace{3pt}
Rotation of\dss the half-space $\rrr_{\qff \leq\trf 0}^{\dff n}$ into\sss the
half-space $\rrr_{\qff \geq\trf 0}^{\dff n}$\sss
leads\sss to a deformation of\dss the second\sss map\sss 
to\sss the map defining\sss the index of\dss $\sigma_{\dff 2}$\nsp.\oss
It\dss follows\sss that\vspace{1.5pt}
\begin{equation}
\label{t-additivity}
\quad
\ti\qff (\trf \widehat{\sigma}\qff)
\off =\off\dff
\ti\qff (\trf \sigma_{\dff 1}\fff,\qff N_{\dff 1} \trf)
\qff +\pff
\ti\qff (\trf \sigma_{\dff 2}\fff,\qff N_{\dff 2} \trf)
\qff,
\end{equation}

\vspace{-12pt}\vspace{1.5pt}
where\sss $N_{\dff 1}\dff,\off N_{\dff 2}$\sss
are\sss the boundary conditions corresponding\dss to\dss
$B^{\dff \sa}_{\dff 1}\dff,\off\qff B^{\dff \sa}_{\dff 2}$\sss
respectively.\oss
In\sss fact,\oss since\sss the symbols\sss
$\sigma_{\dff 1}\fff,\off \sigma_{\dff 2}$\sss
are already\dss bundle-like,\oss
the corresponding\sss topological\dss indices are well\sss
defined\sss without\sss boundary conditions.\oss
As usual,\oss the additivity\dss property\qss (\ref{t-additivity})\qss
is\dss valid also when\sss parameters are present\qss
(and\dss is\dss non-vacuous only\sss with\sss parameters).\oss

We still\dss not\sss used\dss the assumption\sss that\sss $\sigma_{\dff 2}$\sss
is\dss the symbol\sss of\dss a standard operator.\oss
In\sss particular,\oss the additivity\sss property\qss (\ref{t-additivity})\qss
holds when\sss $\sigma_{\dff 2}$\sss is\dss a\sss general\sss symbol\sss
of\dss the same\sss type as\sss $\sigma_{\dff 1}$\nsp.\oss
When\sss $\sigma_{\dff 2}$\sss is\dss the symbol\sss of\dss a standard operator,\oss
$\varepsilon^{\dff +}\fff (\trf \sigma_{\dff 2}\trf)\off =\off 0$\nnsp,\oss
as we saw in\dss Section\qss \ref{symbols-conditions}.\oss
It\sss follows\sss that\sss in\sss this case\dss
$\ti\qff (\trf \sigma_{\dff 2}\fff,\qff N_{\dff 2} \trf)\off =\off 0$\sss
and\dss hence\vspace{1.5pt}
\[
\quad
\ti\qff (\trf \widehat{\sigma}\qff)
\off =\off\dff
\ti\qff (\trf \sigma_{\dff 1}\fff,\qff N_{\dff 1} \trf)
\qff.
\]

\vspace{-12pt}\vspace{1.5pt}
In other\sss words,\oss
the\sss topological\dss index of\dss $P^{\trf \sa}_{\dff 1}$\dss
with\sss the boundary conditions\dss
$B^{\dff \sa}_{\dff 1}$\sss is\dss equal\sss to\sss the\sss topological\dss index
of\trs $\widehat{P}^{\qff \sa}$\dnsp.\oss
To sum up,\oss for\sss families of\dss operators and\dss boundary conditions having\sss
the standard\sss form\sss near\sss the boundary\sss
the determination of\dss both\sss the analytical\dss and\sss the\sss
topological\dss index can\sss be reduced\dss to\sss 
corresponding\sss problem\sss on closed\sss manifolds.\oss

\newpage
\mysection{The\qss index\qss theorem}{index-theorems}

\myuppar{The framework.}
We will\sss speak about\sss various objects parameterized\dss by a\sss
topological\sss space $Z$\nnsp.\oss
Almost\sss all\sss arguments\sss work\sss
when\sss $Z$\sss is\dss compactly\sss generated and\sss paracompact.\oss
For example,\oss every\sss metric space has\sss these properties.\oss
In\sss particular,\oss we will\sss consider manifolds\sss
$X\trf(\trf z\trf)\fff,\off z\qff \in\qff Z$\sss
continuously depending on $z$\sss in\sss the sense\sss that\sss
$X\trf(\trf z\trf)$\sss is\dss the fiber over $z$ of\dss a\sss
locally\sss trivial\dss bundle,\oss as at\sss the end of\trs Section\qss \ref{symbols-conditions}.\oss
In order\sss to avoid cumbersome notations,\oss
we will\sss usually omit\sss the parameter $z$ and speak
about\sss $X$\nnsp,\oss vector bundles and\sss pseudo-differential\sss
operators on $X$\nnsp,\oss etc.\oss having\sss in\sss mind\sss
families parameterized\dss by $z\qff \in\qff Z$\nnsp.\oss

Let\sss $P$\sss be a  pseudo-differential\dss operator in a bundle\sss $E$ over\sss
$X$\sss belonging\dss to\sss the\dss H\"{o}rmander\dss class,\oss
and\dss let\sss $\sigma$\sss be\sss
the symbol\sss of\dss $P$\dnsp.\oss
Suppose\sss that\sss $\sigma$\sss is\dss elliptic and self-adjoint.\oss
Let\sss $N$\sss be an elliptic self-adjoint\dss boundary condition\sss
for\sss $\sigma$\sss in\sss the sense of\trs Section\qss \ref{symbols-conditions}.\oss
We assume\sss that\sss $N$\sss is\dss bundle-like.\oss
Then we can consider $N$ as a subbundle of\dss the restriction\sss
$E\trf |\trf Y$\nnsp,\oss where\sss $Y\off =\off \partial\trf X$\nnsp.\oss
The symbol\sss $\sigma$ of\dss $P$\sss induces a self-adjoint\sss endomorphism\sss
$\Sigma\dff \colon\dff E\trf |\trf Y\qff \ttoo\qff E\trf |\trf Y$\dnsp.\oss
Recall\dss that\sss $E^{\dff +}_{\trf Y}$\sss and\sss $E^{\dff -}_{\trf Y}$\sss
are\sss the subbundles of\dss $E\trf |\trf Y$\sss
generated\dss by\sss the eigenvectors of\dss $\Sigma$\sss with\sss
positive and\sss negative eigenvalues respectively.\oss

Let\dss
$\Pi\dff \colon\dff
E\trf |\trf Y
\qff \ttoo\qff
E\trf |\trf Y$\sss
be\sss the orthogonal\dss projection onto\sss $N$\nnsp,\oss
and\dss let\sss
$\Gamma\off =\off (\trf 1\qff -\qff \Pi\trf)\dff \circ\dff \gamma$\nnsp,\oss
where $\gamma$\sss takes\sss sections of\dss
$E$\sss to\sss their restrictions\sss to $Y$\dnsp.\oss
We may consider\sss $\Gamma$ as an operator\sss taking values\sss
in section of\dss the orthogonal\sss complement\sss $N^{\dff \perp}$
of\dss $N$\sss in\sss $E\trf |\trf Y$\dnsp.\oss
Then\sss the operators\sss
$P\dff \oplus\dff \Gamma$\dnsp,\oss
considered as operators\sss in appropriate\dss Sobolev\dss spaces,\oss
are\dss Fredholm\qss (and continuously depend on $z$\nnsp).\oss
The results of\trs Section\qss \ref{abstract-index}\qss show\sss that\dss
by restricting\sss $P$\sss to\sss $\kernel\dff \Gamma$\sss we get\sss 
an\sss unbounded self-adjoint\sss operator\sss
$P_{\trf \Gamma}$\nsp,\oss
and\dss that\sss the analytical\dss index of\dss $P_{\trf \Gamma}$\sss is\dss well-defined.\oss
We will\sss denote\sss this index\sss by\sss
$\ai\trf(\trf P\halfff,\qff N\trf)$\nnsp.\oss
Since,\oss in\sss fact,\oss everything depends on $z\qff \in\qff Z$\nnsp,\oss
the analytical\dss index\sss
$\ai\trf(\trf P\halfff,\qff N\trf)$\sss
is\dss an element\sss of\dss
$K^{\dff 1}\dff (\trf Z\trf)$\nnsp.\oss
Suppose now\sss that\sss $Z$\sss is\dss compact.\oss
Then\sss the\sss topological\dss index\sss 
$\ti\trf(\trf \sigma\fff,\qff N\trf)
\qff \in\qff 
K^{\dff 1}\dff (\trf Z\trf)$\sss
is\dss also defined.

\mypar{Theorem.}{index-theorem-ext}
\emph{If\trs the bundle\dss $E^{\dff +}_{\trf Y}$ extends\sss to a bundle over\sss $X$\sss
in a manner continuously depending on\dss the parameter\sss $z$\nnsp,\oss then}\qss
$\ai\trf(\trf P\halfff,\qff N\trf)
\off =\off
\ti\trf(\trf \sigma\fff,\qff N\trf)$\nnsp.\oss

\proof
We may assume\sss that\sss the operators\sss
$\mathbf{T}\trf(\trf x_{\dff n}\trf)$\sss from\dss Section\qss \ref{pdo}\qss
do not\sss depend on $x_{\dff n}$\sss for sufficiently small\sss $x_{\dff n}$\nsp.\oss
See\qss \cite{h},\oss the proof\dss of\trs Proposition\qss 20.3.3.\oss
Let\sss $\tau$\sss be\sss the symbol\sss of\dss $T\dff(\trf 0\trf)$\nnsp.\oss
The boundary condition\sss $N$\sss defines an\sss isometric\sss isomorphism\sss 
$\varphi\dff \colon\dff E^{\dff +}_{\trf Y}\qff \ttoo\qff E^{\dff -}_{\trf Y}$\nsp,\oss
and we will\sss identify\sss the bundles\sss
$E^{\dff +}_{\trf Y}$\sss and\sss $E^{\dff -}_{\trf Y}$\sss by $\varphi$\nnsp.\oss
By\sss the results of\trs Sections\qss \ref{boundary-algebra}\qss and\qss \ref{symbols-conditions}\qss
we can deform\sss $\sigma\fff,\qff N$\sss by a deformation\sss preserving all\sss our assumptions\sss
to a normalized\sss new pair\sss $\sigma\fff,\qff N$\nnsp.\oss
This deformation can\sss be\sss lifted\sss to a deformation of\dss $P$\dnsp.\oss
After\sss this deformation\sss $\sigma$\sss is\dss a bundle-like symbol\sss and,\oss moreover,\pss 
$\sigma,\pff N$\sss have a standard\sss form over $Y$\dnsp.\oss
Namely,\vspace{1.5pt}
\[
\quad
\Sigma
\off =\off\dff
\begin{pmatrix}
\off\dff 1 &
0 \off
\vspace{4.5pt} \\
\off\dff 0 &
-\qff 1 \off 
\end{pmatrix}
\quad
\mbox{and}\quad
\tau
\off =\off\dff
\begin{pmatrix}
\off 0 &
i \off
\vspace{4.5pt} \\
\off -\qff i &
0 \off 
\end{pmatrix}
\qff,
\]

\vspace{-12pt}\vspace{1.5pt}
and\sss $N\off =\off \Delta$\sss
with respect\sss to\sss the decomposition\sss
$E\trf |\trf Y
\off =\off 
E^{\dff +}_{\trf Y}\pff \oplus\pff E^{\dff -}_{\trf Y}$.\oss
See\qss (\ref{sigma-rho-tau-st}).\oss
Let\sss us pass\sss to\sss the decomposition\dss
$E\trf |\trf Y
\off =\off 
\Delta
\dff \oplus\dff 
\Delta^{\fff \perp}$\dnsp,\oss
where\sss the subbundles\sss
$\Delta$\sss
and\sss 
$\Delta^{\fff \perp}$
are defined\sss in\sss the obvious way,\oss
and\dss identify\sss these subbundles
with\sss $E^{\dff +}_{\trf Y}$\sss by\sss the projection.\oss
Then\vspace{1.5pt}\vspace{-0.25pt}
\[
\quad
\Sigma
\off =\off\dff
\begin{pmatrix}
\off 0 &
1 \off
\vspace{4.5pt} \\
\off 1 &
0 \off 
\end{pmatrix}
\quad
\mbox{and}\quad
\tau
\off =\off\dff
\begin{pmatrix}
\off 0 &
-\qff i \off
\vspace{4.5pt} \\
\off i &
0 \off 
\end{pmatrix}
\qff,
\]

\vspace{-12pt}\vspace{1.5pt}\vspace{-0.25pt}
and\sss $N$\sss is\dss equal\dss to\sss the first\sss summand\sss
$\Delta\dff \oplus\dff 0$\sss of\dss this decomposition.\oss
Now\sss the symbol\sss $\sigma$\sss is\dss equal\dss near $Y$\sss to an elementary symbol\sss
in\sss the sense of\trs Section\qss \ref{symbols-conditions}.\oss
More precisely,\pss $\sigma$\sss is\dss equal\dss to\sss the elementary symbol\sss $\sigma^{\dff \sa}$\sss
corresponding\sss to\sss $F^{\dff +}\off =\off \Delta$\nnsp,\qss
$F^{\dff -}\off =\off 0$\nnsp,\oss
and\sss $\lambda^{\fff +}\off =\off 1$\dss 
(near $Y$ we don't\sss need\sss the bundle $F$\nsp).\oss
The corresponding\sss boundary condition\sss $N^{\dff \sa}$\sss
is\dss equal\dss to\sss 
$F^{\dff +}\dff \oplus\dff 0
\off =\off
\Delta\dff \oplus\dff 0$\nnsp,\oss
i.e.\qss to\sss $N$\nnsp.\oss
Let\sss us\sss further deform\sss $P$\sss without\sss changing\sss its symbol\dss
to a new operator $P$ equal\sss near $Y$\sss 
to\sss a standard self-adjoint\sss operator $P^{\dff \sa}$\sss
from\dss Section\qss \ref{pdo}.\oss

The bundle $F^{\dff +}$ is\dss canonically\sss isomorphic\sss 
to\sss $E^{\dff +}_{\trf Y}$\sss and\dss hence
can\sss be extended\dss to a bundle\sss $F$\sss over\sss $X$\nnsp.\oss
Therefore\sss the constructions of\trs Section\qss \ref{glueing}\qss apply\sss
to $P$ in\sss the role of\sss $P^{\trf \sa}_{\dff 1}$
and\sss $N$\sss in\sss the role of\dss
the boundary condition\sss $B^{\dff \sa}_{\dff 1}$\sss from\sss that\dss section\qss
(more precisely,\oss one\sss takes as $B^{\dff \sa}_{\dff 1}$\sss the boundary operator defined\dss by $N$\nnsp).\oss 
In\sss more details,\oss a copy of\sss $E$\sss over $X_{\dff 1}$
plays\sss the role of\sss $F_{\dff 1}$\nsp,\pss
and a copy of\sss $F\dff \oplus\dff F$ over\sss $X_{\dff 2}$
plays\sss the role of\dss $F_{\dff 2}$\nsp.\oss 
The constructions of\trs Section\qss \ref{glueing}\pss lead\dss to
an operator\sss $\widehat{P}^{\qff \sa}$\sss on\sss the double\dss 
$\widehat{X}\off =\off X_{\dff 1}\qff \cup\qff X_{\dff 2}$\sss
of\dss $X$\sss such\sss that\vspace{3pt}\vspace{-0.25pt}
\[
\quad
\ai\trf(\trf P\halfff,\qff N\trf)
\off =\off
\ai\trf(\qff \widehat{P}^{\qff \sa}\qff)
\quad
\mbox{and}\quad
\ti\trf(\trf P\halfff,\qff N\trf)
\off =\off
\ti\trf(\qff \widehat{P}^{\qff \sa}\qff)
\qff.
\]

\vspace{-12pt}\vspace{3pt}\vspace{-0.25pt}
But\sss $\widehat{X}$\sss is\dss a closed\sss manifold,\oss
and\sss for closed\sss manifolds our definitions of\dss
the\sss analytical\sss and\dss the\sss topological\dss index
agree with\sss the classical\sss definitions of\qss
Atiyah-Patodi-Singer\qss \cite{aps}.\oss
For\sss the\sss topological\sss index\sss this\dss is\dss obvious,\oss
and\sss for\sss the analytical\dss index\sss this follows\sss from\qss
\cite{i2},\oss Theorem\qss 8.5,\oss
and\sss the fact\sss that\sss on closed\sss manifolds\sss
families of\dss elliptic self-adjoint\dss pseudo-dif\-fer\-en\-tial\sss operators
are\qss \emph{strictly\trs Fredholm}\pss in\sss the sense of\qss \cite{i2}.\oss
This fact\sss follows\sss from\sss the results of\trs Seeley\qss \cite{se1}.\oss
By\sss the results of\trs 
Atiyah,\qss Patodi,\qss and\dss Singer\vspace{1.5pt}
\[
\quad
\ai\trf(\qff \widehat{P}^{\qff \sa}\qff)
\off =\off
\ti\trf(\qff \widehat{P}^{\qff \sa}\qff)
\qff.
\]

\vspace{-12pt}\vspace{1.5pt}
See\qss \cite{aps},\oss Theorem\qss (3.4),\oss
and\sss the discussion\sss below.\oss
It\sss follows\sss that\vspace{1.5pt}
\[
\quad
\ai\trf(\trf P\halfff,\qff N\trf)
\off =\off
\ti\trf(\trf P\halfff,\qff N\trf)
\qff.
\]

\vspace{-12pt}\vspace{1.5pt}
Since\sss the\sss topological\sss and\dss analytical\dss indices are defined\sss
for families,\oss they are\sss invariant\sss under\sss homotopies.\oss
Hence our deformations of\dss $P\fff,\pff N$\sss do not\sss change neither\sss the\sss topological\sss
nor\sss the analytical\dss index.\oss
Since we were implicitly\sss talking about\dss families,\oss 
the\sss theorem\sss follows.\oss  \eproof

\myuppar{The index\sss theorem\sss for self-adjoint\sss operators on closed\sss manifolds.}
The\sss index\sss theorem\sss for\sss the families of\dss self-adjoint\sss
operators on closed\sss manifolds\sss is\dss stated\dss by\sss
Atiyah,\qss Patodi,\qss and\dss Singer\qss \cite{aps}\qss only\sss
for families of\dss operators parameterized\sss by a compact\sss space 
on a\sss fixed\sss manifold\sss $X$\nnsp.\oss
The assumption\sss that\sss the manifold\sss $X$\sss is\dss fixed\dss
is\dss not\sss essential,\oss and\dss the\sss tools developed\dss by\trs
Atiyah\dss and\dss Singer\qss \cite{as4}\qss allow\sss to extend\dss
this\sss index\sss theorem\sss to\sss the families of\dss closed\sss manifolds
parameterized\sss by a compact\sss space.\oss

A natural\sss way\sss to 
deal\sss with\sss more\sss general\sss spaces of\dss parameters\dss
is\dss to use\sss author's definition of\dss the analytical\dss index\qss \cite{i2}.\oss 
The proof\trs in\qss \cite{aps}\qss is\dss based\sss on\sss 
the\dss Atiyah-Singer\qss \cite{as}\qss correspondence between\sss 
self-adjoint\trs Fredholm\dss operators and\sss loops of\trs Fredholm\dss operators,\oss
which\dss is,\oss in\sss fact,\oss build\sss into\sss the 
definition of\dss the analytical\dss index\sss in\qss \cite{aps}.\oss
If\dss the definition of\pss \cite{i2}\qss is\dss used,\oss this\sss tool\sss should\sss be
replaced\sss by\dss Theorem\qss 8.5\qss from\qss \cite{i2}.\oss
The rest\dss is\dss routine.\oss

\mypar{Corollary.}{index-theorem-two}
\emph{\dnsp$\ai\trf(\trf P\halfff,\qff N\trf)
\off -\off
\ti\trf(\trf \sigma\fff,\qff N\trf)$\sss
is\dss an\sss element\sss of\dss order\sss $2$\sss in\sss $K^{\dff 1}\dff (\trf Z\trf)$\nnsp.}

\proof
Let\sss us consider\sss the bundle\sss $E\dff \oplus\dff E$\sss over\sss $X$\sss
and\sss the operator\sss $P\dff \oplus\dff P$\sss in\sss this bundle\sss
together with\sss the boundary condition\sss $N\dff \oplus\dff N$\nnsp.\oss
Clearly,\vspace{3pt}
\[
\quad
\ai\trf(\trf P\dff \oplus\dff P\halfff,\qff N\dff \oplus\dff N\trf)
\off =\off
2\dff \ai\trf(\trf P\halfff,\qff N\trf)
\quad
\mbox{and}\quad
\ti\trf(\trf P\dff \oplus\dff P\halfff,\qff N\dff \oplus\dff N\trf)
\off =\off
2\dff \ti\trf(\trf P\halfff,\qff N\trf)
\qff.
\]

\vspace{-12pt}\vspace{3pt}
Also,\pss 
$(\trf E\dff \oplus\dff E\trf)^{\dff +}_{\trf Y}
\off =\off
E^{\dff +}_{\trf Y}\trf \oplus\dff E^{\dff +}_{\trf Y}$\nsp.\oss
The boundary condition\sss $N$\sss induces an\sss isomorphism\sss between\sss
$E^{\dff +}_{\trf Y}$ and\sss $E^{\dff -}_{\trf Y}$\nsp,\oss
and\dss therefore\sss 
$(\trf E\dff \oplus\trf E\trf)^{\dff +}_{\trf Y}$\dss
is\dss isomorphic\sss to\sss
$E^{\dff +}_{\trf Y}\trf \oplus\dff E^{\dff -}_{\trf Y}
\off =\off
E\trf |\trf Y$\dnsp.\oss
Since\sss $E\trf |\trf Y$ obviously
extends\sss to $X$\nnsp,\oss Theorem\qss \ref{index-theorem-ext}\qss
applies\sss to\sss $P\dff \oplus\dff P\halfff,\off N\dff \oplus\dff N$\nnsp.\oss
The\sss corollary\sss follows.\oss  \eproof

\myuppar{Remark.}
This corollary\sss also follows from\sss the next\sss theorem.\oss

\mypar{Theorem.}{index-theorem-full}
\emph{\dnsp$\ai\trf(\trf P\halfff,\qff N\trf)
\off =\off
\ti\trf(\trf \sigma\fff,\qff N\trf)$\nnsp.\oss}

\proof
As in\sss the proof\dss of\trs Theorem\qss \ref{index-theorem-ext}\qss
we may assume\sss that\sss the operators\sss
$\mathbf{T}\trf(\trf x_{\dff n}\trf)$\sss from\dss Section\qss \ref{pdo}\qss
do not\sss depend on $x_{\dff n}$\sss for sufficiently small\sss $x_{\dff n}$\dss
(and all\sss values of\dss the parameter $z$\nnsp).\oss
Then\sss the approximation\sss procedures of\qss Section\qss \ref{mult-pdo}\qss apply.\oss
Also,\oss after deforming $P$\sss if\dss necessary,\oss 
we may assume\sss that\sss the pair\sss $\sigma\fff,\qff N$\sss is\dss normalized.\oss

Let\sss $V$\sss be a closed\dss manifold,\pss 
$E\fff'$\sss be a vector bundle over $V$\dnsp,\oss
and\sss $Q$ be a pseudo-differential\sss operator of\dss order $1$ in\sss $E\fff'$\dnsp.\oss
We need only\sss the case when\sss
$V\fff,\pff E\fff'$ and\sss $Q$ do not\sss depend on $z$\nnsp.\oss
Let\sss us denote by $q$\sss the symbol\sss of\dss $Q$\sss
and\sss assume\sss that\sss $q$\sss is\dss an\sss isometric isomorphism.\oss 
This,\oss together with\sss the normalization assumption,\oss
will\sss allow us\sss to apply\trs Theorem\qss \ref{t-index-product}.\oss 
Let\sss $Q^{\fff *}$\sss be an operator\sss formally adjoint\sss to $Q$\nnsp.\oss
Then\sss $Q^{\fff *}$\sss and $Q$ are adjoint\sss to each other
as closed\sss unbounded operators in\sss the\dss Sobolev\dss space\sss 
$H_{\trf 0}\dff(\trf X^{\fff \circ}\fff,\qff E\fff'\trf)$\nnsp.\oss
Suppose further\sss that\sss 
$\kernel\fff Q^{\fff *}\off =\off 0$\sss
and\dss $\dim\dff \kernel\fff Q\off =\off 1$\nnsp.\oss
This will\sss allow us\sss to apply\trs Theorem\qss \ref{mult-of-index}.\oss
Clearly,\oss the analytical\dss index of\dss $Q$\sss is\dss equal\dss to $1$\nnsp.\oss
Hence\sss its\sss topological\dss \dss index\sss 
$\ti\trf(\trf q\trf)$\sss is\dss also equal\dss to $1$\nnsp.\oss

The product\sss $P\pff \omult_{\fff 1}\qff Q$\sss is\dss 
an operator\sss acting\sss in\sss the bundle\sss
$(\trf E\dff \boxtimes\dff E\fff' \trf)
\qff \oplus\qff
(\trf E\dff \boxtimes\dff E\fff' \trf)$ 
and\dss belonging\sss to\sss the class $\overline{\mathcal{P}^{\trf 1}}$\dnsp.\oss
The symbol\sss of\dss $P\pff \omult_{\fff 1}\qff Q$\sss is\dss equal\dss to
$\sigma\pff \omult_{\fff 1}\qff q$\nnsp.\oss
Of\dss course,\oss the operator $P$\sss implicitly depends on\sss $z$\nnsp.\oss 
We will\sss denote\sss it\sss by $P\dff(\trf z\trf)$\sss when\sss there\dss is\dss a need\dss 
to stress\sss the dependence of\dss $P$ on\sss $z\qff \in\qff Z$\nnsp,\oss
and use similar notations for other objects.\oss
Let\sss $A\off =\off P_{\dff \Gamma}$\nsp.\oss
By\trs Theorem\qss \ref{properties-product}\qss the product\sss
$A\pff \omult_{\fff 1}\qff Q$\sss is\trs Fredholm\dss and self-adjoint.\oss
Moreover,\oss by\trs Lemma\qss \ref{product-fredholm}\qss not\sss only\sss
the individual\sss operators\sss 
$A\dff(\trf z\trf)\pff \omult_{\fff 1}\qff Q$\sss
are\dss Fredholm,\oss but\sss the family
$A\dff(\trf z\trf)\pff \omult_{\fff 1}\qff Q\dff,\off z\qff \in\qff Z$\dss
is\dss also\dss Fredholm.\oss
By\trs Theorem\qss \ref{mult-of-index}\qss the analytical\dss index of\dss this family\dss
is\dss equal\dss to\sss the analytical\dss index of\dss the family\sss
$A\dff(\trf z\trf)\dff,\off z\qff \in\qff Z$\nnsp.\oss

As we saw\sss in\dss Section\qss \ref{mult-pdo},\oss
the domains of\dss the unbounded operators $A\dff(\trf z\trf)\pff \omult_{\fff 1}\qff Q$\sss 
are equal\dss to\vspace{3pt}
\[
\quad
\mathcal{D}\dff(\trf z\trf)
\off =\off\dff
\kernel\fff 
\left(\qff \Gamma\dff(\trf z\trf)\qff \boxtimes\qff E\fff'\qff\right)
\dff \oplus\dff
\left(\qff \Gamma\dff(\trf z\trf)\qff \boxtimes\qff E\fff'\qff\right)
\pff. 
\]

\vspace{-12pt}\vspace{3pt}
This means\sss that\sss the operator\sss $A\dff(\trf z\trf)\pff \omult_{\fff 1}\qff Q$\sss
is\dss the unbounded self-adjoint\sss operator associated with
$P\dff(\trf z\trf)\pff \omult_{\fff 1}\qff Q$\sss
and\sss the boundary condition\vspace{3pt}
\begin{equation}
\label{bound-product}
\quad
\mathcal{N}\dff(\trf z\trf)
\off =\off\dff
\left(\qff N\dff(\trf z\trf)\qff \boxtimes\qff E\fff'\qff\right)
\dff \oplus\dff
\left(\qff N\dff(\trf z\trf)\qff \boxtimes\qff E\fff'\qff\right)
\pff. 
\end{equation}

\vspace{-12pt}\vspace{3pt}
The operators $P\dff(\trf z\trf)\pff \omult_{\fff 1}\qff Q$\sss do not\dss
belong\sss to\sss the\dss H\"{o}rmander\dss class.\oss
Nevertheless,\oss as we saw\sss in\dss Section\qss \ref{mult-pdo},\oss the family\sss 
$A\dff(\trf z\trf)\pff \omult_{\fff 1}\qff Q\dff,\off z\qff \in\qff Z$\trs
is\qss Fredholm\dss homotopic\sss to\sss the family of\dss 
restrictions\sss to\sss $\mathcal{D}\trf(\trf z\trf)$\sss of\dss some operators\sss
$[\qff P\trf(\trf z\trf)
\pff \omult_{\fff 1}\qff 
Q \trf]_{\trf \varepsilon\trf(\trf z\trf)}\dff,\off 
z\qff \in\qff Z$\sss
belonging\sss to\sss the\dss H\"{o}rmander\dss class
and such\dss that\sss their symbols approximate\sss
the symbols of\dss operators\sss $P\trf(\trf z\trf)\pff \omult_{\fff 1}\qff Q$\nnsp.\oss 
If\dss the approximation\dss is\dss good enough,\oss
then\sss the boundary condition\qss (\ref{bound-product})\qss is\dss an elliptic
boundary condition\sss not\sss only\sss for\sss
$P\trf(\trf z\trf)\pff \omult_{\fff 1}\qff Q$\nnsp,\oss 
but\sss also for\sss
$[\qff P\trf(\trf z\trf)
\pff \omult_{\fff 1}\qff 
Q \trf]_{\trf \varepsilon\trf(\trf z\trf)}$\nsp.\oss
The\dss Fredholm\dss homotopy shows\sss that\sss
the analytical\dss index of\dss the family\sss 
$[\qff P\trf(\trf z\trf)
\pff \omult_{\fff 1}\qff 
Q \trf]_{\trf \varepsilon\trf(\trf z\trf)}\dff,\off 
z\qff \in\qff Z$\sss
with\sss the boundary conditions\qss (\ref{bound-product})\qss is\dss
equal\dss to\sss the the analytical\dss index of\dss the family\sss of\dss operators
$A\dff(\trf z\trf)\pff \omult_{\fff 1}\qff Q\dff,\off z\qff \in\qff Z$\dss
and\dss hence\sss to\sss the the analytical\dss index of\dss the family\sss
$A\dff(\trf z\trf)\fff,\off z\qff \in\qff Z$\nnsp.\oss
The\sss latter\dss is\dss nothing\sss else but\sss the analytical\dss index\sss
$\ai\trf(\trf P\halfff,\qff N\trf)$\sss of\dss $(\trf P\fff,\qff N\trf)$\nnsp.\oss

Since\sss the symbols of\dss operators\sss 
$[\qff P\trf(\trf z\trf)
\pff \omult_{\fff 1}\qff 
Q \trf]_{\trf \varepsilon\trf(\trf z\trf)}$\sss
approximate\sss the symbols of\dss operators\sss 
$P\trf(\trf z\trf)\pff \omult_{\fff 1}\qff Q$\nnsp,\oss 
the\sss topological\dss index of\dss the family\sss
$[\qff P\trf(\trf z\trf)
\pff \omult_{\fff 1}\qff 
Q \trf]_{\trf \varepsilon\trf(\trf z\trf)}\dff,\off 
z\qff \in\qff Z$\sss
with\sss the boundary conditions\qss (\ref{bound-product})\qss is\dss
is\dss equal\dss to\sss the\sss topological\dss index of\dss the family\sss
$P\trf(\trf z\trf)\pff \omult_{\fff 1}\qff Q\dff,\off z\qff \in\qff Z$\sss
with\sss the same boundary conditions,\oss
i.e.\qss to\sss $\ti\trf(\trf \sigma\qff \omult_{\fff 1}\qff q\halfff,\qff \mathcal{N}\trf)$\nnsp.\oss
Since\sss $\ti\trf(\trf q\trf)\off =\off 1$\nnsp,\oss
Theorem\qss \ref{t-index-product}\qss implies\sss that\sss
$\ti\trf(\trf \sigma\pff \omult_{\fff 1}\qff q\halfff,\qff \mathcal{N}\trf)
\off =\off 
\ti\trf(\trf \sigma\fff,\qff N\trf)$\nnsp.\oss

Therefore,\oss in order\sss to prove\sss that\sss
$\ai\trf(\trf P\halfff,\qff N\trf)
\off =\off
\ti\trf(\trf P\halfff,\qff N\trf)$\sss
it\dss is\dss sufficient\sss to prove\sss that\sss
the analytical\dss index of\dss the family\sss 
$[\qff P\trf(\trf z\trf)
\pff \omult_{\fff 1}\qff 
Q \trf]_{\trf \varepsilon\trf(\trf z\trf)}\dff,\off 
z\qff \in\qff Z$\sss
with\sss the boundary conditions\qss (\ref{bound-product})\qss
is\dss equal\dss to\sss its\sss topological\dss index.\oss
This will\dss follow\sss from\dss Theorem\qss \ref{index-theorem-ext}\qss
once we verify\sss its assumption.\oss
Recall\dss that\sss
$E\trf|\trf Y\off =\off E^{\dff +}_{\trf Y}\qff \oplus\qff E^{\dff -}_{\trf Y}$
and\dss the bundle $N$\sss is\dss isomorphic\sss to both\sss
$E^{\dff +}_{\trf Y}$ and\sss $E^{\dff -}_{\trf Y}$\nnsp.\oss
It\sss follows\sss that\sss $N\dff \oplus\dff N$\sss
is\dss isomorphic\sss to $E\trf|\qff Y$\sss and\dss hence extends\sss to $X$\nnsp.\oss
Also,\vspace{3pt}
\[
\quad
\left(\qff \bigl(\trf E\dff \boxtimes\dff E\fff' \trf\bigr)
\qff \oplus\qff
\left(\trf E\dff \boxtimes\dff E\fff' \trf\right)\qff\right)^{\dff +}_{\qff Y}
\qff
\]

\vspace{-12pt}\vspace{3pt}
is\dss isomorphic\sss to\sss $\mathcal{N}$\sss and\dss hence\sss to\sss
$(\trf N\dff \oplus\dff N\trf)\qff \boxtimes\qff E\fff'$\dnsp.\oss
Since $N\dff \oplus\dff N$ extends\sss to $X$\nnsp,\oss
this implies\sss that\sss $\mathcal{N}$\sss extends\sss to\sss $X\dff \times\dff V$\dnsp.\oss
Therefore\dss Theorem\qss \ref{index-theorem-ext}\qss indeed applies\sss to
$[\qff P\trf(\trf z\trf)
\pff \omult_{\fff 1}\qff 
Q \trf]_{\trf \varepsilon\trf(\trf z\trf)}\dff,\off 
z\qff \in\qff Z$\sss
with\sss the boundary conditions\sss 
$\mathcal{N}\dff(\trf z\trf)\dff,\off 
z\qff \in\qff Z$\nnsp.\oss
The\sss theorem\sss follows.\oss  \eproof

\newpage
\mysection{Self-adjoint\qss Fredholm\qss relations}{relations}

\myuppar{Linear\sss relations.}
Let\sss $H$\sss be a separable\dss Hilbert\dss space,\oss
perhaps finitely dimensional.\oss
A\qss \emph{linear\dss relation}\qss on\sss $H$\sss is\dss simply a linear subspace\sss
$\mathcal{B}\qff \subset\qff H\dff \oplus\dff H$\nnsp.\oss
We will\sss consider mostly\qss \emph{closed}\pss relations,\oss
i.e.\qss relations $\mathcal{B}$ closed as a subset\sss of\dss $H\dff \oplus\dff H$\nnsp.\oss
Linear\sss relations should\sss be\sss thought\sss as generalizations of\trs
linear\sss operators $B\dff \colon\dff H\qff \ttoo\qff H$\nnsp,\oss
with\sss graph\vspace{3pt}
\[
\quad
\mathcal{G}\dff(\trf B\trf)
\off =\off
\bigl\{\qff (\trf x\fff,\qff B\trf(\trf x\trf)\trf)
\pff \bigl|\qff\fff
x\qff \in\qff H \qff\bigr\}
\off \subset\off
H\dff \oplus\dff H
\]

\vspace{-12pt}\vspace{3pt}
corresponding\sss to $B$\nnsp.\oss
We will\sss usually\sss identify operators $B$ with\sss their graphs $\mathcal{G}\dff(\trf B\trf)$\nnsp.\oss
The usual\sss operations naturally\sss extend\sss from operators\sss to relations.\oss
For\sss relations
$\mathcal{B}\fff,\pff \mathcal{C}$\dss
let\sss\vspace{3pt}
\[
\quad
\kernel\fff \mathcal{B}
\off =\off
\{\pff a\pff |\qff (\trf a\fff,\qff 0\trf)\qff \in\qff \mathcal{B} \pff\}
\qff,
\quad
\]

\vspace{-33pt}
\[
\quad
\mathcal{B}^{\dff -\dff 1}
\off =\off
\{\pff (\trf a\fff,\qff b\trf)\pff |\qff (\trf b\fff,\qff a\trf)\qff \in\qff \mathcal{B} \pff\}
\qff,\quad
\mbox{and}\quad
\]

\vspace{-33pt}
\[
\quad 
\mathcal{B}\qff +\qff \mathcal{C}
\off =\off
\{\pff (\trf a\fff,\qff b\qff +\qff c\trf)\pff |\qff (\trf a\fff,\qff b\trf)\qff \in\qff \mathcal{B}
\off\qff
\mbox{and}\off\qff
(\trf a\fff,\qff c\trf)\qff \in\qff \mathcal{C}\pff\}
\pff.
\]

\vspace{-12pt}\vspace{3pt}
The difference\sss $\mathcal{B}\qff -\qff \mathcal{C}$\sss is\dss defined\sss in\sss the same way.\oss
Note\sss that,\oss in\sss contrast\sss with operators,\oss the inverse\sss
$\mathcal{B}^{\dff -\dff 1}$\sss of\dss a relation\sss $\mathcal{B}$\sss is\dss
always well\sss defined as a relation.\oss

For a closed\dss linear\sss relation\sss 
$\mathcal{B}\qff \subset\qff H\dff \oplus\dff H$\sss 
let\sss
$\mathcal{B}_{\dff \infty}\off =\off \mathcal{B}\qff \cap\qff (\trf 0\dff \oplus\dff H \trf)$
and\sss $\mathcal{B}_{\dff s}$ be\sss the orthogonal\sss complement\sss of\dss
$\mathcal{B}_{\dff \infty}$\sss in\sss $\mathcal{B}$\dss
(with respect\sss to\sss the usual\sss scalar\sss product\fff),\oss i.e.\dss
$\mathcal{B}\off =\off \mathcal{B}\qff \ominus\qff \mathcal{B}_{\dff \infty}$\nsp.\oss
Then\sss 
$\mathcal{B}
\off =\off
\mathcal{B}_{\dff s}\qff \oplus\qff \mathcal{B}_{\dff \infty}$\sss 
and\sss $\mathcal{B}_{\dff s}$\sss
is\dss equal\dss to\sss the graph $\mathcal{G}\dff(\trf B\trf)$ of\dss
a closed\dss but\sss not\sss necessarily\sss densely defined\sss operator 
$B\off =\off O\trf(\trf \mathcal{B}\trf)$ called\dss the\qss
\emph{operator\dss part}\qss of\dss $\mathcal{B}$\nnsp.\oss
The space\sss $\mult\fff(\trf \mathcal{B}\trf)$\sss such\sss that\sss
$\mathcal{B}_{\dff \infty}\off =\off 0\trf \oplus\dff \mult\fff(\trf \mathcal{B}\trf)$
is\dss called\dss the\qss \emph{multivalued\dss part}\oss  of\dss $\mathcal{B}$\nnsp.\oss

For a complex number\sss $\lambda\qff \in\qff \ccc$\sss we will\sss denote by $\lambda$\sss
also\sss the operator\sss $\lambda\dff \id_{\dff H}$ and\sss the corresponding\sss relation,\oss
when\sss there\dss is\dss no danger of\dss confusion.\oss
The eigenspace of\dss $\mathcal{B}$\sss with\sss the eigenvalue\sss $\lambda$\sss
is\dss defined as $\kernel\fff (\trf \mathcal{B}\qff -\qff \lambda\trf)$\nnsp.\oss
If\dss $B\off =\off O\trf(\trf \mathcal{B}\trf)$\sss 
is\dss the operator\sss part\sss of\dss $\mathcal{B}$\dnsp,\oss
then\sss the eigenspace 
$\kernel\fff (\trf \mathcal{B}\qff -\qff \lambda\trf)$\sss
is\dss equal\dss to\sss the usual\sss eigenspace 
$\kernel\fff (\trf B\qff -\qff \lambda\trf)$\sss
of\dss $B$\nnsp.\oss
The resolvent\sss set\sss $\rho\trf(\trf \mathcal{B} \trf)$ of\dss $\mathcal{B}$\sss is\dss
defined as\sss the set\sss of\dss numbers $\lambda\qff \in\qff \ccc$\sss
such\sss that\sss $(\trf \mathcal{B}\qff -\qff \lambda\trf)^{\dff -\qff 1}$\sss
is\dss the graph of\dss a bounded operator\sss $K\qff \ttoo\qff K$\nnsp.\oss
The resolvent\sss set\sss $\rho\trf(\trf \mathcal{B} \trf)$\sss is\dss open.\oss
See\qss \cite{bhs},\oss Theorem\qss 1.2.6.\oss
If\dss $\mathcal{B}\off =\off \mathcal{G}\dff(\trf B\trf)$\sss is\dss the graph of\dss $B$\nnsp,\oss
then\sss $\rho\trf(\trf \mathcal{B} \trf)$\sss is\dss 
the usual\dss resolvent\sss set\sss $\rho\trf(\trf B \trf)$\sss of\dss $B$\nnsp.

\myuppar{Self-adjoint\dss relations.}
Let\sss $[\trf \bullet\fff,\qff \bullet\trf]$\sss
be\sss the scalar product\sss on\sss
$H\dff \oplus\dff H$\sss defined\dss by\sss the formula\vspace{3pt}
\[
\quad
[\trf (\trf u\fff,\qff v\trf)\fff,\qff (\trf a\fff,\qff b\trf)\trf]
\off =\off
i\trf \sco{\dff u\fff,\qff b\dff}
\qff -\qff
i\trf \sco{\dff v\fff,\qff a\dff}
\qff.
\]

\vspace{-12pt}\vspace{3pt}
The scalar\sss product\sss $[\trf \bullet\fff,\qff \bullet\trf]$\sss
is\dss Hermitian,\oss i.e.\qss
$[\trf x\fff,\qff y\trf]
\off =\off
\overline{[\trf y\fff,\qff x\trf]}$\nnsp,\oss
but\dss is\dss not\sss positive definite.\oss
It\sss turns\sss $H\dff \oplus\dff H$\sss 
into a\qss \emph{Krein\dss space}.\oss
The relation\qss \emph{adjoint}\pss to a relation 
$\mathcal{B}\qff \subset\qff H\dff \oplus\dff H$\sss
is\dss defined as\sss the orthogonal\sss complement\sss $\mathcal{B}^{\fff *}$ of\dss $\mathcal{B}$
with\sss respect\sss to $[\trf \bullet\fff,\qff \bullet\trf]$\nnsp.\oss
Equivalently,\vspace{3pt}\vspace{-1.25pt}
\[
\quad
\mathcal{B}^{\dff *}
\off =\off
\left\{\pff (\trf a\fff,\qff b\trf)\qff \in\qff H\dff \oplus\dff H
\qff \left|\pff
\sco{\dff u\fff,\qff b\dff}
\off =\off
\sco{\dff v\fff,\qff a\dff}
\off\qff
\mbox{for every}\off\qff
(\trf u\fff,\qff v\trf)\qff \in\qff \mathcal{B} \right.
\qff\right\}
\off.
\]

\vspace{-12pt}\vspace{3pt}\vspace{-1.25pt}
A relation $\mathcal{B}$\sss is\dss said\dss to be\qss \emph{self-adjoint}\pss if\dss
$\mathcal{B}^{\fff *}\off =\off \mathcal{B}$\dnsp.\oss
Equivalently,\pss $\mathcal{B}$\sss is\qss \emph{self-adjoint}\pss if\dss
$\mathcal{B}$\sss is\dss equal\dss to its orthogonal\sss complement\sss
with respect\sss to $[\trf \bullet\fff,\qff \bullet\trf]$\nnsp,\oss
i.e.\qss if\dss $\mathcal{B}$\sss is\qss \emph{lagrangian}\qss
with respect\sss to $[\trf \bullet\fff,\qff \bullet\trf]$\nnsp.\oss
If\dss $\mathcal{B}\off =\off \mathcal{G}\dff(\trf B\trf)$\sss is\dss the graph\sss 
of\dss a\sss linear operator\sss
$B\dff \colon\dff H\qff \ttoo\qff H$\nnsp,\oss then 
$\mathcal{B}^{\fff *}\off =\off \mathcal{G}\dff(\trf B^{\dff *}\trf)$\sss
is\dss the graph of\dss the adjoint\sss operator\sss $B^{\dff *}$\nsp\dnsp.\oss
In\sss particular,\pss
$\mathcal{B}\off =\off \mathcal{G}\dff(\trf B\trf)$\sss is\dss
self-adjoint\dss if\dss and\dss only\trs if\dss $B$\sss is\dss self-adjoint.\oss
Clearly,\oss self-adjoint\sss relations are closed.\oss

In\sss fact,\oss self-adjoint\sss relations provide only a minor\sss
generalization of\dss self-adjoint\sss operators.\oss
Let\sss $H_{\dff B}$\sss
be\sss the closure of\dss the domain $\mathcal{D}\dss(\trf B\trf)$ of\dss $B$\nnsp.\oss
If\dss $\mathcal{B}$\sss is\dss self-adjoint,\oss then\sss $B$\sss
leaves $H_{\dff B}$ invariant\sss and\sss can\sss be considered as
a self-adjoint\sss closed densely defined  
operator\sss in\sss $H_{\dff B}$\nsp.\oss
In\sss this case\sss the multivalued\sss part\sss $\mult\fff(\trf \mathcal{B}\trf)$\sss
is\dss equal\dss to\sss the orthogonal\sss complement\sss 
$H_{\dff B}^{\dff \perp}\off =\off H\dff \ominus\dff H_{\dff B}$\nsp.\oss
See\sss \cite{s},\oss Proposition\qss 14.2.\oss
For example,\oss for\sss the relation $0\dff \oplus\dff H$\sss the operator part\dss 
is\dss the zero operator in\sss the zero subspace $0\qff \subset\qff H$
and\sss the multivalued\sss part\dss is\dss $H$\nnsp.\oss
For\sss $H\dff \oplus\dff 0$ the operator part\dss 
is\dss the zero operator in $H$
and\sss the multivalued\sss part\dss is $0$\nnsp.\oss

\myuppar{Cayley\trs transform.}
Self-adjoint\sss relations
$\mathcal{B}\qff \subset\qff H\dff \oplus\dff H$
can be described\sss in\sss terms of\dss unitary operators
$V\dff \colon\dff 
H\qff \ttoo\qff H$\nnsp.\oss
Namely,\oss $\mathcal{B}$ is\dss self-adjoint\trs 
if\dss and\dss only\trs if\vspace{3pt}\vspace{-1.25pt}
\[
\quad
\mathcal{B}
\off =\off
\bigl\{\pff (\trf a\fff,\qff b\trf)\qff \in\qff H\dff \oplus\dff H
\pff \bigl|\pff
(\trf 1\qff -\qff V\trf)\dff b
\off =\off
i\trf
(\trf 1\qff +\qff V\trf)\dff a 
\pff\bigr\}
\]

\vspace{-12pt}\vspace{3pt}\vspace{-1.25pt}
for some unitary operator\sss $V\dff \colon\dff H\qff \ttoo\qff H$\nnsp.\oss
The operator\sss 
$V\off =\off V\trf(\trf \mathcal{B}\trf)$\sss 
is\dss uniquely determined\dss by $\mathcal{B}$ and\dss is\dss called\dss the\qss
\emph{Cayley\dss transform}\pss of\dss $\mathcal{B}$\dnsp.\oss
The\dss Cayley\dss transform\sss $V$\sss 
is\dss equal\sss to\sss the identity on\sss
$\mult\fff(\trf \mathcal{B}\trf)$\sss and\sss to
$(\trf B\qff -\qff i\trf)\dff (\trf B\qff +\qff i\trf)^{\dff -\dff 1}$
on\sss $H\dff \ominus\dff \mult\fff(\trf \mathcal{B}\trf)$\nnsp,\oss
where\sss
$B\off =\off O\trf(\trf \mathcal{B}\trf)$\nnsp.\oss
In other\sss words,\pss $V$\sss induces a unitary\sss operator\sss
$H_{\dff B}\qff \ttoo\qff H_{\dff B}$\nsp,\oss
and\sss this unitary operator\dss is\dss equal\dss to\sss the\dss
Cayley\trs transform of\dss
$B\dff \colon\dff H_{\dff B}\qff \ttoo\qff H_{\dff B}$\nsp.\oss
Conversely,\oss a unitary operator\sss
$V\dff \colon\dff H\qff \ttoo\qff H$\sss
defines a unique self-adjoint\sss relation\sss 
$\mathcal{B}\off =\off \mathcal{B}\dff(\trf V\trf)$\nnsp,\oss
called\dss the\qss
\emph{Cayley\dss transform}\pss of\dss $V$\dnsp,\oss
such\sss that\sss 
$\mult\fff(\trf \mathcal{B}\trf) 
\off =\off
\kernel\fff (\trf 1\qff -\qff V\trf)$ 
and\sss the operator part\sss 
$B\off =\off O\trf(\trf \mathcal{B}\trf)$ acts in\sss 
$\kernel\fff (\trf 1\qff -\qff V\trf)^{\dff \perp}$\sss as
$i\trf
(\trf 1\qff +\qff V\trf)\dff (\trf 1\qff -\qff V\trf)^{\dff -\dff 1}$\dnsp.\oss
As expected,\oss 
$V\off =\off V\trf(\trf \mathcal{B}\trf)$\sss
if\dss and\dss only\trs if\dss
$\mathcal{B}\off =\off \mathcal{B}\dff(\trf V\trf)$\nnsp.\oss

\myuppar{Fredholm\dss relations and\dss families.}
Recall\dss that\sss a self-adjoint\sss closed densely defined operator $B$\sss is\dss
Fredholm\trs if\dss and\sss only\trs if\dss $0$ does not\sss belong\sss to
the essential\sss spectrum of\dss $B$\nnsp,\oss
i.e.\qss the image of\dss the spectral\sss projection\sss
$P_{\trf [\trf -\qff \varepsilon\fff,\qff \varepsilon\trf]}\qff(\trf B\trf)$\sss
is\dss finitely dimensional\sss for some $\varepsilon\qff >\qff 0$\nnsp.\oss
In\sss terms of\dss the\dss Cayley\dss transform\sss
$V\off =\off V\trf(\trf B\trf)$\sss this means\sss that\sss 
the image of\dss the spectral\sss projection\sss
$P_{\trf D}\dff(\trf V\trf)$\sss
is\dss finitely dimensional\sss for some neighborhood $D\qff \subset\qff \ccc$ of\dss $-\qff 1$\nnsp.\oss
Guided\dss by\sss this equivalence,\oss we will\sss say\sss that\sss
a self-adjoint\sss closed\sss relation $\mathcal{B}$\sss is\qss \emph{Fredholm}\pss
if\dss its\dss Cayley\dss transform\sss
$V\off =\off V\trf(\trf \mathcal{B}\trf)$\sss
has\sss this property.\oss
Clearly,\pss $\mathcal{B}$\sss is\trs Fredholm\trs if\dss and\sss only\trs if\dss
the operator\sss part\sss
$B\off =\off O\trf(\trf \mathcal{B}\trf)$\sss is\dss a self-adjoint\sss
closed\sss densely defined\trs Fredholm\dss operator\sss
$H_{\dff B}\qff \ttoo\qff H_{\dff B}$\nsp.\oss

The definition of\trs Fredholm\dss families of\dss self-adjoint\sss operators
from\qss \cite{i2}\qss can be adapted\dss to families of\dss self-adjoint\sss relations
in at\sss least\sss two equivalent\sss ways.\oss
Let\sss
$\mathcal{B}_{\fff w}\dff,\pff w\qff \in\qff W$\dnsp,\oss
where $W$\sss is\dss a\sss topological\sss space,\oss
be a family of\dss self-adjoint\sss relations\sss
$\mathcal{B}_{\fff w}\qff \subset\qff H\dff \oplus\dff H$\nnsp.\oss
For\sss $w\qff \in\qff W$\sss let\sss
$V_w\off =\off V\trf(\trf \mathcal{B}_{\fff w}\trf)$\sss
and\sss
$B_{\fff w}\off =\off O\trf(\trf \mathcal{B}_{\fff w}\trf)$\nnsp.\oss
One can define\sss the\dss Fredholm\dss property of\dss the family\sss
$\mathcal{B}_{\fff w}\dff,\pff w\qff \in\qff W$\sss
either in\sss terms of\dss the family\sss
$B_{\fff w}\dff,\pff w\qff \in\qff W$
of\dss operator parts,\oss or\sss in\sss terms of\dss the family\sss
$V_w\dff,\pff w\qff \in\qff W$\sss 
of\trs Cayley\dss transforms.\oss

In\sss the first\sss approach
we say\sss that\sss the family\sss
$\mathcal{B}_{\fff w}\dff,\pff w\qff \in\qff W$\sss
is\qss \emph{Fredholm}\pss if\dss the family\sss\vspace{3pt}
\[
\quad
B_{\fff w}\dff \colon\dff
H_{\trf B_{\fff w}}\qff \ttoo\qff H_{\trf B_{\fff w}}\dff,\quad
w\qff \in\qff W
\]

\vspace{-12pt}\vspace{3pt}
is\trs Fredholm\dss in\sss the sense of\pss \cite{i2}.\oss
Of\dss course,\oss it\sss was implicitly assumed\sss in\qss \cite{i2}\qss 
that\sss operators in a\sss family are self-adjoint\sss in\sss the same\dss
Hilbert\dss space or\sss in\sss the fibers of\dss a\dss Hilbert\dss bundle.\oss
But\sss the definition,\oss its main properties,\oss
and even\sss the definition\qss (construction)\qss 
of\dss the analytical\dss index do not\sss depend on\sss this\sss
implicit\sss assumption.\oss

In\sss the second\sss approach we
we say\sss that\sss the family\sss
$\mathcal{B}_{\fff w}\dff,\pff w\qff \in\qff W$\sss
is\qss \emph{Fredholm}\pss if\dss the family\sss
$V_w\dff,\pff w\qff \in\qff W$\sss has\sss the following\sss property.\oss
For every\sss $w\qff \in\qff W$\sss there exists a neighborhood\sss $U_{\fff w}$\sss
of\dss $w$ in\sss $W$\sss and a neighborhood\sss $D_{\fff w}$\sss of\dss $-\qff 1$\sss in $\ccc$
such\sss that\sss for\sss $y\qff \in\qff U_{\fff w}$\sss the subspace\vspace{3pt}
\[
\quad
I_{\dff y}
\off =\off
\image\fff P_{\trf D_{\fff w}}\dff\left(\qff V_y\qff\right)
\]

\vspace{-12pt}\vspace{3pt}
is\dss finitely dimensional,\oss continuously depends on\sss
$y\qff \in\qff U_{\fff w}$\nsp,\oss
and\sss the operators\sss
$I_{\dff y}\qff \ttoo\qff I_{\dff y}$\sss induced\sss by\sss operators $V_y$ 
also continuously depend on $y$\nnsp.\oss
This definition\dss is\dss
equivalent\sss to\sss the first\sss one
because\sss $V_w$\sss is\dss equal\sss to\sss
the identity on\sss the orthogonal\sss complement\sss of\dss $H_{\trf B_{\fff w}}$
and\dss hence\sss this complement\sss does not\sss affect\sss the spectral\dss properties
of\dss $V_w$\sss near $-\qff 1$\nnsp.\oss
The definition of\dss the analytical\dss index can be restated\sss 
in\sss terms of\dss $V_w\dff,\pff w\qff \in\qff W$\nnsp.\oss
It\dss is\dss sufficient\sss to notice\sss that\sss every eigenspace\sss 
$\kernel\fff (\trf \mathcal{B}_{\fff w}\qff -\qff \varepsilon\trf)$\sss
is\dss equal\sss to\sss the eigenspace\sss
$\kernel\fff (\trf V_w\qff -\qff \lambda\trf)$\nnsp,\oss
where 
$\lambda
\off =\off 
(\trf \varepsilon\qff -\qff i\trf)\dff 
(\trf \varepsilon\qff +\qff i\trf)^{\dff -\dff 1}$\dnsp,\oss
and\sss the sign of\dss $\varepsilon$\sss is\dss the same as\sss the sign of\sss 
$\image\fff \lambda$\nnsp.\oss

\myuppar{Continuous\sss families.}
The\dss Fredholm\dss property of\dss a family of\dss self-adjoint\sss operators\dss 
is\dss usually deduced\sss from\sss its continuity\sss in some\sss topology.\oss
The\sss topology of\dss the uniform convergence in\sss the resolvent\sss sense
seems\sss to be\sss the most\sss convenient\sss one.\oss
The same\sss tool\dss is\dss available for\sss the relations.\oss
Let\sss us equip\sss the set\sss 
of\dss self-adjoint\sss closed\sss relations\sss in $H$ with\sss
the\sss topology\sss defined\sss by\sss the norm\sss topology on\sss the space
of\dss orthogonal\sss projections in\sss $H\dff \oplus\dff H$\nnsp.\oss
We will\sss say\sss that\sss a family of\trs Fredholm\dss relations
$\mathcal{B}_{\fff w}\dff,\pff w\qff \in\qff W$\sss is\qss \emph{continuous}\pss
if\dss the map\sss $w\off \longmapsto\off \mathcal{B}_{\fff w}$\sss
is\dss continuous with respect\sss to\sss this\sss topology.\oss
If\dss for every $w$ the  
relation $\mathcal{B}_{\fff w}$\sss
is\dss equal\dss to\sss the graph $\mathcal{G}\dff(\trf B_{\fff w}\trf)$
of\dss a self-adjoint\sss operator $B_{\fff w}$\nsp,\oss
then\sss $\mathcal{B}_{\fff w}\dff,\pff w\qff \in\qff W$\sss is\dss continuous\sss
if\dss and\dss only\trs if\dss 
$B_{\fff w}\dff,\pff w\qff \in\qff W$\sss is\dss continuous\sss
in\sss topology of\dss the uniform convergence in\sss the resolvent\sss sense.

The\dss Cayley\dss transform establishes a homeomorphism between\sss the space\sss
of\dss self-adjoint\sss closed\sss relations in $H$
and\sss the space of\dss unitary operators in $H$ with\sss the norm\sss topology.\oss
Hence a family
$\mathcal{B}_{\fff w}\dff,\pff w\qff \in\qff W$\sss is\dss continuous\sss
if\dss and\dss only\trs if\dss the family 
$V\trf(\trf \mathcal{B}_{\fff w}\trf)\dff,\pff w\qff \in\qff W$\sss is\dss continuous.\oss
The continuity\sss properties of\dss spectral\dss projections imply\sss
that\sss continuous families of\trs Fredholm\dss relations are\dss Fredholm.\oss
Hence\sss the analytical\dss index of\dss such families\dss is\dss well\sss defined.\oss

\myuppar{Fredholm-unitary\sss and\sss related operators.}
Let\sss us say\sss that\sss a unitary operator\sss
$V\dff \colon\dff H\qff \ttoo\qff H$\sss is\qss
\emph{Fredholm-unitary}\pss if\dss
$\image\fff P_{\trf D}\dff(\trf V\trf)$\sss
is\dss finitely dimensional\sss for some neighborhood $D\qff \subset\qff \ccc$ of\dss $-\qff 1$\nnsp.\oss
Equivalently,\pss $V$\sss is\qss \emph{Fredholm-unitary}\pss if\dss
the operator\sss induced\sss by $V$ on\sss the orthogonal\sss complement\sss 
in $H$ of\dss the kernel\sss
$\kernel\fff (\trf V\qff -\qff 1\trf)$\sss is\dss equal\dss to\sss the\dss
Cayley\dss transform of\dss some self-adjoint\trs Fredholm\dss operator\qss
(perhaps,\oss unbounded\fff).\oss
Let\sss $U\fred$\sss be\sss the space of\trs Fredholm-unitary\sss 
operators with\sss the norm\sss topology.\oss

The space\sss $U\fred$\sss is\dss not\sss a\sss group under\sss the composition.\oss
For example,\oss the multiplication\sss by $i$\sss is\trs Fredholm-unitary,\oss
but\sss its square\dss is\dss equal\dss to\sss the\sss the multiplication\sss by $-\qff 1$\sss 
and\dss hence\dss is\dss not\trs Fredholm-unitary.\oss
Let\sss us introduce some subspaces of\dss $U\fred$\sss
homotopy equivalent\sss to $U\fred$ and\sss closed under\sss the composition.\oss
Let\sss $U\ffin$ and\sss $U\comp$
be\sss the groups of\dss unitary\sss operators\sss $H\qff \ttoo\qff H$\sss
differing\sss from\sss the identity 
$\id_{\dff H}$\sss by operators of\dss
finite rank and\dss by compact\sss operators respectively.\oss
Also,\oss let\sss 
$U\dff(\trf \infty\trf)
\off =\off
\bigcup_{\qff n}\fff U\dff(\trf n\trf)$\sss be\sss the most\sss
classical\dss infinite unitary\sss group,\oss
defined\sss in\sss terms of\dss some basis of\dss $H$\nnsp.\oss
Let\sss us equip\sss $U\dff(\trf \infty\trf)$ and\sss $U\ffin$\sss
with\sss the direct\dss limit\sss topologies\qss
(see\qss \cite{i1}\qss for\sss the details)\qss
and\sss $U\comp$ with\sss the norm\sss topology.\oss

\mypar{Theorem.}{unitary-spaces}
\emph{The inclusions}\vspace{3pt}
\[
\quad
U\dff(\trf \infty\trf)
\qff \ttoo\qff
U\ffin
\qff \ttoo\qff
U\comp
\qff \ttoo\qff
U\fred
\]

\vspace{-12pt}\vspace{3pt}
\emph{are homotopy equivalences,\oss
and\sss there\dss is\dss a canonical\dss homotopy equivalence between\sss $U\fred$
and\dss the classifying space\dss $\num{\hat{\mathcal{S}}}$\sss
of\trs the category\dss $\hat{\mathcal{S}}$\sss from\pss 
\textup{\cite{i1}},\oss \textup{\cite{i2}}.\oss}

\proof
The proof\trs is\dss based on\sss methods of\pss \cite{i1}.\oss
Let\sss us begin with constructing a homotopy equivalence between
$U\fred$ and\sss $\num{\hat{\mathcal{S}}}$\nnsp.

For every $\varepsilon\qff >\qff 0$\sss 
let\sss us denote  by\sss $S\dff(\trf \varepsilon\trf)$\sss 
and\sss $D\dff(\trf \varepsilon\trf)$\nnsp,\oss
respectively,\oss the circle and\sss the disc in $\ccc$ 
with\sss the center $-\qff 1$ and\sss the radius $\varepsilon$\nnsp.\oss
Let\sss us define an\qss 
\emph{enhanced\trs Fredholm-unitary\dss operator}\pss
as a pair\sss $(\trf V,\pff \varepsilon\trf)$ such\dss that\sss $V$\sss
is\dss a\dss Fredholm-unitary\sss operator,\pss
$S\dff(\trf \varepsilon\trf)$\sss is\dss disjoint\sss from\sss the spectrum of\dss $V$\dnsp,\oss
and\sss $\image\fff P_{\trf D\dff(\trf \varepsilon\trf)}\dff(\trf V\trf)$\sss
is\dss finitely dimensional.\oss
Let\sss us equip\sss the set\sss $\mathcal{U}\fred$ of\dss enhanced\trs Fredholm-unitary\dss operators\dss
by\sss the\sss topology determined\dss by\sss the norm\sss topology on $U\fred$\sss
and\dss the discrete\sss topology on\sss the set\sss $\rrr_{\qff >\dff 0}$\sss of\dss
parameters $\varepsilon$\nnsp.\oss

The usual\sss order on\sss $\rrr_{\qff >\dff 0}$\sss defines a structure of\dss a\sss
topological\sss category on\sss $\mathcal{U}\fred$\dnsp.\oss
This category\sss has\sss $\mathcal{U}\fred$ as\sss the space of\dss objects 
and a single morphism\sss
$(\trf V,\pff \varepsilon\trf)
\qff \ttoo\qff
(\trf V\fff',\pff \varepsilon'\trf)$
when\sss $V\off =\off V\fff'$ and\sss $\varepsilon\qff \leq\qff \varepsilon'$\nnsp.\oss
There\dss is\dss the obvious forgetting\sss functor\sss
$\hat{\varphi}\dff \colon\dff 
\mathcal{U}\fred\qff \ttoo\qff U\fred$\dnsp,\oss
where\sss $U\fred$\sss is\dss considered as a\sss topological\sss category\sss
having only\sss identity\sss morphisms.\oss
It\sss induces a map of\dss classifying spaces\sss
$\num{\hat{\varphi}}\dff \colon\dff 
\num{\mathcal{U}\fred}\qff \ttoo\qff \num{U\fred}
\off =\off
U\fred$\dnsp.\oss
This map\dss is\dss a homotopy equivalence.\oss
The proof\dss is\dss completely similar\sss to\sss the proof\dss of\trs
Theorem\qss 9.1\qss in\qss \cite{i1}.\oss
Let\sss
$\num{\hat{\varphi}}^{\dff -\dff 1}\dff \colon\dff 
U\fred\qff \ttoo\qff \num{\mathcal{U}\fred}$\sss
be a homotopy\sss inverse of\dss $\num{\hat{\varphi}}$\nnsp.\oss

Replacing\sss the signs of\dss eigenvalues by\sss the signs of\dss 
their\sss imaginary\sss parts 
one can define a functor\sss
$\hat{\pi}\dff \colon\dff
\mathcal{U}\fred\qff \ttoo\qff \hat{\mathcal{S}}$\sss
in exactly\sss the same way as\sss the functor\sss
$\hat{\pi}\dff \colon\dff
\hat{\mathcal{E}}\qff \ttoo\qff \hat{\mathcal{S}}$\sss
from\qss \cite{i1}.\oss
The induced\sss map 
$\num{\hat{\pi}}\dff \colon\dff 
\num{\mathcal{U}\fred}\qff \ttoo\qff \num{\hat{\mathcal{S}}}$
is\dss a homotopy equivalence.\oss
This can\sss be proved\sss in\sss the same way as\dss
Theorem\qss 9.7\qss from\qss \cite{i1},\oss
once we know\sss the contractibility of\dss the spaces\vspace{1.5pt}
\[
\quad
U\fred\dff(\trf \varepsilon\trf)
\off =\off
\left\{\qff \left.
V\qff \in\qff U\fred
\pff \right|\qff\fff
\image\fff P_{\trf D\dff(\trf \varepsilon\trf)}\dff(\trf V\trf)
\off =\off
0
\pff\right\}
\pff
\]

\vspace{-12pt}\vspace{1.5pt}
with $\varepsilon\qff >\qff 0$\nnsp.\oss
This\dss is\dss an analogue of\trs Proposition\qss 8.2\qss from\qss \cite{i1}.\oss
In contrast\sss with\sss that\sss proposition,\oss the proof\dss
of\dss the contractibility of\dss $U\fred\dff(\trf \varepsilon\trf)$ 
does not\sss depend on\sss the\dss Kuiper's\trs theorem.\oss
Indeed,\oss an obvious spectral\sss deformation contracts\sss
$U\fred\dff(\trf \varepsilon\trf)$\sss to\sss the identity map\sss $H\qff \ttoo\qff H$\nnsp.\oss
Hence\sss $\num{\hat{\pi}}$\sss is\dss a homotopy equivalence and\dss the composition\sss
$\num{\hat{\pi}}\dff \circ\dff \num{\hat{\varphi}}^{\dff -\dff 1}\dff \colon\dff 
U\fred\qff \ttoo\qff \num{\hat{\mathcal{S}}}$\sss
is\dss also a homotopy equivalences.\oss
This proves\sss the second statement\sss of\dss the\sss theorem.\oss

Let\sss us prove\sss the first\sss statement.\oss
Let\sss
$\mathcal{U}\dff(\trf \infty\trf)\dff,\off\qff
\mathcal{U}\ffin,\off\qff
\mathcal{U}\comp$\sss
be\sss the\sss topological\sss categories of\dss enhanced operators 
defined as\sss $\mathcal{U}\fred$\dnsp,\oss
but\sss with operators belonging only\sss to\sss
$U\dff(\trf \infty\trf)\dff,\off\qff
U\ffin$\dnsp,\pss
$U\comp$\sss
respectively.\oss
For each of\dss these categories\sss there are obvious analogues
of\trs functors $\hat{\varphi}$ and\sss $\hat{\pi}$\nnsp,\oss
and\sss the diagram of\dss induced\sss maps\vspace{1.5pt}
\[
\quad
\begin{tikzcd}[column sep=tri, row sep=rboom]
U\dff(\trf \infty\trf)
\arrow[rr]
&
&
U\ffin
\arrow[rr]
&
&
U\comp
\arrow[rr]
&
&
U\fred
\\
\protect{\num{\mathcal{U}\dff(\trf \infty\trf)}}
\arrow[u]
\arrow[rrrd]
\arrow[rr]
&
&
\protect{\num{\mathcal{U}\ffin}}
\arrow[u]
\arrow[rr]
\arrow[rd]
&
&
\protect{\num{\mathcal{U}\comp}}
\arrow[rr]
\arrow[u]
\arrow[ld]
&
&
\protect{\num{\mathcal{U}\fred}}
\arrow[u]
\arrow[llld]
\\
&
&
&
\protect{\num{\hat{\mathcal{S}}}}
&
&
&
\end{tikzcd}
\]

\vspace{-12pt}\vspace{1.5pt}
is\dss commutative.\oss
The upper vertical\sss arrows are homotopy equivalences by\sss the same reasons as before.\oss
In order\sss to prove\sss that\sss arrows with\sss the\sss target\sss $\num{\hat{\mathcal{S}}}$
are homotopy equivalences,\oss it\dss is\dss sufficient\sss to prove\sss that\sss
the spaces\sss 
$U\dff(\trf \infty\trf)\dff(\trf \varepsilon\trf)\dff,\off\qff
U\ffin\dff(\trf \varepsilon\trf)\dff,\off\qff
U\comp\dff(\trf \varepsilon\trf)$\nnsp,\oss
defined similarly\sss to
$U\fred\dff(\trf \varepsilon\trf)$\nnsp,\oss
are contractible.\oss
But\sss the spectral\sss deformation contracting\sss
$U\fred\dff(\trf \varepsilon\trf)$\sss to\sss $\id_{\dff H}$\sss
leaves each of\dss these spaces invariant\sss and\dss hence contracts\sss them also.\oss
It\sss follows\sss that\sss all\sss vertical\sss and slanted arrows 
in\sss the above diagram are homotopy equivalences.\oss
Together with\sss the commutativity of\dss the diagram\sss this implies\sss
that\sss the horizontal\sss arrows are homotopy equivalences.\oss
This proves\sss that\sss our inclusions are homotopy equivalences.\oss  \eproof

\myuppar{Remark.}
The fact\sss that\sss $U\dff(\trf \infty\trf)\qff \ttoo\qff U\comp$\sss 
is\dss a homotopy equivalence\dss is\dss implicit\sss in\qss \cite{as}.\oss
The methods of\qss \cite{as}\qss also allow\sss to prove\sss that\sss
$U\comp\qff \ttoo\qff U\fred$\sss 
is\dss a homotopy equivalence.\oss
As\sss M.\dss Prokhorova\dss pointed out\sss to\sss the author,\oss
this was done by\dss P.\dss Kirk\sss and\dss M.\dss Lesch\qss \cite{kl}.\oss
See\qss \cite{kl},\oss Lemma\qss 6.1.\oss
The fact\sss that\sss 
$U\dff(\trf \infty\trf)\qff \ttoo\qff U\ffin$\sss 
is\dss a homotopy equivalence\dss is\dss proved\sss in\qss \cite{i1}.\oss 
See\qss \cite{i1},\oss 
Corollary\qss 12.10.\oss
Probably,\oss it\sss was known\sss long\sss before\qss \cite{i1}.\oss

\myuppar{The analytical\dss index of\dss continuous families of\dss self-adjoint\sss relations.}
Let\sss $\mathcal{B}_{\fff w}\dff,\pff w\qff \in\qff W$\sss
be a continuous family of\dss self-adjoint\trs Fredholm\dss relations,\oss
and\dss let\sss $V_w\dff,\pff w\qff \in\qff W$\sss
be\sss the family of\dss their\dss Cayley\dss transforms\sss
$V_w\off =\off V\trf(\trf \mathcal{B}_{\fff w}\trf)$\nnsp.\oss
Then\sss $\mathcal{B}_{\fff w}\dff,\pff w\qff \in\qff W$\sss
is\dss a\dss Fredholm\dss family,\oss
and,\oss as we already\sss mentioned,\oss 
one can define\sss its analytical\dss index\sss following\qss \cite{i2}.\oss
We can also use\sss $V_w\dff,\pff w\qff \in\qff W$\sss
instead of\dss $\mathcal{B}_{\fff w}\dff,\pff w\qff \in\qff W$\sss
in such definition.\oss
Both constructions\sss lead\dss to an\qss \emph{index\dss map}\qss
$W\qff \ttoo\qff \num{\hat{\mathcal{S}}}$\nnsp,\oss
well\sss defined\sss up\sss to homotopy.\oss
Moreover,\oss both constructions\sss lead\dss to\sss the same index\sss map.\oss
The\qss \emph{analytical\dss index}\qss of\dss
$\mathcal{B}_{\fff w}\dff,\pff w\qff \in\qff W$\sss
is\dss defined as\sss the homotopy class of\dss the index\sss map.\oss
Since\sss $\num{\hat{\mathcal{S}}}$\sss is\dss the classifying space for\sss the
$K^{\dff 1}${\dnsp}-theory,\oss the analytical\dss index\sss may\sss be
considered as an element\sss of\dss the abelian\sss group\sss $K^{\dff 1}\dff(\trf W\trf)$\nnsp.\oss

At\sss the same\sss time\sss the family of\trs Cayley\dss transforms\sss
$V_w\dff,\pff w\qff \in\qff W$\sss defines a continuous map\sss
$W\qff \ttoo\qff U\fred$\dnsp.\oss
Together with\sss the homotopy equivalence
$U\fred\qff \ttoo\qff \num{\hat{\mathcal{S}}}$\dss
from\dss Theorem\qss \ref{unitary-spaces}\qss this map defines another map\sss
$W\qff \ttoo\qff \num{\hat{\mathcal{S}}}$\nnsp.\oss
The constructions of\dss the index\sss map and of\dss the homotopy equivalence\sss
$U\fred\qff \ttoo\qff \num{\hat{\mathcal{S}}}$\dss match each other.\oss
In\sss fact,\oss this homotopy equivalence\dss is\dss the index\sss map
of\dss the\sss tautological\dss family parameterized\dss by\sss $U\fred$\dnsp.\oss
Therefore\sss the analytical\dss index of\dss the family\sss
$\mathcal{B}_{\fff w}\dff,\pff w\qff \in\qff W$\sss
can\sss be also defined as\sss the homotopy class of\dss the
corresponding\sss map\sss $W\qff \ttoo\qff U\fred$\dnsp.\oss\vspace{-0.125pt}

The group operation in $K^{\dff 1}\dff(\trf W\trf)$ can be defined\sss in\sss
terms of\dss classifying spaces.\oss
The most\sss natural\sss approach\dss is\dss based on\sss using\sss
the orthogonal\sss direct\sss sum $H\dff \oplus\dff H$ and some\sss isomorphism\sss
$H\dff \oplus\dff H\qff \ttoo\qff H$\nnsp.\oss
Let\sss $\mathbf{U}$\sss be any of\dss the spaces\sss
$U\dff(\trf \infty\trf)\dff,\off\qff
U\ffin,\off\qff
U\comp,\off\qff
U\fred,\off\qff
\num{\hat{\mathcal{S}}}$\nnsp.\oss
Each of\dss them\sss implicitly\sss depends on\sss the\dss Hilbert\dss space $H$\nnsp,\oss
and\sss the isomorphism\sss
$H\dff \oplus\dff H\qff \ttoo\qff H$\sss
leads\sss to maps\sss
$\mathbf{U}\dff \times\dff \mathbf{U}\qff \ttoo\qff \mathbf{U}$\nnsp.\oss
In\sss the case of\dss
$\mathbf{U}\off =\off U\dff(\trf \infty\trf)$\sss
one needs\sss to choose an\sss isomorphism\sss
$H\dff \oplus\dff H\qff \ttoo\qff H$\sss
respecting\sss the bases of\trs Hilbert\dss spaces in a suitable sense.\oss
Up\sss to homotopy\sss the maps\sss
$\mathbf{U}\dff \times\dff \mathbf{U}\qff \ttoo\qff \mathbf{U}$\sss
do not\sss depend on\sss the choice of\dss isomorphism\sss
$H\dff \oplus\dff H\qff \ttoo\qff H$\sss
and define a binary operation on\sss the set\sss of\dss homotopy classes\sss
$W\qff \ttoo\qff \mathbf{U}$\nnsp.\oss
This operation\dss is\dss nothing else but\sss the group operation\sss
in\sss $K^{\dff 1}\dff(\trf W\trf)$\nnsp.\oss\vspace{-0.125pt}

The spaces\sss 
$U\dff(\trf \infty\trf)\dff,\off\qff
U\ffin,\off\qff
U\comp$\sss
are\sss topological\dss groups,\oss
and\dss if\dss $\mathbf{U}$\sss is\dss equal\dss to one of\dss these spaces,\oss
the group operation\sss also\sss leads\sss to a binary operation
on\sss the set\sss of\dss homotopy classes\sss
$W\qff \ttoo\qff \mathbf{U}$\nnsp.\oss
A well\dss known argument,\oss 
based on homotopies of\trs $2\fff \times\fff 2$ matrices,\oss 
shows\sss that\sss this operation\dss is\dss the same
as\sss the one defined\sss in\sss the previous paragraph.\oss 
While\sss $U\fred$\sss is\dss not\sss a group,\oss if\trs
$\mathbf{U}
\off =\off 
U\dff(\trf \infty\trf)\dff,\off\off
U\ffin,\off\off
U\comp$\dnsp,\oss
then\sss the composition of\dss unitary operators defines
an action of\dss $\mathbf{U}$ on\sss $U\fred$\dnsp.\oss
In\sss turn,\oss this action defines an action of\dss maps\sss
$W\qff \ttoo\qff \mathbf{U}$\sss on\sss maps\sss
$W\qff \ttoo\qff U\fred$\dnsp,\oss
which we will\sss denote by\sss
$(\trf f\fff,\qff g\trf)\off \longmapsto\off f\fff \cdot\dff g$\nnsp.\oss\vspace{-0.125pt}

\mypar{Lemma.}{u-action}
\emph{Suppose\sss that\qss
$\mathbf{U}
\off =\off 
U\dff(\trf \infty\trf)\dff,\off\qff
U\ffin$ or\sss
$U\comp$\dnsp.\oss
Let\qss $a\fff,\qff b\qff \in\qff K^{\dff 1}\dff(\trf W\trf)$\sss
be elements represented\dss by\sss the maps\sss
$f\dff \colon\dff W\qff \ttoo\qff U\fred$\sss
and\sss
$g\dff \colon\dff W\qff \ttoo\qff \mathbf{U}$\sss
respectively.\oss
Then\sss the element\sss 
$a\qff +\qff b\qff \in\qff K^{\dff 1}\dff(\trf W\trf)$\sss 
is\dss represented\dss by\sss the map\sss
$f\fff \cdot\dff g$\nnsp.\oss}\vspace{-0.125pt}

\proof
By\trs Theorem\qss \ref{unitary-spaces}\qss the map $f$\sss is\dss homotopic\sss
to a map\sss $f\fff'\dff \colon\dff W\qff \ttoo\qff \mathbf{U}$\nnsp.\oss
The map\sss $f\fff \cdot\dff g$\sss is\dss homotopic\sss to\sss
$f\fff'\fff \cdot\dff g$\sss and\sss hence represents\sss 
the same element\sss of\dss $K^{\dff 1}\dff(\trf W\trf)$\nnsp.\oss
By\sss the remarks preceding\sss the\sss lemma,\oss
this element\dss is\dss equal\dss to\sss $a\qff +\qff b$\nnsp.\oss  \eproof

\newpage
\mysection{Boundary\qss triplets\qss and\qss the\qss analytical\qss index}{boundary-triplets}

\myuppar{Boundary\sss triplets.}
Let\sss us review\dss the definitions and\sss the basic facts of\dss the\sss theory 
of\pss \emph{boundary\sss triplets}\qss following\trs Schm\"{u}gden\qss \cite{s},\oss Chapter\qss 14.\oss
Let\sss $T$\sss be a symmetric densely defined operator in a\dss Hilbert\dss space $H$\nnsp,\oss
and\dss let\sss $\mathcal{D}\off =\off \mathcal{D}\dff(\trf T^{\dff *}\trf)$\sss
be\sss the domain of\dss its adjoint\sss $T^{\dff *}$\dnsp.\oss
A\qss \emph{boundary\sss triplet}\oss for\sss $T^{\dff *}$\sss is\dss a\sss triple\sss
$(\qff K\dff,\off \Gamma_{\fff 0}\dff,\off \Gamma_{\fff 1}\trf)$
consisting of\dss a\dss Hilbert\dss space\sss $K$ and\sss two\sss linear\sss maps\sss
$\Gamma_{\fff 0}\dff,\qff \Gamma_{\fff 1}\dff \colon\dff 
\mathcal{D}\qff \ttoo\qff K$\sss
such\sss that\sss the map\sss
$\Gamma_{\fff 0}\dff \oplus\dff \Gamma_{\fff 1}\dff \colon\dff 
\mathcal{D}\qff \ttoo\qff K\dff \oplus\dff K$\sss
is\dss surjective and\vspace{3pt}
\begin{equation}
\label{triplet}
\quad
\sco{\trf T^{\dff *} x\fff,\qff y\trf}
\qff -\qff
\sco{\trf x\fff,\qff T^{\dff *} y\trf}
\off =\off
\sco{\trf \Gamma_{\fff 1}\dff x\fff,\pff \Gamma_{\fff 0}\dff y\trf}
\qff -\qff
\sco{\trf \Gamma_{\fff 0}\dff x\fff,\pff \Gamma_{\fff 1}\dff y\trf}
\end{equation}

\vspace{-12pt}\vspace{3pt}
for every\sss $x\fff,\qff y\qff \in\qff \mathcal{D}$\dnsp.\oss
The identity\qss (\ref{triplet})\qss is\dss another abstract\sss version 
of\dss the\dss Lagrange\dss identity\sss and\trs Green\dss formulas.\oss
Cf.\qss Section\qss \ref{abstract-index}.\oss
Sometimes\sss it\dss is\dss convenient\sss to rewrite\qss (\ref{triplet})\qss
in\sss terms of\dss the\sss linear\sss maps\sss
$\Gamma_{+}\off =\off \Gamma_{\fff 1}\qff +\qff i\trf \Gamma_{\fff 0}$\sss
and\sss
$\Gamma_{-}\off =\off \Gamma_{\fff 1}\qff -\qff i\trf \Gamma_{\fff 0}$\nsp.\oss
Namely,\oss (\ref{triplet})\qss is\dss equivalent\sss to\vspace{3pt}
\[
\quad
2\dff i\trf \sco{\trf T^{\dff *} x\fff,\qff y\trf}
\qff -\qff
2\dff i\trf \sco{\trf x\fff,\qff T^{\dff *} y\trf}
\off =\off
\sco{\trf \Gamma_{-}\dff x\fff,\pff \Gamma_{-}\dff y\trf}
\qff -\qff
\sco{\trf \Gamma_{+}\dff x\fff,\pff \Gamma_{+}\dff y\trf}
\off.
\]

\vspace{-12pt}\vspace{3pt}
The surjectivity of\dss $\Gamma_{\fff 0}\dff \oplus\dff \Gamma_{\fff 1}$\sss
is\dss equivalent\sss to\sss the surjectivity of\dss
$\Gamma_{-}\dff \oplus\dff \Gamma_{+}$\nsp.\oss
A boundary\sss triplet\sss for $T^{\dff *}$ exists\sss if\dss and\dss only\trs if\dss
the operator\sss 
$T$ admits a self-adjoint\sss extension.\oss

\myuppar{Boundary\sss triplets and self-adjoint\sss extensions.}
Let\sss $T$ be a symmetric densely defined operator in $H$
and\dss let\sss $(\qff K\dff,\off \Gamma_{\fff 0}\dff,\off \Gamma_{\dff 1}\trf)$
be a boundary\sss triplet\sss for $T^{\dff *}$\nsp\dnsp.\oss
Every\sss linear\sss relation
$\mathcal{B}\qff \subset\qff K\dff \oplus\dff K$\sss
leads\sss to an extension\sss $T_{\fff \mathcal{B}}$\sss of\dss
$T$\dnsp,\oss defined as\sss 
the restriction of\dss $T^{\dff *}$\sss to\sss
the domain\vspace{3pt}
\[
\quad
\mathcal{D}\dff\left(\qff T_{\fff \mathcal{B}} \qff\right)
\off =\off
\left\{\qff \left.
x\qff \in\qff \mathcal{D}\dff(\trf T^{\dff *}\trf)
\pff \right|\qff
(\trf \Gamma_{\fff 0}\dff x\fff,\pff \Gamma_{\fff 1}\dff x\trf)
\qff \in\qff
\mathcal{B} 
\pff\right\} 
\off.
\]

\vspace{-12pt}\vspace{3pt}
It\sss turns out\sss that\sss $T_{\fff \mathcal{B}}$\sss is\dss self-adjoint\trs
if\dss and\dss only\trs if\dss $\mathcal{B}$\sss is\dss a self-adjoint\sss relation.\oss
Moreover,\oss every self-adjoint\sss extension of\dss $T$\sss is\dss equal\dss to
$T_{\fff \mathcal{B}}$\sss for a unique self-adjoint\sss relation $\mathcal{B}$\dnsp.\oss

The extension\sss $T_{\fff \mathcal{B}}$\sss 
can be described\sss in\sss terms of\dss  
$B\off =\off O\trf(\trf \mathcal{B}\trf)$ 
and\sss $K_{\dff B}$ as follows.\oss
Let\sss $P_{\dff B}$\sss  
be\sss the orthogonal\dss projection
$K\qff \ttoo\qff K_{\dff B}$\nnsp.\oss
Then\sss $T_{\fff \mathcal{B}}$\sss is\dss equal\dss to\sss 
the restriction $T_{\dff B}$ of\dss $T^{\dff *}$\sss to\vspace{3pt}
\[
\quad
\mathcal{D}\dff\left(\qff T_{\dff B} \qff\right)
\off =\off
\left\{\qff \left.
x\qff \in\qff \mathcal{D}\dff(\trf T^{\dff *}\trf)
\pff \right|\pff
\Gamma_{\fff 0}\dff x\qff \in\qff \mathcal{D}\dff(\trf B\trf)
\off\qff
\mbox{and}\off\qff
B\dff \circ\dff \Gamma_{\fff 0}\dff x
\off =\off
P_{\dff B}\dff \circ\dff \Gamma_{\fff 1}\dff x\trf)
\pff\right\} 
\off.
\]

\vspace{-12pt}\vspace{3pt}
Of\dss course,\oss this means simply\sss that\sss
$\mathcal{D}\dff\left(\qff T_{\fff \mathcal{B}} \qff\right)
\off =\off
\mathcal{D}\dff\left(\qff T_{\dff B} \qff\right)$\sss
and\sss follows from\sss
$\mathcal{B}
\off =\off
\mathcal{B}_{\dff s}\qff \oplus\qff \mathcal{B}_{\dff \infty}$\nsp.\oss
The extension $T_{\fff \mathcal{B}}$\sss can be also described\sss in\sss terms of\dss
the\dss Cayley\dss transform\sss $V\trf(\trf \mathcal{B}\trf)$\nnsp.\oss
For a unitary operator 
$V\dff \colon\dff 
K\qff \ttoo\qff K$\dss
let\sss $T^{\dff V}$\sss be\sss the restriction of\dss $T^{\dff *}$\sss to\vspace{1.5pt}
\[
\quad
\left\{\qff \left.
x\qff \in\qff \mathcal{D}\dff(\trf T^{\dff *}\trf)
\pff \right|\pff
V\dff \circ\dff \Gamma_{+}\dff x\off =\off \Gamma_{-}\dff x
\pff\right\} 
\off =\off
\mathcal{D}\dff\left(\qff T^{\dff V} \qff\right)
\off.
\]

\vspace{-12pt}\vspace{1.5pt}
Then\sss $T^{\dff V}\off =\off T_{\fff \mathcal{B}}$\dss if\dss and\dss only\trs if\dss
$V\off =\off V\trf(\trf \mathcal{B}\trf)$\sss is\dss the\dss Cayley\dss transform of\dss $\mathcal{B}$\dnsp.\oss
There are\sss two special\sss extensions associated\sss with\sss the boundary\sss triplet,\oss
namely,\oss the extensions\sss $T_{\dff 0}$ and\sss $T_{\fff 1}$ corresponding\sss
to\sss the relations $0\dff \oplus\dff K$\sss and\sss $K\dff \oplus\dff 0$ respectively.\oss
Clearly,\oss the domains of\dss $T_{\dff 0}$ and\sss $T_{\fff 1}$ are\sss
$\kernel\fff \Gamma_{\fff 0}$ and\sss $\kernel\fff \Gamma_{\dff 1}$ respectively.\oss\vspace{-0.625pt}

\myuppar{Gamma\sss fields and\dss Weyl\dss functions.}
Let\sss $T$\sss be a symmetric densely defined operator in\sss  $H$\sss 
and\sss
$(\qff K\dff,\off \Gamma_{\fff 0}\dff,\off \Gamma_{\dff 1}\trf)$
be a boundary\sss triplet\sss for $T^{\dff *}$\dnsp.\oss
For\sss $\lambda\qff \in\qff \ccc$\sss let
$\mathcal{N}_{\trf \lambda}
\off =\off 
\kernel\dff (\trf T^{\dff *}\qff -\qff \lambda\trf)$\nnsp.\oss
Let\sss us equip\sss  
$\mathcal{D}\off =\off \mathcal{D}\dff(\trf T^{\dff *}\trf)$\sss
with\sss the\sss graph\sss norm\sss
$\norm{x}_{\trf T^{\dff *}}
\off =\off
\norm{x}\qff +\qff \norm{T^{\dff *}\fff x}$\nnsp,\oss
where\sss $\norm{\bullet}$\sss is\dss the norm in\sss $H$\nnsp.\oss
Clearly,\oss on subspaces\sss 
$\mathcal{N}_{\trf \lambda}$\sss the graph\sss 
norm\sss $\norm{\bullet}_{\trf T^{\dff *}}$
is\dss equivalent\sss to\sss $\norm{\bullet}$\nnsp.\oss
Recall\dss that\sss $T_{\dff 0}$\sss is\dss the extension of\dss $T$\sss with\sss the domain\sss
$\kernel\fff \Gamma_{\fff 0}$\nsp.\oss
It\sss turns out\sss that\sss the maps\sss
$\Gamma_{\fff 0}\dff,\qff \Gamma_{\fff 1}\dff \colon\dff 
\mathcal{D}\qff \ttoo\qff K$\sss
are continuous,\oss and\dss if\sss 
$\lambda\qff \in\qff \rho\trf(\trf T_{\dff 0}\trf)$\nnsp,\oss
then\sss $\Gamma_{\fff 0}$\sss induces a\sss topological\sss isomorphisms\sss
$\mathcal{N}_{\trf \lambda}\qff \ttoo\qff K$\nnsp.\oss
See\qss \cite{s},\oss Lemma\qss 14.13.\oss

For every\sss $z\qff \in\qff \rho\trf(\trf T_{\dff 0}\trf)$\sss let\sss
$\bm{\gamma}\qff(\trf z\trf)
\dff \colon\dff
K\qff \ttoo\qff \mathcal{N}_{\trf z}$\sss 
be\sss the inverse of\dss the map\sss
$\mathcal{N}_{\trf z}\qff \ttoo\qff K$\sss
induced\sss by\sss $\Gamma_{\fff 0}$\nsp,\oss
and\dss let\sss
$M\trf(\trf z\trf)
\off =\off 
\Gamma_{\dff 1}\dff \circ\qff \bm{\gamma}\trf(\trf z\trf)
\dff \colon\dff
K\qff \ttoo\qff K$\nnsp.\oss
Clearly,\pss $\bm{\gamma}\qff(\trf z\trf)$ and\sss $M\trf(\trf z\trf)$
are bounded operators.\oss
The maps\sss
$z\off \longmapsto\off \bm{\gamma}\qff(\trf z\trf)$
and\sss
$z\off \longmapsto\off M\trf(\trf z\trf)$
defined on\sss $\rho\trf(\trf T_{\dff 0}\trf)\qff \subset\qff \ccc$\sss
are called,\oss respectively,\oss the\qss \emph{gamma\dss field}\qss and\sss the\qss
\emph{Weyl\qss function}\pss of\dss the operator $T_{\dff 0}$ 
associated with\sss the boundary\sss triplet\sss
$(\qff K\dff,\off \Gamma_{\fff 0}\dff,\off \Gamma_{\dff 1}\trf)$\nnsp.\oss
They\sss are holomorphic functions on\sss $\rho\trf(\trf T_{\dff 0}\trf)$\nnsp.\oss
Let\sss us state some basic properties of\dss
the gamma\sss field\sss and\sss the\dss Weyl\dss function.\oss
First,\vspace{0.5pt}
\[
\quad
M\trf(\trf z\trf)^{\fff *}
\off =\off\dff
M\trf(\pff \overline{z}\pff)
\quad
\]

\vspace{-12pt}\vspace{0.5pt}
for every\sss $z\qff \in\qff \rho\trf(\trf T_{\dff 0}\trf)$\nnsp.\oss
Next,\oss suppose\sss that\sss $\mathcal{B}$\sss is\dss a closed\sss relation\sss in\sss $K$\nnsp.\oss
Then\sss
$z\qff \in\qff \rho\trf(\trf T_{\fff \mathcal{B}}\trf)$\sss
if\trs and\dss only\trs if\trs
$0\qff \in\qff \rho\trf(\trf \mathcal{B}\qff -\qff M\trf(\trf z\trf)\trf)$\nnsp.\oss
If\dss $z\qff \in\qff \rho\trf(\trf T_{\dff 0}\trf)$\nnsp,\oss
then\vspace{2.25pt}
\[
\quad
\mathcal{B}\qff -\qff M\trf(\trf z\trf)
\off =\off
\bigl\{\pff \left. 
\left(\trf a\fff,\pff b\qff -\qff M\trf(\trf z\trf)\dff(\trf a\trf) \trf\bigr)
\qff \in\qff K\dff \oplus\dff K
\pff \right|\pff
(\trf a\fff,\qff b\trf)\qff \in\qff \mathcal{B} \right.
\qff\bigr\}
\off
\]

\vspace{-12pt}\vspace{2.25pt}
is\dss also a closed\sss relation,\oss
and\sss $\bm{\gamma}\qff(\trf z\trf)$\sss induces an\sss isomorphism\vspace{2.25pt}
\begin{equation*}
\quad
\kernel\dff
(\qff
\mathcal{B}\qff -\qff M\trf(\trf z\trf)
\qff)
\qff \ttoo\qff
\kernel\dff
(\qff
T_{\fff \mathcal{B}}\qff -\qff z
\qff)
\pff.
\end{equation*}

\vspace{-12pt}\vspace{2.25pt}
See\qss \cite{s},\oss Propositions\qss 14.15\qss and\qss 14.17.\oss
If\dss 
$z\fff,\qff w\qff \in\qff \rho\trf(\trf T_{\dff 0}\trf)$\nnsp,\oss
then\vspace{2.25pt}
\begin{equation}
\label{gamma-field}
\quad
\bm{\gamma}\qff(\trf w\trf)
\off =\off
(\trf T_{\dff 0}\qff -\qff z\trf)\qff
(\trf T_{\dff 0}\qff -\qff w\trf)^{\fff -\dff 1}\qff
\bm{\gamma}\qff(\trf z\trf)
\quad
\mbox{and}\quad
\end{equation}

\vspace{-34.5pt}
\begin{equation}
\label{weyl-function}
\quad
M\trf(\trf w\trf)\qff -\qff M\trf(\trf z\trf) 
\off =\off
(\trf w\qff -\qff z\trf)\qff
\bm{\gamma}\qff(\qff \overline{z}\qff)^{\fff *}\qff
\bm{\gamma}\qff(\trf w\trf)
\pff.
\end{equation}

\vspace{-12pt}\vspace{2.25pt}
See\qss \cite{s},\oss Propositions\qss 14.14\dff({\fff}iv{\fff})\qss
and\qss 14.15\dff({\fff}iii{\fff}).\oss

\myuppar{The\dss Krein--Naimark\dss resolvent\dss formula.}
This\dss is\dss the main\sss result\sss of\dss the\sss theory.\oss
It\sss compares\sss the resolvents of\dss $T_{\fff \mathcal{B}}$ and\sss $T_{\dff 0}$\nsp.\oss 
Namely,\oss if\dss
$z\qff \in\qff \rho\trf(\trf T_{\dff 0}\trf)\qff \cap\qff \rho\trf(\trf T_{\fff \mathcal{B}}\trf)$\nnsp,\oss
then\vspace{2.25pt}
\[
\quad
(\qff
T_{\fff \mathcal{B}}\qff -\qff z
\qff)^{\fff -\dff 1}
\pff -\pff
(\qff
T_{\dff 0}\qff -\qff z
\qff)^{\fff -\dff 1}
\off =\off\dff
\bm{\gamma}\qff(\trf z\trf)\qff
(\qff
\mathcal{B}\qff -\qff M\trf(\trf z\trf)
\qff)^{\fff -\dff 1}\qff
\bm{\gamma}\qff(\qff \overline{z}\qff)^{\fff *}
\pff.
\]

\vspace{-12pt}\vspace{2.25pt}
See\qss \cite{s},\oss the formula\qss (14.43)\qss in\dss Theorem\qss 14.18.\oss

\myuppar{A standard example.}
Given a self-adjoint\sss extension $A$ of\dss $T$\sss one can construct\sss
boundary\sss triplets\sss 
$(\qff K\dff,\off \Gamma_{\fff 0}\dff,\off \Gamma_{\dff 1}\trf)$
such\sss that\sss $A\off =\off T_{\dff 0}$\nsp.\oss
Let\sss us\sss fix\sss $\mu\qff \in\qff \rho\trf(\trf A\trf)$\nnsp.\oss
For every\sss $x\qff \in\qff \mathcal{D}\dff(\trf T^{\dff *}\trf)$\sss
there exists unique vectors\sss
$x_{\dff T}\qff \in\qff \mathcal{D}\dff(\trf T\trf)$\sss
and\sss
$x_{\dff 0}\dff,\qff x_{\dff 1}
\qff \in\qff
\mathcal{N}_{\qff \overline{\mu}}
\off =\off 
\kernel\dff (\trf T^{\dff *}\qff -\qff \overline{\mu}\qff)$\sss
such\sss that\vspace{3pt}
\[
\quad
x
\off =\off
x_{\dff T}
\pff +\pff
A\trf(\trf A\qff -\qff \mu\trf)^{\fff -\dff 1}\dff x_{\dff 0}
\pff +\pff
(\trf A\qff -\qff \mu\trf)^{\fff -\dff 1}\dff x_{\dff 1}
\pff.
\]

\vspace{-12pt}\vspace{3pt}
See\qss \cite{s},\oss Proposition\qss 14.11{\fff}({\fff}i{\fff})\qss
and\dss Example\qss 14.5.\oss
Let\sss 
$K\off =\off \mathcal{N}_{\qff \overline{\mu}}$
and  
$\Gamma_{\fff 0}\dff,\pff \Gamma_{\dff 1}$
be defined\dss
by\sss the rules\sss
$\Gamma_{\fff 0}\dff x\off =\off x_{\dff 0}$,\qss
$\Gamma_{\dff 1}\dff x\off =\off x_{\dff 1}$.\oss
Then\sss
$(\qff K\dff,\off \Gamma_{\fff 0}\dff,\off \Gamma_{\dff 1}\trf)$
is\dss a boundary\sss triplet\sss for $T^{\dff *}$
such\sss that\sss $A\off =\off T_{\dff 0}$\nsp,\oss
which we will\sss call\dss the\qss \emph{boundary $\mu$\dnsp-triplet}.\oss
See\qss \cite{s},\oss Example\qss 14.5.\oss
We need a corollary of\dss the\dss Krein--Naimark\dss formula for\sss the boundary\sss $\mu$\dnsp-triplet,\oss
namely\sss the identity\qss (\ref{cayley-trans-mult})\qss below.\oss
This identity can be also deduced directly\sss from\sss the\dss von\dss Neumann\dss theory of\dss extensions.\oss
See\qss \cite{i4},\oss Section\qss 2\qss for\sss such a proof.\oss

For\sss the boundary\sss $\mu$\dnsp-triplet\sss one can compute\sss the gamma-field\sss
$\bm{\gamma}\qff(\trf z\trf)$\sss and\sss the\dss Weyl\sss function\sss 
$M\trf(\trf z\trf)$\nnsp.\oss
Let\sss us begin\sss with\sss
$\bm{\gamma}\qff(\qff \overline{\mu}\qff)$
and\sss
$M\trf(\qff \overline{\mu}\qff)$\nnsp.\oss
Suppose\sss that\sss $v\qff \in\qff K\off =\off \mathcal{N}_{\qff \overline{\mu}}$.\oss
Then\vspace{3pt}
\[
\quad
v\off =\off
(\trf A\qff -\qff \mu\trf)\qff
(\trf A\qff -\qff \mu\trf)^{\fff -\dff 1}\dff v
\off =\off
A\trf(\trf A\qff -\qff \mu\trf)^{\fff -\dff 1}\dff v
\qff -\qff
\mu\trf
(\trf A\qff -\qff \mu\trf)^{\fff -\dff 1}\dff v
\]

\vspace{-33pt}
\[
\quad
\phantom{v
\off }
=\off
A\trf(\trf A\qff -\qff \mu\trf)^{\fff -\dff 1}\dff v
\qff +\qff
(\trf A\qff -\qff \mu\trf)^{\fff -\dff 1}\dff (\trf -\qff \mu\dff v\trf)
\off.
\]

\vspace{-12pt}\vspace{3pt}
It\sss follows\sss that\sss\vspace{3pt}
\[
\quad
\Gamma_{\fff 0}\dff v\off =\off v_{\dff 0}\off =\off v
\quad
\mbox{and}\quad\dff
\Gamma_{\dff 1}\dff v\off =\off v_{\dff 1}\off =\off -\qff \mu\dff v
\off.
\]

\vspace{-12pt}\vspace{3pt}
Hence\sss
$\bm{\gamma}\qff(\qff \overline{\mu}\qff)\dff v\off =\off v$\sss
and\sss\vspace{3pt}
\[
\quad
M\trf(\qff \overline{\mu}\qff)\dff v
\off =\off 
\Gamma_{\dff 1}\dff \circ\qff \bm{\gamma}\trf(\qff \overline{\mu}\qff)\dff v
\off =\off
-\qff \mu\dff v 
\pff.
\]

\vspace{-12pt}\vspace{3pt}
In other\sss terms,\pss
$\bm{\gamma}\qff(\qff \overline{\mu}\qff)$\sss
is\dss the inclusion\sss $K\qff \ttoo\qff H$\sss
and\sss
$M\trf(\qff \overline{\mu}\qff)\off =\off -\qff \mu\dff \id_{\dff K}$\nsp.\oss
By\sss applying\sss the equality\qss (\ref{gamma-field})\qss to\sss
$z\fff,\off \overline{\mu}$\sss in\sss the roles of\dss $w\fff,\qff z$\sss
respectively,\oss we see\sss that\sss for\sss 
$z\qff \in\qff \rho\trf(\trf A\trf)$\vspace{3pt}
\[
\quad
\bm{\gamma}\qff(\trf z\trf)
\off =\off \left.
(\trf A\qff -\qff \overline{\mu}\qff)\qff
(\trf A\qff -\qff z\trf)^{\fff -\dff 1}\qff 
\right|\trf K
\pff.
\]

\vspace{-12pt}\vspace{3pt}
Let $P_{\dff K}\dff \colon\dff H\qff \ttoo\qff K$ be\sss the orthogonal\dss projection.\oss
Since $A$ is\dss self-adjoint,\oss it\sss follows\sss that\vspace{3pt}
\[
\quad
\bm{\gamma}\qff(\qff \overline{z}\qff)^{\fff *}
\off =\off 
P_{\dff K}\trf
(\trf A\qff -\qff \mu\trf)\qff
(\trf A\qff -\qff z\trf)^{\fff -\dff 1}\qff 
\pff
\]

\vspace{-12pt}\vspace{3pt}
for\sss $z\qff \in\qff \rho\trf(\trf A\trf)$\nnsp.\oss
Therefore\sss the\dss Krein--Naimark\dss formula\sss takes\sss the form\vspace{3pt}
\[
\quad
(\qff
T_{\fff \mathcal{B}}\qff -\qff z
\qff)^{\fff -\dff 1}
\pff -\pff
(\qff
A\qff -\qff z
\qff)^{\fff -\dff 1}
\]

\vspace{-33pt}
\[
\quad
=\off\dff 
(\trf A\qff -\qff \overline{\mu}\qff)\qff
(\trf A\qff -\qff z\trf)^{\fff -\dff 1}\qff
(\qff
\mathcal{B}\qff -\qff M\trf(\trf z\trf)
\qff)^{\fff -\dff 1}\qff
P_{\dff K}\trf
(\trf A\qff -\qff \mu\trf)\qff
(\trf A\qff -\qff z\qff)^{\fff -\dff 1}
\pff.
\]

\vspace{-12pt}\vspace{3pt}
Now\sss the equality\qss (\ref{weyl-function})\qss implies\sss that\vspace{4.5pt}
\[
\quad
M\trf(\trf z\trf) 
\off =\off\dff 
M\trf(\qff \overline{\mu}\qff)
\off +\off
(\trf z\qff -\qff \overline{\mu}\qff)\qff
\bm{\gamma}\qff(\qff \mu\qff)^{\fff *}\qff
\bm{\gamma}\qff(\qff \overline{z}\qff)
\] 

\vspace{-31.5pt}
\[
\quad
\phantom{M\trf(\trf z\trf) 
\off }
=\off\dff \left.
-\qff \mu
\off +\off
(\trf z\qff -\qff \overline{\mu}\qff)\qff
P_{\dff K}\trf
(\trf A\qff -\qff \mu\trf)\qff
(\trf A\qff -\qff \overline{\mu}\qff)^{\fff -\dff 1}\qff
(\trf A\qff -\qff \overline{\mu}\qff)\qff
(\trf A\qff -\qff z\trf)^{\fff -\dff 1}\qff 
\right|\trf K
\] 

\vspace{-31.5pt}
\[
\quad
\phantom{M\trf(\trf z\trf) 
\off }
=\off\dff \left.
-\qff \mu
\off +\off
(\trf z\qff -\qff \overline{\mu}\qff)\qff
P_{\dff K}\trf
(\trf A\qff -\qff \mu\trf)\qff
(\trf A\qff -\qff z\trf)^{\fff -\dff 1}\qff
\right|\trf K
\pff 
\] 

\vspace{-31.5pt}
\[
\quad
\phantom{M\trf(\trf z\trf) 
\off }
=\off\dff \left.
-\qff \mu
\off +\off
(\trf z\qff -\qff \overline{\mu}\qff)\qff
P_{\dff K}\trf
(\trf A\qff -\qff z\qff +\qff z\qff -\qff \mu\trf)\qff
(\trf A\qff -\qff z\trf)^{\fff -\dff 1}\qff
\right|\trf K
\pff 
\] 

\vspace{-31.5pt}
\[
\quad
\phantom{M\trf(\trf z\trf) 
\off }
=\off\dff \left.
-\qff \mu
\off +\off
(\trf z\qff -\qff \overline{\mu}\qff)\qff
\off +\off
(\trf z\qff -\qff \overline{\mu}\qff)\qff
(\trf z\qff -\qff \mu\trf)\qff
P_{\dff K}\trf
(\trf A\qff -\qff z\trf)^{\fff -\dff 1}\qff 
\right|\trf K
\]

\vspace{-31.5pt}
\[
\quad
\phantom{M\trf(\trf z\trf) 
\off }
=\off\dff \left.
(\trf z\qff -\qff 2\trf \re \mu\qff)
\off +\off
(\trf z\qff -\qff \overline{\mu}\qff)\qff
(\trf z\qff -\qff \mu\trf)\qff
P_{\dff K}\trf
(\trf A\qff -\qff z\trf)^{\fff -\dff 1}\qff
\right|\trf K
\pff.
\] 

\vspace{-12pt}\vspace{4.5pt}
In\sss the special\sss case of\dss $\mu\off =\off -\qff i$\sss we get\vspace{4.5pt}
\[
\quad
M\trf(\trf z\trf)
\off =\off\dff \left.
z\off +\off
(\trf z\qff +\qff i\qff)\qff
(\trf z\qff -\qff i\trf)\qff
P_{\dff K}\trf
(\trf A\qff -\qff z\trf)^{\fff -\dff 1}\qff
\right|\trf K
\pff 
\] 

\vspace{-31.5pt}
\[
\quad
\phantom{M\trf(\trf z\trf)
\off }
=\off\dff \left.
z\off +\off
(\trf z^{\dff 2}\qff +\qff 1\qff)\qff
P_{\dff K}\trf
(\trf A\qff -\qff z\trf)^{\fff -\dff 1}\qff
\right|\trf K
\pff, 
\] 

\vspace{-12pt}\vspace{4.5pt}
and\sss for\sss $z\off =\off -\qff i$\sss we\sss get\sss
$M\trf(\trf -\qff i\trf)
\off =\off
-\qff i$\nnsp.\oss
Hence\sss the\dss Krein--Naimark\trs formula\sss implies\sss
that\vspace{3pt}
\[
\quad
(\qff
T_{\fff \mathcal{B}}\qff +\qff i
\qff)^{\fff -\dff 1}
\pff -\pff
(\qff
A\qff +\qff i
\qff)^{\fff -\dff 1}
\]

\vspace{-33pt}
\[
\quad
=\off\dff 
(\trf A\qff -\qff i\qff)\qff
(\trf A\qff +\qff i\trf)^{\fff -\dff 1}\qff
(\qff
\mathcal{B}\qff +\qff i
\qff)^{\fff -\dff 1}\qff
P_{\dff K}\trf
(\trf A\qff +\qff i\trf)\qff
(\trf A\qff +\qff i\qff)^{\fff -\dff 1}\qff
\pff
\]

\vspace{-33pt}
\[
\quad
=\off\dff 
(\trf A\qff -\qff i\qff)\qff
(\trf A\qff +\qff i\trf)^{\fff -\dff 1}\qff
(\qff
\mathcal{B}\qff +\qff i
\qff)^{\fff -\dff 1}\qff
P_{\dff K}
\pff.
\]

\vspace{-12pt}\vspace{3pt}
Let\sss us rewrite\sss this formula\sss in\sss terms of\dss the\dss Cayley\dss transforms of\dss
$T_{\fff \mathcal{B}}\dff,\off A$\nnsp,\oss and $\mathcal{B}$\dnsp.\oss
As\dss is\dss well\dss known,\oss if\dss $B$\sss is\dss an operator\sss 
$K\qff \ttoo\qff K$\nnsp,\oss
then\sss 
$V\trf(\trf B\trf)
\off =\off 
1\qff -\qff 2\dff i\trf (\trf B\qff +\qff i\trf)^{\dff -\dff 1}$\nsp\dnsp.\oss
The above description of\dss $V\trf(\trf \mathcal{B}\trf)$
for self-adjoint\sss relations $\mathcal{B}$ implies\sss that\sss also\vspace{3pt}
\[
\quad
V\trf(\trf \mathcal{B}\trf)
\off =\off
1\qff - \qff 2\dff i\trf (\trf \mathcal{B}\qff +\qff i\trf)^{\dff -\dff 1}
\]

\vspace{-12pt}\vspace{3pt}
for such relations $\mathcal{B}$\dnsp.\oss
By combining\sss this with\sss the\sss last\sss equality\sss 
we see\sss that\vspace{3pt}
\[
\quad
V\trf(\trf T_{\fff \mathcal{B}}\trf)
\qff -\qff 
V\trf(\trf A\trf)
\off =\off
V\trf(\trf A\trf)\qff
\left(\qff
-\qff 2\dff i\trf (\trf \mathcal{B}\qff +\qff i\trf)^{\dff -\dff 1}
\qff\right)\qff
P_{\dff K}
\]

\vspace{-34.0pt}
\[
\quad
\phantom{V\trf(\trf T_{\fff \mathcal{B}}\trf)
\qff -\qff 
V\trf(\trf A\trf)
\off }
=\off
V\trf(\trf A\trf)\qff \bigl(\trf V\trf(\trf \mathcal{B}\trf)\qff -\qff 1\trf\bigr)\qff
P_{\dff K}
\pff.
\]

\vspace{-12pt}\vspace{-3pt}
It\sss follows\sss that\vspace{2.25pt}\vspace{-0.125pt}
\[
\quad
V\trf(\trf A\trf)^{\dff -\dff 1}\qff
V\trf(\trf T_{\fff \mathcal{B}}\trf)
\qff -\qff 
1
\off =\off
\bigl(\trf V\trf(\trf \mathcal{B}\trf)\qff -\qff 1\trf\bigr)\qff
P_{\dff K}
\quad\qff
\mbox{and}\quad
\]

\vspace{-34.5pt}
\[
\quad
V\trf(\trf A\trf)^{\dff -\dff 1}\qff
V\trf(\trf T_{\fff \mathcal{B}}\trf)
\off =\off
V\trf(\trf \mathcal{B}\trf)\qff
P_{\dff K}
\pff +\pff
\left(\trf
1\qff -\qff P_{\dff K}
\trf\right)
\pff.
\]

\vspace{-12pt}\vspace{2.25pt}
The operator\sss in\sss the right\sss hand side\dss is\dss 
the\dss Cayley\trs transform of\dss $\mathcal{B}$ considered as a relation\sss
not\sss in\sss the subspace $K\qff \subset\qff H$\nnsp,\oss but\sss in\sss $H$\nnsp.\oss
Let\sss us denote\sss this relation\sss by\sss $\mathcal{B}^{\dff H}$\dnsp.\oss
In other\sss terms,\pss $\mathcal{B}^{\dff H}$\sss is\dss the image of\dss $\mathcal{B}$
under\sss the inclusion\sss 
$K\dff \oplus\dff K\qff \ttoo\qff H\dff \oplus\dff H$\nnsp.\oss
Finally,\oss we get\vspace{3pt}
\[
\quad
V\trf(\trf A\trf)^{\dff -\dff 1}\qff
V\trf(\trf T_{\fff \mathcal{B}}\trf)
\off =\off
V\trf\left(\trf \mathcal{B}^{\dff H}\trf\right)
\quad
\mbox{and}\quad\dff
\]

\vspace{-34.5pt}
\begin{equation}
\label{cayley-trans-mult}
\quad
V\trf(\trf T_{\fff \mathcal{B}}\trf)
\off =\off
V\trf(\trf A\trf)\qff V\trf\left(\trf \mathcal{B}^{\dff H}\trf\right)
\qff.
\end{equation}

\vspace{-12pt}\vspace{3pt}
This\dss is\trs just\sss a more general\sss form of\trs 
Theorem\qss 14.20\qss form\qss \cite{s}.\oss

\mypar{Theorem.}{spectral-triplet-index}
\emph{Let $A$ be a self-adjoint\sss extension of\qss $T$
and\sss $(\qff K\dff,\qff \Gamma_{\fff 0}\dff,\qff \Gamma_{\fff 1}\trf)$
be\sss the corresponding\dss boundary $\mu$\dnsp-triplet\sss with\sss
$\mu\off =\off -\qff i$\nnsp.\oss
Let\sss $\mathcal{B}_{\fff w}\dff,\pff w\qff \in\qff W$\sss be a continuous family
of\dss self-adjoint\dss relations on\dss $K$\nnsp.\oss
If\qss the operator\sss $A$\sss is\trs Fredholm and\dss
$K$ is\dss finitely dimensional,\oss
then\sss the family}\vspace{0pt}
\[
\quad
T_{\fff \mathcal{B}_{\fff w}}\dff,\off w\qff \in\qff W
\]

\vspace{-12pt}\vspace{0pt}
\emph{is\qss Fredholm\qss and\dss its analytical\dss index\dss is\dss equal\dss 
to\sss the analytical\dss index of\dss 
$\mathcal{B}_{\fff w}\dff,\pff w\qff \in\qff W$\nnsp.}

\proof
Since $\mathcal{B}_{\fff w}\dff,\pff w\qff \in\qff W$\sss is\dss a continuous family,\oss
the map\sss 
$w\off \longmapsto\off V\trf\left(\trf \mathcal{B}_{\fff w}\trf\right)$\sss
is\dss continuous as a map from\sss $W$\sss to\sss 
the unitary\sss group $U\trf(\trf K\trf)$\sss of\dss $K$\dnsp.\oss
It\sss follows\sss that\sss the map\vspace{1.5pt}
\[
\quad
g\dff \colon\dff
w
\off \longmapsto\off 
V\trf\left(\trf \mathcal{B}_{\fff w}^{\dff H}\trf\right)
\]

\vspace{-12pt}\vspace{1.5pt}
is\dss continuous and maps\sss $W$\sss into\sss the space $U\fred$ associated\sss with $H$\nnsp.\oss
Since $A$\sss is\trs Fredholm,\pss  
$V\trf(\trf A\trf)$\sss is\trs Fredholm-unitary,\oss
i.e.\dss $V\trf(\trf A\trf)\qff \in\qff U\fred$\dnsp.\oss
Since $K$\sss is\dss finitely dimensional,\oss
the map\vspace{1.5pt}
\[
\quad
w
\off \longmapsto\off 
V\trf(\trf A\trf)\qff V\trf\left(\trf \mathcal{B}_{\fff w}^{\dff H}\trf\right)
\]

\vspace{-12pt}\vspace{1.5pt}
takes values in $U\fred$\dnsp.\oss
Cf.\dss the discussion preceding\dss Lemma\qss \ref{u-action}.\oss
Clearly,\oss it\dss is\dss continuous.\oss
By\sss applying\sss the equality\qss (\ref{cayley-trans-mult})\qss to
$\mathcal{B}\off =\off \mathcal{B}_{\fff w}\dff,\pff w\qff \in\qff W$\nnsp,\oss
we now see\sss that\vspace{1.5pt}
\[
\quad
w
\off \longmapsto\off
V\trf(\trf T_{\fff \mathcal{B}_{\fff w}}\trf)
\off =\off
V\trf(\trf A\trf)\qff V\trf\left(\trf \mathcal{B}_{\fff w}^{\dff H}\trf\right)
\]

\vspace{-12pt}\vspace{1.5pt}
is\dss a continuous map\sss
$W\qff \ttoo\qff U\fred$\dnsp.\oss
It\sss follows\sss that\sss the family\sss
$T_{\fff \mathcal{B}_{\fff w}}\dff,\off w\qff \in\qff W$\sss
is\trs Fredholm.\oss
Without\sss any\sss loss of\dss generality we can assume\sss that\sss
$U\trf(\trf K\trf)\qff \subset\qff U\dff(\trf \infty\trf)$\nnsp.\oss
Let\sss $f\dff \colon\dff W\qff \ttoo\qff U\fred$\sss be\sss the constant\sss map\sss 
$w
\off \longmapsto\off
V\trf\left(\trf A\trf\right)$\nnsp.\oss
The analytical\dss index of\dss the corresponding\sss family\dss is\dss zero.\oss
It\dss remains\sss to apply\trs Lemma\qss \ref{u-action}\qss to $f$ and $g$\nnsp.\oss  \eproof

\myuppar{Uniqueness.}
When\sss a boundary\sss triplet\sss for $T^{\dff *}$ exists,\oss
then\sss it\dss is\dss essentially\sss unique.\oss
Suppose\sss that\sss 
$(\qff K\dff,\off \Gamma_{\fff 0}\dff,\off \Gamma_{\fff 1}\trf)$
and\sss
$(\qff K\fff'\dff,\off \Gamma\fff'_{0}\dff,\off \Gamma\fff'_{1}\trf)$
are\sss two boundary\sss triplets for $T^{\dff *}$\dnsp.\oss
Let\vspace{4.5pt}
\[
\quad
\Gamma_{\bullet}
\off =\off
\Gamma_{\fff 0}\dff \oplus\dff \Gamma_{\fff 1}
\dff \colon\dff
\mathcal{D}\dff(\trf T^{\dff *}\trf)
\qff \ttoo\qff
K\dff \oplus\dff K
\quad
\mbox{and}\quad
\Gamma\fff'_{\bullet}
\off =\off
\Gamma\fff'_{0}\dff \oplus\dff \Gamma\fff'_{1}
\dff \colon\dff
\mathcal{D}\dff(\trf T^{\dff *}\trf)
\qff \ttoo\qff
K\fff'\dff \oplus\dff K\fff'
\qff,
\]

\vspace{-12pt}\vspace{4.5pt}
and\dss let\sss $[\trf \bullet\fff,\qff \bullet\trf]$\nnsp,\qss
$[\trf \bullet\fff,\qff \bullet\trf]\fff'$
be\sss the\dss Hermitian\sss scalar\sss products on\sss
$K\dff \oplus\dff K$\nnsp,\qss $K\fff'\dff \oplus\dff K\fff'$
respectively defined\sss in\dss Section\qss \ref{relations}.\oss
Then\sss there\dss is\dss a\sss bounded\dss isomorphisms\sss
$U\dff \colon\dff 
K\dff \oplus\dff K
\qff \ttoo\qff
K\fff'\dff \oplus\dff K\fff'$\sss
such\sss that\sss
$\Gamma\fff'_{\bullet}
\off =\off
U\dff \circ\dff \Gamma_{\bullet}$
and\sss $U$\sss is\dss unitary\sss with\sss respect\sss to\sss the products\sss
$[\trf \bullet\fff,\qff \bullet\trf]$
and\sss $[\trf \bullet\fff,\qff \bullet\trf]\fff'$\dnsp,\oss
i.e.\vspace{3pt}
\[
\quad
[\trf u\fff,\qff v\qff]
\off =\off
[\trf U\dff u\fff,\qff U\dff v\qff]\fff'
\]

\vspace{-12pt}\vspace{3pt}
for every\sss $u\fff,\qff v\qff \in\qff K$\nnsp.\oss
Moreover,\oss such an\sss isomorphism\sss $U$\sss is\dss unique.\oss
See\qss \cite{bhs},\oss Section\qss 1.8\qss and\trs Proposition\qss 2.5.1.\oss
Since $U$\sss is\dss bounded\sss isomorphism,\oss $U$\sss takes
closed\sss relations on $K$\sss to closed\sss relations on $K\fff'$\dnsp.\oss
Since\sss the isomorphism\sss $U$\sss is\dss unitary\sss with\sss respect\sss to
$[\trf \bullet\fff,\qff \bullet\trf]$
and\sss $[\trf \bullet\fff,\qff \bullet\trf]\fff'$\dnsp,\oss
it\dss takes self-adjoint\sss relations on\sss $K$\sss to
self-adjoint\sss relations on\sss $K\fff'$\dnsp.\oss
Clearly,\oss if\dss $\mathcal{B}$\sss is\dss a self-adjoint\sss relation on\sss $K$\nnsp,\oss
then\sss the extension $T_{\fff \mathcal{B}}$ defined\dss by\sss
$(\qff K\dff,\off \Gamma_{\fff 0}\dff,\off \Gamma_{\fff 1}\trf)$
and $\mathcal{B}$\sss is\dss equal\sss to\sss the extension\sss
$T_{\fff U\trf(\trf \mathcal{B}\trf)}$ defined\dss by\sss
$(\qff K\fff'\dff,\off \Gamma\fff'_{0}\dff,\off \Gamma\fff'_{1}\trf)$
and $U\trf(\trf \mathcal{B}\trf)$\nnsp.\oss

\myuppar{Operators related\dss to abstract\dss boundary\sss problems.}
Suppose\sss that\sss we are in\sss the framework of\trs Section\qss \ref{abstract-index}.\oss
In\sss particular,\pss $A$ satisfies\sss the\sss abstract\dss Lagrange\dss identity\qss (\ref{lagrange}),\pss 
$A\fff,\qff \Pi$\sss is\dss a self-adjoint\dss elliptic regular\sss boundary\sss problem,\oss
and\sss all\sss assumptions of\trs Section\qss \ref{abstract-index}\qss hold.\oss
Suppose\sss that\sss $H^{\dff \partial}\off =\off H\dff \oplus\dff H$\sss
for some\dss Hilbert\dss space $H$\nnsp,\oss that $\Pi$\sss is\dss the projection 
$H\dff \oplus\dff H\qff \ttoo\qff H\dff \oplus\dff 0$\nnsp,\oss 
and\dss that\sss with respect\sss to\sss the decomposition\sss 
$H^{\dff \partial}\off =\off H\dff \oplus\dff H$\sss
the operator\sss $\Sigma$\sss has\sss the form\vspace{1.5pt}
\[
\quad
\Sigma
\off =\off\dff
\begin{pmatrix}
\off 0 &
1 \qff\off
\vspace{4.5pt} \\
\off\dff 1 &
0 \qff\off 
\end{pmatrix}
\off.
\]

\vspace{-12pt}\vspace{1.5pt}
Let\sss $K$\sss be a finitely dimensional\sss subspace of\dss 
$\kernel\fff \Pi_{\dff 1/2}\qff \subset\qff H_{\dff 1/2}^{\dff \partial}$
and\sss let\sss 
$K^{\dff \perp}$ be its orthogonal\sss complement\sss in\sss 
$\kernel\fff \Pi\qff \subset\qff H^{\dff \partial}$\dnsp.\oss
Then\sss 
$L\off =\off K^{\dff \perp}\dff \cap\qff \kernel\fff \Pi_{\dff 1/2}$\sss
is\dss a closed subspace of\dss $\kernel\fff \Pi_{\dff 1/2}$ such\sss that\sss
$\kernel\fff \Pi_{\dff 1/2}\off =\off L\qff +\qff K$\sss
and\sss $L\dff \cap\dff K\off =\off 0$\nnsp.\oss
Let\sss 
$D
\off =\off
\gamma^{\dff -\dff 1}\dff(\trf \Sigma\trf(\trf L\trf) \trf)$\nnsp.\oss
The subspace $D$ contains $\kernel\fff \gamma$ and\dss hence\dss is\dss dense in $H_{\trf 0}$\nsp.\oss
Let\sss $A_{\trf D}\dff \colon\dff D\qff \ttoo\qff H_{\trf 0}$\sss be\sss
the restriction of\sss $A$\sss to $D$\nnsp.\oss

\mypar{Lemma.}{adjoint-ad}
\emph{The adjoint\sss operator\dss $A_{\trf D}^{\fff *}$\sss is\dss 
equal\dss to\sss the restriction of\qss $A$\sss to}\vspace{3pt}
\[
\quad
\mathbb{D}
\off =\off\dff
\gamma^{\dff -\dff 1}\trf\bigl(\qff
(\trf \image\fff \Pi_{\dff 1/2}\trf)\qff +\qff K
\qff\bigr)
\off =\off
\bigl(\qff
\kernel\fff \Pi_{\dff 1/2}\dff \circ\dff \gamma
\qff\bigr)
\off +\off 
\gamma^{\dff -\dff 1}\trf(\trf K\trf)
\pff.
\]

\vspace{-12pt}\vspace{3pt}
\proof
If\dss $w\qff \in\qff \mathbb{D}$ and\sss $u\qff \in\qff D$\nnsp,\oss
then\vspace{3pt}
\[
\quad
\sco{\dff A\dff u\fff,\qff w \dff}_{\dff 0}
\qff -\qff
\sco{\dff u\dff,\qff A\dff w \dff}_{\dff 0}
\off =\off
\sco{\dff i\trf \Sigma\dff \gamma\dff u\dff,\qff \gamma\dff w \dff}_{\dff \partial}
\off =\off
0
\off.
\]

\vspace{-12pt}\vspace{3pt}
It\sss follows\sss that\sss 
$w\qff \in\qff \mathcal{D}\dff(\trf A_{\trf D}^{\fff *} \trf)$
and\sss 
$A_{\trf D}^{\fff *}\dff w\off =\off A\dff w$\nnsp.\oss
Therefore\sss 
$\mathbb{D}\qff \subset\pff \mathcal{D}\dff(\trf A_{\trf D}^{\fff *} \trf)$\sss
and\sss $A_{\trf D}^{\fff *}$\sss is\dss equal\dss to $A$ on\sss $\mathbb{D}$\nnsp.\oss
It\sss remains to prove\sss that\sss
$\mathbb{D}\off =\off \mathcal{D}\dff(\trf A_{\trf D}^{\fff *} \trf)$\nnsp.\oss
Let\sss \vspace{3pt}
\[
\quad
\mathcal{K}_{\dff +}
\off =\off 
\kernel\fff (\trf A_{\trf D}^{\fff *}\qff -\qff i\trf)
\quad
\mbox{and}\quad
\mathcal{K}_{\dff -}
\off =\off 
\kernel\fff (\trf A_{\trf D}^{\fff *}\qff +\qff i\trf)
\off.
\]

\vspace{-12pt}\vspace{3pt}
Then\sss
$\mathcal{D}\dff(\trf A_{\trf D}^{\fff *} \trf)
\off =\off
D\qff +\qff \mathcal{K}_{\dff +}\qff +\qff \mathcal{K}_{\dff -}$\sss
and\dss hence\sss the codimension of\dss $D$\sss in\sss
$\mathcal{D}\dff(\trf A_{\trf D}^{\fff *} \trf)$\sss
is\dss equal\dss to\sss
$\dim\dff \mathcal{K}_{\dff +}\qff +\pff \dim\dff \mathcal{K}_{\dff -}$\nsp.\oss
The restriction $A_{\trf \Gamma}$ of\dss $A$\sss to\sss
$\kernel\fff \Pi_{\dff 1/2}\dff \circ\dff \gamma$\sss
is\dss a self-adjoint\sss operator,\oss
and\dss hence\dss is\dss a self-adjoint\sss extension of\dss $A_{\trf D}$\nsp.\oss
Therefore\sss $A_{\trf \Gamma}$\sss
is\dss equal\sss to\sss the restriction of\dss $A_{\trf D}^{\fff *}$\sss
to\sss $D\qff +\qff \image\fff (\trf 1\qff +\qff U \trf)$\sss
for some isometry\sss
$U\dff \colon\dff 
\mathcal{K}_{\dff +}\qff \ttoo\qff \mathcal{K}_{\dff -}$\nsp.\oss
It\sss follows\sss that\sss the  codimension of\dss $D$\sss in\sss
$\kernel\fff \Pi_{\dff 1/2}\dff \circ\dff \gamma$\sss
is\dss equal\dss to\sss 
$\dim\dff \mathcal{K}_{\dff +}\off =\off \dim\dff \mathcal{K}_{\dff -}$\nsp.\oss
In\sss turn,\sss this\sss implies\sss that\sss the codimension of\dss $D$\sss 
in\sss $\mathbb{D}$\sss is\dss equal\dss to\sss 
$\dim\dff \mathcal{K}_{\dff +}\qff +\pff \dim\dff \mathcal{K}_{\dff -}$\nsp.\oss
Therefore\sss the codimensions of\dss $D$\sss in\sss $\mathbb{D}$
and\sss in\sss
$\mathcal{D}\dff(\trf A_{\trf D}^{\fff *} \trf)$\sss
are\sss the same and\sss finite.\oss
Since\sss 
$\mathbb{D}
\qff \subset\qff 
\mathcal{D}\dff(\trf A_{\trf D}^{\fff *} \trf)$\nnsp,\oss
this implies\sss that\sss
$\mathbb{D}
\off =\off 
\mathcal{D}\dff(\trf A_{\trf D}^{\fff *} \trf)$\nnsp.\oss \eproof

\myuppar{Boundary\sss triplets and extensions related\dss to abstract\dss boundary\sss problems.}
Let\sss us keep\sss the above assumptions about\sss $A\fff,\qff D$\nnsp,\oss etc.\qss
and construct\sss a boundary\sss triplet\sss for $A_{\trf D}^{\fff *}$\nsp.\oss
Let\sss $\pi$\sss  
be\sss the\vspace{-1pt}\\ orthogonal\sss projection
$H^{\dff \partial}\qff \ttoo\qff K$\nnsp.\oss 
If\dss
$u\qff \in\qff \Sigma\trf(\trf K\trf)$\nnsp,\qss
$v\qff \in\qff K$\nnsp,\oss
and\sss
$x\off =\off u\qff +\qff v$\nnsp,\qss
then\vspace{3pt}
\[
\quad
u\off =\off \Sigma\dff \circ\dff \pi\dff \circ\dff \Sigma\trf x
\qff,\quad
v\off =\off \pi\trf x
\pff.
\]

\vspace{-12pt}\vspace{3pt}
Therefore,\oss if\dss
$u\fff,\qff a\qff \in\qff \Sigma\trf(\trf K\trf)$\nnsp,\qss
$v\fff,\qff b\qff \in\qff K$\nnsp,\oss
and\sss
$x\off =\off u\qff +\qff v$\nnsp,\qss
$y\off =\off a\qff +\qff b$\nnsp,\oss
then\vspace{3pt}
\[
\quad
\sco{\dff i\trf \Sigma\trf x\dff,\qff y \dff}_{\dff \partial}
\off =\off
\sco{\dff i\trf \Sigma\dff u\dff,\qff b \dff}_{\dff \partial}
\off +\off
\sco{\dff i\trf \Sigma\dff v\dff,\qff a \dff}_{\dff \partial}
\off =\off
\sco{\dff i\trf \Sigma\dff u\dff,\qff b \dff}_{\dff \partial}
\off -\off
\sco{\dff v\dff,\qff i\trf \Sigma\dff a \dff}_{\dff \partial}
\pff
\]

\vspace{-33pt}
\[
\quad
\phantom{\sco{\dff i\trf \Sigma\dff x\dff,\qff y \dff}_{\dff \partial}
\off }
=\off
\sco{\dff \pi\dff \circ\dff i\trf \Sigma\trf x\dff,\qff \pi\trf y \dff}_{\dff \partial}
\off -\off
\sco{\dff \pi\trf x\dff,\qff \pi\dff \circ\dff i\trf \Sigma\dff y \dff}_{\dff \partial}
\pff,
\]

\vspace{-12pt}\vspace{3pt}
where at\sss the\sss last\sss step we used\dss that\sss 
$\Sigma\dff \circ\dff \Sigma\off =\off \id$\nnsp.\oss
Therefore,\oss if\dss
$x\fff,\qff y\qff \in\qff 
(\trf \image\fff \Pi\trf)\qff +\qff K$\nnsp,\oss
then\vspace{3pt}
\begin{equation}
\label{i-sigma}
\quad
\sco{\dff i\trf \Sigma\trf x\dff,\qff y \dff}_{\dff \partial}
\off =\off
\sco{\dff \pi\dff \circ\dff i\trf \Sigma\trf x\dff,\qff \pi\trf y \dff}_{\dff \partial}
\off -\off
\sco{\dff \pi\trf x\dff,\qff \pi\dff \circ\dff i\trf \Sigma\dff y \dff}_{\dff \partial}
\pff.
\end{equation}

\vspace{-12pt}\vspace{3pt}
Let\sss $\Gamma_{\fff 0}$ and\sss $\Gamma_{\fff 1}$\sss
be\sss the restrictions\sss to\sss $\mathbb{D}$\sss of\dss the maps\vspace{3pt}
\[
\quad
\pi\dff \circ\dff \gamma
\off\qff \mbox{and}\off\qff 
\pi\dff \circ\dff i\trf \Sigma\dff \circ\dff \gamma
\dff\colon\dff
H_{\dff 1}\qff \ttoo\qff K
\]

\vspace{-12pt}\vspace{3pt}
respectively.\oss
Then\sss $\Gamma_{\fff 0}\dff \oplus\dff \Gamma_{\fff 1}$\sss is\dss surjective and\sss
the\sss identity\qss (\ref{lagrange})\qss
together\sss with\qss (\ref{i-sigma})\qss implies\sss that\sss\vspace{3pt}
\[
\quad
\sco{\dff A\dff x\fff,\qff y \dff}_{\dff 0}
\qff -\qff
\sco{\dff x\dff,\qff A\dff y \dff}_{\dff 0}
\off =\off
\sco{\dff \Gamma_{\dff 1}\dff x\dff,\qff \Gamma_{\fff 0}\dff y \dff}_{\dff \partial}
\off -\off
\sco{\dff \Gamma_{\fff 0}\dff x\dff,\qff \Gamma_{\dff 1}\dff y \dff}_{\dff \partial}
\]

\vspace{-12pt}\vspace{3pt}
for every\sss $x\fff,\qff y\qff \in\qff \mathbb{D}$\nnsp.\oss
Since\sss 
$\mathbb{D}\off =\off \mathcal{D}\dff(\trf A_{\trf D}^{\fff *} \trf)$
and\sss $A_{\trf D}^{\fff *}$\sss is\dss the restriction of\sss $A$\nnsp,\oss
this means\sss that\sss the\sss triple\sss
$(\trf K\dff,\off \Gamma_{\fff 0}\dff,\off \Gamma_{\dff 1}\trf)$\sss
is\dss a boundary\sss triplet\sss for $A_{\trf D}^{\fff *}$\nsp.\oss
Clearly,\pss\vspace{3pt}
\[
\quad
\kernel\fff \Gamma_{\fff 0}\off =\off
\gamma^{\dff -\dff 1}\dff(\trf \image\fff \Pi_{\dff 1/2} \trf)
\off =\off 
\kernel\fff \Gamma
\]

\vspace{-12pt}\vspace{3pt}
and\dss hence\sss $(\trf A_{\trf D}\trf)_{\dff 0}\off =\off A_{\trf \Gamma}$\nsp.\oss

A closed\sss relation $\mathcal{C}$ on $K$ defines a closed extension\sss
$(\trf A_{\trf D}\trf)_{\dff \mathcal{C}}$ of\dss $A_{\trf D}$\nsp.\oss
Let\sss us describe\sss this extension\sss in\sss the\sss language of\trs Section\qss \ref{abstract-index}.\oss
The domain of\dss $(\trf A_{\trf D}\trf)_{\dff \mathcal{C}}$\sss is\dss 
the set\sss of\dss all\sss $x\qff \in\qff \mathbb{D}$ such\sss that\sss
$(\trf \Gamma_{\fff 0}\dff x\fff,\pff \Gamma_{\fff 1}\dff x\trf)
\qff \in\qff
\mathcal{C}
\qff \subset\qff
K\dff \oplus\dff K$\nnsp,\oss
i.e.\dss
$(\trf \pi\dff \circ\dff \gamma\dff x\fff,\pff 
\pi\dff \circ\dff i\trf \Sigma\dff \circ\dff \gamma\dff x \trf)
\qff \in\qff
\mathcal{C}$\nnsp,\oss
or,\oss equivalently,\vspace{3pt}
\[
\quad
\gamma\dff x
\off \in\off 
\Sigma\trf(\trf L\trf)\qff +\pff J^{\dff -\dff 1}\dff(\trf \mathcal{C}\trf)
\pff,
\]

\vspace{-12pt}\vspace{3pt}
where\sss $J$\sss is\dss the map 
$\Sigma\trf(\trf K\trf)\dff \oplus\dff K
\qff \ttoo\qff
K\dff \oplus\dff K$\sss
defined\dss by\sss
$J\trf(\trf u\fff,\qff v\trf)
\off =\off
(\trf v\fff,\qff i\trf \Sigma\dff u\trf)$\nnsp.\oss
Recall\dss that\sss $K^{\dff \perp}$ is\dss the orthogonal\sss complement\sss of\dss $K$\sss in\sss 
$\kernel\fff \Pi$\nnsp.\oss
Clearly,\vspace{3pt}
\[
\quad
\Sigma\trf(\trf L\trf)\qff +\pff J^{\dff -\dff 1}\dff(\trf \mathcal{C}\trf)
\off =\off
\left(\qff
\Sigma\trf(\trf K^{\dff \perp}\trf)\qff +\pff J^{\dff -\dff 1}\dff(\trf \mathcal{C}\trf)
\qff\right)
\qff\fff \cap\qff\fff
H_{\dff 1/2}
\pff.
\]

\vspace{-12pt}\vspace{3pt}
Let\sss $\Pi_{\dff \mathcal{C}}$\sss be\sss the orthogonal\sss projection\sss
$H^{\dff \partial}
\qff \ttoo\qff 
\Sigma\trf(\trf K^{\dff \perp}\trf)\qff +\pff J^{\dff -\dff 1}\dff(\trf \mathcal{C}\trf)$\nnsp.\oss
Let\sss us\sss point\sss out\sss for\sss the future applications\sss
that\sss $\Pi_{\dff \mathcal{C}}$\sss is\dss equal\dss to\sss $\Pi$\sss
on\sss the orthogonal\sss complement\sss of\dss
$\Sigma\trf(\trf K\trf)\qff +\qff K$\sss
and\sss hence differs from $\Pi$\sss by an operator\sss of\dss finite rank.\oss
Since\sss $K\qff \subset\qff H_{\dff 1/2}$\nsp,\oss
the projection\sss $\Pi_{\dff \mathcal{C}}$\sss leaves $H_{\dff 1/2}$\sss invariant.\oss
Therefore we can define\sss corresponding operators\sss $\Gamma_{\fff \mathcal{C}}$\sss
and\sss $A_{\trf \Gamma_{\dff \mathcal{C}}}$\nsp.\oss

\mypar{Lemma.}{boundary-problems-triplets}
\emph{The extension\sss $(\trf A_{\trf D}\trf)_{\dff \mathcal{C}}$\sss is\dss
equal\dss to\sss $A_{\trf \Gamma_{\dff \mathcal{C}}}$\nsp.\oss}

\proof
Lemma\qss \ref{adjoint-ad}\qss implies\sss that\sss
$(\trf A_{\trf D}\trf)_{\fff \mathcal{C}}$ 
is\dss equal\dss to\sss the restriction of\dss $A$\sss to\sss 
$\gamma^{\dff -\dff 1}\trf(\trf \image\fff \Pi_{\dff \mathcal{C}}\trf)$\nnsp.\oss
Since\sss $A_{\trf \Gamma_{\dff \mathcal{C}}}$\sss is\dss also\sss equal\dss
to\sss the restriction of\dss $A$\sss to\sss this subspace,\pss
$(\trf A_{\trf D}\trf)_{\dff \mathcal{C}}
\off =\off
A_{\trf \Gamma_{\dff \mathcal{C}}}$\nsp.\oss  \eproof

\mypar{Lemma.}{two-self-adjointness}
\emph{The boundary\sss problem\sss $A\fff,\pff \Pi_{\dff \mathcal{C}}$
is\trs self-adjoint\trs if\trs and\dss only\trs if\qss 
$\mathcal{C}$\sss is\trs self-adjoint.}

\proof
Since $\Pi_{\dff \mathcal{C}}$\sss is\dss an orthogonal\sss projection,\oss
the boundary\sss problem\sss $A\fff,\off \Pi_{\dff \mathcal{C}}$\sss
is\dss self-adjoint\dss if\trs and\dss only\trs if\dss the subspace\sss\vspace{3pt}
\[
\quad 
\image\fff \Pi_{\dff \mathcal{C}}
\off =\off
\Sigma\trf(\trf K^{\dff \perp}\trf)\qff +\pff J^{\dff -\dff 1}\dff(\trf \mathcal{C}\trf)
\]

\vspace{-12pt}\vspace{3pt}
is\dss equal\dss to\sss its orthogonal\sss
complement\sss with respect\sss to\sss
the scalar product $\sco{\trf i\trf \Sigma\trf \bullet\fff,\qff \bullet\trf}$
on $H\dff \oplus\dff H$\nnsp.\oss
Since\sss the orthogonal\sss complement\sss of\dss $K^{\dff \perp}$\sss
in $H\dff \oplus\dff H$\sss is\dss equal\dss to\sss
$\Sigma\trf(\trf K^{\dff \perp}\trf)
\qff +\qff
\Sigma\trf(\trf K\trf)
\qff +\qff
K$\nnsp,\oss
this condition\dss is\dss equivalent\sss to\sss
$J^{\dff -\dff 1}\dff(\trf \mathcal{C}\trf)$\sss
being equal\dss to its orthogonal\sss complement\sss with respect\sss to\sss
the scalar product $\sco{\trf i\trf \Sigma\trf \bullet\fff,\qff \bullet\trf}$ 
on $\Sigma\trf(\trf K\trf)\dff \oplus\dff K$\nnsp.\oss 
At\sss the same\sss time $\mathcal{C}$\sss is\dss self-adjoint\dss
if\trs and\dss only\trs if\dss $\mathcal{C}$\sss is\dss equal\dss to\sss
its orthogonal\sss complement\sss with respect\sss to\sss
the scalar product\sss $[\trf \bullet\fff,\qff \bullet\trf]$\sss 
from\dss Section\qss \ref{relations}\qss with $K$ in\sss the role of\dss $H$\nnsp.\oss
Therefore\sss it\dss is\dss sufficient\sss to prove\sss that\sss $J$\sss
is\dss an\sss isometry\sss between\sss
$\Sigma\trf(\trf K\trf)\dff \oplus\dff K$\sss with\sss the product\sss
$\sco{\trf i\trf \Sigma\trf \bullet\fff,\qff \bullet\trf}$\sss
and\sss $K\dff \oplus\dff K$\sss  with\sss the product\sss
$i\trf [\trf \bullet\fff,\qff \bullet\trf]$

Let\sss
$(\trf u\fff,\qff v\trf)\fff,\off 
(\trf a\fff,\qff b\trf)
\off \in\off
\Sigma\trf(\trf K\trf)\dff \oplus\dff K$\nnsp.\oss
Then\vspace{4.5pt}
\[
\quad
\sco{\qff i\trf \Sigma\trf(\trf u\fff,\qff v\trf)\fff,\off
(\trf a\fff,\qff b\trf)\qff}
\off =\off
\sco{\trf (\trf i\trf \Sigma\dff v\fff,\qff i\trf \Sigma\dff u\trf)\fff,\off 
(\trf a\fff,\qff b\trf)\qff}
\]

\vspace{-33pt}
\[
\quad
\phantom{\sco{\qff i\trf \Sigma\trf(\trf u\fff,\qff v\trf)\fff,\off
(\trf a\fff,\qff b\trf)\qff}
\off }
=\off
i\trf \sco{\trf \Sigma\dff v\fff,\qff a\trf}
\off +\off
i\trf \sco{\trf \Sigma\dff u\fff,\qff b\trf}
\pff.
\]

\vspace{-12pt}\vspace{3pt}
At\sss the same\sss time\vspace{3pt}
\[
\quad
\bigl[\pff
J\trf(\trf u\fff,\qff v\trf)\fff,\off 
J\trf(\trf a\fff,\qff b\trf)
\qff\bigr]
\off =\off
\bigl[\qff
(\trf v\fff,\qff i\trf \Sigma\dff u\trf)\fff,\off 
(\trf b\fff,\qff i\trf \Sigma\dff a\trf)
\qff\bigr]
\off =\off
i\trf\sco{\trf v\fff,\qff i\trf \Sigma\dff a\trf}
\off -\off
i\trf\sco{\trf i\trf \Sigma\dff u\fff,\qff b\trf}
\]

\vspace{-33pt}
\[
\quad\phantom{\bigl[\pff
J\trf(\trf u\fff,\qff v\trf)\fff,\off 
J\trf(\trf a\fff,\qff b\trf)
\qff\bigr]
\off }
=\off
\sco{\trf v\fff,\qff \Sigma\dff a\trf}
\off +\off
\sco{\trf \Sigma\dff u\fff,\qff b\trf}
\off =\off
\sco{\trf \Sigma\dff v\fff,\qff a\trf}
\off +\off
\sco{\trf \Sigma\dff u\fff,\qff b\trf}
\pff.
\]

\vspace{-12pt}\vspace{4.5pt}
Hence\sss 
$\sco{\qff i\trf \Sigma\trf(\trf u\fff,\qff v\trf)\fff,\off
(\trf a\fff,\qff b\trf)\qff}
\off =\off
i\qff \bigl[\pff
J\trf(\trf u\fff,\qff v\trf)\fff,\off 
J\trf(\trf a\fff,\qff b\trf)
\qff\bigr]$\nnsp,\oss
i.e.\qss $J$\sss is\trs indeed\sss an\sss isometry.\oss  \eproof

\myuppar{Families of\dss extensions and\sss the analytical\dss index.}
Let\sss us\sss keep\sss the above assumptions about\sss
$A\fff,\off \Pi\fff,\off \Sigma$\sss and\sss $K$\nnsp.\oss
Let\sss $W$\sss be a\sss topological\sss space and\sss
$\mathcal{C}_{\fff w}\dff,\pff w\qff \in\qff W$\sss
be a continuous family of\dss self-adjoint\sss relations\sss
$\mathcal{C}_{\fff w}\qff \subset\qff K\dff \oplus\dff K$\nnsp.\oss
Since\sss $K$\sss is\dss finitely dimensional,\oss the relations\sss
$\mathcal{C}_{\fff w}$ are automatically\trs Fredholm.\oss
As we pointed out\sss in\dss Section\qss \ref{relations},\oss
together\sss with continuity of\dss the family\sss
$\mathcal{C}_{\fff w}\dff,\pff w\qff \in\qff W$\sss
this\dss implies\sss that\sss this family\dss is\trs Fredholm.\oss
Let\sss 
$A_{\dff w}\off =\off A_{\trf \Gamma_{\dff \mathcal{C}}}$\nsp.\oss
Lemmas\qss \ref{boundary-problems-triplets}\qss and\qss \ref{adjoint-ad}\qss
imply\sss that\sss $A_{\dff w}$\sss is\dss equal\dss to\sss
extension $(\trf A_{\trf D}\trf)_{\dff \mathcal{C}}$\sss of\sss $A_{\trf D}$\sss
defined\sss by $\mathcal{C}_{\dff w}$ and\sss the boundary\sss triplet\sss
$(\trf K\dff,\off \Gamma_{\fff 0}\dff,\off \Gamma_{\dff 1}\trf)$
constructed\sss above and\dss hence\sss 
to\sss the restriction of\sss $A$\sss to\vspace{3pt}
\[
\quad
\mathcal{D}\dff\left(\qff A_{\dff w} \qff\right)
\off =\off\dff
\left\{\qff \left.
x\qff \in\qff \mathbb{D}
\pff \right|\qff
(\trf \Gamma_{\fff 0}\dff x\fff,\pff \Gamma_{\fff 1}\dff x\trf)
\qff \in\qff
\mathcal{B}_{\fff w} 
\pff\right\} 
\qff.
\]

\vspace{-12pt}\vspace{3pt}
\mypar{Theorem.}{relative-index}
\emph{The family\dss
$A_{\dff w}\dff,\pff w\qff \in\qff W$\dss is\trs Fredholm\trs and\dss its\dss
analytical\dss index of\trs the  equal\dss to\sss 
the analytical\dss index of\trs the\dss family\dss
$\mathcal{C}_{\fff w}\dff,\pff w\qff \in\qff W$\nnsp.\oss}

\proof
The operator\sss $A_{\trf \Gamma}$\sss is\dss a self-adjoint\sss extension
of\dss $A_{\trf D}$\nsp.\oss
Let\sss us apply\sss the construction of\dss the standard example\sss
to $A_{\trf D}$ and $A_{\trf \Gamma}$ in\sss the roles of\dss $T$ and $A$ respectively.\oss
Let\sss $(\qff K\fff'\dff,\off \Gamma\fff'_{0}\dff,\off \Gamma\fff'_{1}\trf)$
be\sss the\sss boundary $\mu$\dnsp-triplet\sss with $\mu\off =\off -\qff 1$\nnsp.\oss
The assumptions of\trs Section\qss \ref{abstract-index}\qss ensure\sss that\sss
the operator $A_{\trf \Gamma}$\sss is\trs Fredholm.\oss
The uniqueness properties of\dss boundary\sss triplet\sss imply\sss that\sss
$K\fff'\dff \oplus\dff K\fff'$\sss
is\dss isomorphic\sss to\sss
$K\dff \oplus\dff K$\nnsp.\oss
Since $K$\sss is\dss finitely dimensional,\pss $K\fff'$\sss
is\dss finitely dimensional\sss also.\oss
Therefore\sss the assumptions of\trs Theorem\qss \ref{spectral-triplet-index}\qss
hold\dss for\sss $A_{\trf \Gamma}$ and\sss
$(\qff K\fff'\dff,\off \Gamma\fff'_{0}\dff,\off \Gamma\fff'_{1}\trf)$\sss
in\sss the roles of\dss $A$ and\sss
$(\qff K\dff,\off \Gamma_{\fff 0}\dff,\off \Gamma_{\fff 1}\trf)$\sss
respectively.\oss

Let\sss
$U\dff \colon\dff
K\dff \oplus\dff K
\qff \ttoo\qff
K\fff'\dff \oplus\dff K\fff'$\sss
be\sss the isomorphism\sss provided\dss by\sss the uniqueness property\sss
of\trs boundary\sss triplets,\oss and\dss let\sss
$\mathcal{B}_{\fff w}
\off =\off
U\trf(\trf \mathcal{C}_{\fff w}\trf)$\sss
for every\sss $w\qff \in\qff W$\nnsp.\oss
Then\sss $\mathcal{B}_{\fff w}\dff,\pff w\qff \in\qff W$\sss
is\dss a continuous family of\dss self-adjoint\sss relations on\sss $K\fff'$\dnsp.\oss
Also,\oss for every\sss $w\qff \in\qff W$\sss
the extension\sss $(\trf A_{\trf D}\trf)_{\fff \mathcal{B}_{\fff w}}$\sss of\dss $A_{\trf D}$
defined\dss by\sss the boundary\sss triplet\sss
$(\qff K\fff'\dff,\off \Gamma\fff'_{0}\dff,\off \Gamma\fff'_{1}\trf)$
and\sss $\mathcal{B}_{\fff w}$\sss
is\dss equal\dss to\sss the extension\sss
$A_{\dff w}\off =\off (\trf A_{\trf D}\trf)_{\dff \mathcal{C}_{\fff w}}$\sss
defined\dss by\sss the boundary\sss triplet\sss
$(\trf K\dff,\off \Gamma_{\fff 0}\dff,\off \Gamma_{\dff 1}\trf)$
and\sss $\mathcal{C}_{\fff w}$\nsp.\oss
Hence\sss the families\vspace{4.5pt}
\[
\quad
A_{\dff w}
\off =\off
(\trf A_{\trf D}\trf)_{\dff \mathcal{C}_{\fff w}}\dff,\off w\qff \in\qff W
\quad
\mbox{and}\quad\dff
(\trf A_{\trf D}\trf)_{\dff \mathcal{B}_{\fff w}}\dff,\off w\qff \in\qff W
\]

\vspace{-12pt}\vspace{4.5pt}
are equal,\oss and so are\sss their analytical\dss indices.\oss 
By\trs Theorem\qss \ref{spectral-triplet-index}\qss the family\sss
$(\trf A_{\trf D}\trf)_{\dff \mathcal{B}_{\fff w}}\dff,\off w\qff \in\qff W$\sss
is\trs Fredholm\dss and\sss its analytical\dss index of\trs equal\dss to\sss 
the analytical\dss index of\dss 
$\mathcal{B}_{\fff w}\dff,\off w\qff \in\qff W$\dnsp.\oss
Since\sss the family\sss
$\mathcal{B}_{\fff w}\dff,\pff w\qff \in\qff W$\sss
is\dss equal\sss to\sss the image of\dss the family\sss
$\mathcal{C}_{\fff w}\dff,\pff w\qff \in\qff W$
under\sss $U$\nnsp,\oss
the analytical\dss indices of\dss these\sss two families are equal.\oss
The\sss theorem\sss follows.\oss  \eproof

\newpage
\mysection{Realizations\qss of\pss boundary\qss symbols}{realizations}

\myuppar{General\dss boundary operators in\sss the self-adjoint\sss case.}
As in\dss Section\qss \ref{pdo},\oss
let\sss $P$ be an elliptic self-adjoint\sss operator from\sss
the\dss H\"{o}rmander class,\pss $B_{\trf Y}$\sss be a pseudo-differential\sss
operator of\dss order $0$\sss from $E\trf|\trf Y$\sss to some bundle $G$ over $Y$\dnsp,\oss 
and\sss $N$\sss be\sss the kernel-symbol\sss of\dss $B_{\trf Y}$\nsp.\oss
Suppose\sss that\dss $B\off =\off B_{\trf Y}\dff \circ\trf \gamma$\sss
satisfies\sss the\dss Shapiro-Lopatinskii\sss condition\dss for\dss $P$\dnsp,\oss
and\sss hence\sss $N$\sss satisfies\sss the\dss Shapiro-Lopatinskii\sss 
condition\dss for\sss the symbol\sss $\sigma$ of\dss $P$\dnsp.\oss
Let\vspace{1.5pt}
\[
\quad
\Pi\dff \colon\dff
H_{\trf 0}\dff(\trf Y\fff,\qff E\trf |\trf Y\trf)
\off \ttoo\off
H_{\trf 0}\dff(\trf Y\fff,\qff E\trf |\trf Y\trf)
\]

\vspace{-12pt}\vspace{1.5pt}
be\sss the orthogonal\dss projection onto\sss the kernel\sss of\dss $B_{\trf Y}$\nsp.\oss
By\sss the results of\trs Seeley\dss $\Pi$\sss is\dss a pseudo-differential\sss
operator of\dss order $0$\nnsp.\oss
See\qss \cite{se2},\oss Theorem\qss IV\halfff.6\qss and\dss the proof\dss of\trs
Theorem\qss IV\halfff.7,\oss pp.\qss 252{\fff}--{\fff}257.\oss
In\sss particular,\pss $\Pi$\sss induces a projection\vspace{3pt}
\[
\quad
\Pi_{\dff 1/2}\dff \colon\dff
H_{\dff 1/2}\dff(\trf Y\fff,\qff E\trf |\trf Y\trf)
\off \ttoo\off
H_{\dff 1/2}\dff(\trf Y\fff,\qff E\trf |\trf Y\trf)
\pff.
\]

\vspace{-12pt}\vspace{3pt}
As in\dss Section\qss \ref{pdo},\oss 
we can replace\sss $B_{\trf Y}$\sss by\sss $1\qff -\qff \Pi_{\dff 1/2}$\sss
and\sss $B$\sss by\sss 
$\Gamma\off =\off (\trf 1\qff -\qff \Pi_{\dff 1/2}\trf)\dff \circ\dff \gamma$\sss
without\sss affecting\sss the kernels and\sss the induced unbounded operators.\oss
But\sss assuming\sss that\sss $N$\sss is\dss self-adjoint\dss is\dss not\sss
sufficient\sss to ensure\sss that\sss 
$P_{\dff \Gamma}\off =\off P_{\dff B}$\sss
is\dss self-adjoint.\oss
Moreover,\oss the property of\dss being self-adjoint\dss is\dss not\sss
determined\sss by\sss the principal\sss symbol\sss of\dss $B_{\trf Y}$\nsp.\oss

In order\sss to discuss\sss this in more details,\oss
it\dss is\dss convenient\sss to change\sss the form of\dss boundary\sss conditions.\oss
Seeley\qss \cite{se2}\qss suggested\dss to\sss
consider arbitrary\sss pseudo-differential\sss operators $B$ of\dss order $0$ 
acting in\sss $E\trf |\trf Y$\sss such\sss that\sss $\image\dff B$\sss is\dss
closed\sss in\sss
$H_{\trf 0}\dff(\trf Y\fff,\qff E\trf |\trf Y\trf)$\sss
and an analogue of\dss the\dss Shapiro-Lopatinskii\sss condition\dss holds.\oss
Seeley\dss called\sss such operators $B$\dss 
\emph{well-posed}\pss for $P$\dnsp.\oss
See\qss \cite{se2},\oss Definition\qss VI.3\qss (p.\qss 289).\oss
For such\sss $B$\sss one can define\sss the boundary operator\sss $\Gamma$\sss as
$B\dff \circ\dff \gamma$\sss
and\sss then define\sss the unbounded operator\sss $P_{\dff \Gamma}$\sss exactly as before.\oss
One can always replace $B$\sss by a projection\sss $B\fff'$ without\sss
affecting\sss the kernel\sss $\kernel\fff \Gamma$ and\dss hence\sss $P_{\dff \Gamma}$\nsp.\oss
See\qss \cite{se2},\oss Lemma\qss VI.3\qss (p.\qss 289).\oss
By\sss this reason we will\sss assume\sss that\sss $B$\sss is\dss a projection.\oss
Let\sss $\Pi\off =\off 1\qff -\qff B$\nnsp.\oss

For general\sss elliptic operators\sss well-posed\sss operators define
a strictly\sss larger class of\dss boundary problems\dss
than\sss the classical\dss boundary operators.\oss
But\sss in\sss the self-adjoint\sss case\sss these classes are\sss the same.\oss
Cf.\oss the discussion of\dss boundary conditions in\sss the classical\dss form\sss
in\dss Section\qss \ref{symbols-conditions}.\oss

The results of\trs Section\qss \ref{abstract-index}\qss imply\sss that\sss 
$P_{\dff \Gamma}$\sss is\dss self-adjoint\dss if\dss 
$\Pi\off =\off 1\qff -\qff B$\sss is\dss a self-adjoint\sss projection and\sss
$\Sigma\trf(\trf \image\dff \Pi\trf)\off =\off \kernel\fff \Pi$\dss
(and an elliptic regularity\sss property\sss holds).\oss
It\dss is\dss only natural\dss to ask when\sss there exists a pseudo-differential\sss 
operator\sss $\Pi$\sss of\dss order $0$ with\sss such properties.\oss

\myuppar{Realization of\dss self-adjoint\sss kernel-symbols by self-adjoint\sss operators.}
If\dss $P$\sss is\dss elliptic,\oss then\sss $\Sigma$\sss is\dss 
an automorphism of\dss the bundle $E\trf|\trf Y$\nnsp.\oss
For\sss the rest\sss of\dss this section we will\sss assume\sss that,\oss moreover\halfff,\dss
$\Sigma$\sss is\dss an\sss isometry.\oss
The general\sss case can\sss be reduced\sss to\sss this one by\sss replacing\sss
$\Sigma$\sss
by\sss $\num{\Sigma}^{\dff -\dff 1}\dff \Sigma$\nnsp,\oss the subbundle $N$\sss by\sss
$\num{\Sigma}^{\dff -\dff 1/2}\dff (\trf N\trf)$\nnsp,\oss etc.\oss
Cf.\qss the first\sss homotopy\sss in\dss Section\qss \ref{boundary-algebra}.\oss

Let\sss $N$\sss be an elliptic self-adjoint\dss boundary condition\sss
for\sss the symbol\sss $\sigma$\sss in\sss the sense of\trs Section\qss \ref{symbols-conditions}.\oss
For\sss every\sss $u\qff \in\qff S\dff Y$\sss and\sss 
$y\off =\off \pi\trf(\dff u\trf)$\dss
let\sss
$\bm{\pi}_{\fff u}\dff \colon\dff
E_{\dff y}\qff \ttoo\qff E_{\dff y}$\sss
be\sss the orthogonal\dss projection with\sss the image $N_{\fff u}$\nsp.\oss
These projections define an endomorphism\sss\vspace{3pt}\vspace{-0.2pt}
\[
\quad
\bm{\pi}\dff \colon\dff 
\pi^{\fff *}\dff E\trf |\trf S\fff Y
\qff \ttoo\qff 
\pi^{\fff *}\dff E\trf |\trf S\fff Y
\qff.
\]

\vspace{-12pt}\vspace{3pt}\vspace{-0.2pt}
Let\sss us define a\qss \emph{realization}\pss of\dss the boundary condition $N$\sss
as a pseudo-differential\sss operator of\dss order $0$\sss
from\sss $E\trf |\trf Y$\sss to\sss $E\trf |\trf Y$\sss with\sss the
principal\sss symbol\sss $\bm{\pi}$\dss such\sss that\sss induced operator\vspace{3pt}\vspace{-0.2pt}
\[
\quad
\Pi\dff \colon\dff
H_{\trf 0}\dff(\trf Y\fff,\qff E\trf |\trf Y\trf)
\off \ttoo\off
H_{\trf 0}\dff(\trf Y\fff,\qff E\trf |\trf Y\trf)
\qff
\]

\vspace{-12pt}\vspace{3pt}\vspace{-0.2pt}
is\dss a self-adjoint\dss projection and\dss
$\Sigma\trf(\trf \image\dff \Pi\trf)
\off =\off
\kernel\dff \Pi$\nnsp.\oss
The\sss last\sss condition\dss is\dss the crucial\sss one.\oss
The\sss results of\trs Seeley\qss \cite{se2}\qss imply\sss that\sss 
there exists a pseudo-differential\sss operator $\Pi$ of\dss order $0$
with\sss the principal\sss symbol $\bm{\pi}$ such\sss that\sss
$\Pi$\sss is\dss a projection,\oss
i.e.\dss $\Pi\dff \circ\dff \Pi\off =\off \Pi$\nnsp.\oss
By applying\sss to\sss $\Pi^{\dff *}\dff \Pi$\sss 
an appropriate real\dss function
we will\dss get\sss
a self-adjoint\sss projection\sss with\sss the same symbol.

The boundary condition $N$ also defines for every\sss $u\qff \in\qff S\fff Y$\sss
an\sss isometry
$\varphi_{\fff u}\dff \colon\dff
E_{\dff y}^{\dff +}
\qff \ttoo\qff
E_{\dff y}^{\dff -}$\sss
having\sss $N_{\fff u}$ as\sss its graph.\oss
In\sss turn,\oss these isometries define
an\sss isometric\sss isomorphism of\dss bundles\vspace{3pt}\vspace{-0.2pt}
\[
\quad
\bm{\varphi}\dff \colon\dff 
\pi^{\fff *}\dff E^{\dff +}_{\trf Y}
\qff \ttoo\qff 
\pi^{\fff *}\dff E^{\dff -}_{\trf Y}
\off
\]

\vspace{-12pt}\vspace{3pt}\vspace{-0.2pt}
already used\sss in\dss Section\qss \ref{symbols-conditions}.\oss
An\qss \emph{isometric\dss realization}\pss of\dss $\bm{\varphi}$\sss
is\dss a pseudo-differential\sss operator of\dss order $0$\sss
from\dss $E^{\dff +}_{\trf Y}$\dss to\dss $E^{\dff -}_{\trf Y}$\dss with\sss the
principal\sss symbol\sss $\bm{\varphi}$\dss such\sss that\sss induced operator\vspace{3pt}\vspace{-0.2pt}
\[
\quad
\Phi\dff \colon\dff
H_{\trf 0}\dff(\trf Y\fff,\qff E^{\dff +}_{\trf Y}\trf)
\off \ttoo\off
H_{\trf 0}\dff(\trf Y\fff,\qff E^{\dff -}_{\trf Y}\trf)
\qff
\]

\vspace{-12pt}\vspace{3pt}\vspace{-0.2pt}
is\dss an\sss isometry of\trs Hilbert\dss spaces.\oss

The isometric realizations of\dss $\bm{\varphi}$\sss
are closely related\dss to\sss the realizations of\dss $N$\nnsp.\oss
Let\vspace{3pt}
\[
\quad
H^{\dff \partial}_{\dff +}
\off =\off
H_{\trf 0}\dff\left(\trf Y\fff,\qff E^{\dff +}_{\trf Y}\trf\right)
\qff,\quad
H^{\dff \partial}_{\dff -}
\off =\off
H_{\trf 0}\dff\left(\trf Y\fff,\qff E^{\dff -}_{\trf Y}\trf\right)
\qff,\quad
\mbox{and\dss let}
\]

\vspace{-33pt}
\[
\quad
\pr_{\dff +}\dff \colon\dff
H^{\dff \partial}_{\dff +}\dff \oplus\dff H^{\dff \partial}_{\dff -}
\qff \ttoo\qff
H^{\dff \partial}_{\dff +}
\off,\quad
\pr_{\dff -}\dff \colon\dff
H^{\dff \partial}_{\dff +}\dff \oplus\dff H^{\dff \partial}_{\dff -}
\qff \ttoo\qff
H^{\dff \partial}_{\dff -}
\]

\vspace{-12pt}\vspace{3pt}
be\sss the canonical\dss projections.\oss
Then\dss
$H^{\dff \partial}
\off =\off
H^{\dff \partial}_{\dff +}\dff \oplus\dff H^{\dff \partial}_{\dff -}$\dss
and $\Sigma$ acts as\sss the identity on\sss $H^{\dff \partial}_{\dff +}$\sss
and\dss the minus identity on\sss $H^{\dff \partial}_{\dff -}$.\oss
Suppose\sss that\sss
$\Pi\dff \colon\dff
H^{\dff \partial}\qff \ttoo\qff H^{\dff \partial}$\sss
is\dss a self-adjoint\dss projection such\sss that\sss
$\Sigma\trf(\trf \image\dff \Pi\trf)
\off =\off
\kernel\dff \Pi$\nnsp.\oss
We claim\sss that\sss  
$\image\dff \Pi$\sss
is\dss the graph of\dss an\sss isometry\sss
$\Phi\dff \colon\dff
H^{\dff \partial}_{\dff +}
\qff \ttoo\qff
H^{\dff \partial}_{\dff -}$.\oss

Since\sss
$\Sigma\trf(\qff \image\dff \Pi  \qff)
\off =\off
\kernel\dff \Pi$\nnsp,\oss
the image\sss $\image\dff \Pi$\sss intersects each of\dss the spaces\sss
$H^{\dff \partial}_{\dff +}$\sss and\sss $H^{\dff \partial}_{\dff -}$\sss
by zero.\oss
Suppose\sss that\sss  $v\qff \in\qff H^{\dff \partial}_{\dff +}$\sss
is\dss orthogonal\dss to\sss $\pr_{\dff +}\trf(\trf \image\dff \Pi \trf)$\nnsp.\oss
Clearly,\oss then\sss $v$\sss is\dss orthogonal\dss to\sss $\image\dff \Pi$\nnsp.\oss
Since\sss $\Sigma\trf(\trf v\trf)\off =\off v$\nnsp,\oss
this implies\sss that\sss $v$\sss is\dss orthogonal\sss to\sss $\kernel\dff \Pi$\sss
also,\oss and\sss hence\sss $v\off =\off 0$\nnsp.\oss
It\sss follows\sss that\sss $\pr_{\dff +}$\sss induces a surjective map\sss
$\image\dff \Pi\qff \ttoo\qff H^{\dff \partial}_{\dff +}$.\oss
Similarly,\pss $\pr_{\dff -}$\sss induces a surjective map\sss
$\image\dff \Pi\qff \ttoo\qff H^{\dff \partial}_{\dff -}$.\oss
It\sss follows\sss that\sss
$\image\dff \Pi$\sss
is\dss equal\sss to\sss the graph of\dss a\sss linear map\sss
$\Phi\dff \colon\dff
H^{\dff \partial}_{\dff +}
\qff \ttoo\qff
H^{\dff \partial}_{\dff -}$.\oss
Arguing\sss as in\dss Section\qss \ref{boundary-algebra}\qss
we see\sss that\sss $\Phi$\sss is\dss an\sss isometry\qss
(see\sss the discussion of\trs lagrangian subspaces when $\sigma$\sss is\dss unitary\fff).\oss
This proves our claim.\oss

Let\sss us express\sss $\Pi$\sss in\sss terms of\dss $\Phi$\sss and vice versa.\oss
The orthogonal\sss complement\sss to\sss the graph of\dss $\Phi$\sss
is\dss equal\dss to\sss the graph of\dss $-\qff \Phi$\nnsp.\oss
Since\sss $\Pi$\sss is\dss a self-adjoint\sss projection,\oss it\dss 
is\dss also equal\dss to\sss $\image\dff (\trf 1\qff -\qff \Pi\trf)$\nnsp.\oss
It\sss follows\sss that\sss every\sss
$(\trf x\fff,\qff y\trf)
\qff \in\qff
H^{\dff \partial}_{\dff +}\dff \oplus\dff H^{\dff \partial}_{\dff -}$\sss
admits a unique presentation\vspace{3pt}
\[
\quad
(\trf x\fff,\qff y\trf)
\off =\off
(\trf a\fff,\qff \Phi\dff(\trf a\trf)\trf)
\qff +\qff
(\trf b\fff,\qff -\qff \Phi\dff(\trf b\trf)\trf)
\qff,
\]

\vspace{-12pt}\vspace{3pt}
where\sss $a\fff,\qff b\qff \in\qff H^{\dff \partial}_{\dff +}$.\oss
An easy calculation shows\sss that\vspace{1.5pt}
\[
\quad
a
\off =\off
\frac{x\qff +\qff \Phi^{\dff -\dff 1}\dff(\trf y\trf)}{2}
\qff,\quad
b
\off =\off
\frac{x\qff -\qff \Phi^{\dff -\dff 1}\dff(\trf y\trf)}{2}
\off.
\]

\vspace{-12pt}\vspace{0pt}
This\sss leads\sss to\sss the formulas\vspace{3pt}
\[
\quad
\Pi\trf(\trf x\fff,\qff y\trf)
\off =\off
\left(\qff
\frac{x\qff +\qff \Phi^{\dff -\dff 1}\dff(\trf y\trf)}{2}\qff,\off
\frac{\Phi\dff(\trf x\trf)\qff +\qff y}{2}
\qff\right)
\quad
\mbox{and}\quad
\]

\vspace{-33pt}\vspace{-15pt}
\begin{equation}
\label{pi-phi}
\end{equation}

\vspace{-33pt}\vspace{-15pt}
\[
\quad
\Phi\dff(\trf x\trf)
\off =\off
\pr_{\dff -}\qff
\bigl(\trf
2\trf \Pi\trf (\trf x\fff,\qff 0\trf)
\trf\bigr)
\]

\vspace{-12pt}\vspace{1.5pt}
relating\sss $\Pi$\sss and\sss $\Phi$\nnsp.\oss

Conversely,\oss suppose\sss that\sss 
$\Phi\dff \colon\dff
H^{\dff \partial}_{\dff +}
\qff \ttoo\qff
H^{\dff \partial}_{\dff -}$\sss is\dss an\sss isometry.\oss
Let\sss us define\sss $\Pi$\sss by\sss the above formula.\oss
Direct\sss calculations show\sss that\sss
$\Pi\dff \circ\dff \Pi\off =\off \Pi$\nnsp,\oss
i.e.\dss $\Pi$\sss is\dss a projection,\oss
and\vspace{1.5pt}
\[
\quad
(\trf 1\qff -\qff \Pi \trf)\trf(\trf x\fff,\qff y\trf)
\off =\off
\left(\qff
\frac{x\qff -\qff \Phi^{\dff -\dff 1}\dff(\trf y\trf)}{2}\qff,\off
\frac{y\qff -\qff \Phi\dff(\trf x\trf)}{2}
\qff\right)
\qff
\]

\vspace{-33pt}
\[
\quad
\phantom{(\trf 1\qff -\qff \Pi \trf)\trf(\trf x\fff,\qff y\trf)
\off }
=\off
\Sigma\qff
\left(\qff
\frac{x\qff -\qff \Phi^{\dff -\dff 1}\dff(\trf y\trf)}{2}\qff,\off
\frac{\Phi\dff(\trf x\trf)\qff -\qff y}{2}
\qff\right)
\off =\off
\Sigma\qff
\bigl(\qff
\Pi\trf(\trf x\fff,\qff -\qff y\trf)
\qff\bigr)
\qff.
\]

\vspace{-12pt}\vspace{1.5pt}
Therefore\sss
$\Sigma\trf(\qff \image\dff \Pi  \qff)
\off =\off
\image\dff (\trf 1\qff -\qff \Pi \trf)
\off =\off
\kernel\dff \Pi$.\oss
Another direct\sss calculation,\oss using\sss the assumption\sss
that\sss $\Phi$\sss is\dss an\sss isometry,\oss
shows\sss that\sss
$\image\dff \Pi$\sss is\dss orthogonal\sss to\sss 
$\image\dff (\trf 1\qff -\qff \Pi \trf)$\nnsp.\oss
Since\sss $\Pi$\sss is\dss a projection,\oss
this implies\sss that\sss $\Pi$\sss is\dss a self-adjoint\sss projection.\oss

\mypar{Theorem.}{two-realizations}
\emph{Suppose\sss that\dss $\Pi$\sss and\dss $\Phi$\sss are related as\sss in\pss \textup{(\ref{pi-phi})}.\oss
Then\sss $\Pi$\sss is\dss a realization of\dss the boundary condition\dss $N$\dss if\dss and\dss only\trs if\qss
$\Phi$\sss is\dss an\sss isometric realization of\qss $\bm{\varphi}$\nnsp.\oss}

\proof
The formulas\qss (\ref{pi-phi})\qss imply\sss that\sss $\Pi$\sss
is\dss a pseudo-differential\sss operator of\dss order $0$\sss
if\dss and\dss only\trs if\dss $\Phi$\sss is.\oss
Moreover,\oss if\dss $\Pi$\sss and\sss $\Phi$\sss are
pseudo-differential\sss operators of\dss order $0$\nnsp,\oss
then\sss their\sss principal\sss symbols are related\sss in\sss the same way as\sss
the operators\sss themselves.\oss
It\sss follows\sss that\sss the principal\sss symbol\sss of\dss $\Pi$\sss
is\dss equal\dss to\sss $\bm{\pi}$\sss
if\dss and\dss only\trs if\dss the principal\sss symbol\sss of\dss $\Phi$\sss
is\dss equal\dss to\sss $\bm{\varphi}$\nnsp.\oss
This proves\sss the\sss theorem.\oss  \eproof

\mypar{Theorem.}{sa-index-obstruction}
\emph{There exists an\dss isometric realization of\dss $\bm{\varphi}$\trs
if\trs and\dss only\trs if\trs the analytical\dss index of\dss
pseudo-differential\sss operators of\dss order\sss $0$\sss
with\sss the principal\sss symbol\sss $\bm{\varphi}$\dss is\dss zero.\oss}\vspace{-0.5pt}

\proof
The\qss ``only\trs if''\qss part\dss is\dss trivial.\oss
Let\sss $\Psi_{\fff 0}$\sss be pseudo-differential\sss operator of\dss order $0$ 
with\sss the principal\sss symbol\sss $\bm{\varphi}$\nnsp.\oss
If\dss the analytical\dss index of\dss $\Psi_{\fff 0}$\sss is\sss $0$\nnsp,\oss
then 
$\Psi_1\off =\off \Psi_{\fff 0}\qff +\qff k$\sss 
is\dss an\sss isomorphism for some operator $k$ of\dss finite rank.\oss
The principal\sss symbol\sss of\dss $\Psi_1$\sss 
is\dss equal\dss to\sss that\sss of\dss $\Psi$\nnsp,\oss
i.e.\qss to\sss $\bm{\varphi}$\nnsp.\oss
Let\sss 
$\Phi\off =\off \Psi_1\qff \num{\Psi_1}^{\dff -\dff 1}$\dnsp,\oss
where\sss
$\num{\Psi_1}\off =\off (\trf\Psi_1^{\dff *}\dff \Psi_1\trf)^{\dff 1/2}$\dnsp.\oss
Then\sss $\Phi$\sss is\dss an\sss isometry.\oss
The results of\trs Seeley\qss \cite{se2}\qss imply\sss that\sss
$\num{\Psi_1}$\sss is\dss a pseudo-differential\sss operator\sss
of\dss order $0$\sss with\sss the principal\sss symbol\sss $\num{\bm{\varphi}}$\nnsp.\oss
Since\sss $\bm{\varphi}$\sss is\dss an\sss isometry,\pss
$\num{\bm{\varphi}}$\sss is\dss the identity\sss map
and\dss hence\sss the principal\sss symbol\sss of\dss
$\Phi$\sss is\dss equal\dss to\sss that\sss of\dss $\Psi_1$\sss
and\dss hence\sss to\sss $\bm{\varphi}$\nnsp.\oss
Since\sss $\Phi$\sss is\dss an\sss isometry,\oss
this proves\sss the\qss ``if''\qss part.\oss  \eproof\vspace{-1.5pt}

\myuppar{Realization of\dss self-adjoint\sss kernel-symbols\sss in\sss families.}
Suppose now,\oss as at\sss the end of\trs Sections\qss \ref{symbols-conditions},\oss
that\sss our symbols and operators depend on a parameter\sss $z\qff \in\qff Z$\nnsp,\oss
where $Z$\sss is\dss a\sss topological\sss space.\oss
As in\dss  Sections\qss \ref{symbols-conditions},\oss
we allow\sss also\sss the manifold\sss $X\off =\off X\trf(\trf z\trf)$\sss to vary
among\sss the fibers of\dss a bundle.\oss
Let\sss $Y\trf(\trf z\trf)\off =\off \partial\trf X\trf(\trf z\trf)$\nnsp.\oss
We will\sss assume\sss that\sss $Z$\sss is\dss a compactly\sss generated and\sss paracompact\sss space.\oss
For example,\oss every\sss metric space has\sss these properties.\oss 
Like above,\oss we will\sss consider only\sss the case when\sss the bundle maps\sss $\Sigma\trf(\trf z\trf)$\sss
are\sss isometries.\oss\vspace{-0.5pt}

Let\sss $N\trf(\trf z\trf)\fff,\qff z\qff \in\qff Z$\sss be a continuous family of\dss
self-adjoint\sss elliptic boundary conditions.\oss
Then\sss the families of\dss the corresponding\dss bundle maps\sss
$\bm{\pi}\trf(\trf z\trf)\fff,\qff z\qff \in\qff Z$\sss
and\sss
$\bm{\varphi}\trf(\trf z\trf)\fff,\qff z\qff \in\qff Z$\sss
are defined.\oss
A\qss \emph{realization}\pss of\dss the family\sss $N\trf(\trf z\trf)\fff,\qff z\qff \in\qff Z$\sss
is\dss defined as a continuous family\sss $\Pi\trf(\trf z\trf)\fff,\qff z\qff \in\qff Z$\sss
of\dss pseudo-differential\sss operators\sss of\dss order $0$\sss
such\sss that\sss $\Pi\trf(\trf z\trf)$\sss is\dss a realization of\dss $N\trf(\trf z\trf)$\sss
for every\sss $z\qff \in\qff Z$\nnsp.\oss
The\qss \emph{isometric\dss realizations}\pss of\dss 
$\bm{\varphi}\trf(\trf z\trf)\fff,\qff z\qff \in\qff Z$\sss
are defined\sss in\sss the same manner.\oss\vspace{-0.5pt}

\mypar{Theorem.}{two-realizations-families}
\emph{Suppose\sss that\dss $\Pi\off =\off \Pi\trf(\trf z\trf)$ and\trs 
$\Phi\off =\off \Phi\trf(\trf z\trf)$\sss are related as\sss in\pss \textup{(\ref{pi-phi})}\qss
for every\sss $z\qff \in\qff Z$\nnsp.\oss
Then\qss $\Pi\trf(\trf z\trf)\fff,\qff z\qff \in\qff Z$\qss is\dss a realization of\oss  
$N\trf(\trf z\trf)\fff,\qff z\qff \in\qff Z$\pss if\dss and\dss only\trs if\pss
$\Phi\trf(\trf z\trf)\fff,\qff z\qff \in\qff Z$\qss 
is\dss an\sss isometric realization of\qss 
$\bm{\varphi}\trf(\trf z\trf)\fff,\qff z\qff \in\qff Z$\nnsp.\oss}\vspace{-1.5pt}

\proof
It\dss is\dss sufficient\sss to repeat\sss the proof\dss of\qss Theorem\qss \ref{two-realizations}\qss
with\sss parameters added.\oss  \eproof\vspace{-0.5pt}

\mypar{Corollary.}{existence-realizations-families}
\emph{There exists a realization of\dss the family of\qss
boundary conditions\dss 
$N\trf(\trf z\trf)\fff,\qff z\qff \in\qff Z$\dss 
if\dss and\dss only\trs if\dss there exists
an\sss isometric realization of\qss 
$\bm{\varphi}\trf(\trf z\trf)\fff,\qff z\qff \in\qff Z$\nnsp.\oss}  \eproof\vspace{-0.5pt}

\myuppar{Fredholm\dss families.}
In order\sss to prove an analogue of\trs Theorem\qss \ref{sa-index-obstruction}\qss for
families,\oss we need\sss a general\dss result\sss 
about\sss families of\dss Fredholm\sss operators with analytical\sss index zero,\oss
namely,\oss Theorem\qss \ref{unitary-sections}\qss below.\oss
It\sss depends on\qss \cite{i1},\oss \cite{i2}.\oss
Let\sss $\mathbb{H}$\sss be a\dss Hilbert\dss bundle with\sss
the base\sss $Z$\nnsp,\oss considered as a family\sss
$H_{\dff z}\dff,\qff z\qff \in\qff Z$\sss of\dss Hilbert\dss spaces.\oss
Let\sss $\mathbb{A}$\sss be a\sss family of\qss Fredholm\dss operators\sss
$A_{\dff z}\dff \colon\dff 
H_{\dff z}\qff \ttoo\qff H_{\dff z}\dff,\off 
z\qff \in\qff Z$\nnsp.\oss 
We will\sss assume\sss that\sss $\mathbb{A}$\sss is\dss continuous in an appropriate sense,\oss
for example,\oss that\sss $\mathbb{A}$\sss is\dss
fully\trs Fredholm\dss in\sss
the sense of\oss \textup{\cite{i3}}.\oss
This\dss is\dss the case\sss if\dss $\mathbb{A}$\sss is\dss 
a continuous family of\dss elliptic pseudo-differential\sss
operators by\sss the results of\trs Atiyah\dss and\dss Singer\qss \cite{as4}.\oss
The continuity\sss in\sss the following\sss theorem\sss can be understood\sss
in\sss the same sense.\oss

\mypar{Theorem.}{unitary-sections}
\emph{For every $z\qff \in\qff Z$\dss let\dss
$A_{\dff z}\off =\off U_{\fff z}\qff \num{A_{\dff z}}$\sss
be\sss the polar decomposition of\dss $A_{\dff z}$\nsp.\oss
If\qss the analytical\sss index of\qss the family\qss $\mathbb{A}$\sss is\trs zero,\oss
then\sss there exists a continuous family\sss of\dss isometries
$U\fff'_{\fff z}\dff \colon\dff 
H_{\dff z}\qff \ttoo\qff H_{\dff z}\dff,\off 
z\qff \in\qff Z$\sss
such\sss that\trs
$U\fff'_{\fff z}\qff -\qff U_{\fff z}$\dss
is\dss a\dss finite rank operator\sss for every\sss $z\qff \in\qff Z$\nnsp.\oss}

\proof
Note\sss that,\oss in\sss general,\oss the family\sss
$U_{\fff z}\dff,\qff z\qff \in\qff Z$\sss is\dss not\sss continuous in any\sss
reasonable sense because\sss the dimension of\dss 
the kernels\sss $\kernel\dff U_{\fff z}$\sss may\sss jump.\oss
Suppose\sss that\sss the analytical\sss index of\dss $\mathbb{A}$\sss is\dss zero.\oss
Then\sss for every\sss $z\qff \in\qff Z$\sss
the index of\dss $A_{\dff z}$\sss is\dss zero
and\sss $U_{\fff z}$\sss is\dss a partial\dss isometry\sss with\sss
the same kernel\sss and cokernel\sss as $A_{\dff z}$\nsp.\oss
Since\sss the index of\dss $A_{\dff z}$\sss is\dss zero,\oss
the dimensions of\dss the kernel\sss and cokernel\sss are equal\sss
and\sss hence\sss there exist\dss isometries\sss
$H_{\dff z}\qff \ttoo\qff H_{\dff z}$\sss
equal\sss to\sss $U_{\fff z}$\sss on\sss the orthogonal\sss complement\sss
$H_{\dff z}\dff \ominus\dff \kernel\dff U_{\fff z}$\sss
of\dss the kernel\sss $\kernel\dff U_{\fff z}\off =\off \kernel\dff A_{\dff z}$\nsp.\oss

Let\sss $U^{\dff \fin}\dff(\trf z\trf)$\sss
be\sss the set\sss of\dss isometries\sss 
$H_{\dff z}\qff \ttoo\qff H_{\dff z}$\sss
equal\sss to\sss $U_{\fff z}$\sss on a closed subspace of\dss finite codimension\sss
in\sss $H_{\dff z}\dff \ominus\dff \kernel\dff U_{\fff z}$\nsp,\oss
and\dss let\sss us equip\sss $U^{\dff \fin}\dff(\trf z\trf)$\sss with\sss the norm\sss topology.\oss
Since\sss $\mathbb{A}$\sss is\dss a\dss Fredholm\dss family,\oss the family of\dss spaces\sss 
$U^{\dff \fin}\dff(\trf z\trf)\dff,\off 
z\qff \in\qff Z$\sss
forms a\sss locally\sss trivial\dss bundle over $Z$\sss
having\sss $U^{\dff \fin}\dff(\trf z\trf)$\sss
as\sss the fiber over\sss $z\qff \in\qff Z$\nnsp.\oss
Let\sss us denote\sss this bundle by\vspace{1.5pt}
\[
\quad
\bm{\pi}\ffin\dff(\trf \mathbb{A}\trf)\dff \colon\dff
U^{\dff \fin}\dff(\trf \mathbb{A}\trf)
\qff \ttoo\qff
Z
\qff.
\]

\vspace{-12pt}\vspace{1.5pt}
This\dss is\dss an analogue of\dss the\dss Grassmannian\dss bundle\sss
$\gr\trf(\trf \mathbb{A}\trf)$\sss defined\sss in\qss \cite{i2}\qss for strictly\trs Fredholm\dss
families of\dss self-adjoint\trs Fredholm\dss operators.\oss
There\dss is\dss also a\qss ``universal''\qss analogue\vspace{1.5pt}
\[
\quad
\bm{\pi}\ffin\dff \colon\dff
\mathbf{U}\ffin
\qff \ttoo\qff 
\num{\mathcal{P}{\nsp}\mathcal{S}_{\qff 0}}
\]

\vspace{-12pt}\vspace{1.5pt}
of\dss $\bm{\pi}\ffin\dff(\trf \mathbb{A}\trf)$\nnsp,\oss
where\sss $\num{\mathcal{P}{\nsp}\mathcal{S}_{\qff 0}}$\sss is\dss
a classifying\sss space.\oss
See\qss \cite{i1},\oss Section\qss 15,\oss
where\sss $\mathbf{U}\ffin$\sss is\dss denoted\sss by\sss $\mathbf{U}$\dnsp.\oss
The bundle\sss $\bm{\pi}\ffin\dff(\trf \mathbb{A}\trf)$\sss
is\dss induced\sss from\sss the bundle $\bm{\pi}\ffin$\sss
by a\qss \emph{polarized\dss index\dss map}\dss
$Z\qff \ttoo\qff \num{\mathcal{P}{\nsp}\mathcal{S}_{\qff 0}}$\nnsp.\oss
This\dss is\dss an analogue of\trs Theorem\qss 5.1\qss in\qss \cite{i2}\qss
with a similar\sss proof.\oss

Let\sss us use\sss the definition of\dss the analytical\dss index of\dss families
from\qss \cite{i2},\oss Section\qss 3.\oss
It\dss is\dss equivalent\sss to\sss the classical\sss one when\sss the\sss latter applies.\oss
See\qss \cite{i2},\oss Section\qss 7.\oss
With\sss this definition,\oss the analytical\dss index\dss is\dss zero\sss
if\dss and only\sss if\dss the polarized\dss index\sss map\dss is\dss homotopic\sss
to a constant\sss map\qss 
(compare\qss \cite{i2},\oss the proof\dss of\trs Theorem\qss 5.2).\oss 
Since\dss the analytical\dss index of\dss $\mathbb{A}$\sss is\dss zero,\oss
this implies\sss that\sss $\bm{\pi}\ffin\dff(\trf \mathbb{A}\trf)$\sss is\dss a\sss trivial\dss bundle\sss  
and\sss hence admits a continuous section.\oss
If\dss $z\off \longmapsto\off U\fff'_z$\sss is\dss a section,\oss
then\sss $U\fff'_z$\sss is\dss an\sss isometry\sss 
$H_{\dff z}\qff \ttoo\qff H_{\dff z}$\sss
equal\sss to\sss $U_{\fff z}$\sss on a subspace of\dss finite codimension.\oss
It\sss follows\sss that\sss the rank\sss of\dss
$U\fff'_{\fff z}\qff -\qff U_{\fff z}$\sss
is\dss finite for every\sss $z\qff \in\qff Z$\nnsp.\oss  \eproof

\myuppar{Remark.}
If\dss
the analytical\dss index of\dss $\mathbb{A}$\sss is\dss zero,\oss
then\sss there exists a continuous family of\dss isomorphisms\sss
$A_{\dff z}'\dff \colon\dff 
H_{\dff z}\qff \ttoo\qff H_{\dff z}\dff,\off 
z\qff \in\qff Z$\sss such\sss that\sss
$A_{\dff z}'\qff -\qff A_{\dff z}$\sss has finite rank\sss for every\sss
$z\qff \in\qff Z$\nnsp.\oss

\mypar{Theorem.}{sa-index-obstruction-families}
\emph{The family\sss 
$\bm{\varphi}\trf(\trf z\trf)\fff,\qff z\qff \in\qff Z$\trs
admits an\dss isomeric realization\dss
if\trs and only\trs if\trs its analytical\dss index\dss is\dss zero.\oss}\vspace{-6pt}

\proof
The\qss ``only\dss if''\qss part\dss is\dss trivial.\oss
Let\sss $\Psi \dff(\trf z\trf)\fff,\qff z\qff \in\qff Z$\sss
be some family of\dss pseudo-differential\sss operators of\dss order $0$ 
with\sss the principal\sss symbols\sss 
$\bm{\varphi}\trf(\trf z\trf)\fff,\qff z\qff \in\qff Z$\nnsp,\oss
and suppose\sss that\sss the analytical\dss index of\dss this family\dss is\dss $0$\nnsp.\oss
Then\sss the index of\dss the operator\sss
$\Psi \dff(\trf z\trf)$\sss is\dss zero for every\sss $z\qff \in\qff Z$\nnsp.
Let\sss 
$\Psi \dff(\trf z\trf)
\off =\off 
U\dff(\trf z\trf)\qff \num{\Psi \dff(\trf z\trf)}$\sss
be\sss the polar decomposition of\dss $\Psi \dff(\trf z\trf)$\nsp.\oss
We claim\sss that\sss $U\dff(\trf z\trf)$\sss is\dss a pseudo-differential\sss
operator of\dss order $0$ with\sss the\qss (principal\fff)\qss symbol\sss $\bm{\varphi}\trf(\trf z\trf)$\nnsp.\oss
Since\sss the index of\dss the operator\sss $\Psi \dff(\trf z\trf)$\sss is\dss zero,\oss
there exists an\sss operator $k$ of\dss finite rank\sss such\sss
that\sss $\Psi \dff(\trf z\trf)\qff +\qff k$\sss is\dss invertible.\oss
By\sss the results of\trs Seeley\qss \cite{se2},\oss
the operator\vspace{1.5pt}\vspace{0.625pt}
\[
\quad
\num{\Psi \dff(\trf z\trf)\qff +\qff k}
\off =\off
\left(\qff
\left(\qff 
\Psi \dff(\trf z\trf)^{\fff *}\qff +\qff k^{\dff *}
\qff\right)\dff
\left(\qff 
\Psi \dff(\trf z\trf)\qff +\qff k\vphantom{k^{\dff *}}
\qff\right)
\qff\right)^{\dff 1/2}
\]

\vspace{-12pt}\vspace{1.5pt}\vspace{0.625pt}
is\dss a pseudo-differential\sss operator of\dss order $0$ with\sss 
the symbol\sss $\num{\bm{\varphi}\trf(\trf z\trf)}$\nnsp.\oss
Since\sss $\bm{\varphi}$\sss is\dss unitary,\oss this symbol\dss
is\dss equal\dss to\sss the identity.\oss
Clearly,\pss 
$\num{\Psi \dff(\trf z\trf)\qff +\qff k}$\sss
is\dss equal\sss to\sss $\num{\Psi \trf(\trf z\trf)}$\sss
on\sss the intersection of\dss the kernels\sss 
$\kernel k$\sss
and\sss 
$\kernel k^{\dff *}\dff \circ\qff \Psi \dff(\trf z\trf)$\nnsp,\oss
a subspace of\dss finite codimension.\oss
It\sss follows\sss that\sss $U\dff(\trf z\trf)$\sss
is\dss equal\sss to\sss 
$\Psi \dff(\trf z\trf)\qff
\num{\Psi \dff(\trf z\trf)\qff +\qff k}^{\dff -\dff 1}$\sss
on a subspace of\dss finite codimension.\oss
In\sss turn,\oss this implies\sss that\sss $U\dff(\trf z\trf)$\sss
is\dss a pseudo-differential\sss operator of\dss order $0$ with\sss
the same symbol\sss as\sss
$\Psi \dff(\trf z\trf)\qff
\num{\Psi \dff(\trf z\trf)\qff +\qff k}^{\dff -\dff 1}$\dnsp.\oss
Since\sss the symbol\sss of\dss
$\num{\Psi \dff(\trf z\trf)\qff +\qff k}^{\dff -\dff 1}$\sss
is\dss equal\dss to\sss the identity,\oss
it\sss follows\sss that\sss the symbol\sss of\dss
$U\dff(\trf z\trf)$\sss is\dss equal\dss to\sss $\bm{\varphi}\trf(\trf z\trf)$\nnsp.\oss
This proves our claim.\oss

By\trs Theorem\qss \ref{unitary-sections}\qss 
there exists a continuous family of\dss isometries\sss
$U\fff'\dff(\trf z\trf)\fff,\off z\qff \in\qff Z$\sss
such\sss that\sss the differences\sss
$U\fff'\dff(\trf z\trf)\qff -\qff U\dff(\trf z\trf)$\sss
have finite rank.\oss
This\sss implies\sss that\sss $U\fff'\dff(\trf z\trf)$\sss
is\dss a pseudo-differential\sss operator of\dss order $0$ with\sss 
the same symbol\sss as\sss $U\dff(\trf z\trf)$\nnsp,\oss
i.e.\qss with\sss the symbol\sss $\bm{\varphi}\trf(\trf z\trf)$\nnsp.\oss
Therefore\sss $U\fff'\dff(\trf z\trf)\fff,\off z\qff \in\qff Z$\sss
is\dss an\sss isometric realization of\dss
$\bm{\varphi}\trf(\trf z\trf)\fff,\qff z\qff \in\qff Z$\nnsp.\oss
This proves\sss the\qss ``if''\qss part.\oss  \eproof

\myuppar{The\sss two obstructions.}
The analytical\dss index of\dss a family\sss 
$\Phi\trf(\trf z\trf)\fff,\qff z\qff \in\qff Z$\sss 
of\dss pseudo-differential\sss operators of\dss order $0$
with\sss the symbols\sss $\bm{\varphi}\trf(\trf z\trf)\fff,\qff z\qff \in\qff Z$\sss
depends only on\sss $\bm{\varphi}\trf(\trf z\trf)\fff,\qff z\qff \in\qff Z$\sss
and\dss is\dss called\dss the\qss \emph{analytical\dss index}\pss of\dss the\sss latter.\oss
Let\sss $P\dff(\trf z\trf)\fff,\qff z\qff \in\qff Z$\sss
be a continuous family of\dss self-adjoint\dss pseudo-differential\sss operators of\dss
order $1$\sss belonging\dss to\sss the\dss H\"{o}rmander's\dss class,\oss
and\sss let\sss $N\trf(\trf z\trf)\fff,\qff z\qff \in\qff Z$\sss be
a family of\dss self-adjoint\dss boundary\sss conditions for\sss the symbols
$\sigma\trf(\trf z\trf)$ of\dss these operators.\oss 
The above results show\sss that\sss the analytical\dss index of\dss
the corresponding\dss family\sss
$\bm{\varphi}\trf(\trf z\trf)\fff,\qff z\qff \in\qff Z$\sss
is\dss an obstruction,\oss and\sss the only obstruction\sss
to\sss the existence of\dss realizations of\dss the family\sss
$N\trf(\trf z\trf)\fff,\qff z\qff \in\qff Z$\nnsp.\oss
As in\dss Section\qss \ref{symbols-conditions},\oss
let\sss $V\qff \ttoo\qff Z$\sss be\sss the bundle having\sss
$Y\trf(\trf z\trf)\off =\off \partial\trf X\trf(\trf z\trf)$\sss
as\sss the fiber over\sss $z\qff \in\qff Z$\nnsp.\oss
The families of\dss bundles\sss\vspace{1.5pt}\vspace{0.625pt}
\[
\quad
B\dff Y\trf(\trf z\trf)\dff,\quad
S\fff Y\trf(\trf z\trf)\dff,\quad
E\trf(\trf z\trf)^{\dff +}_{\trf Y\trf(\trf z\trf)}\dff,\quad
\mbox{and}\quad
E\trf(\trf z\trf)^{\dff -}_{\trf Y\trf(\trf z\trf)}
\]

\vspace{-12pt}\vspace{1.5pt}\vspace{0.625pt}
lead\sss to\sss the corresponding\sss bundles over\sss $V$\dnsp,\oss
which we will\sss denote respectively\sss by\vspace{1.5pt}\vspace{0.625pt}
\[
\quad
B\dff V\dff,\quad
S\fff V\dff,\quad
E^{\dff +}_{\qff V}\dff,\quad
\mbox{and}\quad
E^{\dff -}_{\qff V}
\]

\vspace{-12pt}\vspace{1.5pt}\vspace{0.625pt}
Clearly,\oss the family\sss
$\bm{\varphi}\trf(\trf z\trf)\fff,\qff z\qff \in\qff Z$\sss
defines an\sss isometry\sss 
$\bm{\varphi}\dff \colon\dff 
\pi^{\fff *}\dff E^{\dff +}_{\qff V}
\qff \ttoo\qff
\pi^{\fff *}\dff E^{\dff -}_{\qff V}$,\oss
where now\sss $\pi$\sss is\dss the projection\sss
$S\fff V\qff \ttoo\qff V$\dnsp,\oss
and\sss hence defines a class\vspace{1.5pt}\vspace{0.625pt}
\[
\quad
\mathcal{I}\dff(\trf N\trf)
\off \in\off
K^{\dff 0}\dff (\trf B\dff V,\pff S\fff V\trf)
\off =\off
K^{\dff 0}\dff (\trf T\dff V\trf)
\off,
\]

\vspace{-12pt}\vspace{1.5pt}\vspace{0.625pt}
the parametrized\sss version of\dss the class\sss $\mathcal{I}\dff(\trf N\trf)$\sss
from\dss Section\qss \ref{symbols-conditions}.\oss
By\sss the\dss Atiyah--Singer\dss index\sss theorem\sss for\sss families\qss \cite{as4}\qss
the analytical\dss index of\dss the family\sss 
$\bm{\varphi}\trf(\trf z\trf)\fff,\qff z\qff \in\qff Z$\sss
is\dss equal\dss to\sss the forward\sss image of\dss the class\sss
$\mathcal{I}\dff(\trf N\trf)$\sss in\sss $K^{\dff 0}\dff (\trf Z\trf)$\dss
(cf.\qss the definition of\dss the\sss topological\dss index\sss
in\dss Section\qss \ref{symbols-conditions}).\oss
In\sss particular,\oss if\dss $\mathcal{I}\dff(\trf N\trf)$\sss is\dss
equal\dss to zero,\oss then\sss the analytical\dss index of\dss
$\bm{\varphi}\trf(\trf z\trf)\fff,\qff z\qff \in\qff Z$\sss
is\dss also equal\dss to zero.\oss
At\sss the same\sss time\sss the class\sss $\mathcal{I}\dff(\trf N\trf)$\sss
is\dss the obstruction\sss to\sss defining\sss the\sss topological\sss index
of\dss $(\trf \sigma,\pff N\trf)$\nnsp,\oss
as\sss the parametrized\sss version of\trs Proposition\qss \ref{obstruction}\qss shows.\oss
Therefore,\oss if\dss the\sss topological\dss index of\sss $(\trf \sigma,\pff N\trf)$\sss
is\dss defined,\oss then\sss there exists a realization of\dss
the family of\dss boundary conditions\sss
$N\trf(\trf z\trf)\fff,\qff z\qff \in\qff Z$\nnsp.\oss

\myuppar{The non-uniqueness of\dss realizations.}
It\dss turns out\sss that\sss even\sss
when a realization of\dss a family\sss of\dss boundary conditions\sss
$N\trf(\trf z\trf)\fff,\qff z\qff \in\qff Z$\sss
exists,\oss it\dss is\dss usually\sss not\sss unique,\oss
even up\sss to homotopy.\oss
Since\sss the realizations of\dss $N\trf(\trf z\trf)\fff,\qff z\qff \in\qff Z$\sss
are in one-to-one correspondence with\sss the isometric realizations of\dss
$\bm{\varphi}\trf(\trf z\trf)\fff,\qff z\qff \in\qff Z$\nnsp,\oss
in order\sss to prove\sss this,\oss
it\dss is\dss sufficient\sss to prove\sss that\sss the\sss latter are not\sss unique.\oss
Let\sss $\Phi\trf(\trf z\trf)\fff,\qff z\qff \in\qff Z$\sss be an\sss isometric
realization of\dss $\bm{\varphi}\trf(\trf z\trf)\fff,\qff z\qff \in\qff Z$\nnsp.\oss

Let\sss us\sss return\sss to\sss
the proof\dss of\trs Theorem\qss \ref{unitary-sections}\qss
in\sss the\sss case of\dss the family\dss 
$A_{\dff z}\off =\off \Phi\trf(\trf z\trf)\fff,\off z\qff \in\qff Z$\nnsp.\oss
Sections of\trs the bundle\sss
$\bm{\pi}\ffin\dff(\trf \Phi\trf)\dff \colon\dff
U^{\dff \fin}\dff(\trf \Phi\trf)
\qff \ttoo\qff
Z$\dss
define\sss isometric realizations\sss
$\Phi\fff'\trf(\trf z\trf)\fff,\off z\qff \in\qff Z$\sss
such\sss that\sss $\Phi\fff'\trf(\trf z\trf)$\sss is\dss
equal\dss to\sss $\Phi\fff'\trf(\trf z\trf)$ on a closed subspace of\dss finite codimension.\oss
Since an\sss isometric realization exists,\oss
this bundle\dss is\dss trivial,\oss as was explained\sss in\sss the
proof\dss of\trs Theorem\qss \ref{unitary-sections}.\oss
In\sss general,\oss an\sss isometric realization\sss
$\Phi\fff'\trf(\trf z\trf)\fff,\off z\qff \in\qff Z$\sss
will\sss differ\sss from\sss
$\Phi\trf(\trf z\trf)\fff,\off z\qff \in\qff Z$\sss
by\sss compact\sss operators.\oss
This suggests\sss to consider\sss the bundle\sss\vspace{2.5pt}
\[
\quad
\bm{\pi}\comp\dff(\trf \Phi\trf)\dff \colon\dff
U\comp\dff(\trf \Phi\trf)
\qff \ttoo\qff
Z
\]

\vspace{-12pt}\vspace{2.5pt}
having as a fiber over\sss $z\qff \in\qff Z$\sss
the space of\dss isometries $U$ such\sss that\sss
$U\qff -\qff \Phi\trf(\trf z\trf)$\sss is\dss a compact\sss operator.\oss
There\dss is\dss also a\qss ``universal''\qss bundle\vspace{2.5pt}
\[
\quad
\bm{\pi}\comp\dff \colon\dff
\mathbf{U}\comp
\qff \ttoo\qff 
\num{\mathcal{P}{\nsp}\mathcal{S}_{\qff 0}}
\]

\vspace{-12pt}\vspace{2.5pt}
over\sss the same base\sss $\num{\mathcal{P}{\nsp}\mathcal{S}_{\qff 0}}$\sss
as\sss $\bm{\pi}\ffin$\nnsp,\oss and\sss
$\bm{\pi}\comp\dff(\trf \Phi\trf)$\sss is\dss induced\sss 
from\sss $\bm{\pi}\comp$\sss by\sss the same polarized\dss index\sss map.\oss 
Therefore\sss the bundle\sss $\bm{\pi}\comp\dff(\trf \Phi\trf)$\sss
is\dss trivial\dss by\sss the same reason as\sss the bundle $\bm{\pi}\ffin\dff(\trf \Phi\trf)$\nnsp.\oss
The fibers of\dss the bundles\sss $\bm{\pi}\ffin\dff(\trf \Phi\trf)$\sss
and\sss $\bm{\pi}\comp\dff(\trf \Phi\trf)$\sss
are homeomorphic,\oss respectively,\oss
to\sss the space\sss $U\ffin$\sss
of\dss isometries of\dss a\sss fixed\dss Hilbert\dss space equal\dss
to\sss the identity on a closed subspace of\dss finite codimension,\oss
and\sss to\sss the space\sss $U\comp$\sss of\dss isometries differing\sss
from\sss the identity\sss by a compact\sss operator.\oss
By\trs Theorem\qss \ref{unitary-spaces}\qss
both of\dss them are homotopy equivalent\sss
to\sss the infinite unitary\sss group $U\trf(\trf \infty\trf)$\sss
and\dss hence are classifying spaces for\sss $K^{\dff 1}$\dnsp-theory.\oss

Therefore every\sss two\sss isometric realizations define an element\sss of\dss
$K^{\dff 1}\dff(\trf Z\trf)$\sss and,\oss given an\sss isometric realization,\oss
every element\sss of\dss $K^{\dff 1}\dff(\trf Z\trf)$\sss 
results from some other\sss isometric realization.\oss
Therefore\sss the natural\sss data on\sss the\sss level\sss
of\dss principal\sss symbols do not\sss lead\dss to a well\sss defined,\oss
even only\sss up\sss to homotopy,\oss family of\dss self-adjoint\dss boundary\sss problems.\oss
Moreover,\oss as we will\sss see in a moment,\oss changing\sss the realization of\dss
a family of\dss boundary conditions may change\sss the index of\dss the corresponding\sss
family of\dss self-adjoint\sss operators.\oss
But\dss for\sss bundle-like boundary conditions\sss
there are preferred\sss realizations,\oss namely\sss the ones induced\dss
by\dss bundle maps.\oss

\myuppar{Changing\sss realization and\dss the analytical\dss index.}
Let\sss $E\fff,\qff P\halfff,\qss \Sigma\fff,\qff N$\sss etc.\dss be as above\qss
(without\sss parameters).\oss
As before,\oss we will\sss assume\sss that\sss 
the endomorphism $\Sigma$\sss of\dss $E\trf|\trf Y$\sss is\dss an\sss isometry.\oss
We will\sss assume\sss that\sss $N$\sss is\dss bundle-like and\dss treat\sss $N$\sss
as a subbundle of\dss $E\trf|\trf Y$\dnsp.\oss
The symbol\sss $\sigma$\sss of\dss $P$ defines a decomposition\sss
$E\trf|\trf Y\off =\off E^{\dff +}_{\trf Y}\qff \oplus\qff E^{\dff -}_{\trf Y}$\nsp.\oss
Since $N$\sss is\dss bundle-like and\sss $\Sigma$\sss is\dss an\sss isometry,\pss
$N$ defines an\sss isometry\sss 
$\varphi\dff \colon\dff E^{\dff +}_{\trf Y}\qff \ttoo \qff E^{\dff -}_{\trf Y}$\nsp.\oss 
Let\sss $\Delta$ and\sss $\Delta^{\fff \perp}$\sss be\sss the graphs of\dss $\varphi$
and $-\qff \varphi$\sss respectively.\oss
Then\sss $E\trf|\trf Y\off =\off \Delta\dff \oplus\dff \Delta^{\fff \perp}$
and\sss $N\off =\off \Delta\dff \oplus\dff 0$\nnsp.\oss
If\dss we identify $\Delta$ with $\Delta^{\fff \perp}$\sss by\sss the map\sss
$(\trf x\fff,\qff \varphi\trf(\trf x\trf)\trf)
\off \longmapsto\off
(\trf x\fff,\qff -\qff \varphi\trf(\trf x\trf)\trf)$\nnsp,\oss
then $\Sigma$ will\dss take\sss the form\vspace{1.25pt}
\[
\quad
\Sigma
\off =\off\dff
\begin{pmatrix}
\off 0 &
1 \qff\off
\vspace{4.5pt} \\
\off\dff 1 &
0 \qff\off 
\end{pmatrix}
\off.
\]

\vspace{-12pt}\vspace{1.25pt}
The projection\sss 
$E\trf|\trf Y
\off =\off 
\Delta\dff \oplus\dff \Delta^{\fff \perp} 
\qff \ttoo\qff 
\Delta$\sss
induces an orthogonal\sss projection $\Pi$ in\sss 
$H_{\trf 0}\dff(\trf Y\fff,\qff E\trf |\trf Y\trf)$\nnsp.\oss
Let\sss $K$\sss be a finitely dimensional\sss subspace of\dss
the intersection\sss
$\kernel\fff \Pi
\qff \cap\qff
H_{\dff 1/2}\dff(\trf Y\fff,\qff E\trf |\trf Y\trf)$\nnsp.\oss
As we saw in\dss Section\qss \ref{boundary-triplets},\oss
every closed\sss relation\sss $\mathcal{C}\qff \subset\qff K\dff \oplus\dff K$\sss
leads\sss to an orthogonal\sss projection\sss $\Pi_{\dff \mathcal{C}}$\sss 
in\sss $H^{\dff \partial}$\dnsp.\oss
If\sss $\mathcal{C}$\sss is\dss self-adjoint,\oss
then\sss the abstract\sss boundary\sss problem\sss 
$P{},\pff \Pi_{\dff \mathcal{C}}$\sss is\dss self-adjoint\dss 
by\trs Lemma\qss \ref{two-self-adjointness},\oss
and\dss hence\sss the operator\sss 
$P_{\trf \Gamma_{\dff \mathcal{C}}}$\sss defined\dss by\sss
$P{},\pff \Pi_{\dff \mathcal{C}}$\sss is\dss self-adjoint.\oss
As we pointed out\sss before\dss Lemma\qss \ref{boundary-problems-triplets},\oss
$\Pi_{\dff \mathcal{C}}$ is\dss equal\dss to $\Pi$
on\sss the orthogonal\sss complement\sss of\dss
$\Sigma\trf(\trf K\trf)\qff +\qff K$\sss
and\sss differs from $\Pi$\sss by an operator\sss of\dss finite rank.\oss
Hence\sss $\Pi_{\dff \mathcal{C}}$\sss is\dss a pseudo-differential\sss
operator of\dss order $0$ with\sss the same symbol\sss as $\Pi$\nnsp,\oss
and\dss is\dss a realization of\dss $N$\nnsp.\oss
Note\sss that\sss 
the self-adjointness of\sss $\mathcal{C}$\sss
implies\sss that\sss $\image\fff \Pi_{\dff \mathcal{C}}$\sss is\dss the graph
of\dss a map\sss
$H_{\trf 0}\dff(\trf Y\fff,\qff E^{\dff +}_{\trf Y}\trf)
\qff \ttoo\qff
H_{\trf 0}\dff(\trf Y\fff,\qff E^{\dff -}_{\trf Y}\trf)$\nnsp.

Suppose now\sss that\sss $W$\sss be a\sss topological\sss space and\sss
$\mathcal{C}_{\fff w}\dff,\pff w\qff \in\qff W$\sss
be a continuous family of\dss self-adjoint\sss relations\sss
$\mathcal{C}_{\fff w}\qff \subset\qff K\dff \oplus\dff K$\nnsp.\oss
For $w\qff \in\qff W$\dss let\sss\vspace{3pt} 
\[
\quad
\Pi_{\dff w}\off =\off \Pi_{\dff \mathcal{C}_{\fff w}}\qff,\quad
\Gamma_{\fff w}\off =\off \Gamma_{\dff \mathcal{C}_{\fff w}}\qff,\quad
\mbox{and}\quad\dff
P_{\fff w}\off =\off P_{\trf \Gamma_{\fff w}}
\pff.
\]

\vspace{-12pt}\vspace{3pt}
By\trs Theorem\qss \ref{relative-index}\qss
the family\sss $P_{\fff w}\dff,\pff w\qff \in\qff W$\sss
is\dss also\dss Fredholm\dss and\sss its analytical\dss index\dss
is\dss equal\dss to\sss the analytical\dss index of\dss the family\sss
$\mathcal{C}_{\fff w}\dff,\pff w\qff \in\qff W$\nnsp.\oss
The analytical\dss index of\dss $P_{\fff w}\dff,\pff w\qff \in\qff W$\sss
is\dss nothing else as\sss the analytical\dss index of\dss the family\sss
$(\trf P\halfff,\pff \Pi_{\dff w}\trf)$\sss
of\dss well-posed\dss boundary\sss problems in\sss the sense of\trs Seeley.\oss
The operators in\sss this family do not\sss depend on\sss the parameter.\oss
The boundary conditions depend on\sss the parameter $w$\nnsp,\oss
but\sss for all\sss values of\dss $w$ are realizations of\dss $N$\nnsp.\oss
Hence\sss the\sss topological\dss index of\dss this family\dss is\dss
equal\dss to\sss the\sss topological\dss index\sss of\dss the constant\sss family
of\dss boundary problems equal\dss to $(\trf P\halfff,\pff \Pi\trf)$
and\dss hence\dss is\dss equal\dss to zero.\oss

At\sss the same\sss time\sss the analytical\dss index of\dss
$\mathcal{C}_{\fff w}\dff,\pff w\qff \in\qff W$\nnsp,\oss
and\dss hence of\dss
$(\trf P\halfff,\pff \Pi_{\dff w}\trf)$\nnsp,\oss
is\dss quite arbitrary.\oss
For example,\oss if\trs $W$\sss is\dss compact,\oss
then every element\sss of\dss $K^{\dff 1}\dff(\trf W\trf)$\sss
can\sss be realized as\sss the analytical\dss index of\dss a family\sss
$\mathcal{C}_{\fff w}\dff,\pff w\qff \in\qff W$\nnsp.\oss
Indeed,\oss let\sss us consider\sss the family\sss
$V_w\dff, \pff w\qff \in\qff W$
of\dss the\dss Cayley\trs transforms\sss
$V_w\off =\off V\trf(\trf \mathcal{C}_{\fff w} \trf)$\nnsp.\oss
The discussion at\sss the end of\trs Section\qss \ref{relations}\qss
shows\sss that\sss the analytical\dss index of\dss
$\mathcal{C}_{\fff w}\dff,\pff w\qff \in\qff W$
can\sss be identified\sss with\sss the homotopy class of\dss the map\sss
$w\off \longmapsto\off V_w$\sss from $W$\sss to\sss the unitary\sss group of\dss $K$\nnsp.\oss
For compact\sss $W$\sss an element\sss $a\qff \in\qff K^{\dff 1}\dff(\trf W\trf)$\sss
can\sss be realized\sss as such a homotopy class\sss if\dss the dimension of\dss $K$\sss
is\dss sufficiently\sss large.\oss
We see\sss that\sss the analytical\dss index\sss is\dss not\sss determined\dss
by\sss the\sss topological\dss index even\sss when\sss the boundary conditions
differ\sss from\sss the bundle-like ones by operators of\dss finite rank.\oss

\newpage
\mysection{Special\qss boundary\qss conditions}{special-conditions}

\myuppar{Anti-commuting elliptic pairs and\sss special\dss boundary conditions.}
For a major part\sss of\dss this section we will\sss work in\sss the
framework of\trs Section\qss \ref{boundary-algebra}.\oss
We will\sss say\sss that\sss an elliptic\sss pair $\sigma\fff,\qff \tau$\sss is\qss 
\emph{anti-commuting}\pss if\dss $\sigma$ 
anti-commutes with $\tau$ and\sss hence with $\rho\off =\off \sigma^{\dff -\dff 1}\dff \circ\dff \tau$\nnsp.\oss
In\sss this case\sss $\sigma\trf \rho\trf \sigma^{\dff -\dff 1}\off =\off -\qff \rho$\sss
and\sss hence 
$\sigma\trf(\trf \mathcal{L}_{\dff +}\dff(\trf \rho\dff)\trf)
\off =\off
\mathcal{L}_{\dff -}\dff(\trf \rho\dff)$\sss
and\sss
$\sigma\trf(\trf \mathcal{L}_{\dff -}\dff(\trf \rho\dff)\trf)
\off =\off
\mathcal{L}_{\dff +}\dff(\trf \rho\dff)$\nnsp.\oss
We will\sss say\sss that\sss a boundary condition $N$\sss for an elliptic pair $\sigma\fff,\qff \tau$\sss
is\qss \emph{special}\pss if\dss $N$\sss is\dss transverse\sss
not\sss only\sss to\sss 
$\mathcal{L}_{\dff -}\dff(\trf \rho\dff)$\nnsp,\oss
but\sss also\sss to\sss
$\mathcal{L}_{\dff +}\dff(\trf \rho\dff)$\nnsp.\oss

Let\sss $\sigma\fff,\qff \tau$\sss be a self-adjoint\sss anti-commuting elliptic pair
and\sss $N$\sss be a special\sss boundary condition for it.\pss
We will\sss deform\sss $\sigma\fff,\off \rho\fff,\off \tau\fff,\off N$\sss
to a normal\sss form\sss while keeping\sss these properties.\oss 

To begin\sss with,\oss let\sss us apply\sss to $\sigma\fff,\qff \tau$\sss 
the first\sss two deformations from\dss Section\qss \ref{boundary-algebra}.\oss
Clearly,\oss during\sss the first\sss deformation\sss $\sigma_{\dff \alpha}$\sss
remains anti-commuting\sss with\sss $\tau_{\dff \alpha}$\sss and\sss
$N_{\dff \alpha}$\sss stays\sss transverse\sss to\sss
$\mathcal{L}_{\dff -}\dff(\trf \rho_{\dff \alpha}\dff)$\sss 
and\sss by\sss the same reason\sss 
stays\sss transverse\sss to\sss $\mathcal{L}_{\dff +}\dff(\trf \rho_{\dff \alpha}\dff)$\nnsp.\oss
The second deformation does not\sss affect $\sigma$\dnsp,\oss 
the product $[\trf \bullet\fff,\qff \bullet\trf]$\nnsp,\oss 
and\sss the spaces
$\mathcal{L}_{\dff +}\dff(\trf \rho\dff)\fff,\off
\mathcal{L}_{\dff -}\dff(\trf \rho\dff)$\nnsp.\oss
Hence during\sss the second deformation\sss $N$\sss remains\sss transverse\sss to\sss
$\mathcal{L}_{\dff +}\dff(\trf \rho\dff)$\sss
and\sss
$\mathcal{L}_{\dff -}\dff(\trf \rho\dff)$\nnsp.\oss
Since\sss
$\sigma\trf(\trf \mathcal{L}_{\dff +}\dff(\trf \rho\dff)\trf)
\off =\off
\mathcal{L}_{\dff -}\dff(\trf \rho\dff)$\sss
and\sss
$\sigma\trf(\trf \mathcal{L}_{\dff -}\dff(\trf \rho\dff)\trf)
\off =\off
\mathcal{L}_{\dff +}\dff(\trf \rho\dff)$\nnsp,\oss 
the map\sss $\rho_{\dff 0}$\sss anti-com\-mutes with $\sigma$\sss
and\sss hence during\sss the second\sss deformation\sss $\rho\fff,\qff \tau$\sss 
remain anti-commuting with $\sigma$\nnsp.\oss
After\sss the second deformation\sss the spaces\sss
$\mathcal{L}_{\dff +}\dff(\trf \rho\dff)\fff,\off
\mathcal{L}_{\dff -}\dff(\trf \rho\dff)$\sss
are\sss the eigenspaces of\dss $\rho$ with\sss the eigenvalues
$i\fff,\off -\qff i$ respectively.\oss
Let\sss us represent\sss 
$\mathcal{L}_{\dff +}\dff(\trf \rho\dff)$\sss
as\sss the graph of\dss an\sss isometry\sss
$\varphi\dff \colon\dff 
E^{\dff +}\qff \ttoo\qff E^{\dff -}$\dnsp.\oss
Then\sss\vspace{2pt}
\[
\quad
\mathcal{L}_{\dff +}\dff(\trf \rho\dff)
\off =\off
\Delta\dff(\trf \varphi\trf)
\quad
\mbox{and}\quad
\mathcal{L}_{\dff -}\dff(\trf \rho\dff)
\off =\off
\Delta^{\fff \perp}\dff(\trf \varphi\trf)
\qff,
\]

\vspace{-12pt}\vspace{2pt}
and\dss hence\sss $\sigma\fff,\off \rho\fff,\off \tau$\sss
have\sss the standard\sss form\qss (\ref{sigma-rho-tau-phi}).\oss
But\sss now\sss $N$\sss does not\sss have\sss the form\qss (\ref{delta-delta})\qss
and\sss $\sigma\fff,\off \tau\fff,\off N$\sss are not\sss normalized.\oss
Still,\oss already after\sss the first\sss deformation $\sigma$\sss is\dss unitary\sss and\sss hence 
$N^{\dff \perp}\off =\off \sigma\trf(\trf N\trf)$\nnsp.\oss
It\sss follows\sss that\sss $N^{\dff \perp}$\sss
is\dss a\sss lagrangian subspace\sss transverse\sss to both\sss
$\mathcal{L}_{\dff +}\dff(\trf \rho\dff)$\sss
and\sss
$\mathcal{L}_{\dff -}\dff(\trf \rho\dff)
\off =\off
\sigma\trf(\trf \mathcal{L}_{\dff +}\dff(\trf \rho\dff)\trf)$\nnsp.\oss
Equivalently,\pss $N^{\dff \perp}$\sss 
is\dss a special\dss boundary condition\sss for\sss $\sigma\fff,\qff \tau$\nnsp.\oss

Now we are ready\sss for\sss the\sss third and\dss the\sss last\sss step of\dss our deformation.\oss
Let\sss us present\sss $N$\sss as\sss the graph of\dss an\sss isometry\sss
$\psi\dff \colon\dff 
E^{\dff +}\qff \ttoo\qff E^{\dff -}$\dnsp.\oss
Since $\mathcal{L}_{\dff +}\dff(\trf \rho\dff)$\sss is\dss transverse\sss to\sss
$N$\sss and\sss $N^{\dff \perp}$\dnsp,\oss
the isometry\sss $\psi^{\dff -\dff 1}\dff \circ\trf \varphi$\sss
has no eigenvalues equal\dss to $1$ or $-\qff 1$\nnsp.\oss
Therefore we can deform\sss $\varphi$\sss in\sss the class of\dss isometries
with\sss this property\dss to a new\sss isometry\sss $\varphi$\sss such\sss that\sss
$\psi^{\dff -\dff 1}\dff \circ\trf \varphi$\sss has only\sss $i$ and $-\qff i$
as eigenvalues.\oss
This deformation of\dss $\varphi$\sss induces a deformation\sss of\dss
the graph of\dss $\varphi$\nnsp,\oss
i.e.\oss of\dss
$\mathcal{L}_{\dff +}\dff(\trf \rho\dff)$\nnsp.\oss
Simultaneously we deform\sss
$\mathcal{L}_{\dff -}\dff(\trf \rho\dff)$\sss
as\sss the image of\dss
$\mathcal{L}_{\dff +}\dff(\trf \rho\dff)$\sss
under $\sigma$\dnsp.\oss
This deformation of\dss 
$\mathcal{L}_{\dff +}\dff(\trf \rho\dff)$\sss
and\sss
$\mathcal{L}_{\dff -}\dff(\trf \rho\dff)$\sss
defines a deformation of\dss $\rho$ and\dss hence a deformation of\dss 
$\tau\off =\off \sigma\dff \circ\dff \rho$\nnsp.\oss
We do not\sss move $\sigma$ and\sss $N$\sss during\sss
this deformation and\dss hence keep\sss the anti-commutation property
and\sss $N$\sss being\sss a special\dss boundary condition.\oss
This\sss third deformation\sss partially\sss replaces both\sss the\sss third and\sss
the fourth\sss deformations from\dss Section\qss \ref{boundary-algebra}.\oss

\myuppar{Bott\trs periodicity\sss maps.}
Let\sss us\sss recall\sss some constructions of\trs Bott\qss \cite{bott}.\oss
Let\sss $F$\sss be a complex vector space with a\dss Hermitian\dss metric.\oss
As usual,\oss for a subspace\sss $A\qff \subset\qff F$\sss we will\sss denote by\sss $A^{\fff \perp}$ 
its orthogonal\sss complement.\oss
For a\sss linear\sss isomorphism\sss
$a\dff \colon\dff F\qff \ttoo\qff F$\sss let\vspace{3pt}
\[
\quad
\lambda\trf(\trf a\dff,\pff \theta\trf)
\off =\off 
\bigl\{\pff 
\bigl(\qff u\dff \cos\dff \theta/2\fff,\off 
a\trf(\dff u\trf)\trf \sin\dff \theta/2\qff\bigr)
\pff\fff \bigl|\pff\fff
u\qff \in\qff F
\pff\bigr\}
\off \subset\off
F\dff \oplus\dff F
\qff,\quad
\theta\qff \in\qff
[\trf 0\fff,\qff \pi\trf]
\qff.
\]

\vspace{-12pt}\vspace{3pt}
The path
$\lambda\trf(\trf a\dff,\pff \theta\trf)\dff,\off
\theta\qff \in\qff
[\trf 0\fff,\qff \pi\trf]$
connects $F\dff \oplus\dff 0$ with $0\dff \oplus\dff F$\dnsp.\oss
Cf.\qss Bott\qss \cite{bott},\oss (4.4)\qss and\qss (4.5).\oss
Now,\oss let\sss $A$\sss be a vector subspace of\dss $F$\nnsp.\oss
For $\eta\qff \in\qff [\trf 0\fff,\qff \pi\trf]$\sss let\sss us define a map
$f\trf(\trf A\dff,\qff \eta\qff)\dff \colon\dff F\qff \ttoo\qff F$\sss
by\vspace{3pt}
\[
\quad
f\trf(\trf A\dff,\qff \eta\qff)\dff(\trf u\trf)
\off =\off
u\dff e^{\dff i\dff \eta}
\quad\hspace{0.78em}
\mbox{for}\quad
u\qff \in\qff A
\quad
\mbox{and}
\]

\vspace{-36.0pt}
\[
\quad
f\trf(\trf A\dff,\qff \eta\qff)\dff(\trf u\trf)
\off =\off
u\dff e^{\qff -\dff i\dff \eta}
\quad
\mbox{for}\quad
u\qff \in\qff A^{\fff \perp}
\qff.
\]

\vspace{-12pt}\vspace{3pt}
The path\sss
$f\trf(\trf A\dff,\qff \eta\qff)\dff,\off
\theta\qff \in\qff
[\trf 0\fff,\qff \pi\trf]$\sss
connects\sss $\id$\sss with\sss $-\qff \id$\nnsp.\oss
Cf.\qss Bott\qss \cite{bott},\oss (5.1).\oss
One can extend\dss this path by a path\sss
$f\trf(\trf A\dff,\qff \eta\qff)\dff,\off
\eta\qff \in\qff
[\trf \pi\fff,\qff 2\dff \pi\trf]$\sss
returning\sss from\sss $-\qff \id$\sss to $\id$
and actually\dss independent\sss of\dss $A$\nnsp.\oss
One can\sss take\sss
$f\trf(\trf A\dff,\qff \eta\qff)
\off =\off 
\exp\trf (\trf i\trf \eta \qff)$\sss
for\sss
$\eta\qff \in\qff
[\trf \pi\fff,\qff 2\dff \pi\trf]$\nnsp,\oss
as we will\sss do,\oss
although\dss Bott\dss uses\sss the path\sss
$f\dff(\trf A_{\trf 0}\dff,\qff \eta\trf)\dff,\off  
\eta\qff \in\qff [\trf \pi\fff,\qff 2\dff \pi\trf]$\sss
for a fixed subspace $A_{\trf 0}$\nsp.\oss
Bott\qss \cite{bott}\qss assumes\sss that\sss
$\dim\dff A\off =\off \dim\dff A^{\fff \perp}$\dnsp,\oss
but\sss the definitions make sense without\sss this assumption.\oss
These\sss two constructions can\sss be combined as follows.\oss
For a vector subspace\sss $A\qff \subset\qff F$\sss let\vspace{3pt}
\[
\quad
\gamma\trf(\trf A\dff,\pff \eta\dff,\pff \theta\trf)
\off =\off
\lambda\trf\left(\qff 
f\trf(\trf A\dff,\qff \eta\qff)\dff,\pff
\theta\qff\right)
\off \subset\off
F\dff \oplus\dff F
\qff,\quad
\eta\qff \in\qff
[\trf 0\fff,\qff 2\dff \pi\trf]\dff,\off\qff
\theta\qff \in\qff
[\trf 0\fff,\qff \pi\trf]
\qff.
\]

\vspace{-12pt}\vspace{3pt}
More explicitly,\oss for\sss $\eta\qff \in\qff [\trf 0\fff,\qff \pi\trf]$\sss the subspace
$\gamma\trf(\trf A\dff,\pff \eta\dff,\pff \theta\trf)$\sss
is\dss the sum of\dss subspaces\vspace{3pt}
\[
\quad
\left\{\pff 
\left(\qff 
u\dff \cos\dff \theta/2\fff,\off\hspace{0.78em} 
u\dff e^{\dff i\dff \eta}
\trf \sin\dff \theta/2
\qff\right)
\pff\fff \bigl|\pff\dff
u\qff \in\qff A\phantom{^{\fff \perp}}
\pff\right\}
\quad
\mbox{and}\quad
\]

\vspace{-30.0pt}
\[
\quad
\left\{\pff 
\left(\qff 
u\dff \cos\dff \theta/2\fff,\off 
u\dff e^{\qff -\dff i\dff \eta}
\trf \sin\dff \theta/2
\qff\right)
\pff\fff \bigl|\pff\dff
u\qff \in\qff A^{\fff \perp}
\pff\right\}
\qff.
\]

\vspace{-12pt}\vspace{3pt}
The map\sss assigning\sss to a subspace $A$ of\dss $F$\sss
the family\sss
$\gamma\trf(\trf A\dff,\pff \eta\dff,\pff \theta\trf)
\dff,\off
\eta\qff \in\qff
[\trf 0\fff,\qff 2\dff \pi\trf]
\dff,\pff 
\theta\qff \in\qff
[\trf 0\fff,\qff \pi\trf]$\sss
of\dss subspaces of\dss $F\dff \oplus\dff F$\sss is\dss
essentially\sss the\dss Bott\trs periodicity\dss map.\oss
In order\sss to get\dss the\dss Bott\trs periodicity\dss map\sss itself,\pss
one needs\sss to\sss turn\sss the rectangle\sss
$[\trf 0\fff,\qff 2\dff \pi\trf]
\dff \times\dff 
[\trf 0\fff,\qff \pi\trf]$\sss
into a sphere by\sss identifying\sss the sides\sss $0\dff \times\dff [\trf 0\fff,\qff \pi\trf]$\sss
and\sss $2\dff \pi\dff \times\dff [\trf 0\fff,\qff \pi\trf]$\sss
and collapsing each of\dss the sides\sss 
$[\trf 0\fff,\qff 2\dff \pi\trf]\dff \times\dff 0$\sss
and\sss
$[\trf 0\fff,\qff 2\dff \pi\trf]\dff \times\dff \pi$\sss
to a pole of\dss the sphere.\oss

In\dss Bott's\dss construction a subspace\sss $A\qff \subset\qff F$\sss
is\dss encoded\sss by\sss the operator\sss
$f\trf(\trf A\dff,\qff \pi/2\qff)$\sss
equal\sss to\sss the multiplication by $i$ on $A$ and\dss to\sss
the multiplication by $-\qff i$ on $A^{\fff \perp}$\dnsp.\oss
This\dss is\dss a skew-adjoint\sss operator.\oss
For our purposes\sss it\dss is\dss more convenient\sss
to encode $A$ by\sss the self-adjoint\sss operator
equal\sss to\sss the identity on $A$ and\dss to\sss
the minus identity on $A^{\fff \perp}$\dnsp,\oss
i.e.\qss by\sss the operator\sss
$-\qff i\trf f\trf(\trf A\dff,\qff \pi/2\qff)$\nnsp.\oss
This\sss leads\sss to\sss the following\sss modification of\dss the\dss
Bott\trs periodicity\sss map.\oss
Let\vspace{3pt}
\[
\quad
f\fff'\trf(\trf A\dff,\qff \eta\qff)
\off =\off
-\qff i\trf f\trf(\trf A\dff,\qff \eta\qff)
\quad
\mbox{and}\quad
\gamma\fff'\trf(\trf A\dff,\pff \eta\dff,\pff \theta\trf)
\off =\off
\lambda\trf\left(\qff 
f\fff'\trf(\trf A\dff,\qff \eta\qff)\dff,\pff
\theta\qff\right)
\qff.
\]

\vspace{-12pt}\vspace{3pt}
The family\sss 
$f_{\fff t}\trf(\trf A\dff,\qff \eta\qff)
\off =\off
e^{\trf -\trf i\fff t}\dff f\trf(\trf A\dff,\qff \eta\qff)\dff,\off
t\qff \in\qff [\trf 0\fff,\qff \pi/2\trf]$\sss
substituted\sss into\sss
$\lambda\trf(\trf \bullet\dff,\pff \theta\trf)$\sss
leads\sss to a homotopy\sss between\sss the maps assigning\sss to $A$\sss
families\sss $\gamma\trf(\trf A\dff,\pff \eta\dff,\pff \theta\trf)$
and\sss
$\gamma\fff'\trf(\trf A\dff,\pff \eta\dff,\pff \theta\trf)$\sss
respectively.\oss
We will\sss use\sss the\sss latter.\oss
For\sss $\eta\qff \in\qff [\trf 0\fff,\qff \pi\trf]$\sss the subspace
$\gamma\fff'\trf(\trf A\dff,\pff \eta\dff,\pff \theta\trf)$\sss
is\dss the sum of\dss subspaces\vspace{4.5pt}
\[
\quad
\left\{\pff 
\left(\qff 
u\dff \cos\dff \theta/2\fff,\off\hspace{0.78em} 
u\dff e^{\dff i\trf (\trf \eta\qff -\qff \pi/2 \trf)}
\trf \sin\dff \theta/2
\qff\right)
\pff\fff \bigl|\pff\dff
u\qff \in\qff A\phantom{^{\fff \perp}}
\pff\right\}
\quad
\mbox{and}\quad
\]

\vspace{-30.0pt}
\[
\quad
\left\{\pff 
\left(\qff 
u\dff \cos\dff \theta/2\fff,\off 
u\dff e^{\qff -\dff i\trf (\trf \eta\qff +\qff \pi/2 \trf)}
\trf \sin\dff \theta/2
\qff\right)
\pff\fff \bigl|\pff\dff
u\qff \in\qff A^{\fff \perp}
\pff\right\}
\qff.
\]

\vspace{-12pt}\vspace{3pt}\vspace{0.5pt}
\myuppar{The positive eigenspaces of\sss 
$\sigma\dff \cos\dff \theta\qff +\qff \tau\dff \sin\dff \theta$\nnsp.}
For\sss $\sigma\fff,\off \rho\fff,\off \tau$\sss
having\sss the standard\sss form\qss (\ref{sigma-rho-tau-phi}),\oss
the family\sss 
$E^{\dff +}\fff(\trf \theta\trf)\fff,\off
\theta\qff \in\qff
[\trf 0\fff,\qff \pi\trf]$\sss 
of\dss these eigen\-spaces was determined at\sss the end of\trs
Section\qss \ref{boundary-algebra}.\oss
A reparameterized\sss version of\dss this family\sss better\sss
matches\sss the\dss Bott\trs periodicity\sss map.\oss
When $\theta$ runs over\sss $[\trf 0\fff,\qff \pi\trf]$\nnsp,\oss the fraction\sss
$(\trf 1\qff -\qff z\trf)/(\trf 1\qff +\qff z\trf)$\nnsp,\oss
where\sss $z\off =\off \cos\dff \theta\qff +\qff i\dff \sin\dff \theta$\nnsp,\oss 
runs over\sss the negative\sss imaginary\sss half-line\sss
$i\trf \rrr_{\qff \leq\qff 0}$.\oss
The formula\qss (\ref{plus-path})\qss implies\sss that\sss
up\sss to a canonical\dss homotopy\sss the path\sss 
$E^{\dff +}\fff(\trf \theta\trf)\fff,\off
\theta\qff \in\qff
[\trf 0\fff,\qff \pi\trf]$\sss
is\dss the same as\vspace{4.5pt}
\[
\quad
\left\{\pff \left.
\left(\qff u\fff,\off 
-\qff i\trf \varphi\trf(\trf u\trf)\pff 
\frac{\fff \sin\dff \theta/2\fff}{\fff \cos\dff \theta/2\fff}\qff\right) 
\off\qff \right|\off\fff
u\qff \in\qff E^{\dff +}
\pff\right\}
\]

\vspace{-27pt}
\[
\quad
=\off
\bigl\{\pff 
\bigl(\qff u\dff \cos\dff \theta/2\fff,\off 
-\qff i\trf \varphi\trf(\trf u\trf)\dff \sin\dff \theta/2\qff\bigr)
\pff\fff \bigl|\pff\fff
u\qff \in\qff E^{\dff +}
\pff\bigr\}
\dff,\quad
\theta\qff \in\qff
[\trf 0\fff,\qff \pi\trf]
\qff.
\]

\vspace{-12pt}\vspace{4.5pt}
If\trs the spaces\sss $E^{\dff +}$\sss and\sss $E^{\dff -}$\sss are both\sss identified\sss
with a vector space\sss $F$\dnsp,\oss
then\sss this path\dss is\dss equal\dss to\sss the path\sss
$\lambda\trf(\trf -\qff i\trf \varphi\trf,\off \theta\trf)\dff,\off
\theta\qff \in\qff
[\trf 0\fff,\qff \pi\trf]$\nnsp.\oss
This property\sss depends on\sss $\sigma\fff,\off \rho\fff,\off \tau$\sss
having\sss the standard\sss form\qss (\ref{sigma-rho-tau-phi}),\oss
but\sss not\sss on\sss the boundary condition\sss $N$\nnsp.\oss

\myuppar{Comparing\sss $N$ and\sss $N^{\dff \perp}$\dnsp.}
Let\sss us\sss identify\sss $E^{\dff +}$\sss
with\sss $E^{\dff -}$ by\sss $\psi$\nnsp.\oss
Then\dss $\psi\off =\off 1$\nnsp,\qss $N\off =\off \Delta$\nnsp,\oss 
and\dss the eigenvalues of\dss $\varphi$\sss are $i$ or $-\qff i$\nnsp.\oss
We can\sss bring\sss $\sigma\fff,\off \rho\fff,\off \tau\fff,\off N$\sss
into\sss the normal\dss form\sss in\sss the sense of\trs Section\qss \ref{boundary-algebra}\qss
by moving\sss the eigenvalue $i$ of\dss $\varphi$\sss to $1$ clockwise along\sss
the unit\sss circle in\sss $\ccc$\sss and\sss moving\sss the eigenvalue\sss $-\qff i$\sss
counterclockwise\sss to $1$\nnsp.\oss
In\sss more details,\oss let\vspace{3pt}
\[
\quad 
\varphi_{\dff \eta}\trf(\dff u\trf)
\off =\off
\varphi\dff(\dff u\trf)\qff e^{\qff -\dff i\dff \eta}
\quad
\mbox{for}\quad
u\qff \in\qff \mathcal{L}_{\dff +}\dff(\trf \varphi\trf)
\qff,\off
\]

\vspace{-33pt}
\[
\quad 
\varphi_{\dff \eta}\trf(\trf u\trf)
\off =\off
\varphi\dff(\dff u\trf)\qff e^{\dff i\dff \eta\phantom{\qff - }}
\quad
\mbox{for}\quad
u\qff \in\qff \mathcal{L}_{\dff -}\dff(\trf \varphi\trf)
\qff.\off
\]

\vspace{-12pt}\vspace{3pt}
The family\sss $\varphi_{\dff \eta}\dff,\off \eta\qff \in\qff [\trf 0\fff,\qff \pi/2\trf]$\dss
is\dss our deformation of\dss $\varphi$\dnsp.\oss
It\sss deforms\sss $\varphi$\sss into\sss the identity\sss map 
and\dss leads\sss to a deformation of\sss $\rho$ and $\tau$\dnsp,\oss
while $\sigma$ and\sss $N$\sss are not\sss moved.\oss
At\sss the end\sss $\eta\off =\off \pi/2$\sss of\dss this deformation\sss
$\mathcal{L}_{\dff +}\dff(\trf \varphi\trf)
\off =\off
\Delta\dff(\trf \varphi\trf)
\off =\off
\Delta\dff(\trf \id \trf)
\off =\off
N$\sss 
and\sss $\sigma\fff,\off \rho\fff,\off \tau$\sss have\sss the standard\sss
form\qss (\ref{sigma-rho-tau-st}).\oss
Therefore\sss $\sigma\fff,\off \rho\fff,\off \tau\fff,\off N$\sss
are normalized.\oss
The deformation\sss
$\varphi_{\dff -\qff \eta}\dff,\off \eta\qff \in\qff [\trf -\qff \pi/2\fff,\qff 0\trf]$\sss
moves in\sss the opposite direction and connects\sss $\id$\sss with\sss $\varphi$\nnsp.\oss

Similarly,\oss
we can\sss bring\sss $\sigma\fff,\off \rho\fff,\off \tau\fff,\off N^{\dff \perp}$\sss
into\sss the normal\dss form\sss 
by moving\sss the eigenvalue $i$ of\dss $\varphi$\sss to $-\qff 1$ counterclockwise 
and\sss moving\sss the eigenvalue\sss $-\qff i$\sss
clockwise\sss to $-\qff 1$\nnsp.\oss
More precisely,\oss this deformation\dss is\dss
$\varphi_{\dff -\qff \eta}\dff,\off \eta\qff \in\qff [\trf 0\fff,\qff \pi/2\trf]$\nnsp.\oss
Together\sss the deformations\sss
$\varphi_{\dff -\qff \eta}\dff,\off \eta\qff \in\qff [\trf 0\fff,\qff \pi/2\trf]$\sss
and\sss
$\varphi_{\dff -\qff \eta}\dff,\off \eta\qff \in\qff [\trf -\qff \pi/2\fff,\qff 0\trf]$\sss
define a deformation
$\varphi_{\dff -\qff \eta}\dff,\off \eta\qff \in\qff [\trf -\qff \pi/2\fff,\qff \pi/2\trf]$\sss
connecting\sss $\id$\sss with\sss $-\qff \id$\nnsp.\oss
It\dss is\dss convenient\sss to\sss reparameterize\sss this\sss deformation as\sss
$\omega_{\dff \eta}
\off =\off 
\varphi_{\dff -\qff (\trf \eta\qff -\qff \pi/2\trf)}\dff,\off 
\eta\qff \in\qff [\trf 0\fff,\qff \pi\trf]$\nnsp.\oss
Then\vspace{3pt}
\[
\quad 
\omega_{\dff \eta}\trf(\dff u\trf)
\off =\off
\varphi\dff(\dff u\trf)\qff e^{\dff i\trf (\trf \eta\qff -\qff \pi/2 \trf)\phantom{\qff - }}
\off =\off
-\qff i\trf \varphi\dff(\dff u\trf)\qff e^{\dff i\dff \eta}
\off =\off
u\qff e^{\dff i\dff \eta\phantom{\qff - }}
\quad
\mbox{for}\quad
u\qff \in\qff \mathcal{L}_{\dff +}\dff(\trf \varphi\trf)
\qff,
\]

\vspace{-33pt}
\[
\quad 
\omega_{\dff \eta}\trf(\dff u\trf)
\off =\off
\varphi\dff(\dff u\trf)\qff e^{\qff -\dff i\trf (\trf \eta\qff -\qff \pi/2 \trf)}
\off =\off
i\trf \varphi\dff(\dff u\trf)\qff e^{\qff -\dff i\dff \eta}
\off\hspace{0.2em} =\off
u\qff e^{\qff -\dff i\dff \eta}
\quad
\mbox{for}\quad
u\qff \in\qff \mathcal{L}_{\dff -}\dff(\trf \varphi\trf)
\qff.\off
\]

\vspace{-12pt}\vspace{3pt}
The space\sss
$\lambda\trf(\trf -\qff i\trf \omega_{\dff \eta}\trf,\off \theta\trf)$\sss
is\dss equal\dss to\sss the sum of\dss subspaces\vspace{4.5pt}
\[
\quad
\left\{\pff 
\left(\qff 
u\dff \cos\dff \theta/2\fff,\off\off\off
u\qff e^{\dff i\trf (\trf \eta\qff -\qff \pi/2 \trf)}
\trf \sin\dff \theta/2
\qff\right)
\pff\fff \bigl|\pff\dff
u\qff \in\qff \mathcal{L}_{\dff +}\dff(\trf \varphi\trf)
\pff\right\}
\quad
\mbox{and}\quad
\]

\vspace{-30pt}
\[
\quad
\left\{\pff 
\left(\qff 
u\dff \cos\dff \theta/2\fff,\off 
u\qff e^{\qff -\dff i\trf (\trf \eta\qff +\qff \pi/2 \trf)}
\trf \sin\dff \theta/2
\qff\right)
\pff\fff \bigl|\pff\dff
u\qff \in\qff \mathcal{L}_{\dff -}\dff(\trf \varphi\trf)
\pff\right\}
\qff.
\]

\vspace{-12pt}\vspace{4.5pt}
By comparing\sss this with\sss the similar\sss presentation of\dss
$\gamma\fff'\trf(\trf A\dff,\pff \eta\dff,\pff \theta\trf)$\nnsp,\oss
we see\sss that\vspace{3pt}
\begin{equation}
\label{plus-bott}
\quad
\lambda\trf(\trf -\qff i\trf \omega_{\dff \eta}\trf,\off \theta\trf)
\off =\off
\gamma\fff'\trf\left(\trf 
\mathcal{L}_{\dff +}\dff(\trf \varphi\trf)\dff,\pff \eta\dff,\pff \theta
\trf\right)
\qff.
\end{equation}

\vspace{-12pt}\vspace{3pt}
This also follows directly\sss from\sss the definitions.\oss

\myuppar{Canonical\dss homomorphisms\dss
$\beta\dff \colon\dff
K^{\dff i}\dff (\trf S\fff Y\trf)
\qff \ttoo\qff
K^{\dff i}\dff (\trf S\dff X\dff \cup\dff B\dff X_{\trf Y}\trf)$\nnsp,\qss
$i\off =\off 0\fff,\qff 1$\nnsp.}
Now we move\sss to\sss the framework of\trs Section\qss \ref{symbols-conditions}.\oss
Recall\dss that\sss $X$\sss is\dss a compact\sss manifold and\sss
$Y\off =\off \partial\dff X$\nnsp.\oss
Let\sss $I^{\qff 2}$\sss be\sss the rectangle\sss
$[\trf 0\fff,\qff \pi\trf]\dff \times\dff [\trf 0\fff,\qff 2\dff\pi\trf]$\sss
and\dss let\sss $\partial\trf I^{\qff 2}$\sss be\sss its boundary.\oss
The key element\sss in\sss the construction of\dss the homomorphism\sss $\beta$\sss
are\sss the\dss Bott\trs periodicity\sss isomorphisms\vspace{3pt}
\[
\quad
K^{\dff i}\dff (\trf S\fff Y\trf)
\qff \ttoo\qff
K^{\dff i}\dff \left(\qff S\fff Y\dff \times\dff I^{\qff 2},\off 
S\fff Y\dff \times\dff \partial\trf I^{\qff 2} \qff\right)
\qff,
\]

\vspace{-12pt}\vspace{3pt}
where\sss $i\off =\off 0\fff,\qff 1$\nnsp.\oss
Let\sss us\sss take\sss the product\sss
$S\fff Y\dff \times\dff [\trf 0\fff,\qff \pi \trf]$\sss
and collapse\sss $S\fff Y\dff \times\dff 0$\sss and\sss $S\fff Y\dff \times\dff \pi$\sss
into\sss two different\sss copies $Y_{\dff 0}$\sss and\sss $Y_{\dff \pi}$\sss of\dss $Y$\dnsp.\oss
The result\sss $\mathbb{S}\dff Y$ can\sss be\sss identified\sss with $S\dff X_{\trf Y}$\sss
using\sss the parameterization\sss 
$\nu_y\dff \cos\dff \theta
\pff +\qff
u\dff \sin\dff \theta$\sss
of\dss half-circles from\dss Section\qss \ref{symbols-conditions}.\oss
Naturally,\oss we agree\sss that\sss $Y_{\dff 0}$ corresponds\sss to points with\sss
$\theta\off =\off 0$\sss and\sss $Y_{\dff \pi}$\sss to points with\sss $\theta\off =\off \pi$\nnsp.\oss
The quotient\sss map\sss
$S\fff Y\dff \times\dff [\trf 0\fff,\qff \pi \trf]
\qff \ttoo\qff
\mathbb{S}\dff Y$\sss
together with\sss this identification\sss lead\dss to a continuous map\vspace{3pt}
\[
\quad
S\fff Y\dff \times\dff I^{\qff 2}
\off =\off
\bigl(\trf S\fff Y\dff \times\dff [\trf 0\fff,\qff \pi\trf]\trf\bigr)
\dff \times\dff
[\trf 0\fff,\qff 2\dff \pi\trf]
\qff \ttoo\qff
S\dff X_{\trf Y}\dff \times\dff [\trf 0\fff,\qff 2\dff \pi\trf]
\qff.
\]

\vspace{-12pt}\vspace{3pt}
There\dss is\dss also a natural\dss homeomorphism\vspace{3pt}
\[
\quad
S\dff X
\qff \cup\qff
\bigl(\trf S\dff X_{\trf Y}\dff \times\dff [\trf 0\fff,\qff 2\dff \pi\trf]\trf\bigr)
\qff \cup\qff 
\bigl(\trf B\dff X_{\trf Y}\dff \times\dff 2\dff \pi\trf\bigr)
\qff \ttoo\qff
S\dff X\dff \cup\dff B\dff X_{\trf Y}
\qff.
\]

\vspace{-12pt}\vspace{3pt}
The composition of\dss the\sss last\sss two maps\sss 
leads\sss to canonical\dss homomorphisms\vspace{3pt}
\begin{equation}
\label{inclusion}
\quad
K^{\dff i}\dff \left(\qff S\fff Y\dff \times\dff I^{\qff 2},\off 
S\fff Y\dff \times\dff \partial\trf I^{\qff 2} \qff\right)
\qff \ttoo\qff
K^{\dff i}\dff (\qff S\dff X\dff \cup\dff B\dff X_{\trf Y} \qff)
\qff,
\end{equation}

\vspace{-12pt}\vspace{3pt}
where\sss $i\off =\off 0\fff,\qff 1$\nnsp.\oss
The homomorphisms\sss $\beta$\sss are defined as\sss the compositions of\dss these homomorphisms
with\sss the\dss Bott\trs periodicity\sss isomorphisms.\oss
When everything depends on a parameter\sss $z\qff \in\qff Z$\nnsp,\oss
as at\sss the end of\trs Section\qss \ref{symbols-conditions},\oss
this construction\sss leads\sss to\sss canonical\dss homomorphisms\vspace{4.5pt}
\[
\quad
\beta\dff \colon\dff
K^{\dff i}\dff (\trf S\fff V\trf)
\qff \ttoo\qff
K^{\dff i}\dff (\trf S\dff W\dff \cup\dff B\dff W_{\trf V}\trf)
\qff,
\]

\vspace{-12pt}\vspace{4.5pt}
where\sss the spaces involved are as in\dss Section\qss \ref{symbols-conditions}.\oss

\myuppar{Anti-commuting\sss symbols and special\dss boundary conditions.}
In\sss the framework of\trs Section\qss \ref{symbols-conditions},\oss
let\sss $\sigma$\sss be a self-adjoint\sss symbol\sss
such\sss that\sss $\sigma_y$\sss is\dss unitary for every $y\qff \in\qff Y$\sss
and\sss $\sigma$\sss is\qss \emph{anti-commuting}\pss in\sss the sense\sss that\sss
$\sigma_y$\sss anti-commutes with $\tau_u$\sss for every\sss
$u\qff \in\qff S\fff Y$\sss and\sss $y\off =\off \pi\trf(\trf u\trf)$\nnsp.\oss
Let\sss $N$\sss be a self-adjoint\sss elliptic bundle-like boundary condition for $\sigma$\nnsp.\oss
Then\sss the subspace\sss $N_{\fff u}$\nsp,\oss where\sss $u\qff \in\qff S\fff Y$\dnsp,\oss depends only on\sss
$y\off =\off \pi\trf(\dff u\trf)$ and\sss we may denote it\sss by\sss $N_{\fff y}$\nsp.\oss 
We will\sss say\sss that\sss the boundary condition\sss $N$\sss for $\sigma$\sss is\qss
\emph{special}\oss if\dss $N_{\fff y}$\sss is\dss transverse\sss
not\sss only\sss to\sss 
$\mathcal{L}_{\dff -}\dff(\trf \rho_{\fff u}\dff)$\nnsp,\oss
but\sss also\sss to\sss
$\mathcal{L}_{\dff +}\dff(\trf \rho_{\fff u}\dff)$\sss
for every\sss $u\qff \in\qff S\fff Y$\sss and\sss
$y\off =\off \pi\trf(\dff u\trf)$\nnsp.\oss
If\dss $N$\sss is\dss special,\oss the anti-com\-mu\-ta\-tiv\-i\-ty\sss
property of\dss $\sigma$\sss implies\sss that\sss $N^{\dff \perp}$\sss is\dss
also a special\dss boundary condition\sss for $\sigma$\nnsp.\oss

Since\sss $N$\sss and\sss $N^{\dff \perp}$ are bundle-like boundary conditions for $\sigma$\nnsp,\oss
the invariants\sss $e^{\dff +}\fff (\trf \sigma\fff,\qff N\qff)$\sss and\sss
$\varepsilon^{\dff +}\fff (\trf \sigma\fff,\qff N\qff)$\nnsp,\oss
as also $e^{\dff +}\fff (\trf \sigma\fff,\qff N^{\dff \perp}\qff)$\sss
and\sss $\varepsilon^{\dff +}\fff (\trf \sigma\fff,\qff N^{\dff \perp}\qff)$ are well\sss defined.\oss
We would\dss like\sss to see how\sss they differ.\oss
It\sss turns out\sss that,\oss as one may expect,\oss the difference\dss is\dss determined\sss by\sss
the behavior of\dss
$\sigma$\sss and\sss $N$\sss 
at\sss the boundary\sss $Y$\dnsp.\oss

For\sss $u\qff \in\qff S\fff Y$\sss and\sss $y\off =\off \pi\trf(\trf u\trf)$\sss
let\sss $\psi_y$\sss and\sss $\varphi_{\fff u}$\sss be\sss the isometries\sss
$E^{\dff +}_{\dff y}\qff \ttoo\qff E^{\dff -}_{\dff y}$\sss
having as\sss their\sss graphs\sss $N_{\dff y}$\sss and\sss
$\mathcal{L}_{\dff +}\dff(\trf \rho_{\fff u}\dff)$\sss respectively.\oss
Then\dss
$\psi_y^{\dff -\dff 1}\dff \circ\trf \varphi_{\fff u}$\sss
is\dss an\sss isometry\sss $E^{\dff +}_{\dff y}\qff \ttoo\qff E^{\dff +}_{\dff y}$\nsp.\oss
Since\sss $N$\sss is\dss a special\dss boundary condition,\oss
this isometry\sss has\sss no eigenvalues equal\dss to $1$ or $-\qff 1$\nnsp.\oss
For every\sss $u\qff \in\qff S\fff Y$\sss let\sss us consider\sss subspaces\vspace{3pt}
\[
\quad
\mathcal{L}_{\dff +}\qff\left(\trf \psi_y^{\dff -\dff 1}\dff \circ\trf \varphi_{\fff u} \trf\right)
\quad
\mbox{and}\quad
\mathcal{L}_{\dff -}\qff\left(\trf \psi_y^{\dff -\dff 1}\dff \circ\trf \varphi_{\fff u} \trf\right)
\qff
\]

\vspace{-12pt}\vspace{3pt}
and\dss let\sss 
$\bm{\upsilon}_{\fff u}\dff \colon\dff
E^{\dff +}_{\dff y}\qff \ttoo\qff E^{\dff +}_{\dff y}$\sss
be\sss the operator
acting on\sss the first\sss subspace as\sss the identity\sss $\id$\sss
and on\sss the second one as\sss $-\qff \id$\nnsp.\oss
Clearly,\oss the operators $\bm{\upsilon}_{\fff u}$ are self-adjoint\sss
and\sss and\sss invertible,\oss and\dss hence define\sss a self-adjoint\sss
elliptic symbol\sss $\bm{\upsilon}$\nnsp,\oss
the\qss \emph{boundary\sss symbol}\qss $\bm{\upsilon}$\nnsp.\oss
Since\sss $Y$\sss is\dss a closed\sss manifold,\qss $\bm{\upsilon}$ requires no boundary conditions
and\sss hence\sss there are well\sss defined classes\vspace{4.5pt}
\[
\quad
e^{\dff +}\fff (\trf \bm{\upsilon}\trf)
\qff \in\pff 
K^{\dff 0}\dff (\trf S\fff Y\trf)
\quad
\mbox{and}\quad
\varepsilon^{\dff +}\fff (\trf \bm{\upsilon}\trf)
\qff \in\pff 
K^{\dff 1}\dff (\trf B\dff Y\fff,\qff S\fff Y\trf)
\qff.
\]

\vspace{-12pt}\vspace{4.5pt}
When a parameter\sss $z\qff \in\qff Z$\sss is\dss present,\oss
as at\sss the end of\trs Section\qss \ref{symbols-conditions},\oss
we get\sss classes\vspace{4.5pt}
\[
\quad
e^{\dff +}\fff (\trf \bm{\upsilon}\trf)
\qff \in\pff 
K^{\dff 0}\dff (\trf S\fff V\trf)
\quad
\mbox{and}\quad
\varepsilon^{\dff +}\fff (\trf \bm{\upsilon}\trf)
\qff \in\pff 
K^{\dff 1}\dff (\trf B\dff V\fff,\qff S\fff V\trf)
\qff,
\]

\vspace{-12pt}\vspace{4.5pt}
where we used\dss the notations from\dss Section\qss \ref{symbols-conditions}.\oss

\mypar{Theorem.}{basic-difference}
\emph{Under\sss the above assumptions}\qss
$e^{\dff +}\fff (\trf \sigma\fff,\qff N\qff)
\qff -\qff
e^{\dff +}\fff (\trf \sigma\fff,\qff N^{\dff \perp}\qff)
\off =\off\dff
\beta\dff \left(\trf
e^{\dff +}\fff (\trf \bm{\upsilon}\trf)
\qff\right)$\nnsp.\oss

\proof
Let\sss us consider\sss first\sss the case of\dss a single manifold\dss $X$\sss with\sss
the boundary\sss $Y$\dnsp.\oss
For every\sss $u\qff \in\qff S\fff Y$\sss and\sss $y\off =\off \pi\trf(\trf u\trf)$\sss
let\sss $\psi_{\fff y}$\sss and\sss $\varphi_{\fff u}$\sss be\sss the isometries\sss
$E^{\dff +}_{\dff y}\qff \ttoo\qff E^{\dff -}_{\dff y}$\sss
having as\sss their\sss graphs\sss $N_{\dff y}$\sss and\sss
$\mathcal{L}_{\dff +}\dff(\trf \rho_{\fff u}\dff)$\sss respectively.\oss
Let\sss us deform\sss $\sigma_y\dff,\off \rho_{\fff u}\dff,\off \tau_u$\sss
as explained at\sss the beginning of\dss this section
(actually,\oss by\sss the assumptions,\oss the first\sss deformation\dss is\dss not\sss needed\fff).\oss
These deformations can be arranged\sss to continuously depend on $u$\nnsp.\oss
Moreover,\oss they can be extended\sss to a continuous deformation of\dss the symbol\sss $\sigma$\nnsp.\oss 
At\sss the end of\dss these deformations\sss the isometries\sss
$\psi^{\dff -\dff 1}_{\fff y}\dff \circ\trf \varphi_{\fff u}$\sss 
have only\sss $i$ and $-\qff i$ as eigenvalues.\oss
In\sss general,\oss the corresponding eigenspaces depend\sss on $u$\nnsp,\oss
and\dss  the deformed symbol\sss $\sigma$\sss is\dss not\sss bundle-like.\oss
Let\sss us\sss identify\sss
$E^{\dff +}_{\dff y}$\sss with\sss $E^{\dff -}_{\dff y}$\sss
by\sss the isometry\sss $\psi_y$\sss for every\sss $y\qff \in\qff Y$\dnsp.\oss

Now we can\sss further deform\sss $\sigma_y\dff,\off \rho_{\fff u}\dff,\off \tau_u$\sss
and\sss bring\sss $\sigma_y\dff,\off \rho_{\fff u}\dff,\off \tau_u\dff,\oss N_{\dff y}$\sss
into a normal\sss form\sss in\sss the sense of\trs Section\qss \ref{boundary-algebra}\qss
by moving\sss the eigenvalues $i$ and $-\qff i$\sss to $1$ as explained above.\oss
For\sss the deformed symbol\sss 
$\varphi_{\fff u}\off =\off \id$\sss 
for every $u$\sss
and\sss hence\sss the deformed\sss symbol\dss is\dss bundle-like.\oss
Of\dss course,\oss we need\sss to extend\sss these deformations\sss to a deformation
of\dss the whole symbol\sss $\sigma$\sss over\sss $X$\nnsp.\oss
A convenient\sss way\sss to do\sss this\dss is\dss to identify\sss
$S\dff X$\sss with\sss the result\sss of\dss glueing\sss 
$S\dff X_{\trf Y}\dff \times\dff [\trf 0\fff,\qff \pi/2\trf]$\sss
to\sss $S\dff X$\sss along\sss
$S\dff X_{\trf Y}\dff \times\dff 0\off =\off S\dff X_{\trf Y}$\sss
and\sss placing\sss $\sigma$ before\sss this deformation on\sss $S\dff X$\sss
and\sss placing\sss the deformation on\sss 
$S\dff X_{\trf Y}\dff \times\dff [\trf 0\fff,\qff \pi/2\trf]$\nnsp.\oss
At\sss the end\sss we will\sss get\sss a bundle-like symbol,\oss
which we denote by\sss $\bm{\sigma}$\nnsp.\oss
Since\sss $\bm{\sigma}$\sss is\dss bundle-like,\oss
the bundle\sss $\mathbb{E}^{\dff +}\dff(\trf \bm{\sigma}\trf)$\sss
over\sss\vspace{4.5pt}
\[
\quad
S\dff X\dff \cup\dff B\dff X_{\trf Y}
\off =\off
S\dff X
\qff \cup\qff
\bigl(\qff
S\dff X_{\trf Y}\dff \times\dff [\trf 0\fff,\qff \pi/2\trf]
\qff\bigr)
\qff \cup\qff
\bigl(\qff
B\dff X_{\trf Y}\dff \times\dff \pi/2
\qff\bigr)
\]

\vspace{-12pt}\vspace{4.5pt}
is\dss well\sss defined,\oss and\sss the invariant\sss 
$e^{\dff +}\fff (\trf \sigma\fff,\qff N\qff)$\sss
is\dss the class of\dss this bundle.\oss\vspace{1.5pt}

Similarly,\oss we can bring\sss 
$\sigma_y\dff,\off \rho_{\fff u}\dff,\off \tau_u\dff,\oss N_{\dff y}^{\dff \perp}$\sss
to a normal\sss form\sss by moving\sss the eigenvalues $i$ and $-\qff i$\sss to $-\qff 1$\nnsp.\oss 
For\sss the deformed symbol\sss $\varphi_{\fff u}\off =\off -\qff \id$\sss for every $u$\sss
and\sss hence\sss it\dss is\dss bundle-like.\oss
As above,\oss we will\dss place\sss the deformation on\sss
$S\dff X_{\trf Y}\dff \times\dff [\trf 0\fff,\qff \pi/2\trf]$\sss
and\sss get\sss another\sss bundle-like symbol,\oss
which we denote by\sss $\bm{\sigma}^{\dff \perp}$\dnsp.\oss
The invariant\sss
$e^{\dff +}\fff (\trf \sigma\fff,\qff N^{\dff \perp}\qff)$\sss
is\sss the class of\dss the bundle\sss
$\mathbb{E}^{\dff +}\dff(\trf \bm{\sigma}^{\dff \perp}\trf)$\nnsp.\oss
The bundles\sss
$\mathbb{E}^{\dff +}\dff(\trf \bm{\sigma}\trf)$\sss
and\sss
$\mathbb{E}^{\dff +}\dff(\trf \bm{\sigma}^{\dff \perp}\trf)$\sss
differ only over\vspace{4.5pt}
\begin{equation}
\label{cylinder}
\quad
\bigl(\qff
S\dff X_{\trf Y}\dff \times\dff [\trf 0\fff,\qff \pi/2\trf]
\qff\bigr)
\qff \cup\qff
\bigl(\qff
B\dff X_{\trf Y}\dff \times\dff \pi/2
\qff\bigr)
\qff.
\end{equation}

\vspace{-12pt}\vspace{4.5pt}
In order\sss to compare\sss them,\oss 
let\sss us\sss construct\sss an\dss intermediate,\oss a bundle\sss
$\mathbf{E}^{\dff +}\dff(\trf \bm{\sigma}\trf)$ 
over\vspace{4.5pt}
\begin{equation}
\label{long-cylinder}
\quad
\bigl(\qff
S\dff X_{\trf Y}\dff \times\dff [\trf 0\fff,\qff \pi/2\qff +\qff 2\dff \pi\trf]
\qff\bigr)
\qff \cup\qff
\bigl(\qff
B\dff X_{\trf Y}\dff \times\dff [\trf \pi/2\qff +\qff \pi\fff,\qff \pi/2\qff +\qff 2\dff \pi\trf]
\qff\bigr)
\qff.
\end{equation}

\vspace{-12pt}\vspace{4.5pt}
Over\sss the subset\qss (\ref{cylinder})\qss
we set\sss the bundle\sss
$\mathbf{E}^{\dff +}\dff(\trf \bm{\sigma}\trf)$\sss
to be equal\dss to\sss
$\mathbb{E}^{\dff +}\dff(\trf \bm{\sigma}\trf)$\nnsp.\oss
In order\sss to define\sss
$\mathbf{E}^{\dff +}\dff(\trf \bm{\sigma}\trf)$\sss
over\sss
$S\dff X_{\trf Y}
\dff \times\dff 
[\trf \pi/2\fff,\qff \pi/2\qff +\qff 2\dff \pi\trf]$\nnsp,\oss
recall\dss that\sss we identified\sss $S\dff X_{\trf Y}$\sss
with\sss the quotient\sss $\mathbb{S}\dff Y$ of\dss 
$S\fff Y\dff \times\dff [\trf 0\fff,\qff \pi \trf]$\nnsp.\oss
For every\sss $u\qff \in\qff S\fff Y$\sss and\sss 
$\theta\qff \in\qff [\trf 0\fff,\qff \pi\trf]$\nnsp,\qss
$\eta\qff \in\qff [\trf \pi/2\fff,\qff \pi/2\qff +\qff 2\dff \pi\trf]$\nnsp,\oss
let\sss\vspace{3pt}
\[
\quad
\gamma\fff'\trf\left(\trf 
\mathcal{L}_{\dff +}\dff(\trf \varphi_{\fff u}\trf)\dff,\pff \eta\qff -\qff \pi/2\dff,\pff \theta
\trf\right)
\qff
\]

\vspace{-12pt}\vspace{3pt}
be\sss the fiber of\dss 
$\mathbf{E}^{\dff +}\dff(\trf \bm{\sigma}\trf)$ 
over\sss the point\sss represented\dss by\sss
$(\trf u\fff,\qff \theta\fff,\qff \eta \fff\qff)$\nnsp.\oss
Here\sss $\varphi_{\fff u}$\sss refers\sss to\sss the symbol\sss
before\sss the\sss last\sss deformations.\oss
By\sss the equality\qss (\ref{plus-bott})\qss and\sss
the description of\dss the positive subspaces,\oss
for\sss $\eta\off =\off \pi/2$\sss these fibers are\sss the same
as\sss the fibers of\dss
$\mathbb{E}^{\dff +}\dff(\trf \bm{\sigma}\trf)$\nnsp.\oss
In\sss particular,\oss the bundle\sss 
$\mathbf{E}^{\dff +}\dff(\trf \bm{\sigma}\trf)$\sss
is\dss correctly defined over\sss
$S\dff X_{\trf Y}
\dff \times\dff 
[\trf 0\fff,\qff \pi/2\qff +\qff 2\dff \pi\trf]$\nnsp.\oss
For\sss
$\eta
\qff \in\qff 
[\trf \pi/2\qff +\qff \pi\fff,\qff \pi/2\qff +\qff 2\dff \pi\trf]$\sss
these fibers do not\sss depend on $u$\sss and\dss hence\sss 
the family of\dss these fibers can\sss be extended\dss to\sss
$B\dff X_{\trf Y}
\dff \times\dff 
[\trf \pi/2\qff +\qff \pi\fff,\qff \pi/2\qff +\qff 2\dff \pi\trf]$\sss
in\sss the same way as\sss in\dss Section\qss \ref{symbols-conditions}\qss 
the bundles\sss $E^{\dff +}\fff (\trf \sigma\trf)$\sss
were extended\dss to $B\dff X_{\trf Y}$.\oss
The resulting\sss bundle over\sss the space\qss (\ref{long-cylinder})\qss is\dss the promised\sss bundle\sss
$\mathbf{E}^{\dff +}\dff(\trf \bm{\sigma}\trf)$\nnsp.\vspace{-0.125pt}

By\sss the construction,\oss the bundle\sss 
$\mathbf{E}^{\dff +}\dff(\trf \bm{\sigma}\trf)$
over\sss 
$S\dff X_{\trf Y}\dff \times\dff [\trf 0\fff,\qff \pi/2\trf]$\sss
is\dss equal\sss to\sss the bundle\sss 
$\mathbb{E}^{\dff +}\dff(\trf \bm{\sigma}\trf)$\nnsp,\oss
and over\sss
$S\dff X_{\trf Y}\dff \times\dff [\trf \pi/2\fff,\qff \pi\trf]$\sss
is\dss equal\dss to a\qss ``reflected''\qss copy\sss 
of\trs the restriction of\dss the bundle\sss
$\mathbb{E}^{\dff +}\dff(\trf \bm{\sigma}\trf)$\sss 
to\sss
$S\dff X_{\trf Y}\dff \times\dff [\trf 0\fff,\qff \pi/2\trf]$\nnsp.\oss
Hence over\sss
$S\dff X_{\trf Y}\dff \times\dff [\trf 0\fff,\qff \pi\trf]$\sss 
we can deform\sss the bundle\sss
$\mathbf{E}^{\dff +}\dff(\trf \bm{\sigma}\trf)$\sss
to\vspace{4.4pt}
\[
\quad
\left(\qff \left. \mathbb{E}^{\dff +}\dff(\trf \bm{\sigma}\trf)\trf \right|\qff S\dff X_{\trf Y}\qff\right)
\dff \times\dff
[\trf 0\fff,\qff \pi\trf]
\]

\vspace{-12pt}\vspace{4.4pt}
without\sss affecting\dss it\sss over\sss
$S\dff X_{\trf Y}\dff \times\dff 0$\sss
and\sss
$S\dff X_{\trf Y}\dff \times\dff \pi$\nnsp.\oss
Instead of\dss this,\oss we can
simply cut\sss out\sss the subspace 
$S\dff X_{\trf Y}\dff \times\dff [\trf 0\fff,\qff \pi\trf]$
and\sss the bundle over\sss it.\oss
This will\sss replace\sss the interval\sss
$[\trf 0\fff,\qff \pi/2\qff +\qff 2\dff \pi\trf]$\sss
by\sss the shorter\sss interval\sss
$[\trf 0\fff,\qff \pi/2\qff +\qff \pi\trf]$\sss
and\dss result\dss in a bundle over\sss
$S\dff X_{\trf Y}\dff \times\dff [\trf 0\fff,\qff \pi/2\qff +\qff \pi\trf]$\nnsp.\oss
This bundle\dss is\dss equal\dss to\sss
$\mathbb{E}^{\dff +}\dff(\trf \bm{\sigma}^{\dff \perp}\trf)$\sss
over\sss the subspace\qss (\ref{cylinder}).\oss
After a reparameterization\sss replacing\sss the\sss interval\sss
$[\trf 0\fff,\qff \pi/2\qff +\qff \pi\trf]$\sss
by\sss the\sss interval\sss
$[\trf 0\fff,\qff \pi/2\trf]$\nnsp,\oss
the restriction of\dss this new bundle\sss to\vspace{4.4pt}
\[
\quad
\bigl(\qff
S\dff X_{\trf Y}\dff \times\dff [\trf 0\fff,\qff \pi/2\qff +\qff \pi\trf]
\qff\bigr)
\qff \cup\qff
\bigl(\trf
B\dff X_{\trf Y}\dff \times\dff (\trf \pi/2\qff +\qff \pi\trf)
\qff\bigr)
\qff
\]

\vspace{-12pt}\vspace{4.4pt}
may\sss be considered as another extension\sss of\dss
$E^{\dff +}\dff(\trf \bm{\sigma}^{\dff \perp}\trf)$ from\sss
$S\dff X_{\trf Y}\dff \times\dff [\trf 0\fff,\qff \pi/2\trf]$\sss
to\sss $B\dff X_{\trf Y}\dff \times\dff \pi/2$\nnsp.\oss
Since\sss our new bundle\dss is\dss also defined over\sss
$B\dff X_{\trf Y}\dff \times\dff [\trf \pi/2\fff,\qff \pi/2\qff +\qff \pi\trf]$\nnsp,\oss
the\sss two extensions are\sss isomorphic.\oss
By\sss placing\sss back\sss the cut\sss out\sss piece,\oss
we see\sss that\sss the restriction\vspace{4.4pt}
\[
\quad
\widetilde{\mathbb{E}}^{\dff +}\dff(\trf \bm{\sigma}^{\dff \perp}\trf)
\off =\off
\mathbf{E}^{\dff +}\dff(\trf \bm{\sigma}\trf)
\qff \left|\pff
\bigl(\qff
S\dff X_{\trf Y}\dff \times\dff [\trf 0\fff,\qff \pi/2\qff +\qff 2\dff \pi\trf]
\qff\bigr)
\qff \cup\qff
\bigl(\trf
B\dff X_{\trf Y}\dff \times\dff (\trf \pi/2\qff +\qff 2\dff \pi\trf)
\qff\bigr)
\right.
\qff
\]

\vspace{-12pt}\vspace{4.4pt}
is\dss isomorphic\sss to\sss
$\mathbb{E}^{\dff +}\dff(\trf \bm{\sigma}^{\dff \perp}\trf)$ 
after a reparameterization\sss replacing\sss
$[\trf 0\fff,\qff \pi/2\qff +\qff 2\dff \pi\trf]$\sss
by\sss
$[\trf 0\fff,\qff \pi/2\trf]$\nnsp.\oss 

Let\sss us now replace in\sss 
$\mathbf{E}^{\dff +}\dff(\trf \bm{\sigma}\trf)$\sss
the part\sss over\sss
$S\dff X_{\trf Y}\dff \times\dff [\trf \pi/2\fff,\qff \pi/2\qff +\qff 2\dff \pi\trf]$\sss
by\sss the product\sss\vspace{4.4pt}
\[
\quad
\left(\qff \left.
\mathbf{E}^{\dff +}\dff(\trf \bm{\sigma}\trf)\sss 
\dff(\trf \bm{\sigma}\trf)\qff \right|\qff S\dff X_{\trf Y}\dff \times\dff \pi/2
\qff\right)
\dff \times\dff
[\trf \pi/2\fff,\qff \pi/2\qff +\qff 2\dff \pi\trf]
\qff,
\]

\vspace{-12pt}\vspace{4.4pt}
and\dss let\sss
$\widetilde{\mathbb{E}}^{\dff +}\dff(\trf \bm{\sigma}\trf)$\sss
be\sss the resulting\sss bundle.\oss
If\dss we cut\sss out\sss this product\sss from\sss
$\widetilde{\mathbb{E}}^{\dff +}\dff(\trf \bm{\sigma}\trf)$\nnsp,\oss
we will\dss get\sss the bundle\sss
$\mathbb{E}^{\dff +}\dff(\trf \bm{\sigma}\trf)$\nnsp.\oss
By\sss placing\sss back\sss the cut\sss out\sss piece,\oss
we see\sss that\sss
$\widetilde{\mathbb{E}}^{\dff +}\dff(\trf \bm{\sigma}\trf)$\sss
is\dss isomorphic,\oss after a reparameterization\sss replacing\sss
$[\trf 0\fff,\qff \pi/2\qff +\qff 2\dff \pi\trf]$\sss
by\sss
$[\trf 0\fff,\qff \pi/2\trf]$\nnsp,\oss 
to\sss 
$\mathbb{E}^{\dff +}\dff(\trf \bm{\sigma}\trf)$\nnsp.\oss
Denoting\sss by\sss $[\trf F\trf]$\sss the $K$\dnsp-theory class of\dss a vector bundle $F$\dnsp,\oss
we see\sss that\vspace{4.4pt}
\begin{equation}
\label{bott-difference}
\quad
\left[\qff\fff \mathbb{E}^{\dff +}\dff(\trf \bm{\sigma}\trf) \qff\right]
\qff -\qff
\left[\qff\fff \mathbb{E}^{\dff +}\dff(\trf \bm{\sigma}^{\dff \perp}\trf) \qff\right]
\off =\off
\left[\pff \widetilde{\mathbb{E}}^{\dff +}\dff(\trf \bm{\sigma}\trf) \qff\right]
\qff -\qff
\left[\pff \widetilde{\mathbb{E}}^{\dff +}\dff(\trf \bm{\sigma}^{\dff \perp}\qff) \trf\right]
\qff.
\end{equation}

\vspace{-12pt}\vspace{4.4pt}
The bundles\dss
$\widetilde{\mathbb{E}}^{\dff +}\dff(\trf \bm{\sigma}\trf)$
and\dss
$\widetilde{\mathbb{E}}^{\dff +}\dff(\trf \bm{\sigma}^{\dff \perp}\trf)$\sss
differ only over\sss
$S\dff X_{\trf Y}\dff \times\dff [\trf \pi/2\fff,\qff \pi/2\qff +\qff 2\dff \pi\trf]$\nnsp,\oss
and\sss their\sss fibers over\sss over\sss the point\sss represented\dss by\sss
$(\trf u\fff,\qff \theta\fff,\qff \eta \fff\qff)$\sss are\vspace{3pt}
\[
\quad
\lambda\trf(\trf -\qff i\trf \id\trf,\off \theta\trf)
\quad
\mbox{and}\quad\dff
\gamma\fff'\trf\left(\trf 
\mathcal{L}_{\dff +}\dff(\trf \varphi_{\fff u}\trf)\dff,\pff \eta\qff -\qff \pi/2\dff,\pff \theta
\trf\right)
\qff
\]

\vspace{-12pt}\vspace{3pt}
respectively,\oss
and\dss these fibers are equal\dss for\sss $\eta\off =\off \pi/2$\sss
and\sss $\eta\off =\off \pi/2\qff +\qff 2\dff \pi$\nnsp.\oss
Let\sss us rename\sss $\eta\qff -\qff \pi/2$\sss as\sss $\eta$\nnsp.\oss
Then applying\sss the difference construction\sss to\sss two bundles over\sss\vspace{3pt}
\[
\quad
S\dff X_{\trf Y}\dff \times\dff [\trf 0\fff,\qff 2\dff \pi\trf]
\off =\off
S\dff Y\dff \times\dff I^{\qff 2}
\]

\vspace{-12pt}\vspace{3pt}
with\sss these fibers and\sss their identity\sss isomorphism over\sss
$S\dff Y\dff \times\dff \partial\trf I^{\qff 2}$\sss
leads\sss to\sss an element\dss\vspace{3pt}
\[
\quad
D\qff \in\qff
K^{\dff 0}\dff \left(\qff S\fff Y\dff \times\dff I^{\qff 2},\off 
\mathbb{S}\dff Y\dff \times\dff \partial\trf I^{\qff 2} \qff\right)
\qff.
\]

\vspace{-12pt}\vspace{3pt}
By\sss the construction\sss the image of\dss $D$\sss under\sss the homomorphism\qss
(\ref{inclusion})\qss is\dss equal\dss to\sss the difference\qss (\ref{bott-difference})\qss
and\dss hence\sss to\sss
$e^{\dff +}\fff (\trf \sigma\fff,\qff N\qff)
\qff -\qff
e^{\dff +}\fff (\trf \sigma\fff,\qff N^{\dff \perp}\qff)$\nnsp.\oss

At\sss the same\sss time\sss $D$\sss is\dss equal\sss to\sss image under\sss the\dss
Bott\dss periodicity\sss map of\dss the $K$\dnsp-theory class
of\dss the vector bundle\sss $\mathcal{L}_{\dff +}\dff(\trf \bullet\trf)$
over\sss $S\fff Y$\sss
having\sss the vector space\sss
$\mathcal{L}_{\dff +}\dff(\trf \varphi_{\fff u}\trf)$\sss
as\sss the fiber over\sss $u\qff \in\qff S\fff Y$\dnsp.\oss
By\sss the definition,\pss
$\mathcal{L}_{\dff +}\dff(\trf \varphi_{\fff u}\trf)$\sss
is\dss the only\sss eigen\-space of\dss $\bm{\upsilon}_{\fff u}$\sss
with\sss positive eigenvalue.\oss
It\sss follows\sss that\sss 
$\mathcal{L}_{\dff +}\dff(\trf \bullet\trf)
\off =\off
\mathbb{E}^{\dff +}\dff(\trf \bm{\upsilon}\trf)$\nnsp,\oss
and\sss hence\sss $e^{\dff +}\fff (\trf\bm{\upsilon}\trf)$\sss
is\dss equal\dss to\sss the $K$\dnsp-theory class of\dss
$\mathcal{L}_{\dff +}\dff(\trf \bullet\trf)$\nnsp.\oss
This implies\sss that\sss
$e^{\dff +}\fff (\trf \sigma\fff,\qff N\qff)
\qff -\qff
e^{\dff +}\fff (\trf \sigma\fff,\qff N^{\dff \perp}\qff)$\sss
is\dss equal\sss to\sss the result\sss of\dss applying\sss to\sss
$e^{\dff +}\fff (\trf \bm{\upsilon}\trf)$\sss
first\sss the\dss Bott\dss periodicity\sss map and\sss then\sss
the homomorphism\qss (\ref{inclusion}),\oss
i.e.\qss applying\dss $\beta$\nnsp.\oss

This proves\sss the\sss theorem\sss in\sss the case of\dss a single manifold\sss $X$\nnsp.\oss
Including\sss parameters\sss into\sss this proof\trs is\dss a routine matter.\oss  \eproof

\mypar{Corollary.}{k1-difference}
\emph{Under\sss the same assumptions}\qss
$\varepsilon^{\dff +}\fff (\trf \sigma\fff,\qff N\qff)
\qff -\qff
\varepsilon^{\dff +}\fff (\trf \sigma\fff,\qff N^{\dff \perp}\qff)
\off =\off\dff
\beta\dff \left(\trf
\varepsilon^{\dff +}\fff (\trf \bm{\upsilon}\trf)
\qff\right)$\nnsp.\oss

\proof
It\dss sufficient\sss to note\sss that\sss the\dss Bott\trs periodicity\sss
commutes\sss with\sss the coboundary\sss maps in\sss $K$\dnsp-theory.\oss \eproof

\mypar{Theorem.}{t-index-difference}
\emph{Under\sss the same assumptions}\oss
$\ti\qff (\trf \sigma\fff,\qff N\qff)
\qff -\qff
\ti\qff (\trf \sigma\fff,\qff N^{\dff \perp}\qff)
\off =\off\dff
\ti\qff (\trf \bm{\upsilon}\trf)$\nnsp.\oss

\proof
Of\dss course,\oss this\dss is\dss non-trivial\sss only\sss for\sss families\sss
$\sigma\trf(\trf z\trf)\fff,\off N\trf(\trf z\trf)\dff,\off z\qff \in\qff Z$\nnsp.\oss
In\sss the construction of\dss the\sss topological\dss indices of\dss
$(\trf \sigma,\pff N\qff)\dff,\off (\trf \sigma,\pff N^{\dff \perp}\qff)$\sss
and\sss $\bm{\upsilon}$\sss passing\sss to elements\vspace{3pt}
\[
\quad
t\trf(\trf \sigma,\pff N\qff)\dff,\off
t\trf(\trf \sigma,\pff N^{\dff \perp}\qff)\dff,\off
t\trf(\trf \bm{\upsilon}\qff)
\off \in\off
K^{\dff 1}\dff (\trf S^{\dff 2 n}\dff \times\dff Z\trf)
\]

\vspace{-12pt}\vspace{3pt}
absorbs homomorphisms\qss (\ref{inclusion})\qss and\dss they disappear.\oss
The\sss next\sss and\dss last\sss step\sss in\sss the construction of\dss the\sss
topological\dss index\dss is\dss an application of\dss the\dss Bott\trs
periodicity\sss map,\oss or,\oss rather,\oss 
in\sss the context\sss of\dss this section,\oss 
of\dss its\sss inverse.\oss
This step cancels\sss the application of\dss the\dss Bott\trs periodicity\sss
map\sss in\sss the definition of\dss $\beta$\nnsp.\oss
Therefore\sss the\sss theorem\sss follows\sss from\dss
Corollary\qss \ref{k1-difference}.\oss  \eproof

\newpage
\mysection{Dirac-like\qss boundary\qss problems}{odd}

\myuppar{Odd operators.}
Let\sss us return\sss to\sss the framework of\trs Section\qss \ref{abstract-index}\qss
and\sss prove an abstract\sss vanishing\sss theorem\sss
based on\sss ideas of\trs Gorokhovsky\sss and\dss Lesch\qss \cite{gl}.\oss
See\dss Theorem\qss \ref{odd-families}\qss below.\oss 

Suppose\sss that\sss the\dss Hilbert\sss spaces from\dss Section\qss \ref{abstract-index}\qss 
have\sss the form\vspace{3pt}\vspace{-0.4pt}
\[
\quad
H_{\trf 0}
\off =\off
K_{\trf 0}\dff \oplus\dff K_{\trf 0}\qff,\quad
H_{\dff 1}
\off =\off 
K_{\dff 1}\dff \oplus\dff K_{\dff 1}\qff,\quad
H^{\dff \partial}
\off =\off 
K^{\dff \partial}\dff \oplus\dff K^{\dff \partial}\qff,\quad
H_{\dff 1/2}^{\dff \partial}
\off =\off 
K_{\dff 1/2}^{\dff \partial}\dff \oplus\dff K_{\dff 1/2}^{\dff \partial}
\]

\vspace{-12pt}\vspace{3pt}\vspace{-0.4pt}
for some\dss Hilbert\dss spaces\sss
$K_{\trf 0}\dff,\off
K_{\dff 1}\dff,\off
K^{\dff \partial},\off\dff
K_{\dff 1/2}^{\dff \partial}$,\oss
and\dss that\sss $\gamma$\sss is\dss equal\dss to\sss the direct\sss sum
of\dss two copies of\dss a map\sss
$K_{\fff 1}\qff \ttoo\qff K_{\dff 1/2}^{\dff \partial}$,\oss
which we will\sss also denote by\sss $\gamma$\nnsp.\oss
Suppose\sss that\sss with respect\sss to\sss these orthogonal\sss decompositions\sss
the operators $A$ and\sss $\Sigma$\sss have\sss the form\vspace{1.5pt}\vspace{-0.4pt}
\[
\quad
A
\off =\off\dff
\begin{pmatrix}
\off 0 \dff &
B\fff' \off
\vspace{4.5pt} \\
\off\dff B &
0 \off 
\end{pmatrix}
\quad
\mbox{and}\quad
\Sigma
\off =\off\dff
\begin{pmatrix}
\off 0 &
1 \qff\off
\vspace{4.5pt} \\
\off\dff 1 &
0 \qff\off 
\end{pmatrix}
\off,
\]

\vspace{-12pt}\vspace{1.5pt}\vspace{-0.4pt}
for some operators\sss
$B\fff,\qff B\fff'\dff \colon\dff K_{\trf 0}\qff \ttoo\qff K_{\trf 0}$\nsp,\oss
and\dss that\sss the abstract\trs Lagrange\dss identity\qss (\ref{lagrange})\qss holds.\oss 
Suppose\sss that\sss $\Pi$\sss is\dss the orthogonal\sss projection\sss in\sss
$H^{\dff \partial}
\off =\off 
K^{\dff \partial}\dff \oplus\dff K^{\dff \partial}$\sss
onto\sss the graph of\dss a bounded skew-adjoint\sss operator\sss
$g\dff \colon\dff K^{\dff \partial}\qff \ttoo\qff K^{\dff \partial}$\dss
leaving\sss $K_{\dff 1/2}^{\dff \partial}$ invariant.\oss
Since $g$\sss is\dss skew-adjoint,\vspace{3pt}\vspace{-0.4pt}
\[
\quad
\sco{\trf
\Sigma\trf(\trf a\fff,\qff g\trf(\trf a\trf)\qff)\fff,\pff
(\trf b\fff,\qff g\trf(\trf b\trf)\qff) \dff}_{\dff \partial}
\off =\off
\sco{\dff
(\trf g\trf(\trf a\trf)\fff,\qff a\qff)\fff,\pff
(\trf b\fff,\qff g\trf(\trf b\trf)\qff) \dff}_{\dff \partial}
\]

\vspace{-34.5pt}\vspace{-0.4pt}
\[
\quad
=\off
\sco{\dff
g\trf(\trf a\trf)\fff,\pff
b \dff}_{\dff \partial}
\qff +\qff
\sco{\dff
a\fff,\pff
g\trf(\trf b\trf)\dff}_{\dff \partial}
\off =\off
0
\qff
\]

\vspace{-12pt}\vspace{3pt}\vspace{-0.4pt}
for every\sss $a\fff,\qff b\qff \in\qff K^{\dff \partial}$\dnsp.\oss
Hence\sss $\image\dff \Pi$\sss
is\dss orthogonal\dss to\sss $\Sigma\trf(\qff \image\dff \Pi \qff)$ and,\oss equivalently\nnsp,\oss
is\dss contained\sss in\sss $\kernel\dff \Pi$\nnsp.\oss
Conversely,\oss if\dss $(\trf u\fff,\qff v\trf)$\sss is\dss
orthogonal\dss to $\Sigma\trf(\qff \image\dff \Pi \qff)$\nnsp,\oss 
then\vspace{3pt}\vspace{-0.4pt}
\[
\quad
0
\off =\off
\sco{\trf
\Sigma\trf(\trf a\fff,\qff g\trf(\trf a\trf)\qff)\fff,\pff
(\trf u\fff,\qff v\qff) \dff}_{\dff \partial}
\off =\off
i\trf
\sco{\dff
g\trf(\trf a\trf)\fff,\pff
u \dff}
\qff +\qff
i\trf
\sco{\dff
a\fff,\pff
v\dff}_{\dff \partial}
\]

\vspace{-34.5pt}\vspace{-0.4pt}
\[
\quad
=\off
\sco{\dff
a\fff,\pff
g^{\dff *}\trf(\trf u\trf) \dff}_{\dff \partial}
\qff +\qff
\sco{\dff
a\fff,\pff
v\dff}_{\dff \partial}
\off =\off
\sco{\dff
a\fff,\pff
-\qff 
g\trf(\trf u\trf) \dff}_{\dff \partial}
\qff +\qff
\sco{\dff
a\fff,\pff
v\dff}_{\dff \partial}
\off =\off
\sco{\dff
a\fff,\pff
v
\qff -\qff 
g\trf(\trf u\trf) \dff}_{\dff \partial}
\qff
\]

\vspace{-12pt}\vspace{3pt}\vspace{-0.4pt}
for every $a$ and\dss hence\sss $v\off =\off g\trf(\trf u\trf)$\nnsp.\oss
It\sss follows\sss that\sss
$\Sigma\trf(\qff \image\dff \Pi  \qff)
\off =\off
\kernel\dff \Pi$\sss
and\dss hence\sss $A\fff,\pff \Pi$\sss is\dss a self-adjoint\dss boundary problem.\oss
Suppose\sss that\sss
$A\fff,\pff \Pi$\sss satisfies all\sss assumptions of\trs Section\qss \ref{abstract-index}.\oss
Let\sss $\Gamma$\sss be\sss the boundary operator associated\sss with\sss $\Pi$\nnsp.\oss
Finally,\oss let\sss $\varepsilon$\sss be\sss the automorphism of\dss $H_{\trf 0}
\off =\off
K_{\trf 0}\dff \oplus\dff K_{\trf 0}$\sss 
defined\dss by\sss the matrix\vspace{1.5pt}
\[
\quad
\varepsilon
\off =\off\dff
\begin{pmatrix}
\off\dff 1 &
0 \off
\vspace{4.5pt} \\
\off\dff 0 &
-\qff 1 \off 
\end{pmatrix}
\off.
\]

\vspace{-12pt}\vspace{1.5pt}
\mypar{Theorem.}{abstract-vanishing}
\emph{If\qss the operator\sss $i\fff g$\sss is\dss positive definite,\oss
then\sss 
$(\trf A\qff +\qff \varepsilon \trf)_{\trf \Gamma}
\dff \colon\dff 
\kernel\fff \Gamma\qff \ttoo\qff H_{\trf 0}$\qss
is\dss an\dss isomorphism.\oss
If\qss $i\fff g$\sss is\dss negative definite,\oss then\sss
$(\trf A\qff -\qff \varepsilon \trf)_{\trf \Gamma}$\sss is\dss an\dss isomorphism.\oss}

\proof
Let\sss $P\off =\off A\qff +\qff \varepsilon$\nnsp.\oss
Since\sss $\varepsilon$\sss is\dss self-adjoint,\oss
(\ref{lagrange})\sss
implies\sss the\trs Lagrange\dss identity\vspace{3pt}
\[
\quad
\sco{\dff P\dff u\fff,\qff v \dff}_{\dff 0}
\qff -\qff
\sco{\dff u\dff,\qff P\dff v \dff}_{\dff 0}
\off =\off
\sco{\dff i\trf \Sigma\dff \gamma\dff u\dff,\qff \gamma\dff v \dff}_{\dff \partial}
\qff
\]

\vspace{-12pt}\vspace{3pt}
for\sss $P$\dnsp.\oss
Let\sss
$u\off =\off (\trf x\fff,\qff y\trf)$\sss
and\sss 
$v\off =\off (\trf a\fff,\qff b\trf)$\nnsp,\oss
where\sss 
$x\fff,\off y\fff,\off a\fff,\off b
\qff \in\qff K_{\dff 1}$\nsp.\oss
Then\dss\vspace{3pt}
\[
\quad
\bsco{\dff i\trf \Sigma\dff \gamma\dff u\dff,\qff \gamma\dff v \dff}_{\dff \partial}
\off =\off
\bsco{\dff i\trf (\trf \gamma\dff y\fff,\qff \gamma\dff x\trf)\dff,\qff 
(\trf\gamma\dff a\fff,\qff \gamma\dff b\trf) \dff}_{\dff \partial}
\off =\off
i\trf \bsco{\dff \gamma\dff y\dff,\qff 
\gamma\dff a \dff}_{\dff \partial}
\pff +\pff
i\trf \bsco{\dff \gamma\dff x\dff,\qff 
\gamma\dff b \dff}_{\dff \partial}
\off,
\]

\vspace{-33pt}\vspace{-0.5pt}
\[
\quad
\sco{\dff P\dff u\fff,\qff v \dff}_{\dff 0}
\qff -\qff
\sco{\dff u\dff,\qff P\dff v \dff}_{\dff 0}
\off =\off
\sco{\dff B\fff'\fff y\fff,\qff a \dff}_{\dff 0}
\qff +\qff
\sco{\dff B\trf x\fff,\qff b \dff}_{\dff 0}
\qff -\qff
\sco{\dff x\dff,\qff B\fff'\trf b \dff}_{\dff 0}
\qff -\qff
\sco{\dff y\dff,\qff B\trf a \dff}_{\dff 0}
\off.
\]

\vspace{-12pt}\vspace{3pt}\vspace{-0.5pt}
Suppose\sss that\sss $i\fff g$\sss is\dss positive definite.\oss
Let\sss us consider\sss the operator\sss $A\qff +\qff \varepsilon$\sss
with\sss the same boundary conditions\sss $\Pi$\nnsp.\oss
Suppose\sss that\sss
$(\trf A\qff +\qff \varepsilon \trf)\trf(\trf s_{\dff +}\fff,\qff s_{\dff -}\trf)\off =\off 0$\nnsp.\oss
Then\sss
$s_{\dff +}\qff +\pff B\fff'\trf s_{\dff -}\off =\off 0$\sss
and\sss
$B\trf s_{\dff +}\qff -\qff s_{\dff -}\off =\off 0$\nnsp.\oss
Suppose also\sss that\sss 
$s\off =\off (\trf s_{\dff +}\fff,\qff s_{\dff -}\trf)$\sss satisfies\sss
the boundary condition\sss $\Pi$\nnsp,\oss
i.e.\qss $\gamma\dff s_{\dff -}\off =\off g\trf(\trf \gamma\dff s_{\dff +}\trf)$\nnsp,\oss
and apply\sss the\dss above\dss Lagrange\dss identity\sss for\sss $P$\sss to\sss 
$u
\off =\off
(\trf x\fff,\qff y\trf)
\off =\off
(\trf 0\fff,\qff s_{\dff -}\trf)$\sss
and\sss
$v
\off =\off
(\trf a\fff,\qff b\trf)
\off =\off
(\trf s_{\dff +}\fff,\qff 0\trf)$\nnsp.\oss
The\sss left\sss hand side of\dss the\dss Lagrange\dss identity\dss 
is\dss equal\dss to\vspace{3pt}
\[
\quad
\sco{\dff B\fff'\fff y\fff,\qff a \dff}_{\dff 0}
\qff +\qff
\sco{\dff B\trf x\fff,\qff b \dff}_{\dff 0}
\qff -\qff
\sco{\dff x\dff,\qff B\fff'\trf b \dff}_{\dff 0}
\qff -\qff
\sco{\dff y\dff,\qff B\trf a \dff}_{\dff 0}\off
\]

\vspace{-33pt}
\[
\quad
=\off
\sco{\dff B\fff'\fff s_{\dff -}\fff,\qff s_{\dff +} \dff}_{\dff 0}
\qff +\qff
\sco{\dff B\dff 0\fff,\qff 0 \dff}_{\dff 0}
\qff -\qff
\sco{\dff 0\dff,\qff B\trf 0 \dff}_{\dff 0}
\qff -\qff
\sco{\dff s_{\dff -}\dff,\qff B\dff s_{\dff +} \dff}_{\dff 0}
\]

\vspace{-33pt}
\[
\quad
=\off
\sco{\dff B\fff'\dff s_{\dff -}\fff,\qff s_{\dff +} \dff}_{\dff 0}
\qff -\qff
\sco{\dff s_{\dff -}\dff,\qff B\trf s_{\dff +} \dff}_{\dff 0}
\off =\off
-\qff
\sco{\dff s_{\dff +}\fff,\qff s_{\dff +} \dff}_{\dff 0}
\qff -\qff
\sco{\dff s_{\dff -}\dff,\qff s_{\dff -} \dff}_{\dff 0}
\off =\off
-\qff
\sco{\dff s\dff,\qff s \dff}_{\dff 0}
\off,
\]

\vspace{-12pt}\vspace{3pt}
and\dss the right\dss hand side\dss is\dss
equal\dss to\vspace{3pt}
\[
\quad
i\trf 
\sco{\dff \gamma\dff y\dff,\qff 
\gamma\dff a \dff}_{\dff \partial}
\pff +\pff
i\trf \bsco{\dff \gamma\dff x\dff,\qff 
\gamma\dff b \dff}_{\dff \partial}
\off =\off
i\trf \sco{\dff \gamma\dff s_{\dff -}\dff,\qff 
\gamma\dff s_{\dff +} \dff}_{\dff \partial}
\pff +\pff
i\trf \sco{\dff \gamma\trf 0\dff,\qff 
\gamma\trf 0 \dff}_{\dff \partial}
\]

\vspace{-33pt}
\[
\quad
=\off
\sco{\dff i\dff g\trf(\trf \gamma\dff s_{\dff +}\trf)\dff,\qff 
\gamma\dff s_{\dff +} \dff}_{\dff \partial}
\off.
\]

\vspace{-12pt}\vspace{3pt}
Therefore\sss the\dss Lagrange\dss identity\dss takes\sss the form\vspace{3pt}
\[
\quad
-\qff
\sco{\dff s\dff,\qff s \dff}_{\dff 0}
\off =\off
\sco{\dff i\dff g\trf(\trf \gamma\dff s_{\dff +}\trf)\dff,\qff 
\gamma\dff s_{\dff +} \dff}_{\dff \partial}
\off.
\]

\vspace{-12pt}\vspace{3pt}
The\sss left\sss hand side\dss is\sss $\leq\qff 0$\sss
and,\oss since\sss $i\dff g$\sss is\dss positive definite,\oss 
the right\sss hand\sss side\dss is\dss $\geq\qff 0$\nnsp.\oss
Hence\sss
$\bsco{\dff s\dff,\qff s \dff}_{\dff 0}
\off =\off
0$\sss
and\sss $s\off =\off 0$\nnsp.\oss
It\sss follows\sss that\sss 
$\kernel\dff (\trf A\qff +\qff \varepsilon\trf)\dff \oplus\dff \Gamma$\sss
is\dss equal\dss to $0$
and\dss hence\sss
$\kernel\dff (\trf A\qff +\qff \varepsilon\trf)_{\trf \Gamma}$
is\dss also equal\dss to $0$\nnsp.\oss

By\sss the results of\trs Section\qss \ref{abstract-index}\qss the operator\sss
$A_{\trf \Gamma}$\sss is\dss self-adjoint\sss and\trs Fredholm.\oss
Therefore\sss $A_{\trf \Gamma}$\sss is\dss a\dss Fredholm\dss operator of\dss index\sss $0$\nnsp.\oss
Since\sss the inclusion\sss 
$H_{\dff 1}\qff \ttoo\qff H_{\trf 0}$\sss is\dss compact,\oss
the operator\sss $(\trf A\qff +\qff \varepsilon\trf)_{\trf \Gamma}$\sss
is\dss a compact\sss perturbation of\dss $A_{\trf \Gamma}$\nsp.\oss
It\sss follows\sss that\sss $(\trf A\qff +\qff \varepsilon\trf)_{\trf \Gamma}$\sss
is\dss also a\dss Fredholm\dss operator of\dss index $0$\nnsp.\oss
Since\sss its kernel\dss is\sss $0$\nnsp,\oss it\dss is\dss surjective.\oss
Hence $(\trf A\qff +\qff \varepsilon\trf)_{\trf \Gamma}$\sss 
is\dss an\sss isomorphism\sss $\kernel\fff \Gamma\qff \ttoo\qff H_{\trf 0}$\nsp.\oss
This proves\sss the\sss theorem\sss for positive definite\sss $i\fff g$\nnsp.\oss
The proof\dss for\sss negative definite\sss $i\fff g$\sss is\dss
completely similar.\oss  \eproof

\mypar{Theorem.}{odd-families}
\emph{Let\sss $Z$\sss be a\sss topological\sss space.\oss
Suppose\sss that\sss for every\sss $z\qff \in\qff Z$\sss we are given operators\sss
$A\off =\off A\trf(\trf z\trf)$\sss 
and\sss $g\off =\off g\trf(\trf z\trf)$
satisfying\sss the above assumptions and continuously depending on\sss $z$\nnsp.\oss
Then\sss the analytical\dss index of\dss the family\sss
$A\trf(\trf z\trf)_{\trf \Gamma\trf(\trf z\trf)}\dff,\off z\qff \in\qff Z$\sss
is\dss equal\dss to zero.\oss}

\proof
Theorem\qss \ref{abstract-vanishing}\qss implies\sss that\sss the analytical\dss index
of\dss the family\sss
$(\trf A\trf(\trf z\trf)\qff +\qff \varepsilon\qff)_{\trf \Gamma\trf(\trf z\trf)}\dff,\off 
z\qff \in\qff Z$\sss
is\dss equal\dss to zero.\oss
Since\sss
$A\trf(\trf z\trf)_{\trf \Gamma\trf(\trf z\trf)}$\sss
is\dss a compact\sss perturbation of\dss
$(\trf A\trf(\trf z\trf)\qff +\qff \varepsilon\qff)_{\trf \Gamma\trf(\trf z\trf)}$\sss
for every $z$\nnsp,\oss the\sss theorem\sss follows.\oss  \eproof

\myuppar{Odd operators in\sss the\dss H\"{o}rmander\dss class.}
In\sss the framework of\trs Section\qss \ref{pdo},\oss
suppose\sss that\sss $E$\sss has\sss the form\sss
$E\off =\off F\dff \oplus\dff F$\sss for some\dss  
bundle $F$\dnsp.\oss
Let\sss $P$\sss be a  self-adjoint\sss elliptic operator on sections of\dss $E$\sss 
belonging\sss to\sss the\dss H\"{o}rmander\dss class.\oss
Let\sss $\sigma$\sss be\sss the symbol\sss of\dss $P$\sss and\sss $\bm{\Sigma}\trf(\trf x_{\dff n}\trf)$\sss
be\sss the corresponding endomorphism of\dss $E$\sss over\sss the collar.\oss
Suppose\sss that\sss $P$\sss  
has\sss the form\vspace{1.5pt}\vspace{-0.25pt}
\[
\quad
P
\off =\off\dff
\begin{pmatrix}
\off 0 \dff &
B\fff' \off
\vspace{4.5pt} \\
\off\dff B &
0 \off 
\end{pmatrix}
\qff,\quad
\mbox{and\dss that}\quad
\Sigma
\off =\off
\bm{\Sigma}\trf(\trf 0\trf)
\off =\off\dff
\begin{pmatrix}
\off 0 &
1 \qff\off
\vspace{4.5pt} \\
\off\dff 1 &
0 \qff\off 
\end{pmatrix}
\off
\]

\vspace{-12pt}\vspace{1.5pt}\vspace{-0.25pt}
with respect\sss to\sss the decomposition\sss 
$E\off =\off F\dff \oplus\dff F$\dnsp.\oss
Let\sss 
$f\dff \colon\dff
F\trf |\trf Y
\qff \ttoo\qff
F\trf |\trf Y$\sss
be a bundle map,\oss
and\dss let\sss 
$B_{\trf Y}\dff \colon\dff
(\trf F\dff \oplus\dff F\trf)\trf |\trf Y
\qff \ttoo\qff
F\trf |\trf Y$\sss
be defined\sss by\sss 
$B_{\trf Y}\dff(\trf a\fff,\qff b\trf)
\off =\off
b\qff -\qff f\trf(\trf a\trf)$\nnsp.\oss
Then\sss $B_{\trf Y}\dff \circ\dff \gamma$\sss
is\dss a bundle-like boundary operator having\sss 
$N
\off =\off
\{\qff (\trf a\fff,\qff f\trf(\trf a\trf)\trf)\qff |\qff a\qff \in\qff F \qff\}$\sss
as its kernel-symbol.\oss
Since\vspace{3pt}
\[
\quad
\sco{\dff
\Sigma\trf(\trf a\fff,\qff f\trf(\trf a\trf)\qff)\fff,\qff
(\trf b\fff,\qff f\trf(\trf b\trf)\qff) \dff}
\off =\off
\sco{\dff
(\trf f\trf(\trf a\trf)\fff,\qff a\qff)\fff,\qff
(\trf b\fff,\qff f\trf(\trf b\trf)\qff) \dff}
\off =\off
\sco{\dff
f\trf(\trf a\trf)\fff,\qff
b \dff}
\qff +\qff
\sco{\dff
a\fff,\qff
f\trf(\trf b\trf)\dff}
\qff,
\]

\vspace{-12pt}\vspace{3pt}
the kernel-symbol\sss $N$\sss is\dss self-adjoint\dss 
if\dss and\dss only\trs if\dss $f$\sss is\dss skew-adjoint,\oss
or,\oss equivalently,\pss $i\fff f$\sss is\dss self-adjoint.\oss
In\sss general,\oss such a\sss kernel-symbol\sss $N$\sss is\dss not\sss
an elliptic boundary condition\dss for\sss $\sigma$\nnsp.\oss
Let\sss $\tau_u\dff,\off \rho_{\fff u}$\sss be as in\dss Section\qss \ref{pdo}.\oss
Since\sss $P$\sss is\dss odd,\oss the operators $\tau_u\dff,\off \rho_{\fff u}$ 
have\sss the form\vspace{1.5pt}\vspace{-0.25pt}
\[
\quad
\tau_u
\off =\off\dff
\begin{pmatrix}
\off 0  &
\bm{\tau}_u^{\dff *} \dff\off
\vspace{4.5pt} \\
\off \bm{\tau}_u &
0 \dff\off 
\end{pmatrix}
\qff,\qquad
\rho_{\fff u}
\off =\off\dff
\begin{pmatrix}
\off \bm{\tau}_u  &
0 \off
\vspace{4.5pt} \\
\off 0 &
\bm{\tau}_u^{\dff *} \off 
\end{pmatrix}
\off
\]

\vspace{-12pt}\vspace{1.5pt}\vspace{-0.25pt}
for some operators\sss
$\bm{\tau}_u\dff \colon\dff F_{\dff y}\qff \ttoo\qff F_{\dff y}$\nnsp.\oss
We will\sss say\sss that $f$ is\qss \emph{equivariant}\pss if\dss endomorphisms\sss
$f_{\dff y}\dff \colon\dff F_{\dff y}\qff \ttoo\qff F_{\dff y}$
induced\dss by\sss $f$\sss commute with\sss $\bm{\tau}_u$\sss
when $y\off =\off \pi\trf(\trf u\trf)$\nnsp.\oss

\mypar{Lemma.}{graded-sh-l}
\emph{Suppose\sss that $f$ is\dss skew-adjoint\sss
and\sss either\sss $i f$\sss is\dss positive or 
negative definite,\oss
or\sss $f$\dss is\dss equivariant.\oss 
Then\sss $N$\sss is\dss a special\dss boundary condition\trs for\dss $\sigma$\nnsp.\oss}

\proof
Cf.\qss  \cite{gl},\oss Proposition\qss 3.1(a),\oss
and\qss \cite{blz},\oss Proposition\qss 4.3.\oss
Suppose\sss that\sss a vector\dss 
$(\trf a\fff,\qff b\trf)\qff \in\qff F\dff \oplus\dff F$\sss
belongs\sss to\sss $\mathcal{L}_{\dff +}\dff(\trf \rho_{\fff u}\dff)$\nnsp.\oss
Since\sss $\rho_{\fff u}$\sss has\sss the above diagonal\dss form,\oss
in\sss this case also\sss $(\trf a\fff,\qff 0\trf)$ and\sss $(\trf 0\fff,\qff b\trf)$
belong\sss to $\mathcal{L}_{\dff +}\dff(\trf \rho_{\fff u}\dff)$\nnsp.\oss
By\dss Corollary\qss \ref{plus-minus-orth}\qss the vectors\sss
$(\trf a\fff,\qff 0\trf)$ and\sss $(\trf 0\fff,\qff b\trf)$
are orthogonal\sss with respect\sss to\sss the scalar product\sss
$[\trf \bullet\fff,\qff \bullet\trf]
\off =\off
\sco{\dff \bm{\sigma}\dff \bullet\fff,\qff \bullet \dff}$\nnsp.\oss
It\sss follows\sss that\sss
$\sco{\dff b\fff,\qff a \dff}\off =\off 0$\nnsp.\oss
If\dss also
$(\trf a\fff,\qff b\trf)\qff \in\qff N_{\fff u}$\nsp,\oss
then\sss $b\off =\off f_{\dff y}\dff(\trf a\trf)$\sss
and\dss hence
\[
\quad
0
\off =\off
i\trf \sco{\dff b\fff,\qff a \dff}
\off =\off
i\trf \sco{\dff f_{\dff y}\dff(\trf a\trf)\fff,\qff a \dff}
\off =\off
\sco{\dff i\fff f_{\dff y}\dff(\trf a\trf)\fff,\qff a \dff}
\off.
\]

\vspace{-12pt}\vspace{1.5pt}
If\dss $i f$\sss is\dss either positive or 
negative definite,\oss this implies\sss that\sss $a\off =\off b\off =\off 0$\nnsp.\oss
Therefore in\sss this case\sss
$\mathcal{L}_{\dff +}\dff(\trf \rho_{\fff u}\dff)
\qff \cap\qff
N_{\fff u}
\off =\off 0$\nnsp,\oss
and,\oss similarly,\pss
$\mathcal{L}_{\dff -}\dff(\trf \rho_{\fff u}\dff)
\qff \cap\qff
N_{\fff u}
\off =\off 0$\nnsp.\oss
This proves\sss the\sss lemma when\sss $i f$\sss 
is\dss either positive or negative definite.\oss

Suppose now\sss that\sss $f$\sss is\dss equivariant.\oss
Since\sss $f$\sss is\dss skew-adjoint,\oss
the bundle\sss $F$\sss is\dss equal\sss to\sss the direct\sss sum
of\dss two subbundles with fibers\sss
$\mathcal{L}_{\dff +}\dff(\trf f_{\dff y}\dff)$\sss
and\sss
$\mathcal{L}_{\dff -}\dff(\trf f_{\dff y}\dff)$\sss
respectively.\oss
Since $f_{\dff y}$ commutes with\sss $\bm{\tau}_u$\sss
when\sss $y\off =\off \pi\trf(\trf u\trf)$\nnsp,\oss
both\sss these subbundles are\sss invariant\sss
under operators\sss $\bm{\tau}_u$\nsp.\oss
It\sss follows\sss that\sss $\sigma_y\dff,\off \tau_u$\sss
and\sss $N_{\dff u}$\sss
are also direct\sss sums,\oss and\sss that\sss the already proved\sss part\sss
of\dss the\sss lemma applies\sss to both summands.\oss 
The\sss lemma\sss follows.\oss  \eproof

\mypar{Theorem.}{vanishing}
\emph{Suppose\sss that\sss the above framework depends on a parameter\sss $z\qff \in\qff Z$\nnsp.\oss
If\dss the automorphism\sss $i f\trf(\trf z\trf)$\sss 
is\dss positive or negative definite for every\sss $z\qff \in\qff Z$\nnsp,\oss
then\sss the analytical\dss index of\dss the family\sss 
$P\trf(\trf z\trf)\fff,\qff z\qff \in\qff Z$\dss
with\sss the boundary conditions $N\trf(\trf z\trf)\fff,\qff z\qff \in\qff Z$\dss
is\dss equal\dss to zero.\oss}

\proof
Let\sss
$K^{\dff \partial}
\off =\off
H_{\trf 0}\dff(\trf Y\fff,\qff F\trf |\trf Y\trf)$\sss
and\sss
$g\trf(\trf z\trf)\dff \colon\dff K^{\dff \partial}\qff \ttoo\qff K^{\dff \partial}$\sss
be\sss the operator\sss induced\dss by\sss $f\trf(\trf z\trf)$\nnsp.\oss
Clearly,\oss if\dss $i f\trf(\trf z\trf)$\sss 
is\dss positive or negative definite,\oss
then\sss $i\fff g\trf(\trf z\trf)$\sss is\dss also positive or negative definite.\oss 
Therefore\sss the\sss theorem\sss follows from\dss Theorem\qss \ref{odd-families}.\oss  \eproof

\myuppar{Another\sss form of\dss $P\fff,\off N$\nnsp.}
Let\sss us consider\sss the subbundles
\vspace{3pt}
\[
\quad
\Delta_{\trf F}
\off =\off 
\{\qff (\trf a\fff,\qff a\trf)\qff |\qff a\qff \in\qff F \qff\}
\quad
\mbox{and}\quad
\Delta_{\trf F}^{\fff \perp}
\off =\off 
\{\qff (\trf b\fff,\qff -\qff b\trf)\qff |\qff b\qff \in\qff F \qff\}
\]

\vspace{-12pt}\vspace{3pt}
of\dss $E\off =\off F\dff \oplus\dff F$\dnsp.\oss
Clearly,\pss
$\Delta_{\trf F}\off =\off E^{\dff +}$\dnsp,\qss
$\Delta_{\trf F}^{\fff \perp}\off =\off E^{\dff -}$\dnsp,\oss
and we can\sss identify\sss both\sss $\Delta_{\trf F}$\sss and\sss $\Delta_{\trf F}^{\fff \perp}$\sss
with\sss $F$\sss by\sss the projection\sss
$(\trf x\fff,\qff y\trf)
\off \longmapsto\off
x$\nnsp.\oss
Suppose\sss that\sss $f$\sss is\dss skew-adjoint\sss
and\sss write down\sss the corresponding\sss boundary condition\sss 
$N$\sss in\sss terms of\dss the decomposition\sss 
$E\off =\off \Delta_{\trf F}\dff \oplus\dff \Delta_{\trf F}^{\fff \perp}
\off =\off
F\dff \oplus\dff F$\dnsp.\oss
An easy calculation shows\sss that\sss for every\sss
$u\qff \in\qff S\fff Y$\sss and\sss $y\off =\off \pi\trf(\trf u\trf)$\vspace{3pt}
\[
\quad
N_{\dff u}
\off =\off
N_{\dff y}
\off =\off
\left\{\off \left.
\left(\qff 
a\fff,\off \frac{1\qff -\qff f_{\dff y}}{1\qff +\qff f_{\dff y}} \qff(\trf a\trf)
\qff\right)
\off \right|\off
a\qff \in\qff F_{\dff y}
\off\right\}
\qff.\quad
\]

\vspace{-12pt}\vspace{3pt}
Note\sss that\sss $1\qff +\qff f_{\dff y}$\sss is\dss invertible because $f$\sss is\dss skew-adjoint.\oss
By\sss the same reason\vspace{2.75pt}
\[
\quad
\psi_{\fff y}
\off =\off
\frac{1\qff -\qff f_{\dff y}}{1\qff +\qff f_{\dff y}}
\]

\vspace{-12pt}\vspace{2.75pt}
is\dss unitary and\sss invertible.\oss
If\dss we consider\sss $\psi_{\fff y}$\sss as an operator\sss
$E^{\dff +}_{\dff y}\qff \ttoo\qff E^{\dff -}_{\dff y}$\nsp,\oss
then\sss $N_{\dff y}$\sss is\dss the graph of\dss $\psi_{\fff y}$\nsp.\oss 
If\dss $f$\sss is\dss equivariant,\oss
then\sss $N$\sss is\dss a special\dss boundary condition\sss
by\trs Lemma\qss \ref{graded-sh-l}.\oss

\myuppar{Dirac-like boundary problems.}
Suppose\sss that\sss $f$\sss is\dss skew-adjoint\sss and equivariant\sss
and\sss that\sss the operators\sss
$\bm{\tau}_u\dff \colon\dff F\qff \ttoo\qff F$\dss
are skew-adjoint\sss for every\sss $u\qff \in\qff S\fff Y$\dnsp.\oss
In\sss this case\sss the symbol $\sigma$ 
is\dss anti-commuting,\pss $N$\sss is\dss a special\dss boundary condition,\oss 
and\dss hence we are in\sss the situation of\trs
Section\qss \ref{special-conditions}.\oss
One may call\dss the boundary problem defined\dss by such\sss 
$P\fff,\qff f$\qss \emph{Dirac-like}.\oss

Since $f$ is\dss equivariant,\oss
the operators\sss $\bm{\tau}_u$\sss
leave\sss invariant\sss 
$\mathcal{L}_{\dff +}\dff(\trf f_{\dff y} \trf)$\sss
and\sss
$\mathcal{L}_{\dff -}\dff(\trf f_{\dff y} \trf)$\nnsp.\oss
Let\sss\vspace{3pt}
\[
\quad
\bm{\tau}_u^{\dff +}\qff \colon\qff
\mathcal{L}_{\dff +}\dff(\trf f_{\dff y} \trf)
\qff \ttoo\qff
\mathcal{L}_{\dff +}\dff(\trf f_{\dff y} \trf)
\quad
\mbox{and}\quad
\bm{\tau}_u^{\dff -}\qff \colon\qff
\mathcal{L}_{\dff -}\dff(\trf f_{\dff y} \trf)
\qff \ttoo\qff
\mathcal{L}_{\dff -}\dff(\trf f_{\dff y} \trf)
\]

\vspace{-12pt}\vspace{3pt}
be\sss the operators induced\dss by\sss $\bm{\tau}_u$\nsp.\oss
Let\sss
$\mathcal{L}_{\dff +}\dff(\trf f \trf)$\sss
and\sss
$\mathcal{L}_{\dff -}\dff(\trf f \trf)$\sss
be\sss the bundles having\sss
$\mathcal{L}_{\dff +}\dff(\trf f_{\dff y} \trf)$\sss
and\sss
$\mathcal{L}_{\dff -}\dff(\trf f_{\dff y} \trf)$\nnsp,\oss
respectively,\oss
as\sss the fibers over\sss $y\qff \in\qff Y$\dnsp.\oss
Then\vspace{3pt}
\[
\quad
F
\off =\off
\mathcal{L}_{\dff +}\dff(\trf f \trf)
\qff \oplus\qff
\mathcal{L}_{\dff -}\dff(\trf f \trf)
\]

\vspace{-12pt}\vspace{3pt}
and\sss the operators\sss 
$\bm{\tau}_u^{\dff +}\dff,\off \bm{\tau}_u^{\dff -}$\sss
define skew-adjoint\sss symbols 
$\bm{\tau}^{\dff +}\dff,\off \bm{\tau}^{\dff -}$
and self-adjoint\sss symbols\sss
$i\dff \bm{\tau}^{\dff +}\dff,\off i\dff \bm{\tau}^{\dff -}$ 
over $Y$ in\sss the bundles\sss
$\mathcal{L}_{\dff +}\dff(\trf f \trf)\dff,\off 
\mathcal{L}_{\dff -}\dff(\trf f \trf)$\sss
respectively.\oss

\mypar{Theorem.}{dirac-flip}
$\ti\trf(\trf \sigma\fff,\off N\trf)
\off =\off\dff
\ti\trf(\qff i\dff \bm{\tau^{\dff -}}\qff)
\off =\off
\ti\trf(\qff -\qff i\dff \bm{\tau^{\dff +}}\qff)$\nnsp.\oss

\proof
After a spectral\sss deformation of\sss $f$ we may assume\sss that\sss
$i\fff,\qff -\qff i$\sss are\sss the only eigenvalues 
of\dss $f_{\dff y}$\sss for every $y$\nnsp.\oss
Then\sss $i\fff,\qff -\qff i$\sss are also\sss 
the only eigenvalues of\dss $\psi_{\fff y}$ for every $y$\nnsp,\oss
and\vspace{3pt}
\[
\quad
\mathcal{L}_{\dff +}\dff(\trf \psi_{\fff y} \trf)
\off =\off
\mathcal{L}_{\dff -}\dff(\trf f_{\dff y} \trf)
\dff,\quad
\mathcal{L}_{\dff -}\dff(\trf \psi_{\fff y} \trf)
\off =\off
\mathcal{L}_{\dff +}\dff(\trf f_{\dff y} \trf)
\qff.
\]

\vspace{-12pt}\vspace{3pt}
Let\sss 
$f^{\trf +}\dff \colon\dff
F\trf |\trf Y
\qff \ttoo\qff
F\trf |\trf Y$\sss 
be\sss the multiplication\sss by\sss $i$\nnsp,\oss
and\sss $N^{\dff +}$\sss be\sss the corresponding boundary condition.\oss
Since $i f^{\trf +}$ is\dss negative definite,\pss
$\ai\trf(\trf \sigma\fff,\off  N^{\dff +}\trf)\off =\off 0$\sss by\trs
Theorem\qss \ref{vanishing}.\oss
Clearly,\oss the bundle\sss $E^{\dff +}_{\trf Y}$\sss is\dss isomorphic\sss
to\sss $F\trf |\trf Y$\sss and\dss hence extends\sss to $X$\nnsp.\oss
Hence\sss 
$\ti\trf(\trf \sigma\fff,\off  N^{\dff +}\trf)\off =\off 0$\sss by\trs
Theorem\qss \ref{index-theorem-ext},\oss
and\sss in order\sss to prove\sss the first\sss equality\sss
it\dss is\dss sufficient\sss to prove\sss that\vspace{3pt}
\begin{equation}
\label{dirac-difference}
\quad
\varepsilon^{\dff +}\dff \bigl(\trf \sigma\fff,\off N\qff\bigr)
\qff -\qff
\varepsilon^{\dff +}\dff \bigl(\trf \sigma\fff,\off N^{\dff +}\qff\bigr)
\off =\off\dff
\beta\dff \left(\trf
e^{\dff +}\fff (\trf i\dff \bm{\tau^{\dff -}}\trf)
\qff\right)
\qff.
\end{equation}

\vspace{-12pt}\vspace{3pt}
The boundary conditions\sss $N$\sss and\sss $N^{\dff +}$\sss
differ only\sss in\sss the summand\sss
$\mathcal{L}_{\dff -}\dff(\trf f \trf)$\nnsp,\oss
and\sss in\sss this summand\sss they differ in\sss the same way
as\sss the boundary conditions\sss $N$\sss and\sss $N^{\dff \perp}$\sss
in\dss Theorem\qss \ref{basic-difference}\qss and\dss
Corollary\qss \ref{k1-difference}.\oss
In\sss this summand\dss the proof\dss of\trs Theorem\qss \ref{basic-difference}\qss applies
and\sss will\dss imply\qss (\ref{dirac-difference})\qss
once we\sss prove\sss that\sss the sum of\dss the positive eigenspaces of\dss the operators\sss
$i\dff \bm{\tau_u^{\dff -}}$\sss are equal\dss to\sss the positive eigenspaces
of\dss the operators\sss $\bm{\upsilon}_{\fff u}$\sss from\dss Section\qss \ref{special-conditions}\qss
for\sss the summand\sss $\mathcal{L}_{\dff -}\dff(\trf f \trf)$\nnsp.\oss

Since in\dss Section\qss \ref{special-conditions}\qss we worked\sss with\qss 
``non-graded''\qss operators,\oss
we need\sss to pass\sss to\sss the decomposition\sss
$E\off =\off \Delta_{\trf F}\dff \oplus\dff \Delta_{\trf F}^{\fff \perp}
\off =\off
F\dff \oplus\dff F$\dnsp.\oss
We already\sss computed\sss the isometries\sss $\psi_{\fff y}$\sss in\sss this decomposition,\oss
and on\sss $\mathcal{L}_{\dff -}\dff(\trf f \trf)$\sss 
they are equal\dss
to\sss the multiplication\sss by\sss $i$\nnsp.\oss 
At\sss the same\sss time\vspace{1.5pt}
\[
\quad
\rho_{\fff u}
\off =\off\dff
\begin{pmatrix}
\off \bm{\tau}_u  &
0 \off
\vspace{4.5pt} \\
\off 0 &
-\qff \bm{\tau}_u \off 
\end{pmatrix}
\off
\]

\vspace{-12pt}\vspace{1.5pt}
with respect\sss to\sss the original\sss decomposition\sss
$E\off =\off F\dff \oplus\dff F$\dnsp,\oss
and\dss hence\vspace{3pt}
\[
\quad
\mathcal{L}_{\dff +}\dff(\trf \rho_{\fff u}\dff)
\off =\off
\mathcal{L}_{\dff +}\dff(\trf \bm{\tau}_u\dff)
\dff \oplus\trf
\mathcal{L}_{\dff -}\dff(\trf \bm{\tau}_u\dff)
\off \subset\off\dff
F\dff \oplus\dff F
\]

\vspace{-12pt}\vspace{3pt}
with\sss respect\sss to\sss the original\sss decomposition.\oss
It\sss follows\sss that\sss with respect\sss to\sss the decomposition\sss 
$E\off =\off \Delta_{\trf F}\dff \oplus\dff \Delta_{\trf F}^{\fff \perp}
\off =\off
F\dff \oplus\dff F$\sss
the subspace\sss
$\mathcal{L}_{\dff +}\dff(\trf \rho_{\fff u}\dff)$\sss
is\dss the graph of\dss the isometry\sss $\varphi_{\fff u}$\sss
such\sss that\vspace{3pt}
\[
\quad
\varphi_{\fff u}\dff(\trf a\trf)\off =\off a
\quad
\mbox{if}\quad
a\qff \in\qff \mathcal{L}_{\dff +}\dff(\trf \bm{\tau}_u\dff)
\quad
\mbox{and}\quad
\varphi_{\fff u}\dff(\trf a\trf)\off =\off -\qff a
\quad
\mbox{if}\quad
a\qff \in\qff \mathcal{L}_{\dff -}\dff(\trf \bm{\tau}_u\dff)
\qff.
\]

\vspace{-12pt}\vspace{3pt}
It\sss follows\sss that\sss in\dss the summand\sss
$\mathcal{L}_{\dff -}\dff(\trf f \trf)$\vspace{1.5pt}
\[
\quad
\mathcal{L}_{\dff +}\qff\left(\trf \psi_y^{\dff -\dff 1}\dff \circ\trf \varphi_{\fff u} \trf\right)
\off =\off
\mathcal{L}_{\dff -}\dff\left(\trf \bm{\tau}_u^{\dff -}\qff\right)
\qff.
\]

\vspace{-12pt}\vspace{1.5pt}
But\sss $\mathcal{L}_{\dff -}\dff(\trf \bm{\tau}_u^{\dff -}\qff)$\sss
is\dss exactly\sss the sum of\dss the positive eigenspaces of\dss $i\dff \bm{\tau}_u^{\dff -}$\nsp.\oss
By\sss the previous paragraph,\oss this implies\sss the equality\qss (\ref{dirac-difference}),\oss
and\dss hence\sss the equality\sss
$\ti\trf(\trf \sigma\fff,\off N\trf)
\off =\off\dff
\ti\trf(\qff i\dff \bm{\tau^{\dff -}}\qff)$\nnsp.\oss
The proof\dss of\dss the equality\sss
$\ti\trf(\trf \sigma\fff,\off N\trf)
\off =\off\dff
\ti\trf(\qff -\qff i\dff \bm{\tau^{\dff +}}\qff)$\sss
is\dss completely similar.\oss  \eproof

\mypar{Theorem.}{dirac-flip-analytical}
$\ai\trf(\trf \sigma\fff,\off N\trf)
\off =\off\dff
\ai\trf(\qff i\dff \bm{\tau^{\dff -}}\qff)
\off =\off
\ai\trf(\qff -\qff i\dff \bm{\tau^{\dff +}}\qff)$\nnsp.\oss

\proof
It\dss is\dss sufficient\sss to combine\dss Theorem\qss \ref{dirac-flip}\qss
with\dss Theorem\qss \ref{index-theorem-ext}.\oss  \eproof

\vspace*{\bigskipamount}

\begin{flushright}

First\sss version\qss --\qss July\dss 17\fff,\oss 2022

Second\sss version\qss --\qss December\dss 30\fff,\oss 2022

Present\sss version\qss --\qss June\dss 24\fff,\oss 2023

https\halfff:/\!/\!nikolaivivanov.com

E-mail\halfff:\oss nikolai.v.ivanov{\fff}@{\dff}icloud.com,\oss ivanov{\fff}@{\dff}msu.edu

Department\sss of\qss Mathematics,\oss Michigan\sss State\sss University

\end{flushright}

\end{document}